\documentclass[10pt]{amsart}
\usepackage{amssymb,amsmath,amsfonts,amsthm,graphics,mathrsfs,amscd,hyperref}
\usepackage[makeroom]{cancel}
\usepackage[hmargin=1in,vmargin=1in]{geometry}
\usepackage[all,cmtip]{xy}

\theoremstyle{plain}

\newtheorem{definition}[equation]{Definition}

\newtheorem{lemma}[equation]{Lemma}
\newtheorem{proposition}[equation]{Proposition}
\newtheorem{theorem}[equation]{Theorem}

\newtheorem{notation}{Notation}
\newtheorem{remark}[equation]{Remark}

\numberwithin{equation}{subsection}

\title[Theorems of existence and completeness]{On deformations of foliated complex analytic structures I : Theorems of existence and completeness}
\author{Chunghoon Kim}

\begin{document}

\maketitle
\begin{abstract}
We study deformations of foliated complex analytic structures defined by locally free subsheaves of tangent sheaves and locally free subsheaves of cotangent sheaves on the basis of Kodaira-Spencer's deformation theory. We prove theorems of existence and completeness for deformations of foliated complex analytic structures defined by locally free subsheaves of tangent sheaves as analogues of theorem of existence and completeness for deformations of complex analytic structures by Kodaira-Spencer.  We prove theorems of existence and completeness for deformations of foliated complex analytic structures defined by locally free subsheaves of cotangent sheaves. We also prove theorems of existence and completeness for simultaneous deformations  when singular holomorphic foliations on compact complex manifolds are defined by locally free subsheaves of tangent and contangent sheaves at the same time.
\end{abstract}

\tableofcontents

\section{Introduction}

We study deformations of foliated complex analytic structures defined by locally free subsheaves of tangent sheaves and locally free subsheaves of cotangent sheaves on the basis of Kodaira-Spencer's deformation theory.  There are three equivalent definitions of a (singular) holomorphic foliation $\mathcal{F}_0$ on a $n$-dimensional compact complex manifold $M$ \footnote{We adopt the notations from \cite{LPT13}}:
\begin{enumerate}
\item A (singular) holomorphic foliation $\mathcal{F}_0$ of dimension $p$ on $M$ is a coherent subsheaf $\Theta_{\mathcal{F}_0}$ of the tangent sheaf $\Theta_M$ such that $\left[\Theta_{\mathcal{F}_0}, \Theta_{\mathcal{F}_0}\right]\subset \Theta_{\mathcal{F}_0}$ under the Lie bracket, and $\frac{\Theta_M}{\Theta_{\mathcal{F}_0}} $ is a torsion-free $\mathcal{O}_M$-module. The locus of points where $\frac{\Theta_M}{\Theta_{\mathcal{F}_0}}$ is not locally free is called the singular locus $\textnormal{Sing}(\mathcal{F}_0)$ of $\mathcal{F}_0$. The dimension of $\mathcal{F}_0$ is the generic rank of $\Theta_{\mathcal{F}_0}$. We note that the codimension of $\textnormal{Sing}(\mathcal{F}_0)$ is at least 2. In this case, we will denote the foliated complex manifold by $\left(M, \Theta_{\mathcal{F}_0}\right)$. \label{definition1}
\item A (singular) holomorphic foliation $\mathcal{F}_0$ of codimension $q$ on $M$ is a coherent subsheaf $. \mathcal{N}_{\mathcal{F}_0}^*$ of the cotangent sheaf $\Omega_M^1$ such that  $\frac{\Omega_M^1}{\mathcal{N}_{\mathcal{F}_0}^*}$ is a torsion free $\mathcal{O}_M$-module (the locus of points where $\frac{\Omega_M^1}{\mathcal{N}_{\mathcal{F}_0^*}  }$ is not locally free is called the singular locus $\textnormal{Sing}(\mathcal{F}_0)$ of $\mathcal{F}_0)$, and $d \left( \mathcal{N}_{\mathcal{F}_0,x}^* \right) \subset \mathcal{N}_{\mathcal{F}_0,x}^* \bigwedge \Omega_{M,x}^1$ for $x\in M - \textnormal{Sing}(\mathcal{F}_0)$. The codimension of $\mathcal{F}_0$ is the generic rank of $\mathcal{N}_{\mathcal{F}_0}^*$. In this case, we will denote the foliated complex manifold by $\left(M, \mathcal{N}_{\mathcal{F}_0}^* \right)$. \label{definition2}
\item A (singular) holomorphic foliation $\mathcal{F}_0$ is a section $\omega_0 \in H^0\left(M, \bigwedge^q \Omega_M^1\otimes  \mathcal{L}_0\right)$ where $\mathcal{L}_0$ is a line bundle such that the zero set $\textnormal{Sing}(\mathcal{F}_0)$ of $\omega_0$ has codimension at least 2, and for any $x\in M- \textnormal{Sing}(\mathcal{F}_0)$, there is a neighborhood $U_x$ such that $\omega_0$ is locally written as the product of $q$-local holomorphic one-forms on $U_x$, i.e. $\omega= w_x^1\wedge \cdots \wedge w_x^q$ on $U_x$, and satisfies $w_x^1\wedge \cdots \wedge w_x^q\wedge dw_x^\alpha=0$ on $U_x$ for $\alpha=1,...,q$. In this case, we will denote the foliated complex manifold by $(M, \omega_0)$. \label{definition3}
\end{enumerate}
We shall study foliated deformations of $\left(M, \Theta_{\mathcal{F}_0}\right)$ when $\Theta_{\mathcal{F}_0}$ is locally free, and study foliated deformations of $(M, \mathcal{N}_{\mathcal{F}_0}^*)$ when $\mathcal{N}_{\mathcal{F}_0}^*$ is locally free, which is equivalent to foliated deformations of $(M, \omega_0)$ (see Part II).

The author believes that in this introduction it would be helpful to explain how this present work has been started and developed in order to understand the overall picture. While the author was working on the extension of Horikawa's work \cite{Hor76} on deformations of rational maps to deformations of a special kind of Poisson rational maps, he paid attention to the third definition (\ref{definition3}) of singular holomorphic foliations defined by sections, and observed a similarity between deformations of rational maps and deformations of singular holomorphic foliations in terms of sections. He thought that Horikawa's method might be applied to deformations of singular holomorphic foliations. Let us review Horikawa's deformation theory of rational maps and explain how his theory can be applied to foliated deformations in terms of sections. 

Let $M$ be a compact complex manifold and $g:M\to \mathbb{P}_\mathbb{C}^m$ be a rational map into the projective space $\mathbb{P}_\mathbb{C}^m$ which is determined by sections $g^0,..., g^m\in H^0\left(X, \mathcal{L}_0 \right)$ for a line bundle $\mathcal{L}_0$ on $M$. Then deformations of a rational map $g:M\to \mathbb{P}_\mathbb{C}^m$ consist of (i) deformations of the complex structures of $M$ and (ii) deformations of the line bundle $\mathcal{L}_0$ and (iii) deformations of sections $g^0,..., g^m \in H^0\left(M, \mathcal{L}_0 \right)$. As is known, deformations of the pair $\left(M, \mathcal{L}_0\right)$ is controlled by the Atiyah extension $\mathcal{E}_{\mathcal{L}_0}$ associated to the line bundle $\mathcal{L}_0$ defined in the following way: let $\mathcal{U}=\left\{U_i \right\}$ be an open covering of $M$ by coordinate neighborhoods such that the line bundle $\mathcal{L}_0$ is defined by the transition functions $\left\{\psi_{ij} \right\}\in C^1\left(\mathcal{U},  \mathcal{O}_M^* \right)$. Then $\mathcal{E}_L|_{U_i}$ is isomorphic to $\Theta_M|_{U_i} \bigoplus \mathcal{O}_M|_{U_i}$, and $\left(\rho_i, \alpha_i \right)\in \Theta_M|_{U_i}\bigoplus \mathcal{O}_M|_{U_i}$ equals $\left(\rho_j, \alpha_j \right)\in \Theta_M|_{U_j} \bigoplus \mathcal{O}_M|_{U_j}$ if and only if $\rho_i=\rho_j$ and $\alpha_j-\alpha_i=\left[ \rho_i, \log \psi_{ij}\right]$ on $U_i\cap U_j$. Then we have an exact sequence $0\to \mathcal{O}_X \to \mathcal{E}_{\mathcal{L}_0} \to \Theta_M \to 0$. Let $g^\beta \in H^0\left(X, \mathcal{L}_0 \right)$ be represented by locally $g_i^\beta\in \Gamma\left(U_i, \mathcal{O}_M \right)$ on $U_i$ with $g_i^\beta = \psi_{ij} g_j^\beta$ on $U_i\cap U_j$. Horikawa showed that deformations of a rational map $g:M\to \mathbb{P}_\mathbb{C}^m$ is controlled by the complex of sheaves
\begin{align}\label{I6}
\mathcal{E}_{\mathcal{L}_0} \xrightarrow{\tilde{G}} \bigoplus^{m+1} \mathcal{L} \to 0 \to 0 \to \cdots
\end{align}
where $\tilde{G}$ is defined by $\tilde{G}(\rho_i, \alpha_i)=\left( \left[ \rho_i, g_i^0 \right] +g_i^0\alpha_i,..., \left[ \rho_i, g_i^\beta \right] + g_i^\beta \alpha_i ,..., \left[\rho_i, g_i^m \right] + g_i^m \alpha_i \right)$ on $U_i$.

Considering the definition $(\ref{definition3})$ of a singular holomorphic foliation of codimension $q$ defined by a section $\omega_0 \in H^0\left( M, \bigwedge^q \Omega_M^1\otimes \mathcal{L}_0 \right)$ for a line bundle $\mathcal{L}_0$, which is locally decomposable and satisfies the integrability condition, the author thought that deformations of a singular holomorphic foliation $\omega_0$ on $M$ consist of (i) deformations of the complex structures of $M$ and (ii) deformations of the line bundle $\mathcal{L}_0$ and (iii) deformations of a section $\omega_0\in H^0\left(M, \bigwedge^q \Omega_M^1 \otimes \mathcal{L}_0 \right)$ which should be locally decomposable and satisfy the integrability condition. This situation is quietly similar to deformations of a rational map as above.  In our case, we have only one section. Let $\mathcal{U}=\{U_i\}$ be an open covering of $M$ by coordinate neighborhoods such that $\mathcal{L}_0$ is defined by $\{\psi_{ij}\}$ as above. Furthermore, the author assumed that $\omega_0$ is locally decomposable at every point of $M$. In other words, he assumed that $\omega_0$ is represented by $\omega_i= w_i^1 \wedge \cdots \wedge w_i^q\in \Gamma\left( U_i, \bigwedge^q \Omega_M^1 \right)$, where $w_i^\beta\in \Gamma\left(U_i, \Omega_M^1 \right),\beta=1,...q$ with $\omega_i= \psi_{ij} \omega_j$ and $w_i^1 \wedge \cdots \wedge w_i^q \wedge dw_i^\alpha=0,\alpha=1,...,q$  (This is possible when a holomorphic foliation $\mathcal{N}_{\mathcal{F}_0}^*$ is locally free. That's why we assume the local freeness of $\mathcal{N}_{\mathcal{F}_0}^*$ when we deform $\left(M, \mathcal{N}_{\mathcal{F}_0}^* \right)$ in terms of cotangent sheaves in this work). In a similar way to $(\ref{I6})$, the author could construct a complex of sheaves up to the first two-terms
\begin{align}\label{I9}
\mathcal{E}_{\mathcal{L}_0} \xrightarrow{\tilde{G}} \bigwedge^q \Omega_M^1\otimes \mathcal{L}_0 \xrightarrow{G_1\textbf{?}}  B_2 \textbf{?}   \xrightarrow{G_2\textbf{?}}  B_3\textbf{?}  \xrightarrow{G_3\textbf{?}} B_4\textbf{?}  \xrightarrow{G_4\textbf{?}} \cdots
\end{align}
where $\tilde{G}$ is defined by $\tilde{G}\left(\rho_i,\alpha_i \right)= \mathcal{L}_{\rho_i} (\omega_i) + \alpha_i \omega_i $ on $U_i$. Here $\mathcal{L}$ is the Lie derivative.  When we approach some deformation problem, the first thing to do is to identify infinitesimal deformations. An infinitesimal foliated deformation of $\left(M, \omega_0 \right)$ over the ring of the dual number $\mathbb{C}[\epsilon]=\mathbb{C}[t]/(t^2)$ consists of deformations of the pair $(M, \mathcal{L}_0)$ (for the detail, see \cite{Ser06} p.148) which is a collection of (i) $\textnormal{Id}+\epsilon p_{ij}:\Gamma\left( U_{ij}, \mathcal{O}_M\right)\to \Gamma(U_{ij}, \mathcal{O}_M)$ where $p_{ij}\in \Gamma\left(U_{ij}, \Theta_M \right)$ and (ii) transition functions $F_{ij}= \psi_{ij} + \epsilon g_{ij}$ for some $g_{ij}\in \Gamma(U_{ij}, \mathcal{O}_M)$, and deformations of the section $\omega_0$ which is represented by (iii) $\omega_i+ \epsilon v_i=\left(w_i^1+ \epsilon v_i^1\right)\wedge \cdots \wedge \left(w_i^q +\epsilon v_i^q\right)$ for some $v_i^\alpha \in \Gamma\left(U_i, \Omega_M^1\right), \alpha=1,...,p$ $(\iff v_i= \sum_{\alpha=1}^q w_i^1 \wedge \cdots \wedge v_i^\alpha \wedge \cdots \wedge w_i^q)$ with $\omega_i+ \epsilon v_i= \left( \psi_{ij} + \epsilon g_{ij} \right) \left( \omega_j + \epsilon v_j \right) $, and the integrability condition (iv) $\left(w_i^1+\epsilon v_i^1\right)\wedge \cdots \wedge \left(w_i^q+ \epsilon v_i^q \right)\wedge d\left( w_i^\alpha + \epsilon v_i^\alpha \right)=0,\alpha=1,...,q$ ($\iff  \sum_{\alpha=1}^q w_i^1\wedge \cdots \wedge v_i^\alpha \wedge \cdots \wedge  w_i^q \wedge dw_i^\alpha + w_i^1\wedge \cdots \wedge w_i^q \wedge dv_i^\alpha=0 $). Then the infinitesimal foliated deformation over $\mathbb{C}[\epsilon]$ is represented by $\left( \{p_{ij}\}, \left\{\frac{g_{ij}}{\psi_{ij}} \right\}      \right) \bigoplus \left\{ v_i=\sum_{\alpha=1}^q w_i^1 \wedge \cdots \wedge v_i^\alpha \wedge \cdots \wedge w_i^q   \right\}\in  C^1\left( \mathcal{U}, \mathcal{E}_{\mathcal{L}_0} \right)\bigoplus C^0\left(\mathcal{U}, \bigwedge^q \Omega_M^1\otimes \mathcal{L}_0 \right)$, and the condition (iv) gave a clue to define $G_1$ and the third term $B_2$ in $(\ref{I9})$. It turned out that the second term $\bigwedge^q \Omega_M^1\otimes \mathcal{L}_0$ in $(\ref{I9})$ is too big, so that the author defined the subsheaf $\mathcal{S}^1\subset \bigwedge^q \Omega_M^1\otimes \mathcal{L}_0$ locally generated by $w_i^1\wedge \cdots \wedge \widehat{w_i^\alpha} \wedge \cdots \wedge w_i^q ,\alpha=1,..., q$. Using the condition (iv), the author attempted to define $G_1$
\begin{align*}
v_i=\sum_{\alpha=1}^q w_i^1 \wedge \cdots \wedge v_i^\alpha \wedge \cdots \wedge w_i^q  \mapsto  \left( v_i\wedge dw_i^1 + \omega_i \wedge dv_i^1, ...,  v_i \wedge dw_i^q + \omega_i \wedge dv_i^q         \right) \in \bigoplus^q \Gamma\left(U_i, \bigwedge^{q+2}\Omega_M^1 \right)
\end{align*}
This led to define $B_2$:
\begin{align}\label{I8}
\mathcal{E}_{\mathcal{L}_0} \xrightarrow{\tilde{G}} \mathcal{S}^1 \xrightarrow{G_1}  \left( \mathcal{N}_{\mathcal{F}_0}^* \right)^* \otimes \bigwedge^{q+2} \Omega_M^1 \otimes \mathcal{L}_0 \xrightarrow{G_2 \textbf{?}} B_3 \textbf{?} \xrightarrow{G_3 \textbf{?}} B_4 \textbf{?} \xrightarrow{G_4 \textbf{?}} \cdots
\end{align}
We can study foliated deformations of $\left(M, \omega_0\right)$ with $\mathcal{N}_{\mathcal{F}_0}^*$ locally free by using the first three terms $\mathcal{E}_{\mathcal{L}_0} \to \tilde{S}^1\to \left( \mathcal{N}_{\mathcal{F}_0}^* \right)^* \otimes \bigwedge^{q+2} \Omega_M^1\otimes \mathcal{L}_0\to 0\to 0 \to \cdots$. Actually the complex of these three terms is sufficient to study foliated deformations of $\left(M, \omega_0\right)$ or $\left(M, \mathcal{N}_{\mathcal{F}_0}^* \right)$ with $\mathcal{N}_{\mathcal{F}_0}^*$ locally free for compact complex surfaces and compact complex threefolds. Examples we study in Part IV-2 are compact complex surfaces and threefolds. On the other hand, for general case, it seemed to the author that the third term $\left( \mathcal{N}_{\mathcal{F}_0}^* \right)^* \otimes \bigwedge^{q+2} \Omega_M^1\otimes \mathcal{L}_0$ is not only too big, but also there should be some restriction on the second cohomology group of this three-term complex (which represents obstructions) by extending the complex to $B_3, B_4,\cdots$. In other words, this three-term complex did not seem optimal for general case. However every attempt to extend the complex of sheaves $(\ref{I8})$ to $B_3,B_4, \cdots$ turned out to be unnatural. 

The most crucial paper to the author at that time was G\'omez-Mont's work \cite{GM88} on deformations of singular holomorphic foliations where in particular he identifies the first cohomology group of the leaf complex of a compact foliated complex manifold $(M, \Theta_{\mathcal{F}_0})$
\begin{align}\label{I1}
\Theta_{\mathcal{F}_0}^\bullet:\Theta_M  \xrightarrow{D_0} \mathscr{H}om_{\mathcal{O}_M}\left( \Theta_{\mathcal{F}_0}, \frac{\Theta_M}{\Theta_{\mathcal{F}_0}} \right) \xrightarrow{D_1} \mathscr{H}om_{\mathcal{O}_M}\left(\bigwedge^2 \Theta_{\mathcal{F}_0}, \frac{\Theta_M}{\Theta_{\mathcal{F}_0}} \right) \xrightarrow{D_2} \mathscr{H}om_{\mathcal{O}_M}\left( \bigwedge^3 \Theta_{\mathcal{F}_0}, \frac{\Theta_M}{\Theta_{\mathcal{F}_0}} \right) \xrightarrow{D_3} \cdots
\end{align}
with infinitesimal deformations of $\left(M, \Theta_{\mathcal{F}_0}\right)$.\footnote{We will denote the $i$-cohomology group of $\Theta_{\mathcal{F}_0}^\bullet$ by $\mathbb{H}^i\left(M, \Theta_{\mathcal{F}_0}^\bullet \right)$.}  In his arguments, he implicitly assumes that the underlying local flat deformations of $\Theta_{\mathcal{F}_0}$ are the trivial extension, which does not hold for arbitrary coherent sheaves $\Theta_{\mathcal{F}_0}$. However it holds for locally free sheaves $\Theta_{\mathcal{F}_0}$. That's why we assume that $\Theta_{\mathcal{F}_0}$ is locally free in order to fully use the leaf complex $(\ref{I1})$ for the study of deformations of foliated complex analytic structures in terms of  subsheaves of tangent sheaves. The problem at that time was the following: our deformations of $(M, \omega_0)$ as above is described in terms of cotangent sheaves but the leaf complex $(\ref{I1})$ is described in terms of tangent sheaves. The author's initial thought was that they should give equivalent deformations. However there was no way to transform the leaf complex in terms of cotangent sheaf except for regular foliations. For a regular holomorphic foliation $\Theta_{\mathcal{F}_0}$ on $M$, we have isomorphisms $\Theta_{\mathcal{F}_0}\cong \left(\Omega_M^1/\mathcal{N}_{\mathcal{F}_0}^* \right)^*$ and $\mathcal{N}_{\mathcal{F}_0}^* \cong \left( \Theta_M/\Theta_{\mathcal{F}_0} \right)^*$ and we have the isomorphism defined by the dual
{\small{\begin{center}
$\begin{CD}
\Theta_{\mathcal{F}_0}^\bullet:@. \Theta_M @>D_0>> \mathscr{H}om_{\mathcal{O}_M}\left( \Theta_{\mathcal{F}_0}, \frac{\Theta_M}{\Theta_{\mathcal{F}_0}} \right) @>D_1>> \mathscr{H}om_{\mathcal{O}_M}\left( \bigwedge^2 \Theta_{\mathcal{F}_0}, \frac{\Theta_M}{\Theta_{\mathcal{F}_0}}\right) @>D_2>> \mathscr{H}om_{\mathcal{O}_M}\left( \bigwedge^3 \Theta_{\mathcal{F}_0}, \frac{\Theta_M}{\Theta_{\mathcal{F}_0} }\right) @>D_3>> \cdots\\
 @.@V-\textnormal{Id}VV @V\cong VV @V \cong VV @V\cong VV \\
 \tilde{\mathcal{N}}_{\mathcal{F}_0}^{*\bullet} :@.\Theta_M @>E_0>> \mathscr{H}om_{\mathcal{O}_M}\left(\mathcal{N}_{\mathcal{F}_0}^*,  \frac{\Omega_M^1}{ \mathcal{N}_{\mathcal{F}_0}^* }  \right) @>\tilde{E}_1>> \mathscr{H}om_{\mathcal{O}_M}\left(\mathcal{N}_{\mathcal{F}_0}^*, \bigwedge^2 \frac{\Omega_M^1}{\mathcal{N}_{\mathcal{F}_0}^*}  \right) @> \tilde{E}_2>> \mathscr{H}om_{\mathcal{O}_M}\left( \mathcal{N}_{\mathcal{F}_0}^*, \bigwedge^3 \frac{\Omega_M^1}{\mathcal{N}_{\mathcal{F}_0}^*} \right) @> \tilde{E}_3>> \cdots 
\end{CD}$
\end{center}}}
From the isomorphism $\Theta_{\mathcal{F}_0}^{\bullet}\cong \tilde{\mathcal{N}}_{\mathcal{F}_0}^{*\bullet}$ and $D_0, D_1, D_2, \cdots$ of the leaf complex $\Theta_{\mathcal{F}_0}^\bullet$, the author could define ${E}_0, \tilde{E}_1, \tilde{E}_2,\tilde{E}_3, \cdots$ in the following way: let $\mathcal{U}=\{U_i\}$ be an open covering of $M$ by coordinated neighborhoods such that $\mathcal{N}_{\mathcal{F}_0}^*$ on $U_i$ is generated by $w_{i}^1,..., w_{i}^q\in \Gamma(U_i, \Omega_M^1)$. Then we can write $dw_{i}^\alpha= \sum_{\beta=1}^q a_i^{\alpha\beta}\wedge w_{i}^\beta$ for some $a_i^{\alpha\beta}\in \Gamma(U_i, \Omega_M^1)$ for $\alpha=1,...,q$. Then we define $E_0$ and $\tilde{E}_r$ (see Appendix \ref{appendixA} for the detail) by
\begin{align*}
E_0(T_i)\left(w_{i}^\alpha\right) &= \overline{ \mathcal{L}_{T_i}\left(w_{i}^\alpha\right)} ,\,\,\,\,\,\,\textnormal{where}\,\,\,\mathcal{L}= \textnormal{Lie derivative} ,\,\,\, \alpha=1,...,p \\
\tilde{E}_r(\phi_i)\left(w_{i}^\alpha\right)&= \overline{ d\phi_i\left(w_{i}^\alpha\right) - \sum_{\beta=1}^q a_{i}^{\alpha\beta} \wedge \phi_i\left(w_{i}^\beta\right) },\,\,\,\,\,\,\,\,\textnormal{for}\,\,\,r\geq 1,\,\,\,\,\,\phi_i\in \Gamma\left(U_i,  \mathscr{H}om_{\mathcal{O}_M}\left( \mathcal{N}_{\mathcal{F}_0}^*, \bigwedge^r \frac{\Omega_M^1}{\mathcal{N}_{\mathcal{F}_0}^*} \right) \right)
\end{align*}

Later the author understood that the leaf complex $\Theta_{\mathcal{F}_0}^\bullet$ controls deformations of $(M, \Theta_{\mathcal{F}_0})$ which includes flat deformations of underlying locally free sheaf $\Theta_{\mathcal{F}_0}$ without assuming the flatness of deformations of $\mathcal{N}_{\mathcal{F}_0}^*$. As is known, flat deformations of a locally free sheaf $\Theta_{\mathcal{F}_0}$ is controlled by $\mathscr{H}om_{\mathcal{O}_M}\left( \Theta_{\mathcal{F}_0}, \Theta_{\mathcal{F}_0}    \right)$ (\cite{Har10}). This information is contained in $\mathscr{H}om_{\mathcal{O}_M}\left(\Theta_{\mathcal{F}_0}, \frac{\Theta_M}{\Theta_{\mathcal{F}_0}}  \right)$ in the leaf complex $\Theta_{\mathcal{F}_0}^\bullet$ (see Part II).  However there is no information in the leaf complex $\Theta_{\mathcal{F}_0}^\bullet$ on $\mathscr{H}om_{\mathcal{O}_M}\left(\mathcal{N}_{\mathcal{F}_0}^*, \mathcal{N}_{\mathcal{F}_0}^* \right)$ which controls flat deformations of a locally free sheaf $\mathcal{N}_{\mathcal{F}_0}^*$.  Hence we have to distinguish deformations of $(M, \Theta_{\mathcal{F}_0})$ from deformations of $(M, \mathcal{N}_{\mathcal{F}_0}^*)$. As we see as above, the information of  $\mathscr{H}om_{\mathcal{O}_M}\left(\mathcal{N}_{\mathcal{F}_0}^*, \mathcal{N}_{\mathcal{F}_0}^* \right)$ controling flat deformations of a locally free sheaf $\mathcal{N}_{\mathcal{F}_0}^*$ is contained in $\mathscr{H}om_{\mathcal{O}_M}\left( \mathcal{N}_{\mathcal{F}_0}^*, \frac{\Omega_M^1}{\mathcal{N}_{\mathcal{F}_0}^*}  \right)$ of $\tilde{\mathcal{N}}_{\mathcal{F}_0}^{*\bullet}$ (see Part II). It was natural for the author to expect that $\tilde{\mathcal{N}}_{\mathcal{F}_0}^{*\bullet}$ might control deformations of $(M, \mathcal{N}_{\mathcal{F}_0}^*)$ when $\mathcal{N}_{\mathcal{F}_0}^*$ is locally free. However there was a problem: our definition of $\tilde{E}_r (r\geq 1)$ is defined on the basis of the fact that $dw_{i}^\alpha= \sum_{\beta=1}^q a_i^{\alpha\beta}\wedge w_{i}^\beta$ at every point of $M$, which is possible when $\mathcal{N}_{\mathcal{F}_0}^*$ defines a regular foliation. However when $\mathcal{N}_{\mathcal{F}_0}^*$ defines a singular holomorphic foliation on $M$, we cannot guarantee that $dw_{i}^\alpha$ is of the form $ \sum_{\beta=1}^q a_i^{\alpha\beta}\wedge w_{i}^\beta$ around a point in $\textnormal{Sing}(\mathcal{F}_0)$. Hence we have to modify $\tilde{\mathcal{N}}_{\mathcal{F}_0}^{*\bullet}$. We note that $\bigwedge^r \frac{\Omega_M^1}{\mathcal{N}_{\mathcal{F}_0}^*}(r\geq 2)$ are locally free on $M-\textnormal{Sing}(\mathcal{F}_0)$, so that on $M-\textnormal{Sing}(\mathcal{F}_0)$, the complex $\tilde{\mathcal{N}}_{\mathcal{F}_0}^{*\bullet}|_{M-\textnormal{Sing}(\mathcal{F}_0)}$ is well-defined. The author's strategy was to replace $\bigwedge^r \frac{\Omega_M^1}{\mathcal{N}_{\mathcal{F}_0}^*}(r\geq 2)$ by some reflexive sheaves $\tilde{\mathcal{S}}^r(r\geq 2)$ which is isomorphic to $\bigwedge^r \frac{\Omega_M^1}{\mathcal{N}_{\mathcal{F}_0}^*}$ on $M-\textnormal{Sing}(\mathcal{F}_0)$ and construct a complex of sheaves $\mathcal{N}_{\mathcal{F}_0}^{*\bullet}$ which is isomorphic to $\tilde{\mathcal{N}}_{\mathcal{F}_0}^{*\bullet}$ on $M-\textnormal{Sing}(\mathcal{F}_0)$, so that by reflexiveness of $\mathscr{H}om_{\mathcal{O}_M}\left( \mathcal{N}_{\mathcal{F}_0}^*, \tilde{S}^r \right)(r\geq 2)$, we can naturally extend $\tilde{E}_i|_{M-\textnormal{Sing}(\mathcal{F}_0)}$ to $E_i(i\geq 1)$. We define $\tilde{\mathcal{S}}^r$ in the following way: first we note that $\mathcal{N}_{\mathcal{F}_0}^*$ corresponds to a section $\omega_0 \in H^0\left(M, \bigwedge^q \Omega_M^1\otimes \mathcal{L} \right)$, where $\mathcal{L}=\bigwedge^q \left(\mathcal{N}_{\mathcal{F}_0}^*\right)^*$. Then we define $\tilde{\mathcal{S}}^r$ as a subsheaf of $\bigwedge^{q+r}\Omega_M^1\otimes \mathcal{L}$, and any section $\phi\in \Gamma\left(U,\tilde{\mathcal{S}}^r \right)\subset \Gamma\left( U, \bigwedge^{q+r} \Omega_M^1\otimes \mathcal{L}  \right)$ on an open set $U$ of $M$ is defined by the following property that for any $x\in U- \textnormal{Sing}({\mathcal{F}_0})$, there is a neighborhood $U_x$ such that $\phi$ is of the form $\phi= \omega_0 \wedge A_{x}$ for some $A_x \in \Gamma(U_x, \bigwedge^r \Omega_M^1)$. Then we have the dual leaf complex $\mathcal{N}_{\mathcal{F}_0}^{*\bullet}$ constructed from $\tilde{\mathcal{N}}_{\mathcal{F}_0}^{*\bullet}|_{M- \textnormal{Sing}(\mathcal{F}_0)}$:
\begin{align}\label{I2}
\mathcal{N}_{\mathcal{F}_0}^{*\bullet}:\Theta_M \xrightarrow{E_0} \mathscr{H}om_{\mathcal{O}_M}\left( \mathcal{N}_{\mathcal{F}_0}^* , \frac{\Omega_M^1}{\mathcal{N}_{\mathcal{F}_0}^* } \right) \xrightarrow{E_1} \mathscr{H}om_{\mathcal{O}_M}\left( \mathcal{N}_{\mathcal{F}_0}^*, \tilde{\mathcal{S}}^2 \right) \xrightarrow{E_2} \mathscr{H}om_{\mathcal{O}_M}\left( \mathcal{N}_{\mathcal{F}_0}^*, \tilde{\mathcal{S}}^3  \right) \xrightarrow{E_3} \cdots 
\end{align}
We will denote the $i$-th cohomology group of $\mathcal{N}_{\mathcal{F}_0}^{*\bullet}$ by $\mathbb{H}^i\left( M, \mathcal{N}_{\mathcal{F}_0}^{*\bullet} \right)$. We shall show that the dual leaf complex $\mathcal{N}_{\mathcal{F}_0}^{*\bullet}$ controls foliated deformations of $(M, \mathcal{N}_{\mathcal{F}_0}^*)$ in terms of cotangent sheaves when $\mathcal{N}_{\mathcal{F}_0}^*$ is locally free. Later the author understood that there is a notion of foliated modules on $(M, \Theta_{\mathcal{F}_0})$ as an analogue of a notion of Poisson modules on holomorphic Poisson manifolds (see \cite{Pol97}). While the author tried to emulate his previous work on deformations of Poisson invertible sheaves and Poisson vector bundles (\cite{Kim16}) in the context of deformations of foliated invertible sheaves and foliated vector bundles (which is the subject of Part V), he understood that the leaf complex $(\ref{I1})$ and the dual leaf complex $(\ref{I2})$ can be constructed by using the notion of foliated modules in the following way:
\begin{definition}[see Part V and compare \cite{Pol97}, \cite{Kim16}]
Let $\left(M, \Theta_{\mathcal{F}_0}\right)$ be a compact foliated complex manifold. A foliated connection on a $\mathcal{O}_M$-module $G$ is a $\mathbb{C}$-linear bracket $\{-,-\}_G:\Theta_{\mathcal{F}_0} \times G \to G$ such that $\left\{fT, a \right\}_G= f \left\{T,a \right\}_G$ and $\left\{T, fa \right\}_G= \left[T,f \right]a+f\left\{T, a\right\}_G$ where $T\in \Theta_{\mathcal{F}_0}, f\in \mathcal{O}_M,a\in G$. Equivalently, a foliated connection is given by a homomorphism $\nabla_G^0 :G\to \mathscr{H}om_{\mathcal{O}_X}\left(\Theta_{\mathcal{F}_0}, G \right)$ which satisfies the identity
\begin{align*}
\nabla_G^0(fa)(T)=[T, f]a + f\nabla_G^0(a)(T)
\end{align*}
 and the bracket is defined by the formula $\nabla_G^0(a)(T)=\{T,a\}_G$. A foliated connection $\nabla_G^0$ is called flat if it satisfies the Jacobi identity
\begin{align*}
\left\{T_1, \left\{T_2,a \right\}_G \right\}_G = \left\{ \left[T_1,T_2 \right], a \right\}_G+ \left\{T_2, \left\{T_1,a\right\}_G \right\}_G,\,\,\,\,\,\,\,\,\,T_1,T_2\in \Theta_\mathcal{F}, a\in G.
\end{align*}
An $\mathcal{O}_M$-module $G$ equipped with a flat foliated connection is called a foliated $\mathcal{O}_M$-module.  Then we have a complex of sheaves of $G^\bullet$ associated to $\left( G, \{-,-\}_G\right)$ on $\left(M, \Theta_{\mathcal{F}_0}\right)$$:$
\begin{align*}
G^\bullet: G\xrightarrow{\nabla_G^0} \mathscr{H}om_{\mathcal{O}_X}\left(\Theta_{\mathcal{F}_0}, G \right) \xrightarrow{\nabla_G^1} \mathscr{H}om_{\mathcal{O}_X}\left( \bigwedge^2 \Theta_{\mathcal{F}_0}, G  \right) \xrightarrow{\nabla_G^2} \mathscr{H}om_{\mathcal{O}_X}\left( \bigwedge^3 \Theta_{\mathcal{F}_0}, G \right) \xrightarrow{\nabla_G^3}\to \cdots
\end{align*}
\end{definition}
Now assume that $\Theta_{\mathcal{F}_0}$ is locally free. We note that $\frac{\Theta_M}{\Theta_{\mathcal{F}_0}}$ is a foliated $\mathcal{O}_M$-module with the bracket $\{-,-\}_{\frac{\Theta_M}{\Theta_{\mathcal{F}_0}}}: \Theta_{\mathcal{F}_0}\times \frac{\Theta_M}{\Theta_{\mathcal{F}_0}} \to \frac{\Theta_M}{\Theta_{\mathcal{F}_0}}$ defined by $\left\{T, \bar{V}\right\}_{\frac{\Theta_M}{\Theta_{\mathcal{F}_0}}}:=\overline{\left[ T, V  \right]}$ for $T\in \Theta_{\mathcal{F}_0}$ and $\bar{V}$ is the image of $V$ in $\frac{\Theta_M}{\Theta_{\mathcal{F}_0}}$, so that we have a complex of sheaves
\begin{align*}
\frac{\Theta_M}{\Theta_{\mathcal{F}_0} }  \xrightarrow{\tilde{D}_0} \mathscr{H}om_{\mathcal{O}_M}\left( \Theta_{\mathcal{F}_0}, \frac{\Theta_M}{\Theta_{\mathcal{F}_0}} \right) \xrightarrow{D_1} \mathscr{H}om_{\mathcal{O}_M}\left(\bigwedge^2 \Theta_{\mathcal{F}_0}, \frac{\Theta_M}{\Theta_{\mathcal{F}_0}} \right) \xrightarrow{D_2} \mathscr{H}om_{\mathcal{O}_M}\left( \bigwedge^3 \Theta_{\mathcal{F}_0}, \frac{\Theta_M}{\Theta_{\mathcal{F}_0}} \right) \xrightarrow{D_3} \cdots
\end{align*}
Then the leaf complex $(\ref{I1})$ can be constructed by defining $D_0$ to be the composition of the canonical map $\Theta_M\to \frac{\Theta_M}{\Theta_{\mathcal{F}_0}}$ with $-\tilde{D}_0$. 

On the other hand, assume that $\mathcal{N}_{\mathcal{F}_0}^*$ is locally free ($\Theta_{\mathcal{F}_0}$ is not necessarily locally free). We note that there is a foliated module on $\left(\mathcal{N}_{\mathcal{F}_0}^*\right)^*$ with the bracket $\{-,-\}_{\left(\mathcal{N}_{\mathcal{F}_0}^*\right)^*}:\Theta_{\mathcal{F}_0}\times \left( \mathcal{N}_{\mathcal{F}_0}^* \right)^* \to \left( \mathcal{N}_{\mathcal{F}_0}^* \right)^*$ defined by $\left\{ T, \phi \right\}_{\left( \mathcal{N}_{\mathcal{F}_0}^* \right)^*}(w)= \left[ T, \phi(w)\right]- \phi\left( \mathcal{L}_T(w) \right)$ for $T\in \Theta_{\mathcal{F}_0}, \phi \in \mathscr{H}om_{\mathcal{O}_M}\left(\mathcal{N}_{\mathcal{F}_0}^*, \mathcal{O}_M \right)$ and $w\in \mathcal{N}_{\mathcal{F}_0}^*$. Then we have a complex of sheaves 
{\small{\begin{align*}
\mathscr{H}om_{\mathcal{O}_M}\left( \mathcal{N}_{\mathcal{F}_0}^*, \mathcal{O}_M \right) \to  \mathscr{H}om_{\mathcal{O}_M}\left(  \Theta_{\mathcal{F}_0}, \left(\mathcal{N}_{\mathcal{F}_0}^*\right)^* \right) \to \mathscr{H}om_{\mathcal{O}_M}\left( \bigwedge^2 \Theta_{\mathcal{F}_0}, \left(\mathcal{N}_{\mathcal{F}_0}^*\right)^*   \right) \to \mathscr{H}om_{\mathcal{O}_M}\left( \bigwedge^3 \Theta_{\mathcal{F}_0}, \left(\mathcal{N}_{\mathcal{F}_0}^*\right)^* \right) \to \cdots
\end{align*}}}
which induces
{\small{\begin{align}\label{I3}
\mathscr{H}om_{\mathcal{O}_M}\left( \mathcal{N}_{\mathcal{F}_0}^*, \mathcal{O}_X \right) \to \mathscr{H}om_{\mathcal{O}_M}\left( \mathcal{N}_{\mathcal{F}_0}^*, \Theta_{\mathcal{F}_0}^* \right) \to \mathscr{H}om_{\mathcal{O}_M}\left( \mathcal{N}_{\mathcal{F}_0}^*, \left( \bigwedge^2 \Theta_{\mathcal{F}_0}\right)^* \right)\to \mathscr{H}om_{\mathcal{O}_M}\left( \mathcal{N}_{\mathcal{F}_0}^*, \left( \bigwedge^3 \Theta_{\mathcal{F}_0}  \right)^*   \right) \to \cdots
\end{align}}}
We note that we have a map $\Theta_M \xrightarrow{E_0} \mathscr{H}om_{\mathcal{O}_M}\left( \mathcal{N}_{\mathcal{F}_0}^*, \frac{\Omega_M^1}{\mathcal{N}_{\mathcal{F}_0}^*} \right)$ defined by $E_0(T)(\omega)=\overline{\mathcal{L}_T(\omega)}$, where $\mathcal{L}$ is the Lie derivative, and the inclusion  $\frac{\Omega_X^1}{\mathcal{N}_{\mathcal{F}_0}^*}\hookrightarrow \Theta_{\mathcal{F}_0}^*$, so that we have $\Theta_M \xrightarrow{E_0} \mathscr{H}om_{\mathcal{O}_M } \left(\mathcal{N}_{\mathcal{F}_0}^*, \frac{\Omega_X^1}{\mathcal{N}_{\mathcal{F}_0}^*} \right)\to \mathscr{H}om_{\mathcal{O}_X}\left( \mathcal{N}_{\mathcal{F}_0}^*, \Theta_{\mathcal{F}_0}^* \right)$. By combining this with $(\ref{I3})$, we obtain a complex of sheaves
\begin{align}\label{I5}
\Theta_M \xrightarrow{E_0} \mathscr{H}om_{\mathcal{O}_X}\left( \mathcal{N}_{\mathcal{F}_0}^*, \frac{\Omega_M^1}{ \mathcal{N}_{\mathcal{F}_0}^*} \right) \to \mathscr{H}om_{\mathcal{O}_X}\left(  \mathcal{N}_{\mathcal{F}_0}^*, \left( \bigwedge^2 \Theta_{\mathcal{F}_0} \right)^* \right) \to  \mathscr{H}om_{\mathcal{O}_X}\left( \mathcal{N}_{\mathcal{F}_0}^*, \left( \bigwedge^3 \Theta_{\mathcal{F}_0} \right)^* \right) \to \cdots
\end{align}
Since $\left(\bigwedge^r \Theta_{\mathcal{F}_0}\right)^*(r\geq 2)$ is reflexive, and $\left(\bigwedge^r \Theta_{\mathcal{F}_0} \right)^*|_{M-\textnormal{Sing}(\mathcal{F}_0)} \cong 
\left( \bigwedge^r \frac{\Omega_M^1}{\mathcal{N}_{\mathcal{F}_0}^*}\right)|_{M- \textnormal{Sing}(\mathcal{F}_0)}\cong \tilde{S}^r|_{M-\textnormal{Sing}(\mathcal{F}_0)}$, we have $\left(\bigwedge^r \Theta_{\mathcal{F}_0}\right)^* \cong \tilde{\mathcal{S}}^r$. Then we have an isomorphism between the dual leaf complex $\mathcal{N}_{\mathcal{F}_0}^{*\bullet}$ $(\ref{I2})$ and $(\ref{I5})$
{\small{\begin{center}
$\begin{CD}
\Theta_M @>E_0>> \mathscr{H}om_{\mathcal{O}_M}\left( \mathcal{N}_{\mathcal{F}_0}^*,   \frac{\Omega_M^1}{\mathcal{N}_{\mathcal{F}_0}^*}  \right) @>E_1>> \mathscr{H}om_{\mathcal{O}_M}\left( \mathcal{N}_{\mathcal{F}_0}^*,  \tilde{\mathcal{S}}^2 \right) @>E_2>> \mathscr{H}om_{\mathcal{O}_M}\left( \mathcal{N}_{\mathcal{F}_0}^*, \tilde{\mathcal{S}}^3 \right) @>>> \cdots \\
@| @| @V\cong VV @V\cong VV\\
\Theta_M @>E_0>>  \mathscr{H}om_{\mathcal{O}_M}\left( \mathcal{N}_{\mathcal{F}_0}^*, \frac{\Omega_M^1}{\mathcal{N}_{\mathcal{F}_0}^*} \right) @>>> \mathscr{H}om_{\mathcal{O}_M}\left( \mathcal{N}_{\mathcal{F}_0}^*, \left( \bigwedge^2 \Theta_{\mathcal{F}_0} \right)^* \right) @>>> \mathscr{H}om_{\mathcal{O}_M}\left( \mathcal{N}_{\mathcal{F}_0}^* , \left( \bigwedge^3 \Theta_{\mathcal{F}_0} \right)^* \right) @>>> \cdots
\end{CD}$
\end{center}}}

Now let us go back to our initial approach to holomorphic foliated deformations in terms of sections, inspired by Horikawa's work. We recall that we constructed the complex of sheaves $(\ref{I8})$
\begin{align*}
\mathcal{E}_{\mathcal{L}_0} \xrightarrow{\tilde{G}} \mathcal{S}^1  \xrightarrow{G_1}  \left( \mathcal{N}_{\mathcal{F}_0}^* \right)^* \otimes \bigwedge^{q+2}\Omega_M^1\otimes \mathcal{L}_0 \cong \mathscr{H}om_{\mathcal{O}_M}\left( \mathcal{N}_{\mathcal{F}_0}^*  , \bigwedge^{q+2} \Omega_M^1 \otimes \mathcal{L}_0 \right) \to \textbf{?} \to \textbf{?} \to \cdots
\end{align*}
We note that actually the image of $G_1$ lies in $\mathscr{H}om_{\mathcal{O}_M}\left(\mathcal{N}_{\mathcal{F}_0}^*, \tilde{\mathcal{S}}^2 \right)$, so that the dual leaf complex $\mathcal{N}_{\mathcal{F}_0}^{*\bullet}$ enables us to extend the complex of sheaves $(\ref{I8})$ to 
\begin{align*}
\mathcal{E}_{\mathcal{L}_0}^\bullet: \mathcal{E}_{\mathcal{L}_0} \to \tilde{S}^1 \to  \mathscr{H}om_{\mathcal{O}_M}\left( \mathcal{N}_{\mathcal{F}_0}^*, \tilde{\mathcal{S}}^2 \right) \to \mathscr{H}om_{\mathcal{O}_M}\left( \mathcal{N}_{\mathcal{F}_0}^*, \tilde{\mathcal{S}}^3 \right) \to \cdots
\end{align*}
We will denote the $i$-th cohomology group of by $\mathcal{E}_{\mathcal{L}_0}^\bullet$ by $\mathbb{H}^i\left( M, \mathcal{E}_{\mathcal{L}_0}^\bullet \right)$. We note that we have an exact sequence of complexes of sheaves (see Part II)
\begin{center}
$\begin{CD}
@. \cdots @. \cdots @. \cdots \\
@. @AAA @AAA @AAA \\
0 @>>> 0 @>>>  \mathscr{H}om_{\mathcal{O}_M}\left( \mathcal{N}_{\mathcal{F}_0}^*, \tilde{\mathcal{S}}^3  \right) @>>>  \mathscr{H}om_{\mathcal{O}_M}\left( \mathcal{N}_{\mathcal{F}_0}^*, \tilde{\mathcal{S}}^3 \right) @>>> 0 \\
@. @AAA @AAA @AAA \\
0 @>>> 0 @>>> \mathscr{H}om_{\mathcal{O}_M}\left( \mathcal{N}_{\mathcal{F}_0}^*, \tilde{\mathcal{S}}^2 \right) @>>> \mathscr{H}om_{\mathcal{O}_M} \left( \mathcal{N}_{\mathcal{F}_0}^*, \tilde{\mathcal{S}}^2  \right) @>>> 0 \\
@. @AAA @AAA @AAA \\
0 @>>> \mathcal{O}_M @>>> \mathcal{S}^1 @>>> \mathscr{H}om_{\mathcal{O}_M}\left( \mathcal{N}_{\mathcal{F}_0}^*, \frac{\Omega_M^1}{ \mathcal{N}_{\mathcal{F}_0}^* } \right) @>>> 0 \\
@. @AAA @AAA @AAA \\
0 @>>> \mathcal{O}_M @>>> \mathcal{E}_{\mathcal{L}_0} @>>> \Theta_M @>>> 0
\end{CD}$
\end{center}
where $0 \to \mathcal{O}_M \to \mathcal{E}_{\mathcal{L}_0}\to \Theta_M \to 0$ is the Atiyah extension associated to $\mathcal{L}_0$. The above exact sequence implies that
\begin{align*}
\mathbb{H}^i\left( M, \mathcal{E}_{\mathcal{L}_0}^\bullet \right) \cong \mathbb{H}^i\left( M,  \mathcal{N}_{\mathcal{F}_0}^{*\bullet} \right)
\end{align*} 
which shows that foliated deformations of $\left(M, \mathcal{N}_{\mathcal{F}_0}^* \right)$ is equivalent to foliated deformations of $\left(M, \omega_0\right)$. Initially the author proved theorem of existence and completeness for deformations of $\left(M, \mathcal{N}_{\mathcal{F}_0}^* \right)$ with $\mathcal{N}_{\mathcal{F}_0}^*$ locally free by using the complex of sheaves $\mathcal{E}_{\mathcal{L}_0}^\bullet$ by which we can extend the method of Horikawa \cite{Hor76}. However the author decided to use the dual leaf complex $\mathcal{N}_{\mathcal{F}_0}^{*\bullet}$ since $\mathcal{N}_{\mathcal{F}_0}^{*\bullet}$ looks more standard, and is used to define a complex of sheaves controlling simultaneous deformations of $\left(M, \Theta_{\mathcal{F}_0}, \mathcal{N}_{\mathcal{F}_0}^*\right)$ when $\Theta_{\mathcal{F}_0}$ and $\mathcal{N}_{\mathcal{F}_0}^*$ are both locally free as below. However we present infinitesimal foliated deformations of $\left(M, \mathcal{N}_{\mathcal{F}_0}^* \right)$ by using $\mathcal{E}_{\mathcal{L}_0}^\bullet$ over Artin rings in Part II, and study obstructions of codimension $1$-foliated deformations of $\left( M, \mathcal{N}_{\mathcal{F}_0}^* \right)$ by using $\mathcal{E}_{\mathcal{L}_0}^{\bullet}$ in terms of deformations of sections in Part IV-2. This is the essential story about how this present work has been started and developed.

 It is natural to consider simultaneous deformations of $\left(M, \Theta_{\mathcal{F}_0}\right)$ and $\left(M, \mathcal{N}_{\mathcal{F}_0}^*\right)$ which includes flat deformations of $\Theta_{\mathcal{F}_0}$ and $\mathcal{N}_{\mathcal{F}_0}^*$ simultaneously  when $\Theta_{\mathcal{F}_0}$ and $\mathcal{N}_{\mathcal{F}_0}^*$ are both locally free, and ask what kind of complex of sheaves controls simultaneous deformations of $\left( M, \Theta_{\mathcal{F}_0}, \mathcal{N}_{\mathcal{F}_0}^* \right)$. By combining the leaf complex $(\ref{I1})$ and the dual leaf complex $(\ref{I2})$, we shall construct a complex of sheaves $\mathcal{F}_0^\bullet$ (see Appendix \ref{AppendixA4}):
{\small{\begin{center}
$\begin{CD}
\cdots \\
@AF_3AA \\
\mathscr{H}om_{\mathcal{O}_M}\left(\bigwedge^3 \Theta_{\mathcal{F}_0},  \frac{\Theta_M}{\Theta_{\mathcal{F}_0}} \right) \bigoplus \mathscr{H}om_{\mathcal{O}_M}\left( \mathcal{N}_{\mathcal{F}_0}^*, \tilde{\mathcal{S}}^3 \right) \bigoplus \mathscr{H}om_{\mathcal{O}_M}\left( \mathcal{N}_{\mathcal{F}_0}^*, \bigwedge^2 \Theta_{\mathcal{F}_0}^* \right) \\
@AF_2AA \\
\mathscr{H}om_{\mathcal{O}_M}\left( \bigwedge^2 \Theta_{\mathcal{F}_0}, \frac{\Theta_M}{\Theta_{\mathcal{F}_0}} \right) \bigoplus \mathscr{H}om_{\mathcal{O}_M}\left( \mathcal{N}_{\mathcal{F}_0}^*, \tilde{\mathcal{S}}^2 \right) \bigoplus \mathscr{H}om_{\mathcal{O}_M}\left(\mathcal{N}_{\mathcal{F}_0}^*, \Theta_{\mathcal{F}_0}^* \right) \\
@AF_1AA\\
\mathscr{H}om_{\mathcal{O}_M}\left(\Theta_{\mathcal{F}_0}, \frac{\Theta_M}{\Theta_{\mathcal{F}_0}}   \right) \bigoplus \mathscr{H}om_{\mathcal{O}_M}\left( \mathcal{N}_{\mathcal{F}_0}^* , \frac{\Omega_M^1}{\mathcal{N}_{\mathcal{F}_0}^*} \right)\\
@AF_0AA\\
\Theta_M
\end{CD}$
\end{center}}}
and show that the complex of sheaves $\mathcal{F}_0^\bullet$ controls simultaneous foliated deformations of $\left( M, \Theta_{\mathcal{F}_0}, \mathcal{N}_{\mathcal{F}_0}^* \right)$ when $\Theta_{\mathcal{F}_0}$ and $\mathcal{N}_{\mathcal{F}_0}^*$ are both locally free.
 
 In Part II we study deformations of foliated nonsingular varieties defined by locally free subsheaves of tangent sheaves and locally free subsheaves of cotangent sheaves in the language of functors of Artin rings which is the algebraic version of deformations of foliated complex analytic structures. We identify first-order deformations and obstructions.

In Part III we study unfoldings of foliated complex analytic structures defined by locally free subsheaves of cotangent sheaves. We prove theorem of existence for unfoldings of foliated complex analytic structures defined by locally free subsheaves of cotangent sheaves for one-parameter families. The author could not succeed in proving the complete form of theorem of existence for unfoldings and showed only for one-parameter families. We conjecture the complete form of theorem of existence for unfoldings of compact foliated complex manifolds defined by locally free subsheavse of cotangent sheaves. We also prove (the complete form of) theorem of completeness for (global) unfoldings of singular holomorphic foliations of arbitrary codimensions defined by locally free subsheaves of cotanget sheaves which generalizes Suwa's work \cite{Suw81} on the versality theorem for local unfoldings of codimension $1$ foliations.

In Part IV, we study obstructed and unobstructed holomorphic foliated deformations and unfoldings of compact complex foliated manifolds defined by locally free subsheaves of tangent sheaves and cotangent sheaves.

In Part V we define a notion of foliated modules as an analogue of Poisson modules in holomorphic Poisson geometry and study deformations of foliated invertible sheaves and foliated vector bundles on foliated nonsingular varieties defined by locally free subsheaves of tangent sheaves. This is the complete `foliated' version of the author's previous work on deformations of Poisson invertible sheaves and Poisson vector bundles which is based on Sernesi's textbook \cite{Ser06} on deformations of algebraic schemes.

In Part VI we study deformations of a special kind of compact invariant complex submanifolds of foliated complex manifolds defined by locally free subehaves of tangent sheaves. This is the complete `foliated' version of the author's previous work on deformations of compact holomorphic Poisson submanifolds which is based on Kodaira's series of papers \cite{Kod59}, \cite{Kod62}, \cite{Kod63} on deformations of compact complex submanifolds.

In Part VII we study deformations of a special kind of holomorphic maps between foliated complex manifolds defined by locally free subsheaves of tangent sheaves which generalize the special kind of compact invariant complex submanifolds in Part VI. This is the complete `foliated' version of the author's previous work on deformations of holomorphic Poisson maps which is based on Horikawa's series of papers \cite{Hor73}, \cite{Hor74}, \cite{Hor75} on deformations of holomorphic maps.

Codimension $1$ holomorphic foliations defined by $\omega_0 \in H^0\left(M , \Omega_M^1\otimes K_M^{-1} \right)$  on a compact complex threefold $X$ define holomorphic Poisson structures $\Lambda_{\omega_0} \in H^0\left(M, \bigwedge^2 \Theta_M \right)$ on $M$ and so we can also consider their holomorphic Poisson deformations. There is a natural deformation functor from deformations of holomorphic Poisson structures of $\left(M, \Lambda_{\omega_0} \right)$ to deformations of holomorphic foliated structures of $\left(M, \omega_0\right)$ in terms of cotangent sheaves (equivalently in terms of sections). In the language of functors of Artin rings, we have a natural morphism from $\textnormal{Def}_{\left(X, \Lambda_{\omega_0} \right)}$ (see \cite{Kim16}) to $\textnormal{FDef}_{\left(M, \mathcal{N}_{\mathcal{F}_0}^*\cong K_M \right)}$ (see Part II). In Part VIII we compare holomorphic Poisson deformations of compact complex threefolds defined by codimension $1$-foliations and holomorphic foliated deformations of them in terms of locally free subsheaves of cotangent sheaves (equivalently in terms of deformations of sections).

\section{Deformations of foliated complex analytic structures in terms of tangent sheaves}

\begin{definition}[compare \cite{Kod05} p.59 and \cite{GM88} p.61]\label{t3}
Suppose that given a domain $B\subset \mathbb{C}^m$, there is a set $\left\{\left(M_t, \Theta_{\mathcal{F}_t}\right)|t\in B \right\} $ of $n$-dimensional compact $($singularly$)$ foliated complex manifolds $(M_t, \Theta_{\mathcal{F}_t})$ of dimension $p$, so that for each $t\in B$, we have an exact sequence
\begin{align}\label{t1}
0\to \Theta_{\mathcal{F}_t}\to \Theta_{M_t}\to \Theta_{M_t}/\Theta_{\mathcal{F}_t}\to 0
\end{align}
where $\Theta_{M_t}/\Theta_{\mathcal{F}_t}$ is a torsion-free $\mathcal{O}_{M_t}$-module. We say that $\left\{\left(M_t, \Theta_{\mathcal{F}_t}\right)|t\in B\right\}$ is a family of $($singularly$)$ foliated complex manifold or $($singularly$)$ foliated complex analytic family in terms of tangent sheaves if there exist a complex manifold $\mathcal{M}$ and a holomorphic map $\pi:\mathcal{M}\to B$  such that
\begin{enumerate}
\item we have an exact sequence
\begin{align}\label{t2}
0\to \Theta_\mathcal{F} \to \Theta_{\mathcal{M}/B}\to \Theta_{\mathcal{M}/B}/\Theta_{\mathcal{F}}\to 0
\end{align}
where $\Theta_{\mathcal{F}}$ is a coherent subsheaf of the relative tangent sheaf $\Theta_{\mathcal{M}/B}$ over $B$.
\item for each $t\in B$, $\pi^{-1}(t)= \left(M_t, \Theta_{\mathcal{F}_t}\right)$ and the exact sequence $(\ref{t1})$ is induced from $(\ref{t2})$ by restricting to $M_t$
\item the rank of the Jacobian of $\pi$ is equal to $m$ at every point of $\mathcal{M}$.
\item $[\Theta_\mathcal{F},\Theta_\mathcal{F}]\subset \Theta_\mathcal{F}$ under the Lie bracket and $\Theta_{\mathcal{M}/ B}/\Theta_\mathcal{F}$ is flat over $B$.
\end{enumerate}
We will denote the $($singularly$)$ foliated complex analytic family by $(\mathcal{M}, \Theta_\mathcal{F}, B, \pi)$.
\end{definition}

\begin{remark}
Let $\left(\mathcal{M}, \Theta_\mathcal{F}, B, \pi \right)$ with $\Theta_\mathcal{F}$ be a foliated analytic family. Let $\Delta$ be an open set of $B$. Then the restriction $\left(\mathcal{M}_\Delta= \pi^{-1}(\Delta), \Theta_{\mathcal{F}}|_\Delta, \Delta, \pi |_{\mathcal{M}_\Delta} \right)$ is also a foliated complex analytic family in terms of tangent sheaves. We will denote the family by $\left( \mathcal{M}_\Delta, \Theta_{\mathcal{F}_\Delta}, \Delta, \pi\right)$.
\end{remark}

\begin{remark}
When we ignore foliated structures, a foliated $($complex$)$ analytic family $\left( \mathcal{M}, \Theta_\mathcal{F}, B, \pi \right)$ is a complex analytic family $\left( \mathcal{M}, B, \pi \right)$ in the sense of Kodaira-Spencer $($see \textnormal{\cite{Kod05} p.59}$)$.
\end{remark}

\begin{remark}
We note that the flatness of $\Theta_{\mathcal{M}/B}/\Theta_\mathcal{F}$ implies the flatness of $\Theta_\mathcal{F}$ over $B$. In particular if $\Theta_{\mathcal{F}_{t_0}}$ is a locally free subsheaf of $\Theta_{M_{t_0}}$ for some $t_0\in B$, then $\Theta_{\mathcal{F}_t}$ is also locally free in some neighborhood of $t_0$. In the following we shall assume that $\Theta_\mathcal{F}$ is a locally free $\mathcal{O}_\mathcal{M}$-submodule of $\Theta_{\mathcal{M}/B}$ with rank $p$, so that $\Theta_{\mathcal{F}_t}$ is a locally free subsheaf of $\Theta_{M_t}$ for each $t\in B$.
\end{remark}

\begin{remark}
If $\Theta_{\mathcal{F}_{t_0}}$ is a regular foliation on $M_{t_0}$ for some $t_0\in B$, then $\Theta_{\mathcal{F}_t}$ defines a regular foliation on $M_t$ in some neigborhood of $t_0$.
\end{remark}

\begin{remark}\label{ta1}
From \textnormal{Definition \ref{t3}}, we can choose a system of local complex coordinates $\{z_i,...,z_j,...\}, z_j:p\to z_j(p)$, and coordinate polydisks $\mathcal{U}_j$ with respect to $z_j$, satisfying the following conditions\footnote{for the detail, we refer to \cite{Kod05} p.60.}:
\begin{enumerate}
\item $z_j(p)=\left(z_j^1(p),...,z_j^n(p),t_1,...,t_m\right),\left(t_1,...,t_m\right)=\pi(p)$;
\item $\mathcal{U}=\left\{\mathcal{U}_j|j=1,2,...\right\}$ is locally finite.
\end{enumerate}
Then $\left\{p\to \left(z_j^1(p),...,z_j^n(p)\right)|\mathcal{U}_j\cap M_t\ne \emptyset \right\}$ gives a system of local complex coordinates on $M_t$. In terms of these coordinates, $\pi$ is the projection given by $\pi: \left(z_j^1,...,z_j^n, t_1,...,t_m \right)\to \left(t_1,...,t_m\right)$. For $j,k$ with $\mathcal{U}_j\cap \mathcal{U}_k\ne \emptyset$, we denote that coordinate transformation from $z_k$ to $z_j$ by 
\begin{align*}
f_{jk}:\left(z_k^1,...,z_k^n,t \right)\to \left(z_j^1,...,z_j^n,t \right)=f_{jk}\left(z_k^1,...,z_k^n,t \right)
\end{align*}
Thus $f_{jk}$ is given by
\begin{align*}
z_j^\alpha=f_{jk}^\alpha\left(z_k^1,...,z_k^n,t_1,...,t_m\right), \alpha=1,...,n.
\end{align*}
On the other hand, we may assume that the locally free sheaf $\Theta_\mathcal{F}$ is trivialized in the covering $\mathcal{U}$ so that 
\begin{align*}
\Gamma\left(\mathcal{U}_j,\Theta_\mathcal{F}\right)&\cong \bigoplus^p \Gamma\left(\mathcal{U}_j, \mathcal{O}_\mathcal{M}\right)\\
T_j^\alpha(z_j,t) &\mapsto e_j^\alpha=\left(0,...,\overbrace{1}^{\alpha-\textnormal{th}},...0\right),\alpha=1,...,p
\end{align*}
which is generated by vector fields $T_j^\alpha(z_j,t)\in \Gamma\left(\mathcal{U}_j, \Theta_\mathcal{F}\right),\alpha=1,...,p$ of the form
\begin{align}\label{ta3}
T_j^\alpha(z_j,t):=\sum_{\beta=1}^n T_j^{\alpha\beta}(z_j,t)\frac{\partial}{\partial z_j^\beta},\,\,\,\,\,\,\,\,\,\alpha=1,...,p
\end{align}
for some $T_j^{\alpha\beta}(z_j,t)\in \Gamma\left( \mathcal{U}_j , \mathcal{O}_\mathcal{M} \right)$, and satisfies
\begin{align}\label{ta5}
\left[T_j^\alpha(z_j,t), T_j^\beta(z_j,t)\right]=\sum_{\gamma=1}^n g_{j\alpha\beta}^\gamma(z_j,t) T_j^\gamma(z_j,t)
\end{align}
for some $g_{j\alpha\beta}^\gamma(z_j,t)\in \Gamma\left(\mathcal{U}_j, \mathcal{O}_\mathcal{M}\right)$ with $g_{j\alpha\beta}^\gamma(z_j,t)=-g_{j\beta\alpha}^\gamma(z_j,t)$, and for $j,k$ with $\mathcal{U}_j\cap \mathcal{U}_k \ne \emptyset$, we have invertible $p\times p$ matrices $R_{jk}:=\left(r_{jk}^{\alpha\beta}(z_k,t)\right)$ with components  $r_{jk}^{\alpha\beta} (z_k,t)\in \Gamma\left(\mathcal{U}_j\cap \mathcal{U}_k, \mathcal{O}_\mathcal{M}\right)$ such that
\begin{center}
$\left[\begin{matrix}
T_j^1(z_j,t)\\
T_j^2(z_j,t)\\
\cdot\\
\cdot\\
\cdot\\
T_j^p(z_j,t)
\end{matrix}\right]
=\left[\begin{matrix}
r_{jk}^{11}(z_k,t) & \cdots & r_{jk}^{1p}(z_k,t)\\
r_{jk}^{21}(z_k,t) & \cdots & r_{jk}^{2n}(z_k,t)\\
\cdot & \cdots & \cdot \\
\cdot & \cdots & \cdot\\
\cdot & \cdots & \cdot\\
r_{jk}^{p1}(z_k,t) & \cdots & r_{jk}^{pp}(z_k,t)
\end{matrix}\right]
\cdot \left[\begin{matrix}
T_k^1(z_k,t)\\
T_k^2(z_k,t)\\
\cdot\\
\cdot\\
\cdot\\
T_k^p(z_k,t)
\end{matrix}\right]$
\end{center}
and for $i,j,k\in \mathcal{U}_i\cap \mathcal{U}_j\cap \mathcal{U}_k\ne \emptyset$, we have
\begin{align*}
R_{ik}=R_{ij}\cdot R_{jk}\,\,\,\,\,\textnormal{on}\,\,\,\,\,\mathcal{U}_i\cap \mathcal{U}_j\cap\mathcal{U}_k.
\end{align*}

More specifically, $T_j^\alpha(z_j,t)=\sum_{\beta=1}^p r_{jk}^{\alpha\beta}(z_k,t) T_k^\beta(z_k,t)$ implies that
\begin{align*}
\sum_{\gamma=1}^n T_j^{\alpha\gamma}(z_j,t)\frac{\partial}{\partial z_j^\gamma}=\sum_{\beta=1}^p\sum_{\eta=1}^n r_{jk}^{\alpha\beta}(z_k,t) T_k^{\beta\eta}(z_k,t)\frac{\partial}{\partial z_k^\eta},
\end{align*}
so that
\begin{align}\label{ta4}
\sum_{\gamma=1}^n T_j^{\alpha\gamma}(f_{jk}(z_k,t),t)\frac{\partial}{\partial z_j^\gamma}=\sum_{\beta=1}^p\sum_{\eta=1}^n\sum_{\gamma=1}^n r_{jk}^{\alpha\beta}(z_k,t) T_k^{\beta\eta}(z_k,t)\frac{\partial f_{jk}^\gamma(z_k,t)}{\partial z_k^\eta}\frac{\partial}{\partial z_j^\gamma}
\end{align}

\end{remark}

\subsection{Leaf complex controlling deformations of singular holomorphic foliations in terms of locally free subsheaves of tangent sheaf}\

Let $\left(\mathcal{M}, \Theta_\mathcal{F}, B, \pi \right)$ be a foliated analytic family with $\Theta_\mathcal{F}$ locally free as in Definition \ref{t3}, so that $\Theta_{\mathcal{F}_t}$ defines a (singular) holomorphic foliation on a compact complex manifold $M_t$ with $\Theta_{\mathcal{F}_t}$ locally free. Then we have the leaf complex on $M_t$ associated to $\Theta_{\mathcal{F}_t}$ (see Appendix \ref{AppendixA1})
\begin{align*}
\Theta_{\mathcal{F}_t}^\bullet:\Theta_{M_t}\xrightarrow{D_0^t} \mathscr{H}om_{\mathcal{O}_{M_t}}\left(\Theta_{\mathcal{F}_t}, \frac{\Theta_{M_t}}{\Theta_{\mathcal{F}_t}} \right)\xrightarrow{D_1^t} \mathscr{H}om_{\mathcal{O}_{M_t}}\left(\bigwedge^2 \Theta_{\mathcal{F}_t}, \frac{\Theta_{M_t}}{\Theta_{\mathcal{F}_t} }\right) \xrightarrow{D_2^t} \mathscr{H}om_{\mathcal{O}_{M_t}}\left(\bigwedge^3 \Theta_{\mathcal{F}_t}, \frac{\Theta_{M_t}}{\Theta_{\mathcal{F}_t} } \right)\xrightarrow{D_3^t} \cdots
\end{align*}
We will denote the $i$-th cohomology group by $\mathbb{H}^1\left( M_t , \Theta_{\mathcal{F}_t}^\bullet \right)$. We can compute $\mathbb{H}^i\left(M_t, \Theta_{\mathcal{F}_t}^\bullet \right)$  by the following \v Cech resolution of $\Theta_{\mathcal{F}_t}^\bullet$. Here $\delta$ is the \v Cech map and $\mathcal{U}_t= \mathcal{U}\cap M_t=\left\{ U_j^t :=\mathcal{U}_j \cap M_t | j=1,2,... \right\}$ is an open covering of $M_t$:
{\small{\begin{center}
$\begin{CD}
\cdots \\
@AD_2^t AA \\
C^0\left(\mathcal{U}_t, \mathscr{H}om_{\mathcal{O}_{M_t}}\left( \bigwedge^2  \Theta_{\mathcal{F}_t}, \frac{\Theta_{M_t}}{\Theta_{\mathcal{F}_t} } \right) \right) @>-\delta>> \cdots  \\
@AD_1^tAA @AD_1^t AA \\
C^0\left(\mathcal{U}_t, \mathscr{H}om_{\mathcal{O}_{M_t}}\left(\Theta_{\mathcal{F}_t}, \frac{\Theta_{M_t}}{\Theta_{\mathcal{F}_t} } \right) \right)@>\delta>> C^1\left(\mathcal{U}_t, \mathscr{H}om_{\mathcal{O}_{M_t}}\left(\Theta_{\mathcal{F}_t}, \frac{\Theta_{M_t}}{\Theta_{\mathcal{F}_t}}  \right) \right) @>-\delta>> \cdots \\
@AD_0^t AA @AD_0^t AA @AD_0^t AA \\
C^0(\mathcal{U}_t, \Theta_{M_t}) @>-\delta>> C^1(\mathcal{U}_t, \Theta_{M_t}) @>\delta>> C^2(\mathcal{U}_t, \Theta_{M_t}) @>>> \cdots 
\end{CD}$
\end{center}}}

We will  relate the first cohomology group $\mathbb{H}^1\left( M_t, \Theta_{\mathcal{F}_t}^\bullet \right)$ to infinitesimal foliated deformations of $\pi^{-1}(t)=\left(M_t, \Theta_{\mathcal{F}_t} \right)$ in the foliated analytic family $\left(\mathcal{M}, \Theta_\mathcal{F}, \pi, B \right)$ in terms of tangent sheaves.

\subsection{Infinitesimal foliated deformations in terms of tangent sheaves}\label{ta11}\

Let $\left(\mathcal{M}, \Theta_\mathcal{F}, B, \pi  \right)$ be a foliated analytic family with $\Theta_{\mathcal{F}}$ locally free as in Definition \ref{t3}. We keep the notations in Remark \ref{ta1}. By taking the derivative of $(\ref{ta3} )$ with respect to $t$, we set
\begin{align*}
\frac{\partial T_j^\alpha(z_j,t)}{\partial t}:=\sum_{\beta=1}^n \frac{\partial T_j^{\alpha\beta}(z_j,t)}{\partial t}\frac{\partial}{\partial z_j^\beta},\,\,\,\,\,\,\,\,\,\alpha=1,...,p.
\end{align*}
We define an element $\alpha_j(t)$ in $\Gamma\left(U_j^t, \mathscr{H}om_{\mathcal{O}_{M_t}}\left(\Theta_{\mathcal{F}_t}, \frac{\Theta_{M_t}}{\Theta_{\mathcal{F}_t}}  \right) \right)$ in the following way:
\begin{align*}
\alpha_j(t):\Gamma\left(U_j^t, \Theta_{\mathcal{F}_t}\right)&\to \Gamma\left(U_j^t,\frac{\Theta_{M_t}}{\Theta_{\mathcal{F}_t}} \right) \\
                                                         T_j^\alpha(z_j,t)&\mapsto -\overline{\frac{\partial T_j^\alpha(z_j,t)}{\partial t}}
\end{align*}
where $\overline{\frac{T_j^\alpha(z_j,t)}{\partial t}}$ denotes the image of $\frac{\partial T_j^\alpha(z_j,t)}{\partial t}$ in the natural map $\Theta_{M_t}\to \Theta_{M_t}/\Theta_{\mathcal{F}_t}$, and linearly extends to $\Gamma\left( U_j^t, \Theta_{\mathcal{F}_t} \right)$. We also define $\tilde{\alpha}_j(t)\in \Gamma\left(U_j^t, \mathscr{H}om_{\mathcal{O}_{M_t}}\left( \Theta_{\mathcal{F}_t}, \Theta_{M_t} \right) \right)$ in the following way:
\begin{align*}
\tilde{\alpha}_j(t):\Gamma\left(U_j^t, \Theta_{\mathcal{F}_t}\right)&\to \Gamma\left(U_j^t,\Theta_{M_t}  \right) \\
                                                         T_j^\alpha(z_j,t)&\mapsto -\frac{\partial T_j^\alpha(z_j,t)}{\partial t}
\end{align*}
and linearly extends to $\Gamma\left( U_j^t, \Theta_{\mathcal{F}_t} \right)$. Then we have
\begin{proposition}[compare \cite{GM88} p.63]\label{ta2}
{\Small{\begin{align*}
\left(\left\{\theta_{jk}(t):=\sum_{\alpha=1}^n \frac{\partial f_{jk}^\alpha(z_k,t)}{\partial t}\frac{\partial}{\partial z_j^\alpha} \right\}, \left\{ \alpha_j(t):=\left(T_j^\alpha(z_j,t)\mapsto -\overline{\frac{\partial T_j^\alpha(z_j,t)}{\partial t}}\right)_{\alpha=1,...,p} \right\} \right) \in C^1\left(\mathcal{U}^t, \Theta_{M_t} \right) \bigoplus C^0\left( \mathcal{U}_t, \mathscr{H}om_{\mathcal{O}_{M_t}}\left(\Theta_{\mathcal{F}_t}, \frac{\Theta_{M_t}}{\Theta_{\mathcal{F}_t}} \right) \right)
\end{align*}}}
defines a $1$-cocycle in the above \v Cech resolution of $\Theta_{\mathcal{F}_t}^\bullet$ and call its cohomology class in $\mathbb{H}^1\left( M_t, \Theta_{\mathcal{F}_t}^\bullet \right)$ the infinitesimal $($foliated$)$ deformation along $\frac{\partial}{\partial t}$. This expression is independent of the choice of local coordinates.
\end{proposition}

\begin{proof}
First we note that $\delta\left( \left\{ \theta_{jk}(t)\right\} \right)=0$ (see \cite{Kod05} p.201). Second, by taking the derivative of $(\ref{ta5})$ with respect to $t$, 
\begin{align*}
&\left[\frac{\partial T_j^\alpha}{\partial t}, T_j^\beta\right] + \left[ T_j^\alpha, \frac{\partial T_j^\beta}{\partial t} \right] =\sum_{\gamma=1}^p \frac{\partial g_{j\alpha\beta}^\gamma}{\partial t} T_j^\gamma + \sum_{\gamma=1}^p g_{j\alpha\beta}^\gamma \frac{\partial T_j^\gamma}{\partial t} \\
&\Longrightarrow \left[ T_j^\alpha, \tilde{\alpha}_j(t)\left( T_j^\beta \right) \right] - \left[ T_j^\beta, \tilde{\alpha}_j(t)\left( T_j^\alpha \right) \right] - \tilde{\alpha}_j(t)\left(\left[ T_j^\alpha, T_j^\beta \right] \right) = - \sum_{\gamma=1}^p  \frac{\partial g_{j\alpha\beta}^\gamma}{\partial t} T_j^\gamma \in \Theta_{\mathcal{F}_t}
\end{align*}
This implies that $D_1^t\left(\alpha_j(t) \right)=0$. It remains to show that $\delta\left(\{\alpha_j(t)\} \right) + D_0^t \left(\left\{ \theta_{jk} \right\} \right)=0$, equivalently $\left(\alpha_k(t)- \alpha_j(t)\right)(T_j^\alpha) + \overline{ \left[ \theta_{jk}(t) , T_j^\alpha \right]} =0$ for $\alpha=1,...,p$. In fact, by taking the derivative of $(\ref{ta4})$ with respect to $t$ and considering the coefficient of $\frac{\partial}{\partial z_j^\gamma}$, we have
{\small{\begin{align}\label{ta6}
\sum_{\eta=1}^n\frac{\partial T_j^{\alpha\gamma}}{\partial z_j^\eta}\frac{\partial f_{jk}^\eta}{\partial t} +\frac{\partial T_j^{\alpha\gamma}}{\partial t} = \sum_{\beta=1}^p\sum_{\eta=1}^n \frac{\partial r_{jk}^{\alpha\beta}}{\partial t} T_k^{\beta\eta}\frac{\partial f_{jk}^\gamma}{\partial z_k^\eta} +\sum_{\beta=1}^p\sum_{\eta=1}^n r_{jk}^{\alpha\beta} \frac{\partial T_k^{\beta\eta}}{\partial t} \frac{\partial f_{jk}^\gamma}{\partial z_k^\eta} +\sum_{\beta=1}^p\sum_{\eta=1}^n r_{jk}^{\alpha\beta}T_k^{\beta\eta}\frac{\partial}{\partial z_k^\eta}\left(\frac{\partial f_{jk}^\gamma}{\partial t} \right)
\end{align}}}

Since $T_j^\alpha(z_j,t)=\sum_{\beta=1}^p r_{jk}^{\alpha\beta}(z_k,t)T_k^\beta(z_k,t)$, we have
\begin{align}\label{ta7}
\left(\tilde{\alpha}_k(t)- \tilde{\alpha}_j(t)\right)\left(T_j^\alpha\right)&=-\sum_{\beta=1}^p r_{jk}^{\alpha\beta} \frac{\partial T_k^\beta}{\partial t}+\frac{\partial T_j^\alpha}{\partial t}
                                                       =-\sum_{\beta=1}^p\sum_{\eta=1}^n r_{jk}^{\alpha\beta} \frac{\partial T_k^{\beta\eta}}{\partial t}\frac{\partial}{\partial z_k^\eta} +\sum_{\gamma=1}^n \frac{\partial T_j^\alpha}{\partial t}\frac{\partial}{\partial z_j^\gamma}\\
                                                       &=-\sum_{\beta=1}^p\sum_{\eta=1}^n r_{jk}^{\alpha\beta}\frac{\partial T_k^{\beta\eta}}{\partial t}\frac{\partial f_{jk}^\gamma}{\partial z_k^\eta}\frac{\partial}{\partial z_j^\gamma}+\sum_{\gamma=1}^n \frac{\partial T_j^{\alpha\gamma} }{\partial t}\frac{\partial}{\partial z_j^\gamma} \notag
\end{align}

Let us compute
\begin{align}\label{ta8}
\left[ \theta_{jk}(t), T_j^\alpha\right]  = \left[ \sum_{\gamma=1}^n \frac{\partial f_{jk}^\gamma}{\partial t}\frac{\partial}{\partial z_j^\gamma} , \sum_{\beta=1}^n T_j^{\alpha\beta}\frac{\partial}{\partial z_j^\beta} \right] =\sum_{\gamma,\beta=1}^n\frac{\partial f_{jk}^\gamma}{\partial t}\frac{\partial T_j^{\alpha\beta}}{\partial z_j^\gamma}\frac{\partial}{\partial z_j^\beta}-\sum_{\gamma,\beta=1}^n T_j^{\alpha\beta} \frac{\partial}{\partial z_j^\beta}\left( \frac{\partial f_{jk}^\gamma}{\partial t} \right)\frac{\partial}{\partial z_j^\gamma}
\end{align}

Then $(\ref{ta6}), (\ref{ta7})$ and $(\ref{ta8})$ implies that 
\begin{align*}
\left(\tilde{\alpha}_k(t)- \tilde{\alpha}_j(t)\right)\left(T_j^\alpha\right) + \left[ \theta_{jk}(t) ,  T_j^\alpha \right] = \sum_{\gamma=1}^n\sum_{\beta=1}^p \sum_{\eta=1}^n \frac{\partial r_{jk}^{\alpha\beta}}{\partial t} T_k^{\beta \eta} \frac{\partial f_{jk}^\gamma}{\partial z_k^\eta} \frac{\partial}{\partial z_j^\gamma}=\sum_{\beta=1}^p \frac{\partial r_{jk}^{\alpha\beta}}{\partial t} T_k^\beta \in \Theta_{\mathcal{F}_t}
\end{align*}
Hence we have $\alpha_k(t)-\alpha_j(t)+D_0^t\left(\theta_{jk}(t)\right)=0$. Next we show that $\left(\left\{ \theta_{jk}(t)\right\}, \left\{ \alpha_j(t)\right\} \right)$ is independent of the choice of system of local coordinates. We can show that the infinitesimal deformation does not change under the refinement of the open covering (see \cite{Kod05} p.190). Since we can choose a common refinement for two systems of local coordinates, it is sufficient to show that given two local coordinates $x_j=(z_j,t)$ and $u_j=(y_j,t)$ on each $\mathcal{U}_j$, the infinitesimal foliated deformation $\left(\left\{\eta_{jk}(t) \right\}, \left\{ \sigma_j(t) \right\} \right)$ with respect to $\{u_j\}$ coincides with $\left( \left\{\theta_{jk} (t) \right\} , \left\{ \alpha_j (t) \right\} \right)$ with respect to $\{x_j\}$. Let $\Theta_{\mathcal{F}}$ be generated by $V_j^\alpha(y_j,t)=\sum_{\gamma=1}^n V_j^{\alpha\gamma}(y_j,t)\frac{\partial}{\partial y_j^\gamma},\alpha=1,...,p$ in terms of coordinates $u_j=(y_j,t)$ on $\mathcal{U}_j$. Let $\left(y_k,t \right)\to \left( y_j, t \right)= \left(e_{jk}(y_k,t),t \right)$ be the coordinate transformation of $\{u_j\}$ on $\mathcal{U}_j\cap \mathcal{U}_k\ne \emptyset$. Then we have $\eta_{jk}(t)=\sum_{\alpha=1}^n \frac{\partial e_{jk}^\alpha(y_k,t)}{\partial t}\frac{\partial}{\partial y_j^\alpha}, y_k= e_{kj}(y_j,t)$, and 
\begin{align*}
\sigma_j(t):\Gamma\left(U_j^t,  \Theta_{\mathcal{F}_t} \right) \to \Gamma\left( U_j^t, \frac{\Theta_{M_t}}{\Theta_{\mathcal{F}_t}} \right), \,\,\,\,\,& V_j^\alpha \mapsto - \overline{\frac{\partial V_j^\alpha}{\partial t}}:= - \overline{\sum_{\gamma=1}^n \frac{\partial V_j^{\alpha \gamma}(y_j,t)}{\partial t} \frac{\partial}{\partial y_j^\gamma}}\\
\tilde{\sigma}_j(t):\Gamma\left(U_j^t,  \Theta_{\mathcal{F}_t} \right) \to \Gamma\left( U_j^t, \Theta_{M_t} \right),\,\,\,\,\,& V_j^\alpha \mapsto -\frac{\partial V_j^\alpha}{\partial t}:= -\sum_{\gamma=1}^n \frac{\partial V_j^{\alpha \gamma}(y_j,t)}{\partial t} \frac{\partial}{\partial y_j^\gamma}
\end{align*}
We show that $\left(\left\{ \theta_{jk}(t) \right\}, \left\{ \alpha_j(t)\right\} \right)$ is cohomologous to $\left(\left\{ \eta_{jk}(t)\right\}, \left\{ \sigma_j(t)\right\} \right)$. Let $y_j^\alpha=g_j^\alpha\left(z_j^1,..., z_j^n, t \right),\alpha=1,...,n$, define the coordinate transformation from $x_j=(z_j,t)$ to $u_j=(y_j,t)$. Since $\left\{T_j^\alpha\right\}$ and $\left\{V_j^\alpha \right\},\alpha=1,...,n$ are two bases of $\Gamma\left(\mathcal{U}_j, \Theta_\mathcal{F} \right)$, we can write $V_j^\alpha(y_j,t)= \sum_{\beta=1}^p a_j^{\alpha\beta}(z_j,t) T_j^\beta(z_j,t)$ for some $a_j^{\alpha \gamma}(z_j,t)\in \Gamma\left(\mathcal{U}_j, \mathcal{O}_{\mathcal{M}} \right)$, i.e. $\sum_{\gamma=1}^n V_j^{\alpha\gamma}(y_j,t)\frac{\partial}{\partial y_j^\gamma} = \sum_{\beta=1}^p \sum_{\eta=1}^n a_j^{\alpha\beta}(z_j,t) T_j^{\beta \eta}(z_j,t)\frac{\partial}{\partial z_j^\eta}$, so that by setting $g_j(z_j,t)=\left(g_j^1(z_j,t),..., g_j^n(z_j,t) \right)$, we have
\begin{align}\label{ta9}
\sum_{\gamma=1}^n V_j^{\alpha\gamma}\left( g_j(z_j,t),t \right)\frac{\partial}{\partial y_j^\gamma}= \sum_{\beta=1}^p \sum_{\gamma=1}^n\sum_{\eta=1}^n a_j^{\alpha\beta}(z_j,t) T_j^{\beta \eta}(z_j,t)\frac{\partial g_j^\gamma(z_j,t)}{\partial z_j^\eta}\frac{\partial}{\partial y_j^\gamma}
\end{align}

We set $\theta_j(t)=\sum_{\alpha=1}^n \frac{\partial g_j^\alpha(z_j,t)}{\partial t} \frac{\partial}{\partial y_j^\alpha}, y_j^\alpha= g_j^\alpha(z_j,t)$. Then $-\delta\left( \left\{\theta_j(t) \right\} \right)= \left\{ \eta_{jk}(t)- \theta_{jk}(t)\right\}$ (for the detail, see \cite{Kod05} p.191-192), and we claim that $\overline{\left[ \theta_j(t) ,V_j^\alpha \right]}=\sigma_j(t)(V_j^\alpha)- \alpha_j(t)\left( V_j^\alpha \right)$. In fact, by taking the derivative of $(\ref{ta9})$ with respect to $t$, and considering the coefficient of $\frac{\partial}{\partial y_j^\gamma}$, we have
{\small{\begin{align}\label{ta10}
\sum_{\eta=1}^n\frac{\partial V_j^{\alpha\gamma}}{\partial y_j^\eta}\frac{\partial g_{j}^\eta}{\partial t} +\frac{\partial V_j^{\alpha\gamma}}{\partial t} = \sum_{\beta=1}^p\sum_{\eta=1}^n \frac{\partial a_{j}^{\alpha\beta}}{\partial t} T_j^{\beta\eta}\frac{\partial g_{j}^\gamma}{\partial z_j^\eta} +\sum_{\beta=1}^p\sum_{\eta=1}^n a_{j}^{\alpha\beta} \frac{\partial T_j^{\beta\eta}}{\partial t} \frac{\partial g_{j}^\gamma}{\partial z_j^\eta} +\sum_{\beta=1}^p\sum_{\eta=1}^n a_{j}^{\alpha\beta}T_j^{\beta\eta}\frac{\partial}{\partial z_j^\eta}\left(\frac{\partial g_{j}^\gamma}{\partial t} \right)
\end{align}}}
Then since $V_j^\alpha(y_j,t)= \sum_{\beta=1}^p a_j^{\alpha\beta}(z_j,t) T_j^\beta(z_j,t)$, we have
{\small{\begin{align*}
&\tilde{\alpha}_j(t)\left( V_j^\alpha \right) - \tilde{\sigma}_j(t)\left( V_j^\alpha \right) + \left[ \theta_j(t), V_j^\alpha \right] = - \sum_{\beta=1}^p a_j^{\alpha\beta}\frac{\partial T_j^\beta}{\partial t} + \frac{\partial V_j^\alpha}{\partial t} +  \left[ \sum_{\eta=1}^n \frac{\partial g_j^\eta}{\partial t } \frac{\partial}{\partial y_j^\eta} , \sum_{\gamma=1}^n V_j^{\alpha \gamma} \frac{\partial}{\partial y_j^\gamma} \right] \\
&= - \sum_{\beta=1}^p\sum_{\eta,\gamma=1}^n a_j^{\alpha\beta}\frac{\partial T_j^{\beta \eta}}{\partial t} \frac{\partial g_j^\gamma}{\partial z_j^\eta}\frac{\partial}{\partial y_j^\gamma} + \sum_{\gamma=1}^n \frac{\partial V_j^{\alpha \gamma}}{\partial t} \frac{\partial}{\partial y_j^\gamma} + \sum_{\eta,\gamma=1}^n \frac{\partial g_j^\eta}{\partial t} \frac{\partial V_j^{\alpha \gamma}}{\partial y_j^\eta} \frac{\partial}{\partial y_j^\gamma} - \sum_{\eta, \gamma=1}^n  V_j^{\alpha \gamma} \frac{\partial}{\partial y_j^\gamma}\left( \frac{\partial g_j^\eta}{\partial t} \right)\frac{\partial}{\partial y_j^\eta} \\
&= \sum_{\beta=1}^p \sum_{\eta,\gamma=1}^n \frac{\partial a_j^{\alpha \beta}}{\partial t} T_j^{\beta \eta} \frac{\partial g_j^\gamma}{\partial z_j^\eta}\frac{\partial}{\partial z_j^\gamma} = \sum_{\beta=1}^p \frac{\partial a_j^{\alpha \beta}}{\partial t} T_j^\beta \in \Theta_{\mathcal{F}_t}
\end{align*}}}
This implies that $\alpha_j(t)\left( V_j^\alpha \right) - \sigma_j(t)\left( V_j^\alpha \right) + \overline{\left[ \theta_j(t), V_j^\alpha \right]}=0$. Hence $\left(\left\{ \theta_{jk}(t) \right\}, \left\{ \alpha_j(t)\right\} \right)$ is cohomologous to $\left(\left\{ \eta_{jk}(t)\right\}, \left\{ \sigma_j(t)\right\} \right)$.

\end{proof}

\begin{definition}[foliated Kodaira-Spencer map in terms of tangent sheaf]\label{ta12}
Let $(\mathcal{M}, \Theta_\mathcal{F}, B, \pi)$ with $\Theta_{\mathcal{F}}$ locally free be a foliated complex analytic family of deformations of $(M_t, \Theta_{\mathcal{F}_t})=\pi^{-1}(t),t\in B$, where $B$ is a domain of $\mathbb{C}^m$. We keep the notations in the proof of \textnormal{Proposition \ref{ta2}}. For a tangent vector $\frac{\partial}{\partial t}=\sum_{\lambda=1}^m c_\lambda \frac{\partial}{\partial t_\lambda}, c_\lambda \in \mathbb{C}$, of $B$, we put
\begin{align*}
\frac{\partial \Theta_{\mathcal{F}_t}}{\partial t}:=\left\{\left(  T_j^\alpha(z_j,t) \mapsto - \overline{\sum_{\lambda=1}^m c_\lambda \frac{\partial T_j^\alpha(z_j,t)}{\partial t_\lambda}}     \right)_{\alpha=1,...,p} \right\} \in C^0\left( \mathcal{U}_t, \mathscr{H}om_{\mathcal{O}_{M_t}}\left(\Theta_{\mathcal{F}_t}, \frac{\Theta_{M_t}}{\Theta_{\mathcal{F}_t}} \right) \right)
\end{align*}

Then the foliated Kodaira-Spencer map in terms of tangent sheaf at $t$ is defined to be a $\mathbb{C}$-linear map
{\small{\begin{align*}
\varphi_t:T_t(B)&\to \mathbb{H}^1\left(M_t, \Theta_{\mathcal{F}_t}^\bullet \right)\\
                \frac{\partial}{\partial t}&\mapsto \frac{\partial \left(M_t, \Theta_{\mathcal{F}_t}\right)}{\partial t}:=\left( \rho_t\left(\frac{\partial}{\partial t} \right)= \frac{\partial M_t}{\partial t}= \left\{  \sum_{\alpha=1}^n \frac{\partial f_{jk}^\alpha(z_k,t)}{\partial t}\frac{\partial}{\partial z_j^\alpha} \right\}, \left\{ \left(T_j^\alpha(z_j,t)\mapsto -\overline{\frac{\partial T_j^\alpha(z_j,t)}{\partial t}}\right)_{\alpha=1,...,p} \right\} \right) \\
                 &\,\,\,\,\,\,\,\,\,\,\in C^1\left(\mathcal{U}_t, \Theta_{M_t} \right) \bigoplus C^0\left( \mathcal{U}_t, \mathscr{H}om_{\mathcal{O}_{M_t}}\left(\Theta_{\mathcal{F}_t}, \frac{\Theta_{M_t}}{\Theta_{\mathcal{F}_t}} \right) \right)
\end{align*}}}
where $\rho_t:T_t(B) \to H^1\left( M_t, \Theta_{M_t}\right)$ is the Kodaira-Spencer map at $t$ of the underlying complex analytic family $\left( \mathcal{M}, B, \pi \right)$ $($see \cite{Kod05} \textnormal{p.201}$)$.
\end{definition}

\section{Theorem of existence for deformations of foliated complex analytic structures in terms of tangent sheaves}

\subsection{Preliminaries}\label{tt1}\

Let $\left(\mathcal{M}, \Theta_\mathcal{F}, B, \pi \right)$ with $\Theta_\mathcal{F}$ locally free be a foliated analytic family in terms of tangent sheaves, where $B$ is a domain of $\mathbb{C}^m$ containing the origin $0$ as in Definition \ref{t3}. Define $|t|=\max_\lambda |t_\lambda|$ for $t=(t_1,...,t_m)\in \mathbb{C}^m$, and let $\Delta=\{t\in \mathbb{C}^m | |t|<r\}$ the polydisk  of radius $r>0$. If we take a sufficiently  small $\Delta \subset B$, then $\left(\mathcal{M}_\Delta, \Theta_{\mathcal{F}_\Delta} \right)=\pi^{-1}\left(\Delta\right)$ is represented in the form
\begin{align*}
\left(\mathcal{M}_\Delta, \Theta_{\mathcal{F}_\Delta}\right)=\bigcup_j \left(U_j\times \Delta, \Theta_\mathcal{F}|_{U_j\times \Delta}\right)
\end{align*}
We denote a point of $U_j$ by $\xi_j=\left(\xi_j^1,...,\xi_j^n\right)$ and $\Gamma\left(U_j\times \Delta, \Theta_\mathcal{F}\right)$ is generated by $T_j^1\left(\xi_j,t\right),...,T_j^p\left(\xi_j,t\right)$, where
\begin{align*}
T_j^\alpha\left(\xi_j,t\right)=\sum_{\beta=1}^n T_j^{\alpha\beta}\left(\xi_j, t\right) \frac{\partial}{\partial \xi_j^\beta}
\end{align*}
for some $T_j^{\alpha\beta}(\xi_j,t)\in \Gamma\left( U_j\times \Delta , \mathcal{O}_\mathcal{M} \right)$, and satisfies
\begin{align}\label{tai5}
\left[T_j^\alpha(\xi_j,t), T_j^\beta(\xi_j,t)\right]=\sum_{\gamma=1}^n g_{j\alpha\beta}^\gamma(\xi_j,t) T_j^\gamma(\xi_j,t)
\end{align}
for some $g_{j\alpha\beta}^\gamma(\xi_j,t)\in \Gamma\left(U_j\times \Delta,  \mathcal{O}_\mathcal{M}\right)$ with $g_{j\alpha\beta}^\gamma(\xi_j,t)=-g_{j\beta\alpha}^\gamma(\xi_j,t)$.

For simplicity, we assume that $U_j=\left\{\xi_j\in \mathbb{C}^m | |\xi_j|<1\right\}$ where $|\xi_j |=\max_a |\xi_j^a |$. $(\xi_j,t)\in U_j\times \Delta$ and $(\xi_k, t)\in U_k\times \Delta$ are the same point on $\mathcal{M}_\Delta$ if $\xi_j^\alpha=f_{jk}^\alpha(\xi_k,t),\alpha=1,...,n$ where $f_{jk}(\xi_k,t)$ is a holomorphic map of $\xi_k^1,..., \xi_k^n, t_1,...,t_m$, defined on $\left(U_k\times \Delta\right) \cap \left(U_j\times \Delta\right)$, and we have the following relation
\begin{align*}
T_j^\alpha(\xi_j,t)=\sum_{\beta=1}^p r_{jk}^{\alpha\beta}(\xi_k,t) T_k^\beta(\xi_k,t),\,\,\,\,\,\alpha=1,...,p
\end{align*}
for a holomorphic function $r_{jk}^{\alpha\beta}(\xi_k,t)$ on $\left(U_k\times \Delta\right) \cap \left( U_j\times \Delta\right)$. More precisely
\begin{align}
\sum_{\alpha,\delta=1}^n T_j^{\alpha\delta}\left(f_{jk}(\xi_k,t),t \right)\frac{\partial}{\partial \xi_j^\delta } =\sum_{\beta=1}^p \sum_{\gamma,\delta=1}^n r_{jk}^{\alpha\beta}(\xi_k,t) T_k^{\beta\gamma}(\xi_k,t)\frac{\partial f_{jk}^\delta(\xi_j,t)}{\partial \xi_k^\gamma}\frac{\partial}{\partial \xi_j^\delta},\,\,\,\,\,\alpha,\beta=1,...,p
\end{align}
Equivalently
\begin{align} \label{oi5}
T_j^{\alpha\delta}\left(f_{jk}(\xi_k,t),t\right)= \sum_{\beta=1}^p \sum_{\gamma=1}^n r_{jk}^{\alpha\beta}(\xi_k,t) T_k^{\beta\gamma}(\xi_k,t)\frac{\partial f_{jk}^{\delta}(\xi_k,t)}{\partial \xi_k^\gamma}
\end{align}
and on $\left(U_i\times \Delta\right)\cap \left( U_j\times \Delta \right)\cap \left( U_k\times \Delta \right)$, we have 
\begin{align*}
r_{ik}^{\alpha \gamma}(\xi_k,t)= \sum_{\beta=1}^p r_{ij}^{\alpha\beta}\left(f_{jk}(\xi_k,t),t) \right) r_{jk}^{\beta \gamma}(\xi_k,t),\,\,\,\,\,\alpha,\gamma=1,...,p
\end{align*}

By \cite{Kod05} Theorem 2.3, $M_t$ is diffeomorphic to $M_0=\pi^{-1}(0)$ as differentiable manifolds for each $t\in \Delta$. We put $M:=M_0$. By \cite{Kod05} Theorem 2.5, if we take a sufficiently small $\Delta$, there is a diffeomorphism $\Psi$ of $M\times \Delta$ onto $\mathcal{M}_\Delta$ as differentiable manifolds such that $\pi\circ \Psi$ is the projection $M\times \Delta \to \Delta$. Let $z=\left(z^1,..., z^n \right)$ be local coordinates of $M=M_0$. Then we have $\pi\circ \Psi(z,t)=t, t\in \Delta$. For $\Psi(z,t)\in U_j\times \Delta$, put
\begin{align*}
\Psi\left( z, t\right) = \left(\xi_j^1(z,t), ..., \xi_j^n(z,t), t_1,..., t_m \right)
\end{align*}
Then each component $\xi_j^\alpha=\xi_j^\alpha(z,t),\alpha=1,...,n$ is a $C^\infty$ function. If we identify $\mathcal{M}_\Delta =\Psi\left(M\times \Delta \right)$ with $M\times \Delta$ via $\Psi$, $\left(\mathcal{M}_\Delta,  \Theta_{\mathcal{F}_\Delta} \right)$ is considered as a complex manifold with the complex structure defined on the $C^\infty$ manifold $M\times \Delta$ by the system of local complex coordinates on $U_j\times \Delta$
\begin{align}\label{te1}
\left\{ \left(\xi_j,t \right)| j=1 ,2, 3,...   \right\},\,\,\,\,\,(\xi_j,t)=\left( \xi_j^1(z,t), ..., \xi_j^n(z,t), t_1,..., t_m \right)
\end{align}
and holomorphic foliation $\Theta_{\mathcal{F}_\Delta}$ on $U_j\times \Delta$ with respect to the coordinates $(\ref{te1})$ generated by $T_j^\alpha\left(\xi_j,t \right)= \sum_{\beta=1}^n T_j^{\alpha\beta}(\xi_j,t)\frac{\partial}{\partial \xi_j^\beta},\alpha=1,...,p$. 
We note that since $(z^1,..., z^n)$ and $\left(\xi_j^1(z,0), ..., \xi_j^n(z,0)\right)$ are local complex coordinates on $M=M_0$, $\xi_j^\alpha(z,0)$ are holomorphic functions of $z^1,...,z^n,\alpha=1,...,n$. We also note that if we take $\Delta$ sufficiently small, we have
\begin{align}\label{oi3}
\det\left( \frac{\partial \xi_j^\alpha(z,t)}{\partial z^\lambda} \right)_{\alpha, \lambda=1,...,n} \ne 0\,\,\,\,\,\, \, \textnormal{for}\,\,\,t\in \Delta
\end{align}

We put $\mathscr{U}_j=\Psi^{-1}\left(U_j\times \Delta \right)$. Then $\mathscr{U}_j \subset M\times \Delta$ is the domain of $\xi_j^{\alpha}(z,t)$. From $(\ref{oi3})$, we can define a $(0,1)$-form $\varphi_j^\lambda(z,t)= \sum_{v=1}^n \varphi_{jv}^\lambda(z,t) d\bar{z}^v$ for each $\lambda=1,...,n $ such that $\bar{\partial}\xi_j^\alpha =  \sum_{\lambda, v=1}^n \varphi_{jv}^\lambda(z,t)\frac{\partial \xi_j^\alpha}{\partial z^\lambda}d\bar{z}^v$. Then on $\mathscr{U}_j \cap \mathscr{U}_k$, we have $\sum_{\lambda=1}^n \varphi_j^\lambda(z,t)\frac{\partial}{\partial z^\lambda} = \sum_{\lambda=1}^n \varphi_k^\lambda(z,t) \frac{\partial}{\partial z^\lambda}$ (see \cite{Kod05} p.262). Then if for $(z,t)\in \mathscr{U}_j$, we define $\varphi(z,t):=\sum_{\lambda=1}^n \varphi_j^\lambda(z,t)\frac{\partial}{\partial z^\lambda}$, then $\varphi(t):=\varphi(z,t)\in A^{0,1}\left(M, \Theta_M \right)$ for $t\in \Delta$, and satisfies
\begin{align}\label{oi23}
\varphi(0)=0,\,\,\,\,\,\,\,\bar{\partial}\varphi(t)- \frac{1}{2}\left[ \varphi(t), \varphi(t) \right]=0
\end{align}
(see \cite{Kod05} p.263, p.265), We also point out that
\begin{theorem}\label{oi13}
If we take a sufficiently small polydisk $\Delta$ as above, then for $t\in \Delta$, a local $C^\infty$ function $f$ on $M$ is holomorphic with respect to the complex structure $M_t$ if and only if
\begin{align*}
\left( \bar{\partial} - \bar{\partial} \varphi(t)\right)f=0
\end{align*} 
\end{theorem}
\begin{proof}
See \cite{Kod05} Theorem 5.3 p.263
\end{proof}

\begin{theorem}
If we take a sufficiently small polydisk $\Delta$ as above, then for $t\in \Delta$, a $C^\infty$ vector field $T=\sum_{\alpha=1}^n T^\alpha(z)\frac{\partial}{\partial z^\alpha}$ on $M$ defines a holomorphic vector field $\sum_{\beta=1}^n T^\alpha(z)\frac{\partial \xi_j^\beta}{\partial z^\alpha}\frac{\partial}{\partial \xi_j^\beta}$ on $M_t$ with respect to the complex structure $M_t$ induced by $\varphi(t)$ if and only if it satisfies the equation
\begin{align*}
\bar{\partial} T- \left[ \varphi(t), T \right]=0
\end{align*}
\end{theorem}

\begin{proof}
By Theorem $\ref{oi23}$, it suffices to show that $\bar{\partial}T- \left[\varphi(t),T \right]=0$ if and only if for each $\beta=1,...,n$,
\begin{align}
&\bar{\partial}\left(\sum_{\alpha=1}^n T^\alpha(z)\frac{\partial \xi_j^\beta}{\partial z^\alpha} \right)- \left[ \sum_{\lambda, vf=1}^n \varphi_{v}^\lambda  d\bar{z}^v \frac{\partial}{\partial z^\lambda}, \sum_{\alpha=1}^n T^\alpha(z)\frac{\partial \xi_j^\beta}{\partial z^\alpha} \right] \notag \\
&=\sum_{v,\alpha=1}^n \frac{\partial T^\alpha}{\partial \bar{z}^v}\frac{\partial \xi_j^\beta}{\partial z^\alpha}d\bar{z}^v+\sum_{v,\alpha=1}^n T^\alpha \frac{\partial}{\partial \bar{z}^v}\left(\frac{\partial \xi_j^\beta}{\partial z^\alpha} \right) d\bar{z}^v-\sum_{\lambda, v,\alpha=1}^n \varphi_{v}^\lambda \frac{\partial T^\alpha}{\partial z^\lambda}\frac{\partial \xi_j^\beta}{\partial z^\alpha} d\bar{z}^v  - \sum_{\lambda, v,\alpha=1}^n \varphi_v^\lambda T^\alpha \frac{\partial^2 \xi_j^\beta}{\partial z^\lambda \partial z^\alpha} d\bar{z}^v=0 \notag \\
&\iff    \sum_{\alpha=1}^n \frac{\partial T^\alpha}{\partial \bar{z}^v}\frac{\partial \xi_j^\beta}{\partial z^\alpha}+\sum_{\alpha=1}^n T^\alpha \frac{\partial}{\partial z^\alpha}\left(\frac{\partial \xi_j^\beta}{\partial \bar{z}^v} \right)-\sum_{\lambda, \alpha=1}^n \varphi_{v}^\lambda\frac{\partial T^\alpha}{\partial z^\lambda}\frac{\partial \xi_j^\beta}{\partial z^\alpha}-\sum_{\lambda, \alpha=1}^n \varphi_v^\lambda T^\alpha \frac{\partial^2 \xi_j^\beta}{\partial z^\lambda \partial z^\alpha} =0\,\,\,\,\,\,\,\,\text{for each}\,\,\, v,\beta=1,...,n \label{te2}
\end{align}
Since $\frac{\partial \xi_j^\beta}{\partial \bar{z}^v}=\sum_{\lambda=1}^n \frac{\partial \xi_j^\beta}{\partial z^\lambda}\varphi_v^\lambda$, $(\ref{te2})$ is equivalent to that for $v,\beta=1,...,n$,
\begin{align}
& \sum_{\alpha=1}^n \frac{\partial T^\alpha}{\partial \bar{z}^v}\frac{\partial \xi_j^\beta}{\partial z^\alpha}+\sum_{\alpha,\lambda=1}^n T^\alpha \frac{\partial}{\partial z^\alpha}\left(\frac{\partial \xi_j^\beta}{\partial z^\lambda} \varphi_v^\lambda\right)-\sum_{\lambda, \alpha=1}^n \varphi_{v}^\lambda\frac{\partial T^\alpha}{\partial z^\lambda}\frac{\partial \xi_j^\beta}{\partial z^\alpha}-\sum_{\lambda, \alpha=1}^n \varphi_v^\lambda T^\alpha \frac{\partial^2 \xi_j^\beta}{\partial z^\lambda \partial z^\alpha}=0 \notag\\
& \iff  \sum_{\alpha=1}^n \frac{\partial T^\alpha}{\partial \bar{z}^v}\frac{\partial \xi_j^\beta}{\partial z^\alpha}+\sum_{\alpha,\lambda=1}^n T^\lambda \frac{\partial \varphi_v^\alpha}{\partial z^\lambda}\frac{\partial \xi_j^\beta}{\partial z^\alpha}-\sum_{\lambda, \alpha=1}^n \varphi_{v}^\lambda\frac{\partial T^\alpha}{\partial z^\lambda}\frac{\partial \xi_j^\beta}{\partial z^\alpha}=0 \notag \label{ta121}\\
&\iff  \frac{\partial T^\alpha}{\partial \bar{z}^v}+\sum_{\lambda=1}^n T^\lambda \frac{\partial \varphi_v^\alpha}{\partial z^\lambda}-\sum_{\lambda=1}^n \varphi_v^\lambda \frac{\partial T^\alpha}{\partial z^\lambda}=0 \,\,\,\,\,\,\,\, \textnormal{from}\,\,\,(\ref{oi3})
\end{align}

On the other hand, let us compute
\begin{align*}
\bar{\partial} T- \left[\varphi(t),T \right]&=\bar{\partial}\left(\sum_{\alpha=1}^n T^\alpha \frac{\partial}{\partial z_\alpha} \right) -\left[\sum_{\lambda, v=1}^n \varphi_{v}^\lambda d\bar{z}^v\frac{\partial}{\partial z^\lambda}, \sum_{\alpha=1}^n T^\alpha \frac{\partial}{\partial z^\alpha} \right]\\
 & = \sum_{\alpha,v=1}^n \frac{\partial T^\alpha}{\partial \bar{z}^v}d \bar{z}^v\wedge \frac{\partial}{\partial z^\alpha}-\sum_{\lambda, \alpha, v=1}^n d\bar{z}^v\wedge\left( \varphi_v^\lambda \frac{\partial T^\alpha}{\partial z^\lambda}\frac{\partial}{\partial z^\alpha} - T^\alpha \frac{\partial \varphi_v^\lambda}{\partial z^\alpha} \frac{\partial}{\partial z^\lambda} \right)=0
\end{align*}
which is equivalent to (\ref{ta121}).
\end{proof}

With this preparation, we shall prove theorem of existence for deformations of foliated complex analytic structure in terms of tangent sheaves.

\subsection{Theorem of existence of deformations of foliated complex analytic structures in terms of tangent sheaves}

\begin{theorem}[Theorem of existence for deformations of foliated analytic structures in terms of tangent sheaves]\label{tt2}
Let $\left( M , \Theta_{\mathcal{F}_0} \right)$ be a compact foliated complex manifold with $\Theta_{\mathcal{F}_0}$ locally free. Suppose that $\mathbb{H}^2\left( M,  \Theta_{\mathcal{F}_0}^\bullet \right)=0$. Then there exists a foliated analytic family $\left( \mathcal{M}, \Theta_\mathcal{F}, B, \pi \right)$ with $0\in B \subset \mathbb{C}^m$ satisfying the following conditions$:$
\begin{enumerate}
\item $\pi^{-1}(0)=\left( M, \Theta_{\mathcal{F}_0} \right)$
\item The foliated Kodaira-Spencer map $\varphi_0: T_0 B \to \mathbb{H}^1\left( M, \Theta_{\mathcal{F}_0}^\bullet \right)$ in terms of tangent sheaf is an isomorphism.
\end{enumerate}
\end{theorem}

\begin{proof}
We may assume the following:
\begin{enumerate}
\item $M$ is covered by a finite number of coordinate neighborhoods $U_i(i\in I)$ with a system of coordinates $z_i=\left(z_i^1,...,z_i^n\right)$ such that $U_i=\left\{z_i\in \mathbb{C}^n |\max_{1\leq \alpha \leq n} |z_i^\alpha|<1\right\}$ \label{d1}
\item $\Theta_{\mathcal{F}_0}$ is generated by $\left\{T_{0i}^1,...,T_{0i}^p \right\}$ with $T_{0i}^\alpha=\sum_{\beta=1}^q T_{0i}^{\alpha\beta}(z_i) \frac{\partial}{\partial z_i^\alpha}\in \Gamma\left(U_i, \Theta_{\mathcal{F}_0}\right)$ and $T_{0j}^\alpha=\sum_{\beta=1}^p r_{0ij}^{\alpha\beta} T_{0j}^\beta$ on $U_i\cap U_j$ for $r_{0ij}^{\alpha\beta}(z_j)\in \Gamma\left(U_i\cap U_j, \mathcal{O}_M\right),\alpha=1,...,p$.
\item $\left[T_{0i}^\alpha, T_{0i}^\beta \right]= \sum_{\gamma=1}^p g_{0i\alpha\beta}^\gamma T_{0i}^\gamma$ where $g_{i\alpha\beta}^\gamma (z_i) \in \Gamma(U_i, \mathcal{O}_M)$ with $g_{0i\alpha\beta}^\gamma=- g_{0i\beta \alpha}^\gamma $ for $\alpha,\beta=1,...,p$.
\item $z_i$ coincides with $z_j$ if and only if $z_i=f_{ij}(z_j)$.
\end{enumerate}
Then we note that from $\sum_{\gamma=1}^p r_{0ik}^{\alpha\gamma}T_{0k}^\gamma=T_{0i}^\alpha=\sum_{\beta=1}^p r_{0ij}^{\alpha\beta} T_{0j}^\beta=\sum_{\beta,\gamma=1}^p r_{0ij}^{\alpha\beta}r_{0jk}^{\beta\gamma} T_{0k}^\gamma$, we have
\begin{align*}
r_{0ik}^{\alpha\gamma}=\sum_{\beta=1}^p r_{0ij}^{\alpha\beta} r_{0jk}^{\beta\gamma}
\end{align*}
and we note that
{\small{\begin{align*}
\sum_{\gamma,\xi=1}^p g_{0i\alpha\beta}^\gamma r_{0ij}^{\gamma\xi} T_{0j}^\xi&=\sum_{\gamma=1}^p g_{0i\alpha\beta}^\gamma T_{0i}^\gamma=\left[T_{0i}^\alpha, T_{0i}^\beta \right]= \sum_{\xi,\eta=1}^p \left[r_{0ij}^{\alpha\xi} T_{0j}^\xi, r_{0ij}^{\beta \eta}T_{0j}^\eta \right]\\
 &=\sum_{\xi,\eta=1}^p r_{0ij}^{\alpha\xi}r_{0ij}^{\beta\eta}\left[T_{0j}^\xi, T_{0j}^\eta \right] +\sum_{\xi,\eta=1}^p r_{0ij}^{\alpha\xi} T_{0j}^\xi \left(r_{0ij}^{\beta\eta}\right)T_{0j}^\eta-\sum_{\xi,\eta=1}^p r_{0ij}^{\beta\eta} T_{0j}^\eta \left(r_{0ij}^{\alpha\xi} \right) T_{0j}^\xi\\
 &=\sum_{\xi,\eta,\gamma=1}^p r_{0ij}^{\alpha\xi} r_{0ij}^{\beta\eta} g_{0j\xi\eta}^\gamma T_{0j}^\gamma +\sum_{\xi,\eta=1}^p r_{0ij}^{\alpha\xi} T_{0j}^\xi\left(r_{0ij}^{\beta\eta}\right)T_{0j}^\eta-\sum_{\xi,\eta=1}^p r_{0ij}^{\beta\eta} T_{0j}^\eta \left(r_{0ij}^{\alpha\xi}\right) T_{0j}^\xi\\
 &=\sum_{\xi,\eta,\gamma=1}^p r_{0ij}^{\alpha\eta} r_{0ij}^{\beta\gamma} g_{0j\eta\gamma}^\xi T_{0j}^\xi +\sum_{\xi,\eta=1}^p r_{0ij}^{\alpha\eta} T_{0j}^\eta \left(r_{0ij}^{\beta\xi} \right)T_{0j}^\xi-\sum_{\xi,\eta=1}^p r_{0ij}^{\beta\eta} T_{0j}^\eta \left(r_{0ij}^{\alpha\xi}  \right) T_{0j}^\xi 
\end{align*}}}
Hence we have
\begin{align}\label{te3}
\sum_{\xi=1}^p\left(\sum_{\eta,\gamma=1}^p r_{0ij}^{\alpha\eta} r_{0ij}^{\beta \gamma} g_{0j\eta\gamma}^\xi+\sum_{\eta=1}^p r_{0ij}^{\alpha\eta} T_{0j}^\eta\left(r_{0ij}^{\beta\xi}\right) -\sum_{\eta=1}^p r_{0ij}^{\beta\eta} T_{0j}^\eta\left( r_{0ij}^{\alpha\xi}\right)-\sum_{\gamma=1}^p g_{0i\alpha\beta}^\gamma r_{0ij}^{\gamma\xi}\right) T_{0j}^\xi=0
\end{align}
Since $\Theta_{\mathcal{F}_0}$ is locally free, $(\ref{te3})$ is equivalent to
\begin{align}\label{te19}
\sum_{\eta,\gamma=1}^p r_{0ij}^{\alpha\eta} r_{0ij}^{\beta \gamma} g_{0j\eta\gamma}^\xi+T_{0i}^\alpha \left(r_{0ij}^{\beta\xi}\right) -T_{0i}^\beta \left( r_{0ij}^{\alpha\xi}\right)-\sum_{\gamma=1}^p g_{0i\alpha\beta}^\gamma r_{0ij}^{\gamma\xi}=0
\end{align}
On the other hand, we note that
{\small{\begin{align*}
0&=\left[ \left[T_{0i}^\alpha,T_{0i}^\beta\right]-\sum_{\gamma=1}^p g_{0i\alpha\beta}^\gamma T_{0i}^\gamma, T_{0i}^\delta \right]=\left[\left[T_{0i}^\alpha, T_{0i}^\delta\right], T_{0i}^\beta \right] +\left[ T_{0i}^\alpha, \left[ T_{0i}^\beta, T_{0i}^\delta \right]\right] -\sum_{\gamma=1}^p g_{0i\alpha\beta}^\gamma \left[ T_{0i}^\gamma, T_{0i}^\delta\right]+\sum_{\gamma=1}^p T_{0i}^\delta\left( g_{0i\alpha\beta}^\gamma \right) T_{0i}^\gamma\\
&=\left[\sum_{\gamma=1}^p g_{0i\alpha\delta}^\gamma T_{0i}^\gamma, T_{0i}^\beta \right] +\left[T_{0i}^\alpha, \sum_{\gamma=1}^p g_{0i\beta \delta}^\gamma T_{0i}^\gamma \right]-\sum_{\gamma,\xi=1}^p g_{0i\alpha\beta}^\gamma g_{0i\gamma\delta}^\xi T_{0i}^\xi +\sum_{\gamma=1}^p T_{0i}^\delta\left( g_{0i\alpha\beta}^\gamma \right) T_{0i}^\gamma\\
&= \sum_{\gamma,\xi=1}^p g_{0i\alpha\delta}^\gamma g_{0i\gamma \beta}^\xi T_{0i}^\xi -\sum_{\gamma=1}^p T_{0i}^\beta\left( g_{0i\alpha\delta}^\gamma \right) T_{0i}^\gamma +\sum_{\gamma,\xi=1}^p g_{0i\beta\delta}^\gamma g_{0i\alpha\gamma}^\xi T_{0i}^\xi+ \sum_{\gamma=1}^p T_{0i}^\alpha\left( g_{0i\beta\delta}^\gamma \right) T_{0i}^\gamma -\sum_{\gamma,\xi=1}^p g_{0i\alpha\beta}^\gamma g_{0i\gamma\delta}^\xi T_{0i}^\xi +\sum_{\gamma=1}^p T_{0i}^\delta\left( g_{0i\alpha\beta}^\gamma \right) T_{0i}^\gamma\\
&=\sum_{\xi=1}^p\left(\sum_{\gamma=1}^p g_{0i\alpha\delta}^\gamma g_{0i\gamma\beta}^\xi - T_{0i}^\beta\left( g_{0i\alpha\delta}^\xi\right) +\sum_{\gamma=1}^p g_{0i\beta\delta}^\gamma g_{0i\alpha\gamma}^\xi +T_{0i}^\alpha\left( g_{0i\beta\delta}^\xi\right) -\sum_{\gamma=1}^p g_{0i\alpha\beta}^\gamma g_{0i\gamma \delta}^\xi + T_{0i}^\delta\left( g_{0i\alpha\beta}^\xi \right)\right) T_{0i}^\xi
\end{align*}}}

Hence we have
\begin{align}\label{te25}
\sum_{\gamma=1}^p g_{0i\alpha\delta}^\gamma g_{0i\gamma\beta}^\xi - T_{0i}^\beta\left( g_{0i\alpha\delta}^\xi\right) +\sum_{\gamma=1}^p g_{0i\beta\delta}^\gamma g_{0i\alpha\gamma}^\xi +T_{0i}^\alpha\left( g_{0i\beta\delta}^\xi\right) -\sum_{\gamma=1}^p g_{0i\alpha\beta}^\gamma g_{0i\gamma \delta}^\xi + T_{0i}^\delta\left( g_{0i\alpha\beta}^\xi \right)=0
\end{align}

Let $r=\dim_\mathbb{C} \mathbb{H}^1\left(M, \Theta_{\mathcal{F}_0}^\bullet \right)$ and $B=\left\{t= \left(t_1,...,t_r \right)\in \mathbb{C}^r||t|<\epsilon \right\}$ with sufficiently small $\epsilon >0$. Our purpose is to construct
\begin{enumerate}
\item a $C^\infty$ $(0,1)$-form $\varphi(t)$ with coefficients in $\Theta_M$ depending holomorphically on $t$,
\item $C^\infty$-vector fields $T_i^\alpha\left(z_i,t\right)$ on $U_i\times B$ of the form $T_i^\alpha(z_i,t)=\sum_{\beta=1}^n T_i^{\alpha\beta}(z_i,t)\frac{\partial}{\partial z_i^\beta}$ depending holomorphically on $t$,\\
\item $p\times p$ matrices $R_{ij}:=\left(r_{ij}^{\alpha\beta}(z_j,t)\right)$  with components $C^\infty$-function $r_{ij}^{\alpha\beta}(z_j,t)$ on $U_{ij}\times B$ depending holomorphically on $t$, such that
\end{enumerate}
\begin{align}
\varphi(0)=0,\,\,\,\,\,\,\,\,\,\,\,\,T_i^\alpha(z_i,0)=T_{0i}^\alpha,\,\,\,\,\,&\,\,\,\,\,r_{ij}^{\alpha\beta}(z_j,0)=r_{0ij}^{\alpha\beta}(z_j),\,\,\,\,\,\,\,\,\,\,\,\,\,\,g_{i\alpha\beta}^\gamma(z_i,0)=g_{0i\alpha\beta}^\gamma(z_i)\label{t10}
\end{align}
\begin{align}
\bar{\partial}\varphi&-\frac{1}{2}\left[\varphi,\varphi \right]=0, \label{t12}
\end{align}
\begin{align}
\bar{\partial} T_i^\alpha(z_i,t)&- \left[\varphi, T_i^\alpha(z_i,t) \right]=0, \label{ua1}
\end{align}
\begin{align}
\bar{\partial}\left(r_{ij}^{\alpha\beta}(z_i,t)\right)&- \left[\varphi, r_{ij}^{\alpha\beta}(z_j,t)\right]=0, \label{ua2}
\end{align}
\begin{align}
\bar{\partial}\left( g_{i\alpha\beta}^\gamma (z_i,t) \right)&- \left[ \varphi, g_{i\alpha\beta}^\gamma(z_i,t)\right]=0,\label{ua3}
\end{align}
\begin{align}
T_i^\alpha(z_i,t)&=\sum_{\beta=1}^p r_{ij}^{\alpha\beta}(z_j,t) T_j^\beta(z_j,t),\label{ua4}
\end{align}
\begin{align}
\left[T_i^\alpha(z_i,t), T_i^\beta(z_i,t)\right]&=\sum_{\gamma=1}^p g_{i\alpha\beta}^\gamma(z_i,t) T_i^\gamma(z_i,t), \label{ua5}
\end{align}
\begin{align}
R_{ij}R_{jk}=R_{ik}\iff r_{ik}^{\alpha\gamma}(z_k,t)&=\sum_{\beta=1}^p r_{ij}^{\alpha\beta}\left(f_{jk}(z_k),t\right) r_{jk}^{\beta\gamma}(z_k,t)\,\,\,\,\,\,\,\textnormal{on}\,\,\,\,\,\left(U_i\cap U_j\cap U_k \right)\times M \label{t11}
\end{align}
\begin{align}
\sum_{\eta,\gamma=1}^p r_{ij}^{\alpha\eta} r_{ij}^{\beta \gamma} g_{j\eta\gamma}^\xi+T_i^\alpha &\left(r_{ij}^{\beta\xi}\right) -T_i^\beta \left( r_{ij}^{\alpha\xi}\right)-\sum_{\gamma=1}^p g_{i\alpha\beta}^\gamma r_{ij}^{\gamma\xi}=0  \label{te4}
\end{align}
\begin{align}
\sum_{\gamma=1}^p g_{i\alpha\delta}^\gamma g_{i\gamma\beta}^\xi - T_i^\beta\left( g_{i\alpha\delta}^\xi\right)  +\sum_{\gamma=1}^p & g_{i\beta\delta}^\gamma g_{i\alpha\gamma}^\xi +T_i^\alpha\left( g_{i\beta\delta}^\xi\right) -\sum_{\gamma=1}^p g_{i\alpha\beta}^\gamma g_{i\gamma \delta}^\xi + T_i^\delta\left( g_{i\alpha\beta}^\xi \right)=0 \label{te5}
\end{align}

\subsection{Existence of formal solutions}\

We construct solutions of $(\ref{t10})-(\ref{te5})$ which are formal power series in $t$. 

\begin{notation}\label{te27}
We write the power series expansion of a holomorphic function $P(t)$ in $t_1,..,t_m$ defined on a neighborhood of the origin $0$ in the form$:$ $P(t)=P(0)+ P_1(t)+ \cdots + P_\mu(t)+ \cdots $, where each $P_\mu(t)$ denotes a homogenous polynomial of degree $\mu$ in $t_1,..., t_m$. We set $P^\mu(t):=P(0)+ P_1(t)+ \cdots + P_\mu(t)$. For any power series $P(t)$ and $Q(t)$ in $t=(t_1,..., t_m)$, we indicate $P(t)\equiv_\mu Q(t)$ contains no term of degree $\leq \mu$ in $t$. Then $P^\mu=0$ means $P\equiv_\mu 0$.
\end{notation}

We write $\varphi=\sum_{\mu=0}^\infty \varphi_\mu, T_i^\alpha= \sum_{\mu=0}^\infty T_{i|\mu}^\alpha, r_{ij}^{\alpha\beta}= \sum_{\mu=0}^\infty r_{ij|\mu}^{\alpha\beta}$ and $g_{i\alpha\beta}^{\gamma}=\sum_{\mu=0}^\infty g_{i\alpha\beta|\mu}^\gamma$. In view of $(\ref{t10})$, we set $\varphi_0=0, T_{i|0}^\alpha=T_{0i}^\alpha(z_i)$, $r_{ij|0}^{\alpha\beta}(z_i,t)=r_{0ij}^{\alpha\beta}(z_j)$ and $g_{i\alpha\beta|0}^{\gamma}(z_i,t)=g_{0i\alpha\beta}^\gamma(z_i)$. Then $(\ref{t12})-(\ref{te5})$ are equivalent to the following systems of congruences:
\begin{align}
\bar{\partial}\varphi^{\mu}&-\frac{1}{2}\left[\varphi^\mu, \varphi^\mu\right]\equiv_{\mu} 0 \label{t13}
\end{align}
\begin{align}
\bar{\partial} T_i^{\alpha\mu}&-\left[\varphi^\mu, T_i^{\alpha\mu}\right]\equiv_\mu 0 \label{tpc6}
\end{align}
\begin{align}
\bar{\partial} r_{ij}^{\alpha\beta\mu}&- \left[\varphi^\mu, r_{ij}^{\alpha\beta\mu}\right]\equiv_\mu 0 \label{tpc5}
\end{align}
\begin{align}
\bar{\partial} g_{i\alpha\beta}^{\gamma\mu}&- \left[\varphi^\mu, g_{i\alpha\beta}^{\gamma\mu}\right]\equiv_\mu  0 \label{tpc7}
\end{align}
\begin{align}
T_i^{\alpha\mu}&\equiv_\mu \sum_{\beta=1}^p r_{ij}^{\alpha\beta\mu} T_j^{\beta\mu}\label{tpc4}
\end{align}
\begin{align}
\left[T_i^{\alpha\mu}, T_i^{\beta\mu}\right]&\equiv_\mu \sum_{\gamma=1}^p g_{i\alpha\beta}^{\gamma\mu} T_i^{\gamma\mu}  \label{tpc2}
\end{align}
\begin{align}
r_{ik}^{\alpha\gamma \mu} (z_k,t) &\equiv_\mu \sum_{\beta=1}^p r_{ij}^{\alpha\beta \mu}\left(f_{jk},t \right) r_{jk}^{\beta \gamma \mu}(z_k,t) \label{tpc8}
\end{align}
\begin{align}
\sum_{\eta,\gamma=1}^p r_{ij}^{\alpha\eta\mu} r_{ij}^{\beta \gamma\mu} g_{j\eta\gamma}^{\xi\mu}&+T_i^{\alpha\mu} \left(r_{ij}^{\beta\xi\mu}\right) -T_i^{\beta\mu} \left( r_{ij}^{\alpha\xi \mu}\right)-\sum_{\gamma=1}^p g_{i\alpha\beta}^{\gamma \mu} r_{ij}^{\gamma\xi \mu}\equiv_\mu 0 \label{tpc3}
\end{align}
\begin{align}
\sum_{\gamma=1}^p g_{i\alpha\delta}^{\gamma \mu} g_{i\gamma\beta}^{\xi\mu} - T_i^{\beta\mu}\left( g_{i\alpha\delta}^{\xi\mu}\right)& +\sum_{\gamma=1}^p g_{i\beta\delta}^{\gamma\mu} g_{i\alpha\gamma}^{\xi\mu} +T_i^{\alpha\mu}\left( g_{i\beta\delta}^{\xi\mu}\right) -\sum_{\gamma=1}^p g_{i\alpha\beta}^{\gamma\mu} g_{i\gamma \delta}^{\xi\mu} + T_i^{\delta\mu}\left( g_{i\alpha\beta}^{\xi\mu} \right)\equiv_\mu 0 \label{te6}
\end{align}
for $\mu=1,.2,3,\cdots$.

We construct solutions of $\left(\ref{t13} \right)_\mu- (\ref{te6})_\mu$ by induction on $\mu$. We suppose that $\varphi^{\mu-1}, T_i^{\alpha (\mu-1)}, r_{ij}^{\alpha\beta(\mu-1)}$ and $g_{i\alpha\beta}^{\gamma(\mu-1)}$ satisfying $(\ref{t13})_{\mu-1}-(\ref{te6})_{\mu-1}$ are already determined.

We define homogeneous polynomials 
\begin{align}\label{tt3}
\xi_\mu  &\equiv_\mu \bar{\partial} \varphi^{\mu-1} -\frac{1}{2}\left[\varphi^{\mu-1},\varphi^{\mu-1}\right]
\end{align}
\begin{align}
-\Phi_{i|\mu}^\alpha &\equiv_\mu \bar{\partial} T_i^{\alpha (\mu-1)}- \left[\varphi^{\mu-1}, T_i^{\alpha(\mu-1)}\right] \label{tt7}
\end{align}
\begin{align}\label{ds6}
-\sum_{\eta=1}^p \Lambda_{ij|\mu}^{\alpha\eta} r_{0ij}^{\eta\beta} &\equiv_\mu \bar{\partial} r_{ij}^{\alpha\beta(\mu-1)} - \left[ \varphi^{\mu-1}, r_{ij}^{\alpha\beta(\mu-1)} \right]
\end{align}
\begin{align}
-\eta_{i\alpha\beta|\mu}^{\gamma}&\equiv_\mu \bar{\partial} g_{i\alpha\beta}^{\gamma(\mu-1)}- \left[\varphi^{\mu-1}, g_{i\alpha\beta}^{\gamma(\mu-1)} \right]
\end{align}
\begin{align}\label{tt6}
\Gamma_{ij|\mu}^\alpha &\equiv_\mu T_i^{\alpha(\mu-1)}-\sum_{\beta=1}^p r_{ij}^{\alpha\beta(\mu-1)} T_j^{\beta(\mu-1)}
\end{align}
\begin{align}
\Pi_{i|\mu}^{\alpha\beta}&\equiv_\mu \left[T_i^{\alpha(\mu-1)},T_i^{\beta(\mu-1)}\right]-\sum_{\gamma=1}^p g_{i\alpha\beta}^{\gamma(\mu-1)}T_i^{\gamma(\mu-1)} \label{tt121}
\end{align}
\begin{align}
\sum_{\beta=1}^p \lambda_{ijk|\mu}^{\alpha\beta} r_{0ik}^{\beta\gamma}&\equiv_\mu r_{ik}^{\alpha\gamma(\mu-1)} (z_k,t) -\sum_{\beta=1}^p r_{ij}^{\alpha\beta(\mu-1)}(f_{jk},t) r_{jk}^{\beta \gamma(\mu-1)}(z_k,t)
\end{align}
\begin{align}\label{ds47}
\Psi_{ij|\mu}^{\alpha\beta\xi}\equiv_\mu \sum_{\eta,\gamma=1}^p r_{ij}^{\alpha\eta(\mu-1)}& r_{ij}^{\beta \gamma(\mu-1)} g_{j\eta\gamma}^{\xi(\mu-1)}+ T_i^{\alpha(\mu-1)}\left(r_{ij}^{\beta\xi(\mu-1)}\right) -T_i^{\beta(\mu-1)}\left( r_{ij}^{\alpha\xi(\mu-1)}\right)-\sum_{\gamma=1}^p g_{i\alpha\beta}^{\gamma(\mu-1)} r_{ij}^{\gamma\xi(\mu-1)}
\end{align}
{\small{\begin{align}\label{tt4}
G_{i|\mu}^{\alpha\beta \delta}\equiv_\mu \sum_{\gamma=1}^p g_{i\alpha\delta}^{\gamma (\mu-1)} g_{i\gamma\beta}^{\xi(\mu-1)} - T_i^{\beta(\mu-1)}\left( g_{i\alpha\delta}^{\xi(\mu-1)}\right)+\sum_{\gamma=1}^p g_{i\beta\delta}^{\gamma(\mu-1)} g_{i\alpha\gamma}^{\xi(\mu-1)} +T_i^{\alpha(\mu-1)}\left( g_{i\beta\delta}^{\xi(\mu-1)}\right) -\sum_{\gamma=1}^p g_{i\alpha\beta}^{\gamma(\mu-1)} g_{i\gamma \delta}^{\xi(\mu-1)} + T_i^{\delta(\mu-1)}\left( g_{i\alpha\beta}^{\xi(\mu-1)} \right)
\end{align}}}

\begin{lemma}\label{te10}
We have the following equalities$:$
\begin{align}
\bar{\partial}\xi_\mu=0 \label{t19}
\end{align}
\begin{align}
\sum_{\beta,\eta=1}^p r_{0ij}^{\alpha\beta}\lambda_{jkl|\mu}^{\beta\eta}r_{0ji}^{\eta\xi}&-\lambda_{ikl|\mu}^{\alpha\xi}+\lambda_{ijl|\mu}^{\alpha\xi}-\lambda_{ijk|\mu}^{\alpha\xi}=0\label{t20}
\end{align}
\begin{align}
\bar{\partial} &\Phi_{i|\mu}^\alpha=\left[ \xi_\mu,  T_{0i}^\alpha\right] \label{t22}
\end{align}
\begin{align}
\bar{\partial} \Lambda_{ij|\mu}^{\alpha\xi} &=\sum_{\beta=1}^p \left[\xi_\mu, r_{0ij}^{\alpha\beta} \right] r_{0ji}^{\beta\xi} \label{t23}
\end{align}
\begin{align}
\sum_{\eta,\gamma=1}^p r_{0ij}^{\alpha\eta}\Lambda_{jk|\mu}^{\eta\gamma}r_{0ji}^{\gamma\xi}&-\Lambda_{ik|\mu}^{\alpha\xi}+ \Lambda_{ij|\mu}^{\alpha\xi}=\bar{\partial} \lambda_{ijk|\mu}^{\alpha\xi} \label{t24}
\end{align}
\begin{align}
\sum_{\gamma=1}^p r_{0ij}^{\alpha\gamma} \Gamma_{jk|\mu}^\gamma -\Gamma_{ik|\mu}^\alpha  & + \Gamma_{ij|\mu}^\alpha=\sum_{\gamma=1}^p \lambda_{ijk|\mu}^{\alpha\gamma} T_{0i}^\gamma \label{t25}
\end{align}
\begin{align}
\sum_{\beta=1}^p r_{0ij}^{\alpha\beta} \Phi_{j|\mu}^\beta-\Phi_{i|\mu}^\alpha&=\bar{\partial} \Gamma_{ij|\mu}^\alpha  - \sum_{\eta=1}^p \Lambda_{ij|\mu}^{\alpha\eta} T_{0i}^\eta \label{t26} 
\end{align}
\begin{align}
\Pi_{i|\mu}^{\alpha\beta}-\sum_{\eta,\xi=1}^p r_{0ij}^{\alpha\eta} r_{0ij}^{\beta \xi} \Pi_{j|\mu}^{\eta\xi}&= \left[\Gamma_{ij|\mu}^\alpha, T_{0i}^\beta \right]+\left[ T_{0i}^\alpha, \Gamma_{ij|\mu}^\beta \right]-\Gamma_{ij|\mu}\left(\left[T_{0i}^\alpha, T_{0i}^\beta \right] \right) +\sum_{\xi=1}^p \left(\Psi_{ij|\mu}^{\alpha\beta\xi} -\Gamma_{ij|\mu}^\alpha\left(r_{0ij}^{\beta \xi} \right) + \Gamma_{ij|\mu}^\beta \left( r_{0ij}^{\alpha \xi} \right)  \right)T_{0j}^\xi \label{t27}
\end{align}
\begin{align}
\bar{\partial} \Pi_{i|\mu}^{\alpha\beta}= -\left[ \Phi_{i|\mu}^\alpha, T_{0i}^\beta \right]& -\left[ T_{0i}^\alpha, \Phi_{i|\mu}^\beta \right] +\sum_{\gamma=1}^p g_{0i\alpha\beta}^\gamma \Phi_{i|\mu}^\gamma  +\sum_{\gamma=1}^p \eta_{i\alpha\beta|\mu}^\gamma T_{0i}^\gamma \label{t28}
\end{align}
\begin{align} \label{te7}
\left[ \Pi_{i|\mu}^{\alpha \beta}, T_{0i}^\delta \right]=\left[ \Pi_{i|\mu}^{\alpha \delta}, T_{0i}^\beta \right]  +\sum_{\gamma=1}^p g_{0i\alpha\delta}^\gamma \Pi_{i|\mu}^{\gamma \beta}  +\left[ T_{0i}^\alpha, \Pi_{i|\mu}^{\beta\delta} \right]  +\sum_{\gamma=1}^p g_{0i\beta\delta}^\gamma \Pi_{i|\mu}^{\alpha \gamma} - \sum_{\gamma=1}^p g_{0i\alpha\beta}^\gamma \Pi_{i|\mu}^{\gamma \delta} + \sum_{\xi=1}^p G_{i|\mu}^{\alpha \beta \delta \xi} T_{0i}^\xi 
\end{align}
\end{lemma}

\begin{proof}
$(\ref{t19})$ follows from \cite{Kod05} p.273. We prove $(\ref{t20})$. In fact,
{\small{\begin{align*}
&\sum_{\beta=1}^p  \lambda_{ijl|\mu}^{\alpha\beta}r_{0il}^{\beta\gamma}\equiv_\mu  r_{il}^{\alpha\gamma(\mu-1)}(z_l,t)-\sum_{\beta=1}^p r_{ij}^{\alpha\beta(\mu-1)}\left(f_{jl},t\right)r_{jl}^{\beta\gamma(\mu-1)}(z_l,t) \\
&\equiv_\mu \sum_{\beta=1}^p  \lambda_{ikl|\mu}^{\alpha\beta} r_{0il}^{\beta\gamma} +\sum_{\beta=1}^p r_{ik}^{\alpha\beta(\mu-1)}\left(f_{kl},t \right) r_{kl}^{\beta\gamma(\mu-1)}(z_l,t) -\sum_{\beta=1}^p r_{ij}^{\alpha\beta(\mu-1)}\left(f_{jl},t \right)\left(\sum_{\eta=1}^p r_{jk}^{\beta\eta(\mu-1)}\left(f_{kl},t \right)r_{kl}^{\eta \gamma(\mu-1)}(z_l,t) +\sum_{\eta=1}^p \lambda_{jkl|\mu}^{\beta\eta}r_{0jl}^{\eta\gamma}\right)\\
&\equiv_\mu \sum_{\beta=1}^p  \lambda_{ikl|\mu}^{\alpha\beta} r_{0il}^{\beta\gamma} +\sum_{\beta=1}^p r_{ik}^{\alpha\beta(\mu-1)}\left(f_{kl},t\right) r_{kl}^{\beta\gamma(\mu-1)}(z_l,t)-\sum_{\eta=1}^p\left( r_{ik}^{\alpha\eta(\mu-1)}\left(f_{kl},t\right)-\sum_{\beta=1}^p \lambda_{ijk|\mu}^{\alpha\beta}  r_{0ik}^{\beta\eta}\right) r_{kl}^{\eta\gamma(\mu-1)}(z_l,t) -\sum_{\beta,\eta=1}^p r_{0ij}^{\alpha\beta}  \lambda_{jkl|\mu}^{\beta\eta}r_{0jl}^{\eta\gamma}\\
&\equiv_\mu \sum_{\beta=1}^p \lambda_{ikl|\mu}^{\alpha\beta}r_{0il}^{\beta\gamma}+\sum_{\beta=1}^p \lambda_{ijk|\mu}^{\alpha\beta} r_{0il}^{\beta\gamma}-\sum_{\beta,\eta=1}^p r_{0ij}^{\alpha\beta} \lambda_{jkl|\mu}^{\beta\eta} r_{0jl}^{\eta\gamma}
\end{align*}}}

Then we have
\begin{align*}
\sum_{\beta,\eta=1}^p r_{0ij}^{\alpha\beta}\lambda_{jkl|\mu}^{\beta\eta}r_{0jl}^{\eta\gamma} -\sum_{\beta=1}^p  \lambda_{ikl|\mu}^{\alpha\beta}r_{0il}^{\beta\gamma} +\sum_{\beta=1}^p  \lambda_{ijl|\mu}^{\alpha\beta}r_{0il}^{\beta\gamma} -\sum_{\beta=1}^p  \lambda_{ijk|\mu}^{\alpha\beta}r_{0il}^{\beta\gamma}=0\\
\end{align*}
By multiplying $\sum_{\gamma=1}^p r_{0li}^{\gamma \xi} $, we have
\begin{align*}
\sum_{\beta,\eta=1}^p r_{0ij}^{\alpha\beta}\lambda_{jkl|\mu}^{\beta\eta}r_{0ji}^{\eta\xi}-\lambda_{ikl|\mu}^{\alpha\xi}+\lambda_{ijl|\mu}^{\alpha\xi}-\lambda_{ijk|\mu}^{\alpha\xi}=0
\end{align*}

We prove $(\ref{t22})$. In fact,
\begin{align*}
\bar{\partial} \Phi_{i|\mu}^\alpha&\equiv_\mu \bar{\partial}\left(\left[\varphi^{\mu-1} , T_i^{\alpha(\mu-1)}\right]  \right)\equiv_\mu \left[\bar{\partial} \varphi^{\mu-1}, T_i^{\alpha(\mu-1)}\right]-\left[\varphi^{\mu-1}, \bar{\partial} T_i^{\alpha(\mu-1)} \right]\\
                     &\equiv_\mu \left[\xi_\mu+\frac{1}{2}\left[ \varphi^{\mu-1}, \varphi^{\mu-1}\right], T_i^{\alpha(\mu-1)} \right] - \left[ \varphi^{\mu-1}, \left[ \varphi^{\mu-1}, T_i^{\alpha(\mu-1)}\right] -\Phi_{i|\mu}^\alpha \right]
                     \equiv_\mu \left[\xi_\mu, T_{0i}^\alpha \right]
\end{align*}

We prove $(\ref{t23})$. In fact,
\begin{align*}
\sum_{\eta=1}^p \bar{\partial}\Lambda_{ij|\mu}^{\alpha\eta}\cdot r_{0ij}^{\eta\beta}&=\bar{\partial}\left(\sum_{\eta=1}^p \Lambda_{ij|\mu}^{\alpha\eta} r_{0ij}^{\eta\beta} \right)=\bar{\partial} \left(\left[\varphi^{\mu-1}, r_{ij}^{\alpha\beta(\mu-1)} \right]-\bar{\partial} r_{ij}^{\alpha\beta(\mu-1)} \right)=\bar{\partial}\left(\left[\varphi^{\mu-1}, r_{ij}^{\alpha\beta(\mu-1)} \right] \right)\\
              &=\left[\xi_\mu+\frac{1}{2}\left[\varphi^{\mu-1},\varphi^{\mu-1} \right] , r_{ij}^{\alpha\beta(\mu-1)}\right]    -\left[\varphi^{\mu-1},\left[\varphi^{\mu-1}, r_{ij}^{\alpha\beta(\mu-1)} \right] -\sum_{\eta=1}^p \Lambda_{ij|\mu}^{\alpha\eta} r_{0ij}^{\eta\beta} \right]  
              =\left[\xi_\mu, r_{0ij}^{\alpha\beta} \right]                         
\end{align*}
By multiplying $\sum_{\beta=1}^p r_{0ji}^{\beta\xi}$, we have
\begin{align*}
\bar{\partial} \Lambda_{ij|\mu}^{\alpha\xi} =\sum_{\beta=1}^p \left[\xi_\mu, r_{0ij}^{\alpha\beta} \right] r_{0ji}^{\beta\xi}
\end{align*}

We prove $(\ref{t24})$. In fact,
{\small{\begin{align*}
&-\sum_{\eta=1}^p \Lambda_{ik|\mu}^{\alpha\eta}r_{0ik}^{\eta\beta}  \equiv_\mu \bar{\partial} r_{ik}^{\alpha\beta(\mu-1)}(z_k,t)- \left[\varphi^{\mu-1}, r_{ik}^{\alpha\beta(\mu-1)}(z_k,t) \right]\\
                                                                                                      &\equiv_\mu \bar{\partial}\left( \sum_{\eta=1}^p r_{ij}^{\alpha\eta(\mu-1)}\left(f_{jk},t\right)r_{jk}^{\eta\beta(\mu-1)}(z_k,t)+\sum_{\eta=1}^p\lambda_{ijk|\mu}^{\alpha\eta} r_{0ik}^{\eta\beta}\right)-\left[\varphi^{\mu-1},   \sum_{\eta=1}^p r_{ij}^{\alpha\eta(\mu-1)}\left( f_{jk},t\right)r_{jk}^{\eta\beta(\mu-1)}(z_k,t)+\sum_{\eta=1}^p \lambda_{ijk|\mu}^{\alpha\eta} r_{0ik}^{\eta\beta}  \right]\\
                                                                                                      &\equiv_\mu \sum_{\eta=1}^p r_{ij}^{\alpha\eta(\mu-1)}\left( \bar{\partial} r_{jk}^{\eta\beta(\mu-1)}-\left[\varphi^{\mu-1} , r_{jk}^{\eta\beta(\mu-1)}\right]\right) +\sum_{\eta=1}^p\left( \bar{\partial} r_{ij}^{\alpha\eta} -\left[\varphi^{\mu-1} , r_{ij}^{\alpha\eta(\mu-1)}\right]   \right) r_{jk}^{\eta\beta(\mu-1)}+\sum_{\eta=1}^p \bar{\partial} \lambda_{ijk|\mu}^{\alpha\eta}r_{0ik}^{\eta\beta}   \\
                                                                                                      &\equiv_\mu -\sum_{\eta,\gamma=1}^p r_{0ij}^{\alpha\eta}  \Lambda_{jk|\mu}^{\eta\gamma}r_{0jk}^{\gamma\beta} -\sum_{\eta,\gamma=1}^p \Lambda_{ij|\mu}^{\alpha\gamma}r_{0ij}^{\gamma\eta}  r_{0jk}^{\eta\beta} +\sum_{\eta=1}^p \bar{\partial} \lambda_{ijk|\mu}^{\alpha\eta} r_{0ik}^{\eta\beta}
\end{align*}}}
Then we have
\begin{align*}
\sum_{\eta,\gamma=1}^p r_{0ij}^{\alpha\eta} \Lambda_{jk|\mu}^{\eta\gamma}r_{0jk}^{\gamma\beta}  -\sum_{\eta=1}^p  \Lambda_{ik|\mu}^{\alpha\eta}r_{0ik}^{\eta\beta} + \sum_{\gamma=1}^p  \Lambda_{ij|\mu}^{\alpha\gamma} r_{0ik}^{\gamma\beta}= \sum_{\eta=1}^p\bar{\partial} \lambda_{ijk|\mu}^{\alpha\eta}  r_{0ik}^{\eta\beta} 
\end{align*}
By multiplying $\sum_{\beta=1}^p r_{0ki}^{\beta\xi}$, we have
\begin{align*}
\sum_{\eta,\gamma=1}^p r_{0ij}^{\alpha\eta}\Lambda_{jk|\mu}^{\eta\gamma}r_{0ji}^{\gamma\xi}-\Lambda_{ik|\mu}^{\alpha\xi}+ \Lambda_{ij|\mu}^{\alpha\xi}=\bar{\partial} \lambda_{ijk|\mu}^{\alpha\xi}
\end{align*}

We prove $(\ref{t25})$. In fact,
\begin{align*}
\Gamma_{ik|\mu}^\alpha&\equiv_\mu T_i^{\alpha(\mu-1)}-\sum_{\beta=1}^p r_{ik}^{\alpha\beta(\mu-1)} T_k^{\beta(\mu-1)}\equiv_\mu T_i^{\alpha(\mu-1)} -\sum_{\beta=1}^p\left(\sum_{\gamma=1}^p \lambda_{ijk|\mu}^{\alpha\gamma}r_{0ik}^{\gamma\beta} +\sum_{\gamma=1}^p r_{ij}^{\alpha\gamma(\mu-1)} r_{jk}^{\gamma\beta(\mu-1)} \right) T_k^{\beta(\mu-1)}\\
     &\equiv_\mu \sum_{\beta=1}^p r_{ij}^{\alpha\beta(\mu-1)}T_j^{\beta(\mu-1)} + \Gamma_{ij|\mu}^\alpha -\sum_{\gamma=1}^p \lambda_{ijk|\mu}^{\alpha\gamma} T_{0i}^\gamma -\sum_{\gamma=1}^p r_{ij}^{\alpha\gamma(\mu-1)}\left( T_j^{\gamma(\mu-1)}-\Gamma_{jk|\mu}^\gamma \right)\\
     &\equiv_\mu \Gamma_{ij|\mu}^\alpha +\sum_{\gamma=1}^p r_{0ij}^{\alpha\gamma} \Gamma_{jk|\mu}^\gamma -\sum_{\gamma=1}^p \lambda_{ijk|\mu}^{\alpha\gamma} T_{0i}^\gamma
\end{align*}

We prove $(\ref{t26})$. In fact,
\begin{align*}
-\Phi_{i|\mu}^\alpha &\equiv_\mu \bar{\partial} T_i^{\alpha(\mu-1)}- \left[\varphi^{\mu-1}, T_i^{\alpha(\mu-1)} \right]  \\
                    &\equiv_\mu \bar{\partial}\left(\sum_{\beta=1}^p r_{ij}^{\alpha\beta(\mu-1)}T_j^{\beta(\mu-1)} +\Gamma_{ij|\mu}^\alpha \right) -\left[\varphi^{\mu-1}, \sum_{\beta=1}^p r_{ij}^{\alpha\beta(\mu-1)} T_j^{\beta(\mu-1)}+\Gamma_{ij|\mu}^\alpha \right] \\
                    &\equiv_\mu \sum_{\beta=1}^p r_{ij}^{\alpha\beta(\mu-1)}\left( \bar{\partial} T_j^{\beta(\mu-1)}- \left[ \varphi^{\mu-1}, T_j^{\beta(\mu-1)}\right]\right) +\sum_{\beta=1}^p\left(\bar{\partial} r_{ij}^{\alpha\beta(\mu-1)}-\left[ \varphi^{\mu-1}, r_{ij}^{\alpha\beta(\mu-1)} \right] \right) T_j^{\beta(\mu-1)} +\bar{\partial} \Gamma_{ij|\mu}^\alpha \\
                    &\equiv_\mu -\sum_{\beta=1}^p r_{0ij}^{\alpha\beta} \Phi_{j|\mu}^\beta -\sum_{\eta,\beta=1}^p  \Lambda_{ij|\mu}^{\alpha\eta} r_{0ij}^{\eta\beta} T_{0j}^\beta + \bar{\partial} \Gamma_{ij|\mu}^\alpha 
\end{align*}

We prove $(\ref{t27})$. In fact,
{\Small{\begin{align*}
&\Pi_{i|\mu}^{\alpha\beta}\equiv_\mu\left[ T_i^{\alpha(\mu-1)},T_i^{\beta(\mu-1)}\right] -\sum_{\gamma=1}^p g_{i\alpha\beta}^{\gamma(\mu-1)} T_i^{\gamma(\mu-1)}\\
                                       &\equiv_\mu \left[\Gamma_{ij|\mu}^\alpha+\sum_{\eta=1}^p r_{ij}^{\alpha\eta(\mu-1)} T_j^{\eta(\mu-1)}, \Gamma_{ij|\mu}^\beta+\sum_{\xi=1}^p  r_{ij}^{\beta\xi(\mu-1)} T_j^{\xi(\mu-1)} \right] -\sum_{\gamma=1}^p g_{i\alpha\beta}^{\gamma(\mu-1)}\left(\Gamma_{ij|\mu}^\gamma+\sum_{\xi=1}^p r_{ij}^{\gamma\xi(\mu-1)} T_j^{\xi(\mu-1)} \right)\\
                                       &\equiv_\mu \left[\Gamma_{ij|\mu}^\alpha, \sum_{\xi=1}^p r_{0ij}^{\beta\xi} T_{0j}^\xi\right]   + \left[ \sum_{\eta=1}^p r_{0ij}^{\alpha\eta} T_{0j}^\eta, \Gamma_{ij|\mu}^\beta\right]+\sum_{\eta,\xi=1}^p r_{ij}^{\alpha\eta(\mu-1)}r_{ij}^{\beta\xi(\mu-1)}\left[T_j^{\eta(\mu-1)}, T_j^{\xi(\mu-1)}\right] \\
                                       &+\sum_{\eta,\xi=1}^p r_{ij}^{\alpha\eta(\mu-1)} T_j^{\eta(\mu-1)}\left(r_{ij}^{\beta\xi(\mu-1)}\right) T_j^{\xi(\mu-1)}-\sum_{\eta,\xi=1}^p r_{ij}^{\beta\xi(\mu-1)} T_j^{\xi(\mu-1)}\left(r_{ij}^{\alpha\eta(\mu-1)} \right) T_j^{\eta(\mu-1)} -\sum_{\gamma=1}^p g_{i\alpha\beta}^{\gamma(\mu-1)} \Gamma_{ij|\mu}^\gamma 
                                       -\sum_{\gamma,\xi=1}^p g_{i\alpha\beta}^{\gamma(\mu-1)} r_{ij}^{\gamma\xi(\mu-1)} T_j^{\xi(\mu-1)}\\
                                       &\equiv_\mu \left[ \Gamma_{ij|\mu}^\alpha, T_{0i}^\beta \right] + \left[T_{0i}^\alpha , \Gamma_{ij|\mu}^\beta \right] +\sum_{\eta,\xi=1}^p r_{ij}^{\alpha\eta(\mu-1)} r_{ij}^{\beta\xi(\mu-1)}\left( \Pi_{j|\mu}^{\eta\xi}+\sum_{\gamma=1}^p g_{j\eta\xi}^{\gamma(\mu-1)} T_j^{\gamma(\mu-1)} \right)\\
     &+\sum_{\eta,\xi=1}^p r_{ij}^{\alpha\eta(\mu-1)} T_j^{\eta(\mu-1)}\left(r_{ij}^{\beta\xi(\mu-1)}\right) T_j^{\xi(\mu-1)}-\sum_{\eta,\xi=1}^p r_{ij}^{\beta\xi(\mu-1)} T_j^{\xi(\mu-1)}\left(r_{ij}^{\alpha\eta(\mu-1)} \right) T_j^{\eta(\mu-1)} -\sum_{\gamma=1}^p g_{0i\alpha\beta}^{\gamma} \Gamma_{ij|\mu}^\gamma 
               -\sum_{\gamma,\xi=1}^p g_{i\alpha\beta}^{\gamma(\mu-1)} r_{ij}^{\gamma\xi(\mu-1)} T_j^{\xi(\mu-1)}\\
     &\equiv_\mu \left[\Gamma_{ij|\mu}^\alpha, T_{0i}^\beta \right] +\left[ T_{0i}^\alpha, \Gamma_{ij|\mu}^\beta \right]-\sum_{\gamma=1}^p g_{0i\alpha\beta}^\gamma \Gamma_{ij|\mu}^\gamma +\sum_{\eta,\xi=1}^p r_{0ij}^{\alpha\eta} r_{0ij}^{\beta\xi} \Pi_{j|\mu}^{\eta\xi} \\
   &+\sum_{\xi=1}^p\left(\sum_{\eta,\gamma=1}^p r_{ij}^{\alpha\eta(\mu-1)}r_{ij}^{\beta\gamma(\mu-1)} g_{j\eta\gamma}^{\xi(\mu-1)}+\sum_{\eta=1}^p r_{ij}^{\alpha\eta(\mu-1)} T_j^{\eta(\mu-1)}\left(r_{ij}^{\beta\xi(\mu-1)}\right)-\sum_{\eta=1}^p r_{ij}^{\beta\eta(\mu-1)} T_j^{\eta(\mu-1)}\left(r_{ij}^{\alpha\xi(\mu-1)} \right) -\sum_{\gamma=1}^p g_{i\alpha\beta}^{\gamma(\mu-1)} r_{ij}^{\gamma\xi(\mu-1)}  \right) T_j^{\xi(\mu-1)}\\
      &\equiv_\mu \left[\Gamma_{ij|\mu}^\alpha, T_{0i}^\beta \right] +\left[ T_{0i}^\alpha, \Gamma_{ij|\mu}^\beta \right]-\Gamma_{ij|\mu}\left(\left[T_{0i}^\alpha, T_{0i}^\beta\right]\right) +\sum_{\eta,\xi=1}^p r_{0ij}^{\alpha\eta} r_{0ij}^{\beta\xi} \Pi_{j|\mu}^{\eta\xi} \\
    &+\sum_{\xi=1}^p\left(\sum_{\eta,\gamma=1}^p r_{ij}^{\alpha\eta(\mu-1)}r_{ij}^{\beta\gamma(\mu-1)} g_{j\eta\gamma}^{\xi(\mu-1)}+T_i^{\alpha(\mu-1)}\left(r_{ij}^{\beta\xi(\mu-1)}\right)-T_j^{\beta(\mu-1)}\left(r_{ij}^{\alpha\xi(\mu-1)} \right) -\sum_{\gamma=1}^p g_{i\alpha\beta}^{\gamma(\mu-1)} r_{ij}^{\gamma\xi(\mu-1)} \right) T_j^{\xi(\mu-1)} \\
    &-\sum_{\xi=1}^p\Gamma_{ij|\mu}^\alpha\left(r_{0ij}^{\beta \xi} \right)T_{0j}^\xi+\sum_{\xi=1}^p \Gamma_{ij|\mu}^\beta\left(r_{0ij}^{\alpha \xi} \right) T_{0j}^\xi\\
       &\equiv_\mu \left[\Gamma_{ij|\mu}^\alpha, T_{0i}^\beta \right] +\left[ T_{0i}^\alpha, \Gamma_{ij|\mu}^\beta \right]-\Gamma_{ij|\mu}\left(\left[T_{0i}^\alpha, T_{0i}^\beta\right]\right) +\sum_{\eta,\xi=1}^p r_{0ij}^{\alpha\eta} r_{0ij}^{\beta\xi} \Pi_{j|\mu}^{\eta\xi}  +\sum_{\xi=1}^p \Psi_{ij|\mu}^{\alpha\beta\xi} T_{0j}^\xi-\sum_{\xi=1}^p\Gamma_{ij|\mu}^\alpha\left(r_{0ij}^{\beta \xi} \right)T_{0j}^\xi+\sum_{\xi=1}^p \Gamma_{ij|\mu}^\beta\left(r_{0ij}^{\alpha \xi} \right) T_{0j}^\xi  
      \end{align*}}}
      
We prove $(\ref{t28})$. In fact,
\begin{align*}
&\bar{\partial} \Pi_{i|\mu}^{\alpha\beta}\equiv_\mu \bar{\partial}\left( \left[T_i^{\alpha(\mu-1)}, T_i^{\beta(\mu-1)}  \right]-\sum_{\gamma=1}^p g_{i\alpha\beta}^{\gamma(\mu-1)} T_i^{\gamma(\mu-1)}\right)\\
&\equiv_\mu \left[\bar{\partial} T_i^{\alpha(\mu-1)}, T_i^{\beta(\mu-1)}  \right] +\left[ T_i^{\alpha(\mu-1)}, \bar{\partial} T_i^{\beta(\mu-1)}\right] -\sum_{\gamma=1}^p \bar{\partial}\left( g_{i\alpha\beta}^{\gamma(\mu-1)} \right) T_i^{\gamma(\mu-1)} -\sum_{\gamma=1}^p g_{i\alpha\beta}^{\gamma(\mu-1)} \bar{\partial }T_i^{\gamma(\mu-1)}\\
&\equiv_\mu \left[ \left[ \varphi^{\mu-1}, T_i^{\alpha(\mu-1)}   \right]-\Phi_{i|\mu}^{\alpha},  T_i^{\beta(\mu-1)}\right]+\left[T_i^{\alpha(\mu-1)}, \left[\varphi^{\mu-1}, T_i^{\beta(\mu-1)}\right]- \Phi_{i|\mu}^\beta \right]\\
&-\sum_{\gamma=1}^p\left(\left[\varphi^{\mu-1}, g_{i\alpha\beta}^{\gamma(\mu-1)} \right] -\eta_{i\alpha\beta|\mu}^\gamma \right) T_i^{\gamma(\mu-1)}-\sum_{\gamma=1}^p g_{i\alpha\beta}^{\gamma(\mu-1)}\left(\left[\varphi^{\mu-1}, T_i^{\gamma(\mu-1)} \right]- \Phi_{i|\mu}^\gamma \right)\\
&\equiv_\mu -\left[ \Phi_{i|\mu}^\alpha, T_{0i}^\beta \right] -\left[ T_{0i}^\alpha, \Phi_{i|\mu}^\beta \right] +\sum_{\gamma=1}^p g_{0i\alpha\beta}^\gamma \Phi_{i|\mu}^\gamma  +\sum_{\gamma=1}^p \eta_{i\alpha\beta|\mu}^\gamma T_{0i}^\gamma\\
&+\left[\left[ \varphi^{\mu-1}, T_i^{\alpha(\mu-1)}\right], T_i^{\beta(\mu-1)} \right] +\left[ T_i^{\alpha(\mu-1)}, \left[\varphi^{\mu-1}, T_i^{\beta(\mu-1)} \right]\right]
-\sum_{\gamma=1}^p \left[ \varphi^{\mu-1}, g_{i\alpha\beta}^{\gamma(\mu-1)} \right] T_i^{\gamma(\mu-1)}-\sum_{\gamma=1}^p g_{i\alpha\beta}^{\gamma(\mu-1)}\left[\varphi^{\mu-1}, T_i^{\gamma(\mu-1)} \right]\\
&\equiv_\mu -\left[ \Phi_{i|\mu}^\alpha, T_{0i}^\beta \right] -\left[ T_{0i}^\alpha, \Phi_{i|\mu}^\beta \right] +\sum_{\gamma=1}^p g_{0i\alpha\beta}^\gamma \Phi_{i|\mu}^\gamma  +\sum_{\gamma=1}^p \eta_{i\alpha\beta|\mu}^\gamma T_{0i}^\gamma +\left[ \varphi^{\mu-1}, \Pi_{i|\mu}^{\alpha\beta}\right]\\
&\equiv_\mu -\left[ \Phi_{i|\mu}^\alpha, T_{0i}^\beta \right] -\left[ T_{0i}^\alpha, \Phi_{i|\mu}^\beta \right] +\sum_{\gamma=1}^p g_{0i\alpha\beta}^\gamma \Phi_{i|\mu}^\gamma  +\sum_{\gamma=1}^p \eta_{i\alpha\beta|\mu}^\gamma T_{0i}^\gamma
\end{align*}

We prove $(\ref{te7})$. In fact,      
{\small{\begin{align*}
&\left[\Pi_{i|\mu}^{\alpha \beta}, T_{0i}^\delta \right]\equiv_\mu \left[ \left[ T_i^{\alpha(\mu-1)}, T_i^{\beta(\mu-1)}\right],  T_i^{\delta(\mu-1)}\right] -\sum_{\gamma=1}^p\left[ g_{i\alpha\beta}^{\gamma(\mu-1)} T_i^{\gamma(\mu-1)}, T_i^{\delta(\mu-1)}\right]\\
&\equiv_\mu \left[\left[ T_i^{\alpha(\mu-1)}, T_i^{\delta(\mu-1)} \right],  T_i^{\beta(\mu-1)}\right]+ \left[ T_i^{\alpha(\mu-1)}, \left[T_i^{\beta(\mu-1)}, T_i^{\delta(\mu-1)} \right]\right]-\sum_{\gamma=1}^p g_{i\alpha\beta}^{\gamma(\mu-1)}\left[ T_i^{\gamma(\mu-1)}, T_i^{\delta(\mu-1)}\right]+\sum_{\gamma=1}^p T_i^{\delta(\mu-1)}\left(g_{i\alpha\beta}^{\gamma(\mu-1)}\right) T_i^{\gamma(\mu-1)}\\
&\equiv_\mu \left[\Pi_{i|\mu}^{\alpha\delta}, T_{0i}^\beta \right] +\left[\sum_{\gamma=1}^p g_{i\alpha \delta}^{\gamma(\mu-1)} T_i^{\gamma(\mu-1)} , T_i^{\beta(\mu-1)}\right] +\left[ T_{0i}^\alpha, \Pi_{i|\mu}^{\beta \delta}\right] +\left[ T_i^{\alpha(\mu-1)},\sum_{\gamma=1}^p g_{i\beta\delta}^{\gamma(\mu-1)} T_i^{\gamma(\mu-1)} \right]\\
& -\sum_{\gamma=1}^p g_{0i\alpha\beta}^\gamma\Pi_{i|\mu}^{\gamma \delta} -\sum_{\gamma,\xi=1}^p g_{i\alpha\beta}^{\gamma(\mu-1)} g_{i\gamma \delta}^{\xi(\mu-1)} T_i^{\xi(\mu-1)} +\sum_{\gamma=1}^p T_i^{\delta(\mu-1)}\left(g_{i\alpha\beta}^{\gamma(\mu-1)}\right) T_i^{\gamma(\mu-1)}\\
&\equiv_\mu \left[ \Pi_{i|\mu}^{\alpha \delta}, T_{0i}^\beta \right] +\sum_{\gamma=1}^p g_{0i\alpha \delta}^\gamma \Pi_{i|\mu}^{\gamma \beta} +\sum_{\gamma,\xi=1}^p g_{i\alpha \delta}^{\gamma(\mu-1)} g_{i\gamma \beta}^{\xi(\mu-1)} T_i^{\xi(\mu-1)}-\sum_{\gamma=1}^p T_i^{\beta(\mu-1)}\left( g_{i\alpha\delta}^{\gamma(\mu-1)}\right) T_i^{\gamma(\mu-1)}\\
& +\left[ T_{0i}^\alpha, \Pi_{i|\mu}^{\beta \delta} \right] +\sum_{\gamma=1}^p g_{0i\beta\delta}^\gamma \Pi_{i|\mu}^{\alpha \gamma}+\sum_{\gamma,\xi=1}^p g_{i\beta\delta}^{\gamma(\mu-1)}g_{i\alpha \gamma}^{\xi(\mu-1)} T_i^{\xi(\mu-1)} +\sum_{\gamma=1}^p T_i^{\alpha(\mu-1)}\left( g_{i\beta\delta}^{\gamma(\mu-1)} \right) T_i^{\gamma(\mu-1)}\\
& -\sum_{\gamma=1}^p g_{0i\alpha\beta}^\gamma\Pi_{i|\mu}^{\gamma \delta} -\sum_{\gamma,\xi=1}^p g_{i\alpha\beta}^{\gamma(\mu-1)} g_{i\gamma \delta}^{\xi(\mu-1)} T_i^{\xi(\mu-1)} +\sum_{\gamma=1}^p T_i^{\delta(\mu-1)}\left(g_{i\alpha\beta}^{\gamma(\mu-1)}\right) T_i^{\gamma(\mu-1)}\\
&\equiv_\mu \left[ \Pi_{i|\mu}^{\alpha \delta}, T_{0i}^\beta \right] +\sum_{\gamma=1}^p g_{0i\alpha \delta}^\gamma \Pi_{i|\mu}^{\gamma \beta} +\left[ T_{0i}^\alpha, \Pi_{i|\mu}^{\beta \delta} \right] +\sum_{\gamma=1}^p g_{0i\beta\delta}^\gamma \Pi_{i|\mu}^{\alpha \gamma} -\sum_{\gamma=1}^p g_{0i\alpha\beta}^\gamma\Pi_{i|\mu}^{\gamma \delta} + \sum_{\gamma=1}^p G_{i|\mu}^{\alpha\beta \delta\xi} T_{0i}^\xi\\
\end{align*}}}
This completes the proof of Lemma \ref{te10}.
\end{proof}

Our purpose is to construct $\varphi^\mu=\varphi^{\mu-1}+\varphi_\mu, r_{ij}^{\alpha\beta\mu}=r_{ij}^{\alpha\beta(\mu-1)}+r_{ij|\mu}^{\alpha\beta}$, and $T_i^{\alpha\mu}=T_i^{\alpha(\mu-1)}+T_{i|\mu}^\alpha$, and $g_{i\alpha\beta}^{\gamma\mu}=g_{i\alpha\beta}^{\gamma(\mu-1)}+g_{i\alpha\beta | \mu}^\gamma$ satisfying $(\ref{t13})_\mu-(\ref{te6})_\mu$.

\begin{lemma}\label{tt5}
$(\ref{t13})_\mu-(\ref{te6})_\mu$ are equivalent to the following equalities:
\begin{align}
\bar{\partial}\varphi_\mu&=-\xi_\mu \label{te11}\\
\bar{\partial} T_{i|\mu}^\alpha-\left[\varphi_\mu, T_{0i}^\alpha \right]&=\Phi_{i|\mu}^\alpha \label{te12}\\
\bar{\partial} r_{ij|\mu}^{\alpha\beta}-\left[\varphi_\mu, r_{0ij}^{\alpha\beta} \right]&= \sum_{\eta=1}^p \Lambda_{ij|\mu}^{\alpha\eta} r_{0ij}^{\eta\beta} \label{te13}\\
\bar{\partial} g_{i\alpha\beta |\mu}^\gamma -\left[\varphi_\mu, g_{0i\alpha\beta}^{\gamma} \right] &= \eta_{i\alpha\beta|\mu}^\gamma \label{te14}\\
\sum_{\beta=1}^p r_{0ij}^{\alpha\beta} T_{j|\mu}^{\beta}- T_{i|\mu}^\alpha+&\sum_{\beta=1}^p r_{ij|\mu}^{\alpha\beta} T_{0j}^\beta=\Gamma_{ij|\mu}^\alpha \label{te15}\\
 -\left[T_{i|\mu}^\alpha, T_{0i}^\beta\right] - \left[T_{0i}^\alpha, T_{i|\mu}^\beta \right] +\sum_{\gamma=1}^p & g_{0i\alpha\beta}^\gamma T_{i|\mu}^\gamma +\sum_{\gamma=1}^p g_{i\alpha\beta|\mu}^\gamma T_{0i}^{\gamma}=\Pi_{i|\mu}^{\alpha\beta} \label{te16}\\
\sum_{\beta,\xi=1}^p r_{0ij}^{\alpha\beta} r_{jk|\mu}^{\beta\gamma} r_{0ki}^{\gamma\xi}-\sum_{\xi=1}^p &r_{ik|\mu}^{\alpha\gamma} r_{0ki}^{\gamma\xi}+\sum_{\beta=1}^p r_{ij|\mu}^{\alpha\beta} r_{0ji}^{\beta\xi}=\lambda_{ijk|\mu}^{\alpha\xi} \label{te17}
\end{align}
\begin{align}
 -\Psi_{ij|\mu}^{\alpha\beta\xi}&= \sum_{\gamma,\eta=1}^p r_{ij|\mu}^{\alpha\eta} r_{0ij}^{\beta\gamma} g_{0j\eta\gamma}^{\xi} +\sum_{\gamma,\eta=1}^p r_{0ij}^{\alpha\eta} r_{ij|\mu}^{\beta\gamma} g_{0j\eta\gamma}^\xi + \sum_{\gamma,\eta=1}^p r_{0ij}^{\alpha\eta} r_{0ij}^{\beta\gamma} g_{i\eta\gamma|\mu}^\xi  + T_{0i}^\alpha\left(r_{ij|\mu}^{\beta \xi} \right) \label{t31} \\ 
&+ T_{i|\mu}^\alpha \left( r_{0ij}^{\beta\xi}\right) - T_{0i}^\beta\left( r_{ij|\mu}^{\alpha \xi}\right) - T_{i|\mu}^\beta\left( r_{0ij}^{\alpha \xi} \right) -\sum_{\gamma=1}^p g_{0i\alpha\beta}^\gamma r_{ij|\mu}^{\gamma\xi}-\sum_{\gamma=1}^p g_{i\alpha\beta|\mu}^\gamma r_{0ij}^{\gamma\xi}\notag 
\end{align}
{\small{\begin{align}
-G_{i|\mu}^{\alpha\beta \delta\xi}&= \sum_{\gamma=1}^p g_{0i\alpha\delta}^\gamma g_{i\gamma\beta|\mu}^\xi +\sum_{\gamma=1}^p g_{i\alpha\delta|\mu}^\gamma g_{0i\gamma\beta}^{\xi}- T_{0i}^\beta\left( g_{i\alpha\delta|\mu}^\xi\right) -T_{i|\mu}^\beta\left( g_{0i\alpha\delta}^\xi \right) +T_{0i}^\alpha\left(g_{i\beta\delta|\mu}^\xi \right) + T_{i|\mu}^\alpha\left( g_{0i\beta\delta}^\xi \right) \label{t33}\\
&-\sum_{\gamma=1}^p g_{0i\alpha\beta}^\gamma g_{i\gamma\delta|\mu}^\xi -\sum_{\gamma=1}^p g_{i\alpha\beta|\mu}^\gamma g_{0i\gamma\delta}^\xi + T_{0i}^\delta\left(g_{i\alpha\beta|\mu}^\xi \right) + T_{i|\mu}^\delta\left( g_{0i\alpha\beta}^\xi \right) + \sum_{\gamma=1}^p g_{0i\beta\delta}^\gamma g_{i\alpha\gamma |\mu}^\xi + \sum_{\gamma=1}^p g_{i\beta \delta|\mu}^\gamma g_{0i\alpha\gamma}^\xi \notag
\end{align}}}
\end{lemma}

\begin{proof}
$(\ref{te11})$ follows from \cite{Kod05} p.272. We prove $(\ref{te12})$. In fact,
\begin{align*}
\bar{\partial} T_i^{\alpha\mu}-\left[\varphi^\mu, T_i^{\alpha\mu} \right] &\equiv_\mu \bar{\partial}\left( T_i^{\alpha(\mu-1)}+T_{i|\mu}^\alpha \right) -\left[\varphi^{\mu-1}+\varphi_\mu, T_i^{\alpha(\mu-1)}+T_{i|\mu}^\alpha \right]\equiv_\mu -\Phi_{i|\mu}^\alpha + \bar{\partial} T_{i|\mu}^\alpha -\left[\varphi_\mu, T_{0i}^\alpha \right]=0
\end{align*}
We prove $(\ref{te13})$. In fact,
\begin{align*}
\bar{\partial} r_{ij}^{\alpha\beta\mu}-\left[\varphi^\mu, r_{ij}^{\alpha\beta\mu} \right]&\equiv_\mu \bar{\partial}\left(r_{ij}^{\alpha\beta(\mu-1)}+r_{ij|\mu}^{\alpha\beta} \right)-\left[\varphi^{\mu-1}+\varphi_\mu, r_{ij}^{\alpha\beta(\mu-1)}+r_{ij|\mu}^{\alpha\beta} \right]\equiv_\mu  -\sum_{\eta=1}^p \Lambda_{ij|\mu}^{\alpha\eta} r_{0ij}^{\eta\beta} + \bar{\partial} r_{ij|\mu}^{\alpha\beta} -\left[\varphi_\mu, r_{0ij}^{\alpha\beta} \right]=0
\end{align*}
We prove $(\ref{te14})$. In fact, 
\begin{align*}
\bar{\partial} g_{i\alpha\beta}^{\gamma\mu}-\left[\varphi^\mu, g_{i\alpha\beta}^{\gamma\mu} \right] \equiv_\mu \bar{\partial}\left( g_{i\alpha\beta}^{\gamma(\mu-1)}+g_{i\alpha\beta|\mu}^{\gamma}\right) -\left[ \varphi^{\mu-1}+\varphi_\mu, g_{i\alpha\beta}^{\gamma(\mu-1)}+ g_{i\alpha\beta|\mu}^\gamma \right]\equiv_\mu -\eta_{i\alpha\beta|\mu}^\gamma+ \bar{\partial} g_{i\alpha\beta |\mu}^\gamma -\left[\varphi_\mu, g_{0i\alpha\beta}^{\gamma} \right] = 0
\end{align*}
We prove $(\ref{te15})$. In fact,
\begin{align*}
T_i^{\alpha\mu}-\sum_{\beta=1}^p r_{ij}^{\alpha\beta\mu} T_j^{\beta\mu}&\equiv_\mu T_i^{\alpha(\mu-1)}+T_{i|\mu}^{\alpha}-\sum_{\beta=1}^p \left(r_{ij}^{\alpha\beta(\mu-1)}+r_{ij|\mu}^{\alpha\beta} \right) \left(  T_j^{\beta(\mu-1)}+ T_{j|\mu}^\beta \right) \\
                                  &\equiv_\mu \Gamma_{ij|\mu}^\alpha+ T_{i|\mu}^\alpha -\sum_{\beta=1}^p r_{ij|\mu}^{\alpha\beta} T_{0j}^\beta -\sum_{\beta=1}^p r_{0ij}^{\alpha\beta} T_{j|\mu}^\beta =0
\end{align*}
We prove $(\ref{te16})$. In fact,
\begin{align*}
\left[ T_i^{\alpha\mu}, T_i^{\beta\mu}\right]- \sum_{\gamma=1}^p g_{i\alpha\beta}^{\gamma\mu} T_i^{\gamma\mu} &\equiv_\mu \left[ T_i^{\alpha(\mu-1)}+T_{i|\mu}^\alpha, T_i^{\beta(\mu-1)}+T_{i|\mu}^\beta\right] -\sum_{\gamma=1}^p \left(g_{i\alpha\beta}^{\gamma(\mu-1)}+g_{i\alpha\beta|\mu}^\gamma \right)\left(T_i^{\gamma(\mu-1)}+T_{i|\mu}^\gamma \right)\\
&\equiv_\mu \Pi_{i|\mu}^{\alpha\beta} +\left[T_{i|\mu}^\alpha, T_{0i}^\beta\right] + \left[T_{0i}^\alpha, T_{i|\mu}^\beta \right] -\sum_{\gamma=1}^p g_{0i\alpha\beta}^\gamma T_{i|\mu}^\gamma -\sum_{\gamma=1}^p g_{i\alpha\beta|\mu}^\gamma T_{0i}^{\gamma}=0
\end{align*}
We prove $(\ref{te17})$. In fact,
\begin{align*}
r_{ik}^{\alpha\gamma\mu}-\sum_{\beta=1}^p r_{ij}^{\alpha\beta\mu} r_{jk}^{\beta\gamma\mu}&\equiv_\mu \left( r_{ik}^{\alpha\gamma(\mu-1)}+r_{ik|\mu}^{\alpha\gamma}\right)- \sum_{\beta=1}^p \left(r_{ij}^{\alpha\beta(\mu-1)}+r_{ij|\mu}^{\alpha\beta}\right)\left(r_{jk}^{\beta\gamma(\mu-1)}+r_{jk|\mu}^{\beta\gamma} \right)\\
                         &\equiv_\mu \sum_{\beta=1}^p \lambda_{ijk|\mu}^{\alpha\beta} r_{0ik}^{\beta\gamma} +r_{ik|\mu}^{\alpha\gamma}-\sum_{\beta=1}^p r_{0ij}^{\alpha\beta} r_{jk|\mu}^{\beta\gamma} -\sum_{\beta=1}^p r_{ij|\mu}^{\alpha\beta} r_{0jk}^{\beta\gamma}=0
\end{align*}
By multiplying $\sum_{\gamma=1}^p r_{0ki}^{\gamma\xi}$, we have
\begin{align*}
\sum_{\beta,\xi=1}^p r_{0ij}^{\alpha\beta} r_{jk|\mu}^{\beta\gamma} r_{0ki}^{\gamma\xi}-\sum_{\xi=1}^p r_{ik|\mu}^{\alpha\gamma} r_{0ki}^{\gamma\xi}+\sum_{\beta=1}^p r_{ij|\mu}^{\alpha\beta} r_{0ji}^{\beta\xi}=\lambda_{ijk|\mu}^{\alpha\xi}
\end{align*}
We prove $(\ref{t31})$. In fact,
\begin{align*}
&\sum_{\gamma,\eta=1}^p r_{ij}^{\alpha\eta\mu} r_{ij}^{\beta\gamma\mu} g_{j\eta\gamma}^{\xi\mu}+ T_i^{\alpha \mu}\left(r_{ij}^{\beta\xi\mu}\right)- T_j^{\beta \mu}\left( r_{ij}^{\alpha\xi\mu}\right) -\sum_{\gamma=1}^p g_{i\alpha\beta}^{\gamma\mu} r_{ij}^{\gamma\xi\mu} \\
&\equiv_\mu \sum_{\gamma,\eta=1}^p \left(r_{ij}^{\alpha\eta(\mu-1)} +r_{ij|\mu}^{\alpha\eta} \right)\left(r_{ij}^{\beta\gamma(\mu-1)}+r_{ij|\mu}^{\beta\gamma} \right) \left(g_{j\eta\gamma}^{\xi(\mu-1)}+g_{i\eta\gamma|\mu}^\xi \right)  +\left( T_i^{\alpha(\mu-1)}+ T_{i|\mu}^\alpha \right)\left(r_{ij}^{\beta \xi(\mu-1)}+r_{ij|\mu}^{\beta \xi} \right)\\
&-\left(T_j^{\beta(\mu-1)}+ T_{j|\mu}^\beta \right) \left( r_{ij}^{\alpha \xi(\mu-1)}+ r_{ij|\mu}^{\alpha \xi} \right) -\sum_{\gamma=1}^p \left(g_{i\alpha\beta}^{\gamma(\mu-1)}+g_{i\alpha\beta |\mu}^\gamma \right)\left( r_{ij}^{\gamma\xi(\mu-1)}+r_{ij|\mu}^{\gamma\xi}\right)\\
&\equiv_\mu \Psi_{ij|\mu}^{\alpha\beta\xi} +\sum_{\gamma,\eta=1}^p r_{ij|\mu}^{\alpha\eta} r_{0ij}^{\beta\gamma} g_{0j\eta\gamma}^{\xi} +\sum_{\gamma,\eta=1}^p r_{0ij}^{\alpha\eta} r_{ij|\mu}^{\beta\gamma} g_{0j\eta\gamma}^\xi + \sum_{\gamma,\eta=1}^p r_{0ij}^{\alpha\eta} r_{0ij}^{\beta\gamma} g_{i\eta\gamma|\mu}^\xi  + T_{0i}^\alpha\left(r_{ij|\mu}^{\beta \xi} \right) + T_{i|\mu}^\alpha \left( r_{0ij}^{\beta\xi}\right)\\
& - T_{0j}^\beta\left( r_{ij|\mu}^{\alpha \xi}\right) - T_{j|\mu}^\beta\left( r_{0ij}^{\alpha \xi} \right) -\sum_{\gamma=1}^p g_{0i\alpha\beta}^\gamma r_{ij|\mu}^{\gamma\xi}-\sum_{\gamma=1}^p g_{i\alpha\beta|\mu}^\gamma r_{0ij}^{\gamma\xi} = 0
\end{align*}
We prove $(\ref{t33})$. In fact,
{\small{\begin{align*}
&0=\sum_{\gamma=1}^p g_{i\alpha\delta}^{\gamma \mu} g_{i\gamma\beta}^{\xi\mu} - T_i^{\beta \mu}\left( g_{i\alpha\delta}^{\xi \mu} \right) +\sum_{\gamma=1}^p g_{i\beta\delta}^{\gamma \mu} g_{i\alpha\gamma}^{\xi \mu} + T_i^{\alpha\mu}\left( g_{i\beta \delta}^{\xi \mu} \right) - \sum_{\gamma=1}^p g_{i\alpha\beta}^{\gamma \mu} g_{i\gamma \delta}^{\xi \mu}+ T_i^{\delta \mu}\left( g_{i\alpha\beta}^{\xi \mu} \right)\\
&=\sum_{\gamma=1}^p\left( g_{i\alpha\delta}^{\gamma(\mu-1)}+ g_{i\alpha\delta|\mu}^{\gamma}\right) \left( g_{i\gamma\beta}^{\xi(\mu-1)}+ g_{i\gamma\beta|\mu}^{\xi} \right) -\left(T_i^{\beta(\mu-1)}+ T_{i|\mu}^{\beta} \right)\left( g_{i\alpha\delta}^{\xi(\mu-1)}+ g_{i\alpha\delta|\mu}^{\xi} \right)+ \sum_{\gamma=1}^p \left(g_{i\beta \delta}^{\gamma(\mu-1)} + g_{i\beta\delta |\mu}^\gamma \right) \left( g_{i\alpha\gamma}^{\xi(\mu-1)} + g_{i\alpha\gamma |\mu}^{\xi} \right)\\
& + \left( T_i^{\alpha(\mu-1)} + T_{i|\mu}^\alpha \right)\left( g_{i\beta \delta}^{\xi(\mu-1)}+ g_{i\beta\delta|\mu}^\xi\right)-\sum_{\gamma=1}^p \left( g_{i\alpha\beta}^{\gamma(\mu-1)} + g_{i\alpha\beta|\mu}^\gamma \right)\left( g_{i\gamma\delta}^{\xi(\mu-1)}+ g_{i\gamma\delta|\mu}^\xi \right) + \left(T_i^{\delta(\mu-1)}+ T_{i|\mu}^{\delta} \right)\left( g_{i\alpha\beta}^{\xi(\mu-1)} + g_{i\alpha\beta|\mu}^\xi\right)\\
&\equiv_\mu  G_{i|\mu}^{\alpha\beta \delta \xi} + \sum_{\gamma=1}^p g_{0i\alpha\delta}^\gamma g_{i\gamma\beta|\mu}^\xi +\sum_{\gamma=1}^p g_{i\alpha\delta|\mu}^\gamma g_{0i\gamma\beta}^{\xi}- T_{0i}^\beta\left( g_{i\alpha\delta|\mu}^\xi\right) -T_{i|\mu}^\beta\left( g_{0i\alpha\delta}^\xi \right) + \sum_{\gamma=1}^p g_{0i\beta\delta}^\gamma g_{i\alpha\gamma |\mu}^\xi + \sum_{\gamma=1}^p g_{i\beta \delta|\mu}^\gamma g_{0i\alpha\gamma}^\xi\\
& +T_{0i}^\alpha\left(g_{i\beta\delta|\mu}^\xi \right) + T_{i|\mu}^\alpha\left( g_{0i\beta\delta}^\xi \right)-\sum_{\gamma=1}^p g_{0i\alpha\beta}^\gamma g_{i\gamma\delta|\mu}^\xi -\sum_{\gamma=1}^p g_{i\alpha\beta|\mu}^\gamma g_{0i\gamma\delta}^\xi + T_{0i}^\delta\left(g_{i\alpha\beta|\mu}^\xi \right) + T_{i|\mu}^\delta\left( g_{0i\alpha\beta}^\xi \right)
\end{align*}}}
\end{proof}

\begin{lemma}\label{te47}
Under the hypothesis $\mathbb{H}^2\left(M, \Theta_{\mathcal{F}_0}^\bullet \right)=0$, we can find $\varphi_\mu, r_{ij|\mu}^{\alpha\beta}, T_{i|\mu}^\alpha$ and $g_{i\alpha\beta |\mu}^\gamma$ which satisfy $(\ref{te11})-(\ref{t33})$.
\end{lemma}

\begin{proof}
We define $\lambda_{ijk|\mu}\in \Gamma\left( U_{ijk},  \mathcal{A}^{0,0 }\left( \mathscr{H}om_{\mathcal{O}_X}\left( \Theta_{\mathcal{F}_0}, \Theta_{\mathcal{F}_0} \right) \right) \right)$ by
\begin{align*}
\lambda_{ijk|\mu}:\Gamma\left( U_{ijk}, \Theta_{\mathcal{F}_0} \right) &\to \Gamma\left( U_{ijk}, \mathcal{A}^{0,0} \left( \Theta_{\mathcal{F}_0} \right) \right) \\
       T_{0i}^\alpha &\mapsto \sum_{\xi=1}^p \lambda_{ijk|\mu}^{\alpha \xi} T_{0i}^\xi
\end{align*}
and linearly extends to $\Gamma\left( U_{ijk}, \Theta_{\mathcal{F}_0} \right)$. Then by $(\ref{t20})$, we have
\begin{align}\label{t201}
\left(\lambda_{jkl|\mu} - \lambda_{ikl|\mu}+ \lambda_{ijl|\mu}- \lambda_{ijk|\mu}\right)\left( T_{0i}^\alpha \right)= \sum_{\beta, \eta=1}^p r_{0ij}^{\alpha\beta}\lambda_{jkl|\mu}^{\beta \eta} T_{0j}^\eta - \sum_{\xi=1}^p \lambda_{ikl|\mu}^{\alpha \xi} T_{0i}^\xi + \sum_{\xi=1}^p \lambda_{ijl|\mu}^{\alpha \xi} T_{0i}^\xi - \sum_{\xi=1}^p \lambda_{ijk|\mu}^{\alpha \xi} T_{0i}^\xi=0
\end{align}
The we can find $\left\{\lambda_{ij|\mu} \right\} \in C^1\left(\mathcal{U}, \mathcal{A}^{0,0}\left( \mathscr{H}om_{\mathcal{O}_M}\left(\Theta_{\mathcal{F}_0} , \Theta_{\mathcal{F}_0} \right) \right) \right)$ with $\lambda_{ij|\mu}\in \Gamma\left( U_{ij}, \mathcal{A}^{0,0}\left( \mathscr{H}om_{\mathcal{O}_M}\left( \Theta_{\mathcal{F}_0}, \Theta_{\mathcal{F}_0} \right) \right) \right)$ defined by
\begin{align*}
\lambda_{ij|\mu}:\Gamma\left( U_{ij},  \Theta_{\mathcal{F}_0}  \right) &\to \Gamma\left( U_{ij} , \mathcal{A}^{0,0} \left( \Theta_{\mathcal{F}_0} \right)  \right) \\ 
   T_{0i}^\alpha &\mapsto \sum_{\xi=1}^p \lambda_{ij|\mu}^{\alpha \xi} T_{0i}^\xi
\end{align*}
such that $\lambda_{jk|\mu} - \lambda_{ik|\mu} + \lambda_{ij|\mu}= \lambda_{ijk|\mu}$, so that we have
\begin{align} \label{te18}
\sum_{\beta, \eta=1}^p r_{0ij}^{\alpha\beta} \lambda_{jk|\mu}^{\beta\eta} r_{0ji}^{\eta\xi}-\lambda_{ik|\mu}^{\alpha\xi}+\lambda_{ij|\mu}^{\alpha\xi}=\lambda_{ijk|\mu}^{\alpha\xi}
\end{align}
We define $\Lambda_{ij|\mu}\in \Gamma\left( U_{ij}, \mathcal{A}^{0,1}\left( \mathscr{H}om_{\mathcal{O}_M}\left( \Theta_{\mathcal{F}_0}, \Theta_{\mathcal{F}_0} \right) \right) \right)$ by
\begin{align*}
\Lambda_{ij|\mu}:\Gamma\left( U_{ij}, \Theta_{\mathcal{F}_0} \right) &\to \Gamma\left(U_{ij}, \mathcal{A}^{0,1}\left( \Theta_{\mathcal{F}_0} \right) \right) \\
 T_{0i}^\alpha &\mapsto \sum_{\xi=1}^p  \Lambda_{ij|\mu}^{\alpha \xi} T_{0i}^\xi
\end{align*}
Then from $(\ref{te18})$ and $(\ref{t24})$, we can find $\left\{\Lambda_{i|\mu} \right\}\in C^0\left( \mathcal{U}, \mathcal{A}^{0,1}\left(\mathscr{H}om_{\mathcal{O}_M}\left( \Theta_{\mathcal{F}_0}, \Theta_{\mathcal{F}_0} \right)\right)\right)$ with $\Lambda_{i|\mu}\in \Gamma\left( U_i, \mathcal{A}^{0,1}\left( \mathscr{H}om_{\mathcal{O}_M}\left( \Theta_{\mathcal{F}_0}, \Theta_{\mathcal{F}_0} \right) \right) \right)$ defined by
\begin{align*}
\Lambda_{i|\mu}:\Gamma\left( U_i, \Theta_{\mathcal{F}_0} \right) &\to \Gamma\left( U_i, \mathcal{A}^{0,1}\left( \Theta_{\mathcal{F}_0} \right) \right) \\
 T_{0i}^\alpha &\mapsto \sum_{\xi=1}^p \Lambda_{i|\mu}^{\alpha \xi} T_{0i}^\xi
\end{align*}
such that
\begin{align}\label{te21}
\sum_{\beta,\eta=1}^p r_{0ij}^{\alpha\beta}\Lambda_{j|\mu}^{\beta\eta} r_{0ji}^{\eta\xi} -\Lambda_{i|\mu}^{\alpha\xi}&=\Lambda_{ij|\mu}^{\alpha\xi}-\bar{\partial}\lambda_{ij|\mu}^{\alpha\xi}
\end{align}
We define $\Gamma_{ij|\mu}\in \Gamma\left(U_{ij}, \mathcal{A}^{0,0}\left( \mathscr{H}om_{\mathcal{O}_M}\left( \Theta_{\mathcal{F}_0}, \Theta_M \right) \right) \right)$ by
\begin{align*}
\Gamma_{ij|\mu}: \Gamma\left( U_{ij}, \Theta_{\mathcal{F}_0} \right) &\to \Gamma\left( U_{ij},  \mathcal{A}^{0,0} \left( \Theta_M \right) \right) \\
    T_{0i}^\alpha &\mapsto  \Gamma_{ij|\mu}^{\alpha}
\end{align*}
Then from $(\ref{t25})$ and $(\ref{te18})$, we have
\begin{align*}
\sum_{\beta=1}^p r_{0ij}^{\alpha\beta}\left( \Gamma_{jk|\mu}^\gamma -\sum_{\eta=1}^p \lambda_{jk|\mu}^{\beta\eta} T_{0j}^\eta \right) -\left(\Gamma_{ik|\mu}^\alpha-\sum_{\xi=1}^p \lambda_{ik|\mu}^{\alpha\xi} T_{0i}^\xi \right) +\left(\Gamma_{ij|\mu}^\alpha-\sum_{\xi=1}^p \lambda_{ij|\mu}^{\alpha\xi} T_{0i}^\xi \right) = 0.
\end{align*}
Then we can find $\left\{\Gamma_{i|\mu} \right\}\in C^0\left( \mathcal{U}, \mathcal{A}^{0,0}(\mathscr{H}om_{\mathcal{O}_M}\left(\Theta_{\mathcal{F}_0} ,\Theta_M\right)  \right)$ with $\Gamma_{i|\mu}\in \Gamma\left( U_i, \mathcal{A}^{0,0}\left( \mathscr{H}om_{\mathcal{O}_M}\left( \Theta_{\mathcal{F}_0}, \Theta_M \right) \right) \right)$ defined by
\begin{align*}
\Gamma_{i|\mu}: \Gamma\left( U_i, \Theta_{\mathcal{F}_0} \right) &\to \Gamma\left( U_i, \mathcal{A}^{0,0}\left( \Theta_M \right) \right) \\
 T_{0i}^\alpha &\mapsto \sum_{\xi=1}^p \Gamma_{i|\mu}^{\alpha \xi} T_{0i}^\xi
\end{align*}

such that
\begin{align}
\sum_{\beta=1}^p r_{0ij}^{\alpha\beta} \Gamma_{j|\mu}^\beta- \Gamma_{i|\mu}^\alpha =\Gamma_{ij|\mu}^\alpha-\sum_{\xi=1}^p \lambda_{ij|\mu}^{\alpha\xi} T_{0i}^\xi \label{t57}
\end{align}
Then from $(\ref{t27})$ and $(\ref{t57})$ and $(\ref{te19})$, we have
{\small{\begin{align*}
&\Pi_{i|\mu}^{\alpha\beta}-\sum_{\eta,\xi=1}^p r_{0ij}^{\alpha\eta} r_{0ij}^{\beta \xi} \Pi_{j|\mu}^{\eta\xi}= \left[\Gamma_{ij|\mu}^\alpha, T_{0i}^\beta \right]+\left[ T_{0i}^\alpha, \Gamma_{ij|\mu}^\beta \right]-\Gamma_{ij|\mu}\left(\left[T_{0i}^\alpha, T_{0i}^\beta \right] \right) +\sum_{\xi=1}^p \left(\Psi_{ij|\mu}^{\alpha\beta\xi} -\Gamma_{ij|\mu}^\alpha\left(r_{0ij}^{\beta \xi} \right) + \Gamma_{ij|\mu}^\beta \left( r_{0ij}^{\alpha \xi} \right)  \right)T_{0j}^\xi\\
&=\sum_{\xi=1}^p r_{0ij}^{\beta \xi}\left[\Gamma_{ij|\mu}^{\alpha}, T_{0j}^\xi \right] -\sum_{\eta=1}^p r_{0ij}^{\alpha \eta} \left[ \Gamma_{ij|\mu}^\beta,T_{0j}^\eta \right] -\sum_{\gamma=1}^p g_{0i\alpha\beta}^\gamma \Gamma_{ij|\mu}^\gamma +\sum_{\xi=1}^p \Psi_{ij|\mu}^{\alpha \beta \xi} T_{0j}^\xi\\
&=\sum_{\xi,\eta=1}^p r_{0ij}^{\beta \xi}\left[r_{0ij}^{\alpha\eta} \Gamma_{j|\mu}^\eta, T_{0j}^\xi \right] -\sum_{\xi=1}^p r_{0ij}^{\beta \xi}\left[ \Gamma_{i|\mu}^\alpha, T_{0j}^\xi\right] +\sum_{\xi,\eta=1}^p r_{0ij}^{\beta \xi}\left[ \lambda_{ij|\mu}^{\alpha\eta}T_{0i}^\eta, T_{0j}^\xi\right]\\
&-\sum_{\xi,\eta=1}^p r_{0ij}^{\alpha \eta}\left[r_{0ij}^{\beta\xi} \Gamma_{j|\mu}^\xi, T_{0j}^\eta \right]+\sum_{\eta=1}^p r_{0ij}^{\alpha \eta}\left[ \Gamma_{i|\mu}^\beta, T_{0j}^\eta\right] - \sum_{\xi,\eta=1}^p r_{0ij}^{\alpha \eta}\left[ \lambda_{ij|\mu}^{\beta\xi}T_{0i}^\xi, T_{0j}^\eta\right]\\
&-\sum_{\gamma,\xi=1}^p g_{0i\alpha\beta}^\gamma r_{0ij}^{\gamma\xi} \Gamma_{j|\mu}^\xi +\sum_{\gamma=1}^p g_{0i\alpha\beta}^\gamma \Gamma_{i|\mu}^\gamma -\sum_{\gamma,\delta=1}^p g_{0i\alpha\beta}^\gamma \lambda_{ij|\mu}^{\gamma \delta}  T_{0i}^\delta +\sum_{\xi=1}^p \Psi_{ij|\mu}^{\alpha \beta \xi} T_{0j}^\xi\\
&= \sum_{\xi,\eta=1}^p r_{0ij}^{\alpha\eta} r_{0ij}^{\beta \xi}\left[ \Gamma_{j|\mu}^\eta, T_{0j}^\xi \right]-\sum_{\eta=1}^p T_{0i}^\beta\left( r_{0ij}^{\alpha \eta}\right) \Gamma_{j|\mu}^\eta-\left[\Gamma_{i|\mu}^\alpha, T_{0i}^\beta \right]+\sum_{\xi=1}^p  \Gamma_{i|\mu}^\alpha\left( r_{0ij}^{\beta \xi} \right) T_{0j}^\xi +\sum_{\eta=1}^p \left[ \lambda_{ij|\mu}^{\alpha\eta} T_{0i}^\eta, T_{0i}^\beta\right]-\sum_{\xi,\eta=1}^p \lambda_{ij|\mu}^{\alpha \eta} T_{0i}^\eta\left( r_{0ij}^{\beta \xi} \right) T_{0j}^\xi\\
&-\sum_{\xi,\eta=1}^p r_{0ij}^{\alpha\eta} r_{0ij}^{\beta \xi}\left[\Gamma_{j|\mu}^\xi, T_{0j}^\eta \right]+\sum_{\xi=1}^p T_{0i}^\alpha \left( r_{0ij}^{\beta\xi} \right)\Gamma_{j|\mu}^\xi +\left[ \Gamma_{i|\mu}^\beta, T_{0i}^\alpha \right] -\sum_{\eta=1}^p \Gamma_{i|\mu}^\beta\left( r_{0ij}^{\alpha\eta} \right) T_{0j}^\eta -\sum_{\xi=1}^p\left[\lambda_{ij|\mu}^{\beta\xi} T_{0i}^\xi, T_{0i}^\alpha \right]+\sum_{\xi,\eta=1}^p \lambda_{ij|\mu}^{\beta\xi} T_i^\xi\left(r_{0ij}^{\alpha\eta} \right) T_{0j}^\eta\\
&-\sum_{\gamma,\xi=1}^p g_{0i\alpha\beta}^\gamma r_{0ij}^{\gamma\xi} \Gamma_{j|\mu}^\xi +\sum_{\gamma=1}^p g_{0i\alpha\beta}^\gamma \Gamma_{i|\mu}^\gamma -\sum_{\gamma,\delta=1}^p g_{0i\alpha\beta}^\gamma \lambda_{ij|\mu}^{\gamma \delta} T_{0i}^\delta +\sum_{\xi=1}^p \Psi_{ij|\mu}^{\alpha \beta \xi} T_{0j}^\xi\\
&= \sum_{\xi,\eta=1}^p r_{0ij}^{\alpha\eta} r_{0ij}^{\beta \xi}\left[ \Gamma_{j|\mu}^\eta, T_{0j}^\xi \right]-\left[\Gamma_{i|\mu}^\alpha, T_{0i}^\beta \right]+\sum_{\xi=1}^p  \Gamma_{i|\mu}^\alpha\left( r_{0ij}^{\beta \xi} \right) T_{0j}^\xi +\sum_{\eta=1}^p \left[ \lambda_{ij|\mu}^{\alpha\eta} T_{0i}^\eta, T_{0i}^\beta\right]-\sum_{\xi,\eta=1}^p \lambda_{ij|\mu}^{\alpha \eta} T_{0i}^\eta\left( r_{0ij}^{\beta \xi} \right) T_{0j}^\xi  \\
&-\sum_{\xi,\eta=1}^p r_{0ij}^{\alpha\eta} r_{0ij}^{\beta \xi}\left[\Gamma_{j|\mu}^\xi, T_{0j}^\eta \right]+\left[ \Gamma_{i|\mu}^\beta, T_{0i}^\alpha \right] -\sum_{\eta=1}^p \Gamma_{i|\mu}^\beta\left( r_{0ij}^{\alpha\eta} \right) T_{0j}^\eta -\sum_{\xi=1}^p\left[\lambda_{ij|\mu}^{\beta\xi} T_{0i}^\xi, T_{0i}^\alpha \right]+\sum_{\xi,\eta=1}^p \lambda_{ij|\mu}^{\beta\xi} T_{0i}^\xi\left(r_{0ij}^{\alpha\eta} \right) T_{0j}^\eta \\
&-\sum_{\eta,\gamma,\xi=1}^p r_{0ij}^{\alpha\eta} r_{0ij}^{\beta\gamma} g_{0j\eta\gamma}^\xi \Gamma_{j|\mu}^\xi +\sum_{\gamma=1}^p g_{0i\alpha\beta}^\gamma \Gamma_{i|\mu}^\gamma -\sum_{\gamma,\delta=1}^p g_{0i\alpha\beta}^\gamma \lambda_{ij|\mu}^{\gamma \delta}  T_{0i}^\delta +\sum_{\xi=1}^p \Psi_{ij|\mu}^{\alpha \beta \xi} T_{0j}^\xi
\end{align*}}}
Then we see that
{\small{\begin{align*}
&\Pi_{i|\mu}^{\alpha\beta} +\left[ \Gamma_{i|\mu}^\alpha, T_{0i}^\beta \right]  -  \left[ \Gamma_{i|\mu}^\beta, T_{0i}^\alpha \right] -\sum_{\gamma=1}^p g_{0i\alpha\beta}^\gamma \Gamma_{i|\mu}^\gamma       - \sum_{\eta,\xi=1}^p r_{0ij}^{\alpha\eta} r_{0ij}^{\beta\xi}  \left(\Pi_{j|\mu}^{\eta\xi}  +\left[ \Gamma_{j|\mu}^\eta, T_{0j}^\xi\right] -\left[ \Gamma_{j|\mu}^\xi, T_{0j}^\eta \right] -\sum_{\gamma=1}^p g_{0j\eta\xi}^{\gamma} \Gamma_{j|\mu}^\gamma     \right) \in \Gamma\left( U_{ij}, \mathcal{A}^{0,0}\left( \Theta_{\mathcal{F}_0}\right) \right)
\end{align*}}}

Then there exists $\left\{\Psi_{i|\mu} \right\} \in C^0\left( \mathcal{U}, \mathcal{A}^{0,0}\left( \mathscr{H}om_{\mathcal{O}_M}\left( \bigwedge^2 \Theta_{\mathcal{F}_0} , \Theta_{\mathcal{F}_0} \right) \right) \right)$ with $\Psi_{i|\mu}\in \Gamma\left( U_i, \mathcal{A}^{0,0}\left( \mathscr{H}om_{\mathcal{O}_M} \left( \Theta_{\mathcal{F}_0}, \Theta_{\mathcal{F}_0}   \right) \right) \right)$ defined by
\begin{align*}
\Psi_{i|\mu} : \Gamma\left( U_i, \Theta_{\mathcal{F}_0} \right) &\to \Gamma\left( U_i, \mathcal{A}^{0,0}\left( \Theta_{\mathcal{F}_0}  \right) \right) \\
              T_{0i}^\alpha &\mapsto \sum_{\xi=1}^p \Psi_{i|\mu}^{\alpha \xi} T_{0i}^\xi
\end{align*}
such that
\begin{align}\label{ds8}
&\Pi_{i|\mu}^{\alpha\beta} +\left[ \Gamma_{i|\mu}^\alpha, T_{0i}^\beta \right]  -  \left[ \Gamma_{i|\mu}^\beta, T_{0i}^\alpha \right] -\sum_{\gamma=1}^p g_{0i\alpha\beta}^\gamma \Gamma_{i|\mu}^\gamma   +\sum_{\gamma=1}^p \Psi_{i|\mu}^{\alpha\beta \gamma} T_{0i}^\gamma   \\
& - \sum_{\eta,\xi=1}^p r_{0ij}^{\alpha\eta} r_{0ij}^{\beta\xi}  \left(\Pi_{j|\mu}^{\eta\xi}  +\left[ \Gamma_{j|\mu}^\eta, T_{0j}^\xi\right] -\left[ \Gamma_{j|\mu}^\xi, T_{0j}^\eta \right] -\sum_{\gamma=1}^p g_{0j\eta\xi}^{\gamma} \Gamma_{j|\mu}^\gamma  +\sum_{\gamma=1}^p \Psi_{j|\mu}^{\eta \xi \gamma} T_{0j}^\gamma    \right)=0 \notag
\end{align}

This implies that
\begin{align*}
B_\mu :=\left\{B_{i|\mu} \right\} \in A^{0,0}\left(M, \mathscr{H}om_{\mathcal{O}_M}\left(\bigwedge^2 \Theta_{\mathcal{F}_0}, \Theta_M \right) \right)
\end{align*}
where 
\begin{align}\label{te24}
B_{i|\mu}\left(T_{0i}^\alpha \wedge T_{0i}^\beta \right):= \Pi_{i|\mu}^{\alpha\beta} +\left[ \Gamma_{i|\mu}^\alpha, T_{0i}^\beta \right]  -  \left[ \Gamma_{i|\mu}^\beta, T_{0i}^\alpha \right] -\sum_{\gamma=1}^p g_{0i\alpha\beta}^\gamma \Gamma_{i|\mu}^\gamma   +\sum_{\gamma=1}^p \Psi_{i|\mu}^{\alpha\beta \gamma} T_{0i}^\gamma 
\end{align}

We take $\overline{B_\mu}$ be the image in $\frac{A^{0,0}\left( M, \mathscr{H}om_{\mathcal{O}_M}\left(\bigwedge^2 \Theta_{\mathcal{F}_0} ,  \Theta_M \right) \right)}{A^{0,0}\left(M, \mathscr{H}om_{\mathcal{O}_M}\left( \bigwedge^2 \Theta_{\mathcal{F}_0} , \Theta_{\mathcal{F}_0} \right) \right)}$.

On the other hand, we set
\begin{align}\label{te23}
\tilde{\Phi}_{i|\mu}^\alpha := \Phi_{i|\mu}^\alpha-\bar{\partial} \Gamma_{i|\mu}^\alpha +\sum_{\xi=1}^p   \Lambda_{i|\mu}^{\alpha\xi} T_{0i}^\xi
\end{align}
and we define $\tilde{\Phi}_{i|\mu} \in \Gamma\left( U_i, \mathcal{A}^{0,1}\left(\mathscr{H}om_{\mathcal{O}_M}\left( \Theta_{\mathcal{F}_0}, \Theta_M  \right) \right) \right)$ by
\begin{align*}
\tilde{\Phi}_{i|\mu}:\Gamma\left( U_i, \Theta_{\mathcal{F}_0} \right) &\to \Gamma\left( U_i, \mathcal{A}^{0,1}\left( \Theta_M \right)  \right)\\
          T_{0i}^\alpha &\mapsto \tilde{\Phi}_{i|\mu}^\alpha
\end{align*}

We note that from $(\ref{t57}),(\ref{t26})$ and $(\ref{te21})$, we have
\begin{align*}
\sum_{\beta=1}^p r_{0ij}^{\alpha\beta} \tilde{\Phi}_{j|\mu}^\beta- \tilde{\Phi}_{i|\mu}^\alpha &= \sum_{\beta=1}^pr_{0ij}^{\alpha\beta}\left(\Phi_{j|\mu}^\beta-\bar{\partial} \Gamma_{j|\mu}^\beta+\sum_{\gamma=1}^p \Lambda_{j|\mu}^{\beta\gamma} T_{0j}^\gamma \right) -\left( \Phi_{i|\mu}^\alpha-\bar{\partial} \Gamma_{i|\mu}^\alpha +\sum_{\xi=1}^p   \Lambda_{i|\mu}^{\alpha\xi} T_{0i}^\xi \right)\\
&=\bar{\partial} \Gamma_{ij|\mu}^\alpha-\sum_{\eta=1}^p \Lambda_{ij|\mu}^{\alpha\eta} T_{0i}^\eta   -   \bar{\partial}\Gamma_{ij|\mu}^\alpha+\sum_{\xi=1}^p \bar{\partial}\lambda_{ij|\mu}^{\alpha\xi}T_{0i}^\xi +\sum_{\beta,\gamma=1}^p r_{0ij}^{\alpha\beta} \Lambda_{j|\mu}^{\beta\gamma} T_{0j}^\gamma -\sum_{\xi=1}^p \Lambda_{i|\mu}^{\alpha\xi} T_{0i}^\xi \\
&=-\sum_{\xi=1}^p \Lambda_{ij|\mu}^{\alpha\xi} T_{0i}^\xi  +\sum_{\xi=1}^p \bar{\partial}\lambda_{ij|\mu}^{\alpha\xi}T_{0i}^\xi +\sum_{\beta,\eta=1}^p r_{0ij}^{\alpha\beta} \Lambda_{j|\mu}^{\beta\eta} T_{0j}^\eta -\sum_{\xi=1}^p \Lambda_{i|\mu}^{\alpha\xi} T_{0i}^\xi = 0
\end{align*}

This implies that
\begin{align*}
\tilde{\Phi}_\mu:=\left\{ \tilde{\Phi}_{i|\mu} \right\} \in A^{0,1}\left(  M, \mathscr{H}om_{\mathcal{O}_M}\left( \Theta_{\mathcal{F}_0}, \Theta_M     \right)     \right)
\end{align*}
We take $\bar{\Phi}_\mu$ to be the image of $\tilde{\Phi}_\mu$ in $ \frac{A^{0,1}\left(M, \mathscr{H}om_{\mathcal{O}_M}\left(\Theta_{\mathcal{F}_0} , \Theta_M \right) \right)}{A^{0,1}\left(M, \mathscr{H}om_{\mathcal{O}_M} \left( \Theta_{\mathcal{F}_0} , \Theta_{\mathcal{F}_0} \right) \right)}$. Then we claim that 
\begin{align}\label{d2}
\left(\overline{B}_\mu, \overline{\Phi}_\mu, - \xi_\mu \right) \in \frac{A^{0,0}\left( M, \mathscr{H}om_{\mathcal{O}_M}\left( \bigwedge^2 \Theta_{\mathcal{F}_0} , \Theta_M \right) \right)}{A^{0,0}\left(M, \mathscr{H}om_{\mathcal{O}_M} \left( \bigwedge^2 \Theta_{\mathcal{F}_0}  , \Theta_{\mathcal{F}_0} \right) \right)}\bigoplus \frac{A^{0,1}\left( M, \mathscr{H}om_{\mathcal{O}_M}\left(\Theta_{\mathcal{F}_0} ,  \Theta_M \right) \right)}{A^{0,1}\left(M, \mathscr{H}om_{\mathcal{O}_M}\left(  \Theta_{\mathcal{F}_0} , \Theta_{\mathcal{F}_0} \right) \right)} \bigoplus A^{0,2}\left( M, \Theta_M \right)
\end{align}
defines a $2$-cocycle in the following Dolbeault resolution of $\Theta_{\mathcal{F}_0}^\bullet$ (see Appendix \ref{app3}).

\begin{equation}\label{tpc1}
\begin{CD}
\cdots \\
@A\hat{D}_3AA \\
\frac{A^{0,0}\left( \bigwedge^3 \Theta_{\mathcal{F}_0}^*\otimes \Theta_M\right)}{A^{0,0}\left(\bigwedge^3 \Theta_{\mathcal{F}_0}^*\otimes \Theta_{\mathcal{F}_0} \right)}@>-\bar{\partial}>> \cdots\\
@A\hat{D}_2AA @A\hat{D}_2AA \\
\frac{A^{0,0}\left(M,  \bigwedge^2 \Theta_{\mathcal{F}_0}^*\otimes \Theta_M\right)}{A^{0,0}\left(M, \bigwedge^2 \Theta_{\mathcal{F}_0}^* \otimes \Theta_{\mathcal{F}_0} \right)} @>\bar{\partial}>> \frac{A^{0,1}\left( M, \bigwedge^2 \Theta_{\mathcal{F}_0}^*\otimes \Theta_M \right)}{A^{0,1}\left(M, \bigwedge^2 \Theta_{\mathcal{F}_0}^*\otimes \Theta_{\mathcal{F}_0} \right)}@>-\bar{\partial}>> \cdots \\
@A\hat{D}_1AA @A\hat{D}_1AA @A\hat{D}_1AA\\
\frac{A^{0,0}\left(M, \Theta_{\mathcal{F}_0}^*\otimes \Theta_M\right)}{A^{0,0}\left(M, \Theta_{\mathcal{F}_0}^*\otimes \Theta_{\mathcal{F}_0} \right)} @>-\bar{\partial}>> \frac{A^{0,1}\left(M, \Theta_{\mathcal{F}_0}^*\otimes \Theta_M\right)}{A^{0,1}\left(M, \Theta_{\mathcal{F}_0}^*\otimes \Theta_{\mathcal{F}_0} \right)} @>\bar{\partial}>> \frac{A^{0,2}\left(M, \Theta_{\mathcal{F}_0}^*\otimes \Theta_M\right)}{A^{0,2}\left(M, \Theta_{\mathcal{F}_0}^*\otimes \Theta_{\mathcal{F}_0}\right)} @>-\bar{\partial}>>\cdots \\
@A\hat{D}_0AA @A\hat{D}_0AA @A\hat{D}_0AA @A\hat{D}_0AA\\
A^{0,0}(M, \Theta_M) @>\bar{\partial} >> A^{0,1}(M, \Theta_M) @>-\bar{\partial}>> A^{0,2}(M, \Theta_M)  @>\bar{\partial}>> A^{0,3}(M, \Theta_M) @>-\bar{\partial}>> \cdots\\
\end{CD}
\end{equation}

From $(\ref{t19})$, $\bar{\partial}\left(-\xi_\mu\right)=0$. We show that $\bar{\partial}\bar{\Phi}_\mu+ \hat{D}_0\left(- \xi_\mu \right)=0$. In fact, from $(\ref{te23})$ and $(\ref{t22})$ 
\begin{align}\label{d3}
\bar{\partial} \tilde{\Phi}_{i|\mu}^\alpha - \left[ \xi_\mu, T_{0i}^\alpha \right]= \bar{\partial} \Phi_{i|\mu}^\alpha + \sum_{\xi=1}^p \bar{\partial} \Lambda_{i|\mu}^{\alpha \xi} T_{0i}^\xi - \left[ \xi_\mu, T_{0i}^\alpha \right]=\sum_{\xi=1}^p \bar{\partial}\Lambda_{i|\mu}^{\alpha\xi}T_{0i}^\xi \in \Gamma\left(U_i, \mathcal{A}^{0,1}\left(\Theta_{\mathcal{F}_0} \right) \right)
\end{align}
We show that $\bar{\partial} B_\mu + \hat{D}_1\left(\bar{\Phi}_\mu \right)=0$. In fact, from $(\ref{te24}),(\ref{te23})$ and $(\ref{t28})$, we have
{\small{\begin{align*}
&\bar{\partial}B_{i|\mu}\left(T_{0i}^\alpha \wedge T_{0i}^\beta \right) + \left[ T_{0i}^\alpha,\tilde{\Phi}_{i|\mu}^{\beta}\right] - \left[ T_{0i}^\beta, \tilde{\Phi}_{i|\mu}^\alpha \right] - \sum_{\gamma=1}^p g_{0i\alpha\beta}^\gamma \tilde{\Phi}_{i|\mu}^\gamma\\
&=\bar{\partial} \Pi_{i|\mu}^{\alpha\beta} + \bar{\partial}\left[\Gamma_{i|\mu}^\alpha, T_{0i}^\beta \right]  -  \bar{\partial} \left[ \Gamma_{i|\mu}^\beta, T_{0i}^\alpha \right] -\sum_{\gamma=1}^p g_{0i\alpha\beta}^\gamma\bar{\partial} \Gamma_{i|\mu}^\gamma   +\sum_{\gamma=1}^p \bar{\partial}\Psi_{i|\mu}^{\alpha\beta \gamma} T_{0i}^\gamma \\
&+\left[T_{0i}^\alpha, \Phi_{i|\mu}^\beta \right]-\bar{\partial}\left[ T_{0i}^\alpha, \Gamma_{i|\mu}^\beta\right]+\sum_{\xi=1}^p \left[T_{0i}^\alpha, \Lambda_{i|\mu}^{\beta\xi}T_{0i}^\xi \right]  - \left[ T_{0i}^\beta, \Phi_{i|\mu}^\alpha \right] + \bar{\partial}\left[T_{0i}^\beta, \Gamma_{i|\mu}^\alpha \right] - \sum_{\xi=1}^p \left[ T_{0i}^\beta, \Lambda_{i|\mu}^{\alpha \xi} T_{0i}^\xi  \right] \\
&-\sum_{\gamma=1}^p g_{0i\alpha\beta}^\gamma\Phi_{i|\mu}^\gamma+\sum_{\gamma=1}^p g_{0i\alpha\beta}^\gamma \bar{\partial} \Gamma_{i|\mu}^\gamma -\sum_{\gamma,\xi=1}^p g_{0i\alpha\beta}^\gamma \Lambda_{i|\mu}^{\gamma\xi} T_{0i}^\xi\\
&=\sum_{\gamma=1}^p \eta_{i\alpha\beta|\mu}^\gamma T_{0i}^\gamma +\sum_{\gamma=1}^p \bar{\partial} \Psi_{i|\mu}^{\alpha\beta\gamma} T_{0i}^\gamma+\sum_{\xi=1}^p\left[T_{0i}^\alpha, \Lambda_{i|\mu}^{\beta\xi} T_{0i}^\xi \right] + \sum_{\xi=1}^p \left[ T_{0i}^\beta, \Lambda_{i|\mu}^{\alpha\xi}T_{0i}^\xi \right]  -\sum_{\gamma,\xi=1}^p g_{0i\alpha\beta}^\gamma \Lambda_{i|\mu}^{\gamma \xi} T_{0i}^\xi \in \Gamma\left( U_i, \mathcal{A}^{0,1}\left( \Theta_{\mathcal{F}_0} \right) \right)
\end{align*}}}
We show that $\hat{D}_2\left(\bar{B}_\mu \right)=0$. In fact, from $(\ref{te7})$ and $(\ref{te24})$ and $(\ref{te25})$, we have
{\Small{\begin{align*}
&\left[T_{0i}^\alpha, B_{i|\mu}\left(T_{0i}^\beta, T_{0i}^\gamma\right)\right]- \left[T_{0i}^\beta, B_{i|\mu}\left(T_{0i}^\alpha, T_{0i}^\gamma \right) \right]+\left[T_{0i}^\gamma, B_{i|\mu}\left(T_{0i}^\alpha ,  T_{0i}^\beta\right) \right] - B_{i|\mu}\left( \left[T_{0i}^\alpha, T_{0i}^\beta \right], T_{0i}^\gamma \right) + B_{i|\mu}\left(\left[T_{0i}^\alpha, T_{0i}^\gamma \right], T_{0i}^\beta \right) -B_{i|\mu}\left( \left[T_{0i}^\beta, T_{0i}^\gamma  \right], T_{0i}^\alpha \right)\\
&\equiv_{\Theta_{\mathcal{F}_0}} \left[T_{0i}^\alpha,  \Pi_{i|\mu}^{\beta \gamma}+\left[ \Gamma_{i|\mu}^\beta, T_{0i}^\gamma  \right] -\left[ \Gamma_{i|\mu}^\gamma, T_{0i}^\beta \right] - \sum_{\eta=1}^p g_{0i\beta\gamma}^\eta \Gamma_{i|\mu}^\eta \right] - \left[ T_{0i}^\beta, \Pi_{i|\mu}^{\alpha \gamma}+\left[ \Gamma_{i|\mu}^\alpha, T_{0i}^\gamma \right] -\left[ \Gamma_{i|\mu}^\gamma, T_{0i}^\alpha \right] - \sum_{\eta=1}^p g_{0i\alpha\gamma}^\eta \Gamma_{i|\mu}^\eta \right] \\
&+\left[ T_{0i}^\gamma, \Pi_{i|\mu}^{\alpha\beta} + \left[ \Gamma_{i|\mu}^{\alpha}, T_{0i}^\beta   \right] -\left[ \Gamma_{i|\mu}^\beta, T_{0i}^\alpha \right] - \sum_{\eta=1}^p g_{0i\alpha\beta}^\eta \Gamma_{i|\mu}^\eta \right] - \sum_{\eta=1}^p g_{0i\alpha\beta}^\eta \left( \Pi_{i|\mu}^{\eta \gamma} + \left[ \Gamma_{i|\mu}^\eta, T_{0i}^\gamma  \right] -\left[ \Gamma_{i|\mu}^\gamma, T_{0i}^\eta \right] - \sum_{\xi=1}^p g_{0i\eta\gamma}^\xi \Gamma_{i|\mu}^\xi \right)\\
&+\sum_{\eta=1}^p g_{0i\alpha\gamma}^\eta \left( \Pi_{i|\mu}^{\eta \beta} + \left[ \Gamma_{i|\mu}^\eta, T_{0i}^\beta \right] - \left[ \Gamma_{i|\mu}^\beta, T_{0i}^\eta \right]  - \sum_{\xi=1}^p g_{i\eta \beta}^\xi \Gamma_{i|\mu}^\xi \right) - \sum_{\eta=1}^p g_{0i\beta\gamma}^\eta\left( \Pi_{i|\mu}^{\eta \alpha} +\left[ \Gamma_{i|\mu}^\eta, T_{0i}^\alpha \right] - \left[ \Gamma_{i|\mu}^\alpha, T_{0i}^\eta  \right] - \sum_{\xi=1}^p g_{0i\eta\alpha}^\xi \Gamma_{i|\mu}^\xi \right)\\
&= \left[ T_{0i}^\alpha, \Pi_{i|\mu}^{\beta \gamma} \right]+ \left[ T_{0i}^\alpha, \left[ \Gamma_{i|\mu}^\beta, T_{0i}^\gamma \right]  \right] - \left[ T_{0i}^\alpha, \left[ \Gamma_{i|\mu}^\gamma, T_{0i}^\beta \right] \right] - \sum_{\eta=1}^p T_{0i}^\alpha \left( g_{0i\beta\gamma}^\eta \right) \Gamma_{i|\mu}^\eta - \sum_{\eta=1}^p g_{0i\beta\gamma}^\eta \left[  T_{0i}^\alpha, \Gamma_{i|\mu}^\eta \right] \\
&- \left[ T_{0i}^\beta, \Pi_{i|\mu}^{\alpha \gamma} \right] - \left[ T_{0i}^\beta, \left[ \Gamma_{i|\mu}^\alpha, T_{0i}^\gamma \right] \right] + \left[ T_{0i}^\beta, \left[ \Gamma_{i|\mu}^\gamma, T_{0i}^\alpha \right] \right] + \sum_{\eta=1}^p T_{0i}^\beta\left( g_{0i\alpha \gamma}^\eta \right) \Gamma_{i|\mu}^\eta + \sum_{\eta=1}^p g_{0i\alpha \gamma}^\eta\left[ T_{0i}^\beta, \Gamma_{i|\mu}^\eta \right]\\
&+ \left[ T_{0i}^\gamma , \Pi_{i|\mu}^{\alpha\beta } \right] +\left[ T_{0i}^\gamma, \left[ \Gamma_{i|\mu}^\alpha, T_{0i}^\beta \right] \right] - \left[ T_{0i}^\gamma, \left[  \Gamma_{i|\mu}^\beta, T_{0i}^\alpha\right]\right] - \sum_{\eta=1}^p T_{0i}^\gamma \left( g_{0i\alpha\beta}^\eta \right) \Gamma_{i|\mu}^\eta - \sum_{\eta=1}^p g_{0i\alpha\beta}^\eta \left[ T_{0i}^\gamma, \Gamma_{i|\mu}^\eta \right]\\
&-\sum_{\eta=1}^p g_{0i\alpha\beta}^\eta \Pi_{i|\mu}^{\eta\gamma} - \sum_{\eta=1}^p g_{0i\alpha\beta}^\eta\left[ \Gamma_{i|\mu}^\eta, T_{0i}^\gamma \right] + \sum_{\eta=1}^p g_{0i\alpha\beta}^\eta \left[ \Gamma_{i|\mu}^\gamma, T_{0i}^\eta \right] + \sum_{\eta,\xi=1}^p g_{0i\alpha\beta}^\eta g_{0i\eta\gamma}^\xi \Gamma_{i|\mu}^\xi\\
&+ \sum_{\eta=1}^p g_{0i\alpha\gamma}^\eta \Pi_{i|\mu}^{\eta \beta}+\sum_{\eta=1}^p g_{0i\alpha \gamma}^\eta \left[ \Gamma_{i|\mu}^\eta,  T_{0i}^\beta  \right] - \sum_{\eta=1}^p g_{0i\alpha \gamma}^\eta \left[ \Gamma_{i|\mu}^\beta, T_{0i}^\eta  \right] - \sum_{\eta, \xi=1}^p g_{0i\alpha \gamma}^\eta g_{0i\eta \beta}^\xi \Gamma_{i|\mu}^\xi \\
& - \sum_{\eta=1}^p g_{0i\beta \gamma}^\eta \Pi_{i|\mu}^{\eta \alpha} - \sum_{\eta=1}^p g_{0i\beta\gamma}^\eta \left[ \Gamma_{i|\mu}^\eta, T_{0i}^\alpha \right] + \sum_{\eta=1}^p g_{0i\beta\gamma}^\eta \left[ \Gamma_{i|\mu}^\alpha, T_{0i}^\eta \right] + \sum_{\eta, \xi=1}^p g_{0i\beta\gamma}^\eta g_{0i\eta\alpha}^\xi \Gamma_{i|\mu}^\xi\\
&= - \sum_{\xi=1}^p G_{i|\mu}^{\alpha \beta\gamma \xi} T_{0i}^\xi + \left[\Gamma_{i|\mu}^\beta, \left[ T_{0i}^\alpha, T_{0i}^\gamma  \right] \right] -\left[ \Gamma_{i|\mu}^\gamma ,  \left[ T_{0i}^\alpha, T_{0i}^\beta\right] \right] - \left[ \Gamma_{i|\mu}^\alpha, \left[ T_{0i}^\beta, T_{0i}^\gamma \right] \right]\\
&  + \sum_{\eta=1}^p g_{0i\alpha\beta}^\eta \left[ \Gamma_{i|\mu}^\gamma, T_{0i}^\eta\right] - \sum_{\eta=1}^p g_{0i\alpha \gamma}^\eta \left[ \Gamma_{i|\mu}^\beta, T_{0i}^\eta \right] + \sum_{\eta=1}^p g_{0i\beta\gamma}^\eta \left[ \Gamma_{i|\mu}^\alpha, T_{0i}^\eta \right]\\
&=- \sum_{\xi=1}^p G_{i|\mu}^{\alpha\beta \gamma \xi} T_{0i}^\xi + \sum_{\eta=1}^p \left[ \Gamma_{i|\mu}^{\beta}, g_{0i\alpha \gamma}^\eta T_{0i}^\eta \right] - \sum_{\eta=1}^p \left[\Gamma_{i|\mu}^\gamma, g_{0i\alpha\beta}^\eta T_{0i}^\eta \right] - \sum_{\eta=1}^p \left[ \Gamma_{i|\mu}^\alpha , g_{0i\beta\gamma}^\eta T_{0i}^\eta \right] \\
&+ \sum_{\eta=1}^p g_{0i\alpha\beta}^\eta \left[ \Gamma_{i|\mu}^\gamma, T_{0i}^\eta\right] - \sum_{\eta=1}^p g_{0i\alpha \gamma}^\eta \left[ \Gamma_{i|\mu}^\beta, T_{0i}^\eta \right] + \sum_{\eta=1}^p g_{0i\beta\gamma}^\eta \left[ \Gamma_{i|\mu}^\alpha, T_{0i}^\eta \right]\\
&= - \sum_{\xi=1}^p G_{i|\mu}^{\alpha\beta \gamma \xi} T_{0i}^\xi + \sum_{\eta=1}^p \Gamma_{i|\mu}^\beta\left( g_{0i\alpha \gamma}^\eta \right) T_{0i}^\eta - \sum_{\eta=1}^p \Gamma_{i|\mu}^\gamma\left( g_{0i\alpha\beta}^\eta \right) T_{0i}^\eta - \sum_{\eta=1}^p \Gamma_{i|\mu}^\alpha \left( g_{0i\beta\gamma}^\eta \right) T_{0i}^\eta \equiv_{\Theta_{\mathcal{F}_0}} 0
\end{align*}}}
Here $A\equiv_{\Theta_{\mathcal{F}_0}} B$ means that $A-B\in \Gamma\left(U_i, \mathcal{A}^{0,0}\left( \Theta_{\mathcal{F}_0} \right) \right)$.

This proves $\left( \bar{B}_\mu , \bar{\Phi}_\mu, -\xi_\mu \right)$ defines a $2$-cocycle in the above Dolbeualt resolution of $\Theta_{\mathcal{F}_0}^\bullet$. Then by hypothesis $\mathbb{H}^2\left( M, \Theta_{\mathcal{F}_0}^{\bullet} \right)=0$, there exists
\begin{align*}
\left(\overline{T_\mu'}, \varphi_\mu' \right)\in \frac{A^{0,0}\left(M, \mathscr{H}om_{\mathcal{O}_M}\left( \Theta_{\mathcal{F}_0}, \Theta_M \right) \right)}{ A^{0,0}\left( M, \mathscr{H}om_{\mathcal{O}_M}\left(\Theta_{\mathcal{F}_0} , \Theta_{\mathcal{F}_0} \right) \right)} \bigoplus A^{0,1}\left(M, \Theta_M \right)
\end{align*}
such that
\begin{align}
-\bar{\partial} \varphi_\mu'&=-\xi_\mu \label{te30}\\
-\bar{\partial} \overline{T_\mu'}+ \hat{D}_0\left(\varphi_\mu' \right)&=\overline{\Phi}_\mu= \overline{\left\{ T_{0i}^\alpha\mapsto  \Phi_{i|\mu}^\alpha - \bar{\partial} \Gamma_{i|\mu}^\alpha + \sum_{\xi=1}^p \Lambda_{i|\mu}^{\alpha \xi} T_{0i}^\xi  \right\}} \label{te31} \\
\hat{D}_1\left(\overline{T_\mu' }\right)=\overline{B}_\mu&=\overline{\left\{  T_{0i}^\alpha \wedge T_{0i}^\beta \mapsto  \Pi_{i|\mu}^{\alpha\beta}+\left[ \Gamma_{i|\mu}^\alpha, T_i^\beta \right] -\left[ \Gamma_{i|\mu}^\beta, T_i^\alpha  \right]- \sum_{\gamma=1}^p g_{i\alpha\beta}^\gamma \Gamma_{i|\mu}^\gamma + \sum_{\gamma=1}^p \Psi_{i|\mu}^{\alpha\beta\gamma} T_{0i}^\gamma   \right\}} \label{te32}
\end{align}
where $\overline{T'_\mu}$ is the image of the global section $T_\mu'=\left\{ T_{i|\mu}' \right\} \in A^{0,0}\left(M, \mathscr{H}om_{\mathcal{O}_M} \left( \Theta_{\mathcal{F}_0}, \Theta_M   \right)\right)$ in $\frac{A^{0,0}\left(M, \mathscr{H}om_{\mathcal{O}_M}\left( \Theta_{\mathcal{F}_0}, \Theta_M \right) \right)}{ A^{0,0}\left( M, \mathscr{H}om_{\mathcal{O}_M}\left(\Theta_{\mathcal{F}_0} , \Theta_{\mathcal{F}_0} \right) \right)}$ such that $T_{i|\mu}':\Gamma\left(U_i, \Theta_{\mathcal{F}_0} \right)\to \Gamma\left( U_i, \mathcal{A}^{0,0}\left(\Theta_M \right)\right), T_{0i}^\alpha \mapsto T_{i|\mu}'^\alpha,\alpha=1,...,p$. In particular, we have
\begin{align}\label{te38}
T_{i|\mu}'^\alpha= \sum_{\beta=1}^p r_{0ij}^{\alpha\beta} T_{j|\mu}'^\beta
\end{align}

From $(\ref{te30})$, if we take 
\begin{align} \label{te33}
\varphi_\mu:=-\varphi_\mu',
\end{align}
then $(\ref{te11})_\mu$ is satisfied. From $(\ref{te31})$, there exist $\left\{ b_{i|\mu} \right\}\in C^0\left(\mathcal{U}, \mathcal{A}^{0,1}\left( \mathscr{H}om_{\mathcal{O}_M}\left( \Theta_{\mathcal{F}_0}, \Theta_{\mathcal{F}_0} \right) \right) \right)$ defined by $b_{i|\mu}:\Gamma\left( U_i, \Theta_{\mathcal{F}_0} \right)\to \Gamma\left(U_i,  \mathcal{A}^{0,1}\left( \Theta_{\mathcal{F}_0} \right) \right), T_{0i}^\alpha\mapsto \sum_{\xi=1}^p b_{i|\mu}^{\alpha \xi} T_{0i}^\xi$ such that
\begin{align} \label{d46}
\sum_{\beta,\xi=1}^p r_{0ij}^{\alpha\beta} b_{j|\mu}^{\beta \xi} T_{0j}^\xi - \sum_{\xi=1}^p b_{i|\mu}^{\alpha \xi} T_{0i}^\xi= \sum_{\beta=1}^p \left[\varphi_\mu', r_{0ij}^{\alpha\beta} \right] T_{0j}^\beta
\end{align}
and
\begin{align}
 - \bar{\partial} T_{i|\mu}'^\alpha +\left[ \varphi_\mu', T_{0i}^\alpha \right] + \sum_{\xi=1}^p b_{i|\mu}^{\alpha \xi} T_{0i}^\xi &= \Phi_{i|\mu}^\alpha - \bar{\partial} \Gamma_{i|\mu}^\alpha + \sum_{\xi=1}^p \Lambda_{i|\mu}^{\alpha \xi} T_{0i}^\xi \label{t50}
\end{align}
By taking $\bar{\partial}$ on $(\ref{t50})$,  we have from $(\ref{te30})$ and $(\ref{t22})$
\begin{align*}
\left[ \bar{\partial} \varphi_\mu' , T_{0i}^\alpha \right]+ \sum_{\xi=1}^p \bar{\partial} b_{i|\mu}^{\alpha \xi} T_{0i}^\xi = \bar{\partial} \Phi_{i|\mu}^\alpha + \sum_{\xi=1}^p \bar{\partial} \Gamma_{i|\mu}^\alpha T_{0i}^\xi \Longrightarrow \sum_{\xi=1}^p \bar{\partial}\left( b_{i|\mu}^{\alpha \xi} - \Lambda_{i|\mu}^{\alpha \xi} \right) T_{0i}^\xi=0 \Longrightarrow \bar{\partial}\left( b_{i|\mu}^{\alpha \xi} - \Lambda_{i|\mu}^{\alpha \xi} \right)=0,
\end{align*}
so that there exist $c_{i|\mu}^{\alpha \xi}\in \Gamma\left( U_i , \mathcal{A}^{0,0}  \right)$ such that $\bar{\partial} c_{i|\mu}^{\alpha \xi}= b_{i|\mu}^{\alpha \xi} - \Lambda_{i| \mu}^{\alpha \xi}$. Then if we take
\begin{align}\label{te37}
T_{i|\mu}^\alpha:= \Gamma_{i|\mu}^\alpha - T_{i|\mu}'^\alpha + \sum_{\xi=1}^p c_{i|\mu}^{\alpha \xi} T_{0i}^\xi,
\end{align}
then we have from $(\ref{te33})$ and $(\ref{t50})$
\begin{align} \label{te40}
\bar{\partial}\left(\Gamma_{i|\mu}^\alpha- T_{i|\mu}'^\alpha +\sum_{\xi=1}^p c_{i|\mu}^{\alpha \xi} T_{0i}^\xi \right) - \left[ - \varphi_\mu', T_{0i}^\alpha \right] = \Phi_{i|\mu}^\alpha \Longrightarrow \bar{\partial} T_{i|\mu}^\alpha - \left[ \varphi_\mu, T_{0i}^\alpha \right] = \Phi_{i|\mu}^\alpha,
\end{align}
so that $(\ref{te12})_\mu$ is satisfied. We will find $r_{ij|\mu}^{\alpha\beta}$ satisfying $(\ref{te15})_\mu$. We note that from $(\ref{te37})$ and $(\ref{t57})$ and $(\ref{te38})$, 
\begin{align}\label{te39}
&-\sum_{\beta=1}^p r_{0ij}^{\alpha\beta} T_{j|\mu}^\beta + T_{i|\mu}^\alpha + \Gamma_{ij|\mu}^\alpha \\
&= -\sum_{\beta=1}^p r_{0ij}^{\alpha\beta}\left( \Gamma_{j|\mu}^\beta - T_{j|\mu}'^\beta + \sum_{\xi=1}^p c_{j|\mu}^{\beta \xi} T_{0j}^\xi \right) + \left( \Gamma_{i|\mu}^\alpha - T_{i|\mu}'^\alpha + \sum_{\xi=1}^p c_{i|\mu}^{\alpha \xi}  T_{0i}^\xi \right) + \sum_{\beta=1}^p r_{0ij}^{\alpha\beta} \Gamma_{j|\mu}^\beta - \Gamma_{i|\mu}^\alpha + \sum_{\xi=1}^p \lambda_{ij|\mu}^{\alpha \xi} T_{0i}^\xi \notag \\
&= \sum_{\xi=1}^p\left(  \sum_{\eta=1}^p \lambda_{ij|\mu}^{\alpha \eta} r_{0ij}^{\eta\xi} + \sum_{\eta=1}^p c_{i|\mu}^{\alpha \eta}  r_{0ij}^{\eta\xi} - \sum_{\beta=1}^p r_{0ij}^{\alpha\beta} c_{j|\mu}^{\beta\xi} \right) T_{0j}^\xi \notag
\end{align}
If we take
\begin{align}\label{te41}
r_{ij|\mu}^{\alpha \xi} :=  \sum_{\eta=1}^p \lambda_{ij|\mu}^{\alpha \eta} r_{0ij}^{\eta\xi}+ \sum_{\eta=1}^p c_{i|\mu}^{\alpha \eta} r_{0ij}^{\eta\xi} - \sum_{\beta=1}^p r_{0ij}^{\alpha\beta} c_{j|\mu}^{\beta\xi},
\end{align}
then $(\ref{te39})$ implies $(\ref{te15})_\mu$, i.e.
\begin{align}\label{te45}
\sum_{\beta=1}^p r_{0ij}^{\alpha\beta} T_{j|\mu}^{\beta}- T_{i|\mu}^\alpha+&\sum_{\beta=1}^p r_{ij|\mu}^{\alpha\beta} T_{0j}^\beta=\Gamma_{ij|\mu}^\alpha
\end{align}

We check $(\ref{te13})_\mu$. In fact, from $(\ref{t26})$, $(\ref{te40}), (\ref{t57})$ and $(\ref{te41})$, we have
{\small{\begin{align*}
&\sum_{\eta=1}^p \Lambda_{ij|\mu}^{\alpha \eta} T_{0i}^\eta= - \sum_{\beta=1}^p r_{0ij}^{\alpha\beta} \Phi_{j|\mu}^\beta + \Phi_{i|\mu}^\alpha + \bar{\partial} \Gamma_{ij|\mu}^\alpha \\
& = - \sum_{\beta=1}^p r_{0ij}^{\alpha\beta}\left(\bar{\partial} \left(\Gamma_{j|\mu}^\beta - T_{j|\mu}'^\beta + \sum_{\xi=1}^p c_{j|\mu}^{\beta \xi} T_{0j}^\xi \right) - \left[ \varphi_\mu, T_{0j}^\beta\right] \right) + \bar{\partial} \left(\Gamma_{i|\mu}^\alpha - T_{i|\mu}'^\alpha + \sum_{\xi=1}^p c_{i|\mu}^{\alpha \xi} T_{0i}^\xi \right) - \left[ \varphi_\mu , T_{0i}^\alpha \right] \\
&+ \sum_{\beta=1}^p r_{0ij}^{\alpha\beta} \bar{\partial} \Gamma_{j|\mu}^\beta - \bar{\partial} \Gamma_{i|\mu}^\alpha + \sum_{\xi=1}^p \bar{\partial} \lambda_{ij|\mu}^{\alpha \xi} T_{0i}^\xi = \bar{\partial}\left(- \sum_{\beta,\xi=1}^p r_{0ij}^{\alpha\beta} c_{j|\mu}^{\beta \xi} T_{0j}^\xi  + \sum_{\xi=1}^p \bar{\partial} c_{i|\mu}^{\alpha \xi} T_{0i}^\xi + \sum_{\xi=1}^p \lambda_{ij|\mu}^{\alpha \xi} T_{0i}^\xi     \right) - \sum_{\xi=1}^p\left[ \varphi_\mu, r_{0ij}^{\alpha\xi} \right] T_{0j}^\xi\\
&= \sum_{\xi=1}^p \left( \bar{\partial} r_{ij|\mu}^{\alpha \xi} - \left[ \varphi_\mu, r_{0ij}^{\alpha \xi} \right] \right)  T_{0j}^\xi
\end{align*}}}

We check $(\ref{te17})_\mu$. In fact, from $(\ref{te18})$ and $(\ref{te41})$, we have
{\small{\begin{align*}
&\sum_{\beta,\gamma=1}^p r_{0ij}^{\alpha\beta} r_{jk|\mu}^{\beta\gamma} r_{0ki}^{\gamma\xi}-\sum_{\gamma=1}^p r_{ik|\mu}^{\alpha\gamma} r_{0ki}^{\gamma\xi}+\sum_{\beta=1}^p r_{ij|\mu}^{\alpha\beta} r_{0ji}^{\beta\xi}\\
&=\sum_{\beta, \gamma=1}^p r_{0ij}^{\alpha\beta}\left( \sum_{\eta=1}^p \lambda_{jk|\mu}^{\beta\eta} r_{0jk}^{\eta \gamma} +\sum_{\eta=1}^p c_{j|\mu}^{\beta \eta} r_{0jk}^{\eta\gamma} - \sum_{\eta=1}^p r_{0jk}^{\beta\eta} c_{k|\mu}^{\eta \gamma} \right) r_{0ki}^{\gamma \xi} - \sum_{\gamma=1}^p\left( \sum_{\eta=1}^p \lambda_{ik|\mu}^{\alpha \eta} r_{0ik}^{\eta \gamma} + \sum_{\eta=1}^p c_{i|\mu}^{\alpha \eta} r_{0ik}^{\eta\gamma} - \sum_{\eta=1}^p  r_{0ij}^{\alpha \eta} c_{k|\mu}^{\eta \gamma} \right) r_{0ki}^{\gamma \xi} \\
& + \sum_{\beta=1}^p \left( \sum_{\eta=1}^p \lambda_{ij|\mu}^{\alpha \beta} r_{0ij}^{\eta\beta} +\sum_{\eta=1}^p c_{i|\mu}^{\alpha \eta} r_{0ij}^{\eta \beta} - \sum_{\eta=1}^p r_{0ij}^{\alpha \eta} c_{j|\mu}^{\eta\beta}   \right) r_{0ji}^{\beta\xi}  = \lambda_{ijk|\mu}^{\alpha \xi}
\end{align*}}}

On the other hand, from $(\ref{te32})$, there exist $\left\{ W_{i|\mu} \right\} \in C^0\left(\mathcal{U}, \mathcal{A}^{0,0}\left( \mathscr{H}om_{\mathcal{O}_M} \left( \bigwedge^2 \Theta_{\mathcal{F}_0}, \Theta_{\mathcal{F}_0}    \right)    \right) \right)$ defined by $W_{i|\mu}:\Gamma\left( U_i, \bigwedge^2 \Theta_{\mathcal{F}_0} \right)\to \Gamma\left( U_i , \Theta_{\mathcal{F}_0} \right), T_{0i}^\alpha \wedge T_{0i}^\beta \mapsto \sum_{\xi=1}^p W_{i|\mu}^{\alpha \beta \xi} T_{0i}^\xi$ such that
\begin{align} \label{d49}
\sum_{\gamma,\eta, \xi=1}^p r_{0ij}^{\alpha \gamma} r_{0ij}^{\beta\eta} W_{j|\mu}^{\gamma\eta\xi}T_{0j}^\xi - \sum_{\xi=1}^p W_{i|\mu}^{\alpha\beta \xi} T_{0i}^\xi =\sum_{\eta=1}^p \left[  T_{i|\mu}'^\alpha, r_{0ij}^{\beta\eta} \right] T_{0j}^\eta - \sum_{\gamma=1}^p \left[ T_{i|\mu}'^\beta, r_{0ij}^{\alpha\gamma} \right] T_{0j}^\gamma
\end{align}
and
{\small{\begin{align}
&\left[T_{0i}^\alpha, T_{i|\mu}'^\beta \right] - \left[ T_{0i}^\beta, T_{i|\mu}'^\alpha  \right] -\sum_{\gamma=1}^p g_{0i\alpha\beta}^\gamma T_{i|\mu}'^\gamma +\sum_{\gamma=1}^p W_{i|\mu}^{\alpha\beta \gamma }T_{0i}^\gamma=\Pi_{i|\mu}^{\alpha\beta} +\left[ \Gamma_{i|\mu}^\alpha, T_{0i}^\beta \right] -\left[ \Gamma_{i|\mu}^\beta, T_{0i}^\alpha \right] - \sum_{\gamma=1}^p g_{0i\alpha\beta}^\gamma \Gamma_{i|\mu}^\gamma + \sum_{\gamma=1}^p \Psi_{i|\mu}^{\alpha\beta\gamma}  T_{0i}^\gamma \notag \\
&\iff -\left[\Gamma_{i|\mu}^\alpha - T_{i|\mu}'^\alpha, T_{0i}^\beta \right] - \left[ T_{0i}^\alpha , \Gamma_{i|\mu}^\beta -  T_{i|\mu}'^\beta \right] + \sum_{\gamma=1}^p g_{0i\alpha\beta}^\gamma \left( \Gamma_{i|\mu}^\gamma - T_{i|\mu}'^\gamma  \right) + \sum_{\gamma=1}^p \left(  W_{i|\mu}^{\alpha\beta\gamma} - \Psi_{i|\mu}^{\alpha\beta\gamma} \right) T_{0i}^\gamma = \Pi_{i|\mu}^{\alpha\beta} \label{te35}
\end{align}}}
We will find $g_{i\alpha\beta|\mu}^\gamma$ satisfying $(\ref{te16})_\mu$. In fact, from $(\ref{te35})$, we have
\begin{align}\label{te36}
\Pi_{i|\mu}^{\alpha\beta}=&-\left[\Gamma_{i|\mu}^\alpha - T_{i|\mu}'^\alpha + \sum_{\xi=1}^p  c_{i|\mu}^{\alpha \xi} T_{0i}^\xi , T_{0i}^\beta \right]  - \left[ T_{0i}^\alpha, \Gamma_{i|\mu}^\beta- T_{i|\mu}'^\beta +\sum_{\xi=1}^p c_{i|\mu}^{\beta\xi} T_{0i}^\xi \right] + \sum_{\gamma=1}^p g_{0i\alpha\beta}^\gamma\left( \Gamma_{i|\mu}^\gamma- T_{i|\mu}'^\gamma + \sum_{\xi=1}^p c_{i|\mu}^{\gamma \xi} T_{0i}^\xi \right)\\
                   &+\left[ \sum_{\xi=1}^p c_{i|\mu}^{\alpha \xi} T_{0i}^\xi, T_{0i}^\beta \right] + \left[  T_i^\alpha, \sum_{\xi=1}^p c_{i|\mu}^{\beta \xi} T_i^\xi  \right] - \sum_{\gamma, \xi=1}^p g_{0i\alpha\beta}^\gamma c_{i|\mu}^{\gamma \xi} T_{0i}^\xi + \sum_{\xi=1}^p \left( W_{i|\mu}^{\alpha \beta\xi}- \Psi_{i|\mu}^{\alpha\beta\xi} \right)T_{0i}^\xi \notag\\
                   &=- \left[ T_{i|\mu}, T_{0i}^\beta\right] - \left[ T_{0i}^\alpha, T_{i|\mu}^\beta \right] + \sum_{\gamma=1}^p g_{0i\alpha\beta}^\gamma T_{i|\mu}^\gamma \notag \\
                   & +\sum_{\xi=1}^p\left( \sum_{\gamma=1}^p c_{i|\mu}^{\alpha\gamma} g_{0i\gamma \beta}^\xi - T_{0i}^\beta\left(  c_{i|\mu}^{\alpha \xi} \right) + T_{0i}^\alpha \left( c_{i|\mu}^{\beta\xi} \right) + \sum_{\gamma=1}^p c_{i|\mu}^{\beta\gamma} g_{0i\alpha\gamma}^\xi - \sum_{\gamma=1}^p g_{0i\alpha\beta}^\gamma c_{i|\mu}^{\gamma \xi} + W_{i|\mu}^{\alpha\beta\xi}- \Psi_{i|\mu}^{\alpha\beta \xi} \right) T_{0i}^\xi  \notag
\end{align}   
If we take
\begin{align}\label{te42}
g_{i\alpha\beta |\mu}^\xi&:= \sum_{\gamma=1}^p c_{i|\mu}^{\alpha\gamma} g_{0i\gamma \beta}^\xi - T_{0i}^\beta\left(  c_{i|\mu}^{\alpha \xi} \right) + T_{0i}^\alpha \left( c_{i|\mu}^{\beta\xi} \right) + \sum_{\gamma=1}^p c_{i|\mu}^{\beta\gamma} g_{0i\alpha\gamma}^\xi - \sum_{\gamma=1}^p g_{0i\alpha\beta}^\gamma c_{i|\mu}^{\gamma \xi} + W_{i|\mu}^{\alpha\beta\xi}- \Psi_{i|\mu}^{\alpha\beta \xi}, 
\end{align}
from $(\ref{te36})$, we have
\begin{align}\label{te43}
\Pi_{i|\mu}^{\alpha\beta}= - \left[ T_{i|\mu}^\alpha, T_{0i}^\beta\right] - \left[ T_{0i}^\alpha, T_{i|\mu}^\beta \right] + \sum_{\gamma=1}^p g_{0i\alpha\beta}^\gamma T_{i|\mu}^\gamma + \sum_{\xi=1}^p g_{i\alpha\beta|\mu}^\xi T_{0i}^\xi,
\end{align}
so that $(\ref{te16})_\mu$ is satisfied. We check $(\ref{te14})_\mu$. In fact, from $(\ref{t28})$ and $(\ref{te43})$, we have
\begin{align*}
&\sum_{\gamma=1}^p \eta_{i\alpha\beta|\mu}^\gamma T_{0i}^\gamma = \bar{\partial} \Pi_{i|\mu}^{\alpha\beta} + \left[ \Phi_{i|\mu}^\alpha, T_{0i}^\beta \right] - \left[ \Phi_{i|\mu}^\beta , T_{0i}^\alpha \right] - \sum_{\gamma=1}^p g_{0i\alpha\beta}^\gamma \Phi_{i|\mu}^\gamma \\
&= - \left[\bar{\partial} T_{i|\mu}^\alpha, T_{0i}^\beta   \right] + \left[ \bar{\partial} T_{i|\mu}^\beta, T_{0i}^\alpha \right] + \sum_{\gamma=1}^p g_{0i\alpha\beta}^\gamma \bar{\partial} T_{i|\mu}^\gamma + \sum_{\gamma=1}^p \bar{\partial} g_{i\alpha\beta|\mu}^\gamma T_{0i}^\gamma\\
& + \left[ \bar{\partial} T_{i|\mu}^\alpha- \left[ \varphi_\mu , T_{0i}^\alpha \right] , T_{0i}^\beta       \right] - \left[ \bar{\partial} T_{i|\mu}^\beta - \left[ \varphi_\mu, T_{0i}^\beta \right], T_{0i}^\alpha \right] - \sum_{\gamma=1}^p g_{0i\alpha\beta}^\gamma \left( \bar{\partial} T_{i|\mu}^\gamma - \left[ \varphi_\mu, T_{0i}^\gamma \right] \right) \\
&= -  \left[ \left[ \varphi_\mu , T_{0i}^\alpha \right] , T_{0i}^\beta \right] + \left[ \left[ \varphi_\mu , T_{0i}^\beta \right], T_{0i}^\alpha \right] + \sum_{\gamma=1}^p g_{0i\alpha\beta}^\gamma \left[ \varphi_\mu, T_{0i}^\gamma \right] + \sum_{\gamma=1}^p \bar{\partial } g_{i\alpha\beta|\mu}^\gamma T_{0i}^\gamma  = \sum_{\gamma=1}^p \left( \bar{\partial} g_{i\alpha\beta|\mu}^\gamma - \left[ \varphi_\mu, g_{0i\alpha\beta}^\gamma \right] \right) T_{0i}^\gamma ,
\end{align*}
so that $(\ref{te14})_\mu$ is satisfied. It remains to check $(\ref{t31})_\mu$ and $(\ref{t33})_\mu$. Let us check $(\ref{t31})_\mu$. In fact, from $(\ref{t27})$ and $(\ref{te43})$ and $(\ref{te45})$ and $(\ref{te19})$, we have
{\small{\begin{align*}
&-\sum_{\xi=1}^p \Psi_{ij|\mu}^{\alpha\beta \xi} T_{0i}^\xi=-\Pi_{i|\mu}^{\alpha\beta}+\sum_{\eta,\xi=1}^p r_{0ij}^{\alpha\eta} r_{0ij}^{\beta \xi} \Pi_{j|\mu}^{\eta\xi} +\left[\Gamma_{ij|\mu}^\alpha, T_{0i}^\beta \right]+\left[ T_{0i}^\alpha, \Gamma_{ij|\mu}^\beta \right]-\sum_{\gamma=1}^p g_{0i\alpha\beta}^\gamma\Gamma_{ij|\mu}^\gamma+\sum_{\xi=1}^p \left( -\Gamma_{ij|\mu}^\alpha\left(r_{0ij}^{\beta \xi} \right) + \Gamma_{ij|\mu}^\beta \left( r_{0ij}^{\alpha \xi} \right)  \right)T_{0j}^\xi \\
&= -\Pi_{i|\mu}^{\alpha\beta}+\sum_{\eta,\xi=1}^p r_{0ij}^{\alpha\eta} r_{0ij}^{\beta \xi} \Pi_{j|\mu}^{\eta\xi} +\sum_{\xi=1}^p r_{0ij}^{\beta \xi} \left[\Gamma_{ij|\mu}^\alpha, T_{0j}^\xi \right] - \sum_{\xi=1}^p  r_{0ij}^{\alpha \xi}\left[  \Gamma_{ij|\mu}^\beta , T_{0j}^\xi \right]  -\sum_{\gamma=1}^p g_{0i\alpha\beta}^\gamma\Gamma_{ij|\mu}^\gamma \\
&=  \left[ T_{i|\mu}^\alpha, T_{0i}^\beta\right] + \left[ T_{0i}^\alpha, T_{i|\mu}^\beta \right] - \sum_{\gamma=1}^p \cancel{g_{0i\alpha\beta}^\gamma T_{i|\mu}^\gamma} - \sum_{\xi=1}^p g_{i\alpha\beta|\mu}^\xi T_{0i}^\xi + \sum_{\eta, \xi=1}^p r_{0ij}^{\alpha \eta} r_{0ij}^{\beta \xi}\left( - \left[ T_{j|\mu}^\eta, T_{0j}^\xi\right] - \left[ T_{0j}^\eta, T_{j|\mu}^\xi \right] + \sum_{\gamma=1}^p g_{0j\eta\xi}^\gamma T_{j|\mu}^\gamma + \sum_{\gamma=1}^p g_{j\eta\xi|\mu}^\gamma T_{0j}^\gamma \right)\\
&+\sum_{\xi=1}^p r_{0ij}^{\beta \xi}\left( \sum_{\eta=1}^p \left[ r_{0ij}^{\alpha\eta} T_{j|\mu}^{\eta} , T_{0j}^\xi \right]- \left[ T_{i|\mu}^\alpha , T_{0j}^\xi \right]+\sum_{\eta=1}^p \left[r_{ij|\mu}^{\alpha\eta} T_{0j}^\eta, T_{0j}^\xi \right]  \right)- \sum_{\xi=1}^p r_{0ij}^{\alpha \xi} \left(  \sum_{\eta=1}^p \left[r_{0ij}^{\beta\eta} T_{j|\mu}^{\eta}, T_{0j}^\xi \right]- \left[ T_{i|\mu}^\beta , T_{0j}^\xi \right] + \sum_{\eta=1}^p \left[ r_{ij|\mu}^{\beta\eta} T_{0j}^\eta , T_{0j}^\xi \right] \right) \\
& - \sum_{\gamma=1}^p g_{0i\alpha\beta}^\gamma\left(\sum_{\xi=1}^p r_{0ij}^{\gamma\xi} T_{j|\mu}^{\xi}- \cancel{T_{i|\mu}^\gamma}+\sum_{\xi=1}^p r_{ij|\mu}^{\gamma \xi} T_{0j}^\xi \right)\\
&=\sum_{\xi=1}^p T_{i|\mu}^\alpha\left( r_{0ij}^{\beta \xi} \right) T_{0j}^\xi - \sum_{\xi=1}^p T_{i|\mu}^\beta\left( r_{0ij}^{\alpha \xi} \right) T_{0j}^\xi  - \sum_{\xi=1}^p g_{i\alpha\beta|\mu}^\xi T_{0i}^\xi - \sum_{\xi, \eta=1}^p  r_{0ij}^{\beta \xi} T_{0j}^\xi\left(r_{0ij}^{\alpha \eta} \right) T_{j|\mu}^\eta + \sum_{\xi, \eta=1}^p r_{0ij}^{\alpha \xi} T_{0j}^\xi \left( r_{0ij}^{\beta \eta} \right) T_{j|\mu}^\eta + \sum_{\eta, \xi, \gamma=1}^p r_{0ij}^{\alpha \eta} r_{0ij}^{\beta \xi} g_{0j\eta \xi}^\gamma T_{j|\mu}^\gamma \\
&+ \sum_{\eta, \xi, \gamma=1}^p r_{0ij}^{\alpha \eta} r_{0ij}^{\beta\xi} g_{j\eta \xi|\mu}^\gamma T_{0j}^\gamma + \sum_{\xi, \eta, \gamma=1}^p r_{0ij}^{\beta \xi} r_{ij|\mu}^{\alpha \eta} g_{0j\eta \xi}^\gamma T_{0j}^\gamma - \sum_{ \eta=1}^p T_{0i}^\beta  \left( r_{ij|\mu}^{\alpha \eta} \right) T_{0j}^\eta - \sum_{\xi, \eta, \gamma=1}^p r_{0ij}^{\alpha\xi} r_{ij|\mu}^{\beta \eta} g_{0j\eta \xi}^\gamma T_{0j}^\gamma + \sum_{\eta=1}^p T_{0i}^\alpha\left( r_{ij|\mu}^{\beta \eta} \right) T_{0j}^\eta \\
& - \sum_{\gamma, \xi=1}^p g_{0i\alpha\beta}^\gamma r_{0ij}^{\gamma \xi} T_{j|\mu}^\xi - \sum_{\gamma=1}^p g_{i\alpha\beta}^\gamma r_{ij|\mu}^{\gamma \xi} T_{0j}^\xi\\
&=\sum_{\xi=1}^p T_{i|\mu}^\alpha\left( r_{0ij}^{\beta \xi} \right) T_{0j}^\xi - \sum_{\xi=1}^p T_{i|\mu}^\beta\left( r_{0ij}^{\alpha \xi} \right) T_{0j}^\xi  - \sum_{\xi=1}^p g_{i\alpha\beta|\mu}^\xi T_{0i}^\xi + \sum_{\eta, \xi, \gamma=1}^p r_{0ij}^{\alpha \eta} r_{0ij}^{\beta\xi} g_{j\eta \xi|\mu}^\gamma T_{0j}^\gamma + \sum_{\xi, \eta, \gamma=1}^p r_{0ij}^{\beta \xi} r_{ij|\mu}^{\alpha \eta} g_{0j\eta \xi}^\gamma T_{0j}^\gamma \\
& - \sum_{ \eta=1}^p T_{0i}^\beta  \left( r_{ij|\mu}^{\alpha \eta} \right) T_{0j}^\eta - \sum_{\xi, \eta, \gamma=1}^p r_{0ij}^{\alpha\xi} r_{ij|\mu}^{\beta \eta} g_{0j\eta \xi}^\gamma T_{0j}^\gamma + \sum_{\eta=1}^p T_{0i}^\alpha\left( r_{ij|\mu}^{\beta \eta} \right) T_{0j}^\eta - \sum_{\gamma=1}^p g_{i\alpha\beta}^\gamma r_{ij|\mu}^{\gamma \xi} T_{0j}^\xi\
\end{align*}}}

Let us check $(\ref{t33})_\mu$. In fact, from $(\ref{te7})$ and $(\ref{te43})$, we have
{\small{\begin{align} \label{te48}
&- \sum_{\xi=1}^p G_{i|\mu}^{\alpha \beta \delta \xi} T_{0i}^\xi= -\left[ \Pi_{i|\mu}^{\alpha \beta}, T_{0i}^\delta \right]+\left[ \Pi_{i|\mu}^{\alpha \delta}, T_{0i}^\beta \right]  +\sum_{\gamma=1}^p g_{0i\alpha\delta}^\gamma \Pi_{i|\mu}^{\gamma \beta}  +\left[ T_{0i}^\alpha, \Pi_{i|\mu}^{\beta\delta} \right]  +\sum_{\gamma=1}^p g_{0i\beta\delta}^\gamma \Pi_{i|\mu}^{\alpha \gamma} - \sum_{\gamma=1}^p g_{0i\alpha\beta}^\gamma \Pi_{i|\mu}^{\gamma \delta}  \\
&=\left[ \left[T_{i|\mu}^\alpha, T_{0i}^\beta \right], T_{0i}^\delta \right] + \left[\left[T_{0i}^\alpha, T_{i|\mu}^\beta \right], T_{0i}^\delta \right] - \sum_{\gamma=1}^p \left[g_{0i\alpha\beta}^\gamma T_{i|\mu}^\gamma, T_{0i}^\delta \right] - \sum_{\gamma=1}^p \left[ g_{i\alpha\beta |\mu}^\gamma T_{0i}^\gamma, T_{0i}^\delta   \right] \notag  \\
&-\left[ \left[T_{i|\mu}^\alpha, T_{0i}^\delta \right], T_{0i}^\beta \right] - \left[\left[T_{0i}^\alpha, T_{i|\mu}^\delta \right], T_{0i}^\beta \right] + \sum_{\gamma=1}^p \left[g_{0i\alpha\delta}^\gamma T_{i|\mu}^\gamma, T_{0i}^\beta \right] + \sum_{\gamma=1}^p \left[ g_{i\alpha\delta |\mu}^\gamma T_{0i}^\gamma, T_{0i}^\beta   \right] \notag \\
&+\left[ \left[T_{i|\mu}^\beta, T_{0i}^\delta \right], T_{0i}^\alpha \right] + \left[\left[T_{0i}^\beta, T_{i|\mu}^\delta \right], T_{0i}^\alpha \right] - \sum_{\gamma=1}^p \left[g_{0i\beta\delta}^\gamma T_{i|\mu}^\gamma, T_{0i}^\alpha \right] - \sum_{\gamma=1}^p \left[ g_{i\beta\delta |\mu}^\gamma T_{0i}^\gamma, T_{0i}^\alpha   \right]  \notag \\
& +\sum_{\gamma=1}^p g_{0i\alpha\delta}^\gamma \left(   - \left[T_{i|\mu}^\gamma, T_{0i}^\beta \right] - \left[T_{0i}^\gamma, T_{i|\mu}^\beta\right] + \sum_{\eta=1}^p g_{0i\gamma\beta}^\eta T_{i|\mu}^\eta  + \sum_{\eta=1}^p g_{i\gamma\beta|\mu}^\eta T_{0i}^\eta \right) \notag \\
&+\sum_{\gamma=1}^p g_{0i\beta \delta}^\gamma\left(  - \left[T_{i|\mu}^\alpha, T_{0i}^\gamma \right] - \left[T_{0i}^\alpha, T_{i|\mu}^\gamma \right] + \sum_{\eta=1}^p g_{0i\alpha\gamma}^\eta T_{i|\mu}^\eta  + \sum_{\eta=1}^p g_{i\alpha\gamma|\mu}^\eta T_{0i}^\eta   \right) \notag \\
&-\sum_{\gamma=1}^p g_{0i\alpha\beta}^\gamma \left(  - \left[T_{i|\mu}^\gamma, T_{0i}^\delta \right] - \left[T_{0i}^\gamma, T_{i|\mu}^\delta \right] + \sum_{\eta=1}^p g_{0i\gamma\delta}^\eta T_{i|\mu}^\eta  + \sum_{\eta=1}^p g_{i\gamma\delta |\mu}^\eta T_{0i}^\eta  \right) \notag
\end{align}}}
Then from $(\ref{te25})$ and $(\ref{te48})$, we have $(\ref{t33})_\mu$ (we omit the detail). This completes Lemma \ref{te47}.
\end{proof}

It remains to determine $\varphi_1, T_{i|1}^\alpha, r_{ij|1}^{\alpha\beta}$ and $g_{i\alpha\beta|1}^\gamma$ satisfying $(\ref{t13})_1-(\ref{te6})_1$. Given $\dim_\mathbb{C} \mathbb{H}^1\left(M, \Theta_{\mathcal{F}_0}^\bullet \right)=r$, we can find a basis of $\mathbb{H}^1\left(M, \Theta_{\mathcal{F}_0}^\bullet \right)$ by using the Dolbeault resolution of $\Theta_{\mathcal{F}_0}^\bullet$ from $(\ref{tpc1})$ and represent the basis by
\begin{align}\label{st2}
\left( \overline{\psi_\lambda } , \rho_\lambda \right) \in \frac{A^{0,0}\left( M, \mathscr{H}om_{\mathcal{O}_M}\left( \Theta_{\mathcal{F}_0}, \Theta_M \right) \right)}{A^{0,0}\left( M, \mathscr{H}om_{\mathcal{O}_M}\left(\Theta_{\mathcal{F}_0}, \Theta_{\mathcal{F}_0} \right) \right)} \bigoplus A^{0,1}\left(  M, \Theta_M \right),\,\,\,\,\,\,\,\,\lambda=1,...,r
\end{align}
where $\psi_\lambda \in A^{0,0}\left( M, \mathscr{H}om_{\mathcal{O}_M}\left( \Theta_{\mathcal{F}_0}, \Theta_M \right) \right)$ defined by $\psi_\lambda\left(T_{0i}^\alpha \right)= \psi_{i\lambda}^\alpha$ on $U_i$ for $\psi_{i\lambda}^\alpha \in \Gamma\left( U_i, \mathcal{A}^{0,0}\left( \Theta_M \right) \right)$ with $\psi_{i\lambda}^\alpha= \sum_{\beta=1}^p r_{0ij}^{\alpha\beta} \psi_{j\lambda}^\beta$ on $U_{ij}$ such that there exists $\xi_{i\lambda}\in \Gamma\left( U_{i}, \mathcal{A}^{0,1}\left(\mathscr{H}om_{\mathcal{O}_M}\left(\Theta_{\mathcal{F}_0}, \Theta_{\mathcal{F}_0} \right)\right) \right)$ defined by $\xi_{i\lambda}\left(T_{0i}^\alpha \right)=\xi_{i\lambda}^\alpha= \sum_{\beta=1}^p \xi_{i\lambda}^{\alpha\beta} T_{0i}^\beta \in \Gamma\left( U_{ij}, \mathcal{A}^{0,1}\left(\Theta_{\mathcal{F}_0} \right) \right)$  and $c_{i\lambda} \in \Gamma\left(U_i, \mathscr{H}om_{\mathcal{O}_M}\left( \bigwedge^2 \Theta_{\mathcal{F}_0}, \Theta_M \right) \right)$ defined by $c_{i\lambda} \left(T_{0i}^\alpha \wedge T_{0i}^\beta \right) = c_{i\alpha\beta\lambda}^\gamma T_{0i}^\gamma$ satisfying
\begin{align}
&-\bar{\partial}\rho_\lambda=0 \label{st3} \\
\left[ T_{0i}^\alpha, \psi_{i\lambda}^\beta \right] - \left[ T_{0i}^\beta, \psi_{i\lambda}^\alpha \right]& - \sum_{\beta=1}^p g_{0i\alpha\beta}^\gamma \psi_{i\lambda}^\gamma = \sum_{\beta=1}^p c_{i\alpha\beta\lambda}^\gamma T_{0i}^\gamma \label{tp2}\\
-\bar{\partial} \psi_{i\lambda}^\alpha  + & \left[ \rho_\lambda, T_{0i}^\alpha \right] = \sum_{\beta=1}^p \xi_{i\lambda}^{\alpha\beta} T_{0i}^\beta \label{tp3}
\end{align}
We note that since $\bar{\partial}\xi_{i\lambda}^{\alpha\beta}=0$ from $(\ref{tp3})$, there exists $W_{i\lambda}^{\alpha\eta}\in \Gamma\left( U_i, \mathcal{A}^{0,0} \right)$ such that $\bar{\partial} W_{i\lambda}^{\alpha\beta}= \xi_{i\lambda}^{\alpha\beta}$.
We set
\begin{align}
\varphi_1&= \sum_{\lambda=1}^r t_\lambda \rho_\lambda \label{st3} \\
 T_{i}^{\alpha 1}&= T_{0i}^\alpha + \sum_{\lambda=1}^r t_\lambda\left( \psi_{i\lambda}^\alpha + W_{i\lambda}^\alpha \right) := T_{0i}^\alpha + \sum_{\lambda=1}^r t_\lambda\left(\psi_{i\lambda}^\alpha + \sum_{\beta=1}^p W_{i\lambda}^{\alpha\beta} T_{0i}^\beta \right) \label{st4} \\
 r_{ij}^{\alpha\beta 1}&= r_{0ij}^{\alpha\beta} +  \sum_{\lambda=1}^r t_\lambda  \left( \sum_{\gamma=1}^p W_{i\lambda}^{\alpha\gamma} r_{0ij}^{\gamma \beta} - \sum_{\gamma=1}^p r_{0ij}^{\alpha\gamma} W_{j\lambda}^{\gamma  \beta} \right)\label{st5}
\end{align}
Then $(\ref{t13})_1,(\ref{tpc6})_1$, $(\ref{tpc4})_1$ and $(\ref{tpc8})_1$ are satisfied. $(\ref{tpc6})_1$ and $(\ref{tpc4})_1$ implies $(\ref{tpc5})_1$. Let us find $g_{i\alpha\beta}^{\gamma 1}$ satisfying $(\ref{tpc2})_1$. We note that
\begin{align*}
&\left[ T_{0i}^\alpha, T_{0i}^\beta + t_\lambda\psi_{i\lambda}^\beta + t_\lambda W_i^\beta \right] - \left[ T_{0i}^\beta, T_{0i}^\alpha + t_\lambda \psi_{i\lambda}^\alpha + t_\lambda W_i^\alpha \right] - \sum_{\beta=1}^p \left( g_{0i\alpha\beta}^\gamma +  g_{i\alpha\beta \lambda |1 }^{\gamma } \right) \left(T_{0i}^\gamma + t_\lambda \psi_{i\lambda}^\gamma + t_\lambda W_{i\lambda}^\gamma \right)\\
&\equiv_1 t_\lambda \left(\sum_{\beta=1}^p c_{i\alpha\beta\lambda}^\gamma T_{0i}^\gamma + \left[ T_{0i}^\alpha, W_{i\lambda}^\beta \right]  -\left[ T_{0i}^\beta,  W_{i\lambda}^\alpha \right] - \sum_{\gamma=1}^p g_{0i\alpha\beta}^\gamma W_{i\lambda}^\gamma   \right) -   \sum_{\gamma=1}^p  g_{i\alpha\beta \lambda | 1}^\gamma T_{0i}^\gamma \equiv_1 \sum_{\gamma=1}^p \left( t_\lambda C_{i\alpha\beta\lambda}^\gamma- g_{i\alpha \beta \lambda |1}^\gamma \right) T_{0i}^\gamma
\end{align*}
where we write $\sum_{\gamma=1}^p C_{i\alpha\beta \lambda}^\gamma T_{0i}^\gamma:= \sum_{\beta=1}^p c_{i\alpha\beta\lambda}^\gamma T_{0i}^\gamma + \left[ T_{0i}^\alpha, W_{i\lambda}^\beta \right]  -\left[ T_{0i}^\beta,  W_{i\lambda}^\alpha \right] - \sum_{\gamma=1}^p g_{0i\alpha\beta}^\gamma W_{i\lambda}^\gamma $ for $C_{i\alpha\beta\lambda}^\gamma \in \Gamma\left( U_i, \mathcal{A}^{0,0}\right)$, so that we take 
\begin{align}
g_{i\alpha\beta}^{\gamma 1} : = g_{0i\alpha\beta}^\gamma + \sum_{\lambda=1}^r t_\lambda C_{i\alpha\beta\lambda}^\gamma \label{st6}
\end{align}
Then $(\ref{tpc2})_1$ is satisfied. $(\ref{tpc6})_1, (\ref{tpc4})_1$ and $(\ref{tpc2})_1$ implies $(\ref{tpc7})_1, (\ref{tpc3})_1$ and $(\ref{te6})_1$. This completes the inductive construction of $\varphi, T_i^\alpha, r_{ij}^{\alpha\beta}$ and $g_{i\alpha\beta}^\gamma$ satisfying $(\ref{t10})-(\ref{te5})$.

\subsection{Proof of convergence}\

We will prove that $\varphi =\sum_{\mu=1}^\infty \varphi_\mu, T_i^\alpha =\sum_{\mu=0}^\infty T_{i|\mu}^\alpha, r_{ij}^{\alpha\beta}=\sum_{\mu=0}^\infty r_{ij|\mu}^{\alpha\beta}$ and $g_{i\alpha\beta}^\gamma=\sum_{\mu=0}^\infty g_{i\alpha\beta |\mu}^\gamma$ in the preceding subsection converge. Our proof of convergence is based on a combination of the ideas and the methods from \cite{Hor76} and \cite{Hor73} and \cite{Kod05}. Before proceeding to the discussion, we would like to remark that the quotient in the Dolbeault resolution of $\Theta_{\mathcal{F}_0}^\bullet$ makes things complicated. So we will introduce another complex of sheaves which also controls foliated deformations of $\left( M, \Theta_{\mathcal{F}_0} \right)$ but removes the quotient in the degree $1$. We shall define the following complex of sheaves
\begin{align*}
\mathcal{E}_{\Theta_{\mathcal{F}_0}}^\bullet : \mathcal{E}_{\Theta_{\mathcal{F}_0}}\xrightarrow{D_0'} \mathscr{H}om_{\mathcal{O}_M}\left( \Theta_{\mathcal{F}_0}, \Theta_M \right) \xrightarrow{D_1'} \mathscr{H}om_{\mathcal{O}_M}\left(  \bigwedge^2 \Theta_{\mathcal{F}_0}, \frac{\Theta_M}{\Theta_{\mathcal{F}_0}} \right) \xrightarrow{D_2}  \mathscr{H}om_{\mathcal{O}_M}\left( \bigwedge^3 \Theta_{\mathcal{F}_0}, \frac{\Theta_M}{\Theta_{\mathcal{F}_0}} \right) \xrightarrow{D_3} \cdots
\end{align*}
where $\mathcal{E}_{\Theta_{\mathcal{F}_0}}$ is the Atiyah extension defined in the following way: for an open covering of $M$ as in the assumption (\ref{d1}) in the beginning of the proof of Theorem \ref{tt2}, we define the locally free sheaf $\mathcal{E}_{\Theta_{\mathcal{F}_0}}$ locally on $U_i$ by $\mathcal{E}_{\Theta_{\mathcal{F}_0}}|_{U_i}\cong \mathscr{H}om_{\mathcal{O}_M}\left( \Theta_{\mathcal{F}_0} , \Theta_{\mathcal{F}_0} \right)|_{U_i} \bigoplus \Theta_M|_{U_i}$ such that $(\phi_i,X)\in \mathscr{H}om_{\mathcal{O}_M}\left( \Theta_{\mathcal{F}_0}, \Theta_{\mathcal{F}_0} \right)|_{U_i} \bigoplus \Theta_M |_{U_i}$ and $(\phi_j, X)\in \mathscr{H}om_{\mathcal{O}_M}\left( \Theta_{\mathcal{F}_0}, \Theta_{\mathcal{F}_0} \right)|_{U_j} \bigoplus \Theta_M |_{U_j}$ are identified on $U_{ij}$ if $
\phi_i\left(T_{0i}^\alpha\right)=\phi_j\left(T_{0i}^\alpha\right)+\sum_{\beta=1}^p \left[X, r_{0ij}^{\alpha\beta}\right]T_{0j}^\beta$. Then we define $D_0':\mathcal{E}_{\Theta_{\mathcal{F}_0}}\to \mathscr{H}om_{\mathcal{O}_X}\left( \Theta_{\mathcal{F}_0}, \Theta_M \right)$ locally on $U_i$ by $D_0' \left(\left(\phi_i,X\right)\right)\left(T_{0i}^\alpha\right)= \left(\mathcal{L}_X-\phi_i\right)\left(T_{0i}^\alpha\right)$ and then linearly extends to $\Gamma\left( U_i , \Theta_{\mathcal{F}_0} \right)$. We define $D_1'$ by the composition of the natural quotient map $\mathscr{H}om_{\mathcal{O}_M}\left( \Theta_{\mathcal{F}_0}, \Theta_M \right)\to \mathscr{H}om_{\mathcal{O}_M}\left( \Theta_{\mathcal{F}_0}, \frac{\Theta_M }{\Theta_{\mathcal{F}_0}} \right)$ with $D_1$. Then we have 
\begin{align} \label{ttc1}
\mathbb{H}^i\left(  M, \mathcal{E}_{\Theta_{\mathcal{F}_0}}^\bullet \right)\cong \mathbb{H}^i\left( M, \Theta_{\mathcal{F}_0}^\bullet \right),\,\,\,\,\,\,\,\,\,i\geq 0
\end{align}
 (for the detail on the complex of sheaves $\mathcal{E}_{\Theta_{\mathcal{F}_0}}^\bullet$, see Part II). We denote $\mathcal{A}^{0,p}\left( \mathcal{E}_{\Theta_{\mathcal{F}_0}} \right)$ be the sheaf of germs of $C^\infty(0,p)$-forms with coefficients in $\mathcal{E}_{\Theta_{\mathcal{F}_0}}^\bullet$ and denote by $A^{0,p}\left( M, \mathcal{E}_{\Theta_{\mathcal{F}_0}} \right)$ the global section of $\mathcal{A}^{0,p}\left( \mathcal{E}_{\Theta_{\mathcal{F}_0}} \right)$. Then we have the following  Dolbeault resolution of $\mathcal{E}_{\Theta_{\mathcal{F}_0}}^\bullet$:
\begin{equation}
\begin{CD}
\cdots \\
@A\hat{D}_3AA \\
\frac{A^{0,0}\left( M,  \bigwedge^3 \Theta_{\mathcal{F}_0}^*\otimes \Theta_M\right)}{A^{0,0}\left(M, \bigwedge^3 \Theta_{\mathcal{F}_0}^*\otimes \Theta_{\mathcal{F}_0} \right)}@>-\bar{\partial}>> \cdots\\
@A\hat{D}_2AA @A\hat{D}_2AA \\
\frac{A^{0,0}\left(M,  \bigwedge^2 \Theta_{\mathcal{F}_0}^*\otimes \Theta_M\right)}{A^{0,0}\left(M, \bigwedge^2 \Theta_{\mathcal{F}_0}^* \otimes \Theta_{\mathcal{F}_0} \right)} @>\bar{\partial}>> \frac{A^{0,1}\left( M, \bigwedge^2 \Theta_{\mathcal{F}_0}^*\otimes \Theta_M \right)}{A^{0,1}\left(M, \bigwedge^2 \Theta_{\mathcal{F}_0}^*\otimes \Theta_{\mathcal{F}_0} \right)}@>-\bar{\partial}>> \cdots \\
@A\hat{D}_1'AA @A\hat{D}_1'AA @A\hat{D}_1'AA\\
A^{0,0}\left(M, \Theta_{\mathcal{F}_0}^*\otimes \Theta_M\right) @>-\bar{\partial}>> A^{0,1}\left(M, \Theta_{\mathcal{F}_0}^*\otimes \Theta_M\right) @>\bar{\partial}>> A^{0,2}\left(M, \Theta_{\mathcal{F}_0}^*\otimes \Theta_M\right) @>-\bar{\partial}>>\cdots \\
@A\hat{D}_0'AA @A\hat{D}_0'AA @A\hat{D}_0'AA @A\hat{D}_0'AA\\
A^{0,0}\left(M, \mathcal{E}_{\Theta_{\mathcal{F}_0}} \right) @>\bar{\partial} >> A^{0,1}\left(M, \mathcal{E}_{\Theta_{\mathcal{F}_0}}  \right) @>-\bar{\partial}>> A^{0,2}\left(M, \mathcal{E}_{\Theta_{\mathcal{F}_0}} \right)  @>\bar{\partial}>> A^{0,3}(M, \Theta_M) @>-\bar{\partial}>> \cdots\\
\end{CD}
\end{equation}
We reinterpret $(\ref{d2})$ in terms of the above Dolbeault resolution of $\Theta_{\mathcal{F}_0}^\bullet$. From $(\ref{te21})$ we have 
\begin{align}\label{tpc19}
\sum_{\beta,\eta=1}^p r_{0ij}^{\alpha\beta}\bar{\partial}\Lambda_{j|\mu}^{\beta\eta} r_{0ji}^{\eta\xi} - \bar{\partial}\Lambda_{i|\mu}^{\alpha\xi}&=\bar{\partial}\Lambda_{ij|\mu}^{\alpha\xi}
\end{align}
We define $\bar{\partial}\Lambda_{i|\mu}\in \Gamma\left( U_i, \mathcal{A}^{0,2}\left( \mathscr{H}om_{\mathcal{O}_M}\left( \Theta_{\mathcal{F}_0}, \Theta_{\mathcal{F}_0} \right) \right) \right)$ by
\begin{align*}
\bar{\partial}\Lambda_{i|\mu}: \Gamma\left( U_i, \Theta_{\mathcal{F}_0} \right) &\to \Gamma\left( U_i, \mathcal{A}^{0,2}\left(\Theta_{\mathcal{F}_0} \right) \right) \\
 T_{0i}^\alpha &\mapsto  \sum_{\xi=1}^p \bar{\partial} \Lambda_{i|\mu}^{\alpha \xi} T_{0i}^\xi
\end{align*}
Then from $(\ref{t23})$ and $(\ref{tpc19})$, we have
\begin{align*}
\bar{\partial}\Lambda_{j|\mu}\left(T_{0i}^\alpha \right) - \bar{\partial}\Lambda_{i|\mu}\left(T_{0i}^\alpha \right) = \sum_{\xi=1}^p \left[ \xi_\mu, r_{0ij}^{\alpha \xi} \right] T_{0j}^\xi
\end{align*}
This implies that
\begin{align}\label{d28}
\left( -\xi_\mu, \left\{\bar{\partial} \Lambda_{i|\mu}\right\} \right) \in A^{0,2}\left( M, \mathcal{E}_{\Theta_{\mathcal{F}_0}} \right)
\end{align}
Then we see that from $(\ref{d2})$ and $(\ref{d3})$
\begin{align}\label{d13}
\left( \overline{B_\mu}, \tilde{\Phi}_\mu, \left(- \xi_\mu, \left\{ \bar{\partial}\Lambda_{i|\mu} \right\} \right) \right) \in \frac{A^{0,0}\left( M, \mathscr{H}om_{\mathcal{O}_M}\left( \bigwedge^2 \Theta_{\mathcal{F}_0}, \Theta_M \right) \right)}{A^{0,0}\left( M , \mathscr{H}om_{\mathcal{O}_M}\left( \bigwedge^2 \Theta_{\mathcal{F}_0}, \Theta_{\mathcal{F}_0} \right) \right)} \bigoplus A^{0,1}\left( M, \mathscr{H}om_{\mathcal{O}_M}\left(\Theta_{\mathcal{F}_0}, \Theta_M \right)\right) \bigoplus A^{0,2}\left( M , \mathcal{E}_{\Theta_{\mathcal{F}_0}} \right)
\end{align}
defines a $2$-cocycle in the above Dolbeault resolution of $\mathcal{E}_{\Theta_{\mathcal{F}_0}}^\bullet$.

We recall the \"Holder norm on sections of a locally free sheaf (see \cite{Kod05} p.274 and \cite{Hor73} p.388). We will define three \"Holder norms on the sections of $\mathcal{A}^{0,q}\left(\Theta_M \right)$ and $\mathcal{A}^{0,q}\left(\mathscr{H}om_{\mathcal{O}_M}\left(\Theta_{\mathcal{F}_0}, \Theta_M \right) \right)$ and $\mathcal{A}^{0,q}\left( M, \mathcal{E}_{\Theta_{\mathcal{F}_0}} \right)$ and apply to the harmonic theory.

We define the \" Holder norm $|-|_{k+\alpha}$ (k: an integer $\geq 2$, $0<\alpha <1$) for sections of $\mathcal{A}^{0,q}\left( \Theta_M \right)$ as follows: let $\varphi\in \Gamma\left( U_i, \mathcal{A}^{0,q}\left( \Theta_M \right) \right)$ and we write
\begin{align}\label{d6}
\varphi =\sum_{\lambda=1}^n \varphi_i^\lambda \frac{\partial}{\partial z_i^\lambda},\,\,\,\,\,\,\,\,\varphi_i^\lambda= \frac{1}{q!} \sum \varphi_{i\mu_1\cdots \mu_q}^\lambda (z_i)d\bar{z}_i^{\mu_1}\wedge \cdots \wedge d\bar{z}_i^{\mu_q}\in \Gamma\left(U_i, \mathcal{A}^{0,q}\right)
\end{align}
in terms of local coordinates $\left( z_i^1,..., z_i^n \right)$ and let
\begin{align}\label{d7}
\left| \varphi \right|_{k+\alpha}^{U_i} = \sum_{h=0}^k \sup \left| D_i^h \varphi_{i \mu_1\cdots \mu_q }^\lambda(z_i) \right| + \sup \frac{ \left| D_i^k \varphi_{i \mu_1\cdots \mu_q}^\lambda(z_i)- D_i^k \varphi_{\mu_1\cdots \mu_q}^\lambda (y_i) \right| }{\left|z_i-y_i\right|^\alpha}
\end{align}
where the $``\sup$" is extended over all points $z,y\in U_i$, all indices $\lambda, \beta, \gamma, \mu_1,...,\mu_q$, and all partial derivatives $D_i^h, D_i^k$ of order $h, k$ with respect to $z_i^1,..., z_i^n, \bar{z}_i^1,..., \bar{z}_i^n$. For $\varphi \in A^{0,q}\left( M, \Theta_M \right)$, we define
\begin{align}\label{d8}
\left| \varphi \right|_{k+\alpha}=\max_i \left| \varphi \right|_{k+\alpha}^{U_i}
\end{align}
We introduce a Hermitian metric on the fibers of $\Theta_M$ as the vector bundle and define an inner product on $A^{0,q}\left( M, \Theta_M \right)$ and apply a harmonic theory on $\Theta_M$ as in \cite{Kod05} p. 157-p.161.
\begin{align}\label{d9}
&\text{We denote by $\mathfrak{d}$ the adjoint operator of $\bar{\partial}$, and let $\square=\mathfrak{d}\bar{\partial} + \bar{\partial} \mathfrak{d}$ and $G$ the Green's operator}\\
&\text{and $\bold{H}$ the projection onto the space of harmonic forms.} \notag
\end{align}

We define the H\"older norm $|-|_{k+\alpha}$ (k: an integer $\geq 2$, $0<\alpha< 1$) for sections of $\mathcal{A}^{0,q}\left(\mathscr{H}om_{\mathcal{O}_M}\left(\Theta_{\mathcal{F}_0}, \Theta_M \right) \right)$ as follows: let $\phi\in \Gamma\left(U_i, \mathcal{A}^{0,p}\left(\mathscr{H}om_{\mathcal{O}_M}\left(\Theta_{\mathcal{F}_0}, \Theta_M \right) \right) \right)$ and we write 
\begin{align*}
\phi\left(T_{0i}^\beta\right)=\sum_{\gamma=1}^n A_{i}^{\beta\gamma} \frac{\partial}{\partial z_i^\gamma},\,\,\,\,\,\,\,\,\,\,\,A_i^{\beta\gamma}=\frac{1}{q!}\sum A_{i\mu_1\cdots \mu_q}^{\beta\gamma} (z_i)dz_i^{\mu_1}\wedge \cdots \wedge dz_i^{\mu_q} \in \Gamma\left(U_i, \mathcal{A}^{0,q}\right)
\end{align*}
in terms of local coordinates $(z_i^1,..., z_i^n)$ and let
\begin{align}\label{d39}
\left| \phi \right|_{k+\alpha}^{U_i} = \sum_{h=0}^k \sup \left| D_i^h A_{i \mu_1\cdots \mu_q }^{\beta\gamma}(z_i) \right| + \sup \frac{ \left| D_i^kA_{i|\mu_1\cdots \mu_q}^{\beta\gamma}(z_i)- D_i^k A_{i\mu_1\cdots \mu_q}^{\beta\gamma}(y_i) \right| }{\left|z_i-y_i\right|^\alpha}
\end{align}
where the $``\sup$" is extended over all points $z,y\in U_i$, all indices $\beta, \gamma, \mu_1,...,\mu_q$, and all partial derivatives $D_i^h, D_i^k$ of order $h, k$ with respect to $z_i^1,..., z_i^n, \bar{z}_i^1,..., \bar{z}_i^n$. For $\phi\in A^{0,p}\left( M,  \mathscr{H}om_{\mathcal{O}_M}\left(\Theta_{\mathcal{F}_0} , \Theta_M \right)     \right)$, we define
\begin{align*}
\left| \phi \right|_{k+\alpha} = \max_i \left| \phi \right|_{k+\alpha}^{U_i} 
\end{align*}
We introduce a harmonic theory on $\mathscr{H}om_{\mathcal{O}_M}\left( \Theta_{\mathcal{F}_0}, \Theta_M \right)$ as in \cite{Kod05} p.157-p.161.
\begin{align}
\textnormal{We denote by $\mathfrak{d}'$ the adjoint operator of $\bar{\partial}$, and $\square
'=\mathfrak{d}'\bar{\partial} + \bar{\partial} \mathfrak{d}'$ and $G'$ the Green's operator.}
\end{align}

 We define the H\"older norm $|-|_{k+\alpha}$ (k: an integer $\geq 2$, $0<\alpha< 1$) for sections of $\mathcal{A}^{0,q}\left(\mathscr{H}om_{\mathcal{O}_M}\left(\bigwedge^2 \Theta_{\mathcal{F}_0} , \Theta_M \right) \right)$ as follows: let $\phi\in \Gamma\left(U_i, \mathcal{A}^{0,q}\left(\mathscr{H}om_{\mathcal{O}_M}\left( \bigwedge^2 \Theta_{\mathcal{F}_0} , \Theta_M \right) \right) \right)$ and we write 
\begin{align*}
\phi\left(T_i^\beta\wedge T_i^{\gamma}\right)=\sum_{\eta=1}^n A_{i}^{\beta\gamma \eta} \frac{\partial}{\partial z_i^\eta},\,\,\,\,\,\,\,\,\,\,\,A_i^{\beta\gamma\eta}=\frac{1}{q!}\sum A_{i\mu_1\cdots \mu_q}^{\beta\gamma} (z_i)dz_i^{\mu_1}\wedge \cdots \wedge dz_i^{\mu_q} \in \Gamma\left(U_i, \mathcal{A}^{0,q}\right)
\end{align*}
in terms of local coordinates $(z_i^1,..., z_i^n)$ and let
\begin{align*}
\left| \phi \right|_{k+\alpha}^{U_i} = \sum_{h=0}^k \sup \left| D_i^h A_{i \mu_1\cdots \mu_q }^{\beta\gamma \eta}(z_i) \right| + \sup \frac{ \left| D_i^kA_{ i \mu_1\cdots \mu_q}^{\beta\gamma \eta}(z_i)- D_i^k A_{i\mu_1\cdots \mu_q}^{\beta\gamma \eta}(y_i) \right| }{\left|z_i-y_i\right|^\alpha}
\end{align*}
where the $``\sup$" is extended over all points $z,y\in U_i$, all indices $\beta, \gamma, \eta, \mu_1,...,\mu_q$, and all partial derivatives $D_i^h, D_i^k$ of order $h, k$ with respect to $z_i^1,..., z_i^n, \bar{z}_i^1,..., \bar{z}_i^n$. For $\phi\in A^{0,p}\left( M, \mathscr{H}om_{\mathcal{O}_M}\left( \bigwedge^2 \Theta_{\mathcal{F}_0} , \Theta_M  \right)     \right)$, we define
\begin{align*}
\left| \phi \right|_{k+\alpha} = \max_i \left| \phi \right|_{k+\alpha}^{U_i} 
\end{align*}

Let $\mathcal{E}_{\Theta_{\mathcal{F}_0}}$ be the Atiyah extension of $\Theta_{\mathcal{F}_0}$ which is locally of the form 
\begin{align*}
\Gamma\left(U_i, \mathcal{E}_{\Theta_{\mathcal{F}_0}}\right):= \Gamma(U_i, \Theta_M)\bigoplus \Gamma\left(U_i, \mathscr{H}om_{\mathcal{O}_M}\left(\Theta_{\mathcal{F}_0} , \Theta_{\mathcal{F}_0} \right)\right)
\end{align*}
We define the H\"older norm $|-|_{k+\alpha}$ (k: an integer $\geq 2, 0<\alpha <1$) for sections of $\mathcal{A}^{0,q}\left(\mathcal{E}_{\Theta_{\mathcal{F}_0}}\right)$ as follows: let $\xi\in \Gamma(U_i, \mathcal{A}^{0,q}(\mathcal{E}_{\Theta_\mathcal{F}}))$ and we write $\xi=(\varphi_i, \phi_i)$ as
\begin{align}\label{d10}
\varphi_i &=\sum_{\lambda=1}^n \varphi_i^\lambda \frac{\partial}{\partial z_i^\lambda},\,\,\,\,\,\,\,\,\varphi_i^\lambda= \frac{1}{q!} \sum \varphi_{i\mu_1\cdots \mu_q}^\lambda(z_i)d\bar{z}_i^{\mu_1}\wedge \cdots \wedge d\bar{z}_i^{\mu_q}\in \Gamma\left(U_i, \mathcal{A}^{0,q}\right)\\
\phi_i\left(T_{0i}^\beta\right)&=\sum_{\gamma=1}^p A_{i}^{\beta\gamma} T_{0i}^\gamma \,\,\,\,\,\,\,\,\,\,\,A_i^{\beta\gamma}=\frac{1}{q!}\sum A_{i\mu_1\cdots \mu_q}^{\beta\gamma} (z_i)dz_i^{\mu_1}\wedge \cdots \wedge dz_i^{\mu_q} \in \Gamma\left(U_i, \mathcal{A}^{0,q}\right) \notag
\end{align} 
in terms of local coordinates $\left( z_i^1,..., z_i^n \right)$ and let
\begin{align}\label{d11}
\left| \xi \right|_{k+\alpha}^{U_i} = \sum_{h=0}^k \sup \left| D_i^h C_{i \mu_1\cdots \mu_q }(z_i) \right| + \sup \frac{ \left| D_i^k C_{i \mu_1\cdots \mu_q}(z_i)- D_i^k C_{\mu_1\cdots \mu_q}(y_i) \right| }{\left|z_i-y_i\right|^\alpha}
\end{align}
where $C_{i\mu_1\cdots \mu_q}= \varphi_{i \mu_1\cdots \mu_q}^\lambda (z_i) $ or $A_{i\mu_1,...\mu_q}^{\beta \gamma}(z_i)$, and the $``\sup$" is extended over all points $z,y\in U_i$, all indices $\lambda, \beta, \gamma, \mu_1,...,\mu_q$, and all partial derivatives $D_i^h, D_i^k$ of order $h, k$ with respect to $z_i^1,..., z_i^n, \bar{z}_i^1,..., \bar{z}_i^n$. For $\xi \in A^{0,q}\left( M, \mathcal{E}_{\Theta_{\mathcal{F}_0}} \right)$, we define
\begin{align}\label{d12}
\left| \xi \right|_{k+\alpha}=\max_i \left|\xi \right|_{k+\alpha}^{U_i}
\end{align}
In particular, we note that for $\xi=\left( \varphi, \left\{ \phi_i \right\} \right)\in A^{0,q}\left( M, \mathcal{E}_{\Theta_{\mathcal{F}_0}} \right)$, where $\varphi\in A^{0,q}\left(M, \Theta_M \right) $ and $\varphi_i\in \Gamma\left(U_i, \mathcal{A}^{0,q}\left(\mathscr{H}om_{\mathcal{O}_M}\left( \Theta_{\mathcal{F}_0} , \Theta_{\mathcal{F}_0} \right) \right) \right)$, we have from $(\ref{d8})$
\begin{align} \label{ds5}
\left|\xi \right|_{k+\alpha}\ll \left| \varphi \right|_{k+\alpha} + \max_i \left| \phi_i \right|_{k+\alpha}^{U_i}
\end{align}
where $\left| \phi_i \right|_{k+\alpha}^{U_i}$ is defined by the \"Holder norm on sections of $\mathcal{A}^{0,q}\left(\mathscr{H}om_{\mathcal{O}_M}\left(\Theta_{\mathcal{F}_0}, \Theta_{\mathcal{F}_0} \right)\right)$ in terms of $\phi_i\left(T_{0i}^\beta \right)=\sum_{\gamma=1}^p A_i^{\beta \gamma} T_{0i}^\gamma$ as in $(\ref{d10})$ instead of $(\ref{d39})$.

We introduce a harmonic theory on $\mathcal{E}_{\Theta_{\mathcal{F}_0}}$ as in \cite{Kod05} p.157-p.161. 
\begin{align}\label{d14}
\textnormal{We denote by $\tilde{\mathfrak{d}}$ the adjoint operator of $\bar{\partial}$, and $\widetilde{\square}=\tilde{\mathfrak{d}}\bar{\partial} + \bar{\partial} \tilde{\mathfrak{d}}$ and $\tilde{G}$ the Green's operator.}
\end{align}

Consider a formal power series $\phi= \phi(t)=\phi(t_1,..., t_m)=\sum \phi_{v_1\cdots v_m} t_1^{v_1}\cdots t_m^{v_m}$ with coefficients in $A^{0,q}\left( M, \mathcal{G} \right)$ or in $\Gamma\left( U_i, \mathcal{A}^{0,q}\left(\mathcal{G} \right)  \right)$ or in $\Gamma\left(U_{ij}, \mathcal{A}^{0,q}\left( \mathcal{G} \right) \right)$ where $\mathcal{G}=\Theta_M$ or $\Omega_M^1$ or $\mathcal{O}_M$, and a power series $a(t)=\sum a_{v_1\cdots v_m} t_1^{v_1}\cdots t_m^{v_m}, a_{v_1\cdots v_m} \geq 0$. 
\begin{align}\label{d30}
\textnormal{We indicate by $\left|\phi \right|_{k+\alpha} \ll a(t)$ or $\left| \phi \right|_{k+\alpha}^{U_i} \ll a(t)$ that $\left| \phi_{v_1\cdots v_m} \right|_{k+\alpha} \leq a_{v_1\cdots v_m}$   }
\end{align}
In the following we do not indicate explicitly the domain $U_i$ or $U_{ij}$ if no confusion is possible.

We set 
\begin{align}\label{d31}
A(t)=\frac{b}{16c}\sum_{\mu=1}^\infty \frac{1}{\mu^2} c^\mu\left( t_1+ \cdots + t_r \right)^\mu\,\,\,\,\,\,\,\,\,\textnormal{with}\,\,\,b>0\,\,\,\textnormal{and}\,\,\, c>0
\end{align}
Then we have $A(t)^2\ll \frac{b}{c} A(t)$ (see \cite{Kod05} p.279). With this preparation, we will show that for a fixed integer $k\geq 2$ and $0< \alpha<1$, the inductive construction of $\varphi, T_i^\alpha , r_{ij}^{\alpha\beta}$ and $g_{i\alpha\beta}^\gamma$ in the previous subsection can be carried out in such a way that
\begin{align*}
\left|\varphi \right|_{k+\alpha} \ll A(t) \\
\left| T_i^{\alpha }- T_{0i}^{\alpha} \right|_{k+ \alpha}  \ll  A(t) \\
\left| r_{ij}^{\alpha\beta }- r_{0ij}^{\alpha\beta} \right|_{k+ 1+ \alpha} \ll A(t) \\
\left| g_{i\alpha\beta}^{\gamma} - g_{0i\alpha\beta}^\gamma \right|_{ k -1 +  \alpha} \ll A(t)
\end{align*}
Then it suffices to prove that for $\mu=1,2,3,\cdots$,
\begin{align}
\left| \varphi^{\mu}\right|_{k+\alpha}  \ll A(t) \label{d20} \\
\left| T_i^{\alpha\mu}- T_{0i}^{\alpha} \right|_{k+ \alpha}  \ll  A(t) \\
\left| r_{ij}^{\alpha\beta \mu}- r_{0ij}^{\alpha\beta} \right|_{k+ 1+ \alpha} \ll A(t) \\
\left| g_{i\alpha\beta}^{\gamma \mu} - g_{0i\alpha\beta}^\gamma \right|_{ k -1+  \alpha} \ll A(t) \label{d21}
\end{align}
for some proper choice of constants $c>b>0$. We prove $(\ref{d20})_\mu -(\ref{d21})_\mu$ by induction on $\mu$. For $\mu=1$ we have $(\ref{st3})$ and $(\ref{st4})$ and $(\ref{st5})$ and $(\ref{st6})$, and the linear term of $A(t)$ is $\frac{b}{16}\left(t_1+ \cdots  + t_r \right)$. Therefore $(\ref{d20})_1-(\ref{d21})_1$ holds if $b$ is sufficiently large.

Now assume that $(\ref{d20})_{\mu-1}-(\ref{d21})_{\mu-1}$ are satisfied. We will derive $(\ref{d20})_\mu- (\ref{d21})_\mu$. In the following $K_1,K_2, K_3, \cdots$ will denote constants which depend only on $k,\alpha, M, \Theta_{\mathcal{F}_0}$.

From $(\ref{d13})$, we have $\bar{\partial}\left( -\xi_\mu, \left\{\bar{\partial} \Lambda_{i|\mu} \right\} \right)=0$. Then we take $\left(  \tilde{\varphi}_\mu, \left\{ \tilde{b}_{i|\mu} \right\} \right):= \tilde{\mathfrak{d}}\tilde{G}\left(-\xi_\mu, \left\{ \bar{\partial} \Lambda_{i|\mu} \right\} \right)\in A^{0,1}\left( M,  \mathcal{E}_{\Theta_{\mathcal{F}_0}} \right)$ from $(\ref{d14})$, where $\tilde{\varphi}_\mu\in A^{0,1}\left(M, \Theta_M \right)$ and $\tilde{b}_{i|\mu}\in \Gamma\left( U_i, \mathcal{A}^{0,1}\left( \mathscr{H}om_{\mathcal{O}_M}\left( \Theta_{\mathcal{F}_0}, \Theta_{\mathcal{F}_0} \right) \right) \right)$ defined by $\tilde{b}_{i|\mu}:\Gamma\left( U_i, \Theta_{\mathcal{F}_0} \right)\to \Gamma\left( U_i, \mathcal{A}^{0,1}\left( \Theta_{\mathcal{F}_0} \right)\right), T_{0i}^\alpha \mapsto \sum_{\xi=1}^p \tilde{b}_{i|\mu}^{\alpha \xi} T_{0i}^\xi$.  We have
\begin{align}\label{d25}
\left| \left( \tilde{\varphi}_\mu, \left\{ \tilde{b}_{i|\mu}\right\} \right) \right|_{k+\alpha} \leq K_1\left|\left( - \xi_\mu, \left\{ \bar{\partial}\Lambda_{i|\mu}\right\} \right) \right|_{k-1+\alpha}
\end{align}
where $K_1$ is a constant independent of $\left( - \xi_\mu, \left\{ \bar{\partial} \Lambda_{i|\mu} \right\} \right)$. We can prove the inequality $(\ref{d25})$ by the definition of \"Holder norm (see \cite{Kod05} p.270 (5.118)) and the similar type of Lemma to \cite{Kod05} Lemma 5.7 p.276 which is based on the $a\,\,\, priori\,\,\, estimate$ p.275. Since $\tilde{\square}$ from $(\ref{d14})$ is strongly elliptic, we have the $a\,\,\, priori \,\,\, estimate$ (see \cite{Kod05} Theorem 4.3 p.436), so that we can prove the similar type of Lemma. 

From $(\ref{d9})$, we have
\begin{align*}
\tilde{\varphi}_\mu = \bold{H} \tilde{\varphi}_\mu + \left( \bar{\partial} \mathfrak{d} + \mathfrak{d} \bar{\partial} \right) G \tilde{\varphi}_\mu
\end{align*}
It follows that
\begin{align*}
\bar{\partial}\left( \bold{H} \tilde{\varphi}_\mu  +\mathfrak{d}\bar{\partial} G \tilde{\varphi}_\mu \right) &= - \xi_\mu\\
\mathfrak{d}\left( \bold{H} \tilde{\varphi}_\mu + \mathfrak{d} \bar{\partial} G \tilde{\varphi}_\mu \right)&=0
\end{align*}
On the other hand, we define $Q_{ij|\mu} \in \Gamma\left( U_{ij}, \mathcal{A}^{0,0}\left( \mathscr{H}om_{\mathcal{O}_M}\left( \Theta_{\mathcal{F}_0}, \Theta_{\mathcal{F}_0}\right) \right) \right)$ by 
\begin{align*}
Q_{ij|\mu}:\Gamma \left(U_{ij}, \Theta_{\mathcal{F}_0 } \right)&\to \Gamma \left(U_{ij}, \mathcal{A}^{0,0}\left(\Theta_{\mathcal{F}_0} \right) \right)\\
T_{0i}^\alpha &\mapsto  \sum_{\beta=1}^p \left[ \mathfrak{d}G\tilde{\varphi}_\mu, r_{0ij}^{\alpha\beta}\right] T_{0j}^\beta,
\end{align*}
and linearly extends to $\Gamma\left( U_{ij}, \Theta_{\mathcal{F}_0} \right)$. Then $\left\{ Q_{ij|\mu}\right\}\in C^1\left( \mathcal{U}, \mathcal{A}^{0,0}\left( \mathscr{H}om_{\mathcal{O}_M}\left( \Theta_{\mathcal{F}_0}, \Theta_{\mathcal{F}_0} \right) \right) \right)$ defines a $1$-cocycle. Let $\left\{ \rho_i(z) \right\}$ be a partition of unity subordinate to the covering $\mathcal{U}=\left\{ U_i \right\}$. By setting $Q_{i|\mu}=\sum_k \rho_k(z)Q_{ik|\mu}$, we see that  $\left\{Q_{i|\mu} \right\}\in C^0\left( \mathcal{U}, \mathcal{A}^{0,0}\left( \mathscr{H}om_{\mathcal{O}_M}\left(  \Theta_{\mathcal{F}_0}, \Theta_{\mathcal{F}_0}  \right)   \right)  \right)$ such that $Q_{i|\mu} -Q_{j|\mu}=Q_{ij|\mu}$ where
\begin{align*}
Q_{i|\mu}:\Gamma\left(U_i, \Theta_{\mathcal{F}_0} \right) &\to \Gamma\left(U_i, \mathcal{A}^{0,0}\left(\Theta_{\mathcal{F}_0} \right)\right)\\
    T_{0i}^\alpha &\mapsto \sum_{\gamma=1}^p Q_{i|\mu}^{\alpha\gamma} T_{0i}^\gamma, \,\,\,\,\,\,\,\,Q_{i|\mu}^{\alpha\beta}\in \Gamma\left(U_i, \mathcal{A}^{0,0}\right)
\end{align*}
and we have
\begin{align*}
&\sum_{\beta,\gamma=1}^p r_{0ij}^{\alpha\beta} Q_{j|\mu}^{\beta \gamma} T_{0j}^\gamma - \sum_{\gamma=1}^p Q_{i|\mu}^{\alpha \gamma} T_{0i}^\gamma = -\sum_{\beta=1}^p \left[ \mathfrak{d}G\tilde{\varphi}_\mu , r_{0ij}^{\alpha\beta}  \right]T_{0j}^\beta\\
&\Longrightarrow \sum_{\beta,\gamma=1}^p r_{0ij}^{\alpha\beta} \left(\tilde{b}_{j|\mu}^{\beta \gamma}- \bar{\partial} Q_{j|\mu}^{\beta \gamma} \right)T_{0j}^\gamma - \sum_{\gamma=1}^p \left( \tilde{b}_{i|\mu}^{\alpha \gamma}  - \bar{\partial} Q_{i|\mu}^{\alpha \gamma} \right)T_{0i}^\gamma = \sum_{\beta=1}^p \left[ - \tilde{\varphi}_\mu  + \bar{\partial}\mathfrak{d}G\tilde{\varphi}_\mu , r_{0ij}^{\alpha\beta}  \right]T_{0j}^\beta = - \sum_{\beta=1}^p\left[ \bold{H}\tilde{\varphi}_\mu + \mathfrak{d} \bar{\partial} G\tilde{\varphi}_\mu , r_{0ij}^{\alpha\beta} \right] T_{0j}^\beta
\end{align*}

Then we see that
\begin{align}\label{d41}
\left(\varphi_\mu'', \left\{ b_{i|\mu}'' \right\} \right) \in A^{0,1}\left( M, \mathcal{E}_{\Theta_{\mathcal{F}_0}} \right),\,\,\,\,\,\,\,\,\,\, \mathfrak{d}\varphi_\mu''=0,\,\,\,\,\,\,\,\bar{\partial}\left(\varphi_\mu'' , \left\{ b_{i|\mu}''\right\} \right)= \left( - \xi_\mu,  \left\{ \bar{\partial} \Lambda_{i|\mu} \right\} \right)
\end{align}
where $\varphi_\mu'':= \bold{H}\tilde{\varphi}_\mu + \mathfrak{d}\bar{\partial} G \tilde{\varphi}_\mu$ and $b_{i|\mu}''\in \Gamma\left(U_i, \mathcal{A}^{0,1}\left( \mathscr{H}om_{\mathcal{O}_M}\left( \Theta_{\mathcal{F}_0}, \Theta_{\mathcal{F}_0} \right) \right) \right)$ defined by
\begin{align*}
b_{i|\mu}'':\Gamma\left( U_i, \Theta_{\mathcal{F}_0} \right) &\to \Gamma\left( U_i, \mathcal{A}^{0,1}\left( \Theta_{\mathcal{F}_0} \right) \right) \\
 T_{0i}^\alpha &\mapsto \sum_{\beta=1}^p \left( \tilde{b}_{i|\mu}^{\alpha \beta} - \bar{\partial} Q_{i|\mu}^{\alpha \beta} \right) T_{0i}^\beta
\end{align*}

Then by our construction of $\varphi_\mu'' \in A^{0,1}\left( M, \Theta_M \right)$ and $b_{i|\mu}'' \in \Gamma\left( U_i, \mathcal{A}^{0,1}\left(\mathscr{H}om_{\mathcal{O}_M}\left( \Theta_{\mathcal{F}_0}, \Theta_{\mathcal{F}_0} \right) \right) \right)$, we have from $(\ref{d25})$
\begin{align}
&\left| \varphi_\mu'' \right|_{k+\alpha} = \left|\bold{H} \tilde{\varphi}_\mu  + \mathfrak{d} \bar{\partial} G \tilde{\varphi}_\mu \right|_{k+\alpha} = \left| \tilde{\varphi}_\mu  - \bar{\partial}\mathfrak{d} G \tilde{\varphi}_\mu \right|_{k+\alpha} \ll K_2 \left|\tilde{\varphi}_\mu \right|_{k+\alpha} \ll K_1K_2 \left| \left(-\xi_\mu, \left\{ \bar{\partial} \Lambda_{i|\mu} \right\} \right) \right|_{k-1+\alpha}\\
&\left| b_{i|\mu}'' \right|_{k+\alpha}^{U_i}= \left| \tilde{b}_{i|\mu} - \bar{\partial} Q_{i|\mu} \right|_{k+\alpha}^{U_i} \ll \left| \tilde{b}_{i|\mu}\right|_{k+\alpha}^{U_i} + K_3  \left| Q_{i|\mu} \right|_{k+ 1+ \alpha}^{U_i} \ll  \left| \tilde{b}_{i|\mu}\right|_{k+\alpha}^{U_i} + K_3 K_4 \left| \mathfrak{d} G \tilde{\varphi}_\mu \right|_{k+\alpha+1} \ll  K_5\left|\left( - \xi_\mu , \left\{ \bar{\partial} \Lambda_{i|\mu}\right\} \right) \right|_{k-1+\alpha}     \notag
\end{align}
Hence we have
\begin{align} \label{ds3}
\left| \varphi_\mu'' \right|_{k+\alpha},\,\,\,\,\,\, \left| b_{i|\mu}'' \right|_{k+\alpha}^{U_i} \ll K_6 \left|\left( - \xi_\mu , \left\{ \bar{\partial} \Lambda_{i|\mu}\right\} \right) \right|_{k-1+\alpha}
\end{align}
where $K_6$ is a constant independent of $\left( - \xi_\mu, \left\{ \bar{\partial} \Lambda_{i|\mu} \right\} \right)$.

From $(\ref{d13})$, we see that
\begin{align*}
& 0 \in A^{0,2}\left(M , \mathcal{E}_{\Theta_{\mathcal{F}_0}}  \right)\\
& \left\{ T_{0i}^\alpha \mapsto \Phi_{i|\mu}^\alpha - \bar{\partial} \Gamma_{i|\mu}^\alpha +\sum_{\xi=1}^p \Lambda_{i|\mu}^{\alpha\xi} T_{0i}^\xi +\left[ \varphi_\mu'', T_{0i}^\alpha \right]  - \sum_{\xi=1}^p b_{i|\mu}''^{\alpha  \xi} T_{0i}^\xi   \right\}  \in  A^{0,1}\left(M, \mathscr{H}om_{\mathcal{O}_M}\left(\Theta_{\mathcal{F}_0}, \Theta_M \right) \right) \\
&\overline{\left\{ T_{0i}^\alpha \wedge T_{0i}^\beta \mapsto \Pi_{i|\mu}^{\alpha\beta} + \left[ \Gamma_{i|\mu}^\alpha, T_{0i}^\beta \right] - \left[ \Gamma_{i|\mu}^\beta, T_{0i}^\alpha  \right] - \sum_{\gamma=1}^p g_{0i\alpha\beta}^\gamma \Gamma_{i|\mu}^\gamma + \sum_{\gamma=1}^p \Psi_{i|\mu}^{\alpha\beta \gamma} T_{0i}^\gamma    \right\} } \in \frac{ A^{0,0}\left(M, \mathscr{H}om_{\mathcal{O}_M}\left( \bigwedge^2 \Theta_{\mathcal{F}_0}, \Theta_M \right)\right) }{ A^{0,0}\left( M, \mathscr{H}om_{\mathcal{O}_M}\left( \bigwedge^2 \Theta_{\mathcal{F}_0}, \Theta_{\mathcal{F}_0} \right) \right) }
\end{align*}
defines a $2$-cocycle in the above Dolbeault resolution of $\mathcal{E}_{\Theta_{\mathcal{F}_0}}^\bullet$. Since $\mathbb{H}^2\left( M, \Theta_{\mathcal{F}_0}^\bullet \right)=\mathbb{H}^2\left( M, \mathcal{E}_{\Theta_{\mathcal{F}_0}}^\bullet \right)=0$ by the assumption of Theorem \ref{tt2}, there exists $\left(\chi_\mu, \left\{ \eta_i^\chi \right\} \right)\in A^{0,1}\left(M, \mathcal{E}_{\Theta_{\mathcal{F}_0}} \right)$ and $\sigma_\mu \in A^{0,0}\left( M, \mathscr{H}om_{\mathcal{O}_M}\left(\Theta_{\mathcal{F}_0}, \Theta_M \right)\right)$  such that $-\bar{\partial}\left( \chi_\mu, \left\{ \eta_i^{\chi_\mu} \right\} \right)=0, -\bar{\partial} \sigma_\mu + \hat{D}_0'\left(\chi_\mu, \left\{ \eta_i^{\chi_\mu} \right\} \right)=\left\{ T_{0i}^\alpha \mapsto \Phi_{i|\mu}^\alpha - \bar{\partial} \Gamma_{i|\mu}^\alpha +\sum_{\xi=1}^p \Lambda_{i|\mu}^{\alpha\xi} T_{0i}^\xi +\left[ \varphi_\mu'', T_{0i}^\alpha \right]  - \sum_{\xi=1}^p b_{i|\mu}''^{\alpha  \xi} T_{0i}^\xi   \right\}$ and we have $\hat{D}_1'\left( \sigma_\mu\right)=$ $\overline{\left\{ T_{0i}^\alpha \wedge T_{0i}^\beta \mapsto \Pi_{i|\mu}^{\alpha\beta} + \left[ \Gamma_{i|\mu}^\alpha, T_{0i}^\beta \right] - \left[ \Gamma_{i|\mu}^\beta, T_{0i}^\alpha  \right] - \sum_{\gamma=1}^p g_{0i\alpha\beta}^\gamma \Gamma_{i|\mu}^\gamma + \sum_{\gamma=1}^p \Psi_{i|\mu}^{\alpha\beta \gamma} T_{0i}^\gamma    \right\} } $. By using the following Lemma, we will choose appropriate elements $\left(\chi_\mu, \left\{ \eta_i^{\chi_\mu} \right\} \right)$ and $\sigma_\mu$ in a way that the resulting $\varphi, T_i^\alpha, r_{ij}^{\alpha\beta}$ and $g_{i\alpha\beta}^\gamma$ converge.

\begin{lemma}[compar \cite{Hor76} Lemma 3.3]\label{cvt1}
Suppose that $\left(\varphi, \left\{\eta_i \right\}\right)\in A^{0,1}\left(M, \mathcal{E}_{\Theta_{\mathcal{F}_0}}\right)$ where  $\varphi \in A^{0,1}\left(M, \Theta_M\right)$ and $\eta_i\in\Gamma\left(U_i, \mathcal{A}^{0,1}\left(\mathscr{H}om_{\mathcal{O}_M}\left(\Theta_{\mathcal{F}_0}, \Theta_{\mathcal{F}_0} \right)  \right) \right)$ with $\eta_i\left(T_{0i}^\alpha \right)= \eta_j\left( T_{0i}^\alpha\right)+ \sum_{\beta=1}^p \left[ \varphi, r_{0ij}^{\alpha\beta} \right]T_{0j}^\beta$ for $\alpha=1,...,p$, and $\Phi\in A^{0,0}\left(M, \mathscr{H}om_{\mathcal{O}_M}\left( \Theta_{\mathcal{F}_0} , \Theta_M \right)\right)$ and $B\in A^{0,0}\left(M, \mathscr{H}om_{\mathcal{O}_M}\left( \bigwedge^2 \Theta_{\mathcal{F}_0} , \Theta_M \right) \right)$ such that $\left(0,  \Phi + \hat{D}_0' \left(\varphi, \left\{ \eta_i \right\}   \right) , \overline{B }\right)$ defines a $2$-cocycle  in the Dolbeault resolution of $\mathcal{E}_{\Theta_{\mathcal{F}_0}}^\bullet$. Then we can find $\left(\chi, \left\{ \eta_i^\chi\right\}\right) \in A^{0,1}\left( M,  \mathcal{E}_{\Theta_{\mathcal{F}_0}} \right)$ with $\chi \in A^{0,1}(M, \Theta_M)$ and  $\eta_i^\chi \in \Gamma\left(U_i, \mathcal{A}^{0,1}\left(\mathscr{H}om_{\mathcal{O}_M}\left( \Theta_{\mathcal{F}_0}, \Theta_{\mathcal{F}_0} \right) \right)\right)$, and $\sigma \in A^{0,0}\left( M, \mathscr{H}om_{\mathcal{O}_M}\left(\Theta_{\mathcal{F}_0}, \Theta_M \right)     \right)$ and the associated element $\pi_{i\sigma} \in \Gamma\left(U_i, \mathscr{H}om_{\mathcal{O}_M}\left( \bigwedge^2 \Theta_{\mathcal{F}_0}, \Theta_{\mathcal{F}_0} \right)\right)$ in such a way that
\begin{align}
\widetilde{\square} \left(\chi, \left\{ \eta_i^\chi \right\} \right) &=0 \label{tpc13}\\
-\bar{\partial} \sigma\left( T_{0i}^\alpha \right) + \left[\chi, T_{0i}^\alpha \right] -  \eta_i^\chi \left( T_{0i}^\alpha \right) &= \Phi\left( T_{0i}^\alpha \right) + \left[ \varphi , T_{0i}^\alpha \right] - \eta_i\left( T_{0i}^\alpha \right) \label{tpc11}\\
\left[ T_{0i}^\alpha , \sigma\left( T_{0i}^\beta \right) \right]- \left[ T_{0i}^\beta, \sigma\left(T_{0i}^\alpha \right) \right] &- \sum_{\gamma=1}^p g_{0i\alpha\beta}^\gamma \sigma\left(T_{0i}^\gamma \right)+ \pi_{i\sigma}\left( T_{0i}^\alpha \wedge T_{0i}^\beta \right)= B\left( T_{0i}^\alpha \wedge T_{0i}^\beta \right):= B_i^{\alpha\beta} \label{tpc10}\\
\eta_i^\chi(T_{0i}^\alpha)&= \eta_j^\chi \left(T_{0i}^\alpha \right)+ \sum_{\beta=1}^p \left[\chi, r_{0ij}^{\alpha\beta} \right]T_{0j}^\beta \\
\pi_{j\sigma}\left(T_{0i}^\alpha \wedge T_{0i}^\beta \right)- \pi_{i\sigma} \left(T_{0i}^\alpha \wedge T_{0i}^\beta \right)&= \sum_{\eta=1}^p \left[ \sigma\left(T_{0i}^\alpha \right), r_{0ij}^{\alpha\beta } \right] T_{0j}^\eta- \sum_{\eta=1}^p \left[ \sigma\left(T_{0j}^\beta \right), r_{0ij}^{\alpha \eta} \right] T_{0j}^\eta \label{tpc14} \\
\left| \left(\chi, \left\{\eta_i^\chi \right\} \right) \right|_{k+\alpha} \ll K_2&\left( \left|(\varphi,  \left\{ \eta_i \right\})\right|_{k+\alpha} + \left| \Phi \right|_{k-1+\alpha} +  \left| B \right|_{k-1+\alpha}  \right)\\
\left| \sigma \right|_{k+\alpha} \ll K_2&\left( \left| \left(\varphi, \left\{\eta_i \right\}\right)\right|_{k+\alpha} + \left| \Phi \right|_{k-1+\alpha} +  \left| B \right|_{k-1+\alpha}  \right)\\
\left| \pi_{i\sigma} \right|_{k - 1 +\alpha}^{U_i} \ll K_2&\left( \left| \left(\varphi , \left\{\eta_i \right\}\right) \right|_{k+\alpha} + \left| \Phi \right|_{k-1+\alpha} +  \left| B \right|_{k-1+\alpha}  \right)
\end{align}
where $K_2$ is a constant which is independent of $\left(\varphi, \left\{\eta_i \right\} \right), \Phi, B$. We note that $(\ref{tpc11})$ implies $-\bar{\partial} \sigma + \hat{D}_0' \left(\chi, \left\{ \eta_i^\chi \right\} \right)= \Phi + \hat{D}_0' \left(\varphi , \left\{ \eta_i\right\} \right)$ and $(\ref{tpc10})$ implies $\hat{D}_1'\left(\sigma \right)= \overline{B}$.
\end{lemma}

\begin{proof}
For any triple $\zeta=\left( \left( \varphi,  \left\{\eta_i \right\}\right), \Phi, B \right)$ as above, let
\begin{align*}
||\zeta || &= \left| \left(\varphi  , \left\{ \eta_i^\varphi \right\}\right)\right|_{k+\alpha} + \left|\Phi \right|_{k-1+\alpha} + \left|B \right|_{k-1+\alpha} \\
\iota(\zeta)&= \inf \left( \left|\left(\chi, \left\{\eta_i^\chi \right\} \right)\right|_{k+\alpha}  + \left|\sigma \right|_{k+\alpha} +\max_i \left|\pi_{i\sigma} \right|_{k-1+\alpha}^{U_i}  \right)
\end{align*}
where $\inf$ is taken with respect to all solutions $\left( \left(\chi, \left\{\eta_i^\chi \right\}\right),\left( \sigma, \left\{ \pi_{i\sigma}\right\}\right) \right)$ of the equalities $(\ref{tpc13})-(\ref{tpc14})$. It suffices to prove the existence of $K$ such that
\begin{align*}
\iota\left( \zeta \right)\leq K ||\zeta ||\,\,\,\,\,\,\,\,\,\textnormal{for all triple $\zeta$}
\end{align*}
Assume that there is no such constant $K$. Then we can find a sequence $\zeta^{(1)},\zeta^{(2)},\cdots, \zeta^{(v)},\cdots$ of $\zeta^{(v)}=\left( \left( \varphi^{(v)} , \left\{ \eta_i^{(v)} \right\} \right), \Phi^{(v)}, B^{(v)}\right)$ such that
\begin{align*}
\iota\left(\zeta^{(v)} \right)=1\,\,\,\,\,\,\,\,\textnormal{and}\,\,\,\,\,\,\,\,\left|\left| \zeta^{(v)}\right|\right| < \frac{1}{v}
\end{align*}
The first equality implies the existence of $\left(\chi^{(v)}, \left\{ \eta_i^{\chi^{(v)}} \right\} \right)\in A^{0,1}\left( M, \mathcal{E}_{\Theta_{\mathcal{F}_0}} \right)$ and $\sigma^{(v)}\in A^{0,0}\left( M, \mathscr{H}om_{\mathcal{O}_M}\left( \Theta_{\mathcal{F}_0}, \Theta_M \right) \right)$ and $\pi_{i\sigma^{(v)}} \in \Gamma\left( U_i, \mathcal{A}^{0,0}\left( \mathscr{H}om_{\mathcal{O}_M}\left( \bigwedge^2 \Theta_{\mathcal{F}_0}, \Theta_M \right) \right) \right)$ such that
\begin{align}
&\widetilde{\square} \left(\chi^{(v)}, \left\{ \eta_i^{\chi^{(v)}} \right\} \right)=0  \\
-\bar{\partial} \sigma^{(v)}\left( T_{0i}^\alpha \right) + \left[\chi^{(v)}, T_{0i}^\alpha \right]& -  \eta_i^{\chi^{(v)}} \left( T_{0i}^\alpha \right) = \Phi^{(v)}\left( T_{0i}^\alpha \right) + \left[ \varphi^{(v)} , T_{0i}^\alpha \right] - \eta_i^{(v)}\left( T_{0i}^\alpha \right) \label{tpc16} \\
\left[ T_{0i}^\alpha , \sigma^{(v)}\left( T_{0i}^\beta \right) \right]- \left[ T_{0i}^\beta, \sigma^{(v)}\left(T_{0i}^\alpha \right) \right] &- \sum_{\gamma=1}^p g_{0i\alpha\beta}^\gamma \sigma^{(v)}\left(T_{0i}^\gamma \right)+ \pi_{i\sigma^{(v)}}\left( T_{0i}^\alpha \wedge T_{0i}^\beta \right)= B^{(v)}\left( T_{0i}^\alpha \wedge T_{0i}^\beta \right):= B_i^{(v)\alpha\beta} \label{tpc17} \\
\eta_i^{\chi^{(v)}}(T_{0i}^\alpha)&= \eta_j^{\chi^{(v)}} \left(T_{0i}^\alpha \right)+ \sum_{\beta=1}^p \left[\chi^{(v)}, r_{0ij}^{\alpha\beta} \right]T_{0j}^\beta \\
\pi_{j\sigma^{(v)}}\left(T_{0i}^\alpha \wedge T_{0i}^\beta \right)- &\pi_{i\sigma^{(v)}} \left(T_{0i}^\alpha \wedge T_{0i}^\beta \right)= \sum_{\eta=1}^p \left[ \sigma^{(v)}\left(T_{0i}^\alpha \right), r_{0ij}^{\alpha\beta } \right] T_{0j}^\eta- \sum_{\eta=1}^p \left[ \sigma^{(v)}\left(T_{0j}^\beta \right), r_{0ij}^{\alpha \eta} \right] T_{0j}^\eta \label{sd9} \\
\left| \left(\chi^{(v)}, \left\{ \eta_i^{\chi^{(v)}}\right\} \right) \right|_{k+\alpha}&  + \left|\sigma^{(v)} \right|_{k+\alpha} + \max_i \left| \pi_{i\sigma^{(v)}} \right|_{k-1+\alpha}^{U_i} < 2 \label{tpc15}
\end{align}
From $(\ref{tpc15})$, by replacing $\zeta^{(1)}, \zeta^{(2)}, \cdots$ by a suitable subsequence if necessary, we may assume that $\left( \chi, \left\{ \eta_i^\chi\right\} \right)=\lim \left( \chi^{(v)}, \left\{ \eta_i^{\chi^{(v)}} \right\} \right)$ and $\sigma=\lim \sigma^{(v)}$ exists in the norm $|-|_{k}$ and $\pi_{i\sigma}=\lim \pi_{i\sigma^{(v)}}$ exists in the norm $|-|_{k-1}$ by compact embedding for \" Holder norms, so that $\left( \chi, \left\{ \eta_i^{\chi} \right\} \right)$ and $\sigma$ are of class $C^k$ and $\pi_{i\sigma}$ are of class $C^{k-1}$. We note that since $\left(\chi, \left\{\eta_i^{\chi}\right\} \right)$ satisfies an elliptic partial differential equation $\widetilde{\square}\left( \chi, \left\{ \eta_i^\chi \right\} \right)=0$ , $\left( \chi, \left\{ \eta_i^\chi \right\} \right)$ is of $C^\infty$. Moreover, since $\widetilde{\square}$ is strongly elliptic, we have a priori estimate (see \cite{Kod05} Theorem 4.3 p. 436)
\begin{align}\label{ns10}
\left| \left(\chi^{(v)} ,\left\{ \eta_i^{\chi^{(v)}} \right\} \right)- \left( \chi, \left\{\eta_i^\chi \right\} \right) \right|_{k+\alpha} &\leq C_0\left( \left| \widetilde{\square }\left( \left(\chi^{(v)} ,\left\{ \eta_i^{\chi^{(v)}} \right\} \right)- \left( \chi, \left\{\eta_i^\chi \right\} \right)\right)    \right|_{k-2 +\alpha} + \left| \left(\chi^{(v)} ,\left\{ \eta_i^{\chi^{(v)}} \right\} \right)- \left( \chi, \left\{\eta_i^\chi \right\} \right) \right|_0 \right) \\
&= C_0 \left|  \left(\chi^{(v)} ,\left\{ \eta_i^{\chi^{(v)}} \right\} \right)- \left( \chi, \left\{\eta_i^\chi \right\} \right) \right|_0 \notag
\end{align}
for some constant $C_0>0$. Hence $\left(\chi^{(v)}, \left\{ \eta_i^{\chi^{(v)}} \right\} \right)$ converges to $\left( \chi , \left\{ \eta_i^\chi \right\} \right)$ in the norm $|-|_{k+\alpha}$. On the other hand, from $(\ref{tpc16})$ we have $\bar{\partial} \sigma \left( T_{0i}^\alpha \right)=  \left[\chi, T_{0i}^\alpha \right] -  \eta_i^{\chi} \left( T_{0i}^\alpha \right)= \hat{D}_0\left( \chi, \left\{  \eta_i^\chi \right\} \right)$. Setting $\sigma'=\mathfrak{d}' G'\left(\hat{D}_0\left(\chi, \left\{ \eta_i^\chi \right\} \right) \right)$ which is of $C^\infty$, we have $\bar{\partial}\left(\sigma- \sigma' \right)=0$, so that since $\sigma- \sigma'$ is of $C^k$, $\sigma- \sigma' $ is $C^\infty$. Hence $\sigma$ is of $C^\infty$. 

Since $\bar{\partial}$ is elliptic of order $1$, we can apply the following estimate (see \cite{Kod05} p.438)
{\small{\begin{align}\label{ns12}
\left| \sigma^{(v)}\left( T_{0i}^\alpha\right)- \sigma \left(T_{0i}^\alpha\right) \right|_{k+\alpha} & \leq C \left(\left| \bar{\partial}\left(\sigma^{(v)} \left(T_{0i}^\alpha \right) - \sigma \left( T_{0i}^\alpha \right)\right)\right|_{k-1+\alpha} + \left| \sigma^{(v)} \left( T_{0i}^\alpha \right)- \sigma \left(T_{0i}^\alpha \right)\right|_0 \right)  \\
 &=  C  \left( \left|   \left[ \chi^{(v)}- \chi , T_{0i}^\alpha\right]- \left(\eta_i^{(v)}- \eta_i\right)\left(T_{0i}^\alpha \right) - \Phi^{(v)}\left(T_{0i}^\alpha \right) - \left[ \varphi^{(v)}, T_{0i}^\alpha \right]  + \eta_i^{(v)}\left(T_{0i}^\alpha \right)  \right|_{k-1 + \alpha} + \left| \sigma^{(v)}\left( T_{0i}^\alpha \right)- \sigma \left(T_{0i}^\alpha \right) \right|_0            \right) \notag
\end{align}}}
for some constant $C>0$. Since $\left( \chi^{(v)}-\chi, \left\{ \eta_i^{\chi(v)}- \eta_i^\chi \right\} \right), \Phi^{(v)}, \left( \varphi^{(v)}, \left\{ \eta_i^{(v)}\right\} \right)$ converge to $0$ in $|-|_{k-1+\alpha}$ and $\left| \sigma^{(v)}-\sigma \right|_0 \to 0$ by our construction, $\sigma^{(v)}$ converge to $\sigma$ in $|-|_{k+\alpha}$. From $(\ref{tpc17})$, we have $\pi_{i\sigma}\left(T_{0i}^\alpha \wedge T_{0i}^\beta \right)= - \left[ T_{0i}^\alpha, \sigma\left( T_{0i}^\beta \right) \right] + \left[ T_{0i}^\beta , \sigma\left( T_{0i}^\alpha \right)\right] + \sum_{\gamma=1}^p  g_{0i\alpha\beta}^\gamma \sigma\left( T_{0i}^\gamma \right)$, so that $\pi_{i\sigma}$ is of $C^\infty$. From $(\ref{tpc17})$, we have
{\small{\begin{align*}
&\left| \left(\pi_{i\sigma^{(v)}}- \pi_{i\sigma}\right)\left( T_{0i}^\alpha \wedge T_{0i}^\beta \right)    \right|_{k-1+\alpha} \\
&=  \left| -\left[ T_{0i}^\alpha, \left(\sigma^{(v)}- \sigma\right)\left( T_{0i}^\beta \right) \right]   + \left[ T_{0i}^\alpha, \left(\sigma^{(v)}- \sigma\right)\left( T_{0i}^\beta \right) \right] + \sum_{\gamma=1}^p g_{0i\alpha\beta}^\gamma \left(\sigma^{(v)}- \sigma \right)\left( T_{0i}^\gamma \right) + B^{(v)}\left(T_{0i}^\alpha \wedge T_{0i}^\beta \right)\right|_{k-1+\alpha}\\
&\leq C_1\left| \sigma^{(v)} - \sigma \right|_{k+\alpha} + \left| B^{(v)} \right|_{k-1+\alpha}
\end{align*}}}
This implies that $\pi_{i\sigma^{(v)}}$ converges to $\pi_{i\sigma}$ in the norm $|-|_{k-1+\alpha}$. From $(\ref{tpc16})$ and $(\ref{tpc17})$, we have
\begin{align}
-\bar{\partial} \left( \sigma^{(v)}-\sigma\right)\left( T_{0i}^\alpha \right) + \left[\chi^{(v)}- \chi, T_{0i}^\alpha \right]& -  \eta_i^{\chi^{(v)}-\chi} \left( T_{0i}^\alpha \right) = \Phi^{(v)}\left( T_{0i}^\alpha \right) + \left[ \varphi^{(v)} , T_{0i}^\alpha \right] - \eta_i^{(v)}\left( T_{0i}^\alpha \right) \label{sd15}  \\
\left[ T_{0i}^\alpha , \left(\sigma^{(v)}-\sigma\right)\left( T_{0i}^\beta \right) \right]- \left[ T_{0i}^\beta, \left(\sigma^{(v)}-\sigma \right)\left(T_{0i}^\alpha \right) \right] &- \sum_{\gamma=1}^p g_{0i\alpha\beta}^\gamma \left( \sigma^{(v)}-\sigma \right)\left(T_{0i}^\gamma \right)+ \left( \pi_{i\sigma^{(v)}}-\pi_{i\sigma}\right)\left( T_{0i}^\alpha \wedge T_{0i}^\beta \right)= B^{(v)}\left( T_{0i}^\alpha \wedge T_{0i}^\beta \right) \label{sd16}
\end{align}
On the other hand, we have
\begin{align*}
\left|  \left( \chi^{(v)}, \left\{ \eta_i^{\chi^{(v)}} \right\} \right) - \left( \chi , \left\{ \eta_i^\chi \right\} \right)\right|_{k+\alpha} <\frac{1}{4},\,\,\,\,\,\,\,\,\,\left| \sigma^{(v)}- \sigma \right|_{k+\alpha} <\frac{1}{4},\,\,\,\,\,\,\,\,\,\,\, \left| \pi_{i\sigma^{(v)}} -  \pi_{i\sigma} \right|_{k-1+\alpha} < \frac{1}{4}
\end{align*}
for sufficiently large integer $v$. This contradicts to $\iota\left( \zeta \right)=1$. This completes the proof of Lemma \ref{cvt1}.

\end{proof}

In Lemma \ref{cvt1}, we set
\begin{align}\label{ds20}
& \left( \varphi, \left\{ \eta_i \right\} \right):= \left( \varphi_\mu'', \left\{ b_{i|\mu}'' \right\} \right)\in A^{0,1}\left( M, \mathcal{E}_{\Theta_{\mathcal{F}_0}} \right)\\
& \Phi:= \left\{ T_{0i}^\alpha \mapsto \Phi_{i|\mu}^\alpha - \bar{\partial} \Gamma_{i|\mu}^\alpha +\sum_{\xi=1}^p \Lambda_{i|\mu}^{\alpha\xi} T_{0i}^\xi   \right\}  \in  A^{0,1}\left(M, \mathscr{H}om_{\mathcal{O}_M}\left(\Theta_{\mathcal{F}_0}, \Theta_M \right) \right) \notag \\
&B:= \left\{ T_{0i}^\alpha \wedge T_{0i}^\beta \mapsto \Pi_{i|\mu}^{\alpha\beta} + \left[ \Gamma_{i|\mu}^\alpha, T_{0i}^\beta \right] - \left[ \Gamma_{i|\mu}^\beta, T_{0i}^\alpha  \right] - \sum_{\gamma=1}^p g_{0i\alpha\beta}^\gamma \Gamma_{i|\mu}^\gamma + \sum_{\gamma=1}^p \Psi_{i|\mu}^{\alpha\beta \gamma} T_{0i}^\gamma    \right\}  \in A^{0,0}\left(M, \mathscr{H}om_{\mathcal{O}_M}\left( \bigwedge^2 \Theta_{\mathcal{F}_0}, \Theta_M \right)\right)  \notag
\end{align}

We will estimate $\left(\varphi, \left\{ \eta_i\right\} \right)$ and $\Phi$ and $B$. First we will estimate $(\ref{tt3})-(\ref{ds47})$. By induction hypothesis $(\ref{d20})_{\mu-1}- (\ref{d21})_{\mu-1}$, we have
\begin{align}\label{ds22}
\xi_\mu\equiv_\mu  \bar{\partial} \varphi^{\mu-1}-\frac{1}{2}\left[\varphi^{\mu-1}, \varphi^{\mu-1}   \right] \Longrightarrow \left|\xi_\mu \right|_{k-1+\alpha} \leq K_8 \left|\varphi^{\mu-1} \right|_{k+\alpha} \left|\varphi^{\mu-1}\right| _{k+\alpha} \leq \frac{K_8b}{c}A(t)
\end{align}
and we have
\begin{align}
&-\Phi_{i|\mu}^\alpha \equiv_\mu \bar{\partial} T_i^{\alpha(\mu-1)} -\left[\varphi^{\mu-1}, T_i^{\alpha(\mu-1)}   \right]\equiv_\mu \left[ \varphi^{\mu-1}, T_i^{\alpha(\mu-1)}   \right]_\mu \equiv_\mu \left[\varphi^{\mu-1}, T_i^{\alpha(\mu-1)}- T_{0i}^\alpha     \right]_\mu \notag \\
&\Longrightarrow \left| \Phi_{i|\mu}^\alpha \right|_{k-1+\alpha} \leq K_9 \left|\varphi^{\mu-1} \right|_{k+\alpha} \left| T_i^{\alpha(\mu-1)}  -T_{0i}^\alpha \right|_{k+\alpha} \leq  K_9  \frac{b}{c} A(t) \label{ds23} 
\end{align}
and we have
\begin{align}\label{ds11}
&-\sum_{\eta=1}^p \Lambda_{ij|\mu}^{\alpha\eta} r_{0ij}^{\eta \beta} \equiv_\mu \bar{\partial} r_{ij}^{\alpha\beta(\mu-1)} -\left[ \varphi^{\mu-1}, r_{ij}^{\alpha\beta(\mu-1)}    \right] \equiv_\mu -\left[\varphi^{\mu-1}, r_{ij}^{\alpha\beta(\mu-1)}- r_{0ij}^{\alpha\beta} \right]_\mu \notag \\
& \Longrightarrow  \left| \Lambda_{ij|\mu}^{\alpha\beta} \right|_{k+\alpha} \leq K_{10} \left| \varphi^{\mu-1} \right|_{k+\alpha}\left| r_{ij}^{\alpha\beta(\mu-1)} - r_{0ij}^{\alpha\beta}\right|_{k+\alpha+1} \leq K_{10} \frac{b}{c} A(t)
\end{align}
and we have
\begin{align}
&\Gamma_{ij|\mu}^\alpha \equiv_\mu T_i^{\alpha(\mu-1)} - \sum_{\beta=1}^p r_{ij}^{\alpha\beta(\mu-1)} T_j^{\beta(\mu-1)} \equiv_\mu \left[- \sum_{\beta=1}^p \left(r_{ij}^{\alpha\beta(\mu-1)} - r_{0ij}^{\alpha\beta} \right)\left( T_j^{\beta(\mu-1)} - T_{0j}^{\beta}\right)     \right]_\mu \notag \\
& \Longrightarrow \left| \Gamma_{ij|\mu}^\alpha \right|_{k+\alpha} \sum_{\beta=1}^p \left|r_{ij}^{\alpha\beta(\mu-1)} - r_{0ij}^{\alpha\beta} \right|_{k+\alpha} \left| T_j^{\beta(\mu-1)} - T_{0j}^\beta  \right|_{k+\alpha} \leq K_{12} \frac{b}{c} A(t) \label{ds17}
\end{align}
and we have
{\small{\begin{align}
&\Pi_{i|\mu}^{\alpha\beta} \equiv_\mu \left[ T_i^{\alpha(\mu-1)} , T_i^{\beta(\mu-1)}  \right] -\sum_{\beta=1}^p g_{i\alpha\beta}^{\gamma(\mu-1)} T_i^{\gamma(\mu-1)} \equiv_\mu \left[ T_i^{\alpha(\mu-1)} - T_{0i}^\alpha, T_i^{\beta(\mu-1)} - T_{0i}^\beta \right]_\mu - \left[ \sum_{\beta=1}^p \left(g_{i\alpha\beta}^{\gamma(\mu-1)}- g_{0i\alpha\beta}^{\gamma(\mu-1)} \right)\left( T_i^{\gamma(\mu-1)} - T_{0i}^\gamma \right) \right]_\mu \notag \\
&\Longrightarrow \left| \Pi_{i|\mu}^{\alpha\beta} \right|_{k-1+\alpha} \ll K_{13}\left|T_i^{\alpha(\mu-1)} - T_{0i}^\alpha \right|_{k+\alpha} \left| T_i^{\beta(\mu-1)} - T_{0i}^\beta \right|_{k+\alpha} + \sum_{\beta=1}^p \left| g_{i\alpha\beta}^{\gamma(\mu-1)}- g_{0i\alpha\beta}^{\gamma(\mu-1)} \right|_{k-1+\alpha}\left| T_i^{\gamma(\mu-1)} - T_{0i}^\gamma \right|_{k+\alpha} \ll K_{14} \frac{b}{c}A(t) \label{ds19}
\end{align}}}
and we have
\begin{align}\label{ds10}
&\sum_{\beta=1}^p \lambda_{ijk|\mu}^{\alpha\beta} r_{0ik}^{\beta \gamma}\equiv_\mu r_{ik}^{\alpha\gamma(\mu-1)} - \sum_{\beta=1}^p r_{ij}^{\alpha\beta(\mu-1)} r_{jk}^{\beta\gamma(\mu-1)} \equiv_\mu \left[ - \sum_{\beta=1}^p \left( r_{ij}^{\alpha\beta(\mu-1)}  - r_{ij}^{\alpha\beta}\right)\left( r_{jk}^{\beta\gamma(\mu-1)} - r_{0jk}^{\beta \gamma} \right) \right]_\mu \notag \\
&\Longrightarrow \left|\lambda_{ijk|\mu}^{\alpha\beta} \right|_{k+1+\alpha} \ll K_{15} \frac{b}{c}A(t)
\end{align}
and we have
{\small{\begin{align*}
\Psi_{ij|\mu}^{\alpha\beta\xi}&\equiv_\mu  \sum_{\eta, \gamma=1}^p r_{ij}^{\alpha\eta(\mu-1)} r_{ij}^{\beta \gamma(\mu-1)} g_{j\eta \gamma}^{\xi(\mu-1)} +   T_i^{\alpha(\mu-1)} \left( r_{ij}^{\beta \xi(\mu-1)} \right)     -    T_i^{\beta(\mu-1)}\left( r_{ij}^{\alpha \xi(\mu-1)} \right) -  \sum_{\gamma=1}^p g_{i\alpha\beta}^{\gamma(\mu-1)} r_{ij}^{\gamma \xi(\mu-1)}     \\
& \equiv_\mu  \left[ \sum_{\eta, \gamma=1}^p \left( r_{ij}^{\alpha \eta (\mu-1)} - r_{0ij}^{\alpha\eta} \right) \left( r_{ij}^{\beta \gamma(\mu-1)} - r_{0ij}^{\beta \gamma} \right) \left( g_{j\eta \gamma}^{\xi(\mu-1)} - g_{0j\eta\gamma}^\xi \right) + \sum_{\eta, \gamma=1}^p r_{0ij}^{\alpha\eta } \left( r_{ij}^{\beta \gamma(\mu-1)} - r_{0ij}^{\beta \gamma} \right) \left( g_{j\eta \gamma}^{\xi (\mu-1)} - g_{0j\eta \gamma}^\xi \right) \right]_\mu   \\
&+ \left[ \sum_{\eta, \gamma=1}^p \left( r_{ij}^{\alpha \eta (\mu-1)}- r_{0ij}^{\alpha\eta} \right) \left( r_{0ij}^{\beta \gamma} \right) \left( g_{j\eta \gamma}^{\xi(\mu-1)} - g_{0j\eta\gamma}^\xi \right) + \sum_{\eta, \gamma=1}^p \left( r_{ij}^{\alpha \eta (\mu-1)} - r_{0ij}^{\alpha\eta} \right) \left( r_{ij}^{\beta \gamma(\mu-1)} - r_{0ij}^{\beta \gamma} \right) \left(  g_{0j\eta\gamma}^\xi \right)\right]_\mu \\
& + \left[\left( T_i^{\alpha(\mu-1)} - T_{0i}^\alpha\right) \left( r_{ij}^{\beta\xi(\mu-1)} - r_{0ij}^{\beta \xi}  \right)\right]_\mu -  \left[ \left( T_i^{\beta(\mu-1)} - T_{0i}^\beta \right) \left(r_{ij}^{\alpha \xi(\mu-1)} - r_{0ij}^{\alpha \xi} \right)\right]_\mu - \left[ \sum_{\gamma=1}^p \left(g_{i\alpha\beta}^{\gamma (\mu-1)} - g_{0i\alpha\beta}^\gamma \right) \left( r_{ij}^{\gamma \xi(\mu-1)} - r_{0ij}^{\gamma \xi}\right)\right]_\mu
\end{align*}}}
\begin{align}
\Longrightarrow & \left| \Psi_{ij|\mu}^{\alpha\beta \xi}\right|_{k -1 + \alpha} \ll K_{16}\left( \frac{b}{c} \right)^2 A(t) + K_{17} \frac{b}{c} A(t) \ll K_{18}\frac{b}{c} A(t) \label{sd20}
\end{align}

We recall that from $(\ref{te18})$ and $(\ref{te21})$ and $(\ref{t57})$ and $(\ref{ds8})$
\begin{align}
\sum_{\beta, \eta=1}^p r_{0ij}^{\alpha \beta} \lambda_{jk|\mu}^{\beta \eta} r_{0ji}^{\eta \xi} - \lambda_{ik|\mu}^{\alpha \xi} &+ \lambda_{ij|\mu}^{\alpha \xi} =\lambda_{ijk|\mu}^{\alpha \xi} \label{ds13}\\
\sum_{\beta, \eta=1}^p r_{0ij}^{\alpha\beta} \Lambda_{j|\mu}^{\beta \eta} r_{0ij}^{\eta \xi} - \Lambda_{i|\mu}^{\alpha \xi}&= \Lambda_{ij|\mu}^{\alpha \xi} - \bar{\partial} \lambda_{ij|\mu}^{\alpha \xi} \label{ds14}\\
\sum_{\beta=1}^p r_{0ij}^{\alpha\beta} \Gamma_{j|\mu}^\beta - \Gamma_{i|\mu}^\alpha &= \Gamma_{ij|\mu}^\alpha - \sum_{\xi=1}^p \lambda_{ij|\mu}^{\alpha \xi} T_i^\xi \label{ds15}\\
\sum_{\eta, \delta, \xi=1}^p r_{0ij}^{\alpha \eta} r_{0ij}^{\beta \delta} \Psi_{j|\mu}^{\eta \delta \xi} T_{0j}^\xi - \sum_{\gamma=1}^p& \Psi_{i|\mu}^{\alpha\beta \gamma} T_{0i}^\gamma =\Pi_{i|\mu}^{\alpha\beta} + \left[ \Gamma_{i|\mu}^\alpha, T_{0i}^\beta \right] - \left[ \Gamma_{i|\mu}^\beta, T_{0i}^\alpha \right] - \sum_{\gamma=1}^p g_{0i\alpha\beta}^\gamma \Gamma_{i|\mu}^\gamma  \label{ds16} \\
&- \sum_{\eta, \xi=1}^p r_{0ij}^{\alpha \eta} r_{0ij}^{\beta \xi} \left( \Pi_{j|\mu}^{\eta \xi} + \left[ \Gamma_{j|\mu}^\eta, T_{0j}^\xi \right] - \left[ \Gamma_{j|\mu}^\xi, T_{0j}^\eta \right] - \sum_{\gamma=1}^p g_{0j\eta \xi}^\gamma \Gamma_{j|\mu}^\gamma   \right) \notag
\end{align}
By using a partition of unity subordinate to the covering $\left\{ U_i \right\}$, we may choose $\lambda_{ij|\mu}$ from $\lambda_{ijk|\mu}$, and choose $\Lambda_{j|\mu}$ from $\Lambda_{ij|\mu}$ and $\bar{\partial} \lambda_{ij|\mu}$ and choose $\Gamma_{i|\mu}$ from $\Gamma_{ij|\mu}$ and $\lambda_{ij|\mu}$ and choose $\Psi_{i|\mu}$ in a similar way. Hence we may assume that from $(\ref{ds10})$ and $(\ref{ds13})$
\begin{align}\label{ds12}
\left| \lambda_{ij|\mu}\right|_{k+1+\alpha} \ll K_{21}\frac{b}{c}A(t)
\end{align}
and from $(\ref{ds11})$ and $(\ref{ds12})$ and $(\ref{ds14})$ we may assume
\begin{align}\label{ds18}
\left|  \Lambda_{i|\mu} \right|_{k+\alpha} \ll K_{22} \frac{b}{c}A(t)
\end{align}
and from $(\ref{ds17})$ and $(\ref{ds12})$ we may assume 
\begin{align}\label{ds24}
\left| \Gamma_{i|\mu} \right|_{k+\alpha} \ll K_{23}\frac{b}{c}A(t)
\end{align}
and from $(\ref{ds16})$ and $(\ref{ds18})$ and $(\ref{ds19})$ we may assume
\begin{align} \label{ds25}
\left|\Psi_{i|\mu} \right|_{k-1+\alpha} \ll K_{24}\frac{b}{c} A(t)
\end{align}

With this preparation, now we estimate $(\ref{ds20})$ and then apply Lemma \ref{cvt1}. Let us estimate $\left( \varphi_\mu'', \left\{ b_{i|\mu}''\right\} \right)$. First we note that from $(\ref{ds3})$ and $(\ref{ds5})$. 
\begin{align*}
\left|\left( \varphi_\mu'', \left\{ b_{i|\mu}'' \right\}  \right) \right|_{k+\alpha} \ll K_7\left|\left(- \xi_\mu, \left\{ \bar{\partial} \Lambda_{i|\mu} \right\} \right) \right|_{k-1+\alpha}
\end{align*}
Then from $(\ref{ds22})$ and $(\ref{ds18})$, we have
\begin{align}\label{ds31}
\left|\left( \varphi_\mu'', \left\{ b_{i|\mu}'' \right\}  \right) \right|_{k+\alpha} \ll K_{25} \frac{b}{c} A(t)
\end{align}
We estimate  $\Phi:=\left\{ T_{0i}^\alpha \mapsto \Phi_{i|\mu}^\alpha - \bar{\partial} \Gamma_{i|\mu}^\alpha +\sum_{\xi=1}^p \Lambda_{i|\mu}^{\alpha\xi} T_{0i}^\xi   \right\}  \in  A^{0,1}\left(M, \mathscr{H}om_{\mathcal{O}_M}\left(\Theta_{\mathcal{F}_0}, \Theta_M \right) \right)$. From $(\ref{ds23})$ and $(\ref{ds24})$ and $(\ref{ds18})$, we have
\begin{align}\label{sd25}
\left|\Phi \right|_{k-1+\alpha} \ll K_{26} \frac{b}{c} A(t)
\end{align}

We estimate 
\begin{align*}
B:= \left\{ T_{0i}^\alpha \wedge T_{0i}^\beta \mapsto \Pi_{i|\mu}^{\alpha\beta} + \left[ \Gamma_{i|\mu}^\alpha, T_{0i}^\beta \right] - \left[ \Gamma_{i|\mu}^\beta, T_{0i}^\alpha  \right] - \sum_{\gamma=1}^p g_{0i\alpha\beta}^\gamma \Gamma_{i|\mu}^\gamma + \sum_{\gamma=1}^p \Psi_{i|\mu}^{\alpha\beta \gamma} T_{0i}^\gamma    \right\}  \in A^{0,0}\left(M, \mathscr{H}om_{\mathcal{O}_M}\left( \bigwedge^2 \Theta_{\mathcal{F}_0}, \Theta_M \right)\right)
\end{align*}
From $(\ref{ds19})$ and $(\ref{ds24})$ and $(\ref{ds25})$, we have
\begin{align} \label{sd26}
\left| B \right|_{k-1+\alpha} \ll K_{27} \frac{b}{c} A(t)
\end{align}

By Lemma \ref{cvt1}, we can find $\left(\chi_\mu , \left\{ \eta_i^{\chi_\mu} \right\} \right)\in A^{0,1}\left( M,  \mathcal{E}_{\Theta_{\mathcal{F}_0}}\right)$ and $\sigma_\mu \in A^{0,0}\left(M, \mathscr{H}om_{\mathcal{O}_M}\left( \Theta_{\mathcal{F}_0}, \Theta_M \right) \right)$ and $\pi_{i\sigma_\mu} \in \Gamma\left( U_i, \mathscr{H}om_{\mathcal{O}_M}\left( \bigwedge^2 \Theta_{\mathcal{F}_0} , \Theta_{\mathcal{F}_0} \right) \right)$ such that
\begin{align}
\bar{\partial} \chi_\mu =0,&\,\,\,\,\,\,\,\bar{\partial}\eta_i^{\chi_\mu} \left( T_{0i}^\alpha\right)=0 \,\,\,\,\,\left( \widetilde{\square}\left( \chi_\mu , \left\{ \eta_i^{\chi_\mu} \right\} \right)=0 \right) \\
-\bar{\partial} \sigma_\mu \left( T_{0i}^\alpha \right) + \left[ \chi_\mu , T_{0i}^\alpha \right] - \eta_i^{\chi_\mu}\left(T_{0i}^\alpha \right) &=  \Phi_{i|\mu}^\alpha - \bar{\partial} \Gamma_{i|\mu}^\alpha +\sum_{\xi=1}^p \Lambda_{i|\mu}^{\alpha\xi} T_{0i}^\xi +\left[ \varphi_\mu'', T_{0i}^\alpha \right]  - \sum_{\xi=1}^p b_{i|\mu}''^{\alpha  \xi} T_{0i}^\xi \label{sd32}  \\
\left[ T_{0i}^\alpha , \sigma_\mu\left( T_{0i}^\beta \right) \right]- \left[ T_{0i}^\beta, \sigma_\mu\left(T_{0i}^\alpha \right) \right] &- \sum_{\gamma=1}^p g_{0i\alpha\beta}^\gamma \sigma_\mu\left( T_{0i}^\gamma \right) + \pi_{i\sigma_\mu} \left( T_{0i}^\alpha \wedge T_{0i}^\beta \right) \\
= \Pi_{i|\mu}^{\alpha\beta} + \left[ \Gamma_{i|\mu}^\alpha, T_{0i}^\beta \right]& - \left[ \Gamma_{i|\mu}^\beta, T_{0i}^\alpha  \right] - \sum_{\gamma=1}^p g_{0i\alpha\beta}^\gamma \Gamma_{i|\mu}^\gamma + \sum_{\gamma=1}^p \Psi_{i|\mu}^{\alpha\beta \gamma} T_{0i}^\gamma   \notag\\
\pi_{j\sigma_\mu}\left( T_{0i}^\alpha \wedge T_{0i}^\beta \right) - \pi_{i\sigma_\mu}\left( T_{0i}^\alpha \wedge T_{0i}^\beta \right) &=\sum_{\eta=1}^p \left[ \sigma_{i|\mu}^\alpha, r_{0ij}^{\alpha \beta} \right] T_{0j}^\eta - \sum_{\eta=1}^p \left[ \sigma_{i|\mu}^\beta, r_{0ij}^{\alpha\eta} \right] T_{0j}^\eta \label{ds2}\\
\left| \left(\chi_\mu, \left\{ \eta_i^{\chi_\mu}\right\} \right) \right|_{k+\alpha} , \,\,\, &\left| \sigma_\mu\right|_{k+\alpha},\,\,\, \left| \pi_{i\sigma_\mu}\right|_{k-1+\alpha}  \ll K_{28} \frac{b}{c}A(t) \label{ds27}
\end{align}
In view of the equality
\begin{align*}
\chi_\mu=\bold{H} \chi_\mu + (\bar{\partial} \mathfrak{d}+ \mathfrak{d}\bar{\partial}) G \chi_\mu
\end{align*}
we get
\begin{align}
&-\bar{\partial} \sigma_\mu\left(T_{0i}^\alpha\right) +\left[\bold{H}\chi_\mu, T_{0i}^\alpha  \right] +\bar{\partial}\left[ \mathfrak{d}G\chi_\mu, T_{0i}^\alpha  \right] - \eta_i^{\chi_\mu} \left(T_{0i}^\alpha\right) = \Phi_{i|\mu}^\alpha - \bar{\partial} \Gamma_{i|\mu}^\alpha +\sum_{\xi=1}^p \Lambda_{i|\mu}^{\alpha\xi} T_{0i}^\xi +\left[ \varphi_\mu'', T_{0i}^\alpha \right]  - \sum_{\xi=1}^p b_{i|\mu}''^{\alpha  \xi} T_{0i}^\xi \notag \\
&\Longrightarrow -\bar{\partial}\left( \sigma_{i|\mu}^\alpha -\left[\mathfrak{d} G\chi_\mu , T_{0i}^\alpha \right]    \right) +\left[ \bold{H}\chi_\mu - \varphi_\mu'', T_{0i}^\alpha \right] - \sum_{\xi=1}^p \left(\eta_i^{\chi_\mu \alpha\xi}  - b_{i|\mu}''^{\alpha \xi}\right) T_{0i}^\xi= \Phi_{i|\mu}^\alpha - \bar{\partial} \Gamma_{i|\mu}^\alpha + \sum_{\xi=1}^p \Lambda_{i|\mu}^{\alpha \xi} T_{0i}^\xi \label{d42}
\end{align}
On the other hand, we define an element $S_{ij|\mu}\in \Gamma\left( U_{ij}, \mathcal{A}^{0,0}\left( \mathscr{H}om_{\mathcal{O}_M}\left( \Theta_{\mathcal{F}_0}, \Theta_{\mathcal{F}_0} \right) \right) \right)$ by 
\begin{align*}
S_{ij|\mu}:\Gamma \left(U_{ij}, \Theta_{\mathcal{F}_0 } \right)&\to \Gamma \left(U_{ij}, \mathcal{A}^{0,0}\left(\Theta_{\mathcal{F}_0} \right) \right)\\
T_{0i}^\alpha &\mapsto  \sum_{\beta=1}^p \left[ \mathfrak{d}G\chi_\mu, r_{0ij}^{\alpha\beta}\right] T_{0j}^\beta,
\end{align*}
Then $\left\{S_{ij|\mu} \right\}$ is a $1$-cocycle with coefficients in $\mathcal{A}^{0,0}\left( \mathscr{H}om_{\mathcal{O}_M}\left( \Theta_{\mathcal{F}_0} , \Theta_{\mathcal{F}_0} \right) \right)$. Let $\left\{ \rho_i(z) \right\}$ be a partition of unity subordinate to the covering $\mathcal{U}=\left\{ U_i \right\}$. By setting $S_{i|\mu}=\sum_k \rho_k(z)S_{ik|\mu}$, we see that  $\left\{S_{i|\mu} \right\}\in C^0\left( \mathcal{U}, \mathcal{A}^{0,0}\left( \mathscr{H}om_{\mathcal{O}_M}\left(  \Theta_{\mathcal{F}_0}, \Theta_{\mathcal{F}_0}  \right)   \right)  \right)$ such that $S_{i|\mu} -S_{j|\mu}=S_{ij|\mu}$ where
\begin{align*}
S_{i|\mu}:\Gamma\left(U_i, \Theta_{\mathcal{F}_0} \right) &\to \Gamma\left(U_i, \mathcal{A}^{0,0}\left(\Theta_{\mathcal{F}_0} \right)\right)\\
    T_{0i}^\alpha &\mapsto \sum_{\gamma=1}^p S_{i|\mu}^{\alpha\gamma} T_{0i}^\gamma, \,\,\,\,\,\,\,\,S_{i|\mu}^{\alpha\beta}\in \Gamma\left(U_i, \mathcal{A}^{0,0}\right)
\end{align*}
and we have
\begin{align}\label{ds26}
\sum_{\beta,\gamma=1}^p r_{0ij}^{\alpha\beta} S_{j|\mu}^{\beta \gamma} T_{0j}^\gamma - \sum_{\gamma=1}^p S_{i|\mu}^{\alpha \gamma} T_{0i}^\gamma = -\sum_{\beta=1}^p \left[ \mathfrak{d}G\chi_\mu , r_{0ij}^{\alpha\beta}  \right]T_{0j}^\beta
\end{align}
Then we have from $(\ref{ds27})$
\begin{align}\label{ds30}
\left| S_{i|\mu} \right|_{k+ 1+ \alpha} \ll K_{29} \left| \mathfrak{d} G \chi_\mu \right|_{k+1+\alpha} \ll K_{30} \left| \chi_\mu \right|_{k+\alpha} \ll K_{31} \frac{b}{c} A(t)
\end{align}
On the other hand, $(\ref{ds26})$ implies that
\begin{align}\label{d48}
P_i:\Gamma\left(U_i, \Theta_{\mathcal{F}_0} \right) &\to \Gamma\left( U_i, \mathcal{A}^{0,0} \left(\Theta_M \right) \right) \\
 T_{0i}^\alpha & \mapsto \left[ \mathfrak{d}G \chi , T_{0i}^\alpha \right] - S_{i|\mu}\left(T_{0i}^\alpha\right) \notag
\end{align}
defines a global section $P=\left\{P_i \right\}\in A^{0,0}\left(M, \mathscr{H}om_{\mathcal{O}_M}\left( \Theta_{\mathcal{F}_0} , \Theta_M \right)\right)$, so that $\sigma - P\in A^{0,0}\left(M, \mathscr{H}om_{\mathcal{O}_M}\left( \Theta_{\mathcal{F}_0} , \Theta_M  \right)\right)$. Then we have from $(\ref{d42})$
\begin{align} \label{d43}
-\bar{\partial}\left( \sigma_{i|\mu}^\alpha -\left[\mathfrak{d} G\chi_\mu , T_{0i}^\alpha \right]   + \sum_{\xi=1}^p S_{i|\mu}^{\alpha \xi } T_{0i}^\xi \right) +\left[ \bold{H}\chi - \varphi_\mu'', T_{0i}^\alpha \right] - \sum_{\xi=1}^p \left(\eta_i^{\chi \alpha\xi}  - b_i''^{\alpha \xi} - \bar{\partial} S_{i|\mu}^{\alpha \xi}\right) T_{0i}^\xi= \Phi_{i|\mu}^\alpha - \bar{\partial} \Gamma_{i|\mu}^\alpha + \sum_{\xi=1}^p \Lambda_{i|\mu}^{\alpha \xi} T_{0i}^\xi
\end{align}
We note that
\begin{align}\label{d45}
&\sum_{\beta, \xi=1}^p r_{0ij}^{\alpha \beta} \left( \eta_j^{\chi \beta \xi} - b_j''^{\beta \xi} - \bar{\partial} S_{j|\mu}^{\beta \xi} \right) T_{0j}^\xi- \sum_{\xi=1}^p\left( \eta_i^{\chi\alpha \xi} - b_i''^{\alpha \xi} - \bar{\partial} S_{i|\mu}^{\alpha \xi} \right) T_{0i}^\xi\\
&= -\sum_{\beta=1}^p \left[  \chi_\mu - \varphi_\mu'' - \bar{\partial} \mathfrak{d} G\chi_\mu  , r_{0ij}^{\alpha\beta} \right] T_{0j}^\beta = - \sum_{\beta=1}^p \left[ \bold{H}\chi_\mu - \varphi_\mu'' , r_{0ij}^{\alpha\beta} \right] T_{0j}^\beta \notag
\end{align}

If we set
\begin{align}
\varphi_\mu' &:= \bold{H} \chi_\mu - \varphi_\mu'' \label{ds32}\\
T_{i|\mu}'^\alpha &:= \sigma_{i|\mu}^\alpha -\left[\mathfrak{d} G\chi_\mu , T_{0i}^\alpha \right]   + \sum_{\xi=1}^p S_{i|\mu}^{\alpha \xi } T_{0i}^\xi \label{ds33}\\
b_{i|\mu}^{\alpha \xi}  &:= - \eta_i^{\chi \alpha \xi} + b_i''^{\alpha \xi} + \bar{\partial} S_{i|\mu}^{\alpha \xi} \label{ds34}
\end{align}
then we have from $(\ref{ds32})$ and $(\ref{ds31})$ and $(\ref{ds27})$
\begin{align}\label{ds35}
\left| \varphi_\mu'\right|_{k+\alpha} \ll K_{32}\frac{b}{c}A(t)
\end{align}
and we have from $(\ref{ds33})$ and $(\ref{ds30})$ and $(\ref{ds27})$
\begin{align}\label{ds36}
\left| T_{i|\mu}'^\alpha \right|_{k+\alpha} \ll K_{33} \frac{b}{c} A(t)
\end{align}
and we have from $(\ref{ds34})$ and $(\ref{ds27})$ and $(\ref{ds31})$
\begin{align}\label{ds37}
\left|b_{i|\mu}^{\alpha \xi} \right|_{k+\alpha} \ll K_{34}\frac{b}{c} A(t)
\end{align}

From $(\ref{d41})$ we have
\begin{align*}
-\bar{\partial}\varphi_\mu'=-\bar{\partial}\left(\bold{H} \chi_\mu  - \varphi_\mu'' \right) =  \bar{\partial} \varphi_\mu''= - \xi_\mu,\,\,\,\,\,\,\,\, \mathfrak{d} \varphi_\mu'=\mathfrak{d}\left(\bold{H}\chi_\mu - \varphi_\mu'' \right)= -\mathfrak{d}\varphi_\mu''=0
\end{align*}
and so from $(\ref{te33})$, we have $\varphi_\mu:= -\varphi_\mu'$, so that
\begin{align}
\bar{\partial}\varphi_\mu=-\xi_\mu,\,\,\,\,\,\,\,\,\mathfrak{d}\varphi_\mu=0
\end{align}
and $(\ref{d43})$ and $(\ref{d45})$ implies $(\ref{t50})$ and $(\ref{d46})$.
Since $\bar{\partial}\left( b_{i|\mu}^{\alpha \xi} - \Lambda_{i|\mu}^{\alpha \xi}\right)=0$ from $(\ref{d41})$, there exists $c_{i|\mu}^{\alpha \xi}\in \Gamma\left( U_i, \mathcal{A}^{0,0} \right)$ such that $\bar{\partial } c_{i|\mu}^{\alpha \xi}= b_{i|\mu}^{\alpha \xi} - \Lambda_{i|\mu}^\xi$. Then by Potential theoretic Lemma (see \cite{New57}), we have
\begin{align}\label{ds38}
\left| c_{i|\mu}^{\alpha \xi} \right|_{k+1+\alpha}^{U_i} \ll K_{35} \left| b_{i|\mu}^{\alpha \xi} -  \Lambda_{i|\mu}^{\alpha\xi} \right|_{k+\alpha}^{U_i}
\end{align}
and from $(\ref{ds37})$ and $(\ref{ds18})$ and $(\ref{ds38})$, we have
\begin{align}\label{ds45}
\left| c_{i|\mu}^{\alpha \xi} \right|_{k+1+\alpha} \ll K_{36} \frac{b}{c} A(t)
\end{align}

On the other hand, we note that from $(\ref{d48})$
\begin{align}\label{ds1}
&\left[ T_{0i}^\alpha,  P\left(T_{0i}^\beta\right)\right] - \left[ T_{0i}^\beta, P \left(T_{0i}^\alpha\right) \right] - \sum_{\gamma=1}^p g_{0i\alpha\beta}^\gamma P \left(T_{0i}^\gamma\right)\\
&=\left[ T_{0i}^\alpha, \left[ \mathfrak{d}G \chi_\mu , T_{0i}^\beta\right] \right] - \left[ T_{0i}^\alpha, S_{i|\mu}^\beta\right] -\left[ T_{0i}^\beta, \left[ \mathfrak{d}G \chi_\mu , T_{0i}^\alpha \right]  \right] +\left[ T_{0i}^\beta, S_{i|\mu}^\alpha \right] - \sum_{\gamma=1}^p g_{0i\alpha\beta}^\gamma \left[ \mathfrak{d} G\chi_\mu , T_{0i}^\gamma \right] + \sum_{\gamma=1}^p g_{0i\alpha\beta}^\gamma S_{i|\mu}^\gamma \notag \\
&= - \left[ \left[ T_{0i}^\alpha, T_{0i}^\beta\right], \mathfrak{d}G\chi_\mu \right] - \sum_{\gamma=1}^p g_{0i\alpha\beta}^\gamma \left[\mathfrak{d}G\chi_\mu , T_{0i}^\gamma   \right] - \left[ T_{0i}^\alpha, S_{i|\mu}^\beta \right] + \left[ T_{0i}^\beta, S_{i|\mu}^\alpha\right] + \sum_{\gamma=1}^p g_{0i\alpha\beta}^\gamma S_{i|\mu}^\gamma \notag \\
&= \sum_{\gamma=1}^p \left[  \mathfrak{d}G \chi_\mu , g_{0i\alpha\beta}^\gamma \right] T_{0i}^\gamma - \left[ T_{0i}^\alpha, S_{i|\mu}^\beta \right] + \left[ T_{0i}^\beta, S_{i|\mu}^\alpha\right] + \sum_{\gamma=1}^p g_{0i\alpha\beta}^\gamma S_{i|\mu}^\gamma \notag
\end{align}
We define
\begin{align*}
\pi_{iP}:\Gamma\left(U_i,  \bigwedge^2 \Theta_{\mathcal{F}_0} \right)&\to \Gamma\left(U_i,  \mathcal{A}^{0,0}\left( \Theta_{\mathcal{F}_0} \right) \right) \\
T_{0i}^\alpha &\mapsto \pi_{iP}^{\alpha\beta}:=\sum_{\gamma=1}^p \left[  \mathfrak{d}G \chi_\mu , g_{0i\alpha\beta}^\gamma \right] T_{0i}^\gamma - \left[ T_{0i}^\alpha, S_{i|\mu}^\beta \right] + \left[ T_{0i}^\beta, S_{i|\mu}^\alpha\right] + \sum_{\gamma=1}^p g_{0i\alpha\beta}^\gamma S_{i|\mu}^\gamma
\end{align*}
Then we have from $(\ref{ds27})$ and $(\ref{ds30})$
\begin{align}\label{ds28}
\left| \pi_{iP}\right|_{k+\alpha} \ll  K_{37} \frac{b}{c} A(t)
\end{align}
If we take
\begin{align*}
W_{i|\mu}^{\alpha\beta \xi}:= \pi_{i \sigma }^{\alpha\beta \xi}  - \pi_{i P}^{\alpha\beta \xi}
\end{align*}
Then $(\ref{d49})$ is satisfied from $(\ref{ds1})$ and $(\ref{ds2})$, and we have from $(\ref{ds27})$ and $(\ref{ds28})$
\begin{align}\label{ds46}
\left| W_{i|\mu}^{\alpha\beta \xi} \right|_{k-1+\alpha} =\left|\pi_{i\sigma}^{\alpha\beta \xi} - \pi_{iP}^{\alpha \beta\xi} \right|_{k-1+\alpha} \ll K_{38} \frac{b}{c} A(t)
\end{align}

Then we have from $(\ref{te33})$ and $(\ref{te37})$ and $(\ref{te41})$ and $(\ref{te42})$
\begin{align}
\varphi_\mu&= - \varphi_\mu'= \varphi_\mu''- \bold{H} \chi_\mu \label{ds40}\\
T_{i|\mu}^\alpha&= \Gamma_{i|\mu}^\alpha- T_{i|\mu}'^\alpha + \sum_{\xi=1}^p c_{i|\mu}^{\alpha \xi} T_{0i}^\xi \label{ds41}\\
r_{ij|\mu}^{\alpha\beta}&= \sum_{\eta=1}^p \lambda_{ij|\mu}^{\alpha \eta} r_{0ij}^{\eta \xi} + \sum_{\eta=1}^p c_{i|\mu}^{\alpha \eta} r_{0ij}^{\eta \xi} - \sum_{\beta=1}^p r_{0ij}^{\alpha\beta} c_{j|\mu}^{\beta \xi}\label{ds42}\\
g_{i\alpha\beta |\mu}^\xi&= \sum_{\gamma=1}^p c_{i|\mu}^{\alpha\gamma} g_{0i\gamma \beta}^\xi - T_{0i}^\beta\left(  c_{i|\mu}^{\alpha \xi} \right) + T_{0i}^\alpha \left( c_{i|\mu}^{\beta\xi} \right) + \sum_{\gamma=1}^p c_{i|\mu}^{\beta\gamma} g_{0i\alpha\gamma}^\xi - \sum_{\gamma=1}^p g_{0i\alpha\beta}^\gamma c_{i|\mu}^{\gamma \xi} + W_{i|\mu}^{\alpha\beta\xi}- \Psi_{i|\mu}^{\alpha\beta \xi} \label{ds43} 
\end{align}
Then from $(\ref{ds40})$ and $(\ref{ds35})$
\begin{align*}
\left| \varphi_\mu \right|_{k+\alpha} \ll K_{32}\frac{b}{c} A(t)
\end{align*}
and from $(\ref{ds41})$ and $(\ref{ds24})$ and $(\ref{ds36})$ and $(\ref{ds45})$ we have
\begin{align*}
\left| T_{i|\mu}^\alpha \right|_{k+\alpha} \ll K_{39} \frac{b}{c} A(t)
\end{align*}
and from $(\ref{ds42})$ and $(\ref{ds12})$ and $(\ref{ds45})$
\begin{align*}
\left| r_{ij|\mu}^{\alpha\beta} \right|_{k+1+\alpha} \ll K_{40} \frac{b}{c} A(t)
\end{align*}
and from $(\ref{ds43})$ and $(\ref{ds45})$ and $(\ref{ds46})$ and $(\ref{ds25})$
\begin{align*}
\left| g_{i\alpha\beta |\mu}^\xi \right|_{k-1+\alpha} \ll K_{41} \frac{b}{c} A(t)
\end{align*}
Then we can choose $b$ and $c$ satisfying $(\ref{d20})_\mu-(\ref{d21})_\mu$. Then $\varphi(t)$ and $T_i^\alpha$ and  $r_{ij}^{\alpha\beta}$ and $g_{i\alpha\beta}^\gamma$ converge with respect to the \"Holder norm $\left| - \right|_{k+\alpha}, \left|- \right|_{k+\alpha},\left|- \right|_{k+1+\alpha}, \left|- \right|_{k-1+\alpha}$, respectively for some sufficiently small neighborhood $\Delta_\epsilon \subset B$ of $0$. Consequently $\varphi(t)$ is a $C^k$ vector $(0,1)$-form on $M\times \Delta_\epsilon$ and $T_i^\alpha$ is a $C^k$-vector of the form $T_i^\alpha= \sum_{\beta=1}^n T_i^{\alpha\beta}\left(z_i,t \right)\frac{\partial}{\partial z_i^\beta}$ on $U_i\times \Delta_\epsilon$, and $r_{ij}^{\alpha\beta}$ is a $ C^{k+1}$ function on $U_{ij}\times \Delta_\epsilon$ and $g_{i\alpha\beta}^\gamma$ is a $C^{k-1 }$ function on $U_i\times \Delta_\epsilon$. From $(\ref{d41})$ and $(\ref{ds40})$, we have $\mathfrak{d}\varphi_\mu =0$ for $\mu\geq 2$. Hence $\varphi(t)$ is a solution of the quasi-linear partial differential equation of order $2$ 
\begin{align}\label{ds48}
\sum_{\lambda=1}^r \frac{\partial^2}{\partial t_\lambda \partial \bar{t}_\lambda} \varphi(t) +\square \varphi(t) - \mathfrak{d}\left[ \varphi(t), \varphi(t)\right]=\bar{\partial}\mathfrak{d}\varphi_1(t)
\end{align}
We may assume that $(\ref{ds48})$ is a quasi-linear elliptic partial differential equation on $M\times \Delta_\epsilon$ (see \cite{Kod05} p.281). Therefore its solution $\varphi(t)$ is $C^\infty$ on $M\times \Delta_\epsilon$ (see \cite{DN55} Theorem 5 or \cite{Kod05} Appendix \S 8). Then $\varphi(t)$ determines a complex analytic family $\pi: \mathcal{M}\to \Delta_\epsilon$ (see \cite{Kod05} p.281-283). On the other hand, $T_i^\alpha$ and $r_{ij}^{\alpha\beta}$ and $g_{i\alpha\beta}^\gamma$ are holomorphic with respect to the complex structure $\varphi(t)$ by $(\ref{ua1})$ and $(\ref{ua2})$ and $(\ref{ua3})$, respectively. $(\ref{ua4}),(\ref{ua5})$ and $(\ref{t11})$ determines a locally free subsheaf $\Theta_\mathcal{F}$ of $\Theta_{\mathcal{M}/\Delta_\epsilon}$ satisfying the integrability condition. From the definition of linear terms, we infer that the foliated Kodaira-Spencer map $\varphi_0: T_0\left(\Delta_\epsilon \right) \to \mathbb{H}^1\left( M, \Theta_{\mathcal{F}_0}^\bullet \right)$ is bijective. This completes the proof of Theorem \ref{tt2}.

\end{proof}

\section{Theorem of stability for deformations of foliated complex analytic structures in terms of tangent sheaves}

Our proof of theorem of existence of deformations of foliated complex analytic structures in terms of tangent sheaves (Theorem \ref{tt2}) provides the proof of theorem of stability for deformations of foliated complex analytic structures in terms of tangent sheaves (Theorem \ref{stt1}). We define a complex of sheaves which we truncate the $0$-th term of the leaf complex
{\Small{\begin{align*}
\mathscr{H}om_{\mathcal{O}_M}\left( \Theta_{\mathcal{F}_0}, \frac{\Theta_M}{\Theta_{\mathcal{F}_0}} \right)^\bullet : \mathscr{H}om_{\mathcal{O}_M}\left( \Theta_{\mathcal{F}_0}, \frac{\Theta_M}{\Theta_{\mathcal{F}_0}} \right)\xrightarrow{D_1} \mathscr{H}om_{\mathcal{O}_M}\left( \bigwedge^2 \Theta_{\mathcal{F}_0}, \frac{\Theta_M}{\Theta_{\mathcal{F}_0}} \right) \xrightarrow{D_2} \mathscr{H}om_{\mathcal{O}_M}\left( \bigwedge^3  \Theta_{\mathcal{F}_0}, \frac{\Theta_M}{\Theta_{\mathcal{F}_0}} \right) \xrightarrow{D_3} \cdots
\end{align*}}}
We will denote the $i$-th cohomology group of $\mathscr{H}om_{\mathcal{O}_M}\left( \Theta_{\mathcal{F}_0}, \frac{\Theta_M}{\Theta_{\mathcal{F}_0}} \right)^\bullet$ by $\mathbb{H}^i\left( M, \mathscr{H}om_{\mathcal{O}_M}\left( \Theta_{\mathcal{F}_0}, \frac{\Theta_M}{\Theta_{\mathcal{F}_0}} \right)^\bullet \right)$.

\begin{theorem}[Theorem of stability for deformations of foliated complex analytic structures in terms of tangent sheaves] \label{stt1}
Let $\left( M, \Theta_{\mathcal{F}_0} \right)$ be a compact foliated complex manifold with $\Theta_{\mathcal{F}_0}$ locally free. Assume that $\mathbb{H}^1\left( M, \mathscr{H}om_{\mathcal{O}_M}\left( \Theta_{\mathcal{F}_0}, \frac{\Theta_M}{\Theta_{\mathcal{F}_0}} \right)^\bullet \right)=0$. Then for any complex analytic family $\pi:\mathcal{M}\to B$ of deformations of $\pi^{-1}(0)=M, 0\in B$ in the sense of Kodaira-Spencer, there exists an open neighborhood $N\subset B$ of $0$ and a locally free subsheaf $\Theta_{\mathcal{F}}$ of $\Theta_{\frac{\mathcal{M}|_N}{N}}$ such that $\left(\mathcal{M}|_N, \Theta_{\mathcal{F}}\right)$ defines a foliated analytic family $\pi|_N: \left(\mathcal{M}|_N , \Theta_{\mathcal{F}} \right) \to N$ of deformations of $\left(M, \Theta_{\mathcal{F}_0} \right)=\pi^{-1}(0)$ in terms of tangent sheaves.
\end{theorem}

\begin{proof}
We copy the proof of Theorem \ref{tt2}. The difference is that we assume that $\varphi(t)$ (determined by $\mathcal{M}$ over some neighborhood $N'\subset B$ of $0$) is given from the beginning and replace $\varphi^\mu$ by $\varphi(t)$ in the proof of Theorem \ref{tt2}. Accordingly we replace $\xi_\mu$ and $\varphi_\mu$ by $0$. Then the obstructions for solving $(\ref{te12})-(\ref{t33})$ are in $\mathbb{H}^1\left( M, \mathscr{H}om_{\mathcal{O}_M}\left( \Theta_{\mathcal{F}_0} ,\frac{\Theta_M}{\Theta_{\mathcal{F}_0}} \right)^\bullet \right)=0$. Hence we can construct $\Theta_{\mathcal{F}}$ on $\mathcal{M}|_N$ for some neighborhood $N\subset N'$ of $0$. 
\end{proof}

\section{Theorem of completeness of deformations of foliated complex analytic structures in terms of tangent sheaves}

\subsection{Change of parameters (\textnormal{compare \cite{Kod05} p.205})}\label{r5} \

Consider a foliated complex analytic family $\left( \mathcal{M}, \Theta_{\mathcal{F}}, B, \pi \right)$ with $\Theta_{\mathcal{F}}$ locally free of deformations of $\left( M_t, \Theta_{\mathcal{F}_t} \right)=\pi^{-1}(t), t\in B$ where $B$ is a domain of $\mathbb{C}^m$. Let $D$ be a domain of $\mathbb{C}^{m'}$ and $s:u\to t=s(u), u\in D$, a holomorphic map of $D$ into $B$. Then by changing the parameter from $t$ to $u$, we will construct a foliated analytic family $\left\{ \left( M_{s(u)}, \Theta_{\mathcal{F}_{s(u)}} \right) | u \in D \right\}$ on the parameter space $D$ in the following.

Let $\mathcal{M} \times_B D:=\left\{ (p,u)\in \mathcal{M}\times D | \pi(p)=s(u) \right\}$. Then we have the following commutative diagram
\begin{center}
$\begin{CD}
\mathcal{M}\times_B D @>p>> \mathcal{M}\\
@V\pi' VV @VV\pi V \\
D @>s>> B
\end{CD}$
\end{center}
such that $(\mathcal{M}\times_B D, D, \pi')$ is a complex analytic family in the sense of Kodaira-Spencer and $\pi'^{-1}(u)=M_{s(u)}$. On the other hand, $\Pi^*\Theta_{\mathcal{F}}$ is a locally free subsheaf of $\Theta_{\mathcal{M}\times D/D}$, where $\Pi:\mathcal{M}\times D \to \mathcal{M}$ is the natural projection. Since $\mathcal{M}\times_B D$ is a complex submanifold of $\mathcal{M}\times D$, we have an injection $0\to \Theta_{\mathcal{M}\times_B D/D}\hookrightarrow \Theta_{\mathcal{M}\times D/D}|_{\mathcal{M}\times_B D}$ and $\Pi^* \Theta_\mathcal{F}$ is invariant with respect to $\mathcal{M}\times_B D$, the restriction $\Pi^* \Theta_\mathcal{F}|_{\mathcal{M}\times_B D} =p^* \Theta_\mathcal{F}$ is a locally free  subsheaf of $\Theta_{\mathcal{M}\times_B D/D}$, and $\left[ p^* \Theta_\mathcal{F}, p^* \Theta_\mathcal{F} \right] \subset p^* \Theta_\mathcal{F}$. Moreover for each $u\in D$, we have an injection $p^*\Theta_{\mathcal{F}}|_{\pi'^{-1}(u)} = \Theta_{\mathcal{F}_{s(u)}} \hookrightarrow \Theta_{\mathcal{M}\times_B D/D}|_{\pi'^{-1}(u)}= \Theta_{M_{s(u)}}$, so that $\frac{\Theta_{\mathcal{M}\times_B D/D}}{p^* \Theta_{\mathcal{F}}}$ is flat over $D$. This implies that $\left\{ \left( M_{s(u)}, \Theta_{\mathcal{F}_{h(u)}} \right)| u\in D \right\}$ forms a foliated complex analytic family $\left( \mathcal{M}\times_B D, p^* \Theta_\mathcal{F}, D, \pi' \right)$.

\begin{definition}
The foliated complex analytic family $\left( \mathcal{M}\times_B D, p^* \Theta_{\mathcal{F}}, \pi' \right)$ is called the foliated complex analytic family induced from $\left(\mathcal{M}, \Theta_\mathcal{F}, B, \pi \right)$ by the holomorphic map $s:D\to B$.
\end{definition}

To investigate the relation of the infinitesimal foliated deformation of $\left( M_{s(u)}, \Theta_{\mathcal{F}_{s(u)}} \right)$ and that of $\left( M_t, \Theta_{\mathcal{F}_t} \right)$, we assume that $0 \in B, 0 \in D$ and $s(0)=0$. Taking a sufficiently small coordinate polydisk $\Delta$ with $0 \in \Delta \subset B$, we represent $\left(\mathcal{M}_\Delta, \Theta_{\mathcal{F}_\Delta} \right)=\pi^{-1}\left( \Delta \right)$ in the form
\begin{align*}
\left( \mathcal{M}_\Delta, \Theta_{\mathcal{F}_\Delta} \right)= \bigcup \left( U_j\times \Delta, \Theta_{\mathcal{F}}|_{U_j\times \Delta}   \right)
\end{align*}
where $\left(z_j,t \right)\in U_j\times \Delta$ and $\left(z_k,t \right)\in U_k\times \Delta$ are the same point of $\mathcal{M}_\Delta$ if $z_j= f_{jk}\left(z_k,t \right)$, and $\Gamma\left(U_j\times \Delta, \Theta_{\mathcal{F}} \right)$ is generated by $T_j^1\left(z_j,t\right),..., T_j^p \left(z_j,t \right)$, where $T_j^\alpha \left(z_j,t \right)=\sum_{\beta=1}^p T_j^{\alpha\beta}\left(z_j,t \right)\frac{\partial}{\partial z_j^\beta}$ for some $T_j^{\alpha\beta}\left(z_j,t \right)\in \Gamma\left(U_j\times \Delta, \mathcal{O}_\mathcal{M} \right)$ with $T_j^\alpha\left(z_j,t \right)= \sum_{\beta=1}^p r_{jk}^{\alpha\beta}\left(z_k,t \right)T_k^\beta \left(z_k,t \right)$ for some $r_{jk}^{\alpha\beta}\left(z_k,t \right)\in \Gamma\left( \left(U_j\times \Delta \right) \bigcap \left(U_k\times \Delta\right), \mathcal{O}_\mathcal{M} \right)$ for $\alpha=1,...,p$, and $\left[ T_j^\alpha\left(z_j,t\right), T_j^\beta\left(z_j,t\right)\right]= \sum_{\gamma=1}^p g_{j\alpha\beta}^\gamma \left(z_j,t \right) T_j^\gamma\left(z_j,t \right)$ for some $g_{j\alpha\beta}^\gamma\left(z_j,t \right)\in \Gamma \left( U_j\times \Delta, \mathcal{O}_\mathcal{M} \right)$ for $\alpha,\beta=1,...,p$. Take a polydisk $\Delta'$ with $0\in \Delta' \subset D$. Then $\left( \mathcal{M}\times_B D|_{\Delta'}, p^* \Theta_\mathcal{F}|_{\Delta'} \right)$ is represented by the form
\begin{align*}
\left( \mathcal{M}\times_B D|_{\Delta'}, p^* \Theta_\mathcal{F}|_{\Delta'} \right) =\left( U_j\times \Delta', p^* \Theta_\mathcal{F}|_{U_j\times \Delta'} \right)
\end{align*}
where $\left(z_j,u \right)\in U_j\times \Delta'$ and $\left(z_k,u\right)\in U_k\times \Delta'$ are the same point of $\mathcal{M}\times_B D|_{\Delta'}$ if $z_j= f_{jk}\left(z_k, s(u)\right)$, and $\Gamma\left( U_j\times \Delta', p^* \Theta_\mathcal{F}\right)$ is generated by $T_j^1\left(z_j, s(u)\right),..., T_j^p\left(z_j, s(u)\right)$, where $T_j^\alpha \left(z_j, s(u) \right)=\sum_{\beta=1}^p T_j^{\alpha\beta}\left(z_j, s(u)\right)\frac{\partial}{\partial z_j^\beta}$, and $T_j^\alpha \left(z_j, s(u) \right)= \sum_{\beta=1}^p r_{jk}^{\alpha\beta}\left(z_j, s(u)\right) T_j^\beta\left(z_k, s(u)\right)$ for $\alpha=1,...,p$, and $\left[ T_j^\alpha\left(z_j, u(s) \right), T_j^\beta\left(z_j, u(s) \right) \right]= \sum_{\gamma=1}^p g_{j\alpha\beta}^\gamma \left(z_j, u(s) \right) T_j^\gamma \left(z_j, u(s) \right)$. Then we can show that
\begin{theorem}\label{r8}
For any tangent vector $\frac{\partial }{\partial u}=c_1\frac{\partial}{\partial u_1} + \cdots + c_{m'}\frac{\partial}{\partial u_{m'}}\in T_u \left(D \right) $, the infinitesimal foliated deformation of $\left(M_{s(u)}, \Theta_{\mathcal{F}_{s(u)}} \right)$ along $\frac{\partial}{\partial u}$ is given by
\begin{align*}
\frac{\partial \left(M_{s(u)}, \Theta_{\mathcal{F}_{s(u)}}\right)}{\partial u}&=\left(\sum_{\gamma=1}^m \frac{\partial t_\gamma}{\partial u}\frac{\partial M_t}{\partial t_\gamma},\left\{T_j^\alpha(z_j,s(u)) \mapsto - \overline{\sum_{\gamma=1}^m \frac{\partial t_\gamma}{\partial u} \frac{\partial T_j^\alpha}{\partial t_\gamma}} \right\} \right)
            =\sum_{\gamma=1}^m \frac{\partial t_\gamma}{\partial u} \frac{\partial \left(M_t, \Theta_{\mathcal{F}_t} \right) }{\partial t_\gamma}
\end{align*}
\end{theorem}

\subsection{Theorem of completeness of deformations of foliated complex analytic structures in terms of tangent sheaves} \

\begin{definition}
Let $\left( \mathcal{M}, \Theta_\mathcal{F}, B, \pi \right)$ with $\Theta_{\mathcal{F}}$ locally free be a foliated analytic family of compact foliated complex manifolds in terms of tangent sheaves, and $t^0 \in B$. Then $\left( \mathcal{M}, \Theta_\mathcal{F}, B, \pi\right)$ is called complete at $t^0\in B$ if for any foliated complex analytic family $\left(\mathcal{M}' , \Theta_{\mathcal{F}'}, D, \pi' \right)$ in terms of tangent sheaves such that $D$ is a domain of $\mathbb{C}^{m'}$ containing $0$ and that $\pi'^{-1}(0)= \pi^{-1}(t^0)=\left(M, \Theta_{\mathcal{F}_0} \right)$, there are a sufficiently small domain $\Delta$ with $0 \in \Delta \subset D$, and a holomorphic map $s:u\to t=s(u)$ with $s(0)=t^0$ such that $\left( \mathcal{M}'_\Delta, \Theta_{\mathcal{F}'_\Delta}, \Delta, \pi' \right)$ is the foliated complex analytic family induced from $\left( \mathcal{M}, \Theta_\mathcal{F}, B, \pi \right)$ by $s$ where $\left( \mathcal{M}'_\Delta, \Theta_{\mathcal{F}'_\Delta}, \Delta, \pi' \right)$ is the restriction of $\left( \mathcal{M}' , \Theta_{\mathcal{F}'}, D, \pi' \right)$ to $\Delta$.
\end{definition}

We shall prove the following theorem 
\begin{theorem}\label{tc1}
Let $\left(\mathcal{M}, \Theta_\mathcal{F}, B, \pi \right)$ be a foliated analytic family of deformations of a complex foliated complex manifold $(M , \Theta_{\mathcal{F}_0})=\omega^{-1}(0)$ such that $\Theta_{\mathcal{F}_0}$ is locally free, and $B$ is a domain of $\mathbb{C}^r$ containing $0$. If the foliated Kodaira-Spence map $\varphi_0: T_0(B)\to \mathbb{H}^1(M, \Theta_{\mathcal{F}_0}^\bullet)$ is surjective, the foliated analytic family $\left(\mathcal{M}, \Theta_\mathcal{F}, B, \pi \right)$ in terms of tangent sheaves is complete at  $0 \in B$.
\end{theorem}

\begin{proof}
Let $\left(\mathcal{M}, \Theta_\mathcal{F}, B, \pi' \right)$ be a foliated analytic family in terms of tangent sheaves which is represented as in Remark \ref{ta1}.  We keep the notations in Remark \ref{ta1}. Since the problem is local with respect to $B$, we may assume that $B=\left\{ t\in \mathbb{C}^m | |t| <1\right\}$ is a polydisk, and $\mathcal{M}$ is written in the form
\begin{align}
\mathcal{M}=\bigcup_j \,\,\, \mathcal{U}_j,\,\,\,\,\,\,\,\,\mathcal{U}_j=\left\{\left( z_j, t\right) \in \mathbb{C}^n \times B\,\,\, |\,\,\, |z_j| <1 \right\}
\end{align}

Set $U_i:=M \cap \mathcal{U}_i=\left\{ z_i\in \mathbb{C}^n | |z_i|<1 \right\}$. Then $\left(M, \Theta_{\mathcal{F}_0} \right)=\pi^{-1}(0)$ is described in terms of the open covering $\mathcal{U} :=\{U_i\}$ in the following way:  (1) local coordinates $z_i= \left(z_i^1,..., z_i^n \right)$ on $U_i$ with $z_i= f_{0ij}(z_j)= \left( f_{0ij}^1(z_j),..., f_{0ij}^n(z_j) \right)$ where $f_{0ij}(z_j):= f_{ij}(z_j,0)$ and (2) $\Gamma\left(U_i, \Theta_{\mathcal{F}_0} \right)$ is generated by $T_{0i}^\alpha(z_i):= \sum_{\beta=1}^n T_i^{\alpha\beta}(z_i,0)\frac{\partial}{\partial z_i^\beta} \in \Gamma\left(U_i, \Theta_{\mathcal{F}_0} \right)$ with $T_{0i}^\alpha=\sum_{\beta=1}^p r_{0ij}^{\alpha\beta}(z_j) T_{0j}^\beta$ and $\left[ T_{0i}^\alpha, T_{0i}^\beta \right]=\sum_{\gamma=1}^p g_{0i\alpha\beta}^\gamma(z_i) T_{0i}^\gamma$ where $r_{0ij}^{\alpha\beta}(z_j) := r_{ij}^{\alpha\beta}(z_j,0)$ and $g_{0i\alpha\beta}^{\gamma}(z_i) : = g_{i\alpha\beta}^\gamma(z_i,0)$.

Let $\left(\mathcal{M}', \Theta_{\mathcal{F}'}, D, \pi' \right)$ be an another foliated analytic family such that $\pi'^{-1}(0')=\left(M, \Theta_{\mathcal{F}_0}\right)$. We may assume the following:
\begin{enumerate}
\item $D\subset \left\{ u\in \mathbb{C}^{m'} ||u|<1 \right\}$ is a sufficiently small polydisk in $\mathbb{C}^{m'}$ with a system of coordinates $u=(u_1,..., u_{m'})$ centered at $0'$. \label{sc2}
\item $\mathcal{M}'$ is covered by a finite number of coordinate neighborhood $\mathcal{U}_i'=\left\{ \left(\xi_i, u\right) \in \mathbb{C}^n \times D\,\,\,|\,\,\,|\xi_i|<1\right\}$ with a system of coordinate $(\xi_i, u)$ such that $\pi'(\xi_i, u)=u$.
\item $(\xi_i, u)$ coincides with $(\xi_j, u)$ if and only if $\xi_i = f_{ij}'(\xi_i, u)$.
\item $\Theta_{\mathcal{F}'}$ on $\mathcal{U}_i'$ is generated by
\begin{align*}
T_i'^\alpha(\xi_i, u):=\sum_{\beta=1}^n T_i'^{\alpha\beta}(\xi_i, u) \frac{\partial}{\partial \xi_i^\beta},\,\,\,\,\,\,\,\alpha=1,..., p.
\end{align*}
\item On $\mathcal{U}_i'\cap \mathcal{U}_j'\ne \emptyset$, we have $T_i'^{\alpha}(\xi_i, u)= \sum_{\beta=1}^q r_{ij}'^{\alpha\beta}(\xi_j, u) T_j'^\beta(\xi_j, u)$ for $r_{ij}'^{\alpha\beta}(\xi_j, u)\in \Gamma \left(\mathcal{U}_i'\cap \mathcal{U}_j', \mathcal{O}_{\mathcal{M}'} \right)$.
\item $\left[T_i'^\alpha(\xi_i, u), T_i'^\beta(\xi_i, u)    \right]=\sum_{\gamma=1}^p g_{i\alpha\beta}'^\gamma (\xi_i, u) T_i^\gamma(\xi_i, u)$ for $g_{i\alpha\beta}'^\gamma(\xi_i, u)\in \Gamma\left(\mathcal{U}_i', \mathcal{O}_{\mathcal{M}'} \right)$, $\alpha, \beta=1,..., p$.
\item $\pi'^{-1}(0')\cap \mathcal{U}_i'=U_i$, and $\xi_i=z_i$ on $U_i$ and $f_{0ij}(z_j)=f_{0ij}'(\xi_j)$ where $f_{0ij}'(\xi_j) : =f_{ij}'(\xi_j,0)$.
\item Setting $T_{0i}'^{\alpha}(z_i):=T_{i}'^\alpha (z_i, 0)$, i.e.
\begin{align*}
T_{0i}'^\alpha(z_i):= T_i'^\alpha(z_i,0)=\sum_{\beta=1}^n T_i'^{\alpha\beta}(z_i, 0)\frac{\partial}{\partial z_i^\beta},\,\,\,\,\,\,\,\alpha=1,...,p,
\end{align*}
we can find $b_{0i}^{\alpha\beta}(z_i)\in \Gamma\left(U_i, \mathcal{O}_M\right), \alpha,\beta=1,...,p$ such that
\begin{align*}
T_{0i}^\alpha(z_i)= \sum_{\beta=1}^p b_{0i}^{\alpha\beta}(z_i) T_{0i}'^\beta(z_i)
\end{align*}
\item we set $r_{0ij}'^{\alpha\beta}(\xi_j):=r_{ij}'^{\alpha\beta}(\xi_j,0)$ and $g_{0i\alpha\beta}'^\gamma(\xi_i):=g_{i\alpha\beta}^\gamma(\xi_i,0)$. \label{sc3}
\end{enumerate}
We note that since $T_{0i}^\alpha=\sum_{\beta=1}^q r_{0ij}^{\alpha\beta} T_{0j}^\beta$ and $T_{0i}'^{\alpha}=\sum_{\beta=1}^q r_{0ij}'^{\alpha\beta} T_{0j}'^\beta$, we have
\begin{align*}
\sum_{\beta,\gamma=1}^p b_{0i}^{\alpha\beta} r_{0ij}'^{\beta\gamma} T_{0j}'^\gamma = \sum_{\beta=1}^p b_{0i}^{\alpha\beta} T_{0i}'^\beta= \sum_{\beta=1}^p r_{0ij}^{\alpha\beta} T_{0j}^\beta= \sum_{\beta,\gamma=1}^p r_{0ij}^{\alpha\beta} b_{0j}^{\beta \gamma} T_{0j}'^\gamma \Longrightarrow \sum_{\beta=1}^p b_{0i}^{\alpha\beta} r_{0ij}'^{\beta \gamma} = \sum_{\beta=1}^p r_{0ij}^{\alpha\beta} b_{0j}^{\beta \gamma}
\end{align*}
We note that
\begin{align*}
\sum_{\eta=1}^p g_{0i\alpha\beta}^\eta  b_{0i}^{\eta \gamma} T_{0i}'^\gamma &= \sum_{\eta=1}^p g_{0i\alpha\beta}^\eta T_{0i}^\eta =\left[ T_{0i}^\alpha, T_{0i}^\beta \right]= \sum_{\eta, \delta=1}^p \left[ b_{0i}^{\alpha \eta} T_{0i}'^\eta, b_{0i}^{\beta \delta} T_{0i}'^\delta \right] \\
&= \sum_{\eta, \delta, \gamma=1}^p b_{0i}^{\alpha \eta} b_{0i}^{\beta \delta} g_{0i\eta \delta}'^\gamma T_{0i}'^\gamma + \sum_{\eta, \delta=1}^p b_{0i}^{\alpha \eta} T_{0i}'^\eta\left(b_{0i}^{\beta \delta} \right) T_{0i}^\delta - \sum_{\eta, \delta=1}^p b_{0i}^{\beta \delta}  T_{0i}'^\delta\left( b_{0i}^{\alpha \eta} \right) T_{0i}'^\eta
\end{align*}
Then we have
\begin{align*}
\sum_{\eta=1}^p g_{0i\alpha\beta}^\eta b_{0i}^{\eta \gamma} = \sum_{\eta, \delta=1}^p b_{0i}^{\alpha \eta} b_{0i}^{\beta \delta} g_{0i\eta \delta}'^\gamma + \sum_{\eta=1}^p b_{0i}^{\alpha \eta} T_{0i}'^\eta\left( b_{0i}^{\beta \gamma} \right) - \sum_{\delta=1}^p b_{0i}^{\beta \delta} T_{0i}'^\delta \left( b_{0i}^{\alpha \gamma} \right)
\end{align*}

In order to prove Theorem \ref{tc1} it suffices to construct holomorphic functions
\begin{align*}
\varphi_i&:\mathcal{U}_i'\to \mathbb{C}^n \\
s &:D \to \mathbb{C}^{m}
\end{align*}
and a matrix function
\begin{center}
$\left[\begin{matrix}
b_{i}^{11}(\xi_i,u) & b_{i}^{12}(\xi_i,u) & \cdots & b_{i}^{1p}(\xi_i, u)\\
b_{i}^{21}(\xi_i, u) & b_{i}^{22}(\xi_i, u) & \cdots & b_{i}^{2p}(\xi_i, u)\\
\cdot & \cdot & \cdots & \cdot\\
\cdot & \cdot & \cdots & \cdot\\
\cdot & \cdot & \cdots & \cdot\\
b_{i}^{p1}(\xi_i, u) & b_{i}^{p2}(\xi_i, u) & \cdots & b_{i}^{pp}(\xi_i, u)\\
\end{matrix}\right]$
\end{center}
where $b_i^{\alpha\beta}:\mathcal{U}_i'\to \mathbb{C}$ such that
\begin{align}
\varphi_i(\xi_i, 0)=\xi_i, \,\,\,\,\,s(0)&=0,\,\,\,\,\, b_i^{\alpha\beta}(\xi_i,0)= b_{0i}^{\alpha\beta}(\xi_i)\label{tc211}\\
\varphi_i \left(f_{ij}', u \right) &= f_{ij}\left(\varphi_j, s(u)\right) \label{tc2}
\end{align}
\begin{align}
\sum_{\beta=1}^p b_i^{\alpha\beta}\left(f_{ij}' (\xi_j, u), u\right) r_{ij}'^{\beta\gamma}(\xi_j , u) = \sum_{\beta=1}^p r_{ij}^{\alpha\beta} \left (\varphi_j(\xi_j, u), s(u) \right) b_j^{\beta\gamma}(\xi_j, u)\label{tc3}
\end{align}
\begin{align}
\sum_{\eta=1}^p g_{i\alpha\beta}^\eta\left( \varphi_i(\xi_i,u), s(u) \right) b_{i}^{\eta \gamma} = \sum_{\eta, \delta=1}^p b_{i}^{\alpha \eta} b_{i}^{\beta \delta} g_{i \eta \delta}'^\gamma (\xi_i, u) + \sum_{\eta=1}^p b_{i}^{\alpha \eta}  T_{i}'^\eta\left( b_{i}^{\beta \gamma} \right) - \sum_{\delta=1}^p b_{i}^{\beta \delta} T_{0i}'^\delta \left( b_{i}^{\alpha \gamma} \right)\label{tc30}
\end{align}
By setting
\begin{align*}
\left(\varphi_i, s \right): \mathcal{U}_i' \to \mathbb{C}^n \times\mathbb{C}^m,
\end{align*}
we have
\begin{align*}
\left(\varphi_i,s \right)^*(T_i^\alpha)&=\sum_{\beta=1}^n T_i^{\alpha\gamma}\left(\varphi_i(\xi_i, u), s(u)\right)\frac{\partial}{\partial z_i^\gamma} \\
\left(\varphi_i,s\right)_*(T_i'^\alpha)&= \sum_{\beta,\gamma=1}^n T_i'^{\alpha\beta}(\xi_i, u) \frac{\partial \varphi_i^\gamma(\xi_i, u)}{\partial \xi_i^\beta}\frac{\partial}{\partial z_i^\gamma},
\end{align*}
so that the following is satisfied: for $\alpha=1,...,p$
\begin{align}
\sum_{\beta=1}^p \sum_{\gamma,\sigma=1}^n b_i^{\alpha\beta}(\xi_i, u) T_i'^{\beta\sigma}(\xi_i, u)\frac{\partial \varphi_i^\gamma(\xi_i, u)}{\partial \xi_i^\sigma} \frac{\partial}{\partial z_i^\gamma}= \sum_{\gamma=1}^n T_i^{\alpha\gamma}\left(\varphi_i(\xi_i, u), s(u)\right)\frac{\partial}{\partial z_i^\gamma} \label{tc4}
\end{align}

First we prove the existence of formal solution of $(\ref{tc2})-(\ref{tc4})$. We recall Notation \ref{te27}. Then $(\ref{tc2})-(\ref{tc4})$ are equivalent to the following systems of congruences

\begin{align}
\varphi_i^\mu \left(f_{ij}'(\xi_j, u), u \right)& \equiv_\mu f_{ij}\left(\varphi_j^\mu(\xi_j, u), s^\mu(u)\right) \label{tc5}\\
\sum_{\beta=1}^p b_i^{\alpha\beta\mu}\left(f_{ij}'(\xi_j, u),u\right) r_{ij}'^{\beta\gamma}& (\xi_j, u)  \equiv_\mu \sum_{\beta=1}^p  r_{ij}^{\alpha\beta}\left(\varphi_j^\mu(\xi_j, u), s^\mu(u) \right) b_j^{\beta\gamma\mu}(\xi_j, u) \label{tc6} \\
\sum_{\eta=1}^p g_{i\alpha\beta}^\eta\left( \varphi_i^\mu(\xi_i,u), s(u) \right) b_{i}^{\eta \gamma \mu} \equiv_\mu \sum_{\eta, \delta=1}^p& b_{i}^{\alpha \eta \mu} b_{i}^{\beta \delta \mu} g_{i \eta \delta}'^\gamma (\xi_i, u) + \sum_{\eta=1}^p b_{i}^{\alpha \eta \mu}  T_{i}'^\eta\left( b_{i}^{\beta \gamma \mu} \right) - \sum_{\delta=1}^p b_{i}^{\beta \delta \mu} T_{i}'^\delta \left( b_{i}^{\alpha \gamma \mu} \right)\label{tc31}\\
\sum_{\beta=1}^p \sum_{\gamma,\sigma=1}^n b_i^{\alpha\beta \mu}(\xi_i, u) T_i'^{\beta\sigma}(\xi_i, u)&\frac{\partial \varphi_i^{\gamma\mu}(\xi_i, u)}{\partial \xi_i^\sigma} \frac{\partial}{\partial z_i^\gamma} \equiv_\mu  \sum_{\gamma=1}^n T_i^{\alpha\gamma}\left(\varphi_i^\mu (\xi_i, u), s^\mu (u)\right)\frac{\partial}{\partial z_i^\gamma}  \label{tc7}
\end{align}
for each $\mu=1,2,...$.  We shall construct $\varphi_i^\mu, s^\mu$, and $b_i^{\alpha\beta \mu}$ satisfying $(\ref{tc5})_\mu- (\ref{tc7})_\mu$ by induction on $\mu$. We assume that $\varphi_i^{\mu-1}, s^{\mu-1}$ and $b_i^{\alpha\beta(\mu-1)}$ are already determined. Then we define homogenous polynomials $\Gamma_{ij|\mu}^\alpha, B_{ij|\mu}^{\alpha\beta}, G_{i|\mu}^{\alpha \beta \gamma}$, and $\Pi_{i|\mu}^{\alpha \gamma}$ of degree $\mu$ by the following congruences:
{\small{\begin{align}
\Gamma_{ij|\mu}^\alpha &\equiv_\mu \varphi_i^{\alpha(\mu-1)}\left(f_{ij}'(\xi_j, u), u \right)- f_{ij}^{\alpha}\left( \varphi_j^{\mu-1}(\xi_j, u), s^{\mu-1}(u)\right) \label{tc35} \\
\sum_{\beta,\sigma=1}^p B_{ij|\mu}^{\alpha\beta} b_{0i}^{\beta\sigma} r_{0ij}'^{\sigma \gamma} &\equiv_\mu \sum_{\beta=1}^p b_i^{\alpha\beta (\mu-1)}\left(f_{ij}'(\xi_j, u), u\right) r_{ij}'^{\beta\gamma}(\xi_j, u) - \sum_{\beta=1}^p r_{ij}^{\alpha\beta} \left(\varphi_j^{\mu-1}(\xi_j, u\right), s^{\mu-1}(u))b_j^{\beta\gamma(\mu-1)}(\xi_j, u)\label{tc36}\\
G_{i|\mu}^{\alpha \beta \gamma}\equiv_\mu     \sum_{\eta, \delta=1}^p  b_{i}^{\alpha \eta (\mu-1)}& b_{i}^{\beta \delta (\mu-1)} g_{i \eta \delta}'^\gamma + \sum_{\eta=1}^p b_{i}^{\alpha \eta (\mu-1)}  T_{i}'^\eta\left( b_{i}^{\beta \gamma (\mu-1)} \right) - \sum_{\delta=1}^p b_{i}^{\beta \delta (\mu-1)} T_{i}'^\delta \left( b_{i}^{\alpha \gamma (\mu-1)} \right)   -\sum_{\eta=1}^p g_{i\alpha\beta}^\eta\left( \varphi_i^{\mu-1}, s^{\mu-1} \right) b_{i}^{\eta \gamma (\mu-1)}  \label{tc37}    \\
\Pi_{i|\mu}^{\alpha\gamma}&\equiv_\mu T_i^{\alpha\gamma}\left(\varphi_i^{\mu-1}, s^{\mu-1}(u) \right)- \sum_{\beta=1}^p \sum_{\sigma=1}^n b_i^{\alpha\beta (\mu-1)}(\xi_i, u) T_i'^{\beta\sigma}(\xi_i, u)\frac{\partial \varphi_i^{\gamma(\mu-1)}(\xi_i, u)}{\partial \xi_i^\sigma} \label{tc38}
\end{align}}}

We set
\begin{align*}
\Gamma_{ij|\mu}:=\sum_{\alpha=1}^n \Gamma_{ij|\mu}^\alpha \frac{\partial}{\partial z_i^\alpha},\,\,\,\,\,\,\,\,\,\, \Pi_{i|\mu}^\alpha:=\sum_{\gamma=1}^n \Pi_{i|\mu}^{\alpha \gamma} \frac{\partial}{\partial z_i^\gamma}
\end{align*}

\begin{lemma}\label{tc39}
We have the following relations$:$
\begin{align}
&\,\,\,\,\,\,\,\,\Gamma_{jk|\mu}- \Gamma_{ik|\mu} + \Gamma_{ij|\mu}=0\label{tc32}\\
 \sum_{\beta=1}^p r_{0ij}^{\alpha \eta} & \Pi_{j|\mu}^\eta = \Pi_{i|\mu}^\alpha - \left[ \Gamma_{ij|\mu}, T_{0i}^\alpha \right] + \sum_{b=1}^p B_{ij|\mu}^{\alpha b} T_{0i}^b \label{tc33}\\
 \left[ \Pi_{i|\mu}^{\alpha}, T_{0i}^\beta     \right]  - & \left[ \Pi_{i|\mu}^\beta,  T_{0i}^\alpha \right]  - \sum_{\gamma=1}^p g_{0i\alpha\beta}^\gamma \Pi_{i|\mu}^\gamma + \sum_{\gamma=1}^p G_{i|\mu}^{\alpha \beta\gamma} T_{0i}'^\gamma=0 \label{tc34}
\end{align}
\end{lemma}

\begin{proof}
$(\ref{tc32})$ follows from \cite{Kod05} p.292-p.294. We prove $(\ref{tc33})$. In fact, from $(\ref{tc35})$ and $(\ref{tc36})$ and $(\ref{tc38})$, we have
{\small{\begin{align*}
&\sum_{\beta=1}^p r_{0ij}^{\alpha\eta} \Pi_{j|\mu}^\eta\equiv_\mu \sum_{\beta=1}^p r_{ij}^{\alpha\eta}\left(\varphi_j^{\mu-1}, s^{\mu-1}\right) \Pi_{j|\mu }^\eta \\
 & \equiv_\mu \sum_{\eta=1}^p\sum_{\gamma=1}^nr_{ij}^{\alpha\eta}\left(\varphi_j^{\mu-1}, s^{\mu-1}\right)\left( T_j^{\eta\gamma}\left(\varphi_j^{\mu-1} , s^{\mu-1}\right)- \sum_{\beta=1}^p \sum_{\sigma=1}^n b_j^{\eta\beta (\mu-1)} T_j'^{\beta\sigma}\frac{\partial \varphi_j^{\gamma(\mu-1)}}{\partial \xi_j^\sigma}     \right) \frac{\partial}{\partial z_j^\gamma} \\
&\equiv_\mu \sum_{\eta=1}^p\sum_{\gamma,\delta=1}^nr_{ij}^{\alpha\eta}\left(\varphi_j^{\mu-1}, s^{\mu-1} \right)\left( T_j^{\eta\gamma}\left(\varphi_j^{\mu-1} , s^{\mu-1} \right)- \sum_{\beta=1}^p \sum_{\sigma=1}^n b_j^{\eta\beta (\mu-1)} T_j'^{\beta\sigma}\frac{\partial \varphi_j^{\gamma(\mu-1)}}{\partial \xi_j^\sigma}     \right) \frac{\partial f_{0ij}^\delta}{\partial z_j^\gamma} \frac{\partial}{\partial z_i^\delta} \\
&\equiv_\mu  \sum_{\eta=1}^p\sum_{\gamma,\delta=1}^nr_{ij}^{\alpha\eta}\left(\varphi_j^{\mu-1}, s^{\mu-1}\right)\left( T_j^{\eta\gamma}\left( \varphi_j^{\mu-1} , s^{\mu-1}\right)- \sum_{\beta=1}^p \sum_{\sigma=1}^n b_j^{\eta\beta (\mu-1)} T_j'^{\beta\sigma}\frac{\partial \varphi_j^{\gamma(\mu-1)}}{\partial \xi_j^\sigma}     \right) \frac{\partial f_{ij}^\delta}{\partial z_j^\gamma}\left(\varphi_j^{\mu-1}, s^{\mu-1}\right) \frac{\partial}{\partial z_i^\delta} \\
&\equiv_\mu \sum_{\delta=1}^n T_i^{\alpha \delta}\left(f_{ij}\left(\varphi_j^{\mu-1},s^{\mu-1}\right), s^{\mu-1}\right) \frac{\partial}{\partial z_i^\delta}- \sum_{\beta,\eta=1}^p \sum_{\delta, \sigma=1}^n r_{ij}^{\alpha\eta}\left(\varphi_j^{\mu-1}, s^{\mu-1} \right) b_j^{\eta\beta (\mu-1)} T_j'^{\beta\sigma} \frac{\partial f_{ij}^\delta \left(\varphi_j^{\mu-1}, s^{\mu-1}\right)}{\partial \xi_j^\sigma} \frac{\partial}{\partial z_i^\delta} \\
&\equiv_\mu \sum_{\delta=1}^n T_i^{\alpha\delta}\left(\varphi_i^{\mu-1} \left(f_{ij}', u \right)- \Gamma_{ij|\mu} , s^{\mu-1} \right)\frac{\partial}{\partial z_i^\delta}\\
&- \sum_{\beta=1}^p\sum_{\delta, \sigma=1}^n \left( \sum_{\eta=1}^p b_i^{\alpha\eta (\mu-1)}\left(f_{ij}', u \right)r_{ij}'^{\eta\beta} - \sum_{b,c=1}^n B_{ij|\mu}^{\alpha b}b_{0i}^{bc} r_{0ij}'^{c \beta}      \right) T_j'^{\beta\sigma}\frac{\partial \left( \varphi_i^{\delta(\mu-1)}\left(f_{ij}',u\right)- \Gamma_{ij|\mu}^\delta   \right)}{\partial \xi_j^\sigma} \frac{\partial}{\partial z_i^\delta}\\
&\equiv_\mu \sum_{\delta=1}^n \left(T_i^{\alpha\delta}\left(\varphi_i^{\mu-1} \left(f_{ij}', u\right), s^{\mu-1}\right) - \sum_{\gamma=1}^n \frac{\partial T_{0i}^{\alpha\delta}}{\partial z_i^\gamma} \Gamma_{ij|\mu}^\gamma\right) \frac{\partial}{\partial z_i^\delta} - \sum_{\delta, \sigma , \gamma =1}^n\sum_{\eta,\beta =1}^p b_i^{\alpha\eta (\mu-1)}\left(f_{ij}', u \right) r_{ij}'^{\eta\beta} T_j'^{\beta\sigma} \frac{\partial \varphi_i^{\delta(\mu-1)}}{\partial \xi_i^\gamma}\left(f_{ij}',u \right)\frac{\partial f_{ij}'^\gamma}{\partial \xi_j^\sigma}\frac{\partial}{\partial z_i^\delta} \\
&+\sum_{\beta,\eta=1}^p \sum_{\delta, \sigma=1}^n b_{0i}^{\alpha\eta} r_{0ij}'^{\eta\beta} T_{0j}'^{\beta\sigma} \frac{\partial \Gamma_{ij|\mu}^\delta}{\partial z_j^\sigma}\frac{\partial}{\partial z_i^\delta} +\sum_{\beta=1}^p \sum_{\delta,\sigma=1}^n \sum_{b,c=1}^n B_{ij|\mu}^{\alpha b} b_{0i}^{bc} r_{0ij}'^{c\beta} T_{0j}'^{\beta\sigma} \frac{\partial f_{0ij}^\delta }{\partial z_j^\sigma} \frac{\partial}{\partial z_i^\delta}\\
&\equiv_\mu \sum_{\delta=1}^n \left(  T_i^{\alpha\delta}\left(\varphi_i^{\mu-1}\left(f_{ij}', u\right),s^{\mu-1} \right) - \sum_{\delta, \gamma=1}^n \sum_{\eta=1}^p  b_i^{\alpha \eta(\mu-1)}\left(f_{ij}', u\right) T_i'^{\eta \gamma} \frac{\partial \varphi_i^{\delta(\mu-1)}}{\partial \xi_i^\gamma} \left(f_{ij}', u \right)                       \right) \frac{\partial}{\partial z_i^\delta}\\
&- \sum_{\gamma,\delta=1}^n \frac{\partial T_{0i}^{\alpha\delta}}{\partial z_i^\gamma } \Gamma_{ij|\mu}^\gamma \frac{\partial}{\partial z_i^\delta} + \sum_{\sigma, \delta=1}^n T_{0i}^{\alpha \eta}\frac{\partial \Gamma_{ij|\mu}^\delta}{\partial z_i^\eta}\frac{\partial}{\partial z_i^\delta} + \sum_{b =1}^n B_{ij|\mu}^{\alpha b} T_{0i}^b \\
&\equiv_\mu \Pi_{i|\mu}^\alpha + \left[ \sum_{\eta=1}^n T_{0i}^{\alpha\eta} \frac{\partial}{\partial z_i^\eta}, \sum_{\delta=1}^n  \Gamma_{ij|\mu}^\delta  \frac{\partial}{\partial z_i^\delta}    \right]+ \sum_{b=1}^n B_{ij|\mu}^{\alpha b} T_{0i}^b
\end{align*}}}

We prove $(\ref{tc38})$. In fact, since $\left[T_i^\alpha, T_i^\beta \right]= \sum_{\gamma=1}^p g_{i\alpha\beta}^\gamma T_i^\gamma$, we have from $(\ref{tc38})$
{\small{\begin{align*}
&\left[\sum_{\eta=1}^n T_i^{\alpha\eta}(z_i,t)\frac{\partial}{\partial z_i^\eta}  , \sum_{\sigma=1}^n T_i^{\beta\sigma} (z_i, t)  \frac{\partial}{\partial z_i^\sigma}        \right]=  \sum_{\gamma=1}^p \sum_{\eta=1}^n g_{i\alpha\beta}^\gamma(z_i,t) T_i^{\gamma\eta}(z_i,t)\frac{\partial}{\partial z_i^\eta}\\
&\Longrightarrow \sum_{\eta, \sigma=1}^n T_i^{\alpha\eta}(z_i, t)\frac{\partial T_i^{\beta\sigma} (z_i, t)}{\partial z_i^\eta} \frac{\partial}{\partial z_i^\sigma} - \sum_{\eta,\sigma=1}^n T_i^{\beta\sigma}(z_i,t)\frac{\partial T_i^{\alpha\eta}(z_i,t)}{\partial z_i^\sigma} \frac{\partial}{\partial z_i^\eta}=\sum_{\gamma=1}^p \sum_{\eta=1}^n g_{i\alpha\beta}^\gamma(z_i,t) T_i^{\gamma\eta}(z_i,t)\frac{\partial}{\partial z_i^\eta}\\
&\Longrightarrow \sum_{\eta,\sigma=1}^n T_i^{\alpha\eta}\left(\varphi_i^{\mu-1},s^{\mu-1}\right)\frac{\partial T_i^{\beta\sigma}}{\partial z_i^\eta}\left(\varphi_i^{\mu-1},s^{\mu-1}\right) \frac{\partial}{\partial z_i^\sigma} - \sum_{\eta , \sigma=1}^n T_i^{\beta\sigma}\left(\varphi_i^{\mu-1}, s^{\mu-1}\right)\frac{\partial T_i^{\alpha\eta}}{\partial z_i^\sigma}\left( \varphi_i^{\mu-1}, s^{\mu-1}\right)\frac{\partial}{\partial z_i^\eta} \\
&=\sum_{\gamma=1}^p\sum_{\eta=1}^n g_{i\alpha\beta}^\gamma \left(\varphi_i^{\mu-1}, s^{\mu-1}\right) T_i^{\gamma\eta}\left( \varphi_i^{\mu-1}, s^{\mu-1}\right)\frac{\partial}{\partial z_i^\eta}\\
&\Longrightarrow  \sum_{\eta,\sigma=1}^n \left( \Pi_{i|\mu}^{\alpha\eta} + \sum_{c=1}^p \sum_{d=1}^n b_i^{\alpha c(\mu-1)} T_i'^{c d}\frac{\partial \varphi_i^{\eta(\mu-1)}}{\partial \xi_i^d}    \right) \frac{\partial T_i^{\beta\sigma}}{\partial z_i^\eta} \left( \varphi_i^{\mu-1},s^{\mu-1} \right)\frac{\partial}{\partial z_i^\sigma} \\
&-\sum_{\eta, \sigma=1}^n \left( \Pi_{i|\mu}^{\beta\sigma} + \sum_{c=1}^p \sum_{d=1}^n b_i^{\beta c(\mu-1)}T_i'^{cd}\frac{\partial \varphi_i^{\sigma(\mu-1)}}{\partial \xi_i^d}   \right) \frac{\partial T_i^{\alpha\eta}}{\partial z_i^\sigma} \left(\varphi_i^{\mu-1}, s^{\mu-1}\right)\frac{\partial}{\partial z_i^\eta}\\
&\equiv_\mu \sum_{\gamma=1}^p \sum_{\eta=1}^n g_{i\alpha\beta}^{\gamma}\left( \varphi_i^{\mu-1}, s^{\mu-1} \right)\left(  \Pi_{i|\mu}^{\gamma\eta}+ \sum_{c=1}^p \sum_{d=1}^n b_i^{\gamma c(\mu-1)}T_i'^{c d} \frac{\partial \varphi_i^{\eta(\mu-1)}}{\partial \xi_i^d}  \right) \frac{\partial}{\partial z_i^\eta}
\end{align*}}}

Then we have
\begin{align*}
&\sum_{\eta, \sigma=1}^n \Pi_{i|\mu}^{\alpha \eta}\frac{\partial T_{0i}^{\beta\sigma}}{\partial z_i^\eta} \frac{\partial}{\partial z_i^\sigma} + \sum_{c=1}^p\sum_{\sigma, d=1}^n b_i^{\alpha c(\mu-1)} T_i'^{cd}\frac{\partial}{\partial \xi_i^d}\left( T_i^{\beta\sigma}\left(\varphi_i^{\mu-1}, s^{\mu-1}\right)\right) \frac{\partial}{\partial z_i^\sigma}\\
&- \sum_{\eta, \sigma=1}^n \Pi_{i|\mu}^{\beta\sigma}\frac{\partial T_{0i}^{\alpha \eta}}{\partial z_i^\sigma} \frac{\partial}{\partial z_i^\eta} - \sum_{c=1}^p \sum_{\eta, d=1}^n b_i^{\beta c (\mu-1)} T_i'^{cd} \frac{\partial}{\partial \xi_i^d}\left( T_i^{\alpha \eta} \left( \varphi_i^{\mu-1}, s^{\mu-1} \right)  \right) \frac{\partial}{\partial z_i^\eta} \\
&\equiv_\mu \sum_{\gamma=1}^p \sum_{\eta=1}^n g_{0i\alpha\beta}^\gamma \Pi_{i|\mu}^{\gamma\eta} \frac{\partial}{\partial z_i^\eta} + \sum_{\gamma,c=1}^p \sum_{\eta, d=1}^n g_{i\alpha\beta}^\gamma \left( \varphi_i^{\mu-1}, s^{\mu-1}\right) b_i^{\gamma c (\mu-1)} T_i'^{cd}\frac{\partial \varphi_i^{\eta(\mu-1)} }{\partial \xi_i^d}\frac{\partial}{\partial z_i^\eta}
\end{align*} 
Then from $(\ref{tc38})$ and $\left[ T_i'^\alpha, T_i'^\beta \right]=\sum_{\gamma=1}^p g_{i\alpha\beta}'^\gamma T_i'^\gamma$, we have
\begin{align*}
&\sum_{\gamma=1}^p \sum_{\eta=1}^n g_{0i\alpha\beta}^\gamma \Pi_{i|\mu}^{\gamma\eta} \frac{\partial}{\partial z_i^\eta} + \sum_{\gamma,c=1}^p \sum_{\eta, d=1}^n g_{i\alpha\beta}^\gamma\left(\varphi_i^{\mu-1}, s^{\mu-1}\right) b_i^{\gamma c (\mu-1)} T_i'^{cd}\frac{\partial \varphi_i^{\eta(\mu-1)} }{\partial \xi_i^d}\frac{\partial}{\partial z_i^\eta}\\
&\equiv_\mu \sum_{\eta, \sigma=1}^n \Pi_{i|\mu}^{\alpha \eta}\frac{\partial T_{0i}^{\beta\sigma}}{\partial z_i^\eta} \frac{\partial}{\partial z_i^\sigma} + \sum_{c=1}^p\sum_{\sigma, d=1}^n b_i^{\alpha c(\mu-1)} T_i'^{cd}\frac{\partial}{\partial \xi_i^d}\left( \Pi_{i|\mu}^{\beta\sigma} + \sum_{g=1}^p \sum_{h=1}^n b_i^{\beta g(\mu-1)} T_i'^{gh} \frac{\partial \varphi_i^{\sigma(\mu-1)}}{\partial \xi_i^h}\right) \frac{\partial}{\partial z_i^\sigma}\\
&- \sum_{\eta, \sigma=1}^n \Pi_{i|\mu}^{\beta\sigma}\frac{\partial T_{0i}^{\alpha \eta}}{\partial z_i^\sigma} \frac{\partial}{\partial z_i^\eta} - \sum_{c=1}^p \sum_{\eta, d=1}^n b_i^{\beta c (\mu-1)} T_i'^{cd} \frac{\partial}{\partial \xi_i^d}\left( \Pi_{i|\mu}^{\alpha \eta} + \sum_{g=1}^p \sum_{h=1}^n b_i^{\alpha g(\mu -1)} T_i'^{gh} \frac{\partial \varphi_i^{\eta(\mu-1)}}{\partial \xi_i^h} \right) \frac{\partial}{\partial z_i^\eta} \\
&\equiv_\mu \sum_{\eta, \sigma=1}^n \Pi_{i|\mu}^{\alpha \eta}\frac{\partial T_{0i}^{\beta\sigma}}{\partial z_i^\eta} \frac{\partial}{\partial z_i^\sigma} + \sum_{\sigma, d=1}^n T_{0i}^{\alpha d}\frac{ \partial \Pi_i^{\beta\sigma} }{\partial z_i^d}\frac{\partial}{\partial z_i^\sigma} + \sum_{c=1}^p\sum_{\sigma, d=1}^n b_i^{\alpha c(\mu-1)} T_i'^{cd}\frac{\partial}{\partial \xi_i^d}\left( \sum_{g=1}^p \sum_{h=1}^n b_i^{\beta g(\mu-1)} T_i'^{gh} \right) \frac{\partial \varphi_i^{\sigma(\mu-1)}}{\partial \xi_i^h} \frac{\partial}{\partial z_i^\sigma}\\
&- \sum_{\eta, \sigma=1}^n \Pi_{i|\mu}^{\beta\sigma}\frac{\partial T_{0i}^{\alpha \eta}}{\partial z_i^\sigma} \frac{\partial}{\partial z_i^\eta} - \sum_{\eta, d=1}^n T_{0i}^{\beta d}\frac{\partial  \Pi_{i|\mu}^{\alpha \eta}}{\partial z_i^d} \frac{\partial}{\partial z_i^\eta} - \sum_{c=1}^p \sum_{\eta, d=1}^n b_i^{\beta c (\mu-1)} T_i'^{cd} \frac{\partial}{\partial \xi_i^d}\left(  \sum_{g=1}^p \sum_{h=1}^n b_i^{\alpha g(\mu -1)} T_i'^{gh} \right) \frac{\partial \varphi_i^{\eta(\mu-1)}}{\partial \xi_i^h}  \frac{\partial}{\partial z_i^\eta}\\
&\equiv_\mu  \sum_{\eta, \sigma=1}^n \Pi_{i|\mu}^{\alpha \eta}\frac{\partial T_{0i}^{\beta\sigma}}{\partial z_i^\eta} \frac{\partial}{\partial z_i^\sigma} + \sum_{\sigma, d=1}^n T_{0i}^{\alpha d}\frac{ \partial \Pi_i^{\beta\sigma} }{\partial z_i^d}\frac{\partial}{\partial z_i^\sigma} + \sum_{c,g=1}^p \sum_{\sigma , h=1}^n b_i^{\alpha c(\mu-1)} T_i'^c\left( b_i^{\beta g(\mu-1)} \right) T_i'^{gh}\frac{\partial \varphi_i^{\sigma(\mu-1)}}{\partial \xi_i^h} \frac{\partial}{\partial z_i^\sigma}\\
&- \sum_{\eta, \sigma=1}^n \Pi_{i|\mu}^{\beta\sigma}\frac{\partial T_{0i}^{\alpha \eta}}{\partial z_i^\sigma} \frac{\partial}{\partial z_i^\eta} - \sum_{\eta, d=1}^n T_{0i}^{\beta d}\frac{\partial  \Pi_{i|\mu}^{\alpha \eta}}{\partial z_i^d} \frac{\partial}{\partial z_i^\eta}  - \sum_{c,g=1}^p \sum_{\eta , h =1}^n b_i^{\beta c(\mu-1)} T_i'^c\left( b_i^{\alpha g (\mu-1)} \right) T_i'^{gh} \frac{\partial \varphi_i^{\eta(\mu-1)}}{\partial \xi_i^h}\frac{\partial}{\partial z_i^\eta}\\
&+ \sum_{c, g , \sigma =1}^p \sum_{h, \eta=1}^n b_i^{\alpha c(\mu-1)} b_i^{\beta g(\mu-1)} g_{i c g}'^\sigma T_i'^{\sigma h} \frac{\partial \varphi_i^{\eta(\mu-1)}}{\partial \xi_i^h} \frac{\partial}{\partial z^\eta}
\end{align*}
Then from $(\ref{tc37})$, this implies that
\begin{align*}
&\sum_{\eta, \sigma=1}^n \Pi_{i|\mu}^{\alpha \eta}\frac{\partial T_{0i}^{\beta\sigma}}{\partial z_i^\eta} \frac{\partial}{\partial z_i^\sigma} + \sum_{\sigma, d=1}^n T_{0i}^{\alpha d}\frac{ \partial \Pi_i^{\beta\sigma} }{\partial z_i^d}\frac{\partial}{\partial z_i^\sigma}- \sum_{\eta, \sigma=1}^n \Pi_{i|\mu}^{\beta\sigma}\frac{\partial T_{0i}^{\alpha \eta}}{\partial z_i^\sigma} \frac{\partial}{\partial z_i^\eta} - \sum_{\eta, d=1}^n T_{0i}^{\beta d}\frac{\partial  \Pi_{i|\mu}^{\alpha \eta}}{\partial z_i^d} \frac{\partial}{\partial z_i^\eta} \\
&\equiv_\mu \sum_{\gamma=1}^p \sum_{\eta=1}^n g_{0i\alpha\beta}^\gamma \Pi_{i|\mu}^{\gamma\eta} \frac{\partial}{\partial z_i^\eta} - \sum_{c=1}^p \sum_{\eta,d=1}^n G_{i|\mu}^{\alpha\beta c} T_i'^{cd}\frac{\partial \varphi_i^{\eta(\mu-1)} }{\partial \xi_i^d}\frac{\partial}{\partial z_i^\eta}
\end{align*}
which is equivalent to $(\ref{tc38})$. This completes the proof of Lemma \ref{tc39}.

\end{proof}

Our purpose is to determine $\varphi^\mu= \varphi^{\mu-1}+ \varphi_{i|\mu}, s^\mu= s^{\mu-1}+ s_\mu$, and $b_i^{\alpha\beta \mu}=b_i^{\alpha\beta(\mu-1)}+ b_{i|\mu}^{\alpha\beta}$ satisfying $(\ref{tc5})_\mu-(\ref{tc7})_\mu$.

\begin{lemma}\label{tc41}
$(\ref{tc5})_\mu-(\ref{tc7})_\mu$ are equivalent to the following equalities.
\begin{align}
\Gamma_{ij|\mu}&= \varphi_{j|\mu}- \varphi_{i|\mu}+ \sum_{\lambda=1}^m  s_\mu^\lambda \rho_{ij\lambda} \label{tc8}\\
\sum_{\beta, \sigma=1}^p B_{ij|\mu}^{\alpha\beta} b_{0i}^{\beta\sigma} r_{0ij}'^{\sigma \gamma} &= - \sum_{\beta=1}^p b_{i|\mu}^{\alpha \beta} r_{0ij}'^{\beta\gamma} +\sum_{\beta=1}^p \left[  \varphi_{j|\mu} , r_{0ij}^{\alpha\beta} \right] b_{0j}^{\beta\gamma} + \sum_{\beta=1}^p \sum_{\lambda=1}^m \frac{\partial r_{ij}^{\alpha\beta}}{\partial t_\lambda}|_{t=0} s_\mu^\lambda b_{0j}^{\beta\gamma}  + \sum_{\beta=1}^p r_{0ij}^{\alpha\beta} b_{j|\mu}^{\beta\gamma}   \label{tc9} \\
-G_{i|\mu}^{\alpha \beta \gamma} &= \sum_{\eta, \delta=1}^p b_{0i}^{\alpha \eta} b_{i|\mu}^{\beta \delta} g_{0i\eta\delta}'^\gamma + \sum_{\eta, \delta=1}^p b_{i|\mu}^{\alpha \eta} b_{0i}^{\beta \delta} g_{0i\eta \delta}'^\gamma + T_{0i}^\alpha\left(  b_{i|\mu}^{\beta \gamma}\right) + \sum_{\eta=1}^p b_{i|\mu}^{\alpha \eta} T_{0i}'^\eta\left( b_{0i}^{\beta \gamma} \right)- T_{0i}^\beta\left( b_{i|\mu}^{\alpha \gamma} \right) - \sum_{\delta=1}^p b_{i|\mu}^{\beta \delta} T_{0i}'^\delta\left( b_{0i}^{\alpha \gamma} \right) \label{tc40}  \\
&- \sum_{\eta=1}^p\sum_{c=1}^n \frac{\partial g_{0i\alpha\beta}^\eta}{\partial z_i^c}\varphi_{i|\mu}^c b_{0i}^{\eta \gamma} - \sum_{\eta=1}^p \sum_{\lambda=1}^m \frac{\partial g_{i\alpha\beta}^\eta}{\partial t_\lambda}|_{t=0} s_\mu^\lambda b_{0i}^{\eta \gamma} - \sum_{\eta=1}^p g_{0i\alpha\beta}^\eta b_{i|\mu}^{\eta \gamma}      \notag         \\
\Pi_{i|\mu}^\alpha &= - \left[ \varphi_{i|\mu}, T_{0i}^\alpha   \right] + \sum_{\lambda=1}^r s_\mu^\lambda  \alpha_{i\lambda}^\alpha         + \sum_{\beta=1}^p b_{i|\mu}^{\alpha \beta} T_{0i}'^\beta   \label{tc10}
\end{align}
where
\begin{align*}
\varphi_{i|\mu}=\sum_{\alpha=1}^n \varphi_{i|\mu}^\alpha \frac{\partial}{\partial z_i^\alpha}, \,\,\,\,\,\, \rho_{ij\lambda}= \sum_{\alpha=1}^n \frac{\partial f_{ij}^\alpha}{\partial t_\lambda}|_{t=0}\frac{\partial}{\partial z_i^\alpha},\,\,\,\,\,\,\,\,\alpha_{i\lambda}^\alpha = - \frac{\partial T_i^{\alpha}(z_i,t)}{\partial t_\lambda}|_{t=0}
\end{align*}
\end{lemma}

\begin{proof}
$(\ref{tc8})$ follows from \cite{Kod05} p. 290. We prove $(\ref{tc9})$. In fact,
\begin{align*}
&\sum_{\beta=1}^p \left( b_i^{\alpha\beta(\mu-1)}\left(f_{ij}',u\right) + b_{i|\mu}^{\alpha\beta}\left( f_{ij}', u \right)\right) r_{ij}'^{\beta\gamma} - \sum_{\beta=1}^p r_{ij}^{\alpha\beta} \left( \varphi_j^{\mu-1}+ \varphi_{j|\mu}, s^{\mu-1}+ s_\mu \right) \left( b_j^{\beta\gamma(\mu-1)} + b_{j|\mu}^{\beta\gamma}   \right)\\
&\equiv_\mu \sum_{\beta, \sigma=1}^p B_{ij|\mu}^{\alpha\beta} b_{0i}^{\beta\sigma} r_{0ij}'^{\sigma \gamma} + \sum_{\beta=1}^p b_{i|\mu}^{\alpha \beta} r_{0ij}'^{\beta\gamma} - \sum_{\beta=1}^p \sum_{\eta=1}^n \frac{\partial r_{0ij}^{\alpha \beta}}{\partial z_j^\eta} \varphi_{j|\mu}^\eta b_{0j}^{\beta\gamma} - \sum_{\beta=1}^p \sum_{\lambda=1}^r \frac{\partial r_{ij}^{\alpha\beta}}{\partial t_\lambda}|_{t=0} s_\mu^\lambda b_{0j}^{\beta\gamma} - \sum_{\beta=1}^p r_{0ij}^{\alpha\beta} b_{j|\mu}^{\beta\gamma} = 0
\end{align*}

We prove $(\ref{tc40})$. In fact,
\begin{align*}
&\sum_{\eta, \delta=1}^p \left( b_i^{\alpha\eta(\mu-1)} + b_{i|\mu}^{\alpha\eta} \right)\left(  b_i^{\beta \delta(\mu-1)} + b_{i|\mu}^{\beta \delta} \right) g_{i\eta \delta}'^\gamma+ \sum_{\eta=1}^p \left( b_i^{\alpha \eta(\mu-1)} + b_{i|\mu}^{\alpha \eta} \right) T_i'^\eta\left(b_i^{\beta \gamma(\mu-1) } + b_{i|\mu}^{\beta \gamma} \right)\\
&- \sum_{\delta=1}^p \left(b_i^{\beta \delta(\mu-1)} + b_{i|\mu}^{\beta \delta} \right) T_i'^\delta\left( b_i^{\alpha \gamma(\mu-1)} + b_{i|\mu}^{\alpha \gamma} \right) - \sum_{\eta=1}^p g_{i\alpha\beta}^\eta \left( \varphi_i^{\mu-1} + \varphi_{i|\mu}, s^{\mu-1}+ s_\mu \right)\left( b_i^{\eta \gamma(\mu-1)} + b_{i|\mu}^{\eta \gamma} \right)=0
\end{align*}
Then we have
\begin{align*}
- G_{i|\mu}^{\alpha \beta \gamma}&= \sum_{\eta, \delta=1}^p b_{0i}^{\alpha \eta} b_{i|\mu}^{\beta \delta} g_{0i\eta\delta}'^\gamma + \sum_{\eta, \delta=1}^p b_{i|\mu}^{\alpha \eta} b_{0i}^{\beta \delta} g_{0i\eta \delta}^\gamma + T_{0i}^\alpha\left(  b_{i|\mu}^{\beta \gamma}\right) + \sum_{\eta=1}^p b_{i|\mu}^{\alpha \eta} T_{0i}'^\eta\left( b_{0i}^{\beta \gamma} \right)- T_{0i}^\beta\left( b_{i|\mu}^{\alpha \gamma} \right) - \sum_{\delta=1}^p b_{i|\mu}^{\beta \delta} T_{0i}'^\delta\left( b_{0i}^{\alpha \gamma} \right) \\
&- \sum_{\eta=1}^p\sum_{c=1}^n \frac{\partial g_{0i\alpha\beta}^\eta}{\partial z_i^c}\varphi_{i|\mu}^c b_{0i}^{\eta \gamma} - \sum_{\eta=1}^p \sum_{\lambda=1}^m \frac{\partial g_{i\alpha\beta}^\eta}{\partial t_\lambda}|_{t=0} s_\mu^\lambda b_{0i}^{\eta \gamma} - \sum_{\eta=1}^p g_{0i\alpha\beta}^\eta b_{i|\mu}^{\eta \gamma}
\end{align*}

We prove $(\ref{tc10})$. In fact,
\begin{align*}
&T_i^{\alpha \gamma}\left( \varphi_i^{\mu-1} +  \varphi_\mu , s^{\mu-1} + s_\mu \right) - \sum_{\beta=1}^p \sum_{\sigma=1}^n \left(  b_i^{\alpha\beta(\mu-1)} + b_{i|\mu}^{\alpha\beta} \right) T_i'^{\beta\sigma} \frac{\partial \left( \varphi_i^{\gamma(\mu-1)} + \varphi_{i|\mu}^\gamma \right)}{\partial \xi_i^\sigma }\\
&\equiv_\mu \Pi_{i|\mu}^{\alpha \gamma} + \sum_{\eta=1}^n \frac{\partial T_{0i}^{\alpha \gamma}}{\partial z_i^\eta} \varphi_{i|\mu}^\eta + \sum_{\lambda=1}^m \frac{\partial T_i^{\alpha \gamma}}{\partial t_\lambda}|_{t=0} s_\mu^\lambda - \sum_{\beta=1}^p  b_{i|\mu}^{\alpha \beta} T_{0i}'^{\beta\gamma} - \sum_{\beta=1}^p \sum_{\sigma=1}^n b_{0i}^{\alpha \beta} T_{0i}'^{\beta\sigma} \frac{\partial \varphi_{i|\mu}^\gamma}{\partial z_i^\sigma}
\end{align*}
Then by multiplying $\sum_{\gamma=1}^n \frac{\partial}{\partial z_i^\gamma}$,  we have
\begin{align*}
\Pi_{i|\mu}^{\alpha}& \equiv_\mu - \sum_{\eta,\gamma=1}^n \frac{\partial T_{0i}^{\alpha \gamma}}{\partial z_i^\eta}  \varphi_{i|\mu}^\eta \frac{\partial}{\partial z_i^\gamma} - \sum_{\gamma=1}^n \sum_{\lambda=1}^m s_\mu^\lambda \frac{\partial T_i^{\alpha \gamma}}{\partial t_\lambda}|_{t=0}  \frac{\partial}{\partial z_i^\gamma} + \sum_{\beta=1}^p b_{i|\mu}^{\alpha \beta} T_{0i}'^\beta + \sum_{\sigma=1}^n T_{0i}^{\alpha \sigma} \frac{\partial \varphi_{i|\mu}^\gamma}{\partial z_i^\sigma} \frac{\partial}{\partial z_i^\gamma}\\
&\equiv_\mu - \left[ \varphi_{i|\mu}, T_{0i}^\alpha   \right] + \sum_{\lambda=1}^m s_\mu^\lambda  \alpha_{i\lambda}^\alpha         + \sum_{\beta=1}^p b_{i|\mu}^{\alpha \beta} T_{0i}'^\beta
\end{align*}
This completes the proof of Lemma \ref{tc41}.
\end{proof}

We define an element $\overline{\Pi}_{i|\mu} \in \Gamma\left( U_i, \mathscr{H}om_{\mathcal{O}_M}\left( \Theta_{\mathcal{F}_0}, \frac{\Theta_M}{\Theta_{\mathcal{F}_0}} \right) \right)$ defined by
\begin{align} \label{tc213}
\overline{\Pi}_{i|\mu}:\Gamma\left( U_i, \Theta_{\mathcal{F}_0} \right) & \to \Gamma\left( U_i, \frac{\Theta_M}{\Theta_{\mathcal{F}_0}} \right) \\
 T_{0i}^\alpha & \mapsto \overline{\Pi_{i|\mu}^\alpha} \notag
\end{align}
and linearly extends to $\Gamma\left( U_i, \Theta_{\mathcal{F}_0} \right)$. Then $(\ref{tc32})-(\ref{tc34})$ implies that
\begin{align} \label{ttc2}
\left( \left\{  \overline{\Pi}_{i|\mu} \right\}      ,  \{  \Gamma_{ij|\mu} \} \right) \in C^0\left( \mathcal{U}, \mathscr{H}om_{\mathcal{O}_M}\left( \Theta_{\mathcal{F}}, \frac{\Theta_M}{\Theta_{\mathcal{F}_0}} \right)   \right) \bigoplus C^0\left( \mathcal{U}, \Theta_M  \right)
\end{align}
defines a $1$-cocycle in the following \v Cech resolution of $\Theta_{\mathcal{F}_0}^\bullet$:
{\small{\begin{center}
$\begin{CD}
\cdots \\
@AAA \\
C^0\left( \mathcal{U}, \mathscr{H}om_{\mathcal{O}_M}\left( \bigwedge^2 \Theta_{\mathcal{F}_0} , \frac{\Theta_M}{\Theta_{\mathcal{F}_0}}  \right) \right) @>-\delta >> \cdots \\
@AAA @AAA  \\
C^0\left(\mathcal{U}, \mathscr{H}om_{\mathcal{O}_M} \left(\Theta_{\mathcal{F}_0} , \frac{\Theta_M}{\Theta_{\mathcal{F}_0} } \right) \right) @>\delta>>  C^1\left(\mathcal{U}, \mathscr{H}om_{\mathcal{O}_M} \left(\Theta_{\mathcal{F}_0} , \frac{\Theta_M}{\Theta_{\mathcal{F}_0} } \right) \right) @>-\delta >> \cdots \\
@AAA @AAA @AAA \\
C^0(\mathcal{U}, \Theta_M) @>-\delta>> C^1(\mathcal{U}, \Theta_M) @>\delta>> C^2(\mathcal{U}, \Theta_M) @>-\delta>> \cdots
\end{CD}$
\end{center}}}
By the hypothesis that the foliated Kodaira-Spencer map $\varphi_0: T_0(B) \to \mathbb{H}^1\left( M, \Theta_{\mathcal{F}_0}^\bullet \right)$ is surjective, we can find homogeneous polynomial $s_\mu^\lambda$ such that
\begin{align*}
\varphi_0\left( \sum_{\lambda=1}^r s_\mu^\lambda \frac{\partial}{\partial t_\lambda}   \right) = \left(  \left\{  \overline{\Pi}_{i|\mu}  \right\} ,   \left\{ \Gamma_{ij|\mu}  \right\}   \right)
\end{align*}

Since we have
\begin{align*}
\varphi_0\left( \frac{\partial}{\partial t_\lambda} \right) =\left( \{\alpha_{i\lambda}\} =\left\{ T_{0i}^\alpha \mapsto \overline{ -\frac{\partial T_i^\alpha (z_i,t)}{\partial t_\lambda}|_{t=0} }\right\}, \{\rho_{ij\lambda} \} =  \left\{ \sum_{\alpha=1}^n \frac{\partial f_{ij}^\alpha}{\partial t_\lambda}|_{t=0} \frac{\partial}{\partial z_i^\alpha } \right\}   \right),
\end{align*}
there exists $\{\varphi_{i|\mu}\}\in C^0\left(\mathcal{U}, \Theta_M \right)$ and $\left\{ b_{i|\mu}\right\} \in C^0\left( \mathcal{U}, \mathscr{H}om_{\mathcal{O}_M}\left( \Theta_{\mathcal{F}_0}, \Theta_{\mathcal{F}_0} \right) \right)$ where $b_{i|\mu}\in \Gamma\left( U_i, \mathscr{H}om_{\mathcal{O}_M}\left( \Theta_{\mathcal{F}_0}, \Theta_{\mathcal{F}_0} \right) \right)$ defined by $b_{i|\mu}:\Gamma\left( U_i, \Theta_{\mathcal{F}_0} \right)\to \Gamma\left( U_i, \Theta_{\mathcal{F}_0} \right), T_{0i}^\alpha \mapsto \sum_{\gamma=1}^p b_{i|\mu}^{\alpha \gamma} T_{0i}'^\gamma$ such that
\begin{align*}
\varphi_{i|\mu}- \varphi_{j|\mu}&  = -\Gamma_{ij|\mu} + \sum_{\lambda=1}^m s_\mu^\lambda \rho_{ij\lambda}\\
\left[ \varphi_{i|\mu}, T_{0i}^\alpha   \right] &= - \Pi_{i|\mu}^\alpha + \sum_{\lambda=1}^m s_\mu^\lambda \alpha_{i\lambda}^\alpha + \sum_{\gamma=1}^p b_{i|\mu}^{\alpha \gamma} T_{0i}'^\gamma
\end{align*}
Then $(\ref{tc8})$ and $(\ref{tc10})$ holds. We check $(\ref{tc9})$. From $(\ref{tc33})$ and $(\ref{tc8})$ and $(\ref{tc10})$, we have
\begin{align*}
&\sum_{\beta=1}^p r_{0ij}^{\alpha \beta} \left( -\left[ \varphi_{j|\mu}, T_{0j}^\beta \right] + \sum_{\lambda=1}^m s_\mu^\lambda \alpha_{j\lambda}^\beta + \sum_{\gamma=1}^p b_{j|\mu}^{\beta\gamma} T_{0j}'^\gamma   \right) \\
 &=  -\left[ \varphi_{i|\mu}, T_{0i}^\alpha \right] + \sum_{\lambda=1}^m s_\mu^\lambda \alpha_{i\lambda}^\alpha + \sum_{\gamma=1}^p b_{i|\mu}^{\alpha\gamma} T_{0i}'^\gamma +\left[ \varphi_{i|\mu}, T_{0i}^\alpha \right] - \left[ \varphi_{j|\mu}, T_{0i}^\alpha \right] - \sum_{\lambda=1}^m s_\mu^\lambda \left[ \rho_{ij\lambda}, T_{0i}^\alpha   \right] + \sum_{\beta, \sigma, \gamma =1}^p B_{0ij|\mu}^{\alpha  \beta} b_{0i}^{\beta \sigma} r_{0ij}'^{\sigma \gamma} T_{0j}'^\gamma
\end{align*}
which is equivalent to
\begin{align}\label{tc42}
&- \left[ \varphi_{j|\mu}, T_{0i}^\alpha\right] + \sum_{\beta=1}^p \left[ \varphi_{j|\mu}, r_{0ij}^{\alpha \beta} \right] T_{0j}^\beta  + \sum_{\lambda=1}^m s_\mu^\lambda \left( \sum_{\beta=1}^p r_{0ij}^{\alpha \beta} \alpha_{j\lambda}^\beta - \alpha_{i\lambda}^\alpha + \left[ \rho_{ij\lambda}, T_{0i}^\alpha   \right]     \right) + \sum_{\beta,\gamma=1}^p r_{0ij}^{\alpha\beta} b_{j|\mu}^{\beta\gamma} T_{0j}'^\gamma \\
&= \sum_{\gamma=1}^p b_{i|\mu}^{\alpha \gamma} T_{0i}'^\gamma - \left[ \varphi_{j|\mu}, T_{0i}^\alpha \right] + \sum_{\beta, \sigma, \gamma =1}^p B_{0ij|\mu}^{\alpha  \beta} b_{0i}^{\beta \sigma} r_{0ij}'^{\sigma \gamma} T_{0j}'^\gamma \notag
\end{align}
We recall that
\begin{align}\label{tc43}
\sum_{\beta=1}^p r_{0ij}^{\alpha \beta} \alpha_{j\lambda}^\beta - \alpha_{i\lambda}^\alpha + \left[ \rho_{ij\lambda}, T_{0i}^\alpha   \right]    = \sum_{\beta=1}^p \frac{\partial r_{ij}^{\alpha\beta}}{\partial t_\lambda}|_{t=0} T_{0j}^\beta
\end{align}
Then from $(\ref{tc42})$ and $(\ref{tc43})$, we have
\begin{align*}
\sum_{\beta,\gamma=1}^p \left[ \varphi_{j|\mu} , r_{0ij}^{\alpha\beta} \right] b_{0j}^{\beta\gamma} T_{0j}'^\gamma +  \sum_{\lambda=1}^m \sum_{\beta,\gamma=1}^p s_\mu^\lambda \frac{\partial r_{0ij}^{\alpha\beta}}{\partial t_\lambda}|_{t=0} b_{0j}^{\beta\gamma} T_{0j}'^\gamma + \sum_{\beta,\gamma=1}^p r_{0ij}^{\alpha\beta} b_{j|\mu}^{\beta\gamma} T_{0j}'^\gamma = \sum_{\beta,\gamma=1}^p b_{i|\mu}^{\alpha\gamma} r_{0ij}'^{\beta\gamma} T_{0j}'^\gamma + \sum_{\beta, \sigma, \gamma =1}^p B_{0ij|\mu}^{\alpha  \beta} b_{0i}^{\beta \sigma} r_{0ij}'^{\sigma \gamma} T_{0j}'^\gamma
\end{align*}
This implies $(\ref{tc9})$. We check $(\ref{tc40})$. In fact, first we recall that
\begin{align}\label{tc45}
\left[ T_{0i}^\alpha, \alpha_{i\lambda}^\beta \right]- \left[ T_{0i}^\beta, \alpha_{i\lambda}^\alpha \right] - \sum_{\xi=1}^p g_{0i\alpha\beta}^\xi \alpha_{i\lambda}^\xi = - \sum_{\xi=1}^p \frac{\partial g_{i\alpha\beta}^\xi}{\partial t_\lambda}|_{t=0} T_{0i}^\xi
\end{align}
Then from $(\ref{tc10})$ and $(\ref{tc34})$ and $(\ref{tc45})$, we have
{\tiny{\begin{align*}
&- \sum_{\gamma=1}^p G_{i|\mu}^{\alpha \beta \gamma} T_{0i}'^\gamma \\
&=\left[ \sum_{\lambda=1}^m s_\mu^\lambda \alpha_{i\lambda}^\alpha + \sum_{\gamma=1}^p b_{i|\mu}^{\alpha \gamma} T_{0i}'^\gamma- \left[ \varphi_{i| \mu}, T_{0i}^\alpha \right]  , T_{0i}^\beta \right] -\left[ \sum_{\lambda=1}^m s_\mu^\lambda \alpha_{i\lambda}^\beta + \sum_{\gamma=1}^p b_{i|\mu}^{\beta \gamma} T_{0i}'^\gamma- \left[ \varphi_{i| \mu}, T_{0i}^\beta \right] , T_{0i}^\alpha \right]-\sum_{\xi=1}^p g_{0i\alpha\beta}^\xi \left( \sum_{\lambda=1}^m s_\mu^\lambda \alpha_{i\lambda}^\xi + \sum_{\gamma=1}^p b_{i|\mu}^{\xi \gamma} T_{0i}'^\gamma- \left[ \varphi_{i| \mu}, T_{0i}^\xi \right]  \right)\\
&=\sum_{\lambda=1}^m s_\mu^\lambda \left( \left[\alpha_{i\lambda}^\alpha, T_{0i}^\beta \right]  -\left[ \alpha_{i\lambda}^\beta, T_{0i}^\alpha \right] - \sum_{\xi=1}^p g_{0i\alpha\beta}^\xi \alpha_{i\lambda}^\xi   \right) -\sum_{\gamma=1}^p T_{0i}^\beta\left( b_{i|\mu}^{\alpha \gamma} \right) T_{0i}'^\gamma + \sum_{\gamma, \eta=1}^p b_{i|\mu}^{\alpha \gamma}\left[ T_{0i}'^\gamma, b_{0i}^{\beta \eta} T_{0i}'^\eta  \right] + \sum_{\gamma=1}^p T_{0i}^\alpha\left( b_{i|\mu}^{\beta \gamma} \right) T_{0i}'^\gamma - \sum_{\gamma, \eta=1}^p b_{i|\mu}^{\beta \gamma}\left[ T_{0i}'^\gamma, b_{0i}^{\alpha \eta} T_{0i}'^\eta \right] \\
&- \sum_{\xi, \gamma=1}^p g_{0i\alpha\beta}^\xi b_{i|\mu}^{\xi \gamma} T_{0i}'^\gamma - \left[ \left[ \varphi_{i|\mu}, T_{0i}^\alpha \right] , T_{0i}^\beta \right] + \left[ \left[ \varphi_{i|\mu}, T_{0i}^\beta \right], T_{0i}^\alpha \right] + \left[ \varphi_{i|\mu}, \sum_{\xi=1}^p g_{0i\alpha\beta}^\xi T_{0i}^\xi \right] - \sum_{\xi=1}^p \left[ \varphi_{i|\mu}, g_{0i\alpha\beta}^\xi \right] T_{0i}^\xi\\
&= - \sum_{\lambda=1}^m \sum_{\xi=1}^p s_\mu^\lambda \frac{\partial g_{i\alpha\beta}^\xi}{\partial t_\lambda}|_{t=0} T_{0i}^\xi - \sum_{\gamma=1}^p T_{0i}^\beta\left( b_{i|\mu}^{\alpha \gamma} \right) T_{0i}'^\gamma + \sum_{\gamma, \eta=1}^p b_{i|\mu}^{\alpha \gamma} T_{0i}'^\gamma\left( b_{0i}^{\beta \eta}  \right) T_{0i}'^\eta + \sum_{\gamma, \eta, \xi=1}^p b_{i|\mu}^{\alpha \gamma} b_{0i}^{\beta \eta} g_{0i\gamma \eta}'^\xi T_{0i}'^\xi + \sum_{\gamma=1}^p T_{0i}^\alpha \left( b_{i|\mu}^{\beta\gamma} \right) T_{0i}'^\gamma - \sum_{\gamma, \eta=1}^p b_{i|\mu}^{\beta \gamma} T_{0i}'^\gamma\left( b_{0i}^{\alpha \eta} \right) T_{0i}'^\eta \\
&- \sum_{\gamma,\eta, \xi=1}^p b_{i|\mu}^{\beta \gamma} b_{0i}^{\alpha \eta} g_{0i\gamma \eta}'^\xi T_{0i}'^\xi - \sum_{\xi, \gamma=1}^p g_{0i\alpha\beta}^\xi b_{i|\mu}^{\xi \gamma} T_{0i}'^\gamma - \sum_{\xi=1}^p\sum_{c=1}^n \varphi_{i|\mu}^c\frac{\partial  g_{0i\alpha\beta}^\xi}{\partial z_i^c} T_{0i}^\xi
\end{align*}}}
This implies $(\ref{tc40})$. This completes the inductive construction of $\varphi_i^\mu, s^\mu$ and $b_i^{\alpha\beta\mu}$. 

\begin{remark}\label{ttc3}
We recall from $(\ref{ttc1})$ that we have $\mathbb{H}^1\left( M, \Theta_{\mathcal{F}_0}^\bullet \right) \cong \mathbb{H}^1\left(  M, \mathcal{E}_{\Theta_{\mathcal{F}_0}}^\bullet \right)$. We will reinterpret $(\ref{ttc2})$ in $\mathbb{H}^1\left( M, \Theta_{\mathcal{F}_0}^\bullet \right)$ in terms of $\mathbb{H}^1\left( M, \mathcal{E}_{\Theta_{\mathcal{F}_0}}^\bullet \right)$ in the following \v Cech resolution of $\mathcal{E}_{\Theta_{\mathcal{F}_0}}^\bullet$:
{\small{\begin{equation}\label{ttc4}
\begin{CD}
\cdots \\
@AD_2 AA \\
C^0\left( \mathcal{U}, \mathscr{H}om_{\mathcal{O}_M}\left( \bigwedge^2 \Theta_{\mathcal{F}_0}, \frac{\Theta_M}{\Theta_{\mathcal{F}_0}} \right) \right)@>-\delta >> \cdots \\
@AD_1'AA @AD_1'AA\\
C^0\left(\mathcal{U}, \mathscr{H}om_{\mathcal{O}_M}\left( \Theta_{\mathcal{F}_0}, \Theta_M \right) \right) @>\delta>> C^1\left( \mathcal{U}, \mathscr{H}om_{\mathcal{O}_M}\left( \Theta_{\mathcal{F}_0}, \Theta_M \right) \right) @>-\delta>> \cdots \\
@AD_0'AA @AD_0'AA @AD_0'AA \\
C^0\left(\mathcal{U}, \mathcal{E}_{\Theta_{\mathcal{F}_0}} \right) @>-\delta >> C^1\left( \mathcal{U}, \mathcal{E}_{\Theta_{\mathcal{F}_0}} \right) @>\delta>> C^2\left( \mathcal{U}, \mathcal{E}_{\Theta_{\mathcal{F}_0}} \right) @>\delta>> \cdots
\end{CD}
\end{equation}}}
We set $B_{ij|\mu}\in \Gamma\left( U_{ij}, \mathscr{H}om_{\mathcal{O}_M}\left( \Theta_{\mathcal{F}_0}, \Theta_{\mathcal{F}_0} \right) \right)$ by
\begin{align}\label{ssc3}
B_{ij|\mu}:\Gamma\left( U_{ij}, \Theta_{\mathcal{F}_0} \right) &\to \Gamma\left( U_{ij}, \Theta_{\mathcal{F}_0}\right) \\
 T_{0i}^\alpha &\mapsto \sum_{\beta=1}^p B_{ij|\mu}^{\alpha \beta} T_{0i}^\beta \notag
\end{align}
and linearly extends to $\Gamma\left(U_{ij}, \Theta_{\mathcal{F}_0} \right)$. Then from $(\ref{tc33})$ we see that
\begin{align} \label{ttc16}
B_{ji}\left( T_{0i}^\alpha \right)& = \sum_{\beta=1}^p r_{0ij}^{\alpha\beta} B_{ji}\left( T_{0j}^\beta \right) =\sum_{\beta=1}^p r_{0ij}^{\alpha\beta} \left( \sum_{\eta=1}^p r_{0ji}^{\beta \eta} \Pi_{i|\mu}^\eta - \Pi_{j|\mu}^\beta + \left[ \Gamma_{ji}, T_{0j}^\beta \right] \right) \\
& = \Pi_{i|\mu}^\eta - \sum_{\beta=1}^p r_{0ij}^{\alpha\beta} \Pi_{j|\mu}^\beta  -\left[\Gamma_{ij|\mu}, T_{0i}^\alpha \right] + \sum_{\beta=1}^p \left[ \Gamma_{ij|\mu}, r_{0ij}^{\alpha\beta} \right] T_{0j}^\beta = - B_{ij|\mu}\left( T_{0i}^\alpha \right) + \sum_{\beta=1}^p \left[  \Gamma_{ij|\mu}, r_{0ij}^{\alpha\beta} \right] T_{0j}^\beta \notag
\end{align}
This implies that $\left\{ \left( \Gamma_{ij|\mu}, B_{ij|\mu} \right) \right\} \in C^1\left( \mathcal{U},  \mathcal{E}_{\Theta_{\mathcal{F}_0}} \right)$. On the other hand, from $(\ref{tc33})$ we also see that
\begin{align}\label{ttc15}
&\left(B_{ij|\mu}- B_{ik|\mu}+ B_{jk|\mu}\right)\left( T_{0i}^\alpha \right)\\
&= \sum_{\eta=1}^p  r_{0ij}^{\alpha \eta} \Pi_{j|\mu}^\eta - \Pi_{i|\mu}^\alpha + \left[ \Gamma_{ij|\mu}, T_{0i}^\alpha \right] - \sum_{\eta=1}^p  r_{0ik}^{\alpha \eta} \Pi_{k|\mu}^\eta + \Pi_{i|\mu}^\alpha -\left[ \Gamma_{ik|\mu}, T_{0i}^\alpha \right] + \sum_{\beta=1}^p r_{0ij}^{\alpha\beta}\left(\sum_{\eta=1}^p  r_{0jk}^{\alpha \eta} \Pi_{k|\mu}^\eta - \Pi_{j|\mu}^\alpha + \left[ \Gamma_{jk|\mu}, T_{0j}^\alpha \right] \right) \notag \\
&= \left[ \Gamma_{ij|\mu}- \Gamma_{ik|\mu}+ \Gamma_{jk|\mu}, T_{0i}^\alpha \right] - \sum_{\beta=1}^p \left[ \Gamma_{jk|\mu}, r_{0ij}^{\alpha\beta} \right] T_{0j}^\beta =  - \sum_{\beta=1}^p \left[ \Gamma_{jk|\mu}, r_{0ij}^{\alpha\beta} \right] T_{0j}^\beta \notag
\end{align}
This implies that $\delta\left(\left\{ \left( \Gamma_{ij|\mu}, B_{ij|\mu} \right) \right\}\right)=0$. Hence
\begin{align}
\left( \left\{ \Pi_{i|\mu} \right\}  , \left\{  \left(  \Gamma_{ij|\mu}, B_{ij|\mu}  \right)  \right\}  \right) \in C^0\left( \mathcal{U}, \mathscr{H}om_{\mathcal{O}_M}\left( \Theta_{\mathcal{F}_0}, \Theta_M \right) \right) \bigoplus C^1\left( \mathcal{U}, \mathcal{E}_{\Theta_{\mathcal{F}_0}} \right)
\end{align}
defines a $1$-cocycle in the above \v Cech resolution of $\mathcal{E}_{\Theta_{\mathcal{F}_0}}^\bullet$.

\end{remark}

\subsection{Proof of convergence}\

We prove that we can choose appropriate solutions $\varphi_{i|\mu}, s_\mu$ and $b_{i|\mu}^{\alpha\beta}$ satisfying $(\ref{tc8})-(\ref{tc10})$ in each inductive step so that
\begin{align*}
s(u)&=s_1(u)+ s_2(u)+\cdots + s_\mu(u)+ \cdots\\
\varphi_i\left( \xi_i, u\right)&= \xi_i+ \varphi_{i|1}\left(\xi_i, u\right) +  \varphi_{i|2}\left( \xi_i, u \right) + \cdots + \varphi_{i|\mu}(\xi_i, u) + \cdots \\
b_i^{\alpha\beta}(\xi_i, u)&= b_{0i}^{\alpha\beta}(\xi_i) + b_{i|1}^{\alpha\beta}(\xi_i,u) + \cdots + b_{i|\mu}^{\alpha\beta} (\xi_i, u) + \cdots
\end{align*}
converge absolutely and uniformly for $|u|< \epsilon$ if $\epsilon>0$ is sufficiently small. 

Consider a power series 
\begin{align*}
g(u)=\sum g_{v_1,..., v_{m'}} u_1^{v_1} u_2^{v_2} \cdots u_{m'}^{v_{m'}}
\end{align*}
whose coefficients $g_{v_1,..., v_{m'}}$ are vectors, and a power series
\begin{align*}
a(u)=\sum a_{v_1,..., v_{m'}} u_1^{v_1} u_2^{v_2} \cdots u_{m'}^{v_{m'}}
\end{align*}
with non-negative coefficients $a_{v_1\cdots v_{m'}}\geq 0$.
\begin{align}\label{ncc1}
\textnormal{We indicate by writing $g(u)\ll a(u)$ that $\left| g_{v_1,..., v_{m'}} \right| \leq a_{v_1\cdots v_{m'}}$}
\end{align}
We recall that
\begin{align}\label{ncc2}
A(u)= \frac{b}{16c} \sum_{\mu=1}^\infty \frac{1}{\mu^2}c^\mu\left( u_1 + u_2 + \cdots + u_{m'} \right)^\mu \,\,\,\textnormal{with $b,c>0$}\,\,\,\,\, \Longrightarrow A(u)^v \ll \left(\frac{b}{c} \right)^{v-1} A(u),\,\,\,\,\,v=2,3,\cdots
\end{align}

It suffices to prove the estimates $s(u)\ll A(u), \varphi_i(\xi_i,u)-\xi_i\ll A(u)$ and $b_i^{\alpha\beta}(\xi_i, u)- b_{0i}^{\alpha\beta}(\xi_i) \ll A(u)$ for suitable constants $b$ and $c$, equivalently
\begin{align}\label{tp1}
s^\mu(u)\ll A(u),\,\,\,\,\,\,\,\,\,\varphi_i^\mu(\xi_i, u)- \xi_i\ll A(u),\,\,\,\,\,\,\,\,\,\,\,b_i^{\alpha\beta \mu}(\xi_i, u)- b_{0i}^{\alpha\beta}(\xi_i) \ll A(u)
\end{align}
for $\mu=1,2,3,\cdots$. We will prove $(\ref{tp1})$ by induction on $\mu=1,2,3,\cdots$. For $\mu=1$, since the linear term of $A(u)$ is $\frac{b}{16}\left( u_1+ \cdots + u_{r'} \right)$, the estimates $\left(\ref{tp1} \right)_1$ holds if $b$ is sufficiently large. Let $\mu\geq  2$ and assume that the induction $\left(\ref{tp1}\right)_{\mu-1}$ holds for $\mu-1$, i.e.
\begin{align*}
s^{\mu-1}(u)\ll A(u),\,\,\,\,\,\,\,\,\,\varphi_i^{\mu-1}(\xi_i, u)- \xi_i \ll A(u),\,\,\,\,\,\,\,\,\,\,\,b_i^{\alpha\beta(\mu-1)}(\xi_i, u) - b_{0i}^{\alpha\beta}(\xi_i) \ll A(u)
\end{align*}
We will prove that $(\ref{tp1})_\mu$ holds. We estimate $(\ref{tc35}),(\ref{tc36})$ and $(\ref{tc38})$ in the following Lemma.
\begin{lemma}\label{tcc31}
\begin{align}
\Gamma_{ij|\mu}(\xi_i, u) &\ll \left( \frac{K_1}{b}+ \frac{K_2}{c}+ \frac{K_3 b}{c} \right)A(u) \,\,\,\,\,\,\,\,\textnormal{on}\,\,\, U_{ij} \label{tpp1}\\
B_{ij|\mu}^{\alpha\beta} (\xi_i, u) &\ll \left( \frac{K_{27}}{b} + \frac{K_{28}}{c} + \frac{K_{29} b}{c} \right) A(u) \,\,\,\,\,\,\,\,\,\,\textnormal{on}\,\,\, U_{ij} \label{tpp2} \\
\Pi_{i|\mu}^{\alpha\gamma} (\xi_i, u)  & \ll \left( \frac{K_{34} }{c} + \frac{K_{35} b}{c} \right) A(u) \,\,\,\,\,\,\,\textnormal{on}\,\,\, U_i^\delta
\end{align}
where $K_1,K_2,K_3, \cdots$ are constants independent of $\mu$.
\end{lemma}
\begin{proof}
The estimate $(\ref{tpp1})$ follows from \cite{Kod05} p.302. We prove the estimate $(\ref{tpp2})$. In fact, we set
\begin{align}\label{tcc1}
Q_{ij|\mu}^{\alpha \gamma}:= \left[ \sum_{\beta=1}^p b_i^{\alpha\beta (\mu-1)}\left(f_{ij}'(\xi_j, u), u\right) r_{ij}'^{\beta\gamma}(\xi_j, u) \right]_\mu - \left[\sum_{\beta=1}^p r_{ij}^{\alpha\beta} \left(\varphi_j^{\mu-1}(\xi_j, u), s^{\mu-1}(u) \right)b_j^{\beta\gamma(\mu-1)}(\xi_j, u)\right]_\mu
\end{align}
We estimate $Q_{ij|\mu}^{\alpha \gamma}$. We estimate the first part of $Q_{ij|\mu}^{\alpha \gamma}$ $(\ref{tcc1})$:
\begin{align}
&\left[ \sum_{\beta=1}^p b_i^{\alpha\beta (\mu-1)}\left(f_{ij}'(\xi_j, u), u\right) r_{ij}'^{\beta\gamma}(\xi_j, u)  \right]_\mu =\left[ \sum_{\beta=1}^p b_i^{\alpha\beta (\mu-1)}\left(f_{ij}'(\xi_j, u), u\right) \left( r_{ij}'^{\beta \gamma}(\xi_j,u)- r_{0ij}'^{\beta \gamma}(\xi_j) + r_{0ij}'^{\beta \gamma}(\xi_j)     \right)\right]_\mu \notag \\
&=\left[ \sum_{\beta=1}^p b_i^{\alpha\beta (\mu-1)}\left(f_{ij}'(\xi_j, u), u\right) \left( r_{ij}'^{\beta \gamma}(\xi_j,u)- r_{0ij}'^{\beta \gamma}(\xi_j) \right)  \right]_\mu + \left[ \sum_{\beta=1}^p b_i^{\alpha\beta (\mu-1)}\left(f_{ij}'(\xi_j, u), u\right) r_{0ij}'^{\beta \gamma}(\xi_j)\right]_\mu \label{tcc5}
\end{align}
First we estimate 
\begin{align*}
\left[ b_i^{\alpha\beta(\mu-1)} \left( f_{ij}'(\xi_j, u), u \right)\right]_\mu
\end{align*}
In a similar method with \cite{Kod05} p.298-301, we put
\begin{align*}
G(z_i,u)= b_i^{\alpha\beta(\mu-1)}(\xi_i, u)- b_{0i}^{\alpha\beta}(\xi_i)
\end{align*}
then by the induction hypothesis, we have $ G(z_i,u )\ll A(u)$. We recall that $\mathcal{U}=\{U_i\}$ is a finite covering of $M$ where $U_i=\left\{|z_i|\in \mathbb{C}^n |z_i|<1 \right\}$, and $M$ is compact, so that we have 
\begin{align}\label{ncc10}
M= \bigcup_i  U_i^\delta,\,\,\,\,\,\,\,\,\textnormal{where}\,\,\,U_i^\delta=\left\{ z_i\in U_i | |z_i|<1-\delta \right\}\,\,\,\,\,\,\textnormal{for a sufficiently small $\delta$}
\end{align}
Then we have in a similar way with \cite{Kod05} p.298-301
\begin{align}
G\left( f_{ij}(\xi_j,u) , u \right)- G(\xi_i, u) \ll \frac{2^{n+1} b_0}{c\delta} A(u),\,\,\,\,\,b_0>0 \,\,\,\,\,(b> b_0) \,\,\,\,\,\,\textnormal{on}\,\,\,U_i^\delta \cap U_j
\end{align}
namely,
\begin{align}\label{tcc2}
b_i^{\alpha\beta(\mu-1)}\left( f_{ij}(\xi_j, u), u \right) - b_{0i}^{\alpha\beta}\left( f_{ij}(\xi_j, u) \right) - b_i^{\alpha\beta(\mu-1)}(\xi_i, u) + b_{0i}^{\alpha\beta}(\xi_i) \ll \frac{2^{n+1} b_0}{c \delta} A(u)\,\,\,\,\,\,\textnormal{on}\,\,\,U_i^\delta \cap U_j
\end{align}
Since $\tilde{b}_{ij}^{\alpha\beta}(\xi_j, u):=b_{0i}^{\alpha\beta}\left( f_{ij}(\xi_j, u) \right)= b_{0i}^{\alpha\beta}(\xi_i) + \sum_{\mu=1}^\infty \tilde{b}_{ij|\mu}^{\alpha\beta}(\xi_j, u)$ is a holomorphic function, we may assume that
\begin{align}\label{tcc3}
b_{0i}^{\alpha\beta}\left( f_{ij}(\xi_j, u) \right)- b_{0i}^{\alpha\beta}(\xi_i)\ll A_1(s) \ll \frac{b_1}{b}A(u),\,\,\,\,\,\,\,A_1(u)=\frac{b_1}{16 c_1}\sum_{\mu=1}^\infty \frac{c_1^\mu (u_1+ \cdots + u_{m'})}{\mu^2}
\end{align}
holds for $\xi_j\in U_i\cap U_j$ with $b_1 >0$ and $c_1 >0$ and $c>c_1$. Then from $(\ref{tcc2})$ and $(\ref{tcc3})$ we have
\begin{align}\label{tcc4}
b_i^{\alpha\beta(\mu-1)}\left( f_{ij}\left( \xi_j, u\right), u \right) - b_i^{\alpha\beta(\mu-1)}(\xi_i, u) \ll \left(\frac{2^{n+1} b_0}{c\delta} + \frac{b_1}{b} \right) A(u)\,\,\,\,\,\,\,\,\textnormal{on}\,\,\, U_i^\delta \cap U_j
\end{align}
By taking the terms of degree $\mu$ from $(\ref{tcc4})$ and setting $K_4=\frac{2^{n+1}b_0}{\delta}$ and $b_1=K_5$, we have
\begin{align}
\left[ b_i^{\alpha\beta(\mu-1)} \left( f_{ij}'(\xi_j, u), u \right)\right]_\mu \ll \left( \frac{K_4}{c}  + \frac{K_5}{b}\right) A(u)\,\,\,\,\,\,\,\,\textnormal{on $U_i^\delta \cap U_j$}
\end{align}
Then the second term of $(\ref{tcc5})$ is estimated by
\begin{align} \label{tcc9}
 \left[ \sum_{\beta=1}^p b_i^{\alpha\beta (\mu-1)}\left(f_{ij}'(\xi_j, u), u\right) r_{0ij}'^{\beta \gamma}(\xi_j)\right]_\mu \ll K_6 \left( \frac{K_4}{c}  + \frac{K_5}{b} \right) A(u)\,\,\,\,\,\,\,\textnormal{on}\,\,\,U_i^\delta \cap U_j
\end{align}

We may assume $r_{ij}'^{\beta \gamma}(\xi_j, u)- r_{0ij}'^{\beta \gamma}(\xi_j)$,  as in \cite{Kod05} p.298.
\begin{align} \label{tcc6}
r_{ij}'^{\beta \gamma}(\xi_j, u)- r_{0ij}'^{\beta \gamma}(\xi_j)  \ll A_2(u)\ll \frac{b_2}{b}A(u) \,\,\,\,\,\,\,\,\,A_2(u)=\frac{b_2}{16 c_2}\sum_{\mu=1}^\infty \frac{1}{\mu^2} c_2^\mu\left( u_1+ \cdots + u_{m'} \right)^\mu,\,\,\,\,\,\,\,\, b_2>0,\,\,\, c_2 >0 
\end{align}
holds for $z_j\in U_i\cap U_j$. Then we estimate the first term of $(\ref{tcc5})$. We note that 
{\small{\begin{align}\label{tcc7}
& b_i^{\alpha\beta (\mu-1)}\left(f_{ij}'(\xi_j, u), u\right) \left( r_{ij}'^{\beta \gamma}(\xi_j,u)- r_{0ij}'^{\beta \gamma}(\xi_j) \right) \\
&=  \left(b_i^{\alpha\beta (\mu-1)}\left(f_{ij}'(\xi_j, u), u\right)- b_i^{\alpha\beta(\mu-1)}(\xi_i,u) + b_i^{\alpha\beta(\mu-1)}(\xi_i, u) \right) \left( r_{ij}'^{\beta \gamma}(\xi_j,u)- r_{0ij}'^{\beta \gamma}(\xi_j) \right) \notag \\
&=\left(b_i^{\alpha\beta (\mu-1)}\left(f_{ij}'(\xi_j, u), u\right)- b_i^{\alpha\beta(\mu-1)}(\xi_i,u) \right)\left(r_{ij}'^{\beta \gamma}(\xi_j,u)- r_{0ij}'^{\beta \gamma}(\xi_j) \right) + \left(b_i^{\alpha \beta(\mu-1)}(\xi_i, u) - b_{0i}^{\alpha\beta}(\xi_i) \right)\left( r_{ij}'^{\beta \gamma}(\xi_j,u)- r_{0ij}'^{\beta \gamma}(\xi_j)\right) \notag \\
&+ b_{0i}^{\alpha\beta}(\xi_i)\left( r_{ij}'^{\beta \gamma}(\xi_j,u)- r_{0ij}'^{\beta \gamma}(\xi_j) \notag\right)
\end{align}}}
Then from $(\ref{tcc7})$ and $(\ref{tcc4})$ and $(\ref{tcc6})$ and the induction hypothesis, we have
\begin{align}
\left[  b_i^{\alpha\beta (\mu-1)}\left(f_{ij}'(\xi_j, u), u\right) \left( r_{ij}'^{\beta \gamma}(\xi_j,u)- r_{0ij}'^{\beta \gamma}(\xi_j) \right)  \right]_\mu \ll  \left( \frac{K_4}{c} + \frac{K_5}{b}\right)\frac{b_2}{b}A(u)^2 + \frac{b_2}{b} A(u)^2 + K_6 \frac{b_2}{b}A(u) \\
\ll \left(  \left( \frac{K_4}{c} + \frac{K_5}{b} \right) \frac{b_2}{c} + \frac{b_2}{c} + K_6 \frac{b_2}{b}        \right) A(u) \ll \left( \frac{K_7}{c} + \frac{K_8}{b} \right) A(u) \,\,\,\,\,\,\,\textnormal{on}\,\,\, U_i^\delta \cap U_j \notag
\end{align}
Then the first-term of $(\ref{tcc5})$ is estimated by
\begin{align}\label{tcc8}
\left[  \sum_{\beta=1}^p b_i^{\alpha\beta (\mu-1)}\left(f_{ij}'(\xi_j, u), u\right) \left( r_{ij}'^{\beta \gamma}(\xi_j,u)- r_{0ij}'^{\beta \gamma}(\xi_j) \right)  \right]_\mu \ll K_9\left(\frac{K_7}{c}+ \frac{K_8}{b} \right) A(u) \,\,\,\,\,\,\,\textnormal{on}\,\,\,U_i^\delta \cap U_j
\end{align}
Then from $(\ref{tcc9})$ and $(\ref{tcc8})$ the first part of $Q_{ij|\mu}^{\alpha \gamma}$ $(\ref{tcc1})$ is estimated by
\begin{align}\label{tcc15}
\left[ \sum_{\beta=1}^p b_i^{\alpha\beta(\mu-1)} \left( f_{ij}(\xi_j, u), u \right) r_{ij}'^{\beta \gamma}(\xi_j, u) \right]_\mu \ll   \left( \frac{K_{10}}{c} + \frac{ K_{11}}{b} \right) A(u)
\end{align}

We estimate the second part of $Q_{ij|\mu}^{\alpha \gamma}$ $(\ref{tcc1})$.
\begin{align}
&\left[\sum_{\beta=1}^p r_{ij}^{\alpha\beta} \left(\varphi_j^{\mu-1}(\xi_j, u), s^{\mu-1}(u) \right)b_j^{\beta\gamma(\mu-1)}(\xi_j, u)\right]_\mu \notag\\
&=\left[ \sum_{\beta=1}^p \left( r_{ij}^{\alpha\beta}\left( \varphi_j^{\mu-1} (\xi_j, u), s^{\mu-1}(u)\right)  - r_{0ij}^{\alpha\beta} + r_{0ij}^{\alpha\beta} \right) \left( b_j^{\beta \gamma(\mu-1)}(\xi_j, u)- b_{0j}^{\beta \gamma}    + b_{0j}^{\beta \gamma}  \right)  \right]_\mu \notag\\
&=\left[ \sum_{\beta=1}^p \left( r_{ij}^{\alpha\beta}\left( \varphi_j^{\mu-1} (\xi_j, u), s^{\mu-1}(u)\right)  - r_{0ij}^{\alpha\beta} \right) \left( b_j^{\beta \gamma(\mu-1)}(\xi_j, u)- b_{0j}^{\beta \gamma} \right)    \right]_\mu + \left[  \sum_{\beta=1}^p r_{ij}^{\alpha\beta}\left( \varphi_j^{\mu-1}(\xi_j, u), s^{\mu-1}(u) \right) b_{0j}^{\beta \gamma} \right]_\mu \label{tcc11}
\end{align}
We expand $r_{ij}^{\alpha\beta}(\xi_j + y,t)$ into power series in $y_1,..., y_n, t_1,..., t_{m}$. Then we may assume that
\begin{align}
r_{ij}^{\alpha\beta}(\xi_j + y, t)- r_{0ij}^{\alpha\beta}(\xi_j)\ll \sum_{v=1}^\infty a_1^v\left( y_1+ \cdots + y_n + t_1+ \cdots + t_m \right)^v,\,\,\,\,\,\,\,a_1>0
\end{align}
If we set $y=\varphi_j^{\mu-1}(\xi_j, u)- \xi_j$ and $t= s^{\mu-1}(u)$, then since $y\ll A(u)$ and $t\ll A(u)$ by the induction hypothesis, we obtain
\begin{align}\label{tcc10}
r_{ij}^{\alpha\beta}\left( \varphi_j^{\mu-1}, s^{\mu-1}(u) \right) - r_{0ij}^{\alpha\beta}(\xi_j) \ll \sum_{v=1}^\infty a_1^v (n+m)^v A(u)^v 
\end{align}
and by induction hypothesis $b_j^{\beta \gamma(\mu-1  )}(\xi_j, u) - b_{0j}^{\beta \gamma} \ll A(u)$, we have from $(\ref{tcc10})$
\begin{align}
\left(r_{ij}^{\alpha\beta}\left( \varphi_j^{\mu-1} (\xi_j, u), s^{\mu-1}(u)\right)  - r_{0ij}^{\alpha\beta} \right) \left( b_j^{\beta \gamma(\mu-1)}(\xi_j, u)- b_{0j}^{\beta \gamma} \right) & \ll \sum_{v=1}^\infty a_1^v (n+m)^v A(u)^{v+1}  \ll \sum_{v=1}^\infty a_1^v(n+m)^v \left( \frac{b}{c} \right)^v A(u)\\
&\ll  \frac{ba_1(m+n)}{c} \sum_{v=0}^\infty  \left( \frac{b a_1(m+n)}{c}\right)^v A(u)  \notag
\end{align}
By taking a constant $c$ such that
\begin{align}
\frac{b a_1 (m+n) }{c} < \frac{1}{2}
\end{align}
We estimate the first part of $(\ref{tcc11})$ by
\begin{align}\label{tcc12}
\left[ \sum_{\beta=1}^p \left( r_{ij}^{\alpha\beta}\left( \varphi_j^{\mu-1} (\xi_j, u), s^{\mu-1}(u)\right)  - r_{0ij}^{\alpha\beta} \right) \left( b_j^{\beta \gamma(\mu-1)}(\xi_j, u)- b_{0j}^{\beta \gamma} \right)    \right]_\mu \ll  K_{12}\left( \frac{2 b a_1 (m +n )}{c} \right) A(u)
\end{align}
On the other hand, in a similar way with \cite{Kod05} p. 301, we can show that if we take $c$ such that $\frac{ba_0(m+n)}{c}<\frac{1}{2}$ for some constant $a_0>0$,
\begin{align}
\left[ r_{ij}^{\alpha\beta}\left( \varphi^{\mu-1}(\xi_j, u), s^{\mu-1}(u) \right)   \right]_\mu \ll \frac{2 ba_0^2(m+n)^2}{c} A(u)
\end{align}
Then we have
\begin{align}\label{tcc13}
 \left[  \sum_{\beta=1}^p r_{ij}^{\alpha\beta}\left( \varphi_j^{\mu-1}(\xi_j, u), s^{\mu-1}(u) \right) b_{0j}^{\beta \gamma} \right]_\mu  \ll K_{14}\left( \frac{2 ba_0^2(m+n)^2}{c}  \right) A(u)
\end{align}
Then from $(\ref{tcc12})$ and $(\ref{tcc13})$,  the second part of $Q_{ij|\mu}^{\alpha \gamma}$ $(\ref{tcc1})$ is estimated by
\begin{align}\label{tcc14}
\left[\sum_{\beta=1}^p r_{ij}^{\alpha\beta} \left(\varphi_j^{\mu-1}(\xi_j, u), s^{\mu-1}(u) \right)b_j^{\beta\gamma(\mu-1)}(\xi_j, u)\right]_\mu  \ll K_{15}\frac{b}{c} A(u)
\end{align}
Then from $(\ref{tcc15})$ and $(\ref{tcc14})$, $Q_{ij|\mu}^{\alpha \gamma}$ $(\ref{tcc1})$ is estimated by
{\small{\begin{align}
Q_{ij|\mu}^{\alpha \gamma} =  \left[ \sum_{\beta=1}^p b_i^{\alpha\beta (\mu-1)}\left(f_{ij}'(\xi_j, u), u\right) r_{ij}'^{\beta\gamma}(\xi_j, u)  - \sum_{\beta=1}^p r_{ij}^{\alpha\beta} \left(\varphi_j^{\mu-1}(\xi_j, u), s^{\mu-1}(u) \right)b_j^{\beta\gamma(\mu-1)}(\xi_j, u)\right]_\mu \ll \left(  \frac{K_{16}}{b} + \frac{K_{17}}{c} + \frac{K_{15} b}{c} \right) A(u)
\end{align}}}
on $U_i^\delta \cap U_j$. Since $\sum_{\beta, \sigma=1}^p B_{ij|\mu}^{\alpha\beta} b_{0i}^{\beta \sigma} r_{0ij}'^{\sigma \gamma}= Q_{ij|\mu}^{\alpha \gamma}$ from $(\ref{tc36})$, by considering inverse matrix of $b_{0i}^{\beta \sigma}$ and $r_{0ij}'^{\sigma \gamma}$, we have
\begin{align}\label{ttc18}
B_{ij|\mu}^{\alpha \beta} \ll \left(  \frac{K_{18}}{b} + \frac{K_{19}}{c} + \frac{K_{20} b}{c} \right) A(u)\,\,\,\,\,\,\,\textnormal{on}\,\,\, U_i^\delta \cap U_j
\end{align}

Now we estimate $B_{ij|\mu}(\xi_i, u)$ for arbitrary $\xi_i\in U_i\cap U_j$. First we note that $B_{ji}\left(T_{0i}^\alpha\right)= - B_{ij|\mu}\left( T_{0i}^\alpha \right) + \sum_{\beta=1}^p \left[ \Gamma_{ij|\mu}, r_{0ij}^{\alpha\beta} \right] T_{0j}^\beta$ from $(\ref{ttc16})$, equivalently,
\begin{align}
&\sum_{\beta, \gamma=1}^p r_{0ij}^{\alpha\beta} B_{ji|\mu}^{\beta \gamma} T_{0j}^\gamma = - \sum_{\beta=1}^p B_{ij|\mu}^{\alpha\beta} T_{0i}^\beta + \sum_{\gamma=1}^p \left[ \Gamma_{ij|\mu}, r_{0ij}^{\alpha \gamma} \right] T_{0j}^\gamma \notag \\
&\Longrightarrow  \sum_{\beta=1}^p r_{0ij}^{\alpha\beta} B_{ji|\mu}^{\beta \gamma} = - \sum_{\beta=1}^p B_{ij|\mu}^{\alpha\beta} r_{0ij}^{\beta \gamma} + \left[ \Gamma_{ij|\mu}, r_{0ij}^{\alpha\gamma} \right] \label{ttc17}
\end{align}
Then from $(\ref{ttc17})$ and $(\ref{tpp1})$ and $(\ref{ttc18})$ and setting $K'= \frac{K_{1}}{b} + \frac{K_{2}}{c} + \frac{K_{3}b}{c}$ and $K''= \frac{K_{18}}{b} + \frac{K_{19}}{c} + \frac{K_{20}b}{c}$,
\begin{align}\label{ttc19}
B_{ji|\mu}^{\alpha\beta} \ll K_{23} \left( K_{21} K''  + K_{22} K'\right) A(u) \ll \left(\frac{K_{24}}{b} + \frac{K_{25}}{c} + \frac{K_{26}b}{c}\right)A(u) \,\,\,\,\,\,\,\,\,\textnormal{on}\,\,\,\,\,U_i^\delta \cap U_j
\end{align}

We also note that $B_{ij|\mu}^\alpha- B_{ik|\mu}^\alpha + \sum_{\beta=1}^p r_{0ij}^{\alpha\beta}B_{jk|\mu}^\beta+ \sum_{\beta=1}^p\left[ \Gamma_{ij|\mu}, r_{0ij}^{\alpha\beta} \right]T_{0j}^\beta=0$ from $(\ref{ttc15})$, equivalently
\begin{align}
& \sum_{\beta=1}^p B_{ij|\mu}^{\alpha\beta} T_{0i}^\beta - \sum_{\beta=1}^p B_{ik|\mu}^{\alpha\beta} T_{0i}^\beta + \sum_{\beta,\gamma=1}^p r_{0ij}^{\alpha\beta} B_{jk}^{\beta \gamma} T_{0j}^\gamma = - \sum_{\gamma=1}^p \left[\Gamma_{jk|\mu}, r_{0ij}^{\alpha\gamma} \right] T_{0j}^\gamma \notag \\
& \Longrightarrow \sum_{\beta=1}^p B_{ij|\mu}^{\alpha \beta} r_{0ij}^{\beta \gamma} - \sum_{\beta=1}^p B_{ik|\mu}^{\alpha\beta} r_{0ij}^{\beta \gamma} + \sum_{\beta=1}^p r_{0ij}^{\alpha\beta} B_{jk}^{\beta\gamma}= -\left[ \Gamma_{jk|\mu}, r_{0ij}^{\alpha \gamma} \right] \label{ttc20}
\end{align}

 If $\xi_i \in U_i\cap U_j$ and $ \xi_i \notin U_i^\delta$, then $\xi_i\in U_k^\delta$ for some $k\ne i$. Since $\xi_i \in U_i\cap U_j\cap U_k^\delta$, we have $B_{jk|\mu}^{\alpha \gamma} \ll \tilde{K}A(u)$ and $B_{ik|\mu}^{\alpha \gamma}\ll \tilde{K}A(u)$ where $\tilde{K}=\left( \frac{K_{24}}{b} + \frac{K_{25}}{c} + \frac{K_{26}b}{c} \right)$ from $(\ref{ttc19})$. Then from $(\ref{ttc20})$ and $(\ref{tpp1})$, it follows that
 \begin{align}
 B_{ij|\mu}^{\alpha \beta} \ll \left( \frac{K_{27}}{b} + \frac{K_{28}}{c} + \frac{K_{29} b}{c} \right) A(u)\,\,\,\,\,\,\,\,\textnormal{on}\,\,\, U_i \cap U_j
 \end{align}

It remains to estimate $\Pi_{i|\mu}^{\alpha \gamma}$ on $U_i^\delta$ from $(\ref{tc38})$
\begin{align}\label{tpp3}
\left[T_i^{\alpha\gamma}\left(\varphi_i^{\mu-1}, s^{\mu-1}(u) \right)- \sum_{\beta=1}^p \sum_{\sigma=1}^n b_i^{\alpha\beta (\mu-1)}(\xi_i, u) T_i'^{\beta\sigma}(\xi_i, u)\frac{\partial \varphi_i^{\gamma(\mu-1)}(\xi_i, u)}{\partial \xi_i^\sigma} \right]_\mu 
\end{align}
First we estimate the first term of $(\ref{tpp3})$
\begin{align*}
\left[T_i^{\alpha\gamma}\left(\varphi_i^{\mu-1}, s^{\mu-1}(u) \right)\right]_\mu
\end{align*}
In a similar way with \cite{Kod05} p.301, we can show that if we take $c$ such that $\frac{ba_3(m+n)}{c}< \frac{1}{2}$ for some constant $a_3>0$,
\begin{align}\label{tcc26}
\left[T_i^{\alpha\gamma}\left(\varphi_i^{\mu-1}, s^{\mu-1}(u) \right)\right]_\mu \ll \frac{2ba_3^2 (m+n)^2}{c}A(u)
\end{align}

We estimate the second part of $(\ref{tpp3})$. We note that
{\small{\begin{align} \label{tcc20}
&\left[  \sum_{\beta=1}^p \sum_{\sigma=1}^n b_i^{\alpha\beta (\mu-1)}(\xi_i, u) T_i'^{\beta\sigma}(\xi_i, u)\frac{\partial \varphi_i^{\gamma(\mu-1)}(\xi_i, u)}{\partial \xi_i^\sigma} \right]_\mu\\
&=  \left[  \sum_{\beta=1}^p \sum_{\sigma=1}^n \left( b_i^{\alpha\beta (\mu-1)}(\xi_i, u) - b_{0i}^{\alpha \beta} + b_{0i}^{\alpha \beta} \right)\left(T_i'^{\beta\sigma}(\xi_i, u) - T_{0i}'^{\beta \sigma}(\xi_i) + T_{0i}'^{\beta \sigma}(\xi_i) \right)\frac{\partial \left( \varphi_i^{\gamma(\mu-1)}(\xi_i, u) - \xi_i^\gamma + \xi_i^\gamma \right)}{\partial \xi_i^\sigma} \right]_\mu \notag \\
&=\left[  \sum_{\beta=1}^p \sum_{\sigma=1}^n \left( b_i^{\alpha\beta (\mu-1)}(\xi_i, u) - b_{0i}^{\alpha \beta} \right) \left(T_i'^{\beta\sigma}(\xi_i, u)- T_{0i}'^{\beta \sigma} \right) \frac{\partial \left( \varphi_i^{\gamma(\mu-1)}(\xi_i, u)- \xi_i^\gamma \right) }{\partial \xi_i^\sigma}  \right]_\mu\notag \\
&+\left[  \sum_{\beta=1}^p \sum_{\sigma=1}^n b_{0i}^{\alpha \beta}  \left(T_i'^{\beta\sigma}(\xi_i, u)- T_{0i}'^{\beta \sigma} \right) \frac{\partial \left( \varphi_i^{\gamma(\mu-1)}(\xi_i, u)- \xi_i^\gamma \right) }{\partial \xi_i^\sigma}  \right]_\mu + \left[  \sum_{\beta=1}^p \sum_{\sigma=1}^n \left( b_i^{\alpha\beta (\mu-1)}(\xi_i, u) - b_{0i}^{\alpha \beta}\right) T_{0i}'^{\beta\sigma}(\xi_i)\frac{\partial \left( \varphi_i^{\gamma(\mu-1)}(\xi_i, u)- \xi_i^\gamma \right)}{\partial \xi_i^\sigma} \right]_\mu \notag\\
&+\left[ \sum_{\beta=1}^p \left( b_i^{\alpha\beta (\mu-1)}(\xi_i, u)- b_{0i}^{\alpha\beta} \right)\left( T_i'^{\beta \gamma}(\xi_i, u)- T_{0i}'^{\beta \gamma} \right)\right]_\mu \notag
\end{align}}}

First we estimate $\frac{\partial \left(\varphi_i^{\gamma(\mu-1)}(\xi_i, u)- \xi_i^\gamma \right)}{\partial \xi_i^\sigma}$. By induction hypothesis $g_i^\gamma:=\varphi_i^{\gamma (\mu-1)}(\xi_i, u)-\xi_i^\gamma \ll A(u)$. Assume that $\xi_i= \left(\xi_i^1,..., \xi_i^n \right)\in U_i^\delta$. Then
\begin{align*}
\frac{\partial \left( \varphi_i^{\gamma(\mu-1)}(\xi_i, u)- \xi_i^\gamma \right)}{\partial \xi_i^\sigma}=\frac{\partial g_i^{\gamma(\mu-1)}(\xi_i, u)}{\partial \xi_i^\sigma}=\frac{1}{2\pi i} \int_{|z- \xi_i^\gamma |=\delta} \frac{g_i^{\gamma(\mu-1)}(\xi_i^1,..., \overbrace{z}^{\gamma-\textnormal{th}}, ..., \xi_i^n, u)}{(z- \xi_i^\gamma)^2}dz
\end{align*}
so that for $\sigma=1,...,n$, we have
\begin{align}\label{tcc21}
\frac{\partial \varphi_i^{\gamma(\mu-1)}(\xi_i, u)- \xi_i^\gamma }{\partial \xi_i^\sigma} \ll \frac{A(u)}{\delta}
\end{align}
and we may assume that
\begin{align}\label{tcc22}
T_i'^{\beta \gamma} (\xi_i, u)- T_{0i}'^{\beta \gamma}\ll A_4(u) \ll \frac{b_4}{b} A(u),\,\,\,\,\,\,\,\,\,\,A_4(u)=\frac{b_4}{16c_4} \sum_{v=1}^n \frac{c_4^v \left( u_1 + \cdots + u_{m'} \right)^v}{v^2}
\end{align}
holds for $\xi_i \in U_i$ with $b_4>0$ and $c_4>0$ and $c>c_4$. Then from $(\ref{tcc20})$ and $(\ref{tcc21})$ and $(\ref{tcc22})$ and the induction hypothesis, we have
\begin{align}\label{tcc25}
&\left[  \sum_{\beta=1}^p \sum_{\sigma=1}^n b_i^{\alpha\beta (\mu-1)}(\xi_i, u) T_i'^{\beta\sigma}(\xi_i, u)\frac{\partial \varphi_i^{\gamma(\mu-1)}(\xi_i, u)}{\partial \xi_i^\sigma} \right]_\mu\\
&\ll K_{30}\frac{b_4}{b} \frac{1}{\delta}A(u)^3 + K_{31}\frac{b_4}{b}\frac{1}{\delta} A(u)^2 + K_{32} \frac{1}{\delta} A(u)^2 + K_{33}\frac{b_4}{b} A(u)^2 \notag \\
&\ll\left( K_{30} \frac{b_4}{b} \frac{1}{\delta}\frac{b^2}{c^2} + K_{31} \frac{b_4}{b}\frac{1}{\delta} \frac{b}{c}    + K_{32}\frac{1}{\delta} \frac{b}{c} + K_{33} \frac{b_4}{b}\frac{b}{c}  \right) A(u) \ll \left( \frac{K_{34}}{c} + \frac{K_{35} b}{c} \right)A(u) \,\,\,\,\,\,\,\,\textnormal{on}\,\,\,\,\,U_i^\delta \notag
\end{align}
Then from $(\ref{tcc26})$ and $(\ref{tcc25})$ we have
\begin{align}
\Pi_{i|\mu}^{\alpha \gamma}\ll \left( \frac{K_{34}}{c}+ \frac{K_{35}b}{c} \right) A(u)\,\,\,\,\,\,\,\,\textnormal{on}\,\,\, U_i^\delta
\end{align}
This completes the proof of Lemma \ref{tcc31}.

\end{proof}

We recall Remark \ref{ttc3}. For any $\sigma=\left(\Pi, \Gamma, B \right)=\left( \left\{ \Pi_i \right\},  \left\{ \left( \Gamma_{ij}, B_{ij} \right\} \right) \right\} \in C^0\left( \mathcal{U}, \mathscr{H}om_{\mathcal{O}_M}\left( \Theta_{\mathcal{F}_0}, \Theta_M \right) \right) \bigoplus C^1\left( \mathcal{U}, \mathcal{E}_{\Theta_{\mathcal{F}_0}} \right)$ which is a $1$-cocycle in the \v Cech resolution $(\ref{ttc4})$ of $\mathcal{E}_{\Theta_{\mathcal{F}_0}}^\bullet$, where $\Pi_i\in \Gamma\left(U_i, \mathscr{H}om_{\mathcal{O}_M}\left( \Theta_{\mathcal{F}_0}, \Theta_M \right) \right)$ defined by
\begin{align}\label{ssc9}
\Pi_i : \Gamma\left( U_i, \Theta_{\mathcal{F}_0} \right) &\to \Gamma\left( U_i, \Theta_M \right) \\
 T_{0i}^\alpha &\mapsto \Pi_i^\alpha:= \sum_{\gamma=1}^n \Pi_i^{\alpha \gamma} \frac{\partial}{\partial z_i^\gamma}  \notag
\end{align}
and $\Gamma_{ij}\in \Gamma\left( U_{ij}, \Theta_M \right)$ and $B_{ij}\in \Gamma\left( U_{ij}, \mathscr{H}om_{\mathcal{O}_M}\left( \Theta_{\mathcal{F}_0}, \Theta_{\mathcal{F}_0} \right) \right)$ defined by
\begin{align} \label{ssc10}
B_{ij}: \Gamma\left( U_{ij}, \Theta_{\mathcal{F}_0} \right) &\to \Gamma\left( U_{ij} , \Theta_{\mathcal{F}_0} \right), \\
 T_{0i}^\alpha &\mapsto B_{ij}^\alpha:= \sum_{\beta=1}^p B_{ij}^{\alpha \beta} T_{0i}^\beta \notag
\end{align}

we define the norm $||\sigma ||$ by
\begin{align*}
||\sigma|| = \left|\left| \Gamma \right|\right| + \left|\left|B \right|\right| + \left|\left| \Pi \right|\right|
\end{align*}
where
\begin{align*}
\left|\left| \Gamma \right|\right|=\max_{i,j} \sup_{\xi_i\in U_i\cap U_j} \left| \Gamma_{ij} (\xi_i)\right| , \,\,\,\,\,\,\,\,\,\,\,\,\,\,\left|\left| B\right|\right|=\max_{i,j}\max_{\alpha,\beta}\sup_{\xi_i \in U_i\cap U_j}\left|B_{ij}^{\alpha\beta} \right|\,\,\,\,\,\,\,\,\,\,\,\,\,\,\left|\left| \Pi \right| \right| = \max_{i }\max_{\alpha,\gamma}\sup_{\xi_i\in U_i^\delta } \left| \Pi_i^{\alpha\gamma} (\xi_i) \right|
\end{align*}

\begin{lemma}\label{tcc30}
For any triple $\sigma=\left( \Pi,  \left( \Gamma, B \right) \right)=\left( \left\{ \Pi_i \right\}, \left\{ \left( \Gamma_{ij}, B_{ij} \right) \right\}\right)$ which is a $1$-cocycle in the \v Cech resolution of $\mathcal{E}_{\Theta_{\mathcal{F}_0}}^\bullet$, we can find $\varphi_i(\xi_i), s^\lambda$ and $b_i^{\alpha \gamma}(\xi_i)$ satisfying
\begin{align}
\varphi_i- \varphi_j &= - \Gamma_{ij} + \sum_{\lambda=1}^m s^\lambda \rho_{ij\lambda} \label{tcc34}\\
\left[ \varphi_i, T_{0i}^\alpha \right] &= - \Pi_i^\alpha + \sum_{\lambda=1}^m s^\lambda \alpha_{i\lambda}^\alpha + \sum_{\gamma=1}^p b_{i}^{\alpha \gamma} T_{0i}^\gamma \label{tcc35}\\
\left| \varphi_i(\xi_i)\right|\leq K \left| \left| \sigma\right|\right|,&\,\,\,\,\,\,\,\,\,\,\,|s|\leq K ||\sigma|| ,\,\,\,\,\,\,\,\,\,\,\,\, \left|b_i^{\alpha \gamma}(\xi_i)\right| \leq K ||\sigma||
\end{align}
where $K$ is a constant independent of $\sigma$.
\end{lemma}

\begin{proof}
We define
\begin{align}
\iota(\sigma)=\inf \max_{i, \lambda, \alpha, \beta } \left\{  \sup_{\xi_i\in U_i} \left| \varphi_i \right| , \left| s^\lambda \right|,\sup_{\xi_i\in U_i}\left|b_i^{\alpha \beta} \right|  \right\}
\end{align}
where $\inf$ is taken with respect to all the solutions $s^\lambda, \varphi_i(\xi_i)$ and $b_{i}^{\alpha \beta}(\xi_i)$ of $(\ref{tcc34})$ and $(\ref{tcc35})$. It suffices to show that there exists a constant $K$ such that $\iota(\sigma)\leq K \left|\left| \sigma \right|\right|$ for all $1$-cocycle $\sigma$ of $\mathcal{E}_{\Theta_{\mathcal{F}_0}}^\bullet$. Suppose that there is no such constant $K$. Then we can find a sequence $\sigma^{(1)},\sigma^{(2)}, \cdots, \sigma^{(v)},\cdots$ of triple $\sigma^{(v)}=\left( \Pi^{(v)}, \left(\Gamma^{(v)}, B^{(v)} \right) \right)$ of $1$-cocycle of $\mathcal{E}_{\Theta_{\mathcal{F}_0}}^\bullet$ such that $\iota\left( \sigma^{(v)} \right)=1$ and $\left|\left| \sigma^{(v)} \right|\right|< \frac{1}{v}$ such that
\begin{align}
\varphi_i^{(v)}- \varphi_j^{(v)}&= - \Gamma_{ij}^{(v)} + \sum_{\lambda=1}^m s^{\lambda (v)}\rho_{ij\lambda}\label{ttc6} \\
\left[ \varphi_i^{(v)}, T_{0i}^\alpha \right] &= - \Pi_i^{\alpha (v)} + \sum_{\lambda=1}^m s^{\lambda (v)} \alpha_{i\lambda}^\alpha + \sum_{\gamma=1}^p b_{i}^{\alpha \gamma(v)} T_{0i}^\gamma  \label{ttc7}\\
\left| \varphi_i^{(v)} \right|<2,&\,\,\,\,\,\,\,\left| s^{\lambda (v)} \right| <2, \,\,\,\,\,\,\,\,\left| b_i^{\alpha \gamma(v)}\right| <2 \label{ttc5}
\end{align}
Hence by $(\ref{ttc5})$, replacing $\sigma^{(1)},\sigma^{(2)}, \cdots $ by a suitable subsequence, we may assume that $\varphi_i^{(v)}$ converges uniformly on each compact subset of $U_i$ and $s^{(v)}$ converges, and $b_i^{\alpha \gamma(v)}$ converges uniformly on each compact subset of $U_i$. As in \cite{Kod05} p.296, we can deduce that $\varphi_i^{(v)}$ converges uniformly on the whole $U_i$. Now we claim that $b_i^{\alpha \gamma(v)}$ converges uniformly on the whole $U_i$. In fact, first we note that $b_i^{\alpha\gamma(v)}$ converges uniformly on $U_i^\delta$ from $(\ref{ncc10})$. On $U_i\cap U_j^\delta$, $(\ref{ttc6})$ and $(\ref{ttc7})$ implies that
\begin{align}\label{ttc8}
\sum_{\beta, \sigma=1}^p B_{ij}^{\alpha\beta(v)} b_{0i}^{\beta\sigma} r_{0ij}'^{\sigma \gamma} = - \sum_{\beta=1}^p b_{i}^{\alpha \beta(v)} r_{0ij}'^{\beta\gamma} +\sum_{\beta=1}^p \left[  \varphi_{j}^{(v)} , r_{0ij}^{\alpha\beta} \right] b_{0j}^{\beta\gamma} + \sum_{\beta=1}^p \sum_{\lambda=1}^r \frac{\partial r_{ij}^{\alpha\beta}}{\partial t_\lambda}|_{t=0} s^{\lambda (v)}b_{0j}^{\beta\gamma}  + \sum_{\beta=1}^p r_{0ij}^{\alpha\beta} b_{j}^{\beta\gamma(v)} 
\end{align}
Since $s^{\lambda (v)}$ converges, and $\varphi_j$ converges uniformly on $U_j$, and $b_j^{\beta \gamma(v)}$ converges uniformly on $U_j^\delta$, and $\left| B_{ij}^{\alpha\beta(v)}\right| \leq \left|\left| B \right|\right|\to 0$, $(\ref{ttc8})$ implies that $b_i^{\alpha\beta(v)}$ converges uniformly on the whole $U_i$. Then we put $\varphi_i = \lim \varphi_i^{(v)}$ and $s^\lambda= \lim s^{\lambda(v)}$ and $b_i^{\alpha \gamma} =\lim b_i^{\alpha\gamma(v)}$. Then from $(\ref{ttc6})$ and $(\ref{ttc7})$ and $\left| \Gamma_{ij}^{(v)} \right|\to 0$ on $U_{ij}$, and $\left| \Pi_i^{\alpha\gamma(v)}\right|\to 0$ on $U_i^\delta$, we have 
\begin{align}\label{ttc11}
\varphi_i- \varphi_j &= \sum_{\lambda=1}^m s^\lambda \rho_{ij\lambda}\,\,\,\,\,\,\,\,\textnormal{on}\,\,\, U_{ij}\\
\left[ \varphi_i, T_{0i}^\alpha \right] &=  \sum_{\lambda=1}^m s^\lambda \alpha_{i\lambda}^\alpha + \sum_{\gamma=1}^p b_i^{\alpha \gamma} T_{0i}^\gamma \,\,\,\,\,\,\textnormal{on}\,\,\, U_i^\delta \notag
\end{align}
Then by identity theorem, we have
\begin{align}\label{ttc12}
\left[ \varphi_i, T_{0i}^\alpha \right] &=  \sum_{\lambda=1}^m s^\lambda \alpha_{i\lambda}^\alpha + \sum_{\gamma=1}^p b_i^{\alpha \gamma} T_{0i}^\gamma \,\,\,\,\,\,\textnormal{on}\,\,\, U_i
\end{align}

We put $\tilde{\varphi}_i^{(v)} := \varphi_i^{(v)} - \varphi_i$ and $\tilde{s}^{\lambda(v)}:= s^{\lambda(v)}- s^\lambda$ and $\tilde{b}_i^{\alpha \gamma (v)} := b_i^{\alpha \gamma (v)} - b_i^{\alpha \gamma}$. Then for a sufficiently large $v$, we have
\begin{align}
\left| \tilde{\varphi}_i^{(v)} \right| <\frac{1}{2},\,\,\,\,\,\,\,\,\,\left| \tilde{s}^{\lambda (v)} \right|< \frac{1}{2} , \,\,\,\,\,\,\,\,\,\, \left| \tilde{b}_i^{\alpha \gamma(v)} \right| < \frac{1}{2}
\end{align}
while we infer from $(\ref{ttc6}), (\ref{ttc7}),(\ref{ttc11})$ and $(\ref{ttc12})$ that
\begin{align*}
\tilde{\varphi}_i^{(v)}- \tilde{\varphi}_j^{(v)}&= - \Gamma_{ij}^{(v)} + \sum_{\lambda=1}^m \tilde{s}^{\lambda (v)}\rho_{ij\lambda} \\
\left[ \tilde{\varphi}_i^{(v)}, T_{0i}^\alpha \right] &= - \Pi_i^{\alpha (v)} + \sum_{\lambda=1}^m \tilde{s}^{\lambda (v)} \alpha_{i\lambda}^\alpha + \sum_{\gamma=1}^p \tilde{b}_{i}^{\alpha \gamma(v)} T_{0i}^\gamma
\end{align*}
This contradicts to $\iota\left( \sigma \right)=1$.

\end{proof}

Then by Lemma \ref{tcc31} and Lemma \ref{tcc30} we can choose solutions $\varphi_{i|\mu}(\xi_i,u)$ and $s_{\mu}^\lambda (u)$ and $b_{i|\mu}^{\alpha \beta}(\xi_i, u)$ of equations $(\ref{tc8})$ and $(\ref{tc10})$ such that
\begin{align*}
\varphi_{i|\mu}(\xi_i, u)\ll K K^* A(u),\,\,\,\,\,\,\,\,\,\,\,\,\,s_\mu(u)\ll K K^* A(u),\,\,\,\,\,\,\,\,\,\,\,\,\, b_{i|\mu}^{\alpha \beta} \ll K K^* A(u)
\end{align*}
where $K^*=\frac{K_1+ K_{27}}{b}+ \frac{K_2 + K_{28}+ K_{34}}{c} + \frac{( K_3 + K_{29} + K_{35})b}{c}$. We choose $b$ and $c$ in a way that $K K^* <1$. Then we have $\varphi_{i|\mu}(\xi_i, u) \ll A(u), s_\mu (u)\ll A(u)$ and $b_{i|\mu}^{\alpha\beta}(\xi_i, u) \ll A(u)$. Hence the induction $(\ref{tp1})_\mu$ holds for $\mu$. This completes the proof of Theorem \ref{tc1}.

\end{proof}

\section{Deformations of foliated complex analytic structures in terms of cotangent sheaves}

\begin{definition}[compare \cite{Kod05} p.59 and Definition \ref{t3}]\label{n3}
Suppose that given a domain $B\subset \mathbb{C}^m$, there is a set $\left\{\left(M_t, \mathcal{N}_{\mathcal{F}_t}^*\right)|t\in B \right\}$ of $n$-dimensional compact $($singularly$)$ foliated complex manifold $\left(M_t, \mathcal{N}_{\mathcal{F}_t}^*\right)$ of codimension $q$, so that for each $t\in B$, we have an exact sequence
\begin{align}\label{n1}
0\to \mathcal{N}_{\mathcal{F}_t}^*\to \Omega_{M_t}^1\to \Omega_{M_t}^1/\mathcal{N}_{\mathcal{F}_t}^*\to 0
\end{align}
where $\Omega_{M_t}/\mathcal{N}_{\mathcal{F}_t}^*$ is a torsion-free $\mathcal{O}_{M_t}$-module. We say that $\left\{\left(M_t, \mathcal{N}_{\mathcal{F}_t}^*\right)|t\in B\right\}$ is a family of $($singularly$)$ foliated complex manifold or $($singularly$)$ foliated complex analytic family in terms of cotangent sheaves if there exist a complex manifold $\mathcal{M}$ and a holomorphic map $\pi:\mathcal{M}\to B$ such that
\begin{enumerate}
\item we have an exact sequence
\begin{align}\label{n2}
0 \to \mathcal{N}_{\mathcal{F}}^*\to \Omega_{\mathcal{M}/B}^1\to \Omega_{\mathcal{M}/B}^1/\mathcal{N}_{\mathcal{F}}^*\to 0
\end{align}
where $\mathcal{N}_\mathcal{F}^*$ is a coherent subsheaf of the relative cotangent sheaf $\Omega_{\mathcal{M}/B}^1$ over $B$.
\item for each $t\in B$, $\pi^{-1}(t)= \left(M_t, \mathcal{N}_{\mathcal{F}_t }^*\right)$ and the exact sequence $(\ref{n2})$ is induced from $(\ref{n1})$ by restricting to $M_t$.
\item the rank of the Jacobian of $\pi$ is equal to $m$ at every point of $\mathcal{M}$.
\item $\mathcal{N}_{\mathcal{F}}^*$ satisfies the integrability condition$:$ let $\Theta_{\mathcal{F}}$ be the subsheaf of the relative tangent sheaf $\Theta_{\mathcal{M}/B}$ over $B$ which vanishes on $\mathcal{N}_{\mathcal{F}}^*$. Then $d_{\mathcal{M}/B}(\mathcal{N}_\mathcal{F}^*)$ vanishes on $\bigwedge^2 \Theta_\mathcal{F}$, where $d_{\mathcal{M}/B}$ is the relative differential on $\Omega_{\mathcal{M}/B}^1$. By abuse of notation, we will denote $d_{\mathcal{M}/B}$ by $d$. Or equivalently let $S_t=\textnormal{Sing}(\mathcal{F}_t)$ for each $t\in B$. Then we have $d\mathcal{N}_{\mathcal{F}}^*\subset \mathcal{N}_{\mathcal{F}}^*\bigwedge \Omega_{\mathcal{M}/B}^1$ at every point of $\mathcal{M}-\bigcup_{t\in B} S_t$.
\item $\Omega_{\mathcal{M}/B}^1/\mathcal{N}_{\mathcal{F}}^*$ is flat over $B$.
\end{enumerate}
We will denote the $($singularly$)$ foliated complex analytic family by $(\mathcal{M}, \mathcal{N}_\mathcal{F}^*, B, \pi)$.

\end{definition}

\begin{remark}
Let $\left(\mathcal{M}, \mathcal{N}_\mathcal{F}^*, B, \pi \right)$ with $\mathcal{N}_\mathcal{F^*}$ locally free be a foliated analytic family. Let $\Delta$ be an open set of $B$. Then the restriction $\left(\mathcal{M}_\Delta= \pi^{-1}(\Delta), \mathcal{N}_{\mathcal{F}}^*|_\Delta, \Delta, \pi |_{\mathcal{M}_\Delta} \right)$ is also a foliated complex analytic family in terms of cotangent sheaves. We will denote the family by $\left( \mathcal{M}_\Delta, \mathcal{N}_{\mathcal{F}_\Delta}^*, \Delta, \pi\right)$.
\end{remark}

\begin{remark}
When we ignore foliated structures, a foliated $($complex$)$ analytic family $\left( \mathcal{M}, \mathcal{N}_\mathcal{F}^*, B, \pi \right)$ is a complex analytic family $\left( \mathcal{M}, B, \pi \right)$ in the sense of Kodaira-Spencer $($see \textnormal{\cite{Kod05} p.59}$)$.
\end{remark}

\begin{remark}
We elaborate the condition $(4)$. Let $S$ be the  locus of points of $\mathcal{M}$ such that $\frac{\Omega_{\mathcal{M}/B}^1}{\mathcal{N}_\mathcal{F}^*}$ is not locally free. Then $S$ is closed. Since $\frac{\Omega_{\mathcal{M}/B}^1}{\mathcal{N}_\mathcal{F}^*}$ is flat, the locus of points of $\mathcal{M}$ such that $\frac{\Omega_{\mathcal{M}/B}^1}{\mathcal{N}_\mathcal{F}^*}$ is locally free is $\bigcup_{t\in B} M_t-S_t$ since $\frac{\Omega_{\mathcal{M}/B}^1}{\mathcal{N}_\mathcal{F}^*}|_{M_t}= \frac{\Omega_{M_t}^1}{\mathcal{N}_{\mathcal{F}_t}^*}$ is locally free on $M_t-S_t$. This implies that $S=\bigcup_{t\in B} S_t$, and so $\mathcal{M}- \bigcup_{t\in B} S_t$ is open. Consider the restriction of the exact sequence $(\ref{n2})$ to $\mathcal{M}- \bigcup_{t\in B} S_t$ which is an open set, we have the exact sequence
\begin{align*}
0 \to \mathcal{N}_{\mathcal{F}}^* |_{\left(\mathcal{M}- \bigcup_t S_t\right)} \to \Omega_{\mathcal{M}/B}^1 |_{\left(\mathcal{M}- \bigcup_t S_t\right)} \to \frac{\Omega_{\mathcal{M}/B}^1}{\mathcal{N}_{\mathcal{F}}^*}|_{\left(\mathcal{M}-\bigcup_t S_t\right)} \to 0
\end{align*}
Since $\frac{\Omega_{\mathcal{M}/B}^1}{\mathcal{N}_{\mathcal{F}}^*}|_{\left(\mathcal{M}-\bigcup_t S_t\right)}$ is locally free, $\Theta_{\mathcal{F}}|_{\left(\mathcal{M}- \bigcup_t S_t \right)}$ is locally free, and so flat over $B$. In particular, this implies that $\Theta_{\mathcal{F}}|_{M_t- S_t}= \Theta_{\mathcal{F}_t}|_{M_t-S_t}$.  Assume that $d\mathcal{N}_\mathcal{F}^*$ vanishes on $\bigwedge^2 \Theta_\mathcal{F}$ which is equivalent to $\left[ \Theta_\mathcal{F}, \Theta_\mathcal{F}\right]\subset \Theta_\mathcal{F}$. Then $d\mathcal{N}_{\mathcal{F}}^*\subset \mathcal{N}_{\mathcal{F}}^*\bigwedge \Omega_{\mathcal{M}/B}^1$ on $\mathcal{M}- \bigcup_{t\in B} S_t$. Conversely, assume that $d\mathcal{N}_{\mathcal{F}}^*\subset \mathcal{N}_{\mathcal{F}}^*\bigwedge \Omega_{\mathcal{M}/B}^1$ on $\mathcal{M}- \bigcup_{t\in B} S_t$. Then for $T_1,T_2\in \Theta_\mathcal{F}$ and $w\in \mathcal{N}_{\mathcal{F}}^*$, we have $i_{T_1\wedge T_2}(w)=0$ on $\mathcal{M}- \bigcup_{t\in B} S_t$. Then $i_{T_1\wedge T_2}(w)=0$ on $M_t-S_t$, and so it is $0$ on $M_t$. Hence $d\mathcal{N}_\mathcal{F}^*$ vanishes on $\bigwedge^2 \Theta_\mathcal{F}$.
\end{remark}

\begin{remark}
We note that the flatness of $\Omega_{\mathcal{M}/B}^1/\mathcal{N}_\mathcal{F}^*$ implies the flatness of $\mathcal{N}_{\mathcal{F}}^*$ over $B$. In particular, if $\mathcal{N}_{\mathcal{F}_{t_0}}^*$ is a locally free subsheaf of $\Omega_{M_{t_0}}^1$ for some $t_0\in B$, then $\mathcal{N}_{\mathcal{F}_t}^*$ is also locally free in some neighborhood of $t_0$. In the following we shall assume that $\mathcal{N}_{\mathcal{F}}^*$ is a locally free $\mathcal{O}_\mathcal{M}$-submodule of $\Omega_{\mathcal{M}/B}^1$ with rank $q$, so that $\mathcal{N}_{\mathcal{F}_t}^*$ is a locally free subsheaf of $\Omega_{M_t}^1$ for each $t\in B$.
\end{remark}

\begin{remark}
If $\mathcal{N}_{\mathcal{F}_{t_0}}^*$ is a regular foliation on $M_{t_0}$ for some $t_0\in B$, then $\mathcal{N}_{\mathcal{F}_t}^*$ defines a regular foliation in some neighborhood $B_0$ of $t_0$ and the foliated complex analytic family $(\mathcal{M}|_{B_0},\mathcal{N}_{\mathcal{F}}^*|_{B_0},B_0,\pi)$ in terms of cotangent sheaves is also the foliated complex analytic family $\left(\mathcal{M}|_{B_0}, \Theta_{\mathcal{F}}|_{B_0}=\left(\Omega_{\mathcal{M}/B}^1/\mathcal{N}_{\mathcal{F}}^*\right)^*|_{B_0}, B_0, \pi \right)$ in terms of tangent sheaves as in \textnormal{Definition} $\ref{t3}$.
\end{remark}

\begin{remark}\label{n5}
From \textnormal{Definition \ref{n3}}, we can choose a system of local complex coordinates $\{z_i,...,z_j,...\},z_j:p\to z_j(p)$, and coordinate polydisks $\mathcal{U}_j$ with respect to $z_j$, satisfying the following conditions\footnote{for the detail, see \cite{Kod05} p.60}:
\begin{enumerate}
\item $z_j(p)=(z_j^1(p),...,z_j^n(p),t_1,...,t_m),(t_1,...,t_m)=\pi(p)$;
\item $\mathcal{U}=\{\mathcal{U}_j|j=1,2,...\}$ is locally finite.
\end{enumerate}
Then $\left\{p\to \left(z_j^1(p),...,z_j^n(p)\right)|\mathcal{U}_j\cap M_t\ne \emptyset \right\}$ gives a system of local complex coordinates on $M_t$. In terms of these coordinates, $\pi$ is the projection given by $\pi: \left(z_j^1,...,z_j^n,t_1,...,t_m\right)\to \left(t_1,...,t_m\right)$. For $j,k$ with $\mathcal{U}_j\cap \mathcal{U}_k\ne \emptyset$, we denote the coordinate transformation from $z_k$ to $z_j$ by
\begin{align*}
f_{jk}:\left(z_k^1,...,z_k^n,t\right)\to \left(z_j^1,...,z_j^n,t\right)=f_{jk}\left(z_k^1,...,z_k^n,t\right)
\end{align*}
Thus $f_{jk}$ is given by 
\begin{align*}
z_j^\alpha=f_{jk}^\alpha\left(z_k^1,...,z_k^n,t_1,...,t_m\right),\alpha=1,...,n.
\end{align*}

On the other hand, we may assume that the locally free sheaf $\mathcal{N}_\mathcal{F}^*$ is trivialized in terms of the covering $\mathcal{U}$ so that
\begin{align*}
\Gamma\left(\mathcal{U}_j, \mathcal{N}_\mathcal{F}^* \right)&\cong \bigoplus^q \Gamma\left(\mathcal{U}_j, \mathcal{O}_\mathcal{M} \right)\\
w_j^\alpha\left(z_j,t\right)&\mapsto e_j^\alpha=\left(0,...,\overbrace{1}^{\alpha-\textnormal{th}},...,0\right),\alpha=1,...,q
\end{align*}
which is generated by holomorphic $1$-forms $w_j^\alpha(z_j,t)\in \Gamma\left(\mathcal{U}_j,\mathcal{N}_\mathcal{F}^* \right),\alpha=1,...,q$ of the form
\begin{align}\label{n11}
w_j^\alpha(z_j,t):=\sum_{\beta=1}^n w_j^{\alpha\beta}(z_j,t)dz_j^\beta, \,\,\,\,\,\,\alpha=1,...,q
\end{align}
for some $w_j^{\alpha\beta}(z_j,t)\in \Gamma\left(\mathcal{U}_j, \mathcal{O}_\mathcal{M} \right)$. For $(x,a)\in \mathcal{U}_j- \left(S:=\bigcup_{t\in B}\textnormal{Sing}(\mathcal{F}_t) \right)$, there exists a neighborhood $U_{(x,a)}$ of $(x,a)$ such that $U_{(x,a)}\subset  \mathcal{U}_j-S$, and  for $\alpha=1,...,q$,
\begin{align*}
dw_j^\alpha(z_j,t)=\sum_{\beta=1}^q a_{j_{(x,a)}}^{\alpha\beta}(z_j,t) \wedge w_j^\beta(z_j,t)\,\,\,\,\,\,\,\,\,\textnormal{for some}\,\,\,\,\, a_{j_{(x,a)}}^{\alpha\beta}(z_j,t) \in \Gamma\left(U_{(x,a)},\Omega_{\mathcal{M}/B}^1\right)
\end{align*}

Then $w_j^1(z_j,t)\wedge \cdots \wedge w_i^q(z_j,t)\wedge dw_i^\alpha (z_j,t)=0$ on $\mathcal{U}_j \cap (M_t-S_t)$, so that it vanishes on $\mathcal{U}_j\cap M_t$ for each $t\in B$ and so on $\mathcal{U}_j$. Hence we have
\begin{align}\label{n6}
w_j^1(z_j,t)\wedge \cdots \wedge w_j^q(z_j,t)\wedge dw_j^\alpha(z_j,t)=0,\,\,\,\,\,\,\,\alpha=1,...,q,\,\,\,\,\,\textnormal{on}\,\,\,\,\,\mathcal{U}_j.
\end{align}

For $j,k$ with $\mathcal{U}_j\cap \mathcal{U}_k\ne \emptyset$, we have invertible $q\times q$ matrices $H_{jk}:= \left(h_{jk}^{\alpha\beta}(z_j,t)\right)$ with components $h_{jk}^{\alpha\beta}(z_j,t)\in \Gamma\left(\mathcal{U}_j\cap \mathcal{U}_k,\mathcal{O}_\mathcal{M}\right)$ such that
\begin{center}
$\left[\begin{matrix}
w_j^1(z_j,t) \\
w_j^2(z_j,t)\\
\cdot\\
\cdot\\
\cdot\\
w_j^q(z_j,t)
\end{matrix}\right]
=
\left[\begin{matrix}
h_{jk}^{11}(z_j,t) &\cdots & h_{jk}^{1q}(z_j,t)\\
h_{jk}^{21}(z_j,t) & \cdots & h_{jk}^{2q}(z_j,t)\\
\cdot & \cdots &\cdot \\
\cdot & \cdots & \cdot \\
\cdot & \cdots & \cdot \\
h_{jk}^{q1}(z_j,t) & \cdots & h_{jk}^{qq}(z_j,t)
\end{matrix}\right]
\cdot 
\left[\begin{matrix}
w_k^1(z_k,t)\\
w_k^2(z_k,t)\\
\cdot\\
\cdot\\
\cdot\\
w_k^q(z_k,t)
\end{matrix}\right]$
\end{center}
and for $i,j,k\in \mathcal{U}_i\cap \mathcal{U}_j\cap \mathcal{U}_k\ne \emptyset$, we have
\begin{align*}
H_{ik}=H_{ij}\cdot H_{jk}\,\,\,\,\,\textnormal{on}\,\,\,\,\,\mathcal{U}_i\cap \mathcal{U}_j\cap \mathcal{U}_k.
\end{align*}

More specifically, $w_j^\alpha(z_j,t)= \sum_{\beta=1}^p h_{jk}^{\alpha\beta}w_k^\beta(z_k,t)$ implies that
\begin{align}
\sum_{\gamma=1}^n w_j^{\alpha\gamma}(z_j,t)dz_j^\gamma=\sum_{\beta=1}^q\sum_{\eta=1}^n h_{jk}^{\alpha\beta}(z_k,t)w_k^{\beta \eta}(z_k,t)dz_k^\eta
\end{align}
so that
\begin{align}\label{n10}
\sum_{\gamma,\eta=1}^n w_j^{\alpha\gamma}(f_{jk}(z_k,t),t)\frac{\partial f_{jk}^\gamma(z_k,t)}{\partial z_k^\eta} dz_k^\eta=\sum_{\beta=1}^q\sum_{\eta=1}^n h_{jk}^{\alpha\beta}(z_k,t)w_k^{\beta\eta}(z_k,t)dz_k^\eta
\end{align}

\end{remark}

\subsection{Dual leaf complex controlling deformations of singular holomorphic foliations in terms of locally free subsheaves of cotangent sheaf}\

Let $\left( \mathcal{M}, \mathcal{N}_\mathcal{F}^*, B, \pi \right)$ be a foliated complex analytic family with $\mathcal{N}_\mathcal{F}^*$ locally free as in Definition \ref{n3}, so that $\mathcal{N}_{\mathcal{F}_t}^*$ defines a (singular) holomorphic foliation on a compact complex manifold $M_t$ with $\mathcal{N}_{\mathcal{F}_t}^*$ locally free. Then we have the dual leaf complex on $M_t$ associated to $\mathcal{N}_{\mathcal{F}_t}^*$ (see Appendix \ref{AppendixA5})
\begin{align*}
\mathcal{N}_{\mathcal{F}_t}^{* \bullet } :\Theta_{M_t}\xrightarrow{E_0^t} \mathscr{H}om_{\mathcal{O}_{M_t}}\left(\mathcal{N}_{\mathcal{F}_t}^*, \frac{\Omega_{M_t}^1}{\mathcal{N}_{\mathcal{F}_t}^*} \right)\xrightarrow{E_1^t} \mathscr{H}om_{\mathcal{O}_{M_t}}\left( \mathcal{N}_{\mathcal{F}_t}^*, \tilde{\mathcal{S}}_t^2\right) \xrightarrow{E_2^t} \mathscr{H}om_{\mathcal{O}_{M_t}}\left( \mathcal{N}_{\mathcal{F}_t}^*, \tilde{\mathcal{S}}_t^3 \right)\xrightarrow{E_3^t} \cdots
\end{align*}
We will denote the $i$-th cohomology group by $\mathbb{H}^i\left( M_t , \mathcal{N}_{\mathcal{F}_t}^{* \bullet} \right)$. We can compute $\mathbb{H}^i\left(M_t, \mathcal{N}_{\mathcal{F}_t}^{*\bullet} \right)$  by the following \v Cech resolution of $\mathcal{N}_{\mathcal{F}_t}^{*\bullet}$. Here $\delta$ is the \v Cech map and $\mathcal{U}_t= \mathcal{U}\cap M_t=\left\{ U_j^t :=\mathcal{U}_j \cap M_t | j=1,2,... \right\}$ is an open covering of $M_t$:
{\small{\begin{center}
$\begin{CD}
\cdots \\
@AE_2^t AA \\
C^0\left(\mathcal{U}_t, \mathscr{H}om_{\mathcal{O}_{M_t}}\left( \mathcal{N}_{\mathcal{F}_t}^* , \tilde{\mathcal{S}}_t^2 \right) \right) @>-\delta>> \cdots  \\
@AE_1^tAA @AE_1^t AA \\
C^0\left(\mathcal{U}_t, \mathscr{H}om_{\mathcal{O}_{M_t}}\left(\mathcal{N}_{\mathcal{F}_t}^*, \frac{\Omega_{M_t}^1}{\mathcal{N}_{\mathcal{F}_t}^* } \right) \right)@>\delta>> C^1\left(\mathcal{U}_t, \mathscr{H}om_{\mathcal{O}_{M_t}}\left(\mathcal{N}_{\mathcal{F}_t}^*, \frac{\Omega_{M_t}^1}{\mathcal{N}_{\mathcal{F}_t}}  \right) \right) @>-\delta>> \cdots \\
@AE_0^t AA @AE_0^t AA @AE_0^t AA \\
C^0\left(\mathcal{U}_t, \Theta_{M_t}\right) @>-\delta>> C^1\left(\mathcal{U}_t, \Theta_{M_t}\right) @>\delta>> C^2\left(\mathcal{U}_t, \Theta_{M_t}\right) @>>> \cdots 
\end{CD}$
\end{center}}}

We will  relate the first cohomology group $\mathbb{H}^1\left( M_t, \mathcal{N}_{\mathcal{F}_t}^{* \bullet} \right)$ to infinitesimal foliated deformations of $\pi^{-1}(t)=\left(M_t, \mathcal{N}_{\mathcal{F}_t}^{* \bullet} \right)$ in the foliated analytic family $\left(\mathcal{M}, \mathcal{N}_\mathcal{F}^* , \pi, B \right)$ in terms of cotangent sheaves.

\subsection{Infinitesimal foliated deformations in terms of cotangent sheaves}\label{n13}\

Let $\left( \mathcal{M}, \mathcal{N}_\mathcal{F}^*, \pi, B \right)$ be a foliated analytic family with $\mathcal{N}_\mathcal{F}^*$ locally free as in Definition \ref{n3}. We keep the notations in Remark \ref{n5}. By taking the derivative of $(\ref{n11})$ with respect to $t$, we set
\begin{align*}
\frac{\partial w_j^\alpha(z_j,t)}{\partial t}:=\sum_{\beta=1}^n \frac{\partial w_j^{\alpha\beta}(z_j,t)}{\partial t}dz_j^\beta,\,\,\,\,\,\,\,\alpha=1,...,q.
\end{align*}

We define an element $\beta_j(t)$ in $\Gamma\left(U_j^t, \mathscr{H}om_{\mathcal{O}_{M_t}}\left(\mathcal{N}_{\mathcal{F}_t}^*, \frac{\Omega_{M_t}^1}{\mathcal{N}_{\mathcal{F}_t}^*} \right) \right)$ in the following way:
\begin{align*}
\beta_j(t):\Gamma\left(U_j^t, \mathcal{N}_{\mathcal{F}_t}^*\right)&\to  \Gamma\left(U_j^t, \frac{\Omega_{M_t}^1}{\mathcal{N}_{\mathcal{F}_t}^*}\right)\\
                                                                         w_j^\alpha(z_j,t) &\mapsto  -\overline{\frac{\partial w_j^{\alpha}(z_j,t)}{\partial t} }
\end{align*}
where $\overline{\frac{\partial w_j^{\alpha}(z_j,t)}{\partial t} }$ denotes the image of $\frac{\partial w_j^{\alpha}(z_j,t)}{\partial t}$ in the natural map $\Omega_{M_t}^1\to \Omega_{M_t}/\mathcal{N}_{\mathcal{F}_t}^*$, and linearly extends to $\Gamma\left( U_j^t, \mathcal{N}_{\mathcal{F}_t}^* \right)$. We also define an element $\tilde{\beta}_j(t)$ in $\Gamma\left(U_j^t, \mathscr{H}om_{\mathcal{O}_{M_t}}\left(\mathcal{N}_{\mathcal{F}_t}^*, \Omega_{M_t}^1 \right) \right)$ in the following way:
\begin{align*}
\tilde{\beta}_j(t):\Gamma\left(U_j^t, \mathcal{N}_{\mathcal{F}_t}^*\right)&\to  \Gamma\left(U_j^t, \Omega_{M_t}^1\right)\\
                                                                         w_j^\alpha(z_j,t) &\mapsto  -\frac{\partial w_j^{\alpha}(z_j,t)}{\partial t} 
\end{align*}
and linearly extends to $\Gamma\left( U_j^t, \mathcal{N}_{\mathcal{F}_t}^* \right)$. Then we have
\begin{proposition}\label{n14}
{\Small{\begin{align*}
\left(\left\{\theta_{jk}(t):=\sum_{\alpha=1}^n \frac{\partial f_{jk}^\alpha(z_k,t)}{\partial t}\frac{\partial}{\partial z_j^\alpha} \right\}, \left\{ \beta_j(t):=\left(w_j^\alpha(z_j,t)\mapsto -\overline{\frac{\partial w_j^\alpha(z_j,t)}{\partial t}}\right)_{\alpha=1,...,q} \right\} \right) \in C^1\left(\mathcal{U}^t, \Theta_{M_t} \right) \bigoplus C^0\left( \mathcal{U}_t, \mathscr{H}om_{\mathcal{O}_{M_t}}\left(\mathcal{N}_{\mathcal{F}_t}^*, \frac{\Omega_{M_t}^1}{\mathcal{N}_{\mathcal{F}_t}^*} \right) \right)
\end{align*}}}
defines a $1$-cocycle in the above \v Cech resolution of $\mathcal{N}_{\mathcal{F}_t}^{*\bullet}$ and call its cohomology class in $\mathbb{H}^1\left( M_t, \mathcal{N}_{\mathcal{F}_t}^{*\bullet} \right)$ the infinitesimal $($foliated$)$ deformation along $\frac{\partial}{\partial t}$. This expression is independent of the choice of local coordinates.
\end{proposition}

\begin{proof}
First we note that $\delta\left( \left\{ \theta_{jk}(t)\right\} \right)=0$ (see \cite{Kod05} p.201). Second, by taking the derivative of $(\ref{n6})$ with respect to $t$, we have
\begin{align*}
&w_j^1\wedge \cdots \wedge w_j^q \wedge d\left(\frac{\partial w_j^\alpha}{\partial t} \right)+\sum_{\gamma=1}^q w_j^1\wedge \cdots \wedge \overbrace{\frac{\partial w_j^\gamma}{\partial t}}^{\gamma-\textnormal{th}} \wedge \cdots \wedge w_j^q\wedge dw_j^\alpha=0\\
&\Longrightarrow w_j^1\wedge \cdots \wedge w_j^q \wedge d \tilde{\beta}_j^\alpha(t)+\sum_{\gamma=1}^q w_j^1\wedge \cdots \wedge \overbrace{ \tilde{\beta}_j^\gamma(t)}^{\gamma-\textnormal{th}} \wedge \cdots \wedge w_j^q\wedge dw_j^\alpha=0
\end{align*}
This implies that $E_1^t\left( \beta_j(t) \right)=0$. It remains to show that $\delta\left( \left\{ \beta_j(t)\right\} \right) + E_0^t\left(\left\{ \theta_{jk}\right\} \right)=0$, equivalently $\left( \beta_k(t) - \beta_j(t)  \right)(w_j^\alpha) + \overline{\mathcal{L}_{\theta_{jk}(t)}\left(w_j^\alpha \right)}=0$ for $\alpha=1,...,q$. In fact, by taking the derivative of $(\ref{n10})$ with respect to $t$ and considering the coefficient of $dz_k^\eta$, we have
\begin{align*}
\sum_{\gamma,\xi=1}^n \frac{\partial w_j^{\alpha\gamma}}{\partial z_j^\xi}\frac{\partial f_{jk}^\xi}{\partial t}\frac{\partial f_{jk}^\gamma}{\partial z_k^\eta}+\sum_{\gamma=1}^n \frac{\partial w_j^{\alpha\gamma}}{\partial t}\frac{\partial f_{jk}^\gamma}{\partial z_k^\eta} +\sum_{\gamma=1}^n w_j^{\alpha\gamma}\frac{\partial}{\partial z_k^\eta}\left(\frac{\partial f_{jk}^\gamma}{\partial t} \right)=\sum_{\beta=1}^q \frac{\partial h_{jk}^{\alpha\beta}}{\partial t}w_k^{\beta\eta} +\sum_{\beta=1}^q h_{jk}^{\alpha\beta}\frac{\partial w_k^{\beta\eta}}{\partial t},
\end{align*}
so that we have
\begin{align} \label{n7}
\sum_{\gamma,\xi=1}^n \frac{\partial w_j^{\alpha\gamma}}{\partial z_j^\xi}\frac{\partial f_{jk}^\xi}{\partial t}dz_j^\gamma+\sum_{\gamma=1}^n \frac{\partial w_j^{\alpha\gamma}}{\partial t}dz_j^\gamma +\sum_{\gamma,\eta=1}^n w_j^{\alpha\gamma}\frac{\partial}{\partial z_k^\eta}\left(\frac{\partial f_{jk}^\gamma}{\partial t} \right)dz_k^\eta=\sum_{\eta=1}^n\sum_{\beta=1}^q \frac{\partial h_{jk}^{\alpha\beta}}{\partial t}w_k^{\beta\eta}dz_k^\eta +\sum_{\eta=1}^n\sum_{\beta=1}^q h_{jk}^{\alpha\beta}\frac{\partial w_k^{\beta\eta}}{\partial t}dz_k^\eta
\end{align}
Since $w_j^\alpha(z_j,t)= \sum_{\beta=1}^q h_{jk}^{\alpha\beta}(z_k,t) T_k^\beta(z_k,t)$, we have
\begin{align}\label{n8}
\left( \tilde{\beta}_k(t)- \tilde{\beta}_j(t) \right)\left(w_j^\alpha\right)= -\sum_{\beta=1}^q h_{jk}^{\alpha\beta} \frac{\partial w_k^\beta}{\partial t}+\frac{\partial w_j^\alpha}{\partial t}=-\sum_{\beta,\eta=1}^n h_{jk}^{\alpha\beta}\frac{\partial w_k^{\beta\eta}}{\partial t} dz_k^\eta+\sum_{\gamma=1}^n \frac{\partial w_j^{\alpha\gamma}}{\partial t}dz_j^\gamma
\end{align}
Let us compute
\begin{align}\label{n9}
\mathcal{L}_{\theta_{jk}(t)}\left( w_j^\alpha \right)&=\mathcal{L}_{\sum_{\xi=1}^n \frac{\partial f_{jk}^\xi}{\partial t}\frac{\partial}{\partial z_j^\xi}}\left(\sum_{\gamma=1}^n w_j^{\alpha\gamma}dz_j^\gamma \right)=\sum_{\gamma,\xi=1}^n \frac{\partial f_{jk}^\xi}{\partial t}\frac{\partial w_j^{\alpha\gamma}}{\partial z_j^\xi} dz_j^\gamma +\sum_{\gamma=1}^n w_j^{\alpha\gamma}d\left( \frac{\partial f_{jk}^\gamma}{\partial t}\right)\\
 &=\sum_{\gamma,\xi=1}^n \frac{\partial f_{jk}^\xi}{\partial t}\frac{\partial w_j^{\alpha\gamma}}{\partial z_j^\xi} dz_j^\gamma +\sum_{\gamma,\eta=1}^n w_j^{\alpha\gamma}\frac{\partial}{\partial z_k^\eta}\left(\frac{\partial f_{jk}^\gamma}{\partial t} \right) dz_k^\eta \notag
\end{align}
Then $(\ref{n7}), (\ref{n8})$ and $(\ref{n9})$ implies that
\begin{align*}
\left( \tilde{\beta}_k(t)- \tilde{\beta}_j(t) \right)\left(w_j^\alpha\right) + \mathcal{L}_{\theta_{jk}(t)}\left( w_j^\alpha \right)= \sum_{\eta=1}^n \sum_{\beta=1}^q \frac{\partial h_{jk}^{\alpha\beta}}{\partial t} w_k^{\beta \eta} dz_k^\eta \in \mathcal{N}_{\mathcal{F}_t}^*
\end{align*}
Hence we have $\beta_k(t)- \beta_j(t)+ E_0^t \left(\theta_{jk}(t) \right)=0$. Next we show that $\left(\left\{ \theta_{jk}(t)\right\}, \left\{ \beta_j(t)\right\} \right)$ is independent of the choice of system of local coordinates. We can show that the infinitesimal deformation does not change under the refinement of the open covering (see \cite{Kod05} p.190). Since we can choose a common refinement for two systems of local coordinates, it is sufficient to show that given two local coordinates $x_j=(z_j,t)$ and $u_j=(y_j,t)$ on each $\mathcal{U}_j$, the infinitesimal foliated deformation $\left(\left\{\eta_{jk}(t) \right\}, \left\{ \chi_j(t) \right\} \right)$ with respect to $\{u_j\}$ coincides with $\left( \left\{\theta_{jk} (t) \right\} , \left\{ \beta_j (t) \right\} \right)$ with respect to $\{x_j\}$. Let $\mathcal{N}_{\mathcal{F}}^*$ be generated by $W_j^\alpha(y_j,t)=\sum_{\gamma=1}^n W_j^{\alpha\gamma}(y_j,t)dy_j^\gamma,\alpha=1,...,q$ in terms of coordinates $u_j=(y_j,t)$ on $\mathcal{U}_j$. Let $\left(y_k,t \right)\to \left( y_j, t \right)= \left(e_{jk}(y_k,t),t \right)$ be the coordinate transformation of $\{u_j\}$ on $\mathcal{U}_j\cap \mathcal{U}_k\ne \emptyset$. Then we have $\eta_{jk}(t)=\sum_{\alpha=1}^n \frac{\partial e_{jk}^\alpha(y_k,t)}{\partial t}\frac{\partial}{\partial y_j^\alpha}, y_k= e_{kj}(y_j,t)$, and 
\begin{align*}
\chi_j (t):\Gamma\left(U_j^t,  \mathcal{N}_{\mathcal{F}_t}^* \right) \to \Gamma\left( U_j^t, \frac{\Omega_{M_t}^1}{\mathcal{N}_{\mathcal{F}_t}^*} \right), \,\,\,\,\,& W_j^\alpha \mapsto - \overline{\frac{\partial W_j^\alpha}{\partial t}}:= - \overline{\sum_{\gamma=1}^n \frac{\partial W_j^{\alpha \gamma}(y_j,t)}{\partial t} \frac{\partial}{\partial y_j^\gamma}}\\
\tilde{\chi}_j(t):\Gamma\left(U_j^t,  \mathcal{N}_{\mathcal{F}_t}^* \right) \to \Gamma\left( U_j^t, \Omega_{M_t}^1 \right),\,\,\,\,\,& W_j^\alpha \mapsto -\frac{\partial W_j^\alpha}{\partial t}:= -\sum_{\gamma=1}^n \frac{\partial W_j^{\alpha \gamma}(y_j,t)}{\partial t} \frac{\partial}{\partial y_j^\gamma}
\end{align*}
We show that $\left(\left\{ \theta_{jk}(t) \right\}, \left\{ \beta_j(t)\right\} \right)$ is cohomologous to $\left(\left\{ \eta_{jk}(t)\right\}, \left\{ \chi_j(t)\right\} \right)$. Let $y_j^\alpha=g_j^\alpha\left(z_j^1,..., z_j^n, t \right),\alpha=1,...,n$, define the coordinate transformation from $x_j=(z_j,t)$ to $u_j=(y_j,t)$. Since $\left\{w_j^\alpha\right\}$ and $\left\{W_j^\alpha \right\},\alpha=1,...,n$ are two bases of $\Gamma\left(\mathcal{U}_j, \mathcal{N}_\mathcal{F}^* \right)$, we can write $W_j^\alpha(y_j,t)= \sum_{\beta=1}^q b_j^{\alpha\beta}(z_j,t) w_j^\beta(z_j,t)$ for some $b_j^{\alpha \gamma}(z_j,t)\in \Gamma\left(\mathcal{U}_j, \mathcal{O}_{\mathcal{M}} \right)$, i.e. $\sum_{\gamma=1}^n W_j^{\alpha\gamma}(y_j,t)dy_j^\gamma = \sum_{\beta=1}^q \sum_{\eta=1}^n b_j^{\alpha\beta}(z_j,t) w_j^{\beta \eta}(z_j,t)dz_j^\eta$, so that by setting $g_j(z_j,t)=\left(g_j^1(z_j,t),..., g_j^n(z_j,t) \right)$, we have
\begin{align}\label{n12}
\sum_{\gamma,\eta=1}^n W_j^{\alpha\gamma}\left( g_j(z_j,t),t \right)\frac{\partial g_j^\gamma (z_j,t)}{\partial z_j^\eta} d z_j^\eta= \sum_{\beta=1}^q \sum_{\eta=1}^n b_j^{\alpha\beta}(z_j,t) w_j^{\beta \eta}(z_j,t) d z_j^\eta
\end{align}

We set $\theta_j(t)=\sum_{\alpha=1}^n \frac{\partial g_j^\alpha(z_j,t)}{\partial t} \frac{\partial}{\partial y_j^\alpha}, y_j^\alpha= g_j^\alpha(z_j,t)$. Then $-\delta\left( \left\{\theta_j(t) \right\} \right)= \left\{ \eta_{jk}(t)- \theta_{jk}(t)\right\}$ (for the detail, see \cite{Kod05} p.191-192), and we claim that $\overline{\mathcal{L}_{\theta_j(t)} \left(W_j^\alpha\right) }=\chi_j(t)(W_j^\alpha)- \beta_j(t)\left( W_j^\alpha \right)$. In fact, by taking the derivative of $(\ref{n12})$ with respect to $t$, and considering the coefficient of $dz_j^\eta$, we have
\begin{align*}
\sum_{\gamma,\xi=1}^n \frac{\partial W_j^{\alpha\gamma}}{\partial y_j^\xi}\frac{\partial g_j^\xi}{\partial t}\frac{\partial g_{j}^\gamma}{\partial z_j^\eta}+\sum_{\gamma=1}^n \frac{\partial W_j^{\alpha\gamma}}{\partial t}\frac{\partial g_{j}^\gamma}{\partial z_j^\eta} +\sum_{\gamma=1}^n W_j^{\alpha\gamma}\frac{\partial}{\partial z_j^\eta}\left(\frac{\partial g_{j}^\gamma}{\partial t} \right)=\sum_{\beta=1}^q \frac{\partial b_{j}^{\alpha\beta}}{\partial t}w_j^{\beta\eta} +\sum_{\beta=1}^q b_{j}^{\alpha\beta}\frac{\partial w_j^{\beta\eta}}{\partial t},
\end{align*}
so that we have
\begin{align*}
\sum_{\gamma,\xi=1}^n \frac{\partial W_j^{\alpha\gamma}}{\partial y_j^\xi}\frac{\partial g_{j}^\xi}{\partial t}dy_j^\gamma+\sum_{\gamma=1}^n \frac{\partial W_j^{\alpha\gamma}}{\partial t}dy_j^\gamma +\sum_{\gamma,\eta=1}^n W_j^{\alpha\gamma}\frac{\partial}{\partial z_j^\eta}\left(\frac{\partial g_{j}^\gamma}{\partial t} \right)dz_j^\eta=\sum_{\eta=1}^n\sum_{\beta=1}^q \frac{\partial b_{j}^{\alpha\beta}}{\partial t}w_j^{\beta\eta}dz_j^\eta +\sum_{\eta=1}^n\sum_{\beta=1}^q b_{j}^{\alpha\beta}\frac{\partial w_j^{\beta\eta}}{\partial t}dz_j^\eta
\end{align*}
Then since $W_j^\alpha(y_j,t)= \sum_{\beta=1}^q b_j^{\alpha\beta}(z_j,t) w_j^\beta(z_j,t)$, we have
{\small{\begin{align*}
&\tilde{\beta}_j(t)\left( W_j^\alpha \right) - \tilde{\chi}_j(t)\left( W_j^\alpha \right) + \mathcal{L}_{ \theta_j(t)} \left( W_j^\alpha \right) = - \sum_{\beta=1}^q b_j^{\alpha\beta}\frac{\partial w_j^\beta}{\partial t} + \frac{\partial W_j^\alpha}{\partial t} +  \mathcal{L}_{ \sum_{\eta=1}^n \frac{\partial g_j^\eta}{\partial t } \frac{\partial}{\partial y_j^\eta}}\left( \sum_{\gamma=1}^n W_j^{\alpha \gamma}dy_j^\gamma \right) \\
&= - \sum_{\beta=1}^q\sum_{\eta=1}^n b_j^{\alpha\beta}\frac{\partial w_j^{\beta \eta}}{\partial t} d z_j^\eta  + \sum_{\gamma=1}^n \frac{\partial W_j^{\alpha \gamma}}{\partial t} d y_j^\gamma + \sum_{\eta,\gamma=1}^n \frac{\partial g_j^\eta}{\partial t} \frac{\partial W_j^{\alpha \gamma}}{\partial y_j^\eta} d y_j^\gamma + \sum_{\eta, \gamma=1}^n  W_j^{\alpha \gamma} d\left( \frac{\partial g_j^\gamma}{\partial t} \right) \\
&= \sum_{\beta=1}^q \sum_{\eta=1}^n \frac{\partial b_j^{\alpha \beta}}{\partial t} w_j^{\beta \eta} dz_j^\eta  = \sum_{\beta=1}^q \frac{\partial b_j^{\alpha \beta}}{\partial t} w_j^\beta \in \mathcal{N}_{\mathcal{F}_t}^*
\end{align*}}}
This implies that $\beta_j(t)\left( W_j^\alpha \right) - \chi_j(t)\left( W_j^\alpha \right) + \overline{\mathcal{L}_{\theta_j(t)} \left( W_j^\alpha \right)} =0$. Hence $\left(\left\{ \theta_{jk}(t) \right\}, \left\{ \beta_j(t)\right\} \right)$ is cohomologous to $\left(\left\{ \eta_{jk}(t)\right\}, \left\{ \chi_j (t)\right\} \right)$.

\end{proof}

\begin{definition}[foliated Kodaira-Spencer map in terms of cotangent sheaf]\label{n15}
Let $\left(\mathcal{M}, \mathcal{N}_\mathcal{F}^*, B, \pi \right)$ with $\mathcal{N}_{\mathcal{F}}^*$ locally free be a foliated complex analytic family of deformations of $\left(M_t, \mathcal{N}_{\mathcal{F}_t}^*\right)=\pi^{-1}(t),t\in B$, where $B$ is a domain of $\mathbb{C}^m$. We keep the notations in the proof of \textnormal{Proposition \ref{n14}}. For a tangent vector $\frac{\partial}{\partial t}=\sum_{\lambda=1}^m c_\lambda \frac{\partial}{\partial t_\lambda}, c_\lambda \in \mathbb{C}$, of $B$, we put
\begin{align*}
\frac{\partial \mathcal{N}_{\mathcal{F}_t}^*}{\partial t}:=\left\{\left(  w_j^\alpha(z_j,t) \mapsto - \overline{\sum_{\lambda=1}^m c_\lambda \frac{\partial w_j^\alpha(z_j,t)}{\partial t_\lambda}}     \right)_{\alpha=1,...,q} \right\} \in C^0\left( \mathcal{U}_t, \mathscr{H}om_{\mathcal{O}_{M_t}}\left(\mathcal{N}_{\mathcal{F}_t}^*, \frac{\Omega_{M_t}^1}{\mathcal{N}_{\mathcal{F}_t}^*} \right) \right)
\end{align*}

Then the foliated Kodaira-Spencer map in terms of cotangent sheaf at $t$ is defined to be a $\mathbb{C}$-linear map
{\small{\begin{align*}
\varphi_t:T_t(B)&\to \mathbb{H}^1\left(M_t, \mathcal{N}_{\mathcal{F}_t}^{*\bullet} \right)\\
                \frac{\partial}{\partial t}&\mapsto \frac{\partial \left(M_t, \mathcal{N}_{\mathcal{F}_t}^* \right)}{\partial t}:=\left( \rho_t\left(\frac{\partial}{\partial t} \right)= \frac{\partial M_t}{\partial t}= \left\{  \sum_{\alpha=1}^n \frac{\partial f_{jk}^\alpha(z_k,t)}{\partial t}\frac{\partial}{\partial z_j^\alpha} \right\}, \left\{ \left(w_j^\alpha(z_j,t)\mapsto -\overline{\frac{\partial w_j^\alpha(z_j,t)}{\partial t}}\right)_{\alpha=1,...,q} \right\} \right) \\
                 &\,\,\,\,\,\,\,\,\,\,\in C^1\left(\mathcal{U}_t, \Theta_{M_t} \right) \bigoplus C^0\left( \mathcal{U}_t, \mathscr{H}om_{\mathcal{O}_{M_t}}\left(\mathcal{N}_{\mathcal{F}_t}^*, \frac{\Omega_{M_t}^1}{\mathcal{N}_{\mathcal{F}_t}^*} \right) \right)
\end{align*}}}
where $\rho_t:T_t(B) \to H^1\left( M_t, \Theta_{M_t}\right)$ is the Kodaira-Spencer map at $t$ of the underlying complex analytic family $\left( \mathcal{M}, B, \pi \right)$ $($see \cite{Kod05} \textnormal{p.201}$)$.
\end{definition}

\section{Theorem of existence of deformations of foliated complex analytic structures in terms of cotangent sheaves}

\subsection{Preliminaries}\label{nn1}\

Let $\left(\mathcal{M}, \mathcal{N}_\mathcal{F}^* , B, \pi \right)$ with $\mathcal{N}_\mathcal{F}^*$ locally free be a foliated analytic family in terms of cotangent sheaves, where $B$ is a domain of $\mathbb{C}^m$ containing the origin $0$ as in Definition \ref{n3}. Define $|t|=\max_\lambda |t_\lambda|$ for $t=(t_1,...,t_m)\in \mathbb{C}^m$, and let $\Delta=\{t\in \mathbb{C}^m | |t|<r\}$ the polydisk  of radius $r>0$. If we take a sufficiently  small $\Delta \subset B$, then $\left(\mathcal{M}_\Delta, \mathcal{N}_{\mathcal{F}_\Delta}^* \right)=\pi^{-1}\left(\Delta\right)$ is represented in the form
\begin{align*}
\left(\mathcal{M}_\Delta, \mathcal{N}_{\mathcal{F}_\Delta}^* \right)=\bigcup_j \left(U_j\times \Delta, \mathcal{N}_\mathcal{F}^*|_{U_j\times \Delta}\right)
\end{align*}
We denote a point of $U_j$ by $\xi_j=\left(\xi_j^1,...,\xi_j^n\right)$ and $\Gamma\left(U_j\times \Delta, \mathcal{N}_\mathcal{F}^* \right)$ is generated by $w_j^1\left(\xi_j,t\right),...,w_j^q\left(\xi_j,t\right)$, where
\begin{align*}
w_j^\alpha\left(\xi_j,t\right)=\sum_{\beta=1}^n w_j^{\alpha\beta}\left(\xi_j, t\right) d\xi_j^\beta
\end{align*}
for some $w_j^{\alpha\beta}(\xi_j,t)\in \Gamma \left( U_j \times \Delta , \mathcal{O}_\mathcal{M} \right)$, and satisfies
\begin{align}\label{ni6}
w_j^1(\xi_j,t)\wedge \cdots \wedge w_j^q(\xi_j,t)\wedge dw_j^\alpha(\xi_j,t)=0,\,\,\,\,\,\,\,\alpha=1,...,q,\,\,\,\,\,\textnormal{on}\,\,\,\,\, U_j \times \Delta.
\end{align}

For simplicity, we assume that $U_j=\left\{\xi_j\in \mathbb{C}^m | |\xi_j|<1\right\}$ where $|\xi_j|=\max_a |\xi_j^a|$. $(\xi_j,t)\in U_j\times \Delta$ and $(\xi_k, t)\in U_k\times \Delta$ are the same point on $\mathcal{M}_\Delta$ if $\xi_j^\alpha=f_{jk}^\alpha(\xi_k,t),\alpha=1,...,n$ where $f_{jk}(\xi_k,t)$ is a holomorphic map of $\xi_k^1,..., \xi_k^n, t_1,...,t_m$, defined on $(U_k\times \Delta) \cap (U_j\times \Delta)$, and we have the following relation
\begin{align*}
w_j^\alpha(\xi_j,t)=\sum_{\beta=1}^q h_{jk}^{\alpha\beta}(\xi_k,t) w_k^\beta(\xi_k,t),\,\,\,\,\,\alpha=1,...,q
\end{align*}
for a holomorphic function $h_{jk}^{\alpha\beta}(\xi_k,t)$ on $(U_k\times \Delta) \cap ( U_j\times \Delta)$. More precisely
\begin{align*}
\sum_{\gamma,\eta=1}^n w_j^{\alpha\gamma}(f_{jk}(\xi_k,t),t)\frac{\partial f_{jk}^\gamma(\xi_k,t)}{\partial \xi_k^\eta} d\xi_k^\eta=\sum_{\beta=1}^q\sum_{\eta=1}^n h_{jk}^{\alpha\beta}(\xi_k,t)w_k^{\beta\eta}(\xi_k,t)d\xi_k^\eta
\end{align*}
Equivalently
\begin{align}\label{oci1}
\sum_{\gamma=1}^n w_j^{\alpha\gamma}(f_{jk}(\xi_k,t),t)\frac{\partial f_{jk}^\gamma(\xi_k,t)}{\partial \xi_k^\eta}=\sum_{\beta=1}^q h_{jk}^{\alpha\beta}(\xi_k,t)w_k^{\beta\eta}(\xi_k,t)
\end{align}
and on $\left(U_i\times \Delta \right)\cap \left(U_j \times \Delta \right) \cap \left( U_k \times \Delta \right)$, we have
\begin{align}
h_{ik}^{\alpha \gamma}(\xi_k,t)= \sum_{\beta=1}^q h_{ij}^{\alpha\beta}\left(f_{jk}(\xi_k,t),t\right) h_{jk}^{\beta \gamma}(\xi_k,t),\,\,\,\,\,\alpha,\gamma=1,...,q
\end{align}

By \cite{Kod05} Theorem 2.3, $M_t$ is diffeomorphic to $M_0=\pi^{-1}(0)$ as differentiable manifolds for each $t\in \Delta$. We put $M:=M_0$. By \cite{Kod05} Theorem 2.5, if we take a sufficiently small $\Delta$, there is a diffeomorphism $\Psi$ of $M\times \Delta$ onto $\mathcal{M}_\Delta$ as differentiable manifolds such that $\pi\circ \Psi$ is the projection $M\times \Delta \to \Delta$. Let $z=\left(z^1,..., z^n \right)$ be local coordinates of $M=M_0$. Then we have $\pi\circ \Psi(z,t)=t, t\in \Delta$. For $\Psi(z,t)\in U_j\times \Delta$, put
\begin{align*}
\Psi\left( z, t\right) = \left(\xi_j^1(z,t), ..., \xi_j^n(z,t), t_1,..., t_m \right)
\end{align*}
Then each component $\xi_j^\alpha=\xi_j^\alpha(z,t),\alpha=1,...,n$ is a $C^\infty$ function. If we identify $\mathcal{M}_\Delta =\Psi\left(M\times \Delta \right)$ with $M\times \Delta$ via $\Psi$, $\left(\mathcal{M}_\Delta,  \mathcal{N}_{\mathcal{F}_\Delta}^* \right)$ is considered as a complex manifold with the complex structure defined on the $C^\infty$ manifold $M\times \Delta$ by the system of local complex coordinates on $U_j\times \Delta$
\begin{align}\label{ne1}
\left\{ \left(\xi_j,t \right)| j=1 ,2, 3,...   \right\},\,\,\,\,\,(\xi_j,t)=\left( \xi_j^1(z,t), ..., \xi_j^n(z,t), t_1,..., t_m \right)
\end{align}
and holomorphic foliation $\mathcal{N}_{\mathcal{F}_\Delta}^*$ on $U_j\times \Delta$ with respect to the coordinates $(\ref{ne1})$ generated by $w_j^\alpha\left(\xi_j,t \right)= \sum_{\beta=1}^n w_j^{\alpha\beta}(\xi_j,t)d\xi_j^\beta ,\alpha=1,...,q$. 
We note that since $(z^1,..., z^n)$ and $\left(\xi_j^1(z,0), ..., \xi_j^n(z,0)\right)$ are local complex coordinates on $M=M_0$, $\xi_j^\alpha(z,0)$ are holomorphic functions of $z^1,...,z^n,\alpha=1,...,n$. We also note that if we take $\Delta$ sufficiently small, we have
\begin{align}\label{ne2}
\det\left( \frac{\partial \xi_j^\alpha(z,t)}{\partial z^\lambda} \right)_{\alpha, \lambda=1,...,n} \ne 0\,\,\,\,\,\, \, \textnormal{for}\,\,\,t\in \Delta
\end{align}

We put $\mathscr{U}_j=\Psi^{-1}\left(U_j\times \Delta \right)$. Then $\mathscr{U}_j \subset M\times \Delta$ is the domain of $\xi_j^{\alpha}(z,t)$. From $(\ref{ne2})$, we can define a $(0,1)$-form $\varphi_j^\lambda(z,t)= \sum_{v=1}^n \varphi_{jv}^\lambda(z,t) d\bar{z}^v$ for each $\lambda=1,...,n $ such that $\bar{\partial}\xi_j^\alpha =  \sum_{\lambda, v =1}^n \varphi_{jv}^\lambda(z,t)\frac{\partial \xi_j^\alpha}{\partial z^\lambda}d\bar{z}^v$. Then on $\mathscr{U}_j \cap \mathscr{U}_k$, we have $\sum_{\lambda=1}^n \varphi_j^\lambda(z,t)\frac{\partial}{\partial z^\lambda} = \sum_{\lambda=1}^n \varphi_k^\lambda(z,t) \frac{\partial}{\partial z^\lambda}$ (see \cite{Kod05} p.262). Then if for $(z,t)\in \mathscr{U}_j$, we define $\varphi(z,t):=\sum_{\lambda=1}^n \varphi_j^\lambda(z,t)\frac{\partial}{\partial z^\lambda}$, then $\varphi(t):=\varphi(z,t)\in A^{0,1}\left(M, \Theta_M \right)$ for $t\in \Delta$, and satisfies
\begin{align}\label{ne3}
\varphi(0)=0,\,\,\,\,\,\,\,\bar{\partial}\varphi(t)- \frac{1}{2}\left[ \varphi(t), \varphi(t) \right]=0
\end{align}
(see \cite{Kod05} p.263, p.265), We also point out that
\begin{theorem}\label{ne4}
If we take a sufficiently small polydisk $\Delta$ as above, then for $t\in \Delta$, a local $C^\infty$ function $f$ on $M$ is holomorphic with respect to the complex structure $M_t$ if and only if
\begin{align*}
\left( \bar{\partial} - \bar{\partial} \varphi(t)\right)f=0
\end{align*} 
\end{theorem}
\begin{proof}
See \cite{Kod05} Theorem 5.3 p.263
\end{proof}

We note that in the view of $(\ref{ne2})$ we let $\left(\psi_j^{\alpha\lambda}(z,t)\right)_{\alpha,\lambda=1,...,n}$ be the invertible matrix of $\left( \frac{\partial \xi_j^\alpha(z,t)}{\partial z_\lambda}\right)_{\alpha,\lambda=1,...,n}$, so that we have
\begin{align}\label{ca1}
\sum_{\gamma=1}^n \frac{\partial \xi_j^\eta(z,t)}{\partial z^\gamma} \psi_j^{\gamma\beta}(z_j,t)=\delta_\eta^\beta,\,\,\,\,\,\,\,\,\,\,\,\,\,\,\,\,\,\,\,\,\,\, \sum_{\gamma=1}^n \psi_j^{\eta \gamma}(z,t)\frac{\partial \xi_j^{\gamma}(z,t)}{\partial z^\beta} = \delta_\eta^\beta
\end{align}
We note that  $w_j^\alpha(\xi_j,t)=\sum_{\beta=1}^q w_j^{\alpha\beta}(\xi_j,t)d\xi_j^\beta$ implies that
\begin{align}\label{ne5}
w_j^\alpha \left(\xi_j,t\right)=\sum_{\beta,\gamma=1}^n w_j^{\alpha\beta}\left(\xi_j(z,t),t \right)\frac{\partial \xi_j^\beta}{\partial z^\gamma}dz^\gamma +\sum_{\beta,\gamma=1}^n w_j^{\alpha\beta}\left(\xi_j(z,t),t \right)\frac{\partial \xi_j^\beta}{\partial \bar{z}^\gamma} d\bar{z}^\gamma
\end{align}
From the first term $\sum_{\beta,\gamma=1}^n w_j^{\alpha\beta}(\xi_j,t)\frac{\partial \xi_j^\beta}{\partial z^\gamma} dz^\gamma$ in $(\ref{ne5})$, and from $(\ref{ca1})$, we can recover $w_j^\alpha(\xi_j,t)$ since we have $w_j^{\alpha\eta}=\sum_{\beta,\gamma=1}^n w_j^{\alpha\beta}\frac{\partial \xi_j^\beta}{\partial z^\gamma}\psi_j^{\gamma\eta}$. Keeping this in mind, we prove

\begin{theorem}
If we take a sufficiently small polydisk $\Delta$ as above, then for $t\in \Delta$, a $C^\infty$ $1$-form $w(z,t)=\sum_{\alpha=1}^n w^{\alpha}(z,t)dz^\alpha$ defines a holomorphic $1$-form $\sum_{\alpha,\beta=1}^n w^{\alpha}(z,t)\psi_j^{\alpha\beta}(z,t)d\xi_j^\beta$ with respect to the complex structure $M_t$ induced by $\varphi(t)$ if and only if it satisfies the equation
\begin{align*}
\bar{\partial}w(z,t)- \mathcal{L}_{\varphi(t)}\left(w(z,t)\right)=0
\end{align*}
\end{theorem}
For the definition of $\mathcal{L}_{\varphi(t)}$, see Appendix \ref{fb1}.

\begin{proof}
By Theorem \ref{ne4}, $\sum_{\alpha,\beta=1}^n w^\alpha(z,t)\psi_j^{\alpha\beta}(z,t)d\xi_j^\beta$ is holomorphic with respect to $\varphi(t)$  if and only if for each $\beta=1,...,n$, we have
\begin{align*}
0&=\bar{\partial} \left( \sum_{\alpha=1}^n w^\alpha \psi_j^{\alpha\beta}\right)-\left[\sum_{v,\lambda=1}^n \varphi_v^\lambda d\bar{z}^v \frac{\partial}{\partial z^\lambda}, \sum_{\alpha=1}^n w^\alpha \psi_j^{\alpha\beta} \right]\\
 &=\sum_{\alpha,v=1}^n \frac{\partial w^\alpha}{\partial \bar{z}^v} \psi_j^{\alpha\beta}d\bar{z}^v + \sum_{\alpha,v=1}^n w^\alpha \frac{\partial \psi_j^{\alpha\beta}}{\partial \bar{z}^v}  d\bar{z}^v -\sum_{v,\lambda,\alpha=1}^n \varphi_v^\lambda\frac{\partial w^\alpha}{\partial z^\lambda} \psi_j^{\alpha\beta} d\bar{z}^v-\sum_{v,\lambda,\alpha=1}^n \varphi_v^\lambda w^\alpha \frac{\partial \psi_j^{\alpha\beta}}{\partial z^\lambda}d\bar{z}^v
\end{align*}
which is equivalent to that for each $\beta,v =1,...,n$,
\begin{align} \label{ne6}
\sum_{\alpha=1}^n \frac{\partial w^\alpha}{\partial \bar{z}^v}\psi_j^{\alpha\beta} + \sum_{\alpha=1}^n w^\alpha \frac{\partial \psi_j^{\alpha\beta}}{\partial \bar{z}^v}-\sum_{\alpha,\lambda=1}^n \varphi_v^\lambda \frac{\partial w^\alpha}{\partial z^\lambda} \psi_j^{\alpha\beta} - \sum_{\alpha,\lambda=1}^n \varphi_v^\lambda w^\alpha \frac{\partial \psi_j^{\alpha\beta}}{\partial z^\lambda}=0
\end{align}

By multiplying $\sum_{\gamma=1}^n \frac{\partial \xi_j^\beta}{\partial z^\gamma}$ on the right hand side of $(\ref{ne6})$, we have
\begin{align}
\sum_{\alpha,\beta=1}^n \frac{\partial w^\alpha}{\partial \bar{z}^v}&\psi_j^{\alpha\beta}\frac{\partial \xi_j^\beta}{\partial z^\gamma} + \sum_{\alpha,\beta=1}^n w^\alpha \frac{\partial \psi_j^{\alpha\beta}}{\partial \bar{z}^v} \frac{\partial \xi_j^\beta}{\partial z^\gamma} -\sum_{\alpha,\lambda,\beta=1}^n \varphi_v^\lambda \frac{\partial w^\alpha}{\partial z^\lambda} \psi_j^{\alpha\beta}\frac{\partial \xi_j^\beta}{\partial z^\gamma} - \sum_{\alpha,\lambda,\beta=1}^n \varphi_v^\lambda w^\alpha \frac{\partial \psi_j^{\alpha\beta}}{\partial z^\lambda}\frac{\partial \xi_j^\beta}{\partial z^\gamma} =0 \notag \\
\iff &\frac{\partial w^\gamma}{\partial \bar{z}^v} +\sum_{\alpha,\beta=1}^n w^\alpha \frac{\partial \psi_j^{\alpha\beta}}{\partial \bar{z}^v}\frac{\partial \xi_j^\beta}{\partial z^\gamma} -\sum_{\lambda=1}^n \varphi_v^\lambda\frac{\partial w^\gamma}{\partial z^\lambda} -\sum_{\alpha,\lambda,\beta=1}^n w^\alpha \varphi_v^\lambda \frac{\partial \psi_j^{\alpha\beta}}{\partial z^\lambda} \frac{\partial \xi_j^\beta}{\partial z^\gamma}=0 \label{ne7}
\end{align}

On the other hand, we have
\begin{align}
&\bar{\partial}\left(\sum_{\alpha=1}^n w^\alpha dz^\alpha \right)- \mathcal{L}_{\sum_{v,\lambda=1}^n \varphi_v^\lambda d\bar{z}^v \frac{\partial}{\partial z^\lambda}}\left(\sum_{\alpha=1}^n w^\alpha dz^\alpha \right)=\bar{\partial}\left( \sum_{\alpha=1}^n w^\alpha dz^\alpha \right) -\sum_{v=1}^n d\bar{z}^v\wedge \mathcal{L}_{\sum_{\lambda=1}^n \varphi_v^\lambda \frac{\partial}{\partial z^\lambda}}\left( \sum_{\alpha=1}^n w^\alpha dz^\alpha \right) \notag \\
&=\sum_{\alpha, v=1}^n \frac{\partial w^\alpha}{\partial \bar{z}^v} d\bar{z}^v \wedge dz^\alpha-\sum_{v=1}^n d\bar{z}^v \wedge \left(\sum_{\lambda ,\alpha=1}^n \varphi_v^\lambda \frac{\partial w^\alpha}{\partial z^\lambda} dz^\alpha + \sum_{\alpha,\xi=1}^n w^\alpha \frac{\partial \varphi_v^\alpha}{\partial z^\xi} dz^\xi \right) \notag \\
&=\sum_{v,\alpha=1}^n \left(\frac{\partial w^\alpha}{\partial \bar{z}^v}-\sum_{\lambda =1}^n  \varphi_v^\lambda \frac{\partial w^\alpha}{\partial z^\lambda}  -\sum_{\lambda =1}^n w^\lambda \frac{\partial \varphi_v^\lambda}{\partial z^\alpha}   \right) d\bar{z}^v \wedge dz^\alpha = \sum_{v,\gamma=1}^n \left(\frac{\partial w^\gamma}{\partial \bar{z}^v}-\sum_{\lambda =1}^n  \varphi_v^\lambda \frac{\partial w^\gamma}{\partial z^\lambda}  -\sum_{\alpha =1}^n w^\alpha \frac{\partial \varphi_v^\alpha}{\partial z^\gamma}   \right) d\bar{z}^v \wedge dz^\gamma = 0 \notag  \\
&\iff \frac{\partial w^\gamma}{\partial \bar{z}^v}-\sum_{\lambda =1}^n  \varphi_v^\lambda \frac{\partial w^\gamma}{\partial z^\lambda}  -\sum_{\alpha =1}^n w^\alpha \frac{\partial \varphi_v^\alpha}{\partial z^\gamma} =0 \label{ne8}
\end{align}
In order for $(\ref{ne7})$ and $(\ref{ne8})$ to be equivalent, it is sufficient to show that
\begin{align}\label{ne9}
\frac{\partial \varphi_v^\alpha}{\partial z^\gamma}= - \sum_{\beta=1}^n \frac{\partial \psi_j^{\alpha\beta}}{ \partial \bar{z}^v}\frac{\partial \xi_j^\beta}{\partial z^\gamma}+ \sum_{\lambda, \beta=1}^n \varphi_v^\lambda  \frac{\partial \psi_j^{\alpha\beta}}{\partial z^\lambda}\frac{\partial \xi_j^\beta}{\partial z^\gamma}
\end{align}
In fact, we note that since $\varphi_v^\alpha= \sum_{\beta=1}^n \psi_j^{\alpha\beta} \frac{\partial \xi_j^\beta}{\partial \bar{z}^v}$, we have 
\begin{align} \label{ne8}
\frac{\partial \varphi_v^\alpha}{ \partial z^\gamma} &= \sum_{\beta=1}^n \frac{\partial \psi_j^{\alpha\beta}}{\partial z^\gamma} \frac{\partial \xi_j^\beta}{\partial \bar{z}^v}+ \sum_{\beta=1}^n \psi_j^{\alpha\beta}\frac{\partial^2 \xi_j^\beta}{\partial z^\gamma \partial \bar{z}^v} = \sum_{\beta, \lambda=1}^n \frac{\partial \psi_j^{\alpha\beta}}{\partial z^\gamma } \varphi_v^\lambda \frac{\partial \xi_j^\beta}{\partial z^\lambda} + \sum_{\beta=1}^n \psi_j^{\alpha\beta} \frac{\partial^2 \xi_j^\beta}{\partial z^\gamma \partial \bar{z}^v}
\end{align}
On the other hand, we note that
\begin{align*}
&\frac{\partial}{\partial \bar{z}^v}\left( \sum_{\beta=1}^n \psi_j^{\alpha\beta} \frac{\partial \xi_j^\beta}{\partial z^\gamma} \right) = 0 \Longrightarrow  \sum_{\beta=1}^n \frac{\partial \psi_j^{\alpha\beta}}{\partial \bar{z}^v}\frac{\partial \xi_j^\beta}{\partial z^\gamma} = - \sum_{\beta=1}^n \psi_j^{\alpha\beta} \frac{\partial \xi_j^\beta}{\partial \bar{z}^v \partial z^\gamma} \\
&\frac{\partial}{\partial z^\lambda}\left( \sum_{\beta=1}^n \psi_j^{\alpha\beta} \frac{\partial \xi_j^\beta}{\partial z^\gamma} \right) = 0 \Longrightarrow \sum_{\beta=1}^n \frac{\partial \psi_j^{\alpha\beta}}{\partial z^\lambda }\frac{\partial \xi_j^\beta}{\partial z^\gamma} = - \sum_{\beta=1}^n \psi_j^{\alpha\beta} \frac{\partial \xi_j^\beta}{\partial z^\lambda \partial z^\gamma}
\end{align*}
Then $(\ref{ne8})$ is equivalent to that
\begin{align*}
\frac{\partial \varphi_v^\alpha}{\partial z^\gamma}= - \sum_{\beta, \lambda=1}^n \varphi_v^\lambda \psi_j^{\alpha\beta} \frac{\partial^2 \xi_j^\beta }{\partial z^\gamma \partial z^\lambda} - \sum_{\beta=1}^n \frac{\partial \psi_j^{\alpha\beta}}{\partial \bar{z}^v} \frac{\partial \xi_j^\beta}{\partial z^\gamma} = \sum_{\lambda, \beta=1}^n \varphi_v^\lambda \frac{\partial \psi_j^{\alpha\beta}}{\partial z^\lambda}\frac{\partial \xi_j^\beta}{\partial z^\gamma  } - \sum_{\beta=1}^n \frac{\partial \psi_j^{\alpha\beta}}{\partial \bar{z}^v} \frac{\partial \xi_j^\beta}{\partial z^\gamma}
\end{align*}
which implies $(\ref{ne9})$.

\end{proof}

With this preparation, we shall prove theorem of existence for deformations of foliated complex analytic structures in terms of cotangent sheaves.

\subsection{Theorem of existence for deformations of foliated complex analytic structures in terms of cotangent sheaves}

\begin{theorem}[Theorem of existence for deformations of foliated complex analytic structures in terms of cotangent sheaves] \label{nn2}
Let $\left( M , \mathcal{N}_{\mathcal{F}_0}^* \right)$ be a compact foliated complex manifold with $\mathcal{N}_{\mathcal{F}_0}^*$ locally free. Suppose that $\mathbb{H}^2\left( M,  \mathcal{N}_{\mathcal{F}_0}^{ * \bullet} \right)=0$. Then there exists a foliated analytic family $\left( \mathcal{M}, \mathcal{N}_\mathcal{F}^* , B, \pi \right)$ with $0\in B \subset \mathbb{C}^m$ satisfying the following conditions$:$
\begin{enumerate}
\item $\pi^{-1}(0)=\left( M, \mathcal{N}_{\mathcal{F}_0}^* \right)$
\item The foliated Kodaira-Spencer map $\varphi_0: T_0 B \to \mathbb{H}^1\left( M, \mathcal{N}_{\mathcal{F}_0}^{* \bullet} \right)$ in terms of cotangent sheaf is an isomorphism.
\end{enumerate}
\end{theorem}

\begin{proof}
We may assume tha following:
\begin{enumerate}
\item $M$ is covered by a finite number of coordinate neighborhoods $U_i(i\in I)$ with a system of coordinates $z_i=\left(z_i^1,...,z_i^n\right)$ such that $U_i=\left\{z_i\in \mathbb{C}^n|\max_{1\leq \alpha \leq n} | z_i^\alpha| <1 \right\}$ \label{nnc10}
\item $\mathcal{N}_{\mathcal{F}_0}^*$ is generated by $\left\{w_{0i}^1,...,w_{0i}^q \right\}$ with $w_{0i}^\alpha\in \Gamma\left(U_i, \mathcal{N}_{\mathcal{F}_0}^*\right)$ and $w_{0j}^\alpha=\sum_{\beta=1}^p h_{0ij}^{\alpha\beta} w_{0j}^\beta$ on $U_i\cap U_j$ for $h_{0ij}^{\alpha\beta}(z_j)\in \Gamma\left( U_{ij}, \mathcal{O}_M \right),\alpha=1,...,q$.
\item $w_{0i}^1\wedge \cdots \wedge w_{0i}^q\wedge dw_{0i}^\alpha=0$ for $\alpha=1,...,q$. 
\item $z_i$ coincides with $z_j$ if and only if $z_i=f_{ij}(z_j)$.
\item for $x\in U_i-S$, where $S=\textnormal{Sing}\left(\mathcal{F}_0\right)$, there is a neighborhood $U_{i_x} \subset U_i-S$ of $x$ such that $dw_{0i}^\alpha =\sum_{\beta=1}^q a_{0i_x}^{\alpha\beta}  \wedge w_{0i}^\beta$ on $U_{i_x}$ for some $a_{0i_x}^{\alpha\beta} \in \Gamma\left(U_{i_x}, \Omega_M^1 \right)$. \label{nnc11}
\end{enumerate}
Then we note that from $\sum_{\gamma=1}^q h_{0ik}^{\alpha\gamma}w_{0k}^\gamma=w_{0i}^\alpha=\sum_{\beta=1}^q h_{0ij}^{\alpha\beta} w_{0j}^\beta=\sum_{\beta,\gamma=1}^q h_{0ij}^{\alpha\beta} h_{0jk}^{\beta\gamma} w_{0k}^\gamma$, we have
\begin{align*}
h_{0ik}^{\alpha\gamma}=\sum_{\beta=1}^q h_{0ij}^{\alpha\beta} h_{0jk}^{\beta\gamma}
\end{align*}

Let $r=\dim_\mathbb{C} \mathbb{H}^1\left(M, \mathcal{N}_{\mathcal{F}_0}^{*\bullet} \right)$ and $B= \left\{t= \left(t_1,...,t_r \right)\in \mathbb{C}^r||t|<\epsilon \right\}$ with sufficiently small $\epsilon > 0$. Our purpose is to construct
\begin{enumerate}
\item a $C^\infty$ $(0,1)$-form $\varphi(t)$ with coefficients in $\Theta_M$ depending holomorphically on $t$,
\item $C^\infty$-differential forms $w_i^\alpha(z_i,t)$ on $U_i\times B$ of the form $w_i^\alpha(z_i,t)=\sum_{\beta=1}^q w_i^{\alpha\beta}(z_i,t) dz_i^\beta$ depending holomorphically on $t$,
\item $q\times q$ matrices $H_{ij}:=\left( h_{ij}^{\alpha\beta}(z_j,t) \right)$ with components $C^\infty$ function $h_{ij}^{\alpha\beta}(z_j,t)$ on $U_{ij}\times B$ depending holomorphically on $t$ 
\item If $x=(x,0)\in U_i-S$ is a regular point of $\mathcal{F}_0$, then $(z_i,t)$ is also regular for some neighborhood of $(x,0)$, so that we will construct $w_i^\alpha(z_i,t)$ satisfying $dw_i^\alpha(z_i,t)=\sum_{\beta=1}^q a_{i_x}^{\alpha\beta}(z_i,t)\wedge w_i^\beta(z_i,t)$ on $U_{i_x}\times B_{i_x}$ where $ a_{i_x}^{\alpha\beta}(z_i,t)=\sum_{\gamma=1}^n a_{i_x}^{\alpha\beta\gamma}(z_i,t)dz_i^\gamma$ on $U_{i_x}\times B_{i_x}\subset U_i\times B$, depending holomorphically on $t \in B_{i_x}$,  such that
\end{enumerate}
\begin{align}
\varphi(0)=0,\,\,\,\,\,\,\,\,\,\,\, w_i^\alpha(z_i,&0)=w_{0i}^\alpha,\,\,\,\,\,\,\,\,\, h_{ij}^{\alpha\beta}(z_j,0)= h_{0ij}^{\alpha\beta}(z_j),\,\,\,\,\,\,\,\,a_{i_x}^{\alpha\beta}(z_i,0)=a_{0i_x}^{\alpha\beta}(z_i)\,\,\,\,\,\textnormal{on $U_{i_x}$}  \label{n20}\\
&\bar{\partial}- \frac{1}{2}\left[\varphi,\varphi \right]=0 \label{n22}\\
\bar{\partial} w_i^\alpha(z_i,t)&-\mathcal{L}_\varphi\left( w_i^\alpha(z_i,t)\right)=0 \label{ns54}\\
\bar{\partial} h_{ij}^{\alpha\beta}(z_j,t)&-\left[\varphi, h_{ij}^{\alpha\beta} (z_j, t) \right]=0 \label{ns55}\\
w_i^\alpha(z_i,t)&=\sum_{\beta=1}^q h_{ij}^{\alpha\beta}(z_j,t) w_j^\beta(z_j,t) \label{n212}\\
w_i^1(z_i,t)\wedge \cdots  \wedge& w_i^q(z_i,t)\wedge \partial w_i^\alpha(z_i,t)=0 \label{nnc5} \\
H_{ij}H_{jk}=H_{ik} \iff &h_{ik}^{\alpha\gamma}(z_k,t)=\sum_{\beta=1}^q h_{ij}^{\alpha\beta} \left(f_{jk}(z_k,t) , t \right) h_{jk}^{\beta\gamma}(z_k,t)\,\,\,\,\,\,\textnormal{on}\,\,\,\,\,\,(U_i\cap U_j\cap U_k)\times B \label{n21} \\
\partial w_i^\alpha(z_i,t)=\sum_{\beta=1}^q &a_{i_x}^{\alpha\beta}(z_i,t)\wedge w_i^\beta(z_i,t),\,\,\,\,\,\,\,  \textnormal{on $U_{i_x}\times B_{i_x},\,\,\,\,$ where $x=(x,0)\in U_{i_x}\subset U_i-S$}\label{n211}
\end{align}

\subsection{Existence of formal solutions}\

We construct solutions of $(\ref{n20})-(\ref{n211})$ which are formal power series in $t$. We recall Notation \ref{te27}. We write $\varphi=\sum_{\mu=0}^\infty \varphi_\mu, h_{ij}^{\alpha\beta \mu}=\sum_{\mu=0}^\infty h_{ij|\mu}^{\alpha\beta}, w_i^\alpha= \sum_{\mu=0}^\infty w_{i|\mu}^\alpha$ and $a_{i_x}^{\alpha\beta}= \sum_{\mu=0}^\infty a_{i_x|\mu}^{\alpha\beta}$. In view of $(\ref{n20})$, we set $\varphi_0=0, w_{i|0}^\alpha= w_{0i}^\alpha(z_i), h_{ij|0}^{\alpha\beta}(z_i,t)= h_{0ij}^{\alpha\beta}(z_j)$ and $a_{i_x|0}^{\alpha\beta}= a_{0i_x}^{\alpha\beta}$.  Then $(\ref{n22})$-$(\ref{n211})$ are equivalent to the following system of congruences
\begin{align}
\bar{\partial}\varphi^\mu& -\frac{1}{2}\left[ \varphi^\mu, \varphi^\mu \right]\equiv_\mu 0 \label{n23}\\
\bar{\partial} w_i^\alpha(z_i,t)&- \mathcal{L}_{\varphi^\mu}\left(w_i^\alpha(z_i,t) \right)\equiv_\mu 0 \label{nnc4}\\
\bar{\partial} h_{ij}^{\alpha\beta\mu}-&\left[\varphi^\mu, h_{ij}^{\alpha\beta \mu} \right]\equiv_\mu 0 \label{nnc8} \\
w_i^{\alpha\mu}&\equiv_\mu \sum_{\beta=1}^q h_{ij}^{\alpha\beta\mu} w_j^{\beta \mu} \label{n213}\\
w_i^{1\mu}\wedge \cdots&  \wedge w_i^{q\mu}\wedge \partial w_i^{\alpha\mu}\equiv_\mu 0\\
h_{ik}^{\alpha\gamma\mu}(z_k,t)&\equiv_\mu \sum_{\beta=1}^q h_{ij}^{\alpha\beta \mu}\left(f_{jk},t \right) h_{jk}^{\beta\gamma \mu}(z_k,t) \label{n24} \\
\partial w_i^{\alpha\mu}\equiv_\mu \sum_{\beta=1}^qa_{i_x}^{\alpha\beta\mu}&(z_i,t)\wedge w_i^{\beta\mu}(z_i,t)\,\,\,\,\,\,\,\,\textnormal{on $U_{i_x}\subset U_i-S$} \label{n241}
\end{align}
We define homogeneous polynomials
\begin{align}
\xi_\mu &\equiv_\mu \bar{\partial} \varphi^{\mu-1}-\frac{1}{2}\left[\varphi^{\mu-1}, \varphi^{\mu-1} \right] \label{nn3}
\end{align}
\begin{align}
-A_{i|\mu}^\alpha &\equiv_\mu \bar{\partial} w_i^{\alpha(\mu-1)}- L_{\varphi^{\mu-1}}\left( w_i^{\alpha(\mu-1)}\right) \label{nn7}
\end{align}
\begin{align}
-\sum_{\eta=1}^p B_{ij|\mu}^{\alpha\eta} h_{0ij}^{\eta \beta}&\equiv_\mu \bar{\partial} h_{ij}^{\alpha\beta(\mu-1)}-\left[\varphi^{\mu-1}, h_{ij}^{\alpha\beta(\mu-1)} \right]
\end{align}
\begin{align}
C_{ij|\mu}^\alpha &\equiv_\mu w_i^{\alpha(\mu-1)}-\sum_{\beta=1}^q h_{ij}^{\alpha\beta(\mu-1)} w_j^{\beta(\mu-1)} \label{nn6}
\end{align}
\begin{align}
\sum_{\beta=1}^q D_{ijk|\mu}^{\alpha\beta} h_{0ik}^{\beta\gamma}&\equiv_\mu h_{ik}^{\alpha\gamma(\mu-1)}-\sum_{\beta=1}^q r_{ij}^{\alpha\beta(\mu-1)}(f_{jk},t) r_{jk}^{\beta\gamma(\mu-1)}(z_k,t)\label{n215}
\end{align}
\begin{align}
E_{i|\mu}^\alpha&\equiv_\mu w_i^{1(\mu-1)}\wedge \cdots \wedge w_i^{q(\mu-1)}\wedge \partial w_i^{\alpha(\mu-1)} \label{ns14}
\end{align}
\begin{align}
L_{i_x|\mu}^\alpha&\equiv_\mu \partial w_i^{\alpha(\mu-1)}-\sum_{\beta=1}^q a_{i_x}^{\alpha\beta(\mu-1)}\wedge w_i^{\beta(\mu-1)} \,\,\,\,\,\,\,\,\,\textnormal{on $U_{i_x}\subset U_i- S$} \label{nn4}
\end{align}

\begin{lemma} \label{ne10}
We have the following equalities$:$
\begin{align}
\bar{\partial}\xi_\mu=0 \label{n25}
\end{align}
\begin{align}
\sum_{\beta,\eta=1}^q h_{0ij}^{\alpha\beta} D_{jkl|\mu}^{\beta \eta} h_{0ji}^{\eta \xi}-D_{ijk|\mu}^{\alpha\xi}+D_{ijl|\mu}^{\alpha \xi} - D_{ijk|\mu}^{\alpha\xi}=0\label{n26}
\end{align}
\begin{align}
\bar{\partial} A_{i|\mu}^\alpha= \mathcal{L}_{\xi_\mu}\left(w_{0i}^\alpha\right) \label{n29}
\end{align}
\begin{align}
\bar{\partial} B_{ij|\mu}^{\alpha\xi}=\sum_{\beta=1}^p \left[\xi_\mu,h_{0ij}^{\alpha\beta} \right] h_{0ji}^{\beta \xi}\label{n291}
\end{align}
\begin{align}
\sum_{\eta,\gamma=1}^q h_{0ij}^{\alpha\eta} B_{jk|\mu}^{\eta\gamma} h_{0ji}^{\gamma \xi} -B_{ik|\mu}^{\alpha \xi}+B_{ij|\mu}^{\alpha \xi}=\bar{\partial} D_{ijk|\mu}^{\alpha \xi} \label{n30}
\end{align}
\begin{align}
\sum_{\gamma=1}^q h_{0ij}^{\alpha\gamma} C_{jk|\mu}^\gamma-C_{ik|\mu}^\alpha+C_{ij|\mu}^\alpha=\sum_{\gamma=1}^q D_{ijk|\mu}^{\alpha\gamma} w_{0i}^\gamma \label{n31}
\end{align}
\begin{align}
\sum_{\beta=1}^q h_{0ij}^{\alpha\beta} A_{j|\mu}^\beta-A_{i|\mu}^\alpha=\bar{\partial} C_{ij|\mu}^\alpha-\sum_{\eta=1}^q B_{ij|\mu}^{\alpha\eta} w_{0i}^\eta\label{n28}
\end{align}
\begin{align}
\bar{\partial} E_{i|\mu}^\alpha =  \sum_{\beta=1}^q (-1)^\beta w_{0i}^1\wedge \cdots \wedge  \overbrace{ A_{i|\mu}^\beta}^{\beta- \textnormal{th}} \wedge \cdots \wedge w_{0i}^q \wedge dw_{0i}^\alpha + (-1)^q w_{0i}^1 \wedge \cdots \wedge w_{0i}^q \wedge \partial A_{i|\mu}^\alpha \label{n407}
\end{align}
\begin{align}
E_{i|\mu}^\alpha-\sum_{\beta=1}^q h_{0ij}^{\alpha\beta} \det\left( h_{0ij}^{\gamma\eta}\right)E_{j|\mu}^\beta&= \left(\sum_{\beta=1}^q w_{0i}^1\wedge \cdots \wedge \overbrace{C_{ij|\mu}^\beta}^{\beta-\textnormal{th}} \wedge \cdots \wedge w_{0i}^q \right) \wedge dw_{0i}^\alpha+ w_{0i}^1\wedge \cdots \wedge w_{0i}^q \wedge \partial C_{ij|\mu}^\alpha       \label{n417}
\end{align}
\begin{align}
E_{i|\mu}^\alpha&=w_{0i}\wedge L_{i_x|\mu}^\alpha= w_{0i}^1 \wedge \cdots \wedge w_{0i}^q \wedge L_{i_x|\mu}^\alpha \,\,\,\,\,\,\,\,\textnormal{on $U_{i_x}\subset U_i-S$} \label{n53}
\end{align}
\end{lemma}

\begin{proof}
$(\ref{n25})$ follows from \cite{Kod05} p.273. $(\ref{n26})$ follows from the same way with $(\ref{t20})$. We prove $(\ref{n29})$. In fact, from Appendix \ref{fb1},
\begin{align*}
\bar{\partial} A_{i|\mu}^\alpha&\equiv_\mu \bar{\partial}\left( \mathcal{L}_{\varphi^{\mu-1}}\left( w_i^{\alpha(\mu-1)}\right)\right)= \mathcal{L}_{\bar{\partial}\varphi^{\mu-1}}\left( w_i^{\alpha(\mu-1)}\right) - \mathcal{L}_{\varphi^{\mu-1}}\left( \bar{\partial} w_i^{\alpha(\mu-1)}\right)\\
&\equiv_\mu \mathcal{L}_{\xi_\mu+\frac{1}{2}\left[\varphi^{\mu-1},\varphi^{\mu-1} \right]}\left( w_i^{\alpha(\mu-1)} \right) - \mathcal{L}_{\varphi^{\mu-1}}\left( \mathcal{L}_{\varphi^{\mu-1}}\left(w_i^{\alpha(\mu-1)} \right) -A_{i|\mu}^\alpha  \right)\\
&= \mathcal{L}_{\xi_\mu}\left( w_{0i}^\alpha \right) + \frac{1}{2} \mathcal{L}_{\left[\varphi^{\mu-1}, \varphi^{\mu-1}\right]}\left(w_i^{\alpha(\mu-1)} \right) - \mathcal{L}_{\varphi^{\mu-1}}\mathcal{L}_{\varphi^{\mu-1}}\left(w_i^{\alpha(\mu-1)} \right)= \mathcal{L}_{\xi_\mu}\left(w_{0i}^\alpha \right)
\end{align*}
$(\ref{n291})$ follows from the same way with $(\ref{t23})$. $(\ref{n30})$ follows from the same way with $(\ref{t24})$. $(\ref{n31})$ follows from the same way with $(\ref{t25})$. $(\ref{n28})$ follows from the same way with $(\ref{t26})$. We prove $(\ref{n407})$. In fact,
{\small{\begin{align*}
\bar{\partial} E_{i|\mu}^\alpha&\equiv_\mu \sum_{\beta=1}^q (-1)^{\beta -1} w_i^{1(\mu-1)}\wedge \cdots \wedge \bar{\partial} w_i^{\beta(\mu-1)}\wedge \cdots \wedge w_i^{q(\mu-1)}\wedge \partial w_i^{\alpha(\mu-1)} -(-1)^qw_i^{1(\mu-1)}\wedge \cdots \wedge w_i^{q(\mu-1)}\wedge \partial \left(\bar{\partial}w_i^{\alpha(\mu-1)}\right)\\
&\equiv_\mu \sum_{\beta=1}^q (-1)^{\beta-1} w_i^{1(\mu-1)}\wedge \cdots \wedge \left(\mathcal{L}_{\varphi^{\mu-1}}\left(w_i^{\beta(\mu-1)} \right)-A_{i|\mu}^\beta \right)\wedge \cdots \wedge w_i^{q(\mu-1)}\wedge \partial w_i^{\alpha(\mu-1)}\\
&-(-1)^{q}w_i^{1(\mu-1)}\wedge \cdots \wedge w_i^{q(\mu-1)}\wedge \partial \left( \mathcal{L}_{\varphi^{\mu-1}}\left( w_i^{\alpha(\mu-1)}\right)-A_{i|\mu}^\alpha\right)\\
&\equiv_\mu \sum_{\beta=1}^q (-1)^\beta w_{0i}^1\wedge \cdots \wedge A_{i|\mu}^\beta \wedge \cdots \wedge w_{0i}^q\wedge dw_{0i}^\alpha +(-1)^q w_{0i}^1\wedge \cdots \wedge w_{0i}^q\wedge \partial A_{i|\mu}^\alpha + \mathcal{L}_{\varphi^{\mu-1}}\left( w_i^{1(\mu-1)}\wedge \cdots \wedge w_i^{q(\mu-1)}\wedge \partial w_i^{\alpha(\mu-1)}\right)\\
&\equiv_\mu  \sum_{\beta=1}^q (-1)^\beta w_{0i}^1\wedge \cdots \wedge A_{i|\mu}^\beta \wedge \cdots \wedge w_{0i}^q\wedge dw_{0i}^\alpha +(-1)^q w_{0i}^1\wedge \cdots \wedge w_{0i}^q\wedge dA_{i|\mu}^\alpha
\end{align*}}}

We prove $(\ref{n417})$. In fact,
{\small{\begin{align*}
&E_{i|\mu}^\alpha-\sum_{\beta=1}^q h_{0ij}^{\alpha\beta} \det\left( h_{0ij}^{\gamma\eta}\right)E_{j|\mu}^\beta\\
&\equiv_\mu w_i^{1(\mu-1)}\wedge \cdots \wedge w_i^{q(\mu-1)}\wedge \partial w_i^{\alpha(\mu-1)}-\sum_{\beta=1}^qh_{ij}^{\alpha\beta (\mu-1)} \det\left(h_{ij}^{\gamma\eta(\mu-1)}\right) w_j^{1(\mu-1)}\wedge \cdots \wedge w_j^{q(\mu-1)}\wedge \partial w_j^{\beta(\mu-1)}\\
&\equiv_\mu \left(C_{ij|\mu}^1+\sum_{\gamma_1=1}^q h_{ij}^{1\gamma_1(\mu-1)}w_j^{\gamma_1(\mu-1)} \right)\wedge\cdots \wedge \left( C_{ij|\mu}^q+\sum_{\gamma_q=1}^q h_{ij}^{q\gamma_q(\mu-1)} w_j^{\gamma_q(\mu-1)} \right)\wedge \partial \left(C_{ij|\mu}^\alpha+\sum_{\beta=1}^q h_{ij}^{\alpha\beta(\mu-1)} w_j^{\beta(\mu-1)} \right)\\
&-\sum_{\beta=1}^qh_{ij}^{\alpha\beta (\mu-1)} \det\left(h_{ij}^{\gamma\eta(\mu-1)}\right) w_j^{1(\mu-1)}\wedge \cdots \wedge w_j^{q(\mu-1)}\wedge \partial w_j^{\beta(\mu-1)}\\
&\equiv_\mu \left(\sum_{\beta=1}^q w_{0i}^1\wedge \cdots \wedge \overbrace{C_{ij|\mu}^\beta}^{\beta-\textnormal{th}} \wedge \cdots \wedge w_{0i}^q \right) \wedge dw_{0i}^\alpha+ w_{0i}^1\wedge \cdots \wedge w_{0i}^q \wedge \partial C_{ij|\mu}^\alpha 
\end{align*}}}

We prove $(\ref{n53})$. In fact,
We note that
\begin{align*}
E_{i|\mu}^\alpha&\equiv_\mu w_i^{1(\mu-1)}\wedge \cdots \wedge w_i^{q(\mu-1)}\wedge \partial w_i^{\alpha(\mu-1)}\\
 &\equiv_\mu w_i^{1(\mu-1)}\wedge \cdots \wedge w_i^{q(\mu-1)} \wedge \left(\sum_{\gamma=1}^q a_{i_x}^{\alpha\gamma(\mu-1)}\wedge w_i^{\gamma(\mu-1)}+L_{i_x|\mu}^\alpha \right)\equiv_\mu w_{0i}^1 \wedge \cdots \wedge w_{0i}^q \wedge L_{i_x|\mu}^\alpha
\end{align*}
This completes the proof of Lemma \ref{ne10}.

\end{proof}

Our purpose is to construct $\varphi^\mu=\varphi^{\mu-1}+ \varphi_\mu, h_{ij}^{\alpha\beta \mu}= h_{ij}^{\alpha\beta (\mu-1)} + h_{ij|\mu}^{\alpha\beta}, w_i^{\alpha \mu}=w_i^{\alpha (\mu-1)}+ w_{i|\mu}^\alpha$, and $a_{i_x}^{\alpha\beta \mu}= a_{i_x}^{\alpha\beta(\mu-1)}+ a_{i_x|\mu}^{\alpha\beta}$ which satisfy $(\ref{n23})_\mu-(\ref{n241})_\mu$.

\begin{lemma}\label{nn5}
$(\ref{n23})_\mu-(\ref{n241})_\mu$ are equivalent to the following equalities, respectively:
\begin{align}
\bar{\partial} \varphi_\mu&=- \xi_\mu \label{n42}\\
A_{i|\mu}^\alpha &= \bar{\partial} w_{i|\mu}^\alpha- \mathcal{L}_{\varphi_\mu}\left( w_{0i}^\alpha \right)       \label{n43}  \\
\sum_{\eta=1}^q B_{ij|\mu}^{\alpha\eta} h_{0ij}^{\eta \beta}&= \bar{\partial} h_{ij|\mu}^{\alpha\beta}- \left[ \varphi_\mu,  h_{0ij}^{\alpha\beta} \right]     \label{n45} \\ 
 C_{ij|\mu}^\alpha & = -w_{i|\mu}^\alpha +\sum_{\beta=1}^q h_{0ij}^{\alpha\beta} w_{j|\mu}^\beta + \sum_{\beta=1}^q h_{ij|\mu}^{\alpha\beta} w_{0j}^\beta  \label{n46} \\
 -E_{i|\mu}^\alpha &= \sum_{\beta=1}^q w_{0i}^1\wedge \cdots \wedge  \overbrace{ w_{i|\mu}^\beta}^{\beta-\textnormal{th}} \wedge \cdots \wedge w_{0i}^q \wedge d w_{0i}^\alpha + w_{0i}^1\wedge \cdots \wedge w_{0i}^q \wedge \partial w_{i|\mu}^\alpha \label{n47} \\
 D_{ijk|\mu}^{\alpha \xi} &=\sum_{\beta, \gamma =1}^q h_{0ij}^{\alpha\beta} h_{jk|\mu}^{\beta \gamma} h_{0ki}^{\gamma \xi}- \sum_{\gamma=1}^q h_{0ik|\mu}^{\alpha \gamma} h_{0ki}^{\gamma \xi}+ \sum_{\beta=1}^q h_{ij|\mu}^{\alpha\beta} r_{0ji}^{\beta \xi} \label{nn121} \\
L_{i_x|\mu}^\alpha &= - \partial w_{i|\mu}^\alpha + \sum_{\beta=1}^q a_{0i_x}^{\alpha\beta}\wedge w_{i|\mu}^\beta + \sum_{\beta=1}^q a_{i_x|\mu}^{\alpha\beta} \wedge w_{0i}^\beta \label{n48}
\end{align}
\end{lemma}

\begin{proof}
$(\ref{n42})$ follows from \cite{Kod05} p.272. We prove $(\ref{n43})$. In fact,
\begin{align*}
0&\equiv_\mu \bar{\partial}\left(w_i^{\alpha (\mu-1)}+ w_{i|\mu}^\alpha \right) - \mathcal{L}_{ \varphi^{\mu-1} +\varphi_\mu}\left( w_i^{\alpha (\mu-1)}+ w_{i|\mu}^\alpha \right)\\
 &\equiv_\mu \bar{\partial} w_i^{\alpha(\mu-1)} - \mathcal{L}_{\varphi^{\mu-1}}\left(w_i^{\alpha(\mu-1)}\right) + \bar{\partial} w_{i|\mu}^\alpha- \mathcal{L}_{\varphi_\mu}(w_{0i}^\alpha) \iff A_{i|\mu}^\alpha = \bar{\partial} w_{i|\mu}^\alpha- \mathcal{L}_{\varphi_\mu}\left(w_{0i}^\alpha\right)
\end{align*}
$(\ref{n45})$ follows from the same way with $(\ref{te13})$. $(\ref{n46})$ follows from the same way with $(\ref{te15})$. We prove $(\ref{n47})$. In fact,
\begin{align*}
0\equiv_\mu&\left(w_i^{1(\mu-1)}+ w_{i|\mu}^1\right)\wedge \cdots \wedge \left(w_i^{q(\mu-1)}+ w_{i|\mu}^q  \right)  \wedge \partial \left( w_i^{\alpha(\mu-1)}+ w_{i|\mu}^\alpha  \right) \\
&\equiv_\mu E_{i|\mu}^\alpha + \sum_{\beta=1}^q w_{0i}^1\wedge \cdots \wedge  \overbrace{ w_{i|\mu}^\beta}^{\beta-\textnormal{th}} \wedge \cdots \wedge w_{0i}^q \wedge d w_{0i}^\alpha + w_{0i}^1\wedge \cdots \wedge w_{0i}^q \wedge \partial w_{i|\mu}^\alpha
\end{align*}
$(\ref{nn121})$ follows from the same way with $(\ref{te17})$. We prove $(\ref{n48})$. In fact,
\begin{align*}
&\partial\left(w_i^{\alpha(\mu-1)}+ w_{i|\mu}^\alpha \right)\equiv_\mu \sum_{\beta=1}^q \left( a_{i_x}^{\alpha\beta(\mu-1)}+ a_{i_x|\mu}^{\alpha\beta}   \right) \wedge \left( w_i^{\beta(\mu-1)}+ w_{i|\mu}^\beta    \right)\\
&  \iff  L_{i_x|\mu}^\alpha= - \partial w_{i|\mu}^\alpha + \sum_{\beta=1}^q a_{0i_x}^{\alpha\beta}\wedge w_{i|\mu}^\beta + \sum_{\beta=1}^q a_{i_x|\mu}^{\alpha\beta} \wedge w_{0i}^\beta
\end{align*}

\end{proof}

\begin{lemma}\label{ne25}
Under the hypothesis $\mathbb{H}^2\left( M , \mathcal{N}_{\mathcal{F}_0}^{*\bullet} \right)=0$,  we can find $\varphi_\mu,  h_{ij|\mu}^{\alpha\beta}, w_{i|\mu}^\alpha$ and $a_{i_x|\mu}^{\alpha\beta}$ which satisfy $(\ref{n42})- (\ref{n48})$.
\end{lemma}

\begin{proof}
We define $D_{ijk|\mu}\in \Gamma\left( U_{ijk}, \mathcal{A}^{0,0}\left( \mathscr{H}om_{\mathcal{O}_M}\left( \mathcal{N}_{\mathcal{F}_0}^*, \mathcal{N}_{\mathcal{F}_0}^* \right) \right) \right)$ by
\begin{align*}
D_{ijk|\mu}:\Gamma\left( U_{ijk},  \mathcal{N}_{\mathcal{F}_0}^* \right) &\to \Gamma\left( U_{ijk}, \mathcal{A}^{0,0}\left( \mathcal{N}_{\mathcal{F}_0}^* \right) \right)\\
     w_{0i}^\alpha &\mapsto \sum_{\xi=1}^q D_{ijk|\mu}^{\alpha \xi} w_{0i}^\xi
\end{align*}
and linearly extends to $\Gamma\left( U_{ijk}, \mathcal{N}_{\mathcal{F}_0}^* \right)$. Then by $(\ref{n26})$, as in the same way with $(\ref{t201})$, we can find $\left\{D_{ij|\mu} \right\}\in C^1\left( \mathcal{U}, \mathcal{A}^{0,0}\left(\mathscr{H}om_{\mathcal{O}_M}\left( \mathcal{N}_{\mathcal{F}_0}^*, \mathcal{N}_{\mathcal{F}_0}^* \right)   \right) \right)$
with $D_{ij|\mu}\in \Gamma\left( U_{ij}, \mathcal{A}^{0,0}\left( \mathscr{H}om_{\mathcal{O}_M}\left( \mathcal{N}_{\mathcal{F}_0}^*, \mathcal{N}_{\mathcal{F}_0}^* \right) \right) \right)$ defined by
\begin{align*}
D_{ij|\mu}:\Gamma\left( U_{ij} , \mathcal{N}_{\mathcal{F}_0}^* \right) &\to \Gamma\left( U_{ij} , \mathcal{A}^{0,0}\left( \mathcal{N}_{\mathcal{F}_0}^* \right) \right) \\
 w_{0i}^\alpha &\mapsto \sum_{\xi=1}^p D_{ij|\mu}^{\alpha \xi} w_{0i}^\xi
\end{align*}
such that
\begin{align}\label{ne11}
 \sum_{\beta,\eta=1}^q h_{0ij}^{\alpha\beta} D_{jk|\mu}^{\beta \eta} h_{0ji}^{\eta\xi}-D_{ik|\mu}^{\alpha \xi}+ D_{ij|\mu}^{\alpha \xi}=D_{ijk|\mu}^{\alpha\xi}
\end{align}
We define $B_{ij|\mu}\in \Gamma\left(U_{ij}, \mathcal{A}^{0,1}\left(\mathscr{H}om_{\mathcal{O}_M}\left( \mathcal{N}_{\mathcal{F}_0}^* , \mathcal{N}_{\mathcal{F}_0}^* \right) \right) \right)$by
\begin{align*}
B_{ij|\mu}:\Gamma\left( U_{ij}, \mathcal{N}_{\mathcal{F}_0}^* \right) &\to \Gamma\left(U_{ij}, \mathcal{A}^{0,1}\left( \mathcal{N}_{\mathcal{F}_0}^* \right) \right) \\
 w_{0i}^\alpha &\mapsto \sum_{\xi=1}^q B_{ij|\mu}^{\alpha \xi} w_{0i}^\xi
\end{align*}
Then from $(\ref{n30})$ and $(\ref{ne11})$, we can find $\left\{ B_{i|\mu} \right\}\in C^0\left(\mathcal{U}, \mathcal{A}^{0,1}\left( \mathscr{H}om_{\mathcal{O}_M}\left( \mathcal{N}_{\mathcal{F}_0}^*, \mathcal{N}_{\mathcal{F}_0}^* \right) \right) \right)$ with $B_{i|\mu}\in \Gamma\left(U_i, \mathcal{A}^{0,1}\left( \mathscr{H}om_{\mathcal{O}_M}\left(\mathcal{N}_{\mathcal{F}_0}^*, \mathcal{N}_{\mathcal{F}_0}^* \right)\right) \right)$ defined by
\begin{align*}
B_{i|\mu} :\Gamma\left( U_i, \mathcal{N}_{\mathcal{F}_0}^* \right) &\to \Gamma\left( U_i, \mathcal{A}^{0,1}\left( \mathcal{N}_{\mathcal{F}_0}^* \right) \right)\\
 w_{0i}^\alpha &\mapsto \sum_{\xi=1}^q B_{i|\mu}^{\alpha \xi} w_{0i}^\xi
\end{align*}
such that
\begin{align}
\sum_{\eta,\gamma=1}^q h_{0ij}^{\alpha \eta} B_{j|\mu}^{\eta \gamma}h_{0ji}^{\gamma \xi} - B_{i|\mu}^{\alpha \xi}=B_{ij|\mu}^{\alpha \xi}-\bar{\partial} D_{ij|\mu}^{\alpha \xi} \label{n50}
\end{align}

We define $C_{ij|\mu}\in \Gamma\left( U_{ij}, \mathcal{A}^{0,0}\left( \mathscr{H}om_{\mathcal{O}_M}\left( \mathcal{N}_{\mathcal{F}_0}^*, \Omega_M^1 \right)\right) \right)$ by
\begin{align*}
C_{ij|\mu}: \Gamma\left( U_{ij}, \mathcal{N}_{\mathcal{F}_0}^* \right) &\to \Gamma\left( U_{ij} , \mathcal{N}_{\mathcal{F}_0}^* \right) \\
   w_{0i}^\alpha  &\mapsto  C_{ij|\mu}^\alpha
\end{align*}
Then from $(\ref{n31})$ and $(\ref{ne11})$, we have
\begin{align*}
\sum_{\gamma=1}^q h_{0ij}^{\alpha\gamma}\left( C_{jk|\mu}^\gamma - \sum_{\beta=1}^q D_{jk|\mu}^{\gamma \beta} w_{0j}^\beta \right)-\left( C_{ik|\mu}^\alpha -\sum_{\beta=1}^q D_{ik|\mu}^{\alpha\beta} w_{0i}^\beta \right)+\left( C_{ij|\mu}^\alpha -\sum_{\beta=1}^q D_{ij|\mu}^{\alpha\beta} w_{0i}^\beta \right)=0
\end{align*}
This implies that we can find $\left\{ C_{i|\mu} \right\}\in C^0\left( \mathcal{U}, \mathcal{A}^{0,0}\left( \mathscr{H}om_{\mathcal{O}_M}\left( \mathcal{N}_{\mathcal{F}_0}^*, \Omega_M^1 \right) \right) \right)$ with $C_{i|\mu}\in \Gamma\left(U_i, \mathcal{A}^{0,0}\left( \mathscr{H}om_{\mathcal{O}_M} \left( \mathcal{N}_{\mathcal{F}_0}^*, \Omega_M^1   \right) \right)  \right)$ defined by
\begin{align*}
C_{i|\mu} : \Gamma\left( U_i, \mathcal{N}_{\mathcal{F}_0}^* \right) &\to \Gamma\left( U_i, \mathcal{A}^{0,0}\left( \mathcal{N}_{\mathcal{F}_0}^* \right) \right) \\
 w_{0i}^\alpha &\mapsto \sum_{\xi=1}^q C_{i|\mu}^{\alpha \xi} w_{0i}^\xi
\end{align*}
such that
\begin{align}
\sum_{\beta=1}^q h_{0ij}^{\alpha\beta} C_{j|\mu}^\beta - C_{i|\mu}^\alpha= C_{ij|\mu}^\alpha - \sum_{\beta=1}^q  D_{ij|\mu}^{\alpha\beta} w_{0i}^\beta \label{n52}
\end{align}
We note that from $(\ref{n417})$ and $(\ref{n52})$, we have
{\Small{\begin{align*}
&E_{i|\mu}^\alpha- \sum_{\beta=1}^q h_{0ij}^{\alpha\beta} \det\left( h_{0ij}\right) E_{j|\mu}^\beta\\
&= \sum_{\beta=1}^q w_{0i}^1\wedge \cdots \wedge C_{ij|\mu}^\beta \wedge \cdots \wedge w_{0i}^q \wedge d w_{0i}^\alpha + w_{0i}^1 \wedge \cdots \wedge w_{0}^q \wedge  \partial C_{ij|\mu}^\alpha\\
&= \sum_{\beta=1}^q w_{0i}^1 \wedge \cdots \wedge  \left( \sum_{\eta=1}^q h_{0ij}^{\beta\eta} C_{j|\mu}^\eta- C_{i|\mu}^\beta+\sum_{\eta=1}^q D_{ij|\mu}^{\beta \eta} w_{0i}^\eta      \right) \wedge \cdots \wedge w_{0i}^q \wedge  d w_{0i}^\alpha + w_{0i}^1 \wedge \cdots \wedge w_{0i}^q \wedge \partial \left( \sum_{\beta=1}^q h_{0ij}^{\alpha\beta} C_{j|\mu}^\beta- C_{i|\mu}^\alpha + \sum_{\beta=1}^q D_{ij|\mu}^{\alpha\beta} w_{0i}^\beta      \right)\\
&= \sum_{\beta=1}^q w_{0i}^1 \wedge \cdots \wedge \left(\sum_{\eta=1}^q h_{0ij}^{\beta \eta} C_{j|\mu}^\eta - C_{i|\mu}^\beta \right) \wedge \cdots \wedge w_{0i}^q \wedge dw_{0i}^\alpha + w_{0i}^1 \wedge \cdots \wedge w_{0i}^q \wedge \partial \left( \sum_{\beta=1}^q h_{0ij}^{\alpha\beta} C_{j|\mu}^\beta - C_{i|\mu}^\alpha \right)\\
&= - \sum_{\beta=1}^q w_{0i}^1 \wedge \cdots \wedge C_{i|\mu}^\beta \wedge \cdots \wedge w_{0i}^q \wedge d w_{0i}^\alpha - w_{0i}^1\wedge \cdots \wedge w_{0i}^q \wedge \partial C_{i|\mu}^\alpha \\
& + \sum_{\beta=1}^q \left( \sum_{\eta_1=1}^q h_{0ij}^{1 \eta_1} w_{0j}^{\eta_1} \right)\wedge \cdots \wedge \left(\sum_{\eta=1}^q h_{0ij}^{\beta \eta} C_{j|\mu}^\eta \right) \wedge \cdots \wedge \left( \sum_{\eta_q=1}^q h_{0ij}^{q \eta_q} w_{0j}^{\eta_q} \right) \wedge \partial \left(\sum_{\gamma=1}^q h_{ij}^{\alpha \gamma} w_j^\gamma \right) +  w_i^1 \wedge \cdots \wedge w_i^q \wedge d\left( \sum_{\beta=1}^q h_{ij}^{\alpha\beta} C_{j|\mu}^\beta \right)\\
&= - \sum_{\beta=1}^q w_{0i}^1 \wedge \cdots \wedge C_{i|\mu}^\beta \wedge \cdots \wedge w_{0i}^q \wedge d w_{0i}^\alpha - w_{0i}^1\wedge \cdots \wedge w_{0i}^q \wedge \partial C_{i|\mu}^\alpha\\
&+\sum_{\beta=1}^q\det\left(h_{0ij}\right) w_{0j}^1\wedge \cdots \wedge C_{j|\mu}^\beta \wedge \cdots \wedge w_{0j}^q \wedge \left( \sum_{\gamma=1}^q h_{0ij}^{\alpha\gamma} dw_i^\gamma+  \sum_{\gamma=1}^q d h_{0ij}^{\alpha\gamma} \wedge w_{0j}^\gamma \right)\\
& + \det\left(h_{0ij}\right)w_{0j}^1\wedge \cdots \wedge w_{0j}^q \wedge \left( \sum_{\beta=1}^q h_{0ij}^{\alpha\beta} \partial C_{j|\mu}^\beta + \sum_{\beta=1}^q dh_{0ij}^{\alpha\beta} \wedge C_{j|\mu}^\beta  \right) \\
&= - \sum_{\beta=1}^q w_{0i}^1 \wedge \cdots \wedge C_{i|\mu}^\beta \wedge \cdots \wedge w_{0i}^q \wedge d w_{0i}^\alpha - w_{0i}^1\wedge \cdots \wedge w_{0i}^q \wedge \partial C_{i|\mu}^\alpha\\
&+ \sum_{\beta,\gamma=1}^q h_{0ij}^{\alpha\beta} \det\left(h_{0ij}\right) w_{0j}^1\wedge \cdots \wedge C_{j|\mu}^\gamma \wedge \cdots \wedge w_{0j}^q \wedge dw_{0j}^\beta + \sum_{\beta=1}^q h_{0ij}^{\alpha\beta} \det\left(h_{0ij}\right) w_{0j}^1 \wedge \cdots w_{0j}^q \wedge \partial C_{j|\mu}^\beta\\
&+ \sum_{\beta=1}^q \det\left(h_{0ij}\right) w_{0j}^1 \wedge \cdots \wedge C_{j|\mu}^\beta \wedge \cdots \wedge w_{0j}^q \wedge dh_{0ij}^{\alpha\beta} \wedge w_{0j}^\beta +\sum_{\beta=1}^q \det\left(h_{0ij}\right) w_{0j}^1 \wedge \cdots \wedge w_{0j}^q \wedge d h_{0ij}^{\alpha\beta} \wedge C_{j|\mu}^\beta\\
&= - \sum_{\beta=1}^q w_{0i}^1 \wedge \cdots \wedge C_{i|\mu}^\beta \wedge \cdots \wedge w_{0i}^q \wedge d w_{0i}^\alpha - w_{0i}^1\wedge \cdots \wedge w_{0i}^q \wedge \partial C_{i|\mu}^\alpha\\
&+\sum_{\beta,\gamma=1}^q h_{0ij}^{\alpha\beta} \det\left( h_{0ij} \right) w_{0j}^1\wedge \cdots \wedge C_{j|\mu}^\gamma \wedge \cdots \wedge w_{0j}^q \wedge dw_{0j}^\beta + \sum_{\beta=1}^q h_{0ij}^{\alpha\beta} \det\left(h_{0ij}\right) w_{0j}^1 \wedge \cdots w_{0j}^q \wedge \partial C_{j|\mu}^\beta,
\end{align*}}}
so that we have
\begin{align*}
&E_{i|\mu}^\alpha +\sum_{\beta=1}^q w_{0i}^1\wedge \cdots \wedge C_{i|\mu}^\beta \wedge w_{0i}^q \wedge dw_{0i}^\alpha + w_{0i}^1 \wedge \cdots \wedge w_{0i}^q \wedge \partial C_{i|\mu}^\alpha\\
&=\sum_{\beta=1}^q h_{0ij}^{\alpha\beta}\det\left(h_{0ij}\right)\left( E_{j|\mu}^\beta + \sum_{\gamma=1}^q w_{0j}^1 \wedge \cdots \wedge C_{j|\mu}^\gamma\wedge \cdots \wedge w_{0j}^q \wedge dw_{0j}^\beta + w_{0j}^1 \wedge \cdots \wedge w_{0j}^q \wedge \partial C_{j|\mu}^\beta \right)
\end{align*}
This implies that from $(\ref{n53})$
\begin{align*}
I_\mu:=\left\{ I_{i|\mu}\right\} \in A^{0,0}\left(M-S, \mathscr{H}om_{\mathcal{O}_M}\left( \mathcal{N}_{\mathcal{F}_0}^*, \tilde{\mathcal{S}}^2 \right) \right)
\end{align*}
where
\begin{align} \label{ne15}
I_{i|\mu}\left(w_{0i}^\alpha \right)= I_{i|\mu}^\alpha:= E_{i|\mu}^\alpha  +\sum_{\beta=1}^q w_{0i}^1\wedge \cdots \wedge C_{i|\mu}^\beta \wedge w_{0i}^q \wedge dw_{0i}^\alpha + w_{0i}^1 \wedge \cdots \wedge w_{0i}^q \wedge \partial C_{i|\mu}^\alpha\,\,\,\,\,\,\textnormal{on}\,\,\,U_i-S
\end{align}
On the other hand, from $(\ref{n28})$ and $(\ref{n52})$ and $(\ref{n50})$, we have
\begin{align*}
&\sum_{\beta=1}^q h_{0ij}^{\alpha\beta} A_{j|\mu}^\beta- A_{i|\mu}^\alpha= \bar{\partial} C_{ij|\mu}^\alpha - \sum_{\eta=1}^q B_{ij|\mu}^{\alpha\eta} w_{0i}^\eta\\
                                  &= \sum_{\beta=1}^q h_{0ij}^{\alpha\beta} \bar{\partial} C_{j|\mu}^\beta -\bar{\partial} C_{i|\mu}^\alpha   +\sum_{\beta=1}^q \bar{\partial} D_{ij|\mu}^{\alpha\beta} w_{0i}^\beta- \sum_{\xi=1}^q \left( \sum_{\eta,\gamma=1}^q h_{0ij}^{\alpha\eta} B_{j|\mu}^{\eta\gamma} h_{0ji}^{\gamma \xi} - B_{i|\mu}^{\alpha \xi}+\bar{\partial} D_{ij|\mu}^{\alpha \xi} \right) w_{0i}^\xi  \\
                                  &=\sum_{\beta=1}^q h_{0ij}^{\alpha\beta}\left( \bar{\partial} C_{j|\mu}^\beta -\sum_{\gamma=1}^q  B_{j|\mu}^{\beta\gamma} w_{0j}^\gamma \right) -\left(  \bar{\partial} C_{i|\mu}^\alpha  -\sum_{\gamma=1}^q B_{i|\mu}^{\alpha \gamma} w_{0i}^\gamma \right)
\end{align*}
Then we have
\begin{align}\label{ne13}
\sum_{\beta=1}^q h_{0ij}^{\alpha\beta}\left( A_{j|\mu}^\beta- \bar{\partial} C_{j|\mu}^\beta +\sum_{\gamma=1}^q B_{j|\mu}^{\beta\gamma} w_{0j}^\gamma \right)-\left( A_{i|\mu}^\alpha - \bar{\partial} C_{i|\mu}^\alpha + \sum_{\gamma=1}^q B_{i|\mu}^{\alpha\gamma} w_{0i}^\gamma \right)=0
\end{align}
We define $\phi_{i|\mu}\in \Gamma\left(U_i, \mathcal{A}^{0,1}\left(\mathscr{H}om_{\mathcal{O}_M}\left( \mathcal{N}_{\mathcal{F}_0}, \Omega_M^1 \right) \right)  \right)$ by
\begin{align}
\phi_{i|\mu}: \Gamma\left( U_i, \mathcal{N}_{\mathcal{F}_0}^* \right) &\to \Gamma\left( U_i, \mathcal{A}^{0,1}\left( \Omega_M^1 \right) \right) \notag \\
  w_{0i}^\alpha &\mapsto \phi_{i|\mu}^\alpha:=  A_{i|\mu}^\alpha- \bar{\partial} C_{i|\mu}^\alpha + \sum_{\gamma=1}^q B_{i|\mu}^{\alpha \gamma} w_{0i}^\gamma \label{ne14}
\end{align}

Then $(\ref{ne13})$ implies that 
\begin{align*}
\phi_\mu:=\left\{ \phi_{i|\mu} \right\} \in A^{0,1}(X, \mathscr{H}om_{\mathcal{O}_X}\left( \mathcal{N}_{\mathcal{F}_0}^*, \Omega_X^1 \right))
\end{align*}
We take $\bar{\phi}_\mu$ to be the image of $\phi_\mu$ in $\frac{A^{0,1}\left( M, \mathscr{H}om_{\mathcal{O}_M}\left( \mathcal{N}_{\mathcal{F}_0}^*, \Omega_M^1 \right) \right)}{A^{0,1}\left( M, \mathscr{H}om_{\mathcal{O}_M}\left( \mathcal{N}_{\mathcal{F}_0}^*, \mathcal{N}_{\mathcal{F}_0}^* \right) \right)}$. Then we claim that
\begin{align}
\left( I_\mu , \bar{\phi}_\mu, - \xi_\mu \right)\in A^{0,0}\left(M-S, \mathscr{H}om_{\mathcal{O}_M}\left( \mathcal{N}_{\mathcal{F}_0}^*, \tilde{\mathcal{S}}^2 \right) \right) \bigoplus \frac{A^{0,1}\left( M, \mathscr{H}om_{\mathcal{O}_M}\left( \mathcal{N}_{\mathcal{F}_0}^*, \Omega_M^1 \right) \right)}{A^{0,1}\left( M, \mathscr{H}om_{\mathcal{O}_M}\left( \mathcal{N}_{\mathcal{F}_0}^*, \mathcal{N}_{\mathcal{F}_0}^*  \right) \right)} \bigoplus A^{0,2}\left( M, \Theta_M \right) \label{nnc13}
\end{align}
defines a $2$-cocycle in the following Dolbeaut type bicomplex associated to $\mathcal{N}_{\mathcal{F}_0}^{*\bullet}$ (see Appendix \ref{app3}).
{\small{\begin{equation}\label{nnc1}
\begin{CD}
\cdots  \\
@A\hat{E}_3AA   \\
A^{0,0}\left( M-S, \left(\mathcal{N}_{\mathcal{F}_0}^*\right)^*\otimes \tilde{\mathcal{S}}^3 \right)@>-(-1)^q \bar{\partial}>> \cdots \\
@A\hat{E}_2AA @A\hat{E}_2AA \\
A^{0,0}\left( M-S, \left(\mathcal{N}_{\mathcal{F}_0}^*\right)^*\otimes \tilde{\mathcal{S}}^2\right) @>(-1)^q\bar{\partial}>> A^{0,1}\left(M-S, \left(\mathcal{N}_{\mathcal{F}_0}^*\right)^*\otimes \tilde{\mathcal{S}}^2  \right) @>-(-1)^q\bar{\partial}>> \cdots  \\
@A\hat{E}_1AA @A \hat{E}_1AA @A \hat{E}_1AA \\
\frac{A^{0,0}\left(M, \left(\mathcal{N}_{\mathcal{F}_0}^* \right)^* \otimes \Omega_M^1 \right)}{A^{0,0}\left(M,\left(\mathcal{N}_{\mathcal{F}_0}^* \right)^*\otimes \mathcal{N}_{\mathcal{F}_0}^*\right)} @>-\bar{\partial}>>\frac{A^{0,1}\left(M, \left(\mathcal{N}_{\mathcal{F}_0}^*\right)^*\otimes \Omega_M^1 \right)}{A^{0,1}\left(M,\left(\mathcal{N}_{\mathcal{F}_0}^* \right)^*\otimes \mathcal{N}_{\mathcal{F}_0}^*\right)}@>\bar{\partial}>> \frac{A^{0,2}\left(M, \left(\mathcal{N}_{\mathcal{F}_0}^*\right)^*\otimes \Omega_M^1 \right)}{A^{0,2}\left(M,\left(\mathcal{N}_{\mathcal{F}_0}^* \right)^*\otimes \mathcal{N}_{\mathcal{F}_0}^*\right)}@>-\bar{\partial}>> \cdots \\
@A\hat{E}_0AA @A\hat{E}_0AA @A\hat{E}_0AA @A\hat{E}_0AA \\
A^{0,0}\left(M, \Theta_M\right) @>\bar{\partial}>> A^{0,1}\left(M, \Theta_M\right) @>-\bar{\partial}>> A^{0,2}\left(M, \Theta_M \right) @>\bar{\partial}>> A^{0,3}\left(M, \Theta_M\right) @>-\bar{\partial}>> \cdots
\end{CD}
\end{equation}}}

From $(\ref{n25})$, $\bar{\partial}\left( -\xi_\mu \right)=0$. We show that $\bar{\partial}\bar{\phi}_\mu+ \hat{E}_0\left( -\xi_\mu\right)=0$. In fact, from $(\ref{ne14})$ and $(\ref{n29})$, we have
\begin{align*}
\bar{\partial} \phi_{i|\mu}^\alpha+ \mathcal{L}_{-\xi_\mu}\left(w_{0i}^\alpha \right) = \bar{\partial}A_{i|\mu}^\alpha + \sum_{\gamma=1}^q \bar{\partial} B_{i|\mu}^{\alpha \gamma} w_{0i}^\gamma - \bar{\partial} A_{i|\mu}^\alpha = \sum_{\gamma=1}^q \bar{\partial } B_{i|\mu}^{\alpha \gamma} w_{0i}^\gamma \in \Gamma\left( U_i, \mathcal{A}^{0,1}\left( \mathcal{N}_{\mathcal{F}_0}^* \right) \right)
\end{align*}
We show that $(-1)^q\bar{\partial} I_\mu+ \hat{E}_1\left(\bar{\phi}_\mu \right)=0$. In fact, from $(\ref{n407})$ and $(\ref{ne14})$ and $(\ref{ne15})$, we have
{\small{\begin{align*}
&(-1)^q\bar{\partial}I_\mu^\alpha+ \hat{E}_1\left(\phi_{i|\mu}\right)\left(w_{0i}^\alpha \right)\\
&=(-1)^q\bar{\partial}\left( E_{i|\mu}^\alpha +\sum_{\beta=1}^q w_{0i}^1 \wedge \cdots \wedge C_{i|\mu}^\beta \wedge \cdots \wedge w_{0i}^q \wedge dw_{0i}^\alpha+ w_{0i}^1\wedge \cdots \wedge w_{0i}^q \wedge \partial C_{i|\mu}^\alpha    \right)\\
&+ \sum_{\beta=1}^q (-1)^{q-\beta+1} w_{0i}^1 \wedge \cdots \wedge \left(A_{i|\mu}^\beta- \bar{\partial} C_{i|\mu}^\beta +\sum_{\gamma=1}^q B_{i|\mu}^{\beta\gamma} w_{0i}^\gamma \right) \wedge \cdots \wedge w_{0i}^q \wedge dw_{0i}^\alpha - w_{0i}^1 \wedge \cdots \wedge w_{0i}^q \wedge \partial \left( A_{i|\mu}^\alpha -\bar{\partial} C_{i|\mu}^\alpha + \sum_{\gamma=1}^q B_{i|\mu}^{\alpha \gamma} w_{0i}^\gamma \right)\\
&= (-1)^q\sum_{\beta=1}^q (-1)^\beta w_{0i}^1 \wedge \cdots \wedge  A_{i|\mu}^\beta \wedge \cdots \wedge w_i^q \wedge dw_{0i}^\alpha + (-1)^{2q} w_{0i}^1 \wedge \cdots \wedge w_i^q \wedge  \partial A_{i|\mu}^\alpha \\
&+ \sum_{\beta=1}^q(-1)^{q+ \beta-1} w_{0i}^1 \wedge \cdots \wedge \bar{\partial} C_{i|\mu}^\beta \wedge \cdots \wedge w_{0i}^q \wedge dw_{0i}^\alpha +(-1)^{2q+1} w_{0i}^1 \wedge \cdots \wedge w_{0i}^q \partial \left( \bar{\partial} C_{i|\mu}^\alpha \right) \\
&+ \sum_{\beta=1}^q (-1)^{q-\beta+1} w_{0i}^1 \wedge \cdots \wedge \left(A_{i|\mu}^\beta- \bar{\partial} C_{i|\mu}^\beta\right) \wedge \cdots \wedge w_{0i}^q \wedge dw_{0i}^\alpha - w_{0i}^1 \wedge \cdots \wedge w_{0i}^q \wedge \partial \left( A_{i|\mu}^\alpha -\bar{\partial} C_{i|\mu}^\alpha \right) = 0
\end{align*}}}

We show that $\hat{E}_2\left(I_\mu \right)=0$. From $(\ref{ne15})$, it is sufficient to show that $\hat{E}_2\left( E_{i|\mu} \right)=0$ on $U_{i_x}$. First we note that
\begin{align}\label{ne16}
\partial \left(E_{i|\mu}^\alpha\right)&\equiv_\mu \partial \left( w_i^{1(\mu-1)}\wedge \cdots \wedge w_i^{q(\mu-1)}\wedge \partial w_i^{\alpha(\mu-1)} \right)\\
 &\equiv_\mu \sum_{\beta=1}^q (-1)^{\beta-1} w_i^{1(\mu-1)}\wedge \cdots \wedge \overbrace{ \partial w_i^{\beta(\mu-1)}}^{\beta-\textnormal{th}} \wedge \cdots \wedge w_i^{q(\mu-1)}\wedge \partial w_i^{\alpha(\mu-1)} \notag \\
 &\equiv_\mu \sum_{\beta=1}^q(-1)^{\beta-1}  w_i^{1(\mu-1)}\wedge \cdots \wedge \partial w_i^{\beta(\mu-1)} \wedge \cdots \wedge w_i^{q(\mu-1)}\wedge \left(\sum_{\gamma=1}^q a_{i_x}^{\alpha\gamma(\mu-1)} \wedge w_i^{\gamma(\mu-1)} + L_{i_x|\mu}^\alpha\right) \notag \\
  &\equiv_\mu \sum_{\beta=1}^q a_{i_x}^{\alpha\beta(\mu-1)}\wedge w_i^{1(\mu-1)}\wedge  \cdots \wedge w_i^{q(\mu-1)} \wedge \partial  w_i^{\beta(\mu-1)}+ d\left(w_{0i}^1\wedge \cdots \wedge w_{0i}^q \right) \wedge L_{i_x|\mu}^\alpha \notag \\
  &\equiv_\mu \sum_{\beta=1}^q a_{0i_x}^{\alpha\beta} \wedge E_{i|\mu}^\beta + d\left(w_{0i}^1 \wedge \cdots \wedge w_{0i}^q \right)\wedge L_{i_x|\mu}^\alpha  \notag
\end{align}

From $(\ref{n53})$, $E_{i|\mu}^\alpha= w_{0i}^1 \wedge \cdots \wedge w_{0i}^q \wedge L_{i_x|\mu}^\alpha$ on $U_{i_x}$. Then from $(\ref{ne16})$, we have
\begin{align*}
&\hat{E}_2\left(E_{i|\mu} \right)\left(w_{0i}^\alpha \right)=w_{0i}^1 \wedge \cdots \wedge w_{0i}^q \wedge \left( \partial L_{i_x|\mu}^\alpha - \sum_{\beta=1}^q a_{0i_x}^{\alpha\beta} \wedge L_{i_x|\mu}^\beta \right) \\
& = (-1)^q \partial \left( w_{0i}^1 \wedge \cdots \wedge w_{0i}^q \wedge  L_{i_x|\mu}^\alpha \right) -(-1)^q d\left( w_{0i}^1 \wedge \cdots \wedge w_{0i}^q \right) \wedge L_{i_x|\mu}^\alpha - (-1)^q \sum_{\beta=1}^q a_{0i_x}^{\alpha\beta} \wedge w_{0i}^1 \wedge \cdots \wedge w_{0i}^q \wedge L_{i_x|\mu}^\beta \\
&= (-1)^q \partial E_{i|\mu}^\alpha - (-1)^q \left(\partial E_{i|\mu}^\alpha - \sum_{\beta=1}^q a_{0i_x}^{\alpha\beta} \wedge E_{i|\mu}^\beta \right) - (-1)^q a_{0i_x}^{\alpha\beta} \wedge E_{i|\mu}^\beta = 0
\end{align*}
This proves $\left( I_\mu, \bar{\phi}_\mu , - \xi_\mu \right)$ defines a $2$-cocycle in the Dolbeault type bicomplex associated to $\mathcal{N}_{\mathcal{F}_0}^{*\bullet}$. Then by hypothesis $\mathbb{H}^2\left( M, \mathcal{N}_{\mathcal{F}_0}^{*\bullet} \right)=0$, there exists
\begin{align*}
\left(\overline{w'_\mu}, \varphi'_\mu \right) \in \frac{A^{0,0}\left( M, \mathscr{H}om_{\mathcal{O}_M}\left( \mathcal{N}_{\mathcal{F}_0}^* ,  \Omega_M^1 \right) \right)}{A^{0,0}\left( M,  \mathscr{H}om_{\mathcal{O}_M}\left(  \mathcal{N}_{\mathcal{F}_0}^* ,  \mathcal{N}_{\mathcal{F}_0}^* \right) \right)} \bigoplus A^{0,1}\left(M, \Theta_M \right)
\end{align*}
such that
\begin{align}
-\bar{\partial}\varphi_\mu'&= - \xi_\mu \label{n54}\\
\overline{-\bar{\partial} w_\mu'}+ \hat{E}_0 \left(\varphi_\mu' \right)&= \bar{\phi}_\mu=\overline{ \left\{ w_{0i}^\alpha \mapsto A_{i|\mu}^\alpha - \bar{\partial} C_{i|\mu}^\alpha + \sum_{\gamma=1}^q B_{i|\mu}^{\alpha\gamma} w_{0i}^\gamma    \right\} }  \label{n55}\\
\hat{E}_1 \left(w_\mu' \right)= \psi_\mu=& \left\{ w_{0i}^\alpha \mapsto E_{i|\mu}^\alpha+ \sum_{\beta=1}^q w_{0i}^1\wedge \cdots \wedge  \overbrace{C_{i|\mu}^\beta}^{\beta-\textnormal{th}} \wedge \cdots \wedge w_{0i}^q\wedge dw_{0i}^\alpha + w_{0i}^1\wedge \cdots \wedge w_{0i}^q \wedge \partial C_{i|\mu}^\alpha     \right\} \label{n56}
\end{align}
where $\overline{w_\mu'}$ is the image of the global section $w_\mu'=\left\{w_{i|\mu}' \right\}\in A^{0,0}\left( M, \mathscr{H}om_{\mathcal{O}_M}\left(  \mathcal{N}_{\mathcal{F}_0}^*, \Omega_M^1   \right) \right)$ in $\frac{A^{0,0}\left( M, \mathscr{H}om_{\mathcal{O}_M}\left( \mathcal{N}_{\mathcal{F}_0}^* ,  \Omega_M^1 \right) \right)}{A^{0,0}\left( M,  \mathscr{H}om_{\mathcal{O}_M}\left(  \mathcal{N}_{\mathcal{F}_0}^* ,  \mathcal{N}_{\mathcal{F}_0}^* \right) \right)}$ such that $w_{i|\mu}' : \Gamma\left( U_i, \mathcal{N}_{\mathcal{F}_0}^*\right) \to \Gamma\left( U_i, \mathcal{A}^{0,0}\left( \Omega_M^1 \right)\right), w_{0i}^\alpha \mapsto w_{i|\mu}'^\alpha,\alpha=1,...,q$. In particular, we have
\begin{align}\label{ne17}
w_{i|\mu}'^\alpha= \sum_{\beta=1}^q h_{0ij}^{\alpha\beta} w_{j|\mu}'^\beta
\end{align}
From $(\ref{n54})$, if we take
\begin{align}\label{ne18}
\varphi_\mu: = - \varphi_\mu',
\end{align}
the $(\ref{n42})_\mu$ is satisfied. From $(\ref{n55})$, there exist $\left\{ V_{i|\mu} \right\}\in C^0\left( \mathcal{U}, \mathcal{A}^{0,1}\left(\mathscr{H}om_{\mathcal{O}_M}\left( \mathcal{N}_{\mathcal{F}_0}^*, \mathcal{N}_{\mathcal{F}_0}^* \right) \right) \right)$ defined by $V_{i|\mu}:\Gamma\left( U_i, \mathcal{N}_{\mathcal{F}_0}^* \right) \to \Gamma\left( U_i, \mathcal{A}^{0,1}\left( \mathcal{N}_{\mathcal{F}_0}^* \right) \right), w_{0i}^\alpha \mapsto \sum_{\xi=1}^q V_{i|\mu}^{\alpha \xi} w_{0i}^\xi$ such that
\begin{align}\label{ne23}
\sum_{\beta,\gamma=1}^q h_{0ij}^{\alpha\beta} V_{j|\mu}^{\beta\gamma} w_{0j}^\gamma - \sum_{\beta=1}^q V_{i|\mu}^{\alpha\beta} w_{0i}^\beta= \sum_{\beta=1}^q \left[ \varphi_\mu', h_{0ij}^{\alpha\beta} \right] w_{0j}^\beta
\end{align}
and
\begin{align}
-\bar{\partial} w_{i|\mu}'^\alpha+ \mathcal{L}_{\varphi_\mu'} \left(w_{0i}^\alpha \right)+ \sum_{\beta=1}^q V_{i|\mu}^{\alpha\beta} w_{0i}^\beta= A_{i | \mu}^\alpha- \bar{\partial} C_{i|\mu}^\alpha + \sum_{\gamma=1}^q B_{i|\mu}^{\alpha \gamma} w_{0i}^\gamma \label{n51}
\end{align}
Since $\bar{\partial} A_{i|\mu}^\alpha= \mathcal{L}_{\xi_\mu} \left(w_{0i}^\alpha \right)$ from $(\ref{n43})$,  by taking $\bar{\partial}$ on $(\ref{n51})$ we have
\begin{align*}
\mathcal{L}_{\xi_\mu} \left(w_{0i}^\alpha\right) +\sum_{\beta=1}^q \bar{\partial} V_{i|\mu}^{\alpha\beta}  w_{0i}^\beta=\bar{\partial} A_{i|\mu}^\alpha + \sum_{\gamma=1}^q \bar{\partial} B_{i|\mu}^{\alpha\gamma} w_{0i}^\gamma \iff \sum_{\beta=1}^q \bar{\partial}\left( V_{i|\mu}^{\alpha\beta}- B_{i|\mu}^{\alpha\beta} \right)  w_{0i}^\beta =0 \iff \bar{\partial}\left(V_{i|\mu}^{\alpha\beta}-B_{i|\mu}^{\alpha\beta} \right)=0
\end{align*}
Hence there exists $P_{i|\mu}^{\alpha\beta}\in \Gamma\left(U_i, \mathscr{A}^{0,0}\right)$ such that 
\begin{align}\label{ne21}
\bar{\partial} P_{i|\mu}^{\alpha\beta}=V_{i|\mu}^{\alpha\beta}-B_{i|\mu}^{\alpha\beta}.
\end{align}
 Then if we take
\begin{align}\label{ne19}
w_{i|\mu}^\alpha:= C_{i|\mu}^\alpha- w_{i|\mu}'^\alpha+ \sum_{\beta=1}^q P_{i|\mu}^{\alpha\beta} w_{0i}^\beta,
\end{align}
then we have from $(\ref{ne18})$ and $(\ref{n51})$
\begin{align*}
\bar{\partial}\left(C_{i|\mu}^\alpha- w_{i|\mu}'^{\alpha}  +  \sum_{\beta=1}^q P_{i|\mu}^{\alpha\beta} w_i^\beta    \right) -\mathcal{L}_{-\varphi_\mu'}\left(w_i^\alpha \right)= A_{i|\mu}^\alpha \Longrightarrow \bar{\partial} w_{i|\mu}^\alpha - \mathcal{L}_{\varphi_\mu}\left( w_{0i}^\alpha \right)= A_{i|\mu}^\alpha,
\end{align*}
so that $(\ref{n43})_\mu$ is satisfied. We find $h_{ij|\mu}^{\alpha\beta}$ satisfying $(\ref{n46})_\mu$. We note that from $(\ref{ne19})$ and $(\ref{n52})$ and $(\ref{ne17})$, 
\begin{align}\label{ne20}
&-\sum_{\beta=1}^q h_{0ij}^{\alpha\beta} w_{j|\mu}^\beta + w_{i|\mu}^\alpha + C_{ij|\mu}^\alpha\\
&= -  \sum_{\beta=1}^q h_{0ij}^{\alpha\beta}\left( C_{j|\mu}^\beta - w_{j|\mu}'^\beta + \sum_{\xi=1}^q P_{j|\mu}^{\beta \xi} w_{0j}^\xi   \right) + \left( C_{i|\mu}^\alpha - w_{i|\mu}'^\alpha + \sum_{\xi=1}^q P_{i|\mu}^{\alpha \xi} w_{0i}^\xi \right) + \sum_{\beta=1}^q h_{0ij}^{\alpha\beta} C_{j|\mu}^\beta - C_{i|\mu}^\alpha + \sum_{\gamma=1}^q D_{ij|\mu}^{\alpha \gamma} w_{0i}^\gamma \notag \\
&=\sum_{\beta=1}^p \left( \sum_{\gamma=1}^q D_{ij|\mu}^{\alpha\gamma} h_{0ij}^{\gamma \beta} + \sum_{\gamma=1}^q P_{i|\mu}^{\alpha\gamma} h_{0ij}^{\gamma \beta} - \sum_{\gamma=1}^q h_{0ij}^{\alpha\gamma} P_{j|\mu}^{\gamma\beta} \right) w_{0j}^\beta \notag
\end{align}
If we take 
\begin{align}\label{ne22}
h_{ij|\mu}^{\alpha \beta} : = \sum_{\gamma=1}^q D_{ij|\mu}^{\alpha\gamma} h_{0ij}^{\gamma \beta} + \sum_{\gamma=1}^q P_{i|\mu}^{\alpha\gamma} h_{0ij}^{\gamma \beta} - \sum_{\gamma=1}^q h_{0ij}^{\alpha\gamma} P_{j|\mu}^{\gamma\beta},
\end{align}
then $(\ref{ne20})$ implies $(\ref{n46})_\mu$. We check $(\ref{n45})_\mu$. In fact, from
$(\ref{ne22}), (\ref{n50}), (\ref{ne21})$ and $(\ref{ne23})$,
{\small{\begin{align*} 
 &\bar{\partial} h_{ij|\mu}^{\alpha\beta}= \bar{\partial}\left(\sum_{\gamma=1}^q D_{ij|\mu}^{\alpha\gamma} h_{0ij}^{\gamma \beta}  + \sum_{\gamma=1}^q P_{i|\mu}^{\alpha\gamma} h_{0ij}^{\gamma\beta} - \sum_{\gamma=1}^q h_{0ij}^{\alpha\gamma} P_{j|\mu}^{\gamma\beta}  \right) =\sum_{\gamma=1}^q \left( \bar{\partial} D_{ij|\mu}^{\alpha\gamma}\right) h_{0ij}^{\gamma \beta}  + \sum_{\gamma=1}^q \left( \bar{\partial} P_{i|\mu}^{\alpha\gamma} \right) h_{0ij}^{\gamma\beta} - \sum_{\gamma=1}^q h_{0ij}^{\alpha\gamma} \bar{\partial} P_{j|\mu}^{\gamma\beta}  \\
 &= \sum_{\gamma=1}^q \left( B_{ij|\mu}^{\alpha\gamma}- \sum_{\eta,\xi=1}^q h_{0ij}^{\alpha \eta} B_{j|\mu}^{\eta \xi}h_{0ji}^{\xi \gamma} + B_{i|\mu}^{\alpha \gamma} \right) h_{0ij}^{\gamma \beta} + \sum_{\gamma=1}^q \left( V_{i|\mu}^{\alpha\gamma}- B_{i|\mu}^{\alpha\gamma} \right) h_{0ij}^{\gamma\beta} - \sum_{\gamma=1}^q h_{0ij}^{\alpha\gamma} \left( V_{j|\mu}^{\gamma\beta} - B_{j|\mu}^{\gamma\beta} \right)\\
 &= \sum_{\gamma=1}^q B_{ij|\mu}^{\alpha\gamma} h_{0ij}^{\gamma\beta} + \sum_{\gamma=1}^q V_{i|\mu}^{\alpha\gamma} h_{0ij}^{\gamma\beta} - \sum_{\gamma=1}^q h_{0ij}^{\alpha\gamma} V_{j|\mu}^{\gamma\beta} = \sum_{\gamma=1}^q B_{ij|\mu}^{\alpha\gamma} h_{0ij}^{\gamma\beta} - \left[ \varphi_\mu', h_{0ij}^{\alpha\beta} \right] =\sum_{\gamma=1}^q B_{ij|\mu}^{\alpha\gamma} h_{0ij}^{\gamma\beta} + \left[\varphi_\mu, h_{0ij}^{\alpha\beta} \right]
\end{align*}}}
We check $(\ref{nn121})_\mu$. In fact, from $(\ref{ne11})$ and $(\ref{ne22})$, we have
{\small{\begin{align*}
&\sum_{\beta,\gamma=1}^q h_{0ij}^{\alpha\beta} h_{jk|\mu}^{\beta\gamma} h_{0ki}^{\gamma\xi}-\sum_{\gamma=1}^q h_{ik|\mu}^{\alpha\gamma} h_{0ki}^{\gamma\xi}+\sum_{\beta=1}^q h_{ij|\mu}^{\alpha\beta} h_{0ji}^{\beta\xi}\\
&=\sum_{\beta, \gamma=1}^q h_{0ij}^{\alpha\beta}\left( \sum_{\eta=1}^q D_{jk|\mu}^{\beta\eta} h_{0jk}^{\eta \gamma} +\sum_{\eta=1}^q P_{j|\mu}^{\beta \eta} h_{0jk}^{\eta\gamma} - \sum_{\eta=1}^q h_{0jk}^{\beta\eta} P_{k|\mu}^{\eta \gamma} \right) h_{0ki}^{\gamma \xi} - \sum_{\gamma=1}^p\left( \sum_{\eta=1}^q D_{ik|\mu}^{\alpha \eta} h_{0ik}^{\eta \gamma} + \sum_{\eta=1}^q P_{i|\mu}^{\alpha \eta} h_{0ik}^{\eta\gamma} - \sum_{\eta=1}^q  h_{0ij}^{\alpha \eta} P_{k|\mu}^{\eta \gamma} \right) h_{0ki}^{\gamma \xi} \\
& + \sum_{\beta=1}^q \left( \sum_{\eta=1}^q D_{ij|\mu}^{\alpha \beta} h_{0ij}^{\eta\beta} +\sum_{\eta=1}^q P_{i|\mu}^{\alpha \eta} h_{0ij}^{\eta \beta} - \sum_{\eta=1}^q h_{0ij}^{\alpha \eta} P_{j|\mu}^{\eta\beta}   \right) h_{0ji}^{\beta\xi}  = D_{ijk|\mu}^{\alpha \xi}
\end{align*}}}

On the other hand, we find $a_{i_x|\mu}^{\alpha\beta}$ satisfying $(\ref{n48})_\mu$:
\begin{align*}
L_{i_x|\mu}^\alpha &= - \partial w_{i|\mu}^\alpha + \sum_{\beta=1}^q a_{0 i_x}^{\alpha\beta}\wedge w_{i|\mu}^\beta + \sum_{\beta=1}^q a_{i_x|\mu}^{\alpha\beta} \wedge w_{0i}^\beta 
\end{align*}
In fact, from $(\ref{n53})$ and $(\ref{ne19})$ and $(\ref{n56})$,
\begin{align*}
&w_{0i}^1\wedge \cdots \wedge w_{0i}^q \wedge \left(L_{i_x|\mu}^\alpha + \partial w_{i|\mu}^\alpha -  \sum_{\beta=1}^q a_{0 i_x}^{\alpha\beta} \wedge w_{i|\mu}^\beta\right)\\
&=w_{0i}^1\wedge \cdots \wedge w_{0i}^q \wedge \left( L_{i_x|\mu}^\alpha + \partial \left(C_{i|\mu}^\alpha- w_{i|\mu}'^\alpha+ \sum_{\beta=1}^q P_{i|\mu}^{\alpha\beta} w_{0i}^\beta \right) -\sum_{\beta=1}^q a_{0 i_x}^{\alpha\beta} \wedge \left( C_{i|\mu}^\beta- w_{i|\mu}'^\beta + \sum_{\gamma=1}^q P_{j|\mu}^{\beta\gamma} w_{0i}^\gamma  \right) \right) \\
&= w_{0i}^1 \wedge \cdots \wedge w_{0i}^q \wedge L_{i_x|\mu}^\alpha + w_{0i}^1 \wedge \cdots \wedge w_{0i}^q \wedge \partial C_{i|\mu}^\alpha - w_{0i}^1 \wedge \cdots \wedge w_{0i}^q \wedge \partial w_{i|\mu}'^\alpha - w_{0i}^1 \wedge \cdots \wedge w_{0i}^q \wedge \left(\sum_{\beta=1}^q a_{0i_x}^{\alpha\beta}\wedge C_{i|\mu}^\beta \right)\\
&+ w_{0i}^1\wedge \cdots \wedge w_{0i}^q \wedge \left( \sum_{\beta=1}^q a_{0i_x}^{\alpha\beta} \wedge w_{i|\mu}'^\beta   \right) = w_i^1\wedge \cdots \wedge w_i^q \wedge L_{i_x|\mu}^\alpha - E_{i|\mu}^\alpha =0
\end{align*}
Since $w_{0i}^1,...,w_{0i}^q$ defines a regular foliation on $U_{i_x} \subset U_i-S$, we can find $a_{i_x|\mu}^{\alpha\beta}$ satisfying $(\ref{n48})_\mu$. It remains to check $(\ref{n47})_\mu$. In fact, from $(\ref{n56})$ and $(\ref{ne19})$,
\begin{align}\label{ne26}
&\sum_{\beta=1}^q w_{0i}^1\wedge \cdots \wedge  \overbrace{ w_{i|\mu}^\beta}^{\beta-\textnormal{th}} \wedge \cdots \wedge w_{0i}^q \wedge d w_{0i}^\alpha + w_{0i}^1\wedge \cdots \wedge w_{0i}^q \wedge \partial w_{i|\mu}^\alpha\\
&=\sum_{\beta=1}^q w_{0i}^1\wedge \cdots \wedge  \overbrace{ \left( C_{i|\mu}^\beta -w_{i|\mu}'^\beta+ \sum_{\gamma=1}^q P_{i|\mu}^{\beta\gamma} w_{0i}^\gamma \right) }^{\beta-\textnormal{th}} \wedge \cdots \wedge w_{0i}^q \wedge d w_{0i}^\alpha + w_{0i}^1\wedge \cdots \wedge w_{0i}^q \wedge \partial \left( C_{i|\mu}^\alpha - w_{i|\mu}'^\alpha + \sum_{\beta=1}^q P_{i|\mu}^{\alpha\beta} w_{0i}^\beta \right) \notag\\
&= -E_{i|\mu}^\alpha \notag
\end{align}
By the definition of $\hat{E}_2$, $(\ref{ne26})$ implies that 
\begin{align} \label{ne27}
E_{i|\mu}^\alpha +  \sum_{\beta=1}^q w_{0i}^1\wedge \cdots \wedge  \overbrace{ w_{i|\mu}^\beta}^{\beta-\textnormal{th}} \wedge \cdots \wedge w_{0i}^q \wedge d w_{0i}^\alpha + w_{0i}^1\wedge \cdots \wedge w_{0i}^q \wedge \partial w_{i|\mu}^\alpha=0 \,\,\,\,\,\,\,\textnormal{on}\,\,\, U_i- S
\end{align}
We note that from $(\ref{n43})$ and $(\ref{n40})$
\begin{align*}
&\bar{\partial} \left( E_{i|\mu}^\alpha +  \sum_{\beta=1}^q w_{0i}^1\wedge \cdots \wedge  \overbrace{ w_{i|\mu}^\beta}^{\beta-\textnormal{th}} \wedge \cdots \wedge w_{0i}^q \wedge d w_{0i}^\alpha + w_{0i}^1\wedge \cdots \wedge w_{0i}^q \wedge \partial w_{i|\mu}^\alpha \right)\\
&= \sum_{\beta=1}^q (-1)^\beta w_{0i}^1\wedge \cdots \wedge  \overbrace{ A_{i|\mu}^\beta}^{\beta- \textnormal{th}} \wedge \cdots \wedge w_{0i}^q \wedge dw_{0i}^\alpha + (-1)^q w_{0i}^1 \wedge \cdots \wedge w_{0i}^q \wedge \partial A_{i|\mu}^\alpha \\
&+ \sum_{\beta=1}^q (-1)^{\beta-1}w_{0i}^1 \wedge \cdots \wedge \left( A_{i|\mu}^\beta + \mathcal{L}_{\varphi_\mu}\left(w_{0i}^\beta \right) \right) \wedge \cdots \wedge w_{0i}^q \wedge dw_{0i}^\alpha + (-1)^{q+1} w_{0i}^1 \wedge \cdots \wedge w_{0i}^q \wedge \partial \left( A_{i|\mu}^\alpha + \mathcal{L}_{\varphi_\mu}\left( w_{0i}^\alpha \right) \right)\\
&= \sum_{\beta=1}^q (-1)^{\beta-1} w_{0i}^1 \wedge \cdots \wedge \mathcal{L}_{\varphi}\left(w_{0i}^\beta \right) \wedge \cdots \wedge w_{0i}^q \wedge dw_{0i}^\alpha + (-1)^{q+1} w_{0i}^1 \wedge \cdots \wedge w_{0i}^q \wedge \partial\left(\mathcal{L}_{\varphi_\mu}\left(w_{0i}^\alpha \right) \right)\\
&= \mathcal{L}_{\varphi_\mu} \left( w_{0i}^1 \wedge \cdots \wedge w_{0i}^q \wedge dw_{0i}^\alpha \right) = 0
\end{align*}
This implies that $E_{i|\mu}^\alpha +  \sum_{\beta=1}^q w_{0i}^1\wedge \cdots \wedge  \overbrace{ w_{i|\mu}^\beta}^{\beta-\textnormal{th}} \wedge \cdots \wedge w_{0i}^q \wedge d w_{0i}^\alpha + w_{0i}^1\wedge \cdots \wedge w_{0i}^q \wedge \partial w_{i|\mu}^\alpha$ is holomorphic on $U_i$. Hence $(\ref{ne27})$ implies that $E_{i|\mu}^\alpha +  \sum_{\beta=1}^q w_{0i}^1\wedge \cdots \wedge  \overbrace{ w_{i|\mu}^\beta}^{\beta-\textnormal{th}} \wedge \cdots \wedge w_{0i}^q \wedge d w_{0i}^\alpha + w_{0i}^1\wedge \cdots \wedge w_{0i}^q \wedge \partial w_{i|\mu}^\alpha=0$ on $U_i$, so that $(\ref{n47})_\mu$ is satisfied. This completes Lemma \ref{ne25}.

\end{proof}

It remains to determine $\varphi_1, w_{i|1}^\alpha, h_{ij|1}^{\alpha\beta}$ and $a_{i_x|1}^{\alpha\beta}$ satisfying $(\ref{n23})_1- (\ref{n241})_1$. Given $\dim_\mathbb{C} \mathbb{H}^1\left(  M, \mathcal{N}_{\mathcal{F}_0}^{*\bullet} \right)=r$, we can find a basis of $\mathbb{H}^1\left( M, \mathcal{N}_{\mathcal{F}_0}^{*\bullet} \right)$ by using the Dolbeault type bicomplex associated to $\mathcal{N}_{\mathcal{F}_0}^{*\bullet}$ from $(\ref{nnc1})$ and represent the basis by
\begin{align}
\left(  \overline{\kappa_\lambda} , \rho_\lambda \right) \in \frac{A^{0,0}\left( M, \mathscr{H}om_{\mathcal{O}_M}\left( \mathcal{N}_{\mathcal{F}_0}^*, \Omega_M^1 \right) \right)}{A^{0,0}\left(  M, \mathscr{H}om_{\mathcal{O}_M}\left( \mathcal{N}_{\mathcal{F}_0}^*, \mathcal{N}_{\mathcal{F}_0}^* \right)\right)} \bigoplus A^{0,1}\left( M, \Theta_M \right),\,\,\,\,\,\,\lambda=1,...,r \label{st7}
\end{align}
where $\kappa_\lambda \in A^{0,0}\left(M, \mathscr{H}om_{\mathcal{O}_M}\left( \mathcal{N}_{\mathcal{F}_0}^*, \Omega_M^1 \right) \right)$ defined by $\kappa_\lambda \left( w_{0i}^\alpha \right)= \kappa_{i\lambda}^\alpha$ on $U_i$ for $\kappa_{i\lambda}^\alpha\in \Gamma\left( U_i, \mathcal{A}^{0,0}\left( \Omega_M^1 \right) \right)$ with $\kappa_{i\lambda}^\alpha= \sum_{\beta=1}^q h_{0ij}^{\alpha\beta} \kappa_{j\lambda}^\beta$ on $U_{ij}$ such that there exists $\theta_{i\lambda}\in \Gamma\left( U_i, \mathcal{A}^{0,1}\left( \mathscr{H}om_{\mathcal{O}_M}\left( \mathcal{N}_{\mathcal{F}_0}^*, \mathcal{N}_{\mathcal{F}_0}^* \right)\right) \right)$ defined by $\theta_{i\lambda}\left( w_{0i}^\alpha \right)= \theta_{i\lambda}^\alpha= \sum_{\beta=1}^q \theta_{i\lambda}^{\alpha\beta} w_{0i}^\beta \in \Gamma\left( U_i, \mathcal{A}^{0,1}\left(\mathcal{N}_{\mathcal{F}_0}^* \right)\right)$ satisfying
\begin{align}
-\bar{\partial}\rho_\lambda & =0 \\
\sum_{\beta=1}^q w_{0i}^1\wedge \cdots \wedge \kappa_{i\lambda}^\alpha \wedge \cdots \wedge w_{0i}^q & + w_{0i}^1 \wedge \cdots \wedge w_{0i}^q \wedge \partial \kappa_{i\lambda}^\alpha =0 \label{nnc6}  \\
-\bar{\partial}\kappa_{i\lambda}^\alpha +\mathcal{L}_{\rho_\lambda}\left(w_{0i}^\alpha \right) &= \sum_{\beta=1}^q \theta_{i\lambda}^{\alpha\beta} w_{0i}^\beta \label{nnc3}
\end{align}

We note that since $\bar{\partial} \theta_{i\lambda}^{\alpha\beta}=0$ from $(\ref{nnc3})$, there exists $\mathfrak{K}_{i\lambda}^{\alpha \beta} \in \Gamma\left(U_i, \mathcal{A}^{0,0} \right)$ such that $\bar{\partial} \mathfrak{K}_{i\lambda}^{\alpha\beta}= \theta_{i\lambda}^{\alpha\beta}$. We set
\begin{align}
\varphi_1&= \sum_{\lambda=1}^r t_\lambda \rho_\lambda \label{ns3}\\
w_i^{\alpha 1}&= w_{0i}^\alpha + \sum_{\lambda=1}^r t_\lambda \left( \kappa_{i\lambda}^\alpha + \mathfrak{K}_{i\lambda}^\alpha \right) := w_{0i}^\alpha + \sum_{\lambda=1}^r t_\lambda \left( \kappa_{i\lambda}^\alpha + \sum_{\beta=1}^q  \mathfrak{K}_{i\lambda}^{\alpha\beta} w_{0i}^\beta \right) \label{st8}\\
h_{ij}^{\alpha\beta 1} & = h_{0ij}^{\alpha\beta} + \sum_{\lambda=1}^r t_\lambda \left( \sum_{\gamma=1}^q \mathfrak{K}_{i\lambda}^{\alpha \gamma} h_{0ij}^{\gamma \beta} - \sum_{\gamma=1}^q h_{0ij}^{\alpha \gamma} \mathfrak{K}_{j\lambda}^{\gamma \beta}\right)  \label{st9}
\end{align}
Then $(\ref{n23})_1,(\ref{nnc4})_1,(\ref{n213})_1, (\ref{nnc5})_1$ and $(\ref{n24})_1$ are satisfied. $(\ref{nnc4})_1$ and $(\ref{n213})_1$ implies $(\ref{nnc8})_1$.

It remains to find $a_{i_x}^{\alpha\beta 1}$. In fact, from $(\ref{nnc6})$, we have on $U_{i_x}$
\begin{align}
\partial \kappa_{i\lambda}^\alpha - \sum_{\beta=1}^q a_{0i_x}^{\alpha\beta} \wedge \kappa_{i\lambda}^\beta = \sum_{\beta=1}^q b_{i_x \lambda }^{\alpha\beta} \wedge w_{0i}^\beta \,\,\,\,\,\,\,\,\,\textnormal{for some}\,\,\, b_{i_x\lambda}^{\alpha\beta}\in \Gamma\left( U_i, \mathcal{A}^{0,0}\left(\Omega_M^1 \right) \right) \label{st10}
\end{align}
Then we have
\begin{align*}
\partial w_i^{\alpha 1} &= d w_{0i}^\alpha + \sum_{\lambda=1}^r t_\lambda \left( \partial  \kappa_{i\lambda}^\alpha  + \sum_{\beta=1}^q \partial \mathfrak{K}_{i\lambda}^{\alpha\beta} \wedge w_{0i}^\beta + \sum_{\beta=1}^q \mathfrak{K}_{i\lambda}^{\alpha\beta} d w_{0i}^\beta \right) \\
&\equiv_1 \sum_{\beta=1}^q \left( a_{0i_x}^{\alpha\beta} + \sum_{\lambda=1}^r t_\lambda \left( b_{i_x\lambda }^{\alpha\beta}  + \partial \mathfrak{K}_{i\lambda}^{\alpha\beta} - \sum_{\eta=1}^q a_{0i_x}^{\alpha \eta} \mathfrak{K}_{i\lambda}^{\eta \beta}   + \sum_{\eta=1}^q \mathfrak{K}_{i\lambda}^{\alpha\eta} a_{0i_x}^{\eta \beta}  \right) \right) \wedge \left( w_{0i}^\beta + \sum_{\lambda=1}^r t_\lambda \left( \kappa_{i\lambda}^\beta + \sum_{\gamma=1}^q \mathfrak{K}_{i\lambda}^{\beta \gamma} w_{0i}^\gamma \right) \right),
\end{align*}
so that if we take
\begin{align*}
a_{i_x}^{\alpha\beta 1} : = a_{0i_x}^{\alpha\beta} + \sum_{\lambda=1}^r t_\lambda \left( b_{i_x\lambda }^{\alpha\beta}  + \partial \mathfrak{K}_{i\lambda}^{\alpha\beta} - \sum_{\eta=1}^q a_{0i_x}^{\alpha \eta} \mathfrak{K}_{i\lambda}^{\eta \beta}   + \sum_{\eta=1}^q \mathfrak{K}_{i\lambda}^{\alpha\eta} a_{0i_x}^{\eta \beta}  \right),
\end{align*}
then $(\ref{n241})_1$ is satisfied. This completes the inductive construction of $\varphi, w_i^\alpha, h_{ij}^{\alpha\beta}$ and $a_{i_x}^{\alpha\beta}$ satisfying $(\ref{n22})-(\ref{n211})$.

\subsection{Proof of convergence}\

We will prove that $\varphi=\sum_{\mu=1}^\infty \varphi_\mu, w_i^\alpha= \sum_{\mu=0}^\infty w_{i|\mu}^\alpha$, and $h_{ij}^{\alpha\beta}= \sum_{\mu=0}^\infty h_{ij|\mu}^{\alpha\beta}$ converge. We note that once we show $w_i^\alpha$ converges on $U_i\times B$ for some neighborhood $B$ of $0$, then the condition $w_i^1(z_i,t)\wedge \cdots \wedge w_i^q(w_i,t)\wedge \wedge \partial w_i^\alpha(z_i,t)=0$ implies that since $w_{0i}^1,..., w_{0i}^q$ defines a regular foliation on $U_{i_x}$ for $x\in U_i-S$, there exists a neighborhood $B_{i_x}\subset B$ of $0$ such that $w_i^1(z_i,t),..., w_i^q (z_i,t)$ defines a regular foliation on $U_{i_x}$ for $t\in B_{i_x}$, so that we can write $\partial w_i^\alpha(z_i,t)= \sum_{\beta=1}^q a_{i_x}^{\alpha\beta}(z_i,t) \wedge w_i^\beta(z_i,t)$ for some $a_{i_x}^{\alpha\beta}(z_i,t)$ on $U_{i_x}\times B_{i_x}\subset U_{i_x}\times B$. When we apply $\bar{\partial}\left(- \right) - \mathcal{L}_{\varphi}\left( - \right)$ on $ \partial w_i^\alpha- \sum_{\beta=1}^q a_{i_x}^{\alpha\beta}  \wedge w_i^\beta=0$, then the condition $\bar{\partial} w_i^\alpha- \mathcal{L}_\varphi\left(w_i^\alpha \right)=0$ implies that $\bar{\partial}a_{i_x}^{\alpha\beta}- \left[ \varphi, a_{i_x}^{\alpha\beta}\right]=0$ on $U_{i_x}\times B_{i_x}$.

Before proceeding to the discussion, we introduce another complex of sheaves which also controls foliated deformations of $\left( M, \mathcal{N}_{\mathcal{F}_0}^* \right)$ but removes the quotient in the degree $1$. We shall define the following complex of sheaves
\begin{align*}
\mathcal{E}_{\mathcal{N}_{\mathcal{F}_0}^*}^\bullet: \mathcal{E}_{\mathcal{N}_{\mathcal{F}_0}^*} \xrightarrow{E_0'} \mathscr{H}om_{\mathcal{O}_M}\left( \mathcal{N}_{\mathcal{F}_0}^*, \Omega_M^1 \right) \xrightarrow{E_1'} \mathscr{H}om_{\mathcal{O}_M}\left( \mathcal{N}_{\mathcal{F}_0}^*, \tilde{\mathcal{S}}^2 \right) \xrightarrow{E_2} \mathscr{H}om_{\mathcal{O}_M}\left( \mathcal{N}_{\mathcal{F}_0}^*, \tilde{\mathcal{S}}^3 \right) \xrightarrow{E_3} \cdots
\end{align*}
where $\mathcal{E}_{\mathcal{N}_{\mathcal{F}_0}^*}$ is the Aityah extension defined in the following way: we recall $(\ref{nnc10})-(\ref{nnc11})$ in the beginning of the proof of Theorem \ref{nn2}. We define $\mathcal{E}_{\mathcal{N}_{\mathcal{F}_0}^*}$ locally on $U_i$ by $\mathcal{E}_{\mathcal{N}_{\mathcal{F}_0}^*}|_{U_i}\cong \mathscr{H}om_{\mathcal{O}_M}\left( \mathcal{N}_{\mathcal{F}_0}^* , \mathcal{N}_{\mathcal{F}_0}^* \right)|_{U_i} \bigoplus \Theta_M|_{U_i}$ such that $(\phi_i,X)\in \mathscr{H}om_{\mathcal{O}_X}\left( \mathcal{N}_{\mathcal{F}_0}^* , \mathcal{N}_{\mathcal{F}_0}^* \right)|_{U_i} \bigoplus \Theta_X |_{U_i}$ and $(\phi_j, X)\in \mathscr{H}om_{\mathcal{O}_X}\left( \mathcal{N}_{\mathcal{F}_0}^* , \mathcal{N}_{\mathcal{F}_0}^* \right)|_{U_j} \bigoplus \Theta_X|_{U_j}$ are identified on $U_{ij}$ if $\phi_i\left(w_{0i}^\alpha\right)=\phi_j\left(w_{0i}^\alpha\right)+\sum_{\beta=1}^q \left[X, h_{0ij}^{\alpha\beta}\right] w_{0j}^\beta$. Then we define $E_0':\mathcal{E}_{\mathcal{N}_{\mathcal{F}_0}^*}\to \mathscr{H}om_{\mathcal{O}_X}\left( \mathcal{N}_{\mathcal{F}_0}^*, \Omega_M^1 \right)$ locally on $U_i$ by $E_0' \left(\left(\phi_i,X\right)\right)\left(w_{0i}^\alpha\right)= \left(\mathcal{L}_X-\phi_i\right)\left(w_{0i}^\alpha\right)$ and then linearly extends to $\Gamma\left( U_i , \mathcal{N}_{\mathcal{F}_0}^* \right)$. We define $E_1'$ by the composition of the natural quotient map $\mathscr{H}om_{\mathcal{O}_M}\left( \mathcal{N}_{\mathcal{F}_0}^*, \Omega_M^1  \right)\to \mathscr{H}om_{\mathcal{O}_M}\left( \mathcal{N}_{\mathcal{F}_0}^*, \frac{\Omega_M^1}{ \mathcal{N}_{\mathcal{F}_0}^* } \right)$. Then we have 
\begin{align}\label{uuc17}
\mathbb{H}^i\left(  M, \mathcal{E}_{\mathcal{N}_{\mathcal{F}_0}^*}^\bullet \right)\cong \mathbb{H}^i\left( M, \mathcal{N}_{\mathcal{F}_0}^{ * \bullet } \right),\,\,\,\,\,\,\,\,i\geq 0
\end{align}
(for the detail on the complex of sheaves $\mathcal{E}_{\mathcal{N}_{\mathcal{F}_0}^*}^{\bullet}$, see Part II). We denote $\mathcal{A}^{0,p}\left( \mathcal{E}_{\mathcal{N}_{\mathcal{F}_0}^*} \right)$ be the sheaf of germs of $C^\infty(0,p)$-forms with coefficients in $\mathcal{E}_{\mathcal{N}_{\mathcal{F}_0}^*}^\bullet$ and denote by $A^{0,p}\left( M, \mathcal{E}_{\mathcal{N}_{\mathcal{F}_0}^*} \right)$ the global section of $\mathcal{A}^{0,p}\left( \mathcal{E}_{\mathcal{N}_{\mathcal{F}_0}^*} \right)$. Then we have the following Dolbeault type bicomplex associated to $\mathcal{E}_{\mathcal{N}_{\mathcal{F}_0}^*}^\bullet$:
{\small{\begin{equation}
\begin{CD}
\cdots  \\
@A\hat{E}_3AA   \\
A^{0,0}\left( M-S, \left(\mathcal{N}_{\mathcal{F}_0}^*\right)^*\otimes \tilde{\mathcal{S}}^3 \right)@>-(-1)^q \bar{\partial}>> \cdots \\
@A\hat{E}_2AA @A\hat{E}_2AA \\
A^{0,0}\left( M-S, \left(\mathcal{N}_{\mathcal{F}_0}^*\right)^*\otimes \tilde{\mathcal{S}}^2\right) @>(-1)^q\bar{\partial}>> A^{0,1}\left(M-S, \left(\mathcal{N}_{\mathcal{F}_0}^*\right)^*\otimes \tilde{\mathcal{S}}^2  \right) @>-(-1)^q\bar{\partial}>> \cdots  \\
@A\hat{E}_1'AA @A \hat{E}_1'AA @A \hat{E}_1'AA \\
A^{0,0}\left(M, \left(\mathcal{N}_{\mathcal{F}_0}^* \right)^* \otimes \Omega_M^1 \right) @>-\bar{\partial}>> A^{0,1}\left(M, \left(\mathcal{N}_{\mathcal{F}_0}^*\right)^*\otimes \Omega_M^1 \right)@>\bar{\partial}>> A^{0,2}\left(M, \left(\mathcal{N}_{\mathcal{F}_0}^*\right)^*\otimes \Omega_M^1 \right)@>-\bar{\partial}>> \cdots \\
@A\hat{E}_0'AA @A\hat{E}_0AA @A\hat{E}_0'AA @A\hat{E}_0'AA \\
A^{0,0}\left(M, \mathcal{E}_{\mathcal{N}_{\mathcal{F}_0}^*} \right) @>\bar{\partial}>> A^{0,1}\left(M, \mathcal{E}_{\mathcal{N}_{\mathcal{F}_0}^*}\right) @>-\bar{\partial}>> A^{0,2}\left(M, \mathcal{E}_{\mathcal{N}_{\mathcal{F}_0}^*} \right) @>\bar{\partial}>> A^{0,3}\left(M, \mathcal{E}_{\mathcal{N}_{\mathcal{F}_0}^*} \right)
\end{CD}
\end{equation}}}
Let us denote the $i$-cohomology group of the complex associated to the above bicomplex by $\textnormal{S}'^i$. Then $\textnormal{S}'^i\cong \textnormal{S}^i (i\geq 0)$ $\textnormal{(\ref{d5})}$ and so $\textnormal{S}'^i\cong \mathbb{H}^i\left( M ,  \mathcal{E}_{\mathcal{N}_{\mathcal{F}_0}^*}^\bullet \right)$ for $i=0,1,2$.

We reinterpret $(\ref{nnc13})$ in terms of the above bicomplex associate to $\mathcal{E}_{\mathcal{N}_{\mathcal{F}_0}^*}^\bullet$. From $(\ref{n50})$, we have
\begin{align}
\sum_{\eta , \gamma =1}^q h_{0ij}^{\alpha\beta} \bar{\partial} B_{j|\mu}^{\eta \gamma} h_{0ji}^{\gamma \xi} - \bar{\partial} B_{i|\mu}^{\alpha \xi} = \bar{\partial} B_{ij|\mu}^{\alpha \xi} \label{nnc12}
\end{align}
We define $\bar{\partial} B_{i|\mu}\in \Gamma\left( U_i, \mathcal{A}^{0,2}\left(\mathscr{H}om_{\mathcal{O}_M}\left( \mathcal{N}_{\mathcal{F}_0}^*, \mathcal{N}_{\mathcal{F}_0}^* \right) \right) \right)$ by 
\begin{align*}
\bar{\partial} B_{i|\mu}:\Gamma\left( U_i, \mathcal{N}_{\mathcal{F}_0}^* \right) &\to \Gamma\left( U_i, \mathcal{A}^{0,2}\left( \mathcal{N}_{\mathcal{F}_0}^* \right)\right)  \\
 w_{0i}^\alpha &\mapsto \sum_{\xi=1}^q \bar{\partial} B_{i|\mu}^{\alpha \xi} w_{0i}^\alpha
\end{align*}

Then from $(\ref{n291})$ and $(\ref{nnc12})$, we have
\begin{align*}
\bar{\partial} B_{j|\mu}\left(w_{0i}^\alpha \right) - \bar{\partial} B_{i|\mu}\left( w_{0i}^\alpha \right)= \sum_{\xi=1}^q \left[ \xi_\mu, h_{0ij}^{\alpha \xi} \right] w_{0j}^\xi
\end{align*}
This implies that
\begin{align}\label{d29}
\left(- \xi_\mu, \left\{ \bar{\partial} B_{i|\mu}\right\} \right)\in A^{0,2}\left( M, \mathcal{E}_{\mathcal{N}_{\mathcal{F}_0}^*} \right)
\end{align}

Then we see that from $(\ref{nnc13})$
\begin{align} \label{ns4}
\left( I_\mu , \phi_\mu, \left( - \xi_\mu , \left\{ \bar{\partial} B_{i|\mu} \right\} \right) \right)\in A^{0,0}\left(M-S, \mathscr{H}om_{\mathcal{O}_M}\left( \mathcal{N}_{\mathcal{F}_0}^*, \tilde{\mathcal{S}}^2 \right) \right) \bigoplus A^{0,1}\left( M, \mathscr{H}om_{\mathcal{O}_M}\left( \mathcal{N}_{\mathcal{F}_0}^*, \Omega_M^1 \right) \right) \bigoplus A^{0,2}\left( M, \mathcal{E}_{\mathcal{N}_{\mathcal{F}_0}^*} \right)
\end{align}
defines a $2$-cocycle in the above Dolbeault type bicomplex associated to $\mathcal{E}_{\mathcal{N}_{\mathcal{F}_0}^*}^\bullet$.

We will define the \"Holder norms on the sections of $\mathcal{A}^{0,q}\left( \Theta_M \right)$ and $\mathcal{A}^{0,q}\left(\mathscr{H}om_{\mathcal{O}_M}\left( \mathcal{N}_{\mathcal{F}_0}^* , \Omega_M^1 \right) \right)$ and $A^{0,q}\left( \mathscr{H}om_{\mathcal{O}_M}\left( \mathcal{N}_{\mathcal{F}_0}^* , \bigwedge^{q+2} \Omega_M^1 \otimes \mathcal{L}_0 \right) \right)$ and $\mathcal{A}^{0,q}\left( M, \mathcal{E}_{\mathcal{N}_{\mathcal{F}_0}^*} \right)$ and apply to the harmonic theory.

We define the \"Holder norm $|-|_{k+\alpha}$ (k: an integer $\geq 2$, $0<\alpha<1$) for section of $\mathcal{A}^{0,q}\left( \Theta_M \right)$ as in $(\ref{d6}),(\ref{d7}),(\ref{d8})$ and introduce a harmonic theory as in $(\ref{d9})$.

We define the \"Holder norm $|-|_{k+\alpha}$ (k: an integer $\geq 2$, $0<\alpha<1$) for section of $\mathcal{A}^{0,q}\left( \mathscr{H}om_{\mathcal{O}_M}  \left( \mathcal{N}_{\mathcal{F}_0}^* , \Omega_M^1 \right) \right)$ as follows: let $\phi\in \Gamma\left(U_i, \mathcal{A}^{0,p}\left(\mathscr{H}om_{\mathcal{O}_M}\left( \mathcal{N}_{\mathcal{F}_0}^*, \Omega_M^1 \right) \right) \right)$ and we write 
\begin{align*}
\phi\left(w_{0i}^\beta\right)=\sum_{\gamma=1}^n A_{i}^{\beta\gamma} dz_i^\gamma ,\,\,\,\,\,\,\,\,\,\,\,A_i^{\beta\gamma}=\frac{1}{q!}\sum A_{i\mu_1\cdots \mu_q}^{\beta\gamma} (z_i)dz_i^{\mu_1}\wedge \cdots \wedge dz_i^{\mu_q} \in \Gamma\left(U_i, \mathcal{A}^{0,q}\right)
\end{align*}
in terms of local coordinates $(z_i^1,..., z_i^n)$ and let
\begin{align*}
\left| \phi \right|_{k+\alpha}^{U_i} = \sum_{h=0}^k \sup \left| D_i^h A_{i \mu_1\cdots \mu_q }^{\beta\gamma}(z_i) \right| + \sup \frac{ \left| D_i^kA_{i \mu_1\cdots \mu_q}^{\beta\gamma}(z_i)- D_i^k A_{i\mu_1\cdots \mu_q}^{\beta\gamma}(y_i) \right| }{\left|z_i-y_i\right|^\alpha}
\end{align*}
where the $``\sup$" is extended over all points $z,y\in U_i$, all indices $\beta, \gamma, \mu_1,...,\mu_q$, and all partial derivatives $D_i^h, D_i^k$ of order $h, k$ with respect to $z_i^1,..., z_i^n, \bar{z}_i^1,..., \bar{z}_i^n$. For $\phi\in A^{0,p}\left( M,  \mathscr{H}om_{\mathcal{O}_M}\left(\Theta_{\mathcal{F}_0} , \Theta_M \right)     \right)$, we define
\begin{align*}
\left| \phi \right|_{k+\alpha} = \max_i \left| \phi \right|_{k+\alpha}^{U_i} 
\end{align*}
We introduce a harmonic theory on $\mathscr{H}om_{\mathcal{O}_M}\left( \mathcal{N}_{\mathcal{F}_0}^*, \Omega_M^1 \right)$.
\begin{align}
\textnormal{We denote by $\mathfrak{d}''$ the adjoint operator of $\bar{\partial}$, and $\square''=\mathfrak{d}''\bar{\partial} + \bar{\partial} \mathfrak{d}''$ and $G''$ the Green's operator.}
\end{align}

We define the \"Holder norm $|-|_{k+\alpha}$ (k: an integer $\geq 2$, $0<\alpha<1$) for section of $\mathcal{A}^{0,r}\left( \mathscr{H}om_{\mathcal{O}_M}  \left( \mathcal{N}_{\mathcal{F}_0}^* , \bigwedge^{q+2} \Omega_M^1 \otimes \mathcal{L}_0 \right) \right)$ as follows: let $\phi\in \Gamma\left(U_i, \mathcal{A}^{0,r}\left(\mathscr{H}om_{\mathcal{O}_M}\left( \mathcal{N}_{\mathcal{F}_0}^* , \bigwedge^{q+2} \Omega_M^1 \otimes \mathcal{L}_0 \right) \right) \right)$ and we write 
\begin{align*}
\phi\left( w_{0i}^\beta \right)=\sum_{\eta=1}^n \sum_{\sigma_1,..., \sigma_{q+2} }A_{i \sigma_1 \cdots \sigma_{q+2}}^{\beta} dz_i^{\sigma_1}\wedge \cdots \wedge dz_i^{\sigma_{q+2}}\,\,\,\,\,\,\,\,\,\,\,A_{i\sigma_1\cdots \sigma_{q+2}}^{\beta}=\frac{1}{r!}\sum A_{i\mu_1\cdots \mu_r}^{\beta\sigma_1\cdots \sigma_{q+2}} (z_i)dz_i^{\mu_1}\wedge \cdots \wedge dz_i^{\mu_r} \in \Gamma\left(U_i, \mathcal{A}^{0,r}\right)
\end{align*}
in terms of local coordinates $(z_i^1,..., z_i^n)$ and let
\begin{align*}
\left| \phi \right|_{k+\alpha}^{U_i} = \sum_{h=0}^k \sup \left| D_i^h A_{i \mu_1\cdots \mu_q }^{\beta \sigma_1\cdots \sigma_r }(z_i) \right| + \sup \frac{ \left| D_i^kA_{ i\mu_1\cdots \mu_r}^{ \beta \sigma_1 \cdots \sigma_{q+2}}(z_i)- D_i^k A_{i \mu_1 \cdots \mu_r}^{\beta \sigma_1 \cdots \sigma_{q+2}}(y_i) \right| }{\left|z_i-y_i\right|^\alpha}
\end{align*}
where the $``\sup$" is extended over all points $z,y\in U_i$, all indices $\beta, \sigma_1,..., \sigma_{q+2}, \mu_1,...,\mu_q$, and all partial derivatives $D_i^h, D_i^k$ of order $h, k$ with respect to $z_i^1,..., z_i^n, \bar{z}_i^1,..., \bar{z}_i^n$. For $\phi\in A^{0,p}\left( M, \mathscr{H}om_{\mathcal{O}_M}\left( \mathcal{N}_{\mathcal{F}_0}^* , \bigwedge^{q+2} \Omega_M^1 \otimes \mathcal{L}_0  \right)     \right)$, we define
\begin{align*}
\left| \phi \right|_{k+\alpha} = \max_i \left| \phi \right|_{k+\alpha}^{U_i} 
\end{align*}

We define the \"Holder norm $|-|_{k+\alpha}$ (k: an integer $\geq 2$, $0<\alpha<1$) for section of $\mathcal{A}^{0,q}\left(\mathcal{E}_{\mathcal{N}_{\mathcal{F}_0}^*} \right)$ in a similar way to $(\ref{d10}),(\ref{d11})$, and $(\ref{d12})$. We introduce a harmonic theory on $\mathcal{E}_{\mathcal{N}_{\mathcal{F}_0}^*}$. 
\begin{align}\label{ns5}
\textnormal{We denote by $\tilde{\mathfrak{d}}$ the adjoint operator of $\bar{\partial}$, and $\widetilde{\square}=\tilde{\mathfrak{d}}\bar{\partial} + \bar{\partial} \tilde{\mathfrak{d}}$ and $\tilde{G}$ the Green's operator on $\mathcal{E}_{\mathcal{N}_{\mathcal{F}_0}^*}$.}
\end{align}
We recall $(\ref{d30})$ and $(\ref{d31})$ as in the proof of convergence in Theorem \ref{tt2}. With this preparation, we will show that for a fixed integer $k\geq 2$ and $0< \alpha<1$, the inductive construction of $\varphi, w_i^\alpha , h_{ij}^{\alpha\beta}$ in the previous subsection can be carried out in such a way that
\begin{align*}
\left|\varphi  \right|_{k+\alpha} \ll A(t)\\
\left| w_i^{\alpha }- w_{0i}^\alpha \right|_{k+\alpha} \ll A(t) \\
\left| h_{ij}^{\alpha\beta }- h_{0ij}^{\alpha\beta} \right|_{k+1+\alpha} \ll A(t)
\end{align*}
Then it suffices to prove that for $\mu=1,2,3,\cdots$
\begin{align}
\left|\varphi^\mu \right|_{k+\alpha} \ll A(t) \label{ns1}\\
\left| w_i^{\alpha \mu}- w_{0i}^\alpha \right|_{k+\alpha} \ll A(t) \\
\left| h_{ij}^{\alpha\beta \mu}- h_{0ij}^{\alpha\beta} \right|_{k+1+\alpha} \ll A(t)\label{ns2}
\end{align}
for some proper choice of constants $c>b>0$. We prove $(\ref{ns1})_\mu-(\ref{ns2})_\mu$ by induction on $\mu$. For $\mu=1$ we have $(\ref{ns3})$ and $(\ref{st8})$ and $(\ref{st9})$, and the linear term of $A(t)$ is $\frac{b}{16}\left(t_1+ \cdots  + t_r \right)$. Therefore $(\ref{ns1})_1-(\ref{ns2})_1$ holds if $b$ is sufficiently large. 

Now assume that $(\ref{ns1})_{\mu-1}-(\ref{ns2})_{\mu-1}$ are satisfied. We will derive $(\ref{ns1})_\mu- (\ref{ns2})_\mu$. In the following $Z_1,Z_2, Z_3,\cdots$ will denote constants which depend only on $k,\alpha, M , \mathcal{N}_{\mathcal{F}_0}^*$.

From $(\ref{ns4})$, we have $\bar{\partial}\left(-\xi_\mu, \left\{ \bar{\partial} B_{i|\mu} \right\} \right)=0$. Then we take $\left( \tilde{\varphi}_\mu, \left\{ \tilde{ \lambda}_{i|\mu} \right\}\right):= \tilde{\mathfrak{d}}\tilde{G} \left( -\xi_\mu,  \left\{\bar{\partial} B_{i|\mu} \right\}\right)\in A^{0,1}\left(M, \mathcal{E}_{\mathcal{N}_{\mathcal{F}_0}^*} \right)$ from $(\ref{ns5})$, where $\tilde{\varphi}_\mu\in A^{0,1}\left(  M, \Theta_M \right)$ and $\tilde{\lambda}_{i|\mu} \in \Gamma\left( U_i, \mathcal{A}^{0,1}\left( \mathscr{H}om_{\mathcal{O}_M}\left( \mathcal{N}_{\mathcal{F}_0}^*, \mathcal{N}_{\mathcal{F}_0}^* \right) \right) \right)$ defined by $\tilde{\lambda}_{i|\mu}:\Gamma\left( U_i, \mathcal{N}_{\mathcal{F}_0}^* \right)\to \Gamma\left( U_i, \mathcal{A}^{0,1}\left( \mathcal{N}_{\mathcal{F}_0}^* \right) \right), w_{0i}^\alpha \mapsto \sum_{\xi=1}^q \tilde{\lambda}_{i|\mu}^{\alpha \xi} w_{0i}^\xi$. Then
\begin{align}\label{sd5}
\left| \left( \tilde{\varphi}_\mu, \left\{ \tilde{\lambda}_{i|\mu} \right\} \right) \right|_{k+\alpha} \ll Z_1 \left| \left(- \xi_\mu, \left\{\bar{\partial} B_{i|\mu} \right\} \right) \right|_{k-1+\alpha}
\end{align}
as in the same way with $(\ref{d25})$. From $(\ref{d9})$, we have
\begin{align*}
\tilde{\varphi}_\mu = \bold{H} \tilde{\varphi}_\mu + \left( \bar{\partial} \mathfrak{d} + \mathfrak{d} \bar{\partial} \right) G \tilde{\varphi}_\mu
\end{align*}
It follows that
\begin{align*}
\bar{\partial}\left( \bold{H} \tilde{\varphi}_\mu  +\mathfrak{d}\bar{\partial} G \tilde{\varphi}_\mu \right) &= - \xi_\mu\\
\mathfrak{d}\left( \bold{H} \tilde{\varphi}_\mu + \mathfrak{d} \bar{\partial} G \tilde{\varphi}_\mu \right)&=0
\end{align*}
On the other hand, we define $q_{ij|\mu} \in \Gamma\left( U_{ij}, \mathcal{A}^{0,0}\left( \mathscr{H}om_{\mathcal{O}_M}\left( \mathcal{N}_{\mathcal{F}_0}^* , \mathcal{N}_{\mathcal{F}_0}^* \right) \right) \right)$ by 
\begin{align*}
q_{ij|\mu}:\Gamma \left(U_{ij}, \mathcal{N}_{\mathcal{F}_0 }^* \right)&\to \Gamma \left(U_{ij}, \mathcal{A}^{0,0}\left(\mathcal{N}_{\mathcal{F}_0}^* \right) \right)\\
w_{0i}^\alpha &\mapsto  \sum_{\beta=1}^q \left[ \mathfrak{d}G\tilde{\varphi}_\mu, h_{0ij}^{\alpha\beta}\right] w_{0j}^\beta,
\end{align*}
and linearly extends to $\Gamma\left( U_{ij}, \mathcal{N}_{\mathcal{F}_0}^* \right)$. Then $\left\{ q_{ij|\mu}\right\}\in C^1\left( \mathcal{U}, \mathcal{A}^{0,0}\left( \mathscr{H}om_{\mathcal{O}_M}\left( \mathcal{N}_{\mathcal{F}_0}^*, \mathcal{N}_{\mathcal{F}_0}^* \right) \right) \right)$ defines a $1$-cocycle. Let $\left\{ \rho_i(z) \right\}$ be a partition of unity subordinate to the covering $\mathcal{U}=\left\{ U_i \right\}$. By setting $q_{i|\mu}=\sum_k \rho_k(z)q_{ik|\mu}$, we see that  $\left\{q_{i|\mu} \right\}\in C^0\left( \mathcal{U}, \mathcal{A}^{0,0}\left( \mathscr{H}om_{\mathcal{O}_M}\left(  \mathcal{N}_{\mathcal{F}_0}^*, \mathcal{N}_{\mathcal{F}_0}^*  \right)   \right)  \right)$ such that $q_{i|\mu} -q_{j|\mu}=q_{ij|\mu}$ where
\begin{align*}
q_{i|\mu}:\Gamma\left(U_i, \mathcal{N}_{\mathcal{F}_0}^* \right) &\to \Gamma\left(U_i, \mathcal{A}^{0,0}\left(\mathcal{N}_{\mathcal{F}_0}^* \right)\right)\\
    w_{0i}^\alpha &\mapsto \sum_{\gamma=1}^q q_{i|\mu}^{\alpha\gamma} w_{0i}^\gamma, \,\,\,\,\,\,\,\, q_{i|\mu}^{\alpha\beta}\in \Gamma\left(U_i, \mathcal{A}^{0,0}\right)
\end{align*}
Then we see that
\begin{align}\label{ns38}
\left(\varphi_\mu'', \left\{ \lambda_{i|\mu}'' \right\} \right) \in A^{0,1}\left( M, \mathcal{E}_{\mathcal{N}_{\mathcal{F}_0}^*} \right),\,\,\,\,\,\,\,\,\,\, \mathfrak{d}\varphi_\mu''=0,\,\,\,\,\,\,\,\bar{\partial}\left(\varphi_\mu'' , \left\{ \lambda_{i|\mu}''\right\} \right)= \left( - \xi_\mu,  \left\{ \bar{\partial} B_{i|\mu} \right\} \right)
\end{align}
where $\varphi_\mu'':= \bold{H}\tilde{\varphi}_\mu + \mathfrak{d}\bar{\partial} G \tilde{\varphi}_\mu$ and $\lambda_{i|\mu}''\in \Gamma\left(U_i, \mathcal{A}^{0,1}\left( \mathscr{H}om_{\mathcal{O}_M}\left( \mathcal{N}_{\mathcal{F}_0}^*, \mathcal{N}_{\mathcal{F}_0}^* \right) \right) \right)$ defined by
\begin{align*}
\lambda_{i|\mu}'':\Gamma\left( U_i, \mathcal{N}_{\mathcal{F}_0}^* \right) &\to \Gamma\left( U_i, \mathcal{A}^{0,1}\left( \mathcal{N}_{\mathcal{F}_0}^* \right) \right) \\
  w_{0i}^\alpha &\mapsto \sum_{\beta=1}^q \left( \tilde{\lambda}_{i|\mu}^{\alpha \beta} - \bar{\partial} q_{i|\mu}^{\alpha \beta} \right) w_{0i}^\beta
\end{align*}
and we have in a similar way with $(\ref{ds3})$
\begin{align}\label{ns27}
\left| \varphi_\mu'' \right|_{k+\alpha},\,\,\,\,\,\, \left| \lambda_{i|\mu}'' \right|_{k+\alpha}^{U_i} \ll Z_2  \left|\left( - \xi_\mu , \left\{ \bar{\partial} B_{i|\mu}\right\} \right) \right|_{k-1+\alpha}
\end{align}

From $(\ref{ns4})$, we consider
{\small{\begin{align*}
&0\in A^{0,2}\left( M, \mathcal{E}_{\mathcal{N}_{\mathcal{F}_0}^*} \right)\\
&\phi_\mu'=\left\{ w_{0i}^\alpha \mapsto A_{i|\mu}^\alpha - \bar{\partial} C_{i|\mu}^\alpha + \sum_{\gamma=1}^q B_{i|\mu}^{\alpha \gamma} w_{0i}^\gamma + \mathcal{L}_{\varphi_\mu''}\left( w_{0i}^\alpha \right) - \sum_{\xi=1}^q \lambda_{i|\mu}''^{\alpha \xi} w_{0i}^\xi \right\} \in A^{0,1}\left(M, \mathscr{H}om_{\mathcal{O}_M}\left(  \mathcal{N}_{\mathcal{F}_0}^*, \Omega_M^1 \right) \right) \\
&\psi_\mu=\left\{ w_{0i}^\alpha \mapsto  E_{i|\mu}^\alpha + \sum_{\beta=1}^q w_{0i}^1 \wedge \cdots \wedge C_{i|\mu}^\beta \wedge \cdots \wedge w_{0i}^q \wedge dw_{0i}^\alpha + w_{0i}^1 \wedge \cdots \wedge w_{0i}^q \wedge  \partial C_{i|\mu}^\alpha \right\}  \in A^{0,0}\left( M, \mathscr{H}om_{\mathcal{O}_M}\left( \mathcal{N}_{\mathcal{F}_0}^*, \bigwedge^{q+2} \Omega_M^1 \otimes \mathcal{L}_0 \right) \right)
\end{align*}}}
Then $\psi_\mu |_{M-S}\in A^{0,0}\left( M-S, \mathscr{H}om_{\mathcal{O}_M}\left( \mathcal{N}_{\mathcal{F}_0}^*, \tilde{\mathcal{S}}^2 \right) \right)$ and $\left( 0, \phi_\mu', \psi_\mu |_{M-S} \right)$ defines a $2$-cocycle in the Dolbeault type bicomplex associated to $\mathcal{E}_{\mathcal{N}_{\mathcal{F}_0}^*}^{\bullet}$. Since $\mathbb{H}^2\left( M, \mathcal{N}_{\mathcal{F}_0}^{*\bullet} \right)= \mathbb{H}^2\left( M, \mathcal{E}_{\mathcal{N}_{\mathcal{F}_0}^*}^\bullet \right)=0$ by the assumption of Theorem \ref{nn2}, there exists $\left(\chi_\mu, \left\{ \delta_i^{\chi_\mu} \right\} \right)\in A^{0,1}\left( M, \mathcal{E}_{\mathcal{N}_{\mathcal{F}_0}^*} \right)$ and $\Sigma_\mu \in A^{0,0}\left( M, \mathscr{H}om_{\mathcal{O}_M}\left( \mathcal{N}_{\mathcal{F}_0}^*, \Omega_M^1  \right) \right)$ such that $-\bar{\partial}\left(\chi_\mu, \left\{ \delta_i^{\chi_\mu}\right\} \right)=0, -\bar{\partial} \Sigma_\mu + \hat{E}_0'\left(\chi_\mu, \left\{ \delta_i^{\chi_\mu} \right\} \right)=\phi_\mu'$ and $\hat{E}_1'\left( \Sigma_\mu \right)= \psi_\mu |_{M-S}$. By using the following Lemma, we will choose appropriate $\left(\chi_\mu, \left\{ \delta_i^{\chi_\mu}  \right\} \right)$ and $\Sigma_\mu$ in a way that $\varphi, w_i^\alpha$ and $h_{ij}^{\alpha\beta}$ converge.

We define $E_1^\sharp: \mathscr{H}om_{\mathcal{O}_M}\left( \mathcal{N}_{\mathcal{F}_0}^*, \Omega_M^1 \right) \to \mathscr{H}om_{\mathcal{O}_M}\left( \mathcal{N}_{\mathcal{F}_0}^*, \bigwedge^{q+2} \Omega_M^1 \otimes \mathcal{L}_0 \right)$ in the following way: for $\phi\in \Gamma\left( U_i, \mathscr{H}om_{\mathcal{O}_M}\left( \mathcal{N}_{\mathcal{F}_0}^*, \Omega_M^1 \right) \right)$, we have
\begin{align*}
E_1^\sharp\left(\phi \right): \Gamma\left(U_i, \mathcal{N}_{\mathcal{F}_0}^*\right)&\to \Gamma\left( U_i, \bigwedge^{q+2} \Omega_M^1 \otimes \mathcal{L}_0 \right) \\ 
                       w_{0i}^\alpha &\mapsto \sum_{\beta=1}^q w_{0i}^1 \wedge \cdots \wedge \phi\left( w_{0i}^\beta \right) \wedge \cdots \wedge w_{0i}^q \wedge dw_i^\alpha + w_{0i}^1 \wedge \cdots \wedge w_{0i}^q \wedge d \phi\left(w_{0i}^\alpha \right)
\end{align*}
This induces
\begin{align*}
E_1^\sharp: A^{0,p}\left( M, \mathscr{H}om_{\mathcal{O}_M}\left( \mathcal{N}_{\mathcal{F}_0}^*, \Omega_M^1 \right) \right) \to A^{0,p}\left( M, \mathscr{H}om_{\mathcal{O}_M}\left( \mathcal{N}_{\mathcal{F}_0}^*, \bigwedge^{q+2} \Omega_M^1\otimes \mathcal{L}_0 \right) \right)
\end{align*}
Then $\hat{E}_1':A^{0,p}\left(M, \mathscr{H}om_{\mathcal{O}_M}\left( \mathcal{N}_{\mathcal{F}_0}^*, \Omega_M^1 \right)\right)\to A^{0,p}\left( M-S, \mathscr{H}om_{\mathcal{O}_M}\left( \mathcal{N}_{\mathcal{F}_0}^*, \tilde{\mathcal{S}}^2 \right)\right)$ is the composition of $E_1^\sharp$ and the restriction on $M-S$.

\begin{lemma}\label{nnc2}
Suppose that $\left(\varphi, \left\{ \delta_i\right\}\right) \in A^{0,1}\left(M, \mathcal{E}_{\mathcal{N}_{\mathcal{F}_0}^*} \right)$ where $\varphi \in A^{0,1}\left( M, \Theta_M \right)$ and $\delta_i\in \Gamma\left( U_i, \mathscr{H}om_{\mathcal{O}_M}\left( \mathcal{N}_{\mathcal{F}_0}^*, \mathcal{N}_{\mathcal{F}_0}^* \right) \right)$ with $\delta_i\left( w_{0i}^\alpha \right)= \delta_j \left( w_{0i}^\alpha \right)+ \sum_{\beta=1}^q \left[ \varphi, h_{0ij}^{\alpha\beta} \right] w_{0j}^\beta$ for $\alpha=1,...,q$, and $V\in A^{0,1}\left( M, \mathscr{H}om_{\mathcal{O}_M}\left( \mathcal{N}_{\mathcal{F}_0}^*, \Omega_M^1 \right) \right)$ and $I\in A^{0,0}\left(M, \mathscr{H}om_{\mathcal{O}_M}\left( \mathcal{N}_{\mathcal{F}_0}^*, \bigwedge^{q+2} \Omega_M^1 \otimes \mathcal{L}_0\right) \right)$ such that  $I|_{M-S}\in A^{0,0}\left( M-S, \mathscr{H}om_{\mathcal{O}_M}\left( \mathcal{N}_{\mathcal{F}_0}^* , \tilde{\mathcal{S}}^2 \right) \right)$ and $\left( 0, V + \hat{E}_0'\left(\varphi, \left\{ \delta_i \right\} \right), I|_{M-S} \right)$ defines a $2$-cocycle in the Dolbeault type bicomplex associated to $\mathcal{E}_{\mathcal{N}_{\mathcal{F}_0}^*}^{\bullet}$, and moreover $(-1)^q \bar{\partial} I + E_1^\sharp\left( V \right)=0$. Then we can find $\left(\chi, \left\{\delta_i^\chi \right\} \right) \in A^{0,1}\left( M, \mathcal{E}_{\mathcal{N}_{\mathcal{F}_0}^*} \right)$, and $\Sigma\in A^{0,0}\left( M, \mathscr{H}om_{\mathcal{O}_M}\left( \mathcal{N}_{\mathcal{F}_0}^*, \Omega_M^1 \right) \right)$ in such a way that
\begin{align}
&\widetilde{\square}\left( \chi, \left\{ \delta_i^\chi \right\} \right)=0 \label{ns6}\\
-\bar{\partial} \Sigma\left(w_{0i}^\alpha \right) +\mathcal{L}_{\chi}\left( w_{0i}^\alpha\right)- & \delta_i^{\chi} \left( w_{0i}^\alpha \right) = V\left( w_{0i}^\alpha \right) + \mathcal{L}_{\varphi}-  \delta_i \left( w_{0i}^\alpha \right) \label{ns11}\\
\delta_i^\chi\left( w_{0i}^\alpha \right)& = \delta_j^\chi\left( w_{0i}^\alpha \right) + \sum_{\beta=1}^q \left[ \chi, h_{0ij}^{\alpha\beta} \right] w_{0j}^\beta \label{dt7}  \\
&E_1^\sharp\left(  \Sigma \right)= I \label{ns7}
\end{align}
\begin{align*}
\left|\left( \chi, \left\{ \delta^\chi \right\} \right) \right|_{k+\alpha} &\ll Z\left( \left|\left( \varphi, \left\{ \delta_i \right\} \right) \right|_{k+\alpha} + \left| V\right|_{k-1+\alpha} + \left| I \right|_{k-1 + \alpha} \right) \\
\left| \Sigma \right|_{k+\alpha}&\ll Z\left( \left|\left( \varphi, \left\{ \delta_i \right\}\right) \right|_{k+\alpha} + \left| V\right|_{k-1+\alpha} + \left| I \right|_{k-1+\alpha} \right) 
\end{align*}
where $Z$ is a constant which is independent of $\left( \varphi, \left\{\delta_i \right\}\right), V, I$

\end{lemma}
\begin{remark}
We may write $E_1^\sharp\left(\Sigma \right)=I$ instead of $\hat{E}_1'\left(\Sigma \right)= I |_{M-S}$. In fact, assume that $\hat{E}_1'\left(\Sigma \right) = I|_{M-S}$ and let $E_1^\sharp\left(\Sigma \right)= I'$. Then $I-I'=0$ on $M-S$. Since $(-1)^q\bar{\partial} \left( I- I'\right)= - E_1^\sharp(V) - (-1)^q \bar{\partial}E_1^\sharp\left(\Sigma \right)= -E_1^\sharp(V)- E_1^\sharp\bar{\partial}\left(\Sigma \right)= -E_1^\sharp\left( V+ \bar{\partial}\Sigma \right)= - E_1^\sharp\hat{E}_0'\left( \chi- \varphi, \left\{\delta_i^\chi - \delta_i \right\} \right)=0$, so that $I-I'$ is holomorphic on $M$, and so $I=I'$ on $M$.
\end{remark}

\begin{proof}
For any triple $\zeta=\left(\left( \varphi , \left\{ \delta_i \right\} \right) , V, I \right)$ as above, let 
\begin{align*}
||\zeta || &= \left| \left(\varphi  , \left\{ \delta_i^\varphi \right\}\right)\right|_{k+\alpha} + \left|V \right|_{k-1+\alpha} + \left| I \right|_{k-1+\alpha} \\
\iota(\zeta)&= \inf \left( \left|\left(\chi, \left\{\delta_i^\chi \right\} \right)\right|_{k+\alpha}  + \left|\Sigma  \right|_{k+\alpha}  \right)
\end{align*}
where $\inf$ is taken with respect to all solutions $\left( \left(\chi, \left\{\eta_i^\chi \right\}\right), \Sigma \right)$ of the equalities $(\ref{ns6})-(\ref{ns7})$. It suffices to prove the existence of $Z$ such that
\begin{align*}
\iota\left( \zeta \right)\leq Z ||\zeta ||\,\,\,\,\,\,\,\,\,\textnormal{for all triple $\zeta$}
\end{align*}
Assume that there is no such constant $Z$. Then we can find a sequence $\zeta^{(1)},\zeta^{(2)},\cdots, \zeta^{(v)},\cdots$ of triple $\zeta^{(v)}=\left( \left( \varphi^{(v)} , \left\{ \eta_i^{(v)} \right\} \right), V^{(v)}, I^{(v)}\right)$ such that
\begin{align*}
\iota\left(\zeta^{(v)} \right)=1\,\,\,\,\,\,\,\,\textnormal{and}\,\,\,\,\,\,\,\,\left|\left| \zeta^{(v)}\right|\right| < \frac{1}{v}
\end{align*}
The first equality implies the existence of $\left(\chi^{(v)}, \left\{ \delta_i^{\chi^{(v)}}  \right\} \right)\in A^{0,1}\left( M, \mathcal{E}_{\mathcal{N}_{\mathcal{F}_0}^*} \right)$ and $\Sigma^{(v)}\in A^{0,0}\left( M, \mathscr{H}om_{\mathcal{O}_M}\left( \mathcal{N}_{\mathcal{F}_0}^*, \Omega_M^1 \right) \right)$ such that
\begin{align}
&\tilde{\square}\left( \chi^{(v)}, \left\{ \delta_i^{\chi^{(v)}}\right\} \right)=0\\
-\bar{\partial} \Sigma^{(v)}\left( w_{0i}^\alpha \right)+ \mathcal{L}_{\chi^{(v)}}\left(w_{0i}^\alpha \right)&- \delta_i^{\chi^{(v)}}\left(  w_{0i}^\alpha \right)= V^{(v)}\left( w_{0i}^\alpha \right) + \mathcal{L}_{\varphi^{(v)}}\left( w_{0i}^\alpha \right) - \delta_i^{(v)}\left( w_{0i}^\alpha \right)\label{sd10}\\
\delta_i^{\chi^{(v)}}\left(w_{0i}^\alpha \right) &= \delta_j^{\chi^{(v)}}\left( w_{0i}^\alpha\right) + \sum_{\beta=1}^q \left[ \chi^{(v)}, h_{0ij}^{\alpha\beta} \right] w_{0j}^\beta \\
&\,\,\,\,\,E_1^\sharp\left( \Sigma^{(v)} \right)= I^{(v)} \label{sd11}\\
&\left| \left( \chi^{(v)}, \left\{ \delta_i^{\chi^{(v)}} \right\} \right)   \right|_{k+\alpha} + \left| \Sigma^{(v)} \right|_{k+\alpha} <2 \label{ns8}
\end{align}
Then $(\ref{ns8})$ implies that we may assume that $ \left( \chi, \left\{ \delta_i^{\chi} \right\}\right)=\lim \left(  \chi^{(v)}, \left\{ \delta_i^{\chi^{(v)}} \right\} \right)$ and $\Sigma= \lim \Sigma^{(v)}$ exist in the norm $|-|_k$, so that $\left( \chi, \left\{ \delta_i^\chi \right\} \right)$ and $\Sigma$ are of class $C^k$. We note that since $\left( \chi, \left\{ \delta_i^{\chi} \right\} \right)$ satisfies an elliptic partial differential equation $\widetilde{\square}\left( \chi, \left\{ \delta_i^\chi \right\} \right)=0$, $\left(\chi, \left\{ \delta_i^\chi \right\} \right)$ is of $C^\infty$. Moreover since $\widetilde{\square}$ is strongly elliptic, by using `a priori estimate' (see \cite{Kod05} Theorem 4.3 p.436) we can show that $\left( \chi^{(v)}, \left\{ \delta_i^{\chi^{(v)}} \right\} \right)$ converges to $\left(\chi, \left\{ \delta_i^\chi \right\} \right)$ in the norm $|-|_{k+\alpha}$ in the same way with $(\ref{ns10})$. On the other hand, from $(\ref{ns11})$ we have $\bar{\partial}\Sigma= \hat{E}_0'\left(\chi, \left\{ \delta_i^\chi\right\} \right)$, so that $\Sigma$ is $C^\infty$. In the same with $(\ref{ns12})$, we can show $\Sigma^{(v)}$ converge to $\Sigma$ in $|-|_{k+\alpha}$. Then from $(\ref{ns11})$ and $(\ref{ns7})$
\begin{align}
-\bar{\partial}\left( \Sigma^{(v)}- \Sigma\right)\left(w_{0i}^\alpha \right) +\mathcal{L}_{\chi^{(v)}- \chi}\left( w_{0i}^\alpha\right)- & \delta_i^{\chi^{(v)}- \chi} \left( w_{0i}^\alpha \right) = V^{(v)}\left( w_{0i}^\alpha \right) + \mathcal{L}_{\varphi^{(v)}}-  \delta_i^{(v)} \left( w_{0i}^\alpha \right) \label{sd13}\\
&E_1^\sharp\left(  \Sigma^{(v)}- \Sigma \right)= I^{(v)} \label{sd14}
\end{align}
On the other hand, we have
\begin{align*}
\left|\left( \chi^{(v)}, \left\{ \delta_i^{\chi^{(v)}} \right\}\right)- \left(\chi, \left\{ \delta_i^\chi \right\} \right) \right|_{k+\alpha} < \frac{1}{4},\,\,\,\,\,\,\,\,\,\left| \Sigma^{(v)} - \Sigma \right|_{k+\alpha}<\frac{1}{4}
\end{align*}
for sufficiently large integer $v$. This contradicts to $\iota\left( \zeta \right)=1$. This completes the proof of Lemma \ref{nnc2}.
\end{proof}

In Lemma \ref{nnc2}, we set
\begin{align}\label{ns26}
\left( \varphi, \left\{ \delta_i\right\}\right)&:=\left( \varphi_\mu'', \left\{ \lambda_{i|\mu}'' \right\} \right)  \\
V&:= \phi_\mu'=\left\{ w_{0i}^\alpha \mapsto  A_{i|\mu}^\alpha - \bar{\partial} C_{i|\mu}^\alpha + \sum_{\gamma=1}^q B_{i|\mu}^{\alpha \gamma} w_{0i}^\gamma \right\}\notag \\
I&:= \psi_\mu=\left\{ w_{0i}^\alpha \mapsto E_{i|\mu}^\alpha + \sum_{\beta=1}^q w_{0i}^1 \wedge \cdots \wedge C_{i|\mu}^\beta \wedge \cdots \wedge w_{0i}^q \wedge dw_{0i}^\alpha + w_{0i}^1 \wedge \cdots \wedge w_{0i}^q \wedge \partial C_{i|\mu}^\alpha \right\}. \notag
\end{align}
We will estimate $\left(\varphi_\mu'', \left\{ \lambda_{i|\mu}''\right\} \right)$ and $V$ and $I$. First we will estimate $(\ref{nn3})-(\ref{ns14})$. By induction hypothesis $(\ref{ns1})_{\mu-1}-(\ref{ns2})_{\mu-1}$, we have
\begin{align}\label{ns15}
\xi_\mu\equiv_\mu  \bar{\partial} \varphi^{\mu-1}-\frac{1}{2}\left[\varphi^{\mu-1}, \varphi^{\mu-1}   \right] \Longrightarrow \left|\xi_\mu \right|_{k-1+\alpha} \leq Z_3 \left|\varphi^{\mu-1} \right|_{k+\alpha} \left|\varphi^{\mu-1}\right| _{k+\alpha} \leq \frac{Z_3b}{c}A(t)
\end{align}
and we have
{\small{\begin{align}\label{ns16}
-A_{i|\mu}^\alpha \equiv_\mu \bar{\partial} &w_i^{\alpha (\mu-1)} - \mathcal{L}_{\varphi^{\mu-1}}\left( w_i^{\alpha(\mu-1)} \right)\equiv_\mu \left[ \mathcal{L}_{\varphi^{\mu-1}}\left(w_i^{\alpha(\mu-1)} - w_{0i}^\alpha \right)        \right]_\mu=\left[  - \partial i_{\varphi^{\mu-1}}\left( w_i^{\alpha(\mu-1)} -w_{0i}^\alpha \right)   + i_{\varphi^{\mu-1}}\partial \left(w_i^{\alpha(\mu-1)}   - w_{0i}^\alpha\right)    \right]_\mu  \notag \\
& \Longrightarrow \left| A_{i|\mu}^\alpha  \right|_{k-1+\alpha}\leq Z_4 \left| \varphi^{\mu-1} \right|_{k+\alpha} \left| w_i^{\alpha(\mu-1)} - w_{0i}^\alpha \right|_{k+\alpha} \ll Z_4\frac{b}{c} A(t)
\end{align}}}
and we have
\begin{align}\label{ns17}
&-\sum_{\eta=1}^q B_{ij|\mu}^{\alpha \eta} h_{0ij}^{\eta \beta} \equiv_\mu \bar{\partial} h_{ij}^{\alpha\beta(\mu-1)} - \left[ \varphi^{\mu-1}, h_{ij}^{\alpha\beta(\mu-1)}   \right] \notag\\
&\Longrightarrow  \left| B_{ij|\mu}^{\alpha \eta} \right|_{k+\alpha} \leq Z_5\left| \varphi^{\mu-1}\right|_{k+\alpha} \left| h_{ij}^{\alpha\beta(\mu-1)}- h_{0ij}^{\alpha\beta}  \right|_{k+1+\alpha} \ll Z_5\frac{b}{c} A(t)
\end{align}
and we have
\begin{align}\label{ns18}
&C_{ij|\mu}^\alpha\equiv_\mu w_i^{\alpha(\mu-1)} - \sum_{\beta=1}^q h_{ij}^{\alpha\beta(\mu-1)} w_j^{\beta(\mu-1)}\equiv_\mu \left[ -\sum_{\beta=1}^q \left(h_{ij}^{\alpha\beta(\mu-1)}- h_{0ij}^{\alpha\beta}\right)\left( w_j^{\beta(\mu-1)} - w_{0j}^\beta  \right)  \right]_\mu  \notag \\
&\Longrightarrow  \left| C_{ij|\mu}^\alpha \right|_{k+\alpha} \ll Z_6\left| h_{ij}^{\alpha\beta(\mu-1)}- h_{0ij}^{\alpha\beta} \right|_{k+\alpha} \left| w_j^{\beta(\mu-1)} - w_{0j}^\beta   \right|_{k+\alpha} \ll \frac{Z_7 b}{c} A(t)
\end{align}
and we have
\begin{align}\label{ns19}
&\sum_{\beta=1}^q D_{ijk|\mu}^{\alpha\beta} h_{0ik}^{\beta \gamma} \equiv_\mu h_{ik}^{\alpha \gamma(\mu-1)} - \sum_{\beta=1}^q h_{ij}^{\alpha\beta(\mu-1)}(f_{jk}, t) h_{jk} ^{\beta\gamma(\mu-1)}(z_k,t) \equiv_\mu -\left[ \sum_{\beta=1}^q\left(h_{ij}^{\alpha\beta(\mu-1)}(f_{jk}, t) - h_{0ij}^{\alpha\beta}   \right) \left( h_{jk}^{\beta \gamma(\mu-1)} - h_{0jk}^{\beta \gamma} \right)   \right]_\mu \notag \\
&\Longrightarrow \left| D_{ijk|\mu}^{\alpha\beta}   \right|_{k+1+\alpha} \leq Z_8\sum_{\beta=1}^q \left| h_{ij}^{\alpha\beta(\mu-1)}- h_{0ij}^{\alpha\beta} \right|_{k+1+\alpha} \left| h_{jk}^{\beta \gamma(\mu-1)} - h_{0jk}^{\beta \gamma} \right|_{k+1+\alpha} \ll \frac{Z_9 b}{c} A(t)
\end{align}
and we have
\begin{align*}
E_{i|\mu}^\alpha& \equiv_\mu w_i^{1(\mu-1)}\wedge \cdots \wedge w_i^{q(\mu-1)}\wedge dw_i^{\alpha(\mu-1)}\\
&= (w_i^{1(\mu-1)}- w_{0i}^1 +w_{0i}^1)\wedge\cdots \wedge (w_i^{q(\mu-1)}- w_{0i}^q + w_{0i}^q)\wedge d \left(  w_i^{\alpha(\mu-1)} - w_{0i}^\alpha + w_{0i}^\alpha   \right)\\
&=(w_i^{1(\mu-1)}- w_{0i}^1)\wedge \cdots \wedge (w_i^{q(\mu-1)}- w_{0i}^q) \wedge d \left( w_i^{\alpha(\mu-1)} - w_{0i}^\alpha \right)\\
&+\sum_{\beta=1}^q (w_i^{1(\mu-1)}- w_{0i}^1)\wedge \cdots \wedge \overbrace{w_{0i}^\beta}^{\beta-\textnormal{th}}\wedge \cdots \wedge (w_i^{q(\mu-1)}- w_{0i}^q) \wedge d\left( w_i^{\alpha(\mu-1)} - w_{0i}^\alpha \right)\\
&+ \cdots \\
&+ \sum w_{0i}^1\wedge \cdots \wedge (w_i^{\beta(\mu-1)}- w_{0i}^\beta) \wedge \cdots \wedge (w_i^{\gamma (\mu-1)} - w_{0i}^\gamma) \wedge \cdots \wedge w_{0i}^q \wedge dw_{0i}^\alpha           \\
&+ \sum_{\beta=1}^q w_{0i}^1\wedge \cdots \wedge \left(w_i^{\beta(\mu-1)} - w_{0i}^\beta \right) \wedge \cdots \wedge w_{0i}^q \wedge d\left( w_i^{\alpha(\mu-1)}- w_{0i}^\alpha \right)\\
&\left(+ \sum_{\beta=1}^q w_{0i}^1\wedge \cdots \wedge \left( w_i^{\beta(\mu-1)} - w_{0i}^\beta \right) \wedge \cdots \wedge w_{0i}^q \wedge dw_{0i}^\alpha + w_{0i}^1\wedge \cdots \wedge w_{0i}^q \wedge d\left(  w_i^{\alpha(\mu-1)} - w_{0i}^\alpha \right) \right)
\end{align*}
Then we have
\begin{align}\label{ns20}
\left| E_{i|\mu}^\alpha \right|_{k-1+\alpha} \ll Z_{10}\left( \frac{b}{c} \right)^q A(t) + Z_{11}\left(\frac{b}{c} \right)^{q-1} A(t) + \cdots + Z_{12} \frac{b}{c}A(t) \ll Z_{13} \frac{b}{c} A(t)
\end{align} 
We recall that from $(\ref{ne11})$ and $(\ref{n50})$ and $(\ref{n52})$
\begin{align} 
\sum_{\beta, \eta=1}^q h_{0ij}^{\alpha \beta} D_{jk|\mu}^{\beta \eta} h_{0ji}^{\eta \xi} - D_{ik|\mu}^{\alpha \xi} + D_{ij|\mu}^{\alpha \xi} &= D_{ijk|\mu}^{\alpha \xi} \label{ns21}\\
\sum_{\eta, \gamma=1}^q h_{0ij}^{\alpha \eta} B_{j|\mu}^{\eta \gamma} h_{0ji}^{\gamma \xi} - B_{i|\mu}^{\alpha \xi} & = B_{ij|\mu}^{\alpha \xi} - \bar{\partial} D_{ij|\mu}^{\alpha \xi} \label{ns22}\\
\sum_{\beta=1}^q h_{0ij}^{\alpha \beta} C_{j|\mu}^\beta - C_{i|\mu}^\alpha &= C_{ij|\mu}^\alpha - \sum_{\beta=1}^q D_{ij|\mu}^{\alpha \beta} w_{0i}^\beta \label{ns23}
\end{align}
By using a partition of unity subordinate to the covering $\left\{ U_i \right\}$, we may choose $D_{jk|\mu}$ from $D_{ijk|\mu}$, and choose $B_{j|\mu}$ from $B_{ij|\mu}$ and $\bar{\partial} D_{ij|\mu}$, and choose $C_{i|\mu}$ from $C_{ij|\mu}$ and $D_{ij|\mu}$. Hence we may assume that from $(\ref{ns19})$ and $(\ref{ns21})$
\begin{align}\label{ns24}
\left| D_{ij|\mu} \right|_{k+1+\alpha} \ll Z_{14}\frac{b}{c} A(t)
\end{align} 
and from $(\ref{ns17})$ and $(\ref{ns24})$ and $(\ref{ns22})$
\begin{align}\label{ns28}
\left| B_{i|\mu} \right|_{k+\alpha} \ll Z_{15} \frac{b}{c} A(t)
\end{align}
and from $(\ref{ns18})$ and $(\ref{ns24})$ and $(\ref{ns23})$
\begin{align}\label{ns29}
\left| C_{i|\mu} \right|_{k+\alpha} \ll Z_{16}\frac{b}{c} A(t)
\end{align}
With this preparation, now we estimate $(\ref{ns26})$ and then apply Lemma \ref{nnc2}. Let us estimate $\left(  \varphi_\mu'', \left\{ \lambda_{i|\mu}'' \right\} \right)$. From $(\ref{ns27})$ and $(\ref{ns15})$ and $(\ref{ns28})$, we have
\begin{align}\label{ns36}
\left| \left( \varphi_\mu'', \left\{ \lambda_{i|\mu}'' \right\} \right) \right|_{k+\alpha} \ll Z_{17}\frac{b}{c} A(t)
\end{align}
We estimate $V:=\phi_\mu'=\left\{ w_{0i}^\alpha \mapsto A_{i|\mu}^\alpha - \bar{\partial} C_{i|\mu}^\alpha + \sum_{\gamma=1}^q B_{i|\mu}^{\alpha \gamma} w_{0i}^\gamma \right\} \in A^{0,1}\left( M, \mathscr{H}om_{\mathcal{O}_M}\left( \mathcal{N}_{\mathcal{F}_0}^*, \Omega_M^1 \right) \right)$. From $(\ref{ns16})$ and $(\ref{ns29})$ and $(\ref{ns28})$
\begin{align}\label{sd27}
\left| V \right|_{k-1+\alpha} \ll Z_{18} \frac{b}{c} A(t)
\end{align}
We estimate
{\Small{\begin{align*}
I:= \psi_\mu=\left\{ w_{0i}^\alpha \mapsto E_{i|\mu}^\alpha + \sum_{\beta=1}^q w_{0i}^1 \wedge \cdots \wedge C_{i|\mu}^\beta \wedge \cdots \wedge w_{0i}^q \wedge dw_{0i}^\alpha + w_{0i}^1 \wedge \cdots \wedge w_{0i}^q \wedge \partial C_{i|\mu}^\alpha \right\} \in A^{0,0}\left( M, \mathscr{H}om_{\mathcal{O}_M}\left( \mathcal{N}_{\mathcal{F}_0}^*, \bigwedge^{q+2} \Omega_M^1 \otimes \mathcal{L}_0 \right) \right)
\end{align*}}}
From $(\ref{ns20})$ and $(\ref{ns29})$, we have
\begin{align}\label{sd28}
\left| I \right|_{k-1+\alpha} \ll Z_{19} \frac{b}{c} A(t)
\end{align}

By Lemma \ref{nnc2}, we can find $\left(\chi_\mu, \left\{ \delta_i^{\chi_\mu} \right\} \right)\in A^{0,1}\left( M, \mathcal{E}_{\mathcal{N}_{\mathcal{F}_0}^*} \right)$ and $\Sigma_\mu \in A^{0,0}\left( M, \mathscr{H}om_{\mathcal{O}_M}\left( \mathcal{N}_{\mathcal{F}_0}^*, \Omega_M^1 \right) \right)$ such that
{\small{\begin{align}
\bar{\partial} \chi_\mu =0, \,\,\,\,\,\bar{\partial}\delta_i^{\chi_\mu}&=0\,\,\,\left( \widetilde{\square} \left( \chi_\mu , \left\{ \delta_i^{\chi_\mu} \right\} \right)=0\right)\\
-\bar{\partial} \Sigma_\mu \left(w_{0i}^\alpha\right) + \mathcal{L}_{\chi_\mu}\left(w_{0i}^\alpha \right) - \delta_i^{\chi_\mu}\left( w_{0i}^\alpha \right)&= A_{i|\mu}^\alpha - \bar{\partial} C_{i|\mu}^\alpha + \sum_{\gamma=1}^q B_{i|\mu}^{\alpha \gamma} w_{0i}^\gamma + \mathcal{L}_{\varphi_\mu''}\left(w_{0i}^\alpha \right) - \lambda_{i|\mu}''\left(w_{0i}^\alpha \right) \label{sd30} \\
\sum_{\beta=1}^q w_{0i}^1 \wedge \cdots \wedge \Sigma_i^\beta \wedge \cdots \wedge w_{0i}^q \wedge dw_{0i}^\alpha + w_{0i}^1\wedge \cdots \wedge w_{0i}^q \wedge \partial \Sigma_{\mu |i}^\alpha&= E_{i|\mu}^\alpha + \sum_{\beta=1}^q w_{0i}^1 \wedge \cdots \wedge C_{i|\mu}^\beta \wedge \cdots \wedge w_{0i}^q \wedge dw_{0i}^\alpha + w_{0i}^1 \wedge \cdots \wedge w_{0i}^q \wedge \partial C_{i|\mu}^\alpha \label{sd31}\\
\left| \left(\chi_\mu, \left\{\delta_i^{\chi_\mu} \right\} \right) \right|_{k+\alpha}& ,\,\,\,\,\,\,\,\left| \Sigma_\mu \right|_{k+\alpha} \ll Z_{20}\frac{b}{c} A(t) \label{ns30}
\end{align}}}

In view of the equality
\begin{align*}
\chi_\mu = \bold{H} \chi_\mu + \left( \bar{\partial}\mathfrak{d}+ \mathfrak{d}\bar{\partial} \right)G\chi_\mu ,
\end{align*}
we get
\begin{align}\label{ns32}
- \bar{\partial}\left(\Sigma_{i | \mu }^\alpha - \mathcal{L}_{\mathfrak{d}G \chi_\mu }\left(w_{0i}^\alpha \right) \right) + \mathcal{L}_{\bold{H}\chi_\mu - \varphi_\mu''}\left(w_{0i}^\alpha \right) - \left( \delta_i^{\chi_\mu} - \lambda_{i|\mu}''  \right)\left( w_{0i}^\alpha \right) = A_{i|\mu}^\alpha - \bar{\partial} C_{i|\mu}^\alpha + \sum_{\gamma=1}^q B_{i|\mu}^{\alpha \gamma} w_{0i}^\gamma 
\end{align}
On the other hand, we note that $\left\{U_{ij|\mu} \right\}$ is a $1$-cocycle with coefficient in $\mathcal{A}^{0,1}\left( \mathscr{H}om_{\mathcal{O}_M}\left( \mathcal{N}_{\mathcal{F}_0}^*, \mathcal{N}_{\mathcal{F}_0}^*   \right) \right)$ where
\begin{align*}
U_{ij|\mu}: \Gamma\left( U_{ij}, \mathcal{N}_{\mathcal{F}_0}^* \right)& \to \Gamma\left( U_{ij}, \mathcal{A}^{0,0}\left( \mathcal{N}_{\mathcal{F}_0}^* \right) \right)\\
w_{0i}^\alpha &\mapsto \sum_{\beta=1}^q \left[ \mathfrak{d} G \chi_\mu,  h_{0ij}^{\alpha\beta} \right] w_{0j}^\beta
\end{align*}

Let $\left\{ \rho_i(z)\right\}$ be a partition of unity subordinate  to the covering $\mathcal{U}=\left\{ U_i \right\}$. By setting $U_{i|\mu}= \sum_k \rho_k(z)U_{ik|\mu}$, we see that $\left\{ U_{i|\mu} \right\}\in C^0\left( \mathcal{U}, \mathscr{H}om_{\mathcal{O}_M}\left( \mathcal{N}_{\mathcal{F}_0}^*, \mathcal{N}_{\mathcal{F}_0}^* \right)\right)$ such that $U_{i|\mu}- U_{j|\mu}= U_{ij|\mu}$ where
\begin{align*}
U_{i|\mu}:\Gamma\left( U_i, \mathcal{N}_{\mathcal{F}_0}^* \right) &\to \Gamma\left( U_i, \mathcal{A}^{0,0}\left( \mathcal{N}_{\mathcal{F}_0}^* \right) \right) \\
      w_{0i}^\alpha &\mapsto \sum_{\gamma=1}^q U_{i|\mu}^{\alpha \gamma} w_{0i}^\gamma,\,\,\,\,\,\,\,\,U_{i|\mu}^{\alpha \gamma} \in \Gamma\left( U_i, \mathcal{A}^{0,0} \right)
\end{align*}
and we have
\begin{align}\label{ns31}
\sum_{\beta, \gamma=1}^q h_{0ij}^{\alpha\beta} U_{j|\mu}^{\beta \gamma} w_{0j}^\gamma - \sum_{\gamma=1}^q U_{i|\mu}^{\alpha \gamma} w_{0i}^\gamma = - \sum_{\beta=1}^q \left[ \mathfrak{d}G \chi_\mu, h_{0ij}^{\alpha\beta} \right] w_{0j}^\beta
\end{align}
Then we have from $(\ref{ns30})$
\begin{align}\label{ns37}
\left| U_{i|\mu} \right|_{k+1+\alpha} \ll Z_{21}\left|\mathfrak{d}G\chi_\mu \right|_{k+1+\alpha} \ll Z_{22} \left| \chi_\mu\right|_{k+\alpha} \ll Z_{23} \frac{b}{c} A(t)
\end{align}
On the other hand, $(\ref{ns31})$ implies that
\begin{align*}
R_i:\Gamma\left( U_i, \mathcal{N}_{\mathcal{F}_0}^* \right) &\to \Gamma\left( U_i, \mathcal{A}^{0,0}\left( \Omega_M^1 \right) \right) \\
 w_{0i}^\alpha &\mapsto \mathcal{L}_{ \mathfrak{d} G \chi_\mu}\left( w_{0i}^\alpha \right) - U_{i|\mu}\left( w_{0i}^\alpha \right)
\end{align*}
defines a global section $R=\left\{R_i \right\} \in A^{0,0}\left( M, \mathscr{H}om_{\mathcal{O}_M}\left( \mathcal{N}_{\mathcal{F}_0}^*, \Omega_M^1 \right)\right)$. Then $\Sigma_\mu - R\in A^{0,0}\left( M, \mathscr{H}om_{\mathcal{O}_M}\left( \mathcal{N}_{\mathcal{F}_0}^*, \Omega_M^1 \right) \right)$. Then we have from $(\ref{ns32})$
\begin{align}\label{n40}
- \bar{\partial}\left(\Sigma_{i|\mu}^\alpha - \mathcal{L}_{\mathfrak{d}G \chi_\mu}\left(w_{0i}^\alpha \right) + \sum_{\xi=1}^q U_{i|\mu}^{\alpha \xi} w_{0i}^\xi \right) + \mathcal{L}_{\bold{H}\chi_\mu - \varphi_\mu''}\left(w_{0i}^\alpha \right) - \left( \delta_i^{\chi_\mu}-  \lambda_{i|\mu}''  - \bar{\partial} U_{i|\mu}\right)\left( w_{0i}^\alpha \right) = A_{i|\mu}^\alpha - \bar{\partial} C_{i|\mu}^\alpha + \sum_{\gamma=1}^q B_{i|\mu}^{\alpha \gamma} w_{0i}^\gamma 
\end{align}
We note that
\begin{align}\label{n41}
&\sum_{\beta, \xi=1}^q h_{0ij}^{\alpha\beta}\left( \delta_j^{\chi_\mu \beta \xi} - \lambda_{j|\mu}''^{\beta \xi} - \bar{\partial} U_{j|\mu}^{\beta \xi}  \right) w_{0j}^\xi - \sum_{\xi=1}^q \left( \delta_i^{\chi_\mu\alpha\xi}  - \lambda_{i|\mu}^{\alpha \xi} - \bar{\partial} U_{i|\mu}^{\alpha \xi}\right) w_{0i}^\xi\\
&= - \sum_{\beta=1}^q \left[ \chi_\mu - \varphi_\mu'' - \bar{\partial} \mathfrak{d} G \chi_\mu,  h_{0ij}^{\alpha\beta} \right] w_{0j}^\beta = - \sum_{\beta=1}^q \left[ \bold{H} \chi_\mu - \varphi_\mu'', h_{0ij}^{\alpha\beta} \right] w_{0j}^\beta  \notag
\end{align}
If we set 
\begin{align}
\varphi_\mu' &:= \bold{H} \chi_\mu - \varphi_\mu'' \label{ns33} \\
w_{i|\mu}'^\alpha &:= \Sigma_{i|\mu}^\alpha - \mathcal{L}_{\mathfrak{d} G \chi_\mu}\left( w_{0i}^\alpha \right) + U_{i|\mu}\left( w_{0i}^\alpha \right) \label{ns34}\\
V_{i|\mu}^{\alpha \xi} &:= - \delta_i^{\chi_\mu \alpha \xi} + \lambda_{i|\mu}''^{\alpha \xi} + \bar{\partial} U_{i|\mu}^{\alpha \xi} \label{ns35}
\end{align}
then we have from $(\ref{ns33})$ and $(\ref{ns30})$ and $(\ref{ns36})$
\begin{align}\label{ns48}
\left| \varphi_\mu' \right|_{k+\alpha}\ll  Z_{24}\frac{b}{c} A(t)
\end{align}
and we have from $(\ref{ns34})$ and $(\ref{ns30})$ and $(\ref{ns37})$
\begin{align}\label{ns51}
\left| w_{i|\mu}'^\alpha \right|_{k+\alpha} \ll Z_{25} \frac{b}{c} A(t)
\end{align}
and we have from $(\ref{ns35})$ and $(\ref{ns30})$ and $(\ref{ns36})$ and $(\ref{ns37})$
\begin{align}\label{ns42}
\left| V_{i|\mu}^{\alpha \xi} \right|_{k+\alpha} \ll Z_{26} \frac{b}{c} A(t)
\end{align}
From $(\ref{ns38})$ we have
\begin{align*}
-\bar{\partial}\varphi_\mu'=-\bar{\partial}\left(\bold{H} \chi_\mu  - \varphi_\mu'' \right) =  \bar{\partial} \varphi_\mu''= - \xi_\mu,\,\,\,\,\,\,\,\, \mathfrak{d} \varphi_\mu'=\mathfrak{d}\left(\bold{H}\chi_\mu - \varphi_\mu'' \right)= -\mathfrak{d}\varphi_\mu''=0
\end{align*}
and so from $(\ref{ne18})$, we have $\varphi_\mu:= - \varphi_\mu'$, so that
\begin{align}\label{ns52}
\bar{\partial}\varphi_\mu=-\xi_\mu,\,\,\,\,\,\,\,\,\mathfrak{d}\varphi_\mu=0
\end{align}
and $(\ref{n40})$ and $(\ref{n41})$ implies $(\ref{n51})$ and $(\ref{ne23})$. Since $\bar{\partial}\left(V_{i|\mu}^{\alpha\beta}- B_{i|\mu}^{\alpha\beta} \right)=0$, there exists $P_{i|\mu}^{\alpha \beta} \in \Gamma\left(U_i, \mathcal{A}^{0,0} \right)$ such that $\bar{\partial} P_{i|\mu}^{\alpha \beta}= V_{i|\mu}^{\alpha \beta} - B_{i|\mu}^{\alpha\beta}$. Then by Potential theoretic Lemma (see \cite{DN55}), we have
\begin{align} \label{ns43}
\left|  P_{i|\mu}^{\alpha \beta} \right|_{k+ 1+\alpha} \ll Z_{27} \left|  V_{i|\mu}^{\alpha\beta} - B_{i|\mu}^{\alpha \beta} \right|_{k+\alpha}
\end{align}
and from $(\ref{ns42})$ and $(\ref{ns28})$ and $(\ref{ns43})$, we have
\begin{align}\label{ns49}
\left| P_{i|\mu}^{\alpha\beta} \right|_{k+1+\alpha} \ll  Z_{28} \frac{b}{c} A(t)
\end{align}

Then we have from $(\ref{ne18})$ and $(\ref{ne19})$ and $(\ref{ne22})$
\begin{align}
\varphi_\mu&= -\varphi_\mu' =   -\bold{H}\chi_\mu + \varphi_\mu'' \label{ns45} \\
w_{i|\mu}^\alpha&= C_{i|\mu}^\alpha - w_{i|\mu}'^\alpha + \sum_{\beta=1}^q P_{i|\mu}^{\alpha\beta} w_{0i}^\beta \label{ns46}\\
h_{ij|\mu}^{\alpha\beta}&= \sum_{\gamma=1}^q D_{ij|\mu}^{\alpha \gamma} h_{0ij}^{\gamma \beta} + \sum_{\gamma=1}^q P_{i|\mu}^{\alpha \gamma} h_{0ij}^{\gamma \beta} - \sum_{\gamma=1}^q h_{0ij}^{\alpha \gamma} P_{j|\mu}^{\gamma \beta} \label{ns47}
\end{align}
Then from $(\ref{ns45})$ and $(\ref{ns48})$
\begin{align}
\left|\varphi_\mu \right|_{k+\alpha} \ll Z_{24}\frac{b}{c} A(t)
\end{align}
and from $(\ref{ns46})$ and $(\ref{ns29})$ and $(\ref{ns51})$ and $(\ref{ns49})$
\begin{align*}
\left| w_{i|\mu}^\alpha \right|_{k+\alpha} \ll Z_{29} \frac{b}{c} A(t)
\end{align*}
and from $(\ref{ns47})$ and $(\ref{ns24})$ and $(\ref{ns49})$
\begin{align*}
\left| h_{ij|\mu}^{\alpha\beta} \right|_{k+1+\alpha} \ll Z_{30} \frac{b}{c} A(t)
\end{align*}
Then we can choose $b$ and $c$ satisfying $(\ref{ns1})_\mu-(\ref{ns2})_\mu$. Then $\varphi(t)$ and $w_i^\alpha$ and  $h_{ij}^{\alpha\beta}$ converge with respect to the \"Holder norm $\left| - \right|_{k+\alpha}, \left|- \right|_{k+\alpha},\left|- \right|_{k+1+\alpha}$, respectively for some sufficiently small neighborhood $\Delta_\epsilon \subset B$ of $0$. Consequently $\varphi(t)$ is a $C^k$ vector $(0,1)$-form on $M\times \Delta_\epsilon$ and $w_i^\alpha$ is a $C^k$-form of the form $w_i^\alpha= \sum_{\beta=1}^n w_i^{\alpha\beta}\left(z_i,t \right)dz_i^\beta $ on $U_i\times \Delta_\epsilon$, and $h_{ij}^{\alpha\beta}$ is a $ C^{k+1}$ function on $U_{ij}\times \Delta_\epsilon$. From $(\ref{ns52})$, we have $\mathfrak{d}\varphi_\mu =0$ for $\mu\geq 2$. Hence $\varphi(t)$ is a solution of the quasi-linear partial differential equation of order $2$ 
\begin{align}\label{ns53}
\sum_{\lambda=1}^r \frac{\partial^2}{\partial t_\lambda \partial \bar{t}_\lambda} \varphi(t) +\square \varphi(t) - \mathfrak{d}\left[ \varphi(t), \varphi(t)\right]=\bar{\partial}\mathfrak{d}\varphi_1(t)
\end{align}
We may assume that $(\ref{ns53})$ is a quasi-linear elliptic partial differential equation on $M\times \Delta_\epsilon$ (see \cite{Kod05} p.281). Therefore its solution $\varphi(t)$ is $C^\infty$ on $M\times \Delta_\epsilon$ (see \cite{DN55} Theorem 5 or \cite{Kod05} Appendix \S 8). Then $\varphi(t)$ determines a complex analytic family $\pi: \mathcal{M}\to \Delta_\epsilon$ (see \cite{Kod05} p.281-283). On the other hand, $w_i^\alpha$ and $h_{ij}^{\alpha\beta}$ are holomorphic with respect to the complex structure $\varphi(t)$ by $(\ref{ns54})$ and $(\ref{ns55})$, respectively. $(\ref{n212}),(\ref{nnc5})$ and $(\ref{n21})$ determines a locally free subsheaf $\mathcal{N}_\mathcal{F}^*$ of $\Omega_{\mathcal{M}/\Delta_\epsilon}^1$ satisfying the integrability condition. From the definition of linear terms, we infer that the foliated Kodaira-Spence map $\varphi_0: T_0\left(\Delta_\epsilon \right)\to \mathbb{H}^1\left( M, \mathcal{N}_{\mathcal{F}_0}^{*\bullet} \right)$ is bijective. This completes the proof of Theorem \ref{nn2}.

\end{proof}

\section{Theorem of stability for deformations of foliated complex analytic structures in terms of cotangent sheaves}

Our proof of theorem of existence of deformations of foliated complex analytic structures in terms of cotangent sheaves (Theorem \ref{nn2}) provides the proof of theorem of stability for deformations of foliated complex analytic structures in terms of cotangent sheaves (Theorem \ref{nst1}). We define a complex of sheaves which we truncate the $0$-th term of the dual leaf complex
{\Small{\begin{align*}
\mathscr{H}om_{\mathcal{O}_M}\left( \mathcal{N}_{\mathcal{F}_0}^*, \frac{\Omega_M^1}{\mathcal{N}_{\mathcal{F}_0}^*} \right)^\bullet : \mathscr{H}om_{\mathcal{O}_M}\left( \mathcal{N}_{\mathcal{F}_0}^*, \frac{\Omega_M^1}{ \mathcal{N}_{\mathcal{F}_0}^* } \right)\xrightarrow{E_1} \mathscr{H}om_{\mathcal{O}_M}\left( \mathcal{N}_{\mathcal{F}_0}^*, \tilde{\mathcal{S}}^2\right) \xrightarrow{E_2} \mathscr{H}om_{\mathcal{O}_M}\left( \mathcal{N}_{\mathcal{F}_0}^*, \tilde{\mathcal{S}}^3 \right) \xrightarrow{E_3} \cdots
\end{align*}}}
We will denote the $i$-th cohomology group of $\mathscr{H}om_{\mathcal{O}_M}\left( \mathcal{N}_{\mathcal{F}_0}^* , \frac{\Omega_M^1}{\mathcal{N}_{\mathcal{F}_0}^*} \right)^\bullet$ by $\mathbb{H}^i\left( M, \mathscr{H}om_{\mathcal{O}_M}\left( \mathcal{N}_{\mathcal{F}_0}^* , \frac{\Omega_M^1}{\mathcal{N}_{\mathcal{F}_0}^*} \right)^\bullet \right)$.

\begin{theorem}[Theorem of stability for deformations of foliated complex analytic structures in terms of cotangent sheaves] \label{nst1}
Let $\left( M, \mathcal{N}_{\mathcal{F}_0}^* \right)$ be a compact foliated complex manifold with $\mathcal{N}_{\mathcal{F}_0}^*$ locally free. Assume that the first cohomology $\mathbb{H}^1\left( M, \mathscr{H}om_{\mathcal{O}_M}\left( \mathcal{N}_{\mathcal{F}_0}^*, \frac{\Omega_M^1}{\mathcal{N}_{\mathcal{F}_0}^*} \right)^\bullet \right)=0$. Then for any complex analytic family $\pi:\mathcal{M}\to B$ of deformations of $\pi^{-1}(0)=M, 0\in B$ in the sense of Kodaira-Spencer, there exists an open neighborhood $N\subset B$ of $0$ and a locally free subsheaf $\mathcal{N}_{\mathcal{F}}^*$ of $\Omega_{\frac{\mathcal{M}|_N}{N}}^1$ such that $\left(\mathcal{M}|_N, \mathcal{N}_{\mathcal{F}}^*\right)$ defines a foliated analytic family $\pi|_N: \left(\mathcal{M}|_N , \mathcal{N}_{\mathcal{F}}^* \right) \to N$ of deformations of $\left(M, \mathcal{N}_{\mathcal{F}_0}^* \right)=\pi^{-1}(0)$ in terms of cotangent sheaves.
\end{theorem}

\begin{proof}
We copy the proof of Theorem \ref{nn2}. The difference is that we assume that $\varphi(t)$ (determined by $\mathcal{M}$ over some neighborhood $N'\subset B$ of $0$) is given from the beginning and replace $\varphi^\mu$ by $\varphi(t)$ in the proof of Theorem \ref{nn2}. Accordingly we replace $\xi_\mu$ and $\varphi_\mu$ by $0$. Then the obstructions for solving $(\ref{nn7})-(\ref{nn4})$ are in $\mathbb{H}^1\left( M, \mathscr{H}om_{\mathcal{O}_M}\left( \mathcal{N}_{\mathcal{F}_0}^* ,\frac{\Omega_M^1}{ \mathcal{N}_{\mathcal{F}_0}^* } \right)^\bullet \right)=0$. Hence we can construct $\mathcal{N}_{\mathcal{F}}^*$ on $\mathcal{M}|_N$ for some neighborhood $N\subset N'$ of $0$. 
\end{proof}

\section{Theorem of completeness of deformations of foliated complex analytic structures in terms of cotangent sheaves}

\subsection{Change of parameters} \label{r6} \

Consider a foliated complex analytic family $\left( \mathcal{M}, \mathcal{N}_\mathcal{F}^*, B, \pi \right)$ with $\mathcal{N}_\mathcal{F}^*$ locally free in terms of cotangent sheaves of deformations of $\left( M_t, \mathcal{N}_{\mathcal{F}_t}^* \right)=\pi^{-1}(t), t\in B$, where $B$ is a domain of $\mathbb{C}^m$. Let $D$ be a domain of $\mathbb{C}^{m'}$ and $s:u\to t= s(u), u\in D$, a holomorphic map of $D$ into $B$. Then by changing the parameter from $t$ to $u$, we will construct a foliated analytic family $\left\{ \left( M_{s(u)}, \mathcal{N}_{\mathcal{F}_{s(u)}}^*  \right) | u \in D  \right\}$ on the parameter space $D$ in the following.

Let $\mathcal{M} \times_B D:=\left\{ \left( p, u \right)\in \mathcal{M}\times D | \pi(p)=s(u) \right\}$. Then we have the following commutative diagram
\begin{center}
$\begin{CD}
\mathcal{M}\times_B D @>p>> \mathcal{M}\\
@V\pi' VV @VV\pi V\\
D @>s>> B
\end{CD}$
\end{center}
such that $\left( \mathcal{M}\times_B D, D, \pi' \right)$ is a complex analytic family in the sense of Kodaira-Spencer and $\pi'^{-}(u)= M_{s(u)}$. On the other hand, $\Pi^*\mathcal{N}_\mathcal{F}^*$ is a locally free subsheaf of $\Omega_{\mathcal{M}\times D}^1$, where $\Pi: \mathcal{M}\times D \to \mathcal{M}$ is the natural projection. From the exact sequence $0 \to \Pi^* \Omega_D^1 \to \Omega_{\mathcal{M}\times D}^1\to \Omega_{\mathcal{M}\times D/D}^1\to 0 $ and $\Pi^* \Omega_D^1 \cap \Pi^* \mathcal{N}_\mathcal{F}^* = 0$, we see that $\Pi^* \mathcal{N}_\mathcal{F}^*$ can be considered to be a locally free susbsheaf of $\Omega_{\mathcal{M}\times D/D}^1$. Since $\mathcal{M}\times_B D$ is a complex submanifold of $\mathcal{M}\times D$, we have an exact sequence $0\to \mathcal{I}/\mathcal{I}^2 \xrightarrow{d} \Omega_{\mathcal{M}\times D}^1|_{\mathcal{M}\times_B D}\to \Omega_{\mathcal{M}\times_B D}^1 \to 0$ where $\mathcal{I}$ is the ideal sehaf of $\mathcal{O}_{\mathcal{M}\times D}$ defining $\mathcal{M}\times_B D$. The ideal sheaf $\mathcal{I}$ is locally defined by the equation $t=s(u)$, i.e. $t_\alpha= s_\alpha\left(u_1,..., u_{m'}\right),\alpha=1,...,m$. Since the image $d\left( \mathcal{I}\right)$ of $\mathcal{I}$ in $\Omega_{\mathcal{M}\times D}^1$ is locally generated by the form $dt_\alpha- \sum_{\beta=1}^{m'} \frac{\partial s_\alpha}{\partial u_\beta} d u_\beta$, so that $d\left( \mathcal{I}\right)\bigcap \pi^* \Omega_D^1=0$. This implies that the above exact sequence induces an exact sequence $0\to \mathcal{I}/\mathcal{I}^2 \xrightarrow{d} \to \Omega_{\mathcal{M}\times D/D}^1|_{\mathcal{M}\times_B D} \to \Omega_{\mathcal{M}\times_B D/D}^1 \to 0$. We have $d\left(\mathcal{I}/\mathcal{I}^2 \right)\cap p^* \mathcal{N}_\mathcal{F}^*=0$. Here we note that $\Pi^* \mathcal{N}_{\mathcal{F}}^*|_{\mathcal{M}\times_B D}= p^* \mathcal{N}_\mathcal{F}^*$. Hence $p^* \mathcal{N}_{\mathcal{F}}^*$ is isomorphically mapped under $\Omega_{\mathcal{M}\times D/D}^1|_{\mathcal{M}\times_B D} \to \Omega_{\mathcal{M}\times_B D/D}^1$, so that $p^* \mathcal{N}_\mathcal{F}^*$ can be considered to be a locally free subsheaf of $\Omega_{\mathcal{M}\times_B D / D}^1$. Let $S$ be the locus of $\mathcal{M}$ such that $\frac{\Omega_\mathcal{M}^1}{\mathcal{N}_\mathcal{F}^*}$ is not locally free, and let $S' =S\times_B D$. Then $p$ induces $p: \left(\mathcal{M}\times_B D \right) - S' \to \mathcal{M}- S$. Since the pullback commutes with the differential and wedge products, $d_{\mathcal{M}/B} \left( \mathcal{N}_\mathcal{F}^* \right) \subset \mathcal{N}_\mathcal{F}^*\wedge \Omega_{\mathcal{M}/B}^1$ on $\mathcal{M}- S$ implies that $d_{\mathcal{M}\times_B D/D}\left(p^* \mathcal{N}_\mathcal{F}^* \right) \subset p^*{\mathcal{N}}_\mathcal{F}^* \wedge \Omega_{\mathcal{M}\times_B D/D}^1 $ on $\left(\mathcal{M}\times_B D \right) - S'$. Moreover each $u\in D$, $p^*\mathcal{N}_\mathcal{F}^*|_{\pi'^{-1}(u)}= \mathcal{N}_{\mathcal{F}_{u(s)}}^*\to \Omega_{\mathcal{M}\times_B D/D}^1|_{\pi'^{-1}(u)}= \Omega_{M_{s(u)}}^1$ is injective, so that $\frac{\Omega_{\mathcal{M}\times_B D/D}^1}{ p^* \mathcal{N}_\mathcal{F}^* }$ is flat over $D$. Let $\pi'':\left(\mathcal{M}\times_B D \right)- S' \to D$ be the restriction of $\pi'$ to $\left( \mathcal{M}\times_B D \right) - S'$. Then for each $u\in D, \frac{\Omega_{\mathcal{M}\times_B D/D}^1}{p^* \mathcal{N}_\mathcal{F}^* }|_{\pi''^{-1}(u)}= \frac{\Omega_{M_{s(u)}}^1}{\mathcal{N}_{\mathcal{F}_{s(u)}}^*}|_{M_{s(u)}- S_{s(u)}}$ is locally free, so that $\frac{\Omega_{\mathcal{M}\times_B D/D}^1}{p^* \mathcal{N}_\mathcal{F}^* }$ is locally free on $\left(\mathcal{M}\times_B D \right) - S'$. This implies that $\left\{ \left( M_{s(u)}, \mathcal{N}_{\mathcal{F}_{s(u)}}^* \right) | u \in D \right\}$ forms a foliated complex analytic family $\left( \mathcal{M} \times_B D, p^* \mathcal{N}_\mathcal{F}^*, D, \pi' \right)$.

\begin{definition}
The foliated complex analytic family $\left( \mathcal{M}\times_B D, p^* \mathcal{N}_\mathcal{F}^*, D , \pi' \right)$ is called the foliated complex analytic family induced from $\left( \mathcal{M}, \mathcal{N}_\mathcal{F}^*, B, \pi \right)$ by the holomorphic map $s : D \to B$.
\end{definition}

To investigate the relation of the infinitesimal deformation of $\left( M_{s(u)}, \mathcal{N}_{\mathcal{F}_{s(u)}}^* \right)$ and that of $\left( M_t, \mathcal{N}_{\mathcal{F}_t}^* \right)$, we assume that $0\in B, 0\in D$ and $s(0)=0$. Taking a sufficiently small coordinate polydisk $\Delta$ with $0\in \Delta \subset B$, we represent $\left(\mathcal{M}_\Delta, \mathcal{N}_{\mathcal{F}_\Delta}^* \right)=\pi^{-1}\left( \Delta \right)$ in the form
\begin{align*}
\left( \mathcal{M}_\Delta , \mathcal{N}_{\mathcal{F}_\Delta}^* \right) = \bigcup \left( U_j \times \Delta, \mathcal{N}_{\mathcal{F}_\Delta}^*|_{U_j\times \Delta} \right)
\end{align*}
where $(z_j, t)\in U_j\times \Delta$ and $(z_k,t)\in U_k\times \Delta$ are the same point of $\mathcal{M}_\Delta$ if $z_j= f_{jk}(z_k,t)$, and $\Gamma\left( U_j\times \Delta, \mathcal{N}_\mathcal{F}^* \right)$ is generated by $w_j^1(z_j,t),..., w_j^q(z_j,t)$, where $w_j^\alpha(z_j,t)= \sum_{\beta=1}^n w_j^{\alpha \beta}(z_j,t) dz_j^\beta$ for some $w_j^{\alpha\beta}(z_j,t)\in \Gamma\left( U_j \times \Delta, \mathcal{O}_\mathcal{M} \right),\alpha=1,...,q$ with $w_j^\alpha(z_j,t)= \sum_{\beta=1}^q h_{jk}^{\alpha\beta}(z_k,t) w_k^\beta (z_k,t)$ for some $h_{jk}^{\alpha\beta}(z_k,t)\in \Gamma\left( \left( U_j \times \Delta\right) \bigcup \left( U_k \times \Delta \right) , \mathcal{O}_\mathcal{M} \right)$ for $\alpha=1,...,q$ and the integrabiity $w_j^1(z_j,t) \wedge \cdots \wedge w_j^q(z_j,t ) \wedge dw_j^\alpha(z_j,t)=0$ for $\alpha=1,...,q$. Take a polydisk $\Delta'$ with $0\in \Delta' \subset D$. Then $\left( \mathcal{M}\times_B D|_{\Delta'}, p^* \mathcal{N}_\mathcal{F}^*|_{\Delta'}, \Delta', \pi' \right)$ is represente by the form
\begin{align*}
\left( \mathcal{M}\times_B D|_{\Delta'}, p^* \mathcal{N}_\mathcal{F}^*|_{\Delta'} \right)=  \bigcup \left( U_j \times \Delta', p^* \mathcal{N}_\mathcal{F}^*|_{U_j\times \Delta'} \right)
\end{align*}
where $(z_j, u)\in U_j \times \Delta'$ and $(z_k, u)\in U_k \times \Delta'$ are the same point of $\mathcal{M} \times_B D|_{\Delta'}$ if $z_j= f_{jk}\left(z_k, s(u)\right)$ and $\Gamma\left( U_j \times \Delta', p^* \mathcal{N}_\mathcal{F}^* \right)$ is generated by $w_j^1\left( z_j, s(u) \right),..., w_j^q\left( z_j, u(s) \right)$, where $w_j^\alpha \left( z_j, u(s) \right)= \sum_{\beta=1}^q w_j^{\alpha\beta}\left( z_j, u(s) \right) dz_j^\beta$, and $w_j^1\left( z_j, u(s) \right) \wedge \cdots \wedge w_j^q\left( z_j u(s) \right) \wedge dw_j^\alpha\left( z_j, u(s)\right)$ for $\alpha=1,...,q$. Then we can show that 
\begin{theorem}\label{r11}
For any tangent vector $\frac{\partial }{\partial u}=c_1\frac{\partial}{\partial u_1} + \cdots + c_{m'}\frac{\partial}{\partial u_{m'}}\in T_u \left(D \right) $, the infinitesimal foliated deformation of $\left(M_{s(u)}, \mathcal{N}_{\mathcal{F}_{s(u)}}^* \right)$ along $\frac{\partial}{\partial u}$ is given by
\begin{align*}
\frac{\partial \left(M_{s(u)},  \mathcal{N}_{\mathcal{F}_{s(u)}}^*\right)}{\partial u}&=\left(\sum_{\gamma=1}^m \frac{\partial t_\gamma}{\partial u}\frac{\partial M_t}{\partial t_\gamma},\left\{w_j^\alpha\left(z_j,s(u)\right) \mapsto -\overline{\sum_{\gamma=1}^m \frac{\partial t_\gamma}{\partial u} \frac{\partial w_j^\alpha}{\partial t_\gamma}} \right\} \right)
           =\sum_{\gamma=1}^m \frac{\partial t_\gamma}{\partial u} \frac{\partial \left(M_t, \mathcal{N}_{\mathcal{F}_t}^* \right) }{\partial t_\gamma}
\end{align*}
\end{theorem}

\subsection{Theorem of completeness for deformations of foliated complex analytic structures in terms of cotangent sheaves}\

\begin{definition}
Let $\left( \mathcal{M}, \mathcal{N}_\mathcal{F}^*, B, \pi \right)$ with $\mathcal{N}_\mathcal{F}^*$ locally free be a foliated analytic family of deformations of compact foliated complex manifolds, and $t^0\in B$. Then $\left( \mathcal{M}, \mathcal{N}_\mathcal{F}^*, B, \pi \right)$ is called complete at $t^0\in B$ if for any foliated complex analytic family $\left( \mathcal{M}', \mathcal{N}_{\mathcal{F}'}^*, D, \pi'  \right)$ in terms of cotangent sheaves such that $D$ is a domain of $\mathbb{C}^{m'}$ containing $0$ and tha $\pi'^{-1}(0)=\pi^{-1}(t^0)=\left( M, \mathcal{N}_{\mathcal{F}_0}^* \right)$, there are a sufficiently small domain $\Delta$ with $0\in \Delta \subset D$, and a holomorphic map $s:u\to t=s(u)$ with $s(0)=t^0$ such that $\left( \mathcal{M}'_\Delta, \mathcal{N}_{\mathcal{F}'_\Delta}^*, \Delta, \pi' \right)$ is the foliated complex analytic family induced from $\left( \mathcal{M}, \mathcal{N}_\mathcal{F}^*, B, \pi \right)$ by $s$ where $\left( \mathcal{M}'_\Delta, \mathcal{N}_{\mathcal{F}'_\Delta}^*, \Delta, \pi' \right)$ is the restriction of $\left( \mathcal{M}', \mathcal{N}_{\mathcal{F}'}^*, D, \pi' \right)$ to $\Delta$.
\end{definition}

\begin{theorem}\label{nc1}
Let $\left(\mathcal{M}, \mathcal{N}_\mathcal{F}^*, B, \pi \right)$ be a foliated analytic family of deformations of a compact foliated complex manifold $\left(M , \mathcal{N}_{\mathcal{F}_0}^* \right)=\omega^{-1}(0)$ such that $\mathcal{N}_{\mathcal{F}_0}^*$ is locally free, and $B$ is a domain of $\mathbb{C}^r$ containing $0$. If the foliated Kodaira-Spence map $\varphi_0: T_0(B)\to \mathbb{H}^1\left(M, \mathcal{N}_{\mathcal{F}_0}^{*\bullet} \right)$ is surjective, the foliated analytic family $(\mathcal{M}, \mathcal{N}_\mathcal{F}^*, B, \pi)$ in terms of cotangent sheaves is complete at $0 \in B$.
\end{theorem}

\begin{proof}
Let $\left(\mathcal{M}, \mathcal{N}_\mathcal{F}^* , B, \pi \right)$ be a foliated analytic family in terms of cotangent sheaves which is represented as in Remark \ref{n5}.  We keep the notations in Remark \ref{n5}. Since the problem is local with respect to $B$, we may assume that $B=\left\{ t\in \mathbb{C}^m ||t|<1 \right\}$ is a polydisk, and $\mathcal{M}$ is written in the form
\begin{align}
\mathcal{M}=\bigcup_j \,\,\, \mathcal{U}_j,\,\,\,\,\,\,\,\,\mathcal{U}_j=\left\{\left( z_j, t\right) \in \mathbb{C}^n \times B\,\,\, |\,\,\, |z_j| <1 \right\}
\end{align}

Set $U_i:=M \cap \mathcal{U}_i=\left\{ z_i\in \mathbb{C}^n | |z_i|<1 \right\}$. Then $\left( M, \mathcal{N}_{\mathcal{F}_0}^* \right)=\pi^{-1}(0)$ is described in terms of covering $\mathcal{U}:=\{U_i\}=\left\{ M \cap \mathcal{U}_i \right\}$ in the following way: (1) local coordinates $z_i= \left(z_i^1,..., z_i^n \right)$ on $U_i$ with $z_i= f_{0ij}(z_j)= \left( f_{0ij}^1(z_j),..., f_{0ij}^n(z_j) \right)$ where $f_{0ij}(z_j):= f_{ij}(z_j,0)$ and (2) $\Gamma\left(U_i, \mathcal{N}_{\mathcal{F}_0}^* \right)$ is generated by $w_{0i}^\alpha(z_i):= \sum_{\beta=1}^n w_i^{\alpha\beta}(z_i,0)dz_i^\beta \in \Gamma\left(U_i, \mathcal{N}_{\mathcal{F}_0}^* \right)$ with $w_{0i}^\alpha= \sum_{\beta=1}^q h_{0ij}^{\alpha\beta}(z_j) w_{0j}^\beta$ where $h_{0ij}^{\alpha\beta}(z_j):=h_{ij}^{\alpha\beta}(z_j,0)$, and $w_{0i}^1\wedge \cdots \wedge w_{0i}^q \wedge dw_{0i}^\alpha=0$.

Let $(\mathcal{M}', \mathcal{N}_{\mathcal{F}'}^*, D, \pi')$ be an another foliated analytic family such that $\pi'^{-1}(0')= \left(M, \mathcal{N}_{\mathcal{F}_0}^* \right)$. We may assume the following:
\begin{enumerate}
\item $D\subset \left\{  u\in \mathbb{C}^{m'}| |u|<1 \right\}$ is a sufficiently small polydisk in $\mathbb{C}^{m'}$ with a system of coordinates $u=(u_1,..., u_{m'})$ centered at $0'$. \label{sc4}
\item $\mathcal{M}'$ is covered by a finite number of coordinate neighborhood $\mathcal{U}_i'=\left\{ \left(\xi_i, u\right) \in \mathbb{C}^n \times D\,\,\,|\,\,\,|\xi_i|<1\right\}$ with a system of coordinate $(\xi_i, u)$ such that $\pi'(\xi_i, u)=u$.
\item $(\xi_i, u)$ coincides with $(\xi_j, u)$ if and only if $\xi_i = f_{ij}'(\xi_i, u)$.
\item $\mathcal{N}_{\mathcal{F}'}^*$ on $\mathcal{U}_i'$ is generated by
\begin{align*}
w_i'^\alpha(\xi_i,u):= \sum_{\beta=1}^n w_i'^{\alpha\beta}(\xi_i, u) d\xi_i^\beta,\,\,\,\,\,\,\,\,\alpha=1,...,q
\end{align*}
\item on $\mathcal{U}_i' \cap \mathcal{U}_j'\ne \emptyset$, we have $w_i'^\alpha(\xi_i,u)= \sum_{\beta=1}^q h_{ij}'(\xi_j, u) w_j'^\beta(\xi_j, u)$ for $h_{ij}'^{\alpha\beta}(\xi_j, u)\in \Gamma\left( \mathcal{U}_i'\cap \mathcal{U}_j', \mathcal{O}_\mathcal{M} \right)$. \label{cs2}
\item $w_i'^1(\xi_i,u) \wedge \cdots \wedge w_i'^q(\xi_i, u) \wedge d w_i'^\alpha(\xi_i, u)=0,\alpha=1,...,q$.
\item $\pi'^{-1}(0)\cap \mathcal{U}_i'=U_i$, and $\xi_i= z_i$ on $U_i$ and $f_{0ij}(z_j)=f_{0ij}'(\xi_j)$ where $f_{0ij}'(\xi_j) : =f_{ij}'(\xi_j,0)$. \label{cs3}
\item Setting $w_{0i}'^\alpha(z_i):= w_i'^\alpha(z_i,0)$, i.e.
\begin{align*}
w_{0i}'^\alpha(z_i):= w_i'^\alpha(z_i,0)=\sum_{\beta=1}^n w_i'^{\alpha\beta}(z_i,0)dz_i^\beta,\,\,\,\,\,\alpha=1,...,q,
\end{align*}
we can find $c_{0i}^{\alpha\beta}(z_i)\in \Gamma\left( U_i, \mathcal{O}_M\right),\alpha,\beta=1,...,q$ such that
\begin{align*}
w_{0i}^\alpha(z_i)= \sum_{\beta=1}^q c_{0i}^{\alpha\beta}(z_i) w_{0i}'^\beta (z_i)
\end{align*}
\item we set $h_{0ij}'^{\alpha\beta}(\xi_j):= h_{ij}'^{\alpha\beta}(\xi_j,0)$. \label{sc5}
\end{enumerate}
We note that since $w_{0i}^\alpha=\sum_{\beta=1}^q h_{0ij}^{\alpha\beta}w_{0j}^\beta$ and $w_{0i}'^\alpha=\sum_{\beta=1}^q h_{0jk}'^{\alpha\beta} w_{0k}'^\beta$, we have
\begin{align*}
 \sum_{\beta,\gamma=1}^q c_{0i}^{\alpha\beta} h_{0ij}'^{\beta\gamma}  w_{0j}'^\gamma=\sum_{\beta=1}^q c_{0i}^{\alpha\beta} w_{0i}'^\beta = w_{0i}^\alpha=\sum_{\beta=1}^q h_{0ij}^{\alpha\beta} w_{0j}^\beta = \sum_{\beta,\gamma=1}^q h_{0ij}^{\alpha\beta} c_{0j}^{\beta\gamma} w_{0j}'^\gamma 
\end{align*}
so that we have
\begin{align*}
 \sum_{\beta=1}^q c_{0i}^{\alpha\beta} h_{0ij}'^{\beta\gamma} =  \sum_{\beta=1}^q h_{0ij}^{\alpha\beta} c_{0j}^{\beta\gamma}
\end{align*}

In order to prove Theorem \ref{nc1} it suffices to construct holomorphic functions
\begin{align*}
\varphi_i&:\mathcal{U}_i'\to \mathbb{C}^n \\
s &:D \to \mathbb{C}^{m}
\end{align*}
and a matrix function
\begin{center}
$\left[\begin{matrix}
c_{i}^{11}(\xi_i,u) & c_{i}^{12}(\xi_i,u) & \cdots & c_{i}^{1q}(\xi_i, u)\\
c_{i}^{21}(\xi_i, u) & c_{i}^{22}(\xi_i, u) & \cdots & c_{i}^{2q}(\xi_i, u)\\
\cdot & \cdot & \cdots & \cdot\\
\cdot & \cdot & \cdots & \cdot\\
\cdot & \cdot & \cdots & \cdot\\
c_{i}^{q1}(\xi_i, u) & c_{i}^{q2}(\xi_i, u) & \cdots & c_{i}^{qq}(\xi_i, u)\\
\end{matrix}\right]$
\end{center}
where $c_i^{\alpha\beta}:\mathcal{U}_i'\to \mathbb{C}$ such that
\begin{align}
\varphi_i(\xi_i, 0)=\xi_i, \,\,\,\,\,s(0)&=0,\,\,\,\,\, c_i^{\alpha\beta}(\xi_i,0)= c_{0i}^{\alpha\beta}(\xi_i) \label{nc211}\\
\varphi_i \left(f_{ij}', u \right) &= f_{ij}\left(\varphi_j, s(u)\right) \label{nc2}
\end{align}
\begin{align}
\sum_{\beta=1}^q c_i^{\alpha\beta}\left(f_{ij}' (\xi_j, u), u\right) h_{ij}'^{\beta\gamma}(\xi_j , u) = \sum_{\beta=1}^q h_{ij}^{\alpha\beta} \left (\varphi_j(\xi_j, u), s(u) \right) c_j^{\beta\gamma}(\xi_j, u)\label{nc3}
\end{align}
By setting
\begin{align*}
\left(\varphi_i, s \right): \mathcal{U}_i' \to \mathbb{C}^n \times\mathbb{C}^m,
\end{align*}
we have
\begin{align*}
&\left(\varphi_i,s \right)^* \left(w_i^\alpha\right)=\left(\varphi_i,s\right)^*\left(\sum_{\beta=1}^q w_i^{\alpha\beta}(z_i,t)dz_i^\beta \right) \\
 &= \sum_{\beta=1}^q w_i^{\alpha\beta}\left(\varphi_i(\xi_i,u),s(u) \right)d \varphi_i^\beta(\xi_i,u)= \sum_{\beta=1}^q \sum_{\gamma=1}^n w_i^{\alpha \beta}\left(\varphi_i(\xi_i, u), s(u)\right)\frac{\partial \varphi_i^\beta}{\partial \xi_i^\gamma}d\xi_i^\gamma,
\end{align*}
so that the following is satisfied: for $\alpha=1,...,q$
\begin{align}
\sum_{\beta=1}^q\sum_{\gamma=1}^n c_i^{\alpha\beta}(\xi_i, u) w_i'^{\beta\gamma}(\xi_i, u)d\xi_i^\gamma &= \sum_{\beta=1}^n \sum_{\gamma=1}^n w_i^{\alpha\beta}\left(\varphi_i(\xi_i, u), s(u) \right)\frac{\partial \varphi_i^\beta}{\partial \xi_i^\gamma}  d\xi_i^\gamma \label{nc4}
\end{align}

First we prove the existence of formal solution of $(\ref{nc2})-(\ref{nc4})$. We recall Notation \ref{te27}. Then $(\ref{nc2})-(\ref{nc4})$ are equivalent to the following systems of congruences
\begin{align}
\varphi_i^\mu \left(f_{ij}'(\xi_j, u), u \right)& \equiv_\mu f_{ij}\left(\varphi_j^\mu(\xi_j, u), s^\mu(u)\right) \label{nc5}\\
\sum_{\beta=1}^q c_i^{\alpha\beta\mu}\left(f_{ij}'(\xi_j, u),u\right) h_{ij}'^{\beta\gamma}& (\xi_j, u)  \equiv_\mu \sum_{\beta=1}^q  h_{ij}^{\alpha\beta}\left(\varphi_j^\mu(\xi_j, u), s^\mu(u) \right) c_j^{\beta\gamma\mu}(\xi_j, u) \label{nc6} \\
\sum_{\beta=1}^q\sum_{\gamma=1}^n c_i^{\alpha\beta \mu}(\xi_i, u) w_i'^{\beta\gamma}(\xi_i, u)d\xi_i^\gamma &= \sum_{\beta=1}^n \sum_{\gamma=1}^n w_i^{\alpha\beta}\left(\varphi_i^\mu (\xi_i, u), s^\mu (u) \right)\frac{\partial \varphi_i^{ \beta \mu}}{\partial \xi_i^\gamma}  d\xi_i^\gamma \label{nc7}
\end{align}
for each $\mu=1,2,...$.  We shall construct $\varphi_i^\mu, s^\mu$, and $c_i^{\alpha\beta \mu}$ satisfying $(\ref{nc5})_\mu- (\ref{nc7})_\mu$ by induction on $\mu$. We assume that $\varphi_i^{\mu-1}, s^{\mu-1}$ and $c_i^{\alpha\beta(\mu-1)}$ are already determined. Then we define homogenous polynomials $\Gamma_{ij|\mu}^\alpha, A_{ij|\mu}^{\alpha\beta}$, and $W_{i|\mu}^{\alpha \gamma}$ of degree $\mu$ by the following congruences:
{\small{\begin{align}
\Gamma_{ij|\mu}^\alpha &\equiv_\mu \varphi_i^{\alpha(\mu-1)}\left(f_{ij}'(\xi_j, u), u \right)- f_{ij}^{\alpha}\left( \varphi_j^{\mu-1}(\xi_j, u), s^{\mu-1}(u)\right) \label{nc8} \\
\sum_{\beta,\sigma=1}^q A_{ij|\mu}^{\alpha\beta} c_{0i}^{\beta\sigma} h_{0ij}'^{\sigma \gamma} &\equiv_\mu \sum_{\beta=1}^q c_i^{\alpha\beta (\mu-1)}\left(f_{ij}'(\xi_j, u), u\right) h_{ij}'^{\beta\gamma}(\xi_j, u) - \sum_{\beta=1}^q h_{ij}^{\alpha\beta} \left(\varphi_j^{\mu-1}(\xi_j, u\right), s^{\mu-1}(u))c_j^{\beta\gamma(\mu-1)}(\xi_j, u)\label{nc9}\\
W_{i|\mu}^{\alpha\gamma}& \equiv_{\mu} \sum_{\beta=1}^q  c_i^{\alpha\beta (\mu-1)}(\xi_i, u) w_i'^{\beta\gamma}(\xi_i, u)- \sum_{\beta=1}^q w_i^{\alpha\beta}\left(\varphi_i^{\mu-1}(\xi_i, u), s^{\mu-1}(u) \right) \frac{\partial \varphi_i^{\beta (\mu-1)}}{\partial \xi_i^\gamma} \label{nc10}
\end{align}}}
We set 
\begin{align*}
\Gamma_{ij|\mu}&= \sum_{\alpha=1}^n \Gamma_{ij|\mu}^{\alpha}\frac{\partial}{\partial z_i^\gamma},\,\,\,\,\,\,\,\,\,\,\,\, W_{i|\mu}^\alpha = \sum_{\gamma=1}^n W_{i|\mu}^{\alpha\gamma} dz_i^\gamma\\
\end{align*}

\begin{lemma}\label{nc11}
We have the following equalities$:$ for $\alpha=1,...,q$,
\begin{align}
&\,\,\,\,\,\,\,\,\Gamma_{jk|\mu}- \Gamma_{ik|\mu} + \Gamma_{ij|\mu}=0\label{nc12}\\
W_{i|\mu}^\alpha &=  \sum_{\beta=1}^q A_{ij|\mu}^\beta w_{0i}^\beta + \sum_{\beta=1}^q h_{0ij}^{\alpha\beta} W_{j|\mu}^\beta- \mathcal{L}_{\Gamma_{ij|\mu}}\left(w_{0i}^\alpha \right) \label{nc13}\\
\sum_{\beta=1}^q w_{0i}^1 \wedge \cdots \wedge &\overbrace{\left(-  W_{i|\mu}^\beta \right)}^{\beta-\textnormal{th}} \wedge \cdots \wedge w_{0i}^q \wedge d w_{0i}^\alpha + w_{0i}^1 \wedge \cdots \wedge w_{0i}^q \wedge d \left( - W_{i|\mu}^\alpha \right)=0 \label{nc14}
\end{align}

\end{lemma}

\begin{proof}
$(\ref{nc12})$ follows from \cite{Kod05} p.292-p.294. We prove $(\ref{nc13})$. In fact, from $(\ref{nc8})-(\ref{nc10})$,
{\small{\begin{align*}
&W_{i|\mu}^\alpha \equiv_\mu \sum_{\gamma=1}^n \left(\sum_{\beta=1}^q c_i^{\alpha\beta (\mu-1)} w_i'^{\beta\gamma}   -\sum_{\beta=1}^n w_i^{\alpha\beta}\left( \varphi_i^{\mu-1}, s^{\mu-1}\right)\frac{\partial \varphi_i^{\beta(\mu-1)}}{\partial \xi_i^\gamma}  \right) dz_i^\gamma \\
                          &\equiv_\mu \sum_{\gamma,\tau=1}^n \left(\sum_{\beta=1}^q c_i^{\alpha\beta (\mu-1)} w_i'^{\beta\gamma}   -\sum_{\beta=1}^n w_i^{\alpha\beta}\left( \varphi_i^{\mu-1}, s^{\mu-1} \right)\frac{\partial \varphi_i^{\beta(\mu-1)}}{\partial \xi_i^\gamma}  \right)\frac{\partial f_{0ij}^\gamma}{\partial z_j^\tau}dz_j^\tau \\
                          &\equiv_\mu \sum_{\gamma,\tau=1}^n \left(\sum_{\beta=1}^q c_i^{\alpha\beta (\mu-1)} w_i'^{\beta\gamma}   -\sum_{\beta=1}^n w_i^{\alpha\beta}\left( \varphi_i^{\mu-1}, s^{\mu-1}\right)\frac{\partial \varphi_i^{\beta(\mu-1)}}{\partial \xi_i^\gamma}  \right)\frac{\partial f_{ij}'^\gamma(\xi_j, u)}{\partial \xi_j^\tau}dz_j^\tau\\
                          &\equiv_\mu \sum_{\gamma, \tau=1}^n\sum_{\beta=1}^q c_i^{\alpha\beta(\mu-1)} w_i'^{\beta\gamma}\left(f_{ij}', u\right)\frac{\partial f_{ij}'^\gamma}{\partial \xi_j^\tau}       dz_j^\tau - \sum_{\tau=1}^n \sum_{\beta=1}^n w_i^{\alpha\beta}\left( \varphi_i^{\mu-1}\left(f_{ij}'(\xi_j, u), u\right), s^{\mu-1}\right) \frac{\partial \varphi_i^{\beta(\mu-1)}\left(f_{ij}', u\right)}{\partial \xi_j^\tau} dz_j^\tau\\
                          &\equiv_\mu \sum_{\tau=1}^n \sum_{\beta, \gamma=1}^q c_i^{\alpha\beta(\mu-1)} h_{ij}'^{\beta\gamma} w_j'^{\gamma\tau} dz_j^\tau - \sum_{\tau=1}^n \sum_{\beta=1}^n w_i^{\alpha\beta}\left(\Gamma_{ij|\mu}+ f_{ij}\left(\varphi_j^{\mu-1},s^{\mu-1}\right), s^{\mu-1}\right) \frac{\partial \left(\Gamma_{ij|\mu}^\beta+ f_{ij}^\beta\left(\varphi_j^{\mu-1}, s^{\mu-1}\right) \right)}{\partial \xi_j^\tau} dz_j^\tau \\
                          &\equiv_\mu\sum_{\tau=1}^n \left(\sum_{\beta, \sigma,\gamma=1}^q A_{ij|\mu}^{\alpha\beta} h_{0ij}^{\beta\sigma} c_{0j}^{\sigma \gamma} + \sum_{\beta=1}^q h_{ij}^{\alpha\beta}\left( \varphi_j^{\mu-1}, s^{\mu-1} \right) c_j^{\beta\gamma (\mu-1)} \right) w_j'^{\gamma \tau} dz_j^\tau \\
                          &- \sum_{\tau=1}^n\sum_{\beta=1}^q w_i^{\alpha\beta}\left(f_{ij}\left( \varphi_j^{\mu-1}, s^{\mu-1} \right), s^{\mu-1}\right)\frac{\partial f_{ij}^\beta \left( \varphi_j^{\mu-1}, s^{\mu-1} \right)}{\partial \xi_j^\tau} dz_j^\tau - \sum_{\beta,\tau=1}^n \frac{\partial w_{0i}^{\alpha\beta}}{\partial z_i^\tau}\Gamma_{ij|\mu}^\tau dz_i^\beta - \sum_{\tau=1}^n \sum_{\beta=1}^qw_{0i}^{\alpha\beta}\frac{\partial \Gamma_{ij|\mu}^\beta}{\partial z_j^\tau}dz_j^\tau\\
                          &\equiv_\mu \sum_{\beta=1}^q A_{ij|\mu}^\beta w_{0i}^\beta + \sum_{\tau=1}^n \sum_{\beta=1}^q h_{ij}^{\alpha\beta} \left( \varphi_j^{\mu-1}, s^{\mu-1} \right)\left( W_{j|\mu}^{\beta\tau} + \sum_{\gamma=1}^n w_j^{\beta\gamma}\left( \varphi_j^{\mu-1}, s^{\mu-1}\right) \frac{\partial \varphi_j^{\gamma(\mu-1)}}{\partial \xi_j^\tau} \right) dz_j^\tau\\
                          &- \sum_{\tau=1}^n\sum_{\beta=1}^q w_i^{\alpha\beta}\left( f_{ij}\left( \varphi_j^{\mu-1}, s^{\mu-1}\right), s^{\mu-1} \right)\frac{\partial f_{ij}^\beta \left( \varphi_j^{\mu-1}, s^{\mu-1}\right)}{\partial \xi_j^\tau} dz_j^\tau - \sum_{\beta,\tau=1}^n \frac{\partial w_{0i}^{\alpha\beta}}{\partial z_i^\tau}\Gamma_{ij|\mu}^\tau dz_i^\beta - \sum_{\tau=1}^n \sum_{\beta=1}^qw_{0i}^{\alpha\beta}\frac{\partial \Gamma_{ij|\mu}^\beta}{\partial z_j^\tau}dz_j^\tau\\
                          &\equiv_\mu \sum_{\beta=1}^q A_{ij|\mu}^\beta w_{0i}^\beta + \sum_{\beta=1}^q h_{0ij}^{\alpha\beta} W_{j|\mu}^\beta- \mathcal{L}_{\sum_{\tau=1}^n \Gamma_{ij|\mu}^\tau \frac{\partial}{\partial z_i^\tau}}\left(\sum_{\beta=1}^n w_{0i}^{\alpha\beta} dz_i^\beta \right)\\
                          &+\sum_{\tau,\gamma =1}^n \sum_{\beta=1}^q h_{ij}^{\alpha\beta}\left( \varphi_j^{\mu-1}, s^{\mu-1}\right)w_j^{\beta\gamma}\left( \varphi_j^{\mu-1}, s^{\mu-1} \right)\frac{\partial \varphi_j^{\gamma(\mu-1)}}{\partial \xi_j^\tau} dz_j^\tau - \sum_{\tau, \gamma=1}^n \sum_{\beta=1}^q w_i^{\alpha\beta}\left(f_{ij}(\varphi_j^{\mu-1}, s^{\mu-1}), s^{\mu-1}\right)\frac{\partial f_{ij}^\beta}{\partial z_j^\gamma}\left(\varphi_j^{\mu-1}, s^{\mu-1}\right) \frac{\partial \varphi_j^{\gamma(\mu-1)}}{\partial \xi_j^\tau} dz_j^\tau \\
                          &\equiv_\mu  \sum_{\beta=1}^q A_{ij|\mu}^\beta w_{0i}^\beta + \sum_{\beta=1}^q h_{0ij}^{\alpha\beta} W_{j|\mu}^\beta- \mathcal{L}_{\sum_{\tau=1}^n \Gamma_{ij|\mu}^\tau \frac{\partial}{\partial z_i^\tau}}\left(\sum_{\beta=1}^n w_{0i}^{\alpha\beta} dz_i^\beta \right) 
\end{align*}}}
We prove $(\ref{nc14})$. In fact, from $w_i^1\wedge \cdots \wedge w_i^q \wedge d w_i^\alpha =0$, i.e.
\begin{align}\label{nc15}
&\left( \sum_{\beta_1=1}^n w_i^{1\beta_1}(z_i,t) dz_i^{\beta_1}   \right) \wedge \cdots \wedge \left(\sum_{\beta_q=1}^n w_i^{q\beta_q}(z_i,t) dz_i^{\beta_q}   \right) \wedge d\left(\sum_{\beta=1}^q w_i^{\alpha\beta}(z_i,t) dz_i^\beta   \right) =0
\end{align} 
By pulling back $(\ref{nc15})$ under $\left( \varphi_i^{\mu-1}, s^{\mu-1}\right)$, we have
{\Small{\begin{align*}
&\left( \sum_{\beta_1=1}^n w_i^{1\beta_1}\left( \varphi_i^{\mu-1},s^{\mu-1} \right) d \varphi_i^{\beta_1(\mu-1)}  \right) \wedge \cdots \wedge \left(\sum_{\beta_q=1}^n w_i^{q\beta_q}\left( \varphi^{\mu-1},s^{\mu-1} \right) d  \varphi_i^{\beta_q (\mu-1)}    \right) \wedge d\left(\sum_{\beta=1}^q w_i^{\alpha\beta}\left( \varphi_i^{\mu-1},s^{\mu-1}\right) d \varphi_i^{\beta(\mu-1)}   \right) =0 
\end{align*}}}
Then we have from $(\ref{nc10})$
\begin{align*}
&\left(\sum_{\beta,\gamma=1}^n c_i^{1\beta(\mu-1)} w_i'^{\beta \gamma}d\xi_i^\gamma   - \sum_{\gamma=1}^n W_{i|\mu}^{1\gamma} dz_i^\gamma      \right) \wedge \cdots  \wedge \left(\sum_{\beta, \gamma=1}^n c_i^{q\beta(\mu-1)} w_i'^{\beta \gamma}d\xi_i^\gamma   - \sum_{\gamma=1}^n W_{i|\mu}^{q\gamma} dz_i^\gamma     \right)  \\
& \wedge d \left(\sum_{\beta,\gamma=1}^n c_i^{\alpha\beta(\mu-1)} w_i'^{\beta \gamma}d\xi_i^\gamma   - \sum_{\gamma=1}^n W_{i|\mu}^{\alpha \gamma} dz_i^\gamma      \right)  = 0 
\end{align*}
In other words, we have
\begin{align} \label{nc16}
\left( \sum_{\beta=1}^n c_i^{1\beta(\mu-1)} w_i'^\beta - W_{i|\mu}^1 \right) \wedge \cdots \wedge \left( \sum_{\beta=1}^n c_i^{q\beta(\mu-1)} w_i'^\beta - W_{i|\mu}^1 \right)  \wedge d \left( \sum_{\beta=1}^n c_i^{\alpha\beta(\mu-1)} w_i'^\beta - W_{i|\mu}^1 \right) =0
\end{align}
Since $w_i'^1 \wedge \cdots \wedge w_i^q \wedge dw_i'^\alpha=0$, $(\ref{nc16})$ implies $(\ref{nc14})$. This completes the proof of Lemma \ref{nc11}.

\end{proof}

Our purpose is to determine $\varphi^\mu= \varphi^{\mu-1}+ \varphi_{i|\mu}, s^\mu= s^{\mu-1}+ s_\mu$, and $c_i^{\alpha\beta(\mu-1)}+ c_{i|\mu}^{\alpha\beta}$ satisfying $(\ref{nc5})_\mu-(\ref{nc7})_\mu$.

\begin{lemma}\label{nc23}
$(\ref{nc5})_\mu-(\ref{nc7})_\mu$ are equivalent to the following equalities$:$
\begin{align}
\Gamma_{ij|\mu}&= \varphi_{j|\mu}- \varphi_{i|\mu}+ \sum_{\lambda=1}^m  s_\mu^\lambda \rho_{ij|\lambda} \label{nc17} \\
\sum_{\beta, \sigma=1}^q A_{ij|\mu}^{\alpha\beta} h_{0ij}^{\beta\sigma} c_{0j}^{\sigma \gamma} &=  -\sum_{\beta=1}^q c_{i|\mu}^{\alpha\beta} h_{0ij}'^{\beta\gamma} + \sum_{\beta=1}^q \left[ \varphi_{j|\mu} , h_{0ij}^{\alpha\beta}  \right]  c_{0j}^{\beta\gamma} + \sum_{\beta=1}^q \sum_{\lambda=1}^m \frac{\partial h_{ij}^{\alpha\beta}}{\partial t_\lambda }|_{t=0} s_\mu^\lambda c_{0j}^{\beta \gamma} + \sum_{\beta=1}^q h_{0ij}^{\beta\gamma} c_{j|\mu}^{\beta\gamma}    \label{nc18}   \\
W_{i|\mu}^{\alpha } &=  - \sum_{\beta=1}^q c_{i|\mu}^{\alpha\beta} w_{0i}'^{\beta} + \mathcal{L}_{\varphi_{i|\mu}}\left( w_{0i}^\alpha  \right) - \sum_{\lambda=1}^r s_\mu^\lambda \beta_{i\lambda}^\alpha   \label{nc19}
\end{align}
where
\begin{align*}
\varphi_{i|\mu}=\sum_{\alpha=1}^n \varphi_{i|\mu}^\alpha \frac{\partial}{\partial z_i^\alpha},\,\,\,\,\,\,\,\rho_{ij\lambda}=\sum_{\alpha=1}^n \frac{\partial f_{ij}^\alpha}{\partial t_\lambda}|_{t=0} \frac{\partial}{\partial z_i^\alpha},\,\,\,\,\,\,\,\,\,\,\,\, \beta_{i\lambda}^\alpha = - \frac{\partial w_i^\alpha (z_i,t) }{\partial t_\lambda}|_{t=0}
\end{align*}

\end{lemma}

\begin{proof}
$(\ref{nc17})$ follows from \cite{Kod05} p.290. We prove $(\ref{nc18})$. In fact,
\begin{align*}
&\sum_{\beta=1}^q \left( c_i^{\alpha\beta (\mu-1)}\left(f_{ij}', u \right) +  c_{i|\mu}^{\alpha\beta}\left(f_{ij}', u \right)        \right) h_{ij}'^{\beta\gamma} - \sum_{\beta=1}^q h_{ij}^{\alpha\beta} \left( \varphi_j^{\mu-1} + \varphi_{j|\mu}, s^{\mu-1}+ s_\mu \right) \left( c_j^{\beta\gamma (\mu-1)} + c_{j|\mu}^{\beta\gamma }   \right)\\
&\equiv_\mu \sum_{\beta,\sigma=1}^q A_{ij|\mu}^{\alpha\beta} h_{0ij}^{\beta\sigma} c_{0j}^{\sigma \gamma} + \sum_{\beta=1}^q c_{i|\mu}^{\alpha\beta} h_{0ij}'^{\beta\gamma} - \sum_{\beta=1}^q\sum_{\eta=1}^n \frac{\partial h_{ij}^{\alpha\beta}}{\partial z_j^\eta} \varphi_{j|\mu}^\eta c_{0j}^{\beta\gamma} - \sum_{\beta=1}^q \sum_{\lambda=1}^m \frac{\partial h_{ij}^{\alpha\beta}}{\partial t_\lambda } s_\mu^\lambda c_{0j}^{\beta \gamma} - \sum_{\beta=1}^q h_{0ij}^{\beta\gamma} c_{j|\mu}^{\beta\gamma}
\end{align*}
We prove $(\ref{nc19})$. In fact, first we note that
\begin{align*}
&\sum_{\beta=1}^q \left( c_i^{\alpha\beta(\mu-1)} + c_{i|\mu}^{\alpha\beta} \right) w_i'^{\beta\gamma}- \sum_{\beta=1}^n w_i^{\alpha\beta}\left( \varphi_i^{\mu-1}+ \varphi_{i|\mu}, s^{\mu-1}+ s_\mu   \right) \frac{\partial \left( \varphi_i^{\beta(\mu-1)} + \varphi_{i|\mu}^\beta \right) }{\partial \xi_i^\gamma}\\
&= W_{i|\mu}^{\alpha\gamma}+ \sum_{\beta=1}^q c_{i|\mu}^{\alpha\beta} w_{0i}'^{\beta\gamma} -  \sum_{\eta=1}^n \frac{\partial w_{0i}^{\alpha\gamma }}{\partial z_i^\eta} \varphi_{i|\mu}^\eta -  \sum_{\lambda=1}^m \frac{\partial w_{i}^{\alpha\gamma }}{\partial t_\lambda} s_\mu^\lambda - \sum_{\beta=1}^n w_{0i}^{\alpha\beta} \frac{\partial \varphi_{i|\mu}^\beta}{\partial z_i^\gamma}
\end{align*}
Then by adding $\sum_{\gamma=1}^n dz_i^\gamma$, we have
\begin{align*}
W_{i|\mu}^\alpha&=  - \sum_{\gamma=1}^n\sum_{\beta=1}^q c_{i|\mu}^{\alpha\beta} w_{0i}'^{\beta\gamma}dz_i^\gamma +  \sum_{\gamma=1}^n \sum_{\eta=1}^n \frac{\partial w_{0i}^{\alpha\gamma }}{\partial z_i^\eta} \varphi_{i|\mu}^\eta dz_i^\gamma +  \sum_{\gamma=1}^n\sum_{\lambda=1}^r s_\mu^\lambda \frac{\partial w_{i}^{\alpha\gamma }}{\partial t_\lambda} dz_i^\gamma + \sum_{\gamma=1}^n \sum_{\beta=1}^n w_{0i}^{\alpha\beta} \frac{\partial \varphi_{i|\mu}^\beta}{\partial z_i^\gamma}dz_i^\gamma\\
 &=  -  \sum_{\beta=1}^q c_{i|\mu}^{\alpha\beta} w_{0i}'^{\beta} + \mathcal{L}_{\sum_{\eta=1}^n \varphi_{i|\mu}^\eta\frac{\partial}{\partial z_i^\eta}}\left(\sum_{\gamma=1}^n w_{0i}^{\alpha \gamma} dz_i^\gamma  \right) +  \sum_{\gamma=1}^n\sum_{\lambda=1}^r s_\mu^\lambda \frac{\partial w_{i}^{\alpha\gamma }}{\partial t_\lambda} dz_i^\gamma 
\end{align*}
This completes the proof of Lemma \ref{nc23}.
\end{proof}

We define an element $\overline{W}_{i|\mu}\in \Gamma\left( U_i, \mathscr{H}om_{\mathcal{O}_M} \left( \mathcal{N}_{\mathcal{F}_0}^*, \frac{\Omega_M^1}{\mathcal{N}_{\mathcal{F}_0}^*}  \right) \right)$ by
\begin{align}
\overline{W}_{i|\mu}: \Gamma\left( U_i, \mathcal{N}_{\mathcal{F}_0}^*  \right) &\to \Gamma\left(  U_i  ,   \frac{\Omega_M^1}{\mathcal{N}_{\mathcal{F}_0}^*}   \right) \label{nc125} \\
                                       w_{0i}^\alpha &\mapsto \overline{W_{i|\mu}^\alpha } \notag
\end{align}
and linearly extends to $\Gamma\left( U_i, \mathcal{N}_{\mathcal{F}_0}^* \right)$. Then $(\ref{nc12})-(\ref{nc14})$ implies that
\begin{align}\label{uuc15}
\left(  \left\{ -  \overline{W}_{i|\mu} \right\}   ,  \left\{  \Gamma_{ij|\mu} \right\} \right) \in C^0\left( \mathcal{U}, \mathscr{H}om_{\mathcal{O}_M}\left( \mathcal{N}_{\mathcal{F}_0}^* ,  \frac{\Omega_M^1}{\mathcal{N}_{\mathcal{F}_0}^*} \right)   \right) \bigoplus C^1 \left( \mathcal{U}, \Theta_M  \right)
\end{align}
defines a $1$-cocycle in the following \v Cech resolution of $\mathcal{N}_\mathcal{F}^{* \bullet}$:
{\small{\begin{center}
$\begin{CD}
\cdots \\
@AAA \\
C^0\left( \mathcal{U}, \mathscr{H}om_{\mathcal{O}_M}\left( \mathcal{N}_\mathcal{F}^*, \tilde{\mathcal{S}}^2 \right) \right) @>-\delta >> \cdots \\
@AAA @AAA  \\
C^0\left(\mathcal{U}, \mathscr{H}om_{\mathcal{O}_M} \left(\mathcal{N}_{\mathcal{F}_0}^* , \frac{\Omega_M^1}{\mathcal{N}_{\mathcal{F}_0}^*} \right) \right) @>\delta>>  C^1\left(\mathcal{U}, \mathscr{H}om_{\mathcal{O}_M} \left( \mathcal{N}_{\mathcal{F}_0}^* , \frac{\Omega_M^1}{\mathcal{N}_{\mathcal{F}_0}^* } \right) \right) @>-\delta >> \cdots \\
@AAA @AAA @AAA \\ 
C^0 \left(\mathcal{U}, \Theta_M \right) @>-\delta>> C^1\left(\mathcal{U}, \Theta_M \right) @>\delta>> C^2 \left(\mathcal{U}, \Theta_M \right) @>-\delta>> \cdots
\end{CD}$
\end{center}}}
By the hypothesis that the foliated Kodaira-Spencer map $\varphi_0: T_0(B) \to \mathbb{H}^1\left( M, \mathcal{N}_\mathcal{F}^{ * \bullet } \right)$ is surjective, we can find homogeneous polynomial $s_\mu^\lambda$ such that
\begin{align*}
\varphi_0\left( \sum_{\lambda=1}^r s_\mu^\lambda \frac{\partial}{\partial t_\lambda}   \right) = \left(  \left\{  -\overline{W}_{i|\mu}    \right\} ,   \left\{ \Gamma_{ij|\mu}  \right\}   \right)
\end{align*}

Since we have
\begin{align*}
\varphi_0\left( \frac{\partial}{\partial t_\lambda} \right) =\left( \{\beta_{i\lambda}\} =\left\{ w_{0i}^\alpha \mapsto  \beta_{i\lambda}^\alpha:= \overline{ -\frac{\partial w_i^\alpha (z_i,t)}{\partial t_\lambda}|_{t=0} }\right\}, \{\rho_{ij\lambda} \} =  \left\{ \sum_{\alpha=1}^n \frac{\partial f_{ij}^\alpha}{\partial t_\lambda}|_{t=0} \frac{\partial}{\partial z_i^\alpha } \right\}   \right),
\end{align*}
there exists $\left\{ \varphi_{i|\mu} \right\}\in C^0 \left(\mathcal{U}, \Theta_M \right)$ and $\left\{ c_{i|\mu} \right\} \in C^0\left( \mathcal{U}, \mathscr{H}om_{\mathcal{O}_M}\left( \mathcal{N}_{\mathcal{F}_0}^*, \mathcal{N}_{\mathcal{F}_0}^* \right) \right)$ where $c_{i|\mu}\in \Gamma\left( U_i, \mathscr{H}om_{\mathcal{O}_M}\left( \mathcal{N}_{\mathcal{F}_0}^*, \mathcal{N}_{\mathcal{F}_0}^* \right) \right)$ defined by $c_{i|\mu}: \Gamma\left( U_i , \mathcal{N}_{\mathcal{F}_0}^* \right) \to \Gamma\left( U_i, \mathcal{N}_{\mathcal{F}_0}^* \right), w_{0i}^\alpha \mapsto \sum_{\gamma=1}^q c_{i|\mu}^{\alpha \gamma} w_{0i}'^\gamma$ such that
\begin{align*} 
\varphi_{i|\mu}-  \varphi_{j|\mu}&  = -\Gamma_{ij|\mu} + \sum_{\lambda=1}^r s_\mu^\lambda \rho_{ij\lambda}\\
\mathcal{L}_{\varphi_{i|\mu}}\left( w_{0i}^\alpha \right) &=  W_{i|\mu}^\alpha + \sum_{\lambda=1}^r s_\mu^\lambda \alpha_{i\lambda}^\alpha + \sum_{\gamma=1}^p c_{i|\mu}^{\alpha \gamma} w_{0i}'^\gamma
\end{align*}
Then $(\ref{nc17})$ and $(\ref{nc19})$ holds. We check $(\ref{nc18})$. In fact, first we note that from $(\ref{nc19})$ and $(\ref{nc13})$, we have
{\small{\begin{align}\label{nc20}
\mathcal{L}_{ \varphi_i|\mu} \left( w_{0i}^\alpha \right)- \sum_{\gamma=1}^q c_{i|\mu}^{\alpha \gamma} w_{0i}'^\gamma- \sum_{\lambda=1}^r s_\mu^\lambda \beta_{i\lambda}^\alpha - \sum_{\beta=1}^q h_{0ij}^{\alpha \beta}\left( \mathcal{L}_{\varphi_j|\mu}\left(w_{0j}^\beta\right)- \sum_{\gamma=1}^q c_{j|\mu}^{\beta\gamma} w_{0j}'^\gamma   - \sum_{\lambda=1}^r s_\mu^\lambda \beta_{j\lambda}^\beta  \right)= \sum_{\beta=1}^q A_{ij|\mu}^{\alpha \beta } w_{0i}^\beta - \mathcal{L}_{\Gamma_{ij|\mu}}(w_{0i}^\alpha)
\end{align}}}
We recall that
\begin{align*}
\mathcal{L}_{\rho_{ij\lambda}} (w_{0i}^\alpha) - \beta_{i\lambda}^\alpha= \sum_{\beta=1}^q \frac{\partial h_{ij}^{\alpha\beta}}{\partial t_\lambda}|_{t=0}  w_{0j}^\beta - \sum_{\beta=1}^q h_{0jk}^{\alpha\beta} \beta_{j\lambda}^\beta,
\end{align*}
so that by multiplying $\sum_{\lambda=1}^m s_\mu^\lambda$, we have
\begin{align}\label{nc21}
\mathcal{L}_{\sum_{\lambda=1}^m s_\mu^\lambda\rho_{ij\lambda}} (w_{0i}^\alpha) - \sum_{\lambda=1}^m s_\mu^\lambda  \beta_{i\lambda}^\alpha= \sum_{\lambda=1}^m \sum_{\beta=1}^q s_\mu^\lambda \frac{\partial h_{ij}^{\alpha\beta}}{\partial t_\lambda}|_{t=0}  w_{0j}^\beta -  \sum_{\lambda=1}^m \sum_{\beta=1}^q s_\mu^\lambda h_{0ij}^{\alpha\beta} \beta_{j\lambda}^\beta
\end{align}
From $(\ref{nc17})$ and $(\ref{nc21})$, we have
\begin{align}\label{nc22}
\mathcal{L}_{\Gamma_{ij|\mu}}(w_{0i}^\alpha) - \mathcal{L}_{\varphi_{j|\mu}}(w_{0i}^\alpha) + \mathcal{L}_{\varphi_{i|\mu}}(w_{0i}^\alpha) - \sum_{\lambda=1}^r s_\mu^\lambda \beta_{i\lambda}^\alpha + \sum_{\lambda=1}^r \sum_{\beta=1}^q s_\mu^\lambda h_{0ij}^{\alpha\beta} \beta_{j\lambda}^\beta = \sum_{\lambda=1}^r \sum_{\beta=1}^q s_\mu^\lambda \frac{\partial h_{ij}^{\alpha\beta}}{\partial t_\lambda}|_{t=0} w_{0j}^\beta
\end{align}
Then from $(\ref{nc20})$ and $(\ref{nc22})$, we have
\begin{align*}
 \sum_{\lambda=1}^m \sum_{\beta=1}^q s_\mu^\lambda \frac{\partial h_{ij}^{\alpha\beta}}{\partial t_\lambda}|_{t=0} w_{0j}^\beta - \sum_{\beta=1}^q c_{i|\mu}^{\alpha \beta} w_{0i}'^\beta + \sum_{\beta=1}^q \left[ \varphi_{j|\mu}, h_{0ij}^{\alpha \beta}\right] w_{0j}^\beta + \sum_{\beta, \gamma=1}^q h_{0ij}^{\alpha\beta} c_{j|\mu}^{\beta\gamma} w_{0j}'^\gamma = \sum_{\beta=1}^q A_{ij|\mu}^{\alpha \beta} w_{0i}^\beta
\end{align*}
Since $w_{0j}^\beta= \sum_{\gamma=1}^q c_{0j}^{\beta\gamma} w_{0j}'^\gamma$, we have
{\small{\begin{align*}
 \sum_{\gamma=1}^q\sum_{\lambda=1}^m \sum_{\beta=1}^q s_\mu^\lambda \frac{\partial h_{ij}^{\alpha\beta}}{\partial t_\lambda}|_{t=0} c_{0j}^{\beta\gamma} w_{0j}'^\gamma - \sum_{\gamma=1}^q \sum_{\beta=1}^q c_{i|\mu}^{\alpha \beta} h_{0ij}'^{\beta\gamma} w_{0j}'^\gamma + \sum_{\gamma=1}^q\sum_{\beta=1}^q \left[ \varphi_{j|\mu}, h_{0ij}^{\alpha \beta}\right] c_{0j}^{\beta\gamma} w_{0j}'^\gamma + \sum_{\beta, \gamma=1}^q h_{0ij}^{\alpha\beta} c_{j|\mu}^{\beta\gamma} w_{0j}'^\gamma = \sum_{\gamma=1}^q \sum_{\beta=1}^q A_{ij|\mu}^{\alpha\beta} h_{0ij}^{\beta\sigma} c_{0j}^{\sigma \gamma} w_{0j}'^\gamma
\end{align*}}}
This implies $(\ref{nc18})$.

\begin{remark} \label{uuc19}
We recall from $(\ref{uuc17})$ that we have $\mathbb{H}^1\left( M, \mathcal{N}_{\mathcal{F}_0}^{*\bullet} \right) \cong \mathbb{H}^1\left( M, \mathcal{E}_{\mathcal{N}_{\mathcal{F}_0}^*}^\bullet \right)$. We will reinterpret $(\ref{uuc15})$ in $\mathbb{H}^1\left( M, \mathcal{N}_{\mathcal{F}_0}^{*\bullet} \right)$ in terms of $\mathbb{H}^1\left( M, \mathcal{E}_{\mathcal{N}_{\mathcal{F}_0}^*}^\bullet  \right)$ in the following \v Cech resolution of $\mathcal{E}_{\mathcal{N}_{\mathcal{F}_0}^*}^\bullet$.
{\small{\begin{equation}\label{uuc18}
\begin{CD}
\cdots\\
@AE_2AA \\
C^0\left( \mathcal{U} , \mathscr{H}om_{\mathcal{O}_M}\left( \mathcal{N}_{\mathcal{F}_0}^*, \tilde{\mathcal{S}}^2 \right) \right) @>-\delta>> \cdots\\
@AE_1'AA @AE_1'AA\\
C^0\left(\mathcal{U},\mathscr{H}om_{\mathcal{O}_M}\left( \mathcal{N}_{\mathcal{F}_0}^*, \Omega_M^1 \right) \right) @>\delta>> C^1\left(\mathcal{U},\mathscr{H}om_{\mathcal{O}_M}\left( \mathcal{N}_{\mathcal{F}_0}^*, \Omega_M^1 \right) \right) @>-\delta>> \cdots \\
@AE_0' AA @AE_0'AA @AE_0'AA \\
C^0\left(\mathcal{U}, \mathcal{E}_{\mathcal{N}_{\mathcal{F}_0}^*}\right) @>-\delta>> C^1\left(\mathcal{U}, \mathcal{E}_{\mathcal{N}_{\mathcal{F}_0}^*}\right)  @>\delta>> C^2\left(\mathcal{U}, \mathcal{E}_{\mathcal{N}_{\mathcal{F}_0}^*}\right) @>-\delta>> \cdots
\end{CD}
\end{equation}}}
We set $A_{ij|\mu}\in \Gamma\left( U_{ij}, \mathscr{H}om_{\mathcal{O}_M}\left( \mathcal{N}_{\mathcal{F}_0}^*, \mathcal{N}_{\mathcal{F}_0}^* \right) \right)$ by
\begin{align}\label{ssc4}
A_{ij|\mu}:\Gamma\left( U_{ij} , \mathcal{N}_{\mathcal{F}_0}^* \right) &\to \Gamma\left( U_{ij} , \mathcal{N}_{\mathcal{F}_0}^* \right) \\
 w_{0i}^\alpha &\mapsto \sum_{\beta=1}^q A_{ij|\mu}^{\alpha\beta} w_{0i}^\beta \notag
\end{align}
and linearly extends to $\Gamma\left( U_{ij}, \mathcal{N}_{\mathcal{F}_0}^* \right)$. Then from $(\ref{nc13})$, we see that
\begin{align*}
A_{ji|\mu}\left( w_{0i}^\alpha \right) & = \sum_{\beta=1}^q h_{0ij}^{\alpha\beta} A_{ji|\mu}^{\alpha\beta} \left( w_{0j}^\beta \right) = \sum_{\beta=1}^q h_{0ij}^{\alpha\beta}\left( W_{j|\mu}^\beta - \sum_{\eta=1}^q h_{0ji}^{\beta \eta} W_{i|\mu}^\eta + \mathcal{L}_{\Gamma_{ji}|\mu}\left( w_{0j}^\beta \right)   \right) \\
&= \sum_{\beta=1}^q h_{0ij}^{\alpha\beta} W_{j|\mu}^\beta - W_{i|\mu}^\alpha - \mathcal{L}_{\Gamma_{ij|\mu}}\left( w_{0i}^\alpha \right) + \sum_{\beta=1}^q \left[ \Gamma_{ij|\mu}, h_{0ij}^{\alpha\beta} \right] w_{0j}^\beta = - A_{ij|\mu}\left( w_{0i}^\alpha \right) + \sum_{\beta=1}^q \left[ \Gamma_{ij|\mu}, h_{0ij}^{\alpha\beta} \right] w_{0j}^\beta
\end{align*}
This implies that $\left\{ \left(\Gamma_{ij|\mu}, A_{ij|\mu}  \right) \right\}\in C^1\left( \mathcal{U}, \mathcal{E}_{\mathcal{N}_{\mathcal{F}_0}^*} \right)$. On the other hand, from $(\ref{nc13})$, we see that
\begin{align*}
&\left( A_{ij|\mu} - A_{ik|\mu} + A_{jk|\mu}   \right)\left( w_{0i}^\alpha \right)\\
&= W_{i|\mu}^\alpha - \sum_{\beta=1}^q h_{0ij}^{\alpha\beta} W_{j|\mu}^\beta + \mathcal{L}_{\Gamma_{ij|\mu}} \left( w_{0i}^\alpha \right) - W_{i|\mu}^\alpha + \sum_{\beta=1}^q h_{0ik}^{\alpha\beta} W_{k|\mu}^\beta - \mathcal{L}_{\Gamma_{ik|\mu}}\left( w_{0i}^\alpha \right) + \sum_{\beta=1}^q h_{0ij}^{\alpha\beta}\left(  W_{j|\mu}^\beta - \sum_{\eta=1}^q h_{0jk}^{\beta \eta} W_{k|\mu}^\eta + \mathcal{L}_{\Gamma_{jk|\mu}}\left( w_{0j}^\beta \right)  \right)\\
& = - \sum_{\eta=1}^q \left[ \Gamma_{jk|\mu} , h_{0ij}^{\alpha\beta} \right] w_{0j}^\beta
\end{align*}
This implies that $\delta\left(\left\{\left( \Gamma_{ij|\mu}, A_{ij|\mu} \right) \right\} \right) =0$. Hence 
\begin{align}
\left( \left\{ - W_{i|\mu}\right\}, \left\{ \left( \Gamma_{ij|\mu}, A_{ij|\mu}   \right) \right\} \right) \in C^0\left( \mathcal{U}, \mathscr{H}om_{\mathcal{O}_M}\left( \mathcal{N}_{\mathcal{F}_0}^*, \Omega_M^1 \right) \right) \bigoplus C^1\left( \mathcal{U}, \mathcal{E}_{\mathcal{N}_{\mathcal{F}_0}^*} \right)
\end{align}
defines a $1$-cocycle in the above \v Cech resolution of $\mathcal{E}_{\mathcal{N}_{\mathcal{F}_0}^*}^\bullet$.
\end{remark}

\subsection{Proof of convergence}\

We prove that we can choose appropriate solutions $\varphi_{i|\mu},s_\mu$, and $c_{i|\mu}^{\alpha\beta}$ satisfying $(\ref{nc17})-(\ref{nc19})$ in each inductive step so that 
\begin{align*}
s(u)&=s_1(u)+ s_2(u)+\cdots + s_\mu(u)+ \cdots\\
\varphi_i\left( \xi_i, u\right)&= \xi_i+ \varphi_{i|1}\left(\xi_i, u\right) +  \varphi_{i|2}\left( \xi_i, u \right) + \cdots + \varphi_{i|\mu}(\xi_i, u) + \cdots \\
c_i^{\alpha\beta}(\xi_i, u)&= c_{0i}^{\alpha\beta}(\xi_i) + c_{i|1}^{\alpha\beta}(\xi_i,u) + \cdots + c_{i|\mu}^{\alpha\beta} (\xi_i, u) + \cdots
\end{align*}
converge absolutely and uniformly for $|u|< \epsilon$ if $\epsilon>0$ is sufficiently small. We recall the notations $(\ref{ncc1})$ and $(\ref{ncc2})$.

It suffices to prove the estimates $s(u)\ll A(u), \varphi_i(\xi_i,u)-\xi_i\ll A(u)$ and $c_i^{\alpha\beta}(\xi_i, u)- c_{0i}^{\alpha\beta} \ll A(u)$ for suitable constants $b$ and $c$, equivalently
\begin{align}\label{ncc3}
s^\mu(u)\ll A(u),\,\,\,\,\,\,\,\,\,\varphi_i^\mu(\xi_i, u)- \xi_i\ll A(u),\,\,\,\,\,\,\,\,\,\,\,c_i^{\alpha\beta \mu}(\xi_i, u)- c_{0i}^{\alpha\beta}(\xi_i) \ll A(u)
\end{align}
for $\mu=1,2,3,\cdots$. We will prove $(\ref{ncc3})$ by induction on $\mu=1,2,3,\cdots$. For $\mu=1$, since the linear term of $A(u)$ is $\frac{b}{16}\left( u_1+ \cdots + u_{r'} \right)$, the estimates $\left(\ref{ncc3} \right)_1$ holds if $b$ is sufficiently large. Let $\mu\geq  2$ and assume that the induction $\left(\ref{ncc3}\right)_{\mu-1}$ holds for $\mu-1$, i.e.
\begin{align*}
s^{\mu-1}(u)\ll A(u),\,\,\,\,\,\,\,\,\,\varphi_i^{\mu-1}(\xi_i, u)- \xi_i \ll A(u),\,\,\,\,\,\,\,\,\,\,\,c_i^{\alpha\beta(\mu-1)}(\xi_i, u) - c_{0i}^{\alpha\beta}(\xi_i) \ll A(u)
\end{align*}
We will prove that $(\ref{ncc3})_\mu$ holds for $\mu$. We estimate $(\ref{nc8}),(\ref{nc9})$ and $(\ref{nc10})$ in the following Lemma.

\begin{lemma}\label{ncc20}
\begin{align}
\Gamma_{ij|\mu}(\xi_i, u) &\ll \left( \frac{K_1}{b}+ \frac{K_2}{c}+ \frac{K_3 b}{c} \right)A(u) \,\,\,\,\,\,\,\,\textnormal{on}\,\,\, U_{ij}
\label{ncc4}\\
A_{ij|\mu}^{\alpha\beta} &\ll  \left( \frac{K_8}{b}+ \frac{K_9}{c}+ \frac{K_{10} b}{c} \right) A(u) \,\,\,\,\,\,\,\,\textnormal{on}\,\,\, U_{ij}    \label{ncc5} \\
W_{i|\mu}^{\alpha \gamma} &\ll  \left( \frac{K_7}{c} + \frac{K_6b}{c} \right) A(u)\,\,\,\,\,\,\,\,\textnormal{on}\,\,\, U_i^\delta            \label{ncc6}
\end{align}
where $K_1,k_2, K_3, \cdots$ are constants independent of $\mu$.
\end{lemma}

\begin{proof}
The estimate $(\ref{ncc4})$ follows from \cite{Kod05} p.302. The estimate $(\ref{ncc5})$ follows from a similar way with $(\ref{tpp2})$ in Lemma \ref{tcc31}. It remains to estimate $W_{i|\mu}^{\alpha \gamma}$ from $(\ref{nc10})$
\begin{align}\label{ncc7}
\left[\sum_{\beta=1}^q  c_i^{\alpha\beta (\mu-1)}(\xi_i, u) w_i'^{\beta\gamma}(\xi_i, u)- \sum_{\beta=1}^q w_i^{\alpha\beta}\left(\varphi_i^{\mu-1}(\xi_i, u), s^{\mu-1}(u) \right) \frac{\partial \varphi_i^{\beta (\mu-1)}}{\partial \xi_i^\gamma} \right]_\mu
\end{align}
First we estimate the first term of $(\ref{ncc7})$. We note that
\begin{align}\label{ncc9}
&\left[ c_i^{\alpha\beta(\mu-1)}(\xi_i, u) w_i'^{\beta \gamma}(\xi_i, u) \right]_\mu =\left[ \left(c_i^{\alpha\beta(\mu-1)}(\xi_i, u)- c_{0i}^{\alpha\beta} + c_{0i}^{\alpha \beta} \right) \left( w_i'^{\beta \gamma}(\xi_i, u)  - w_{0i}'^{\beta \gamma} + w_{0i}'^{\beta \gamma} \right)     \right]_\mu  \\
&= \left[ \left( c_i^{\alpha\beta(\mu-1)}(\xi_i, u) - c_{0i}^{\alpha\beta}\right) \left( w_i'^{\beta \gamma}(\xi_i, u)- w_{0i}'^{\beta \gamma}\right)  \right]_\mu + \left[ c_{0i}^{\alpha\beta}\left(w_i'^{\beta \gamma}(\xi_i, u)- w_{0i}'^{\beta \gamma}  \right) \right]_\mu \notag
\end{align}
Since $w_i'^{\beta \gamma}(\xi_i, u)$ is a holomorphic function, we may assume that
\begin{align}\label{ncc8}
w_i'^{\beta \gamma}(\xi_i, u) - w_{0i}'^{\beta \gamma} \ll A_5(u) \ll \frac{b_5}{b} A(u)\,\,\,\,\,\,\,\,\,\,\,A(u) =\frac{b_5}{16c_5}\sum_{v=1}^\infty \frac{c_5^v \left(u_1+ \cdots + u_{m'} \right)^v}{v^2}
\end{align}
holds for $\xi_i \in U_i$ with $b_5>0$ and $c_5>0$ and $b>b_5$ and $c>c_5$. Then from $(\ref{ncc9})$ and $(\ref{ncc8})$ and the induction hypothesis
\begin{align}\label{ncc13}
\left[ \sum_{\beta=1}^q c_i^{\alpha\beta(\mu-1)}(\xi_i, u) w_i'^{\beta \gamma}(\xi_i, u) \right]_\mu \ll K_4  \frac{b_5}{b} A(u)^2  + K_5\frac{b_5}{b} A(u) \ll \left( \frac{K_6}{c}+ \frac{K_7}{b} \right) A(u)
\end{align}
We estimate the second term of $(\ref{ncc7})$. We expand $w_i^{\alpha\beta}\left(\xi_i+ y, t \right)$ into power series in $y_1,...,y_n, t_1,...,t_m$. Then we may assume that
\begin{align}
w_i^{\alpha\beta}(\xi_i + y, t)- w_{0i}^{\alpha\beta} \ll \sum_{v=1}^\infty a_2^v\left( y_1+ \cdots + y_n + t_1 + \cdots + t_m \right)^v,\,\,\, a_2 >0
\end{align}
If we set $y=\varphi_i^{\mu-1}(\xi_i, u)- \xi_i$ and $t= s^{\mu-1}(u)$, then since $y \ll A(u)$ and $ t \ll A(u)$ by the induction hypothesis, we obtain
\begin{align}
w_i^{\alpha\beta}\left(\varphi_i^{\mu-1}(\xi_i, u), s^{\mu-1}(u) \right) - w_{0i}^\beta \ll \sum_{v=1}^\infty a_2^v (n+m)^v A(u)^v \ll \sum_{v=1}^{\infty} a_2^v (m+n)^v \left(\frac{b}{c} \right)^{v-1}A(u)
\end{align}
and from $(\ref{ncc10})$ and $(\ref{tcc21})$ we have
\begin{align}\label{ncc11}
\frac{\partial \varphi_i^{\beta(\mu-1)}(\xi_i, u)- \xi_i^\beta }{\partial \xi_i^\gamma} \ll \frac{A(u)}{\delta} \,\,\,\,\,\,\textnormal{on}\,\,\, U_i^\delta
\end{align}
We also note that in a similar way with \cite{Kod05} p.301, we can show that if we take $c$ such that $\frac{ba_3(m+n)}{c} < \frac{1}{2}$ for some constant $a_3>0$, we may assume that
\begin{align}
\left[ w_i^{\alpha \gamma}\left( \varphi_i^{\mu-1}(\xi_i, u), s^{\mu-1}(u) \right) \right]_\mu \ll \frac{2ba_3^2(m+n)^2}{c} A(u)
\end{align}
Then we have
\begin{align}
&\left[\sum_{\beta=1}^q w_i^{\alpha\beta}\left(\varphi_i^{\mu-1}(\xi_i, u), s^{\mu-1}(u) \right) \frac{\partial \varphi_i^{\beta (\mu-1)}}{\partial \xi_i^\gamma} \right]_\mu =  \left[\sum_{\beta=1}^q \left( w_i^{\alpha\beta}\left(\varphi_i^{\mu-1}(\xi_i, u), s^{\mu-1}(u) \right) - w_{0i}^{\alpha\beta} + w_{0i}^{\alpha\beta} \right) \frac{\partial \left(\varphi_i^{\beta (\mu-1)} - \xi_i^\beta + \xi_i^\beta \right)}{\partial \xi_i^\gamma} \right]_\mu \\
&=  \left[\sum_{\beta=1}^q \left( w_i^{\alpha\beta}\left(\varphi_i^{\mu-1}(\xi_i, u), s^{\mu-1}(u) \right) - w_{0i}^{\alpha\beta} \right) \frac{\partial \left(\varphi_i^{\beta (\mu-1)} - \xi_i^\beta  \right)}{\partial \xi_i^\gamma} \right]_\mu    + \left[ w_i^{\alpha \gamma}\left( \varphi_i^{\mu-1} (\xi_i, u), s^{\mu-1}(u)\right)\right]_\mu               \notag \\
& \ll K_5\sum_{v=1}^\infty a_2^v (m+n)^v \left( \frac{b}{c} \right)^{v-1}\frac{1}{\delta} A(u)^2  + \frac{2ba_3^2(m+n)^2}{c} A(u)  \ll \frac{K_5}{\delta} \frac{a_2 b (m+n)}{c} \sum_{v=0}^\infty \left( \frac{a_2b(m+n)}{c} \right)^v A(u)  + \frac{2ba_3^2(m+n)^2}{c} A(u)     \notag 
\end{align}
If we take a constant $c$ such that $\frac{ba_2(m+n)}{c}<\frac{1}{2}$, we obtain
\begin{align}\label{ncc12}
\left[\sum_{\beta=1}^q w_i^{\alpha\beta}\left(\varphi_i^{\mu-1}(\xi_i, u), s^{\mu-1}(u) \right) \frac{\partial \varphi_i^{\beta (\mu-1)}}{\partial \xi_i^\gamma} \right]_\mu  \ll 2\frac{K_5}{\delta}\frac{a_2 b (m+n)}{c} A(u)+ \frac{2ba_3^2(m+n)^2}{c} A(u) \ll K_6 \frac{b}{c} A(u) \,\,\,\textnormal{on}\,\,\, U_i^\delta
\end{align}

Then from $(\ref{ncc13})$ and $(\ref{ncc12})$ we have
\begin{align}
W_{i|\mu}^{\alpha \gamma} \ll \left(  \frac{K_7}{c}  + K_6\frac{b}{c}\right) A(u) \,\,\,\,\,\,\textnormal{on}\,\,\, U_i^\delta
\end{align}
This completes the proof of Lemma \ref{ncc20}.

\end{proof}

We recall Remark \ref{uuc19}. For any $\sigma=\left( -W, \left(\Gamma, A \right) \right) =\left( \left\{ -W_i \right\}, \left\{ \left(\Gamma_{ij} ,  A_{ij} \right) \right\} \right)\in C^0\left( \mathcal{U}, \mathscr{H}om_{\mathcal{O}_M} \left( \mathcal{N}_{\mathcal{F}_0}^*, \Omega_M^1     \right) \right) \bigoplus C^1\left( \mathcal{U}, \mathcal{E}_{\mathcal{N}_{\mathcal{F}_0}^*} \right)$ which is a $1$-cocylce in the \v Cech resolution $(\ref{uuc18})$ of $\mathcal{E}_{\mathcal{N}_{\mathcal{F}_0}^*}^\bullet$, where $W_i\in \Gamma\left( U_i, \mathscr{H}om_{\mathcal{O}_M}\left( \mathcal{N}_{\mathcal{F}_0}^*, \Omega_M^1 \right) \right)$ defined by
\begin{align} \label{ssc13}
W_i : \Gamma\left( U_i, \mathcal{N}_{\mathcal{F}_0}^*  \right) &\to \Gamma\left( U_i, \Omega_M^1  \right)\\
 w_{0i}^\alpha &\mapsto W_i^\alpha:=\sum_{\gamma=1}^n W_i^{\alpha \gamma} dz_i^\gamma \notag
\end{align}
and $\Gamma_{ij}\in \Gamma\left( U_{ij}, \Theta_M \right)$ and $A_{ij}\in \Gamma\left( U_{ij}, \mathscr{H}om_{\mathcal{O}_M}\left( \mathcal{N}_{\mathcal{F}_0}^*, \mathcal{N}_{\mathcal{F}_0}^*   \right) \right)$ defined by
\begin{align} \label{ssc15}
A_{ij}: \Gamma\left( U_{ij}, \mathcal{N}_{\mathcal{F}_0}^* \right) &\to \Gamma\left( U_{ij} , \mathcal{N}_{\mathcal{F}_0}^* \right) , \\
 w_{0i}^\alpha &\mapsto A_{ij}^\alpha:= \sum_{\beta=1}^q A_{ij}^{\alpha\beta} w_{0i}^\beta \notag
\end{align}

we define the norm $\left| \left| \sigma \right| \right|$
\begin{align*}
\left| \left| \sigma \right| \right| =\left| \left| \Gamma \right| \right| + \left|\left| A \right| \right| + \left|\left| W\right|\right|
\end{align*}
where
\begin{align*}
\left|\left| \Gamma \right|\right|=\max_{i,j} \sup_{z_i\in U_i\cap U_j} \left| \Gamma_{ij} (z_i)\right| , \,\,\,\,\,\,\,\,\,\,\,\,\,\,\left|\left| A \right|\right|=\max_{i,j}\max_{\alpha,\beta}\sup_{z_j \in U_i\cap U_j}\left|A_{ij}^{\alpha\beta} \right|\,\,\,\,\,\,\,\,\,\,\,\,\,\,\left|\left| W \right| \right| = \max_{i }\max_{\alpha,\gamma}\sup_{z_i\in U_i^\delta} \left| W_i^{\alpha\gamma} (z_i) \right|
\end{align*}

\begin{lemma}\label{ncc29}
For any pair $\sigma=\left(-W, \left(\Gamma, A \right) \right)=\left(\left\{ - W_i \right\}, \left\{ \left( \Gamma_{ij} , A_{ij} \right) \right\} \right)$ which is a $1$-cocycle in the \v Cech resolution of $\mathcal{E}_{\mathcal{N}_{\mathcal{F}_0}^*}^\bullet$, we can find $\varphi_i(\xi), s^\lambda$ and $c_i^{\alpha \gamma}(\xi_i)$ satisfying
\begin{align}
\varphi_i- \varphi_j & = - \Gamma_{ij} + \sum_{\lambda=1}^m s^\lambda \rho_{ij\lambda} \label{ncc21}\\
\mathcal{L}_{\varphi_i}\left( w_{0i}^\alpha \right) - \sum_{\beta=1}^q & c_{i|\mu}^{\alpha\beta} w_{0i}'^\beta = W_i^\alpha  - \sum_{\lambda=1}^m s_\mu^\lambda \beta_{i\lambda}^\alpha \\
\left|\varphi_i\left(\xi_i \right) \right| \leq K \left|\left| \sigma \right|\right|,\,\,\,\,\,\,\,\,&\left|s \right| \leq K \left|\left| \sigma \right|\right|,\,\,\,\,\,\,\,\,\,\,\,\left| c_i^{\alpha \gamma}(\xi_i) \right| \leq K \left|\left| \sigma \right|\right| \label{ncc22}
\end{align}
where $K$ is a constant independent of $\sigma$.
\end{lemma}

\begin{proof}
We define
\begin{align}
\iota\left(\sigma \right) = \inf \max_{i,\lambda, \alpha, \beta} \left\{\sup_{\xi_i \in U_i} \left|\varphi_i \right| , \left| s^\lambda \right|, \sup_{\xi_i\in U_i} \left| c_i^{\alpha\beta} \right|   \right\}
\end{align}
where $\inf$ is taken with respect to all the solutions $s^\lambda, \varphi_i(\xi_i), b_i^{\alpha\beta}(\xi_i)$ of $(\ref{ncc21})$ and $(\ref{ncc22})$. It suffices to show that there exists a constant $K$ such that $i(\sigma) \leq K \left|\left| \sigma \right|\right|$ for all $1$-cocycle $\sigma$ of $\mathcal{E}_{\mathcal{N}_{\mathcal{F}_0}^*}^\bullet$. Suppose that there is no such constant $K$. Then we can find a sequence $\sigma^{(1)}, \sigma^{(2)},\cdots, \sigma^{(v)}, \cdots$ of triple $\sigma^{(v)}=\left( - W^{(v)}, \left( \Gamma^{(v)}, A^{(v)} \right) \right)$ of $1$-cocycle of $\mathcal{E}_{\mathcal{N}_{\mathcal{F}_0}^*}^\bullet$ such that $\iota\left( \sigma^{(v)} \right) =1 $ and $\left|\left| \sigma^{(v)} \right|\right| < \frac{1}{v}$ such that

\begin{align}
\varphi_i^{(v)}- \varphi_j^{(v)} & = - \Gamma_{ij}^{(v)} + \sum_{\lambda=1}^m s^{\lambda (v)}\rho_{ij\lambda} \label{ncc24}\\
\mathcal{L}_{\varphi_i^{(v)}}\left( w_{0i}^\alpha \right) - \sum_{\beta=1}^q & c_{i|\mu}^{\alpha\beta(v)} w_{0i}'^\beta = W_i^{\alpha(v)}  - \sum_{\lambda=1}^m s_\mu^{\lambda (v)}\beta_{i\lambda}^\alpha \label{ncc25} \\
\left| \varphi_i^{(v)}  \right| < 2, \,\,\,\,\,\,\,\,\,\,&\left| s^{\lambda(v)} \right| <2,\,\,\,\,\,\,\,\,\,\,\,\left| c_i^{\alpha \gamma(v)}\right| <2 \label{ncc23}
\end{align}
Since $\left| W_i^{(v)} \right|\to 0$ on $U_i^\delta$, taking the limit for $v \to \infty$, we have
\begin{align*}
0=- \sum_{\gamma=1}^n \sum_{\beta=1}^q c_i^{\alpha\beta} w_{0i}'^\beta + \mathcal{L}_{\varphi_i}\left( w_{0i}^\alpha \right) - \sum_{\lambda=1}^m s^{\lambda } \beta_{i\lambda}^\alpha \,\,\,\,\,\,\,\textnormal{on}\,\,\, U_i^\delta
\end{align*}
Hence, by $(\ref{ncc23})$, replacing $\sigma^{(1)},\sigma^{(2)},\cdots$ by a suitable subsequence, we may assume $\varphi_i^{(v)}$ converges uniformly on each compact susbet of $U_i$, and $s^{(v)}$ converges, and $c_i^{\alpha \gamma (v)}$ converges uniformly on each compact subset of $U_i$. As in \cite{Kod05} p.296, we can deduce that $\varphi_i^{(v)}$ converges uniformly on the whole $U_i$. Now we claim that $c_i^{\alpha \gamma (v)}$ converges uniformly on the whole $U_i$. In fact, first we note that $b_i^{\alpha \gamma (v)}$ converges uniformly on $U_i^\delta$ from $(\ref{ncc10})$. On $U_i\cap U_j^\delta$, $(\ref{ncc24})$ and $(\ref{ncc25})$ imply that
\begin{align}\label{ncc26}
\sum_{\beta, \sigma=1}^q A_{ij}^{\alpha\beta (v)} h_{0ij}^{\beta\sigma} c_{0j}^{\sigma \gamma} =  -\sum_{\beta=1}^q c_{i}^{\alpha\beta (v)} h_{0ij}'^{\beta\gamma} + \sum_{\beta=1}^q \left[ \varphi_{j}^{(v)} , h_{0ij}^{\alpha\beta}  \right]  c_{0j}^{\beta\gamma} + \sum_{\beta=1}^q \sum_{\lambda=1}^m \frac{\partial h_{ij}^{\alpha\beta}}{\partial t_\lambda }|_{t=0} s^{\lambda (v)} c_{0j}^{\beta \gamma} + \sum_{\beta=1}^q h_{0ij}^{\beta\gamma} c_{j}^{\beta\gamma (v)} 
\end{align}
Since $s^{\lambda (v)}$ converges, and $\varphi_j$ converges uniformly on $U_j$, and $c_j^{\beta \gamma(v)}$ converges uniformly on $U_j^\delta$, and $\left| A_{ij}^{\alpha\beta(v)}\right| \leq \left|\left| A \right|\right|\to 0$, $(\ref{ncc26})$ implies that $c_i^{\alpha \beta(v)}$ converges uniformly on the whole $U_i$. Then we put $\varphi_i=\lim \varphi_i^{(v)}$ and $s^\lambda=\lim s^{\lambda (v)}$ and $c_i^{\alpha \gamma}= \lim c_i^{\alpha \gamma (v)}$. Then from $(\ref{ncc24})$ and $(\ref{ncc25})$ and $\left| \Gamma_{ij}^{(v)} \right|\to 0$ on $U_{ij}$ and $\left| W_i^{\alpha \gamma (v)}\right|\to 0$ on $U_i^\delta$, we have
\begin{align}
\varphi_i- \varphi_j & =  \sum_{\lambda=1}^m s^{\lambda (v)}\rho_{ij\lambda} \,\,\,\,\,\,\,\,\,\textnormal{on}\,\,\, U_{ij} \label{ncc27} \\
\mathcal{L}_{\varphi_i}\left( w_{0i}^\alpha \right) - \sum_{\beta=1}^q & c_{i}^{\alpha\beta} w_{0i}'^\beta =  - \sum_{\lambda=1}^m s_\mu^{\lambda }\beta_{i\lambda}^\alpha\,\,\,\,\,\,\,\,\,\textnormal{on}\,\,\, U_i^\delta \notag
\end{align}
Then by identity theorem, we obtain
\begin{align}
\mathcal{L}_{\varphi_i}\left( w_{0i}^\alpha \right) - \sum_{\beta=1}^q & c_{i}^{\alpha\beta} w_{0i}'^\beta =  - \sum_{\lambda=1}^m s_\mu^{\lambda }\beta_{i\lambda}^\alpha\,\,\,\,\,\,\,\,\,\textnormal{on}\,\,\, U_i \label{ncc28}
\end{align}
We put $\tilde{s}^{\lambda(v)}= s^{\lambda(v)}- s^\lambda$ and $\tilde{\varphi}_j^{(v)}= \varphi_j^{(v)} - \varphi_j$ and $\tilde{c}_i^{\alpha\beta(v)}=c_i^{\alpha\beta(v)}- c_i^{\alpha\beta}$. Then for a sufficiently small large $v$, we have
\begin{align}
\left| \tilde{\varphi}_i^{(v)} \right| < \frac{1}{2},\,\,\,\,\,\,\,\,\,\left| \tilde{s}^{\lambda (v)} \right| < \frac{1}{2},\,\,\,\,\,\,\,\,\,\,\,\left| c_i^{\alpha \gamma (v)} \right| < \frac{1}{2}
\end{align}
while we infer from $(\ref{ncc24}), (\ref{ncc25}), (\ref{ncc27})$ and $(\ref{ncc28})$ that
\begin{align*}
\tilde{\varphi}_i^{(v)}- \tilde{\varphi}_j^{(v)} & = - \Gamma_{ij}^{(v)} + \sum_{\lambda=1}^m \tilde{s}^{\lambda (v)}\rho_{ij\lambda} \\
\mathcal{L}_{\tilde{\varphi}_i^{(v)}}\left( w_{0i}^\alpha \right) - \sum_{\beta=1}^q & \tilde{c}_{i}^{\alpha\beta(v)} w_{0i}'^\beta = W_i^{\alpha(v)}  - \sum_{\lambda=1}^m \tilde{s}^{\lambda (v)}\beta_{i\lambda}^\alpha 
\end{align*}
This contradicts to $\iota\left( \sigma \right)=1$.
\end{proof}
Then by Lemma \ref{ncc20} and Lemma \ref{ncc29} we can choose solutions $\varphi_{i|\mu}(\xi_i,u)$ and $s_{\mu}^\lambda (u)$ and $c_{i|\mu}^{\alpha \beta}(\xi_i, u)$ of equations $(\ref{nc17})$ and $(\ref{nc19})$ such that
\begin{align*}
\varphi_{i|\mu}(\xi_i, u)\ll K K^* A(u),\,\,\,\,\,\,\,\,\,\,\,\,\,s_\mu(u)\ll K K^* A(u),\,\,\,\,\,\,\,\,\,\,\,\,\, c_{i|\mu}^{\alpha \beta} \ll K K^* A(u)
\end{align*}
where $K^*=\frac{K_1+ K_{8}}{b}+ \frac{K_2 + K_{9}+ K_{7}}{c} + \frac{( K_3 + K_{10} + K_6)b}{c}$. We choose $b$ and $c$ in a way that $K K^* <1$. Then we have $\varphi_{i|\mu}(\xi_i, u) \ll A(u), s_\mu (u)\ll A(u)$ and $c_{i|\mu}^{\alpha\beta}(\xi_i, u) \ll A(u)$. Hence the induction $(\ref{ncc3})_\mu$ holds for $\mu$. This completes the proof of Theorem \ref{nc1}.

\end{proof}

\section{Deformations of foliated complex analytic structures in terms of both tangent and cotangent sheaves} \label{ss1}

\begin{definition}[compare Definition \ref{t3} and Definition \ref{n3}]\label{s1}
Suppose that given a domain $B\subset \mathbb{C}^m$, there is a set $\left\{\left(M_t, \Theta_{\mathcal{F}_t}, \mathcal{N}_{\mathcal{F}_t}^*\right)|t\in B\right\}$ of $n$-dimensional compact $($singularly$)$ foliated complex manifolds $(M_t,\mathcal{F}_t)$ of dimension $p$ and codimension $q=n-p$, so that for each $t\in B$, we have exact sequences
\begin{align}
0\to \Theta_{\mathcal{F}_t}\to \Theta_{M_t}\to \Theta_{M_t}/\Theta_{\mathcal{F}_t}\to 0\label{s2}\\
0\to \mathcal{N}_{\mathcal{F}_t}^*\to \Omega_{M_t}^1\to \Omega_{M_t}^1/\mathcal{N}_{\mathcal{F}_t}^*\to 0 \label{s3}
\end{align}
where $\Theta_{\mathcal{F}_t}=\left(\Omega_{M_t}^1/\mathcal{N}_{\mathcal{F}}^*\right)^*$ and $\mathcal{N}_{\mathcal{F}_t}^*= \left(\Theta_{M_t}/ \Theta_{\mathcal{F}_t}\right)^*$, and $\Theta_{M_t}/\Theta_{\mathcal{F}_t}$ and $\Omega_{M_t}^1/\mathcal{N}_{\mathcal{F}_t}^*$ are torsion-free $\mathcal{O}_{M_t}$-modules. We say that $\left\{\left(M_t, \Theta_{\mathcal{F}_t}, \mathcal{N}_{\mathcal{F}_t}\right)|t\in B\right\}$ is a family of $($singularly$)$ foliated complex manifold or $($singularly$)$ foliated complex analytic family in terms of both tangent and cotangent sheaves, or foliated analytic family for simultaneous deformations if there exists a complex manifold $\mathcal{M}$ and a holomorphic map $\pi:\mathcal{M}\to B$ such that
\begin{enumerate}
\item we have exact sequences
\begin{align}
0 \to \Theta_\mathcal{F}\to \Theta_{\mathcal{M}/B}\to \Theta_{\mathcal{M}/B}/\Theta_\mathcal{F} \to 0 \label{s4}\\
0\to \mathcal{N}_\mathcal{F}^*\to \Omega_{\mathcal{M}/B}^1\to \Omega_{\mathcal{M}/B}^1/ \mathcal{N}_\mathcal{F}^*\to 0 \label{s5}
\end{align}
where $\Theta_\mathcal{F}$ is a coherent subsheaf of the relative tangent sheaf $\Theta_{\mathcal{M}/B}$ over $B$, and $\mathcal{N}_\mathcal{F}^*$ is a coherent subsheaf of the relative cotangent sheaf $\Omega_{\mathcal{M}/B}^1$ over $B$
\item $\Theta_{\mathcal{F}}=(\Omega_{\mathcal{M}/B}^1/\mathcal{N}_\mathcal{F}^*)^*$ and $\mathcal{N}_\mathcal{F}^*=(\Theta_{\mathcal{M}/B}/\Theta_\mathcal{F})^*$.
\item for each $t\in B$, $\pi^{-1}(t)=\left(M_t, \Theta_{\mathcal{F}_t},\mathcal{N}_{\mathcal{F}_t}^*\right)$, and the exact sequence $(\ref{s2})$ is induced from $(\ref{s4})$ and $(\ref{s3})$ is induced from $(\ref{s5})$ by restricting to $M_t$, respectively.
\item the rank of the Jacobian of $\pi$ is equal to $m$ at every point of $\mathcal{M}$.
\item $[\Theta_\mathcal{F},\Theta_\mathcal{F}]\subset \Theta_\mathcal{F}$ under the Lie bracket, and $\Theta_\mathcal{M}/\Theta_\mathcal{F}$ is flat over $B$.
\item $d_{\mathcal{M}/B}(\mathcal{N}_\mathcal{F}^*) \subset \mathcal{N}_\mathcal{F}^*\bigwedge \Omega_{\mathcal{M}/B}^1$ at every point of $\mathcal{M}-\bigcup_{t\in B} S_t$, where $S_t=\textnormal{Sing}(\mathcal{F}_t)$ for $t\in B$, and $d_{\mathcal{M}/B}$ is the relative differential on $\Omega_{\mathcal{M}/B}^1$. By abuse of notation, we will denote $d_{\mathcal{M}/B}$ by $d$. 
\item $\Omega_{\mathcal{M}/B}^1/\mathcal{N}_\mathcal{F}^*$ is flat over $B$.
\end{enumerate}
We will denote the $($singularly$)$ foliated complex analytic family by $( \mathcal{M} , \Theta_\mathcal{F}, \mathcal{N}_\mathcal{F}^*, B, \pi)$ or simply by $(\mathcal{M}, \mathcal{F}, B, \pi)$.
\end{definition}

\begin{remark}
Let $\left(\mathcal{M}, \Theta_\mathcal{F}, \mathcal{N}_\mathcal{F}^* ,  B, \pi \right)$ with $\Theta_\mathcal{F}$ and $\mathcal{N}_\mathcal{F}^*$ locally free be a foliated analytic family. Let $\Delta$ be an open set of $B$. Then the restriction $\left(\mathcal{M}_\Delta= \pi^{-1}(\Delta), \Theta_{\mathcal{F}}|_\Delta , \mathcal{N}_\mathcal{F}^*|_\Delta, \pi |_{\mathcal{M}_\Delta} \right)$ is also a foliated complex analytic family in terms of both tangent sheaves and cotangent sheaves. We will denote the family by $\left( \mathcal{M}_\Delta, \Theta_{\mathcal{F}_\Delta}, \mathcal{N}_{\mathcal{F}_\Delta}^*,  \Delta, \pi\right)$.
\end{remark}

\begin{remark}
If $\mathcal{F}_t$ is a regular foliation on $M_t$ for each $t\in B$, then \textnormal{Definition} $\ref{t3}$, \textnormal{Definition} $\ref{n3}$ and \textnormal{Definition} $\ref{s1}$ are all equivalent.
\end{remark}

\begin{remark}
In the following we shall assume that $\Theta_\mathcal{F}$ and $\mathcal{N}_\mathcal{F}^*$ in \textnormal{Definition} $\ref{s1}$ are both locally free $\mathcal{O}_\mathcal{M}$-modules of rank $p$ and rank $q$, respectively.
\end{remark}

\begin{remark}\label{ss17}
We shall keep the notations in \textnormal{Remark \ref{ta1}} and \textnormal{Remark \ref{n5}}. Then we only need to discuss the condition $(2)$  $\Theta_{\mathcal{F}}=\left(\frac{\Omega_{\mathcal{M}/B}^1}{\mathcal{N}_\mathcal{F}^*}\right)^*$ and $\mathcal{N}_\mathcal{F}^*= \left( \frac{\Theta_{\mathcal{M}/B}}{\Theta_\mathcal{F} } \right)^*$ in \textnormal{Definition \ref{s1}}. The condition $(2)$ implies that  for $\alpha=1,...,p$ and $\beta=1,...,q$, we have
\begin{align}\label{s6}
i_{T_j^\alpha(z_j,t)}\left(w_j^\beta(z_j,t) \right)=i_{\sum_{\gamma=1}^n T_j^{\alpha\gamma}(z_j,t)\frac{\partial}{\partial z_j^\gamma} }\left( \sum_{\eta=1}^n w_j^{\beta\eta}(z_j,t)dz_j^\eta \right)= \sum_{\gamma=1}^n T_j^{\alpha\gamma}(z_j,t) w_j^{\beta \gamma}(z_j,t)=0
\end{align}
\end{remark}

\subsection{Complex controlling simultaneous deformations of singular holomorphic foliations in terms of both locally free subsheaves of tangent sheaves and locally free sheaves of cotangent sheaves}\

Let $\left(\mathcal{M}, \Theta_\mathcal{F}, \mathcal{N}_\mathcal{F}^*, B, \pi \right)$ be a foliated analytic family with $\Theta_{\mathcal{F}}$ and $\mathcal{N}_\mathcal{F}^*$ locally free as in Definition \ref{s1}, so that $\mathcal{F}_t=\left(\Theta_{\mathcal{F}_t}, \mathcal{N}_{\mathcal{F}_t}^* \right)$ defines a (singular) holomorphic foliation on a compact complex manifold $M_t$ with $\Theta_{\mathcal{F}_t}$ and $\mathcal{N}_{\mathcal{F}_t}^*$ locally free. Then we have the complex of sheaves $\mathcal{F}_t^\bullet$ on $M_t$ associated to $\mathcal{F}_t=\left( \Theta_{\mathcal{F}_t}, \mathcal{N}_{\mathcal{F}_t}^* \right)$ (see Appendix \ref{AppendixA4}).
\begin{center}
$\mathcal{F}_t^\bullet: \begin{CD}
\cdots \\
@AF_3^tAA \\
\mathscr{H}om_{\mathcal{O}_{M_t}}\left(\bigwedge^3 \Theta_{\mathcal{F}_t},  \frac{\Theta_{M_t}}{\Theta_{\mathcal{F}_t}} \right) \bigoplus \mathscr{H}om_{\mathcal{O}_{M_t}}\left(\mathcal{N}_{\mathcal{F}_t}^*, \tilde{\mathcal{S}}_t^3 \right) \bigoplus \mathscr{H}om_{\mathcal{O}_{M_t}}\left(\mathcal{N}_{\mathcal{F}_t}^*,\bigwedge^2 \Theta_{\mathcal{F}_t}^* \right)       \\
@AF_2^t AA \\
\mathscr{H}om_{\mathcal{O}_{M_t}}\left(\bigwedge^2 \Theta_{\mathcal{F}_t}, \frac{\Theta_{M_t}}{\Theta_{\mathcal{F}_t} }\right) \bigoplus \mathscr{H}om_{\mathcal{O}_{M_t}}\left(\mathcal{N}_{\mathcal{F}_t}^*, \tilde{\mathcal{S}}_t^2 \right) \bigoplus \mathscr{H}om_{\mathcal{O}_{M_t}}\left(\mathcal{N}_{\mathcal{F}_t}^*, \Theta_{\mathcal{F}_t}^* \right)       \\
@AF_1^t AA \\
\mathscr{H}om_{\mathcal{O}_{M_t}}\left(\Theta_{\mathcal{F}_t}, \frac{\Theta_{M_t}}{\Theta_{\mathcal{F}_t}} \right) \bigoplus \mathscr{H}om_{\mathcal{O}_M}\left( \mathcal{N}_{\mathcal{F}_t}^*, \frac{\Omega_{M_t}^1}{\mathcal{N}_{\mathcal{F}_t}^* } \right)  \\
@AF_0^t AA \\
\Theta_{M_t}
\end{CD}$
\end{center}
We will denote the $i$-th cohomology group by $\mathbb{H}^1\left( M_t,  \mathcal{F}_t^\bullet \right)$. We can compute $\mathbb{H}^i\left( M_t, \mathcal{F}_t^\bullet \right)$ by the following \v Cech resolution of $\mathcal{F}_t^\bullet$. Here $\delta$ is the \v Cech map and $\mathcal{U}_t=\mathcal{U}\cap M_t=\left\{ U_j^t:=\mathcal{U}_j\cap M_t|j=1,2,... \right\}$ is an open covering of $M_t$:
{\tiny{\begin{center}
$\begin{CD}
\cdots \\
@AF_2^tAA \\
C^0\left( \mathcal{U}_t, \left(\bigwedge^2 \Theta_{\mathcal{F}_t}^*\otimes \frac{\Theta_{M_t}}{\Theta_{\mathcal{F}_t}}\right) \bigoplus \left(\left(\mathcal{N}_{\mathcal{F}_t}^* \right)^*\otimes \tilde{\mathcal{S}}_t^2\right) \bigoplus \left( \left( \mathcal{N}_{\mathcal{F}_t}^* \right)^*\otimes \Theta_{\mathcal{F}_t}^* \right)   \right) @>-\delta>> \cdots \\
@AF_1^tAA @AF_1^tAA \\
C^0\left(\mathcal{U}_t, \left(\Theta_{\mathcal{F}_t}^*\otimes \frac{\Theta_{M_t}}{\Theta_{\mathcal{F}_t}}\right) \bigoplus \left(\left(\mathcal{N}_{\mathcal{F}_t}^* \right)^*\otimes \frac{\Omega_{M_t}^1}{\mathcal{N}_{\mathcal{F}_t}^*} \right) \right) @>\delta>> C^1\left(\mathcal{U}_t, \left(\Theta_{\mathcal{F}_t}^*\otimes \frac{\Theta_{M_t}}{\Theta_{\mathcal{F}_t}}\right) \bigoplus \left(\left(\mathcal{N}_{\mathcal{F}_t}^* \right)^*\otimes \frac{\Omega_{M_t}^1}{\mathcal{N}_{\mathcal{F}_t}^*} \right) \right) @>-\delta>> \cdots \\
@AF_0^tAA @AF_0^tAA  @AF_0^tAA \\
C^0\left(\mathcal{U}_t, \Theta_{M_t} \right) @>-\delta >> C^1\left( \mathcal{U}_t, \Theta_{M_t} \right) @>>> C^2\left(\mathcal{U}_t, \Theta_{M_t} \right) @>\delta>> \cdots
\end{CD}$
\end{center}}}

We will relate the first cohomology group $\mathbb{H}^1\left( M_t, \mathcal{F}_t^\bullet \right)$ to infinitesimal foliated deformations of $\pi^{-1}(t)= \left( M_t, \Theta_{\mathcal{F}_t}, \mathcal{N}_{\mathcal{F}_t}^*  \right)$ in the foliated analytic family $\left( \mathcal{M}, \Theta_\mathcal{F}, \mathcal{N}_\mathcal{F}^*, B, \pi  \right)$ in terms of both tangent sheaves and cotangent sheaves.

\subsection{Infinitesimal foliated deformations in terms of both tangent sheaves and cotangent sheaves}

We shall keep the notations in subsection $\ref{ta11}$ and subsection $\ref{n13}$.

\begin{proposition}\label{s7}
{\Small{\begin{align*}
&\left(\left\{\theta_{jk}(t):=\sum_{\alpha=1}^n \frac{\partial f_{jk}^\alpha(z_k,t)}{\partial t}\frac{\partial}{\partial z_j^\alpha} \right\},  \left\{ \alpha_j(t):=\left(T_j^\alpha(z_j,t)\mapsto -\overline{\frac{\partial T_j^\alpha(z_j,t)}{\partial t}}\right)_{\alpha=1,...,p} \right\}  , \left\{ \beta_j(t):=\left(w_j^\alpha(z_j,t)\mapsto -\overline{\frac{\partial w_j^\alpha(z_j,t)}{\partial t}}\right)_{\alpha=1,...,q} \right\} \right) \\&\in C^1\left(\mathcal{U}^t, \Theta_{M_t} \right) \bigoplus C^0\left( \mathcal{U}_t, \mathscr{H}om_{\mathcal{O}_{M_t}}\left(\Theta_{\mathcal{F}_t}, \frac{\Theta_{M_t}}{\Theta_{\mathcal{F}_t}} \right) \right) \bigoplus C^0\left( \mathcal{U}_t, \mathscr{H}om_{\mathcal{O}_{M_t}}\left(\mathcal{N}_{\mathcal{F}_t}^*, \frac{\Omega_{M_t}^1}{\mathcal{N}_{\mathcal{F}_t}^*} \right) \right)
\end{align*}}}
defines a $1$-cocycle in the above \v Cech resolution of $\mathcal{F}_t^\bullet $ and call its cohomology class in $\mathbb{H}^1\left( M_t, \mathcal{F}_t^\bullet \right)$ the infinitesimal $($foliated$)$ deformation along $\frac{\partial}{\partial t}$. This expression is independent of the choice of local coordinates.
\end{proposition}

\begin{proof}
We shall keep the notations in the proof of Proposition \ref{ta2} and Proposition \ref{n14}. Then it remains to show that $i_{T_j^\alpha} \left( \tilde{\beta}_j\left(w_j^\gamma\right) \right) + i_{\tilde{\alpha}_j\left(T_j^\alpha \right) } \left(w_j^\gamma \right)=0$ for $\alpha=1,...,p$ and $\gamma=1,...,q$. In fact, by taking the derivative of $(\ref{s6})$ with respect to $t$, we have
\begin{align} \label{ss18}
i_{T_j^\alpha(z_j,t)}\left( \frac{\partial w_j^\gamma(z_j,t)}{\partial t} \right) + i_{\frac{\partial T_j^\alpha(z_j,t)}{\partial t}}\left( w_j^\gamma (z_j,t) \right) =0 \Longrightarrow -i_{T_j^\alpha}\left(\tilde{\beta}_j\left(w_j^\gamma\right) \right) - i_{\tilde{\alpha}_j\left( T_j^\alpha \right)}\left( w_j^\gamma \right) = 0
\end{align}

\end{proof}

\begin{definition}[foliated Kodaira-Spencer map in terms of both tangent sheaf and cotangent sheaf]\label{s8}
Let $(\mathcal{M}, \Theta_\mathcal{F},  \mathcal{N}_\mathcal{F}^*, B, \pi)$ with $\Theta_\mathcal{F}$ and $\mathcal{N}_{\mathcal{F}}^*$ locally free be a foliated complex analytic family of deformations of $\left(M_t, \Theta_{\mathcal{F}_t}, \mathcal{N}_{\mathcal{F}_t}^*\right)=\pi^{-1}(t),t\in B$, where $B$ is a domain of $\mathbb{C}^m$. We keep the notations in \textnormal{Definition \ref{ta12}} and \textnormal{Definition \ref{n15}}. For a tangent vector $\frac{\partial}{\partial t}=\sum_{\lambda=1}^m c_\lambda \frac{\partial}{\partial t_\lambda}, c_\lambda \in \mathbb{C}$, of $B$, the foliated Kodaira-Spencer map in terms of both tangent and cotangent sheaf at $t$ is defined to be a $\mathbb{C}$-linear map
{\small{\begin{align*}
\varphi_t:T_t(B)&\to \mathbb{H}^1\left(M_t, \mathcal{F}_t^\bullet \right)\\
                \frac{\partial}{\partial t}&\mapsto \frac{\partial \left(M_t, \Theta_{\mathcal{F}_t},  \mathcal{N}_{\mathcal{F}_t}^* \right)}{\partial t}:=\left(\frac{\partial M_t}{\partial t}, \frac{ \partial \Theta_{\mathcal{F}_t}}{\partial t} , \frac{\partial \mathcal{N}_{\mathcal{F}_t}^*}{\partial t} \right) \\
                 &\,\,\,\,\,\,\,\,\,\,\in C^1\left(\mathcal{U}_t, \Theta_{M_t} \right) \bigoplus C^0\left(\mathcal{U}_t, \mathscr{H}om_{\mathcal{O}_{M_t}} \left( \Theta_{\mathcal{F}_t}, \frac{\Theta_{M_t}}{\Theta_{\mathcal{F}_t}} \right) \right) \bigoplus C^0\left( \mathcal{U}_t, \mathscr{H}om_{\mathcal{O}_{M_t}}\left(\mathcal{N}_{\mathcal{F}_t}^*, \frac{\Omega_{M_t}^1}{\mathcal{N}_{\mathcal{F}_t}^*} \right) \right)
\end{align*}}}
\end{definition}

\subsection{Preliminaries} \label{ss5}\

Let $\left(\mathcal{M}, \Theta_\mathcal{F}, \mathcal{N}_\mathcal{F}^* , B, \pi \right)$ with $\Theta_\mathcal{F}$ and $\mathcal{N}_\mathcal{F}^*$ locally free be a foliated analytic family in terms of both tangent and cotangent sheaves, where $B$ is a domain of $\mathbb{C}^m$ containing the origin $0$ as in Definition \ref{s1}. We shall keep the notations in subsection \ref{tt1} and subsection \ref{nn1}. We only need to add the condition  $\Theta_{\mathcal{F}}=\left(\frac{\Omega_{\mathcal{M}/B}^1}{\mathcal{N}_\mathcal{F}^*}\right)^*$ and $\mathcal{N}_\mathcal{F}^*= \left( \frac{\Theta_{\mathcal{M}/B}}{\Theta_\mathcal{F} } \right)^*$ in \textnormal{Definition \ref{s1}}, which implies that  for $\alpha=1,...,p$ and $\beta=1,...,q$, we have
\begin{align}\label{s6}
i_{T_j^\alpha(z_j,t)}\left(w_j^\beta(z_j,t) \right)=i_{\sum_{\gamma=1}^n T_j^{\alpha\gamma}(z_j,t)\frac{\partial}{\partial z_j^\gamma} }\left( \sum_{\eta=1}^n w_j^{\beta\eta}(z_j,t)dz_j^\eta \right)= \sum_{\gamma=1}^n T_j^{\alpha\gamma}(z_j,t) w_j^{\beta \gamma}(z_j,t)=0
\end{align}

With this preparation, we shall prove theorem of existence for deformations of foliated complex analytic structures in terms of both tangent and cotangent sheaves. 

\subsection{Theorem of existence for deformations of foliated complex analytic structures in terms of both tangent and cotangent sheaves}

\begin{theorem}[Theorem of existence for deformations of foliated complex analytic structures in terms of both tangent and cotangent sheaves]\label{ss6}
Let $\left( M , \mathcal{F}_0 \right)$ be a compact foliated complex manifold with both $\Theta_{\mathcal{F}_0}$ and $\mathcal{N}_{\mathcal{F}_0}^*$ locally free. Suppose that $\mathbb{H}^2\left( M, \mathcal{F}_0^\bullet \right)=0$. Then there exists a foliated analytic family $\left( \mathcal{M}, \Theta_{\mathcal{F}}, \mathcal{N}_\mathcal{F}^* , B, \pi \right)$ with $0\in B \subset \mathbb{C}^m$ satisfying the following conditions$:$
\begin{enumerate}
\item $\pi^{-1}(0)=\left( M, \mathcal{F}_0 \right)$
\item The foliated Kodaira-Spencer map $\varphi_0: T_0 B \to \mathbb{H}^1\left( M, \mathcal{F}_0^{\bullet} \right)$ in terms of both tangent and cotangent sheaf is an isomorphism.
\end{enumerate}
\end{theorem}

\begin{proof}
We shall keep the notation in the proof of Theorem \ref{tt2} and Theorem \ref{nn2}.

\subsection{Existence of formal solutions}\

We construct solutions of $(\ref{t10})-(\ref{te5})$ and $(\ref{n20})-(\ref{n211})$ which are power series in $t$. Additional we construct solution of
\begin{align}\label{se1}
i_{T_i^\alpha}\left( w_i^\beta \right)=0,\,\,\,\,\,\,\,\,\,\,\alpha=1,...,p, \,\,\,\beta=1,...,q
\end{align}
 Then $(\ref{t10})-(\ref{te5})$ are equivalent to the systems of congruences $(\ref{t13})_\mu-(\ref{te6})_\mu$ for $\mu=1,2,3,\cdots$ and $(\ref{n22})-(\ref{n211})$ are equivalent to the system of congruences $(\ref{n23})_\mu- (\ref{n241})_\mu$ for $\mu=1,2,3,\cdots$ and additionally $(\ref{se1})$ is equivalent to the system of congruenes
\begin{align}\label{se3}
i_{T_i^{\alpha\mu}}\left( w_i^{\beta\mu} \right)\equiv_\mu 0 \,\,\,\,\,\,\alpha=1,...,p, \,\,\,\beta=1,...,q
\end{align}
for $\mu=1,2,3,\cdots$. We define homogenous polynomials $(\ref{tt3})-(\ref{tt4})$ and $(\ref{nn3})-(\ref{nn4})$, and additionally
\begin{align}\label{se6}
K_{i|\mu}^{\alpha\beta} \equiv_\mu  i_{T_i^{\alpha(\mu-1)}}\left( w_i^{\beta(\mu-1)} \right),\,\,\,\,\,\,\alpha=1,...,p,\,\,\,\beta=1,...,q
\end{align}

\begin{lemma}\label{ss8}
We have the equalities $(\ref{t19})-(\ref{te7})$ and $(\ref{n25})-(\ref{n53})$ and additionally
\begin{align} \label{se2}
K_{i|\mu}^{\alpha\beta} - \sum_{\xi=1}^p \sum_{\gamma=1}^q r_{0ij}^{\alpha \xi} h_{0ij}^{\beta\gamma} K_{j|\mu}^{\xi \gamma}=i_{\Gamma_{ij|\mu}^\alpha}\left( w_{0i}^\beta \right)  + i_{T_{0i}^\alpha}\left(C_{ij|\mu}^\beta \right)
\end{align}
\begin{align}\label{se4}
\bar{\partial} K_{i|\mu}^{\alpha\beta}= - i_{\Phi_{i|\mu}^\alpha}\left( w_{0i}^\beta \right) + i_{T_{0i}^\alpha}\left( A_{i|\mu}^\beta \right)
\end{align}
\end{lemma}

\begin{proof}
$(\ref{t19})-(\ref{te7})$ follows from Lemma \ref{te10}, and $(\ref{n25})-(\ref{n53})$ follows from Lemma \ref{ne10}. We prove $(\ref{se2})$. In fact, from $(\ref{se6})$ and $(\ref{tt6})$ and $(\ref{nn6})$, we have
\begin{align*}
K_{i|\mu}^{\alpha\beta}&\equiv_\mu i_{T_i^{\alpha(\mu-1)}}\left( w_i^{\beta(\mu-1)}\right) =i_{\Gamma_{ij|\mu}^\alpha+\sum_{\xi=1}^p r_{ij}^{\alpha\xi(\mu-1)} T_j^{\xi(\mu-1)} } \left( C_{ij|\mu}^\beta+\sum_{\gamma=1}^q h_{ij}^{\beta\gamma(\mu-1)} w_j^{\gamma(\mu-1)}\right)\\
&\equiv_\mu i_{\Gamma_{ij|\mu}^\alpha}\left(w_{0i}^\beta \right) + i_{T_{0i}^\alpha}\left( C_{ij|\mu}^\beta \right) +\sum_{\xi=1}^p \sum_{\gamma=1}^q r_{ij}^{\alpha\xi(\mu-1)} h_{ij}^{\beta \gamma(\mu-1)} i_{T_j^{\xi(\mu-1)}}\left(w_j^{\gamma(\mu-1)} \right)\\
&\equiv_\mu i_{\Gamma_{ij|\mu}^\alpha}\left( w_{0i}^\beta \right)  + i_{T_{0i}^\alpha}\left(C_{ij|\mu}^\beta \right) +\sum_{\xi=1}^p \sum_{\gamma=1}^q r_{0ij}^{\alpha\xi} h_{0ij}^{\beta \gamma} K_{j|\mu}^{\xi \gamma}
\end{align*}
We prove $(\ref{se4})$. In fact, from $(\ref{se6})$ and $(\ref{tt7})$ and $(\ref{nn7})$, we have
\begin{align*}
\bar{\partial} K_{i|\mu}^{\alpha\beta}&\equiv_\mu i_{\bar{\partial} T_i^{\alpha(\mu-1)}} \left(w_i^{\beta(\mu-1)} \right)- i_{T_i^{\alpha(\mu-1)}}\left(\bar{\partial} w_i^{\beta(\mu-1)} \right)\\
                                                         &\equiv_\mu i_{\left[ \varphi^{\mu-1}, T_i^{\alpha(\mu-1)}\right]-\Phi_{i|\mu}^\alpha}\left( w_i^{\beta(\mu-1)}\right) - i_{T_i^{\alpha(\mu-1)}}\left( \mathcal{L}_{\varphi^{\mu-1}}\left( w_i^{\beta(\mu-1)} \right) - A_{i|\mu}^\beta \right)\\
                                                         &= \mathcal{L}_{\varphi^{\mu-1}}\left( K_{i|\mu}^{\alpha\beta}\right) + i_{T_i^{\alpha(\mu-1)}}\left( \mathcal{L}_{\varphi^{\mu-1}} \left( w_i^{\beta(\mu-1)}  \right) \right) -i_{\Phi_{i|\mu}^\alpha}\left( w_{0i}^\beta \right) - i_{T_i^{\alpha(\mu-1)}}\left( \mathcal{L}_{\varphi^{\mu-1}}\left(w_i^{\beta(\mu-1)} \right) \right)+ i_{T_{0i}^\alpha}\left( A_{i|\mu}^\beta \right) \\
                                                         &=  -i_{\Phi_{i|\mu}^\alpha}\left( w_{0i}^\beta \right) + i_{T_{0i}^\alpha}\left( A_{i|\mu}^\beta \right)
\end{align*}
\end{proof}

Our purpose is to construct $\varphi^\mu=\varphi^{\mu-1}+ \varphi_\mu, r_{ij}^{\alpha\beta\mu}= r_{ij}^{\alpha\beta(\mu-1)}+ r_{ij|\mu}^{\alpha\beta}, T_i^{\alpha\mu}=T_i^{\alpha(\mu-1)}+ T_{i|\mu}^\alpha$,  $g_{i\alpha\beta}^{\gamma \mu}= g_{i\alpha\beta}^{\gamma(\mu-1)}+ g_{i\alpha\beta|\mu}^\gamma$ and $h_{ij}^{\alpha\beta \mu}= h_{ij}^{\alpha\beta (\mu-1)} + h_{ij|\mu}^{\alpha\beta}, w_i^{\alpha \mu}=w_i^{\alpha (\mu-1)}+ w_{i|\mu}^\alpha, a_{i_x}^{\alpha\beta \mu}= a_{i_x}^{\alpha\beta(\mu-1)}+ a_{i_x|\mu}^{\alpha\beta}$ satisfying $(\ref{t13})_\mu-(\ref{te6})_\mu$ and $(\ref{n23})_\mu-(\ref{n241})_\mu$
and additionally $(\ref{se3})_\mu$.

\begin{lemma}\label{ss10}
$(\ref{t13})_\mu-(\ref{te6})_\mu$ are equivalent to $(\ref{te11})-(\ref{t33})$, and $(\ref{n23})_\mu-(\ref{n241})_\mu$ are equivalent to $(\ref{n42}) -(\ref{n48})$ and additionally $(\ref{se3})_\mu$ is equivalent to
\begin{align}\label{se5}
- K_{i|\mu}^{\alpha\beta} = i_{T_{0i}^\alpha}\left( w_{i|\mu}^\beta \right) + i_{T_{i|\mu}^\alpha}\left( w_{0i}^\beta \right)
\end{align}
\end{lemma}

\begin{proof}
$(\ref{te11})-(\ref{t33})$ follows from Lemma \ref{tt5}, and $(\ref{n42})-(\ref{n48})$ follows from Lemma \ref{nn5}. It remains to prove $(\ref{se5})$. In fact, from $(\ref{se6})$, we have
\begin{align*}
i_{T_i^{\alpha(\mu-1)}+ T_{i|\mu}^\alpha }\left( w_i^{\beta(\mu-1)} + w_{i|\mu}^\beta \right)=0 \iff K_{i|\mu}^{\alpha\beta}+ i_{T_{0i}^\alpha}\left(w_{i|\mu}^\beta \right)+ i_{T_{i|\mu}^\alpha}\left(w_{0i}^\beta\right)=0
\end{align*}
\end{proof}

\begin{lemma}\label{ss12}
Under the hypothesis $\mathbb{H}^2\left( M, \mathcal{F}_0^\bullet \right)=0$, we can find $\varphi_\mu, r_{ij|\mu}^{\alpha \beta} , T_{i|\mu}^\alpha, g_{i\alpha\beta|\mu}^\gamma$ and $h_{ij|\mu}^{\alpha\beta}, w_{i|\mu}^\alpha, a_{i_x|\mu}^{\alpha\beta}$ which satisfy $(\ref{te11})-(\ref{t33})$ and $(\ref{n42})-(\ref{n48})$ and additionally $(\ref{se5})$.
\end{lemma}

\begin{proof}
We keep the notations in the proof of Lemma \ref{te47} and Lemma \ref{ne25}.  We note that from $(\ref{t57})$ we have $\sum_{\xi=1}^p r_{0ij}^{\alpha \xi} \Gamma_{j|\mu}^\xi-\Gamma_{i|\mu}^\alpha=\Gamma_{ij|\mu}^\alpha-\sum_{\xi=1}^p \lambda_{ij|\mu}^{\alpha \xi} T_{0i}^\xi$ and from $(\ref{n52})$ we have $\sum_{\beta=1}^q h_{0ij}^{\alpha\beta} C_{j|\mu}^\beta - C_{i|\mu}^\alpha= C_{ij|\mu}^\alpha - \sum_{\beta=1}^q  D_{ij|\mu}^{\alpha\beta} w_{0i}^\beta$. Then from $(\ref{se2})$ we have
\begin{align}  \label{se8}
\sum_{\xi=1}^p \sum_{\gamma=1}^q r_{0ij}^{\alpha\xi} h_{0ij}^{\beta \gamma} \left( K_{j|\mu}^{\xi \gamma} + i_{\Gamma_{j|\mu}^\xi} \left( w_{0j}^\gamma \right) + i_{T_{0j}^\xi}\left( C_{j|\mu}^\gamma \right) \right) - \left(K_{i|\mu}^{\alpha\beta} + i_{\Gamma_{i|\mu}^\alpha}\left( w_{0i}^\beta \right)  + i_{T_{0i}^\alpha}\left( C_{i|\mu}^\beta \right)  \right)=0
\end{align}
We define $Z_{i|\mu}\in \Gamma\left( U_i, \mathcal{A}^{0,0}\left( \mathscr{H}om_{\mathcal{O}_M}\left( \mathcal{N}_{\mathcal{F}_0}^*, \Theta_{\mathcal{F}_0}^* \right) \right) \right)$ by 
\begin{align*}
Z_{i|\mu}:\Gamma\left( U_i, \mathcal{N}_{\mathcal{F}_0}^* \right) &\to \Gamma\left( U_i, \mathcal{A}^{0,0}\left( \Theta_{\mathcal{F}_0}^* \right) \right) \\
 w_{0i}^\beta &\mapsto \left(  T_{0i}^\alpha \mapsto Z_{i|\mu}^\alpha:= K_{i|\mu}^{\alpha\beta} + i_{\Gamma_{i|\mu}^\alpha}\left( w_{0i}^\beta \right)  + i_{T_{0i}^\alpha}\left( C_{i|\mu}^\beta \right)    \right)
\end{align*}
Then $(\ref{se8})$ implies that
\begin{align*}
Z_\mu:= \left\{ Z_{i|\mu}\right\} \in A^{0,0}\left( M, \mathscr{H}om_{\mathcal{O}_M} \left( \mathcal{N}_{\mathcal{F}_0}^*, \Theta_{\mathcal{F}_0}^* \right) \right)
\end{align*}
Then we claim that
{\Small{\begin{align} \label{d27}
 \left(\overline{B}_\mu, I_\mu , Z_\mu, \bar{\Phi}_\mu,  \overline{\phi}_\mu   , - \xi_\mu \right) \in &   \frac{A^{0,0}\left( M, \mathscr{H}om_{\mathcal{O}_M}\left( \bigwedge^2 \Theta_{\mathcal{F}_0} ,  \Theta_M \right) \right)}{A^{0,0}\left( M, \mathscr{H}om_{\mathcal{O}_M}\left( \bigwedge^2 \Theta_\mathcal{F} , \Theta_{\mathcal{F}_0}  \right) \right)}\bigoplus A^{0,0}\left(M- S, \mathscr{H}om_{\mathcal{O}_M}\left( \mathcal{N}_{\mathcal{F}_0}^* , \tilde{S}^2 \right)\right) \bigoplus A^{0,0} \left(M, \mathscr{H}om_{\mathcal{O}_M}\left(\mathcal{N}_{\mathcal{F}_0}^*  , \Theta_{\mathcal{F}_0}^*  \right)  \right) \\
 &\bigoplus \frac{A^{0,1}\left(M, \mathscr{H}om_{\mathcal{O}_M}\left( \Theta_{\mathcal{F}_0} ,  \Theta_M \right)  \right)}{A^{0,1}\left(M, \mathscr{H}om_{\mathcal{O}_M} \left( \Theta_{\mathcal{F}_0} , \Theta_{\mathcal{F}_0} \right)  \right)} \bigoplus \frac{A^{0,1}\left( M, \mathscr{H}om_{\mathcal{O}_M}\left( \mathcal{N}_{\mathcal{F}_0}^*  , \Omega_M^1 \right) \right)}{A^{0,1}\left( M,  \mathscr{H}om_{\mathcal{O}_M} \left( \mathcal{N}_{\mathcal{F}_0}^* , \mathcal{N}_{\mathcal{F}_0}^* \right) \right)} \bigoplus A^{0,2}\left(M, \Theta_M\right) \notag
\end{align}}}
defines a $2$-cocycle in the following Dolbeault type bicomplex associated to $\mathcal{F}_0^\bullet$:

{\Tiny{\begin{equation}\label{st1}
\begin{CD}
\cdots\\
@A\hat{F}_3AA \\
\frac{A^{0,0}\left(M, \bigwedge^3 \Theta_{\mathcal{F}_0}^*\otimes \Theta_M\right)}{A^{0,0}\left(M,\bigwedge^3 \Theta_{\mathcal{F}_0}^*\otimes \Theta_{\mathcal{F}_0}\right)}\bigoplus A^{0,0}\left(M - S, \left(\mathcal{N}_{\mathcal{F}_0}^*\right)^*\otimes \tilde{\mathcal{S}}^3 \right)\bigoplus A^{0,0}\left(M-S, \left(\mathcal{N}_{\mathcal{F}_0}^* \right)^*\otimes \bigwedge^2 \Theta_\mathcal{F}^* \right) @> - \bar{\partial} >> \cdots \\
@A \hat{F}_2AA \\
\frac{A^{0,0}\left(M, \bigwedge^2 \Theta_{\mathcal{F}_0}^*\otimes \Theta_M\right)}{A^{0,0}\left(M,\bigwedge^2 \Theta_{\mathcal{F}_0}^*\otimes \Theta_{\mathcal{F}_0}\right)}\bigoplus A^{0,0}\left(M - S, \left(\mathcal{N}_{\mathcal{F}_0}^*\right)^*\otimes \tilde{\mathcal{S}}^2 \right)\bigoplus A^{0,0}\left(M, \left(\mathcal{N}_{\mathcal{F}_0}^* \right)^*\otimes \Theta_\mathcal{F}^* \right)@>\bar{\partial}>> \cdots \\
@A \hat{F}_1AA @A \hat{F}_1AA \\
\frac{A^{0,0}\left(M, \Theta_{\mathcal{F}_0}^*\otimes \Theta_M\right)}{A^{0,0}\left(M, \Theta_{\mathcal{F}_0}^*\otimes \Theta_{\mathcal{F}_0} \right)}\bigoplus \frac{A^{0,0}\left(M,\left(\mathcal{N}_{\mathcal{F}_0}^*\right)^*\otimes \Omega_M^1\right)}{A^{0,0}\left( M, \left(\mathcal{N}_{\mathcal{F}_0}^* \right)^* \otimes \mathcal{N}_{\mathcal{F}_0}^*\right)} @>-\bar{\partial}>>\frac{A^{0,1}\left(M, \Theta_{\mathcal{F}_0}^*\otimes \Theta_M\right)}{A^{0,1}\left(M, \Theta_{\mathcal{F}_0}^*\otimes \Theta_{\mathcal{F}_0} \right)}\bigoplus \frac{A^{0,1}\left(M,\left(\mathcal{N}_{\mathcal{F}_0}^*\right)^*\otimes \Omega_M^1\right)}{A^{0,1}\left( M, \left(\mathcal{N}_{\mathcal{F}_0}^* \right)^* \otimes \mathcal{N}_{\mathcal{F}_0}^*\right)}@>\bar{\partial}>>\cdots \\
@A\hat{F}_0 AA @A \hat{F}_0 AA  \\
A^{0,0}\left(M, \Theta_M \right) @>\bar{\partial}>> A^{0,1} \left(M, \Theta_M\right) @>-\bar{\partial}>> \cdots 
\end{CD}
\end{equation}}}

From the proof of Lemma \ref{te47} and Lemma \ref{ne25}, it remains to show that $\beta_2\left( \overline{B_\mu }\right) + \gamma_2\left( I_\mu \right)- \hat{E}_1' \left( Z_\mu \right)=0$ on $U_i-S$ and $\bar{\partial} Z_\mu + \hat{\beta}_1\left( \overline{\Phi}_\mu \right)+ \hat{\gamma}_1\left( \overline{\phi_\mu} \right)=0$ on $U_i$. We show that $\beta_2\left( \overline{B_\mu }\right) + \gamma_2\left( \overline{I_\mu} \right)- \hat{E}_1\left( Z_\mu \right)=0$ on $U_i-S$. In fact,  first we note that we can write 
\begin{align*}
i_{T_{0i}^\alpha}\left( dw_{0i}^\gamma \right)= \sum_{\eta=1}^q a_{0i\alpha}^{\gamma\eta} w_{0i}^\eta\,\,\,\,\,\,\,\,\,\textnormal{for some}\,\,\,a_{0i\alpha}^{\gamma \eta} \in \Gamma\left( U_i, \mathcal{O}_M \right)
\end{align*}
and since $dw_{0i}^\gamma= \sum_{\eta=1}^q a_{0i_x}^{\gamma \eta} \wedge w_{0i}^\eta$ on $U_{i_x}$ for $x\in U_i-S$, we see that $i_{T_{0i}^\alpha}\left(a_{0i_x}^{\gamma \eta } \right) = a_{0i\alpha}^{\gamma \eta}$. Then we have on $U_{i_x}$
\begin{align}\label{ss15}
&i_{\Pi_{i|\mu}^{\alpha\beta}}\left(w_{0i}^\gamma\right)\equiv_\mu i_{\left[T_i^{\alpha(\mu-1)},T_i^{\beta(\mu-1)}\right] - \sum_{\eta=1}^p g_{i\alpha\beta}^{\eta (\mu-1)} T_i^{\eta(\mu-1)}}\left( w_i^{\gamma(\mu-1)}\right)\\
&\equiv_\mu \mathcal{L}_{T_i^{\alpha(\mu-1)}} i_{T_i^{\beta(\mu-1)}}\left(w_i^{\gamma(\mu-1)} \right) - i_{T_i^{\beta(\mu-1)}} \mathcal{L}_{T_i^{\alpha(\mu-1)}}\left( w_i^{\gamma(\mu-1)} \right) - \sum_{\eta=1}^p g_{i\alpha\beta}^{\eta(\mu-1)} i_{T_i^{\eta(\mu-1)}}\left( w_i^{\gamma(\mu-1)}  \right) \notag \\
&\equiv_\mu  i_{T_{0i}^\alpha} \left(\partial K_{i|\mu}^{\beta\gamma} \right) - i_{T_i^{\beta(\mu-1)}}\left(i_{T_i^{\alpha(\mu-1)}}\partial + \partial i_{T_i^{\alpha(\mu-1)}} \right)\left( w_i^{\gamma(\mu-1)} \right) - \sum_{\eta=1}^p g_{0i\alpha\beta}^\eta K_{i|\mu}^{\eta \gamma}\notag\\
&\equiv_\mu i_{T_{0i}^\alpha}\left( \partial K_{i|\mu}^{\beta\gamma} \right) -i_{T_i^{\beta(\mu-1)}} i_{T_i^{\alpha(\mu-1)}}\left(\partial w_i^{\gamma(\mu-1)} \right) - i_{T_{0i}^\beta}\left( \partial K_{i|\mu}^{\alpha\gamma} \right) - \sum_{\eta=1}^p  g_{0i\alpha\beta}^\eta K_{i|\mu}^{\eta\gamma} \notag \\
&\equiv_\mu i_{T_{0i}^\alpha}\left( \partial K_{i|\mu}^{\beta\gamma} \right) -i_{T_i^{\beta(\mu-1)}} i_{T_i^{\alpha(\mu-1)}}\left( L_{i_x|\mu}^\gamma + \sum_{\eta=1}^q a_{i_x}^{\gamma\eta (\mu-1)} \wedge w_i^{\eta(\mu-1)}\right) - i_{T_i^\beta}\left( \partial K_{i|\mu}^{\alpha\gamma} \right) - \sum_{\eta=1}^p  g_{0i\alpha\beta}^\eta K_{i|\mu}^{\eta\gamma} \notag\\
&\equiv_\mu i_{T_{0i}^\alpha}\left( \partial K_{i|\mu}^{\beta\gamma} \right) -i_{T_{0i}^\alpha \wedge T_{0i}^\beta} \left( L_{i_x|\mu}^\gamma \right) - \sum_{\eta=1}^q i_{T_i^{\beta(\mu-1)}} \left(i_{T_i^{\alpha(\mu-1)}}\left(a_{i_x}^{\gamma\eta (\mu-1)} \right) w_i^{\eta(\mu-1)} - K_{i|\mu}^{\alpha\eta}  a_{0i_x}^{\gamma\eta}\right) - i_{T_{0i}^\beta}\left( \partial K_{i|\mu}^{\alpha\gamma} \right) - \sum_{\eta=1}^p  g_{0i\alpha\beta}^\eta K_{i|\mu}^{\eta\gamma} \notag\\
&\equiv_\mu i_{T_{0i}^\alpha}\left( \partial K_{i|\mu}^{\beta\gamma} \right) -i_{T_{0i}^\alpha \wedge T_{0i}^\beta} \left( L_{i_x|\mu}^\gamma \right) - \sum_{\eta=1}^q  \left(i_{T_{0i}^{\alpha}}\left(a_{0i_x}^{\gamma\eta} \right) K_{i|\mu}^{\beta\eta} - K_{i|\mu}^{\alpha\eta}  i_{T_{0i}^\beta}\left(a_{i0_x}^{\gamma\eta}\right)\right) - i_{T_{0i}^\beta}\left( \partial K_{i|\mu}^{\alpha\gamma} \right) - \sum_{\eta=1}^p  g_{0i\alpha\beta}^\eta K_{i|\mu}^{\eta\gamma} \notag \\
&\equiv_\mu i_{T_{0i}^\alpha}\left( \partial K_{i|\mu}^{\beta\gamma} \right) -i_{T_{0i}^\alpha \wedge T_{0i}^\beta} \left( L_{i_x|\mu}^\gamma \right) - \sum_{\eta=1}^q  \left( a_{0i\alpha}^{\gamma \eta} K_{i|\mu}^{\beta\eta} -  a_{0i\beta}^{\gamma \eta} K_{i|\mu}^{\alpha\eta} \right) - i_{T_{0i}^\beta}\left( \partial K_{i|\mu}^{\alpha\gamma} \right) - \sum_{\eta=1}^p  g_{0i\alpha\beta}^\eta K_{i|\mu}^{\eta\gamma} \notag
\end{align}
On the other hand, we note that $E_{i|\mu}^\gamma= w_{0i}^1 \wedge \cdots \wedge w_{0i}^q \wedge L_{i_x|\mu}^\gamma$ on $U_{i_x}$. Then we have $i_{T_{0i}^\alpha \wedge T_{0i}^\beta}\left(E_{i|\mu}^\gamma \right)= i_{T_{0i}^\alpha \wedge T_{0i}^\beta}\left(L_{i_x|\mu}^\gamma \right) w_{0i}^1 \wedge \cdots \wedge w_{0i}^q$. Then for another $y\in U_i-S$ we have on $U_{i_x}\cap U_{i_y}$, we have
\begin{align*} 
\left(i_{T_{0i}^\alpha \wedge T_{0i}^\beta} \left(L_{i_x|\mu}^\gamma\right) - i_{T_{0i}^\alpha\wedge T_{0i}^\beta}\left( L_{i_y|\mu}^\gamma \right) \right)w_{0i}^1\wedge \cdots \wedge w_{0i}^q = 0 \Longrightarrow  i_{T_{0i}^\alpha \wedge T_{0i}^\beta} \left(L_{i_x|\mu}^\gamma\right) = i_{T_{0i}^\alpha\wedge T_{0i}^\beta}\left( L_{i_y|\mu}^\gamma \right)
\end{align*}
This implies that $L_{i\alpha\beta |\mu }^\gamma:=\left\{i_{T_i^\alpha \wedge T_i^\beta}\left(L_{i_x|\mu}^\gamma\right)| x\in U_i-S \right\}$ is well-defined on $U_i-S$. Then from $(\ref{ss15})$, we have on $U_i-S$
\begin{align}\label{ss16}
i_{\Pi_{i|\mu}^{\alpha\beta}}\left( w_{0i}^\gamma \right) + L_{i\alpha\beta |\mu}^\gamma = i_{T_{0i}^\alpha}\left( \partial K_{i|\mu}^{\beta\gamma} \right) - \sum_{\eta=1}^q  \left( a_{0i\alpha}^{\gamma \eta} K_{i|\mu}^{\beta\eta} -  a_{0i\beta}^{\gamma \eta} K_{i|\mu}^{\alpha\eta} \right) - i_{T_{0i}^\beta}\left( \partial K_{i|\mu}^{\alpha\gamma} \right) - \sum_{\eta=1}^p  g_{0i\alpha\beta}^\eta K_{i|\mu}^{\eta\gamma}
\end{align}
This implies that $\beta_1\left(  \bar{B}_{i|\mu} \right) +\gamma_1\left(I_{i |\mu} \right)= \hat{E}_1'\left( Z_{i|\mu} \right) $ on $U_i-S$. Next we show that $\bar{\partial} Z_\mu + \hat{\beta}_1\left(\bar{\Phi}_\mu \right)+ \hat{\gamma}_1\left( \overline{\phi}_\mu \right)=0$. In fact, it follows from $(\ref{se4})$. Hence $(\overline{B_\mu}, I_\mu, Z_\mu, \overline{\Phi_\mu}, \overline{\phi_\mu}, -\xi_\mu )$ defines a $2$-cocycle in the above Dolbeault type bicomplex associated to $\mathcal{F}_0^\bullet$. The by hypothesis $\mathbb{H}^2\left( M, \mathcal{F}_0^\bullet \right)=0$, there exists
\begin{align*}
\left( \overline{T_\mu'}, \overline{w_\mu'}, \varphi_\mu' \right) \in \frac{A^{0,0}\left( M,  \mathscr{H}om_{\mathcal{O}_M}\left( \Theta_{\mathcal{F}_0}, \Theta_M \right) \right)}{A^{0,0}\left( M ,  \mathscr{H}om_{\mathcal{O}_M }\left( \Theta_{\mathcal{F}_0}, \Theta_{\mathcal{F}_0} \right) \right)} \bigoplus \frac{A^{0,0}\left( M, \mathscr{H}om_{\mathcal{O}_M}\left(  \mathcal{N}_{\mathcal{F}_0}^*, \Omega_M^1 \right) \right)}{A^{0,0}\left( M, \mathscr{H}om_{\mathcal{O}_M}\left( \mathcal{N}_{\mathcal{F}_0}^*, \mathcal{N}_{\mathcal{F}_0}^* \right) \right)} \bigoplus A^{0,1}\left( M, \Theta_M \right)
\end{align*}
such that $(\ref{te30}), (\ref{te31}), (\ref{te32})$ and $(\ref{n55}),(\ref{n56})$ and additionally 
\begin{align}\label{s13}
\beta_1 \left(\overline{T_\mu'}\right)+  \gamma_1 \left( \overline{ w_\mu' } \right) = Z_\mu =\left\{ w_{0i}^\beta \mapsto \left(  T_{0i}^\alpha \mapsto K_{i|\mu}^{\alpha\beta} + i_{\Gamma_{i|\mu}^\alpha}\left( w_{0i}^\beta \right) + i_{T_{0i}^\alpha}\left( C_{i|\mu}^\beta \right)  \right) \right\}
\end{align}
In the proof of Lemma \ref{te47} and Lemma \ref{ne25}, we have already shown that $\varphi_\mu, T_{i|\mu}, r_{ij|\mu}^{\alpha\beta}$ and $w_{i|\mu}, h_{ij|\mu}^{\alpha\beta}$ and $a_{i_x|\mu}^{\alpha\beta}$ satisfying $(\ref{te11})-(\ref{t33})$ and $(\ref{n42})-(\ref{n48})$. It remains to show $(\ref{se5})$. In fact, we recall that from $(\ref{te37})$ and $(\ref{ne19})$, we have
\begin{align*}
T_{i|\mu}^\alpha&= \Gamma_{i|\mu}^\alpha - T_{i|\mu}'+ \sum_{\xi=1}^p c_{i|\mu}^{\alpha \xi} T_{0i}^\xi \\
w_{i|\mu}^\alpha&=  C_{i|\mu}^\alpha - w_{i|\mu}'^\alpha + \sum_{\beta=1}^q P_{i|\mu}^{\alpha\beta} w_{0i}^\beta
\end{align*}
Then we have
\begin{align*}
K_{i|\mu}^{\alpha\beta}+ i_{T_{0i}^\alpha}\left(w_{i|\mu}^\beta\right)+ i_{T_{i|\mu}^\alpha}\left(w_{0i}^\beta \right)&=K_{i|\mu}^{\alpha\beta} + i_{T_{0i}^\alpha}\left(  C_{i|\mu}^\beta - w_{i|\mu}'^\beta + \sum_{\gamma=1}^q P_{i|\mu}^{\beta\gamma} w_{0i}^\gamma \right) + i_{\Gamma_{i|\mu}^\alpha - T_{i|\mu}'^\alpha+ \sum_{\xi=1}^p c_{i|\mu}^{\alpha \xi} T_{0 i }^\xi }\left( w_{0i}^\beta \right)\\
   &= K_{i|\mu}^{\alpha\beta}+ i_{T_{0i}^\alpha}\left( C_{i | \mu}^\beta \right) + i_{\Gamma_{i|\mu}^\alpha}\left( w_{0i}^\beta \right) - \left( i_{T_{0 i }^\alpha} \left(w_{i|\mu}'^\beta \right) + i_{T_{i|\mu}'^\alpha}\left(  w_{0i}^\beta \right) \right) = 0
\end{align*}
which follows from $(\ref{s13})$. This completes the proof of Lemma \ref{ss12}.

\end{proof}
It remains to determine $\varphi_1, T_{i}^{\alpha 1}, r_{ij}^{\alpha\beta 1}, g_{i\alpha\beta}^{\gamma 1}$ and $w_i^{\alpha 1}, h_{ij}^{\alpha 1}, a_{i_x}^{\alpha\beta 1}$ satisfying $(\ref{t13})_1-(\ref{te6})_1$ and $(\ref{n23})_1-(\ref{n241})_1$ and $(\ref{se3})_1$. Given $\dim_\mathbb{C} \mathbb{H}^1\left( M, \mathcal{F}_0^\bullet \right)=r$, we can find a basis of $\mathbb{H}^1\left( M, \mathcal{F}_0^\bullet \right)$ by using the Dolbeualt type bicomplex associated to $\mathcal{F}_0^\bullet$ from $(\ref{st1})$ and represent the basis
\begin{align*}
\left(\overline{\psi_\lambda},  \overline{\kappa_\lambda} , \rho_\lambda \right) \in \frac{A^{0,0}\left(M, \mathscr{H}om_{\mathcal{O}_M}\left(  \Theta_{\mathcal{F}_0}, \Theta_M \right) \right)}{A^{0,0}\left( M, \mathscr{H}om_{\mathcal{O}_M}\left( \Theta_{\mathcal{F}_0}, \Theta_{\mathcal{F}_0} \right) \right)}\bigoplus \frac{A^{0,0}\left( M, \mathscr{H}om_{\mathcal{O}_M}\left( \mathcal{N}_{\mathcal{F}_0}^*, \Omega_M^1 \right) \right)}{A^{0,0}\left( M, \mathscr{H}om_{\mathcal{O}_M}\left( \mathcal{N}_{\mathcal{F}_0}^*, \mathcal{N}_{\mathcal{F}_0}^*  \right) \right)} \bigoplus A^{0,1}\left( M, \Theta_M \right),\,\,\,\,\,\,\,\,\lambda=1,...,r.
\end{align*}
where $\psi_\lambda$ from $(\ref{st2})$ satisfying $(\ref{st3})-(\ref{tp3})$ and $\kappa_\lambda$ from $(\ref{st7})$ satisfying $(\ref{nnc6})$ and $(\ref{nnc3})$. We have also the equality
\begin{align}\label{st15}
i_{T_{0i}^\alpha}\left( \kappa_{i\lambda}^\beta \right) + i_{\psi_{i\lambda}^\alpha}\left( w_{0i}^\beta \right)=0
\end{align}

Then we set $\varphi_1$ as in $(\ref{st3})$, $T_i^{\alpha1}$ as in $(\ref{st4})$, and $r_{ij}^{\alpha\beta 1}$ as in $(\ref{st5})$, and $g_{i\alpha\beta}^{\gamma 1}$ as in $(\ref{st6})$. We also set $w_i^{\alpha 1}$ as in $(\ref{st8})$, and $h_{ij}^{\alpha\beta 1}$ as in $(\ref{st9})$ and $a_{i_x}^{\alpha\beta 1}$ as in $(\ref{st10})$. Then $(\ref{t13})_1-(\ref{te6})_1$ and $(\ref{n23})_1 - (\ref{n241})_1$ are satisfied. It remains to check $(\ref{se3})_1$. In fact,
\begin{align*}
i_{T_{i}^{\alpha 1}}\left( w_{i}^{\beta 1} \right)= i_{T_{0i}^\alpha + \sum_{\lambda=1}^r t_\lambda \left( \psi_{i\lambda}^\alpha + W_{i\lambda}^\alpha \right)}\left( w_{0i}^\beta + \sum_{\eta=1}^r t_\eta \left( \kappa_{i\lambda}^\beta + \mathfrak{K}_{i\lambda}^\beta \right)\right) \equiv_1  \sum_{\lambda=1}^r t_\lambda \left( i_{\psi_{i\lambda}^\alpha}\left( w_{0i}^\beta \right) + i_{T_{0i}^\alpha}\left( \kappa_{i\lambda}^\beta \right) \right)=0
\end{align*}
This completes the inductive construction of $\varphi, T_i^\alpha, r_{ij}^{\alpha\beta}, g_{i\alpha\beta}^{\gamma}$ and $ w_i^\alpha, h_{ij}^{\alpha\beta}, a_{i_x}^{\alpha\beta}$ satisfying $(\ref{t10})-(\ref{te5}),(\ref{n22})-(\ref{n211})$ and $(\ref{se1})$.

\subsection{Proof of convergence}\

We will prove that $\varphi=\sum_{\mu=1}^\infty \varphi_\mu, T_i^\alpha = \sum_{\mu=0}^\infty T_{i|\mu}^\alpha, r_{ij}^{\alpha\beta}= \sum_{\mu=0}^\infty r_{ij|\mu}^{\alpha\beta}, g_{i\alpha\beta}^\gamma = \sum_{\mu=0}^\infty g_{i\alpha\beta|\mu}^\gamma$, and $w_i^\alpha= \sum_{\mu=0}^\infty w_{i|\mu}^\alpha, h_{ij}^{\alpha\beta}= \sum_{\mu=0}^\infty h_{ij|\mu}^{\alpha\beta}$ converge. Before proceeding discussion, we introduce another complex of sheaves which also controls foliated deformations of $\left( M, \Theta_{\mathcal{F}_0}, \mathcal{N}_{\mathcal{F}_0}^* \right)$ but removes the quotient in the degree $1$. We shall define the following complex of sheaves $\mathcal{E}_{\mathcal{F}_0}^\bullet$:
\begin{equation}\label{d32}
\begin{CD}
\cdots \\
@AF_3AA \\
\mathscr{H}om_{\mathcal{O}_M}\left( \bigwedge^3 \Theta_{\mathcal{F}_0}, \frac{\Theta_M}{\Theta_{\mathcal{F}_0}} \right) \bigoplus \mathscr{H}om_{\mathcal{O}_M}\left( \mathcal{N}_{\mathcal{F}_0}^*, \tilde{\mathcal{S}}^3 \right) \bigoplus \mathscr{H}om_{\mathcal{O}_M}\left( \mathcal{N}_{\mathcal{F}_0}^*, \bigwedge^2 \Theta_{\mathcal{F}_0}^* \right) \\
@AF_2AA \\
\mathscr{H}om_{\mathcal{O}_M}\left( \bigwedge^2 \Theta_{\mathcal{F}_0}, \frac{\Theta_M}{\Theta_{\mathcal{F}_0}} \right) \bigoplus \mathscr{H}om_{\mathcal{O}_M}\left( \mathcal{N}_{\mathcal{F}_0}^*, \tilde{\mathcal{S}}^2 \right) \bigoplus \mathscr{H}om_{\mathcal{O}_M}\left( \mathcal{N}_{\mathcal{F}_0}^*, \Theta_{\mathcal{F}_0}^* \right) \\
@AF_1'AA \\
\mathscr{H}om_{\mathcal{O}_M}\left( \Theta_{\mathcal{F}_0}, \Theta_M \right) \bigoplus \mathscr{H}om_{\mathcal{O}_M}\left( \mathcal{N}_{\mathcal{F}_0}^*, \Omega_M^1 \right) \\
@AF_0'AA \\
\mathcal{E}_{\mathcal{F}_0}
\end{CD}
\end{equation}
where $\mathcal{E}_{\mathcal{F}_0}$ is defined by combining $\mathcal{E}_{\Theta_{\mathcal{F}_0}}$ and $\mathcal{E}_{\mathcal{N}_{\mathcal{F}_0}^*}$ in the following way: for an open covering $\left\{ U_i \right\}$ of $M$ as above, we define $\mathcal{E}_{\mathcal{F}_0}$ locally  on $U_i$ by $\mathcal{E}_{\mathcal{F}_0}|_{U_i} \cong \mathscr{H}om_{\mathcal{O}_M}\left(  \Theta_{\mathcal{F}_0}, \Theta_{\mathcal{F}_0}  \right)|_{U_i}\bigoplus \mathscr{H}om_{\mathcal{O}_M}\left( \mathcal{N}_{\mathcal{F}_0}^*, \mathcal{N}_{\mathcal{F}_0}^* \right)|_{U_i} \bigoplus \Theta_M |_{U_i}$ such that $\left( \phi_i, \psi_i , X  \right)\in \mathscr{H}om_{\mathcal{O}_M}\left(\Theta_{\mathcal{F}_0}, \Theta_{\mathcal{F}_0} \right)|_{U_i}\bigoplus \mathscr{H}om_{\mathcal{O}_M}\left( \mathcal{N}_{\mathcal{F}_0}^*,\mathcal{N}_{\mathcal{F}_0}^* \right)|_{U_i} \bigoplus \Theta_M |_{U_i}$ and $\left( \phi_j, \psi_j , X  \right)\in \mathscr{H}om_{\mathcal{O}_M}\left(\Theta_{\mathcal{F}_0}, \Theta_{\mathcal{F}_0} \right)|_{U_j}$ $\bigoplus \mathscr{H}om_{\mathcal{O}_M}\left( \mathcal{N}_{\mathcal{F}_0}^*,\mathcal{N}_{\mathcal{F}_0}^* \right)|_{U_j} \bigoplus \Theta_M |_{U_j}$ are identified on $U_{ij}$ if $\phi_i\left( T_{0i}^\alpha \right)= \phi_j\left(T_{0i}^\alpha \right)+\sum_{\beta=1}^p \left[X, r_{0ij}^{\alpha\beta} \right] T_{0j}^\beta$ and $\psi_i\left(w_{0i}^\alpha \right)= \psi_j\left( w_{0i}^\alpha \right) + \sum_{\beta=1}^q \left[ X, h_{0ij}^{\alpha\beta}\right] w_{0j}^\beta$. Then we define $F_0':\mathcal{E}_{\mathcal{F}_0}\to \mathscr{H}om_{\mathcal{O}_M}\left( \Theta_{\mathcal{F}_0}, \Theta_M \right) \bigoplus \mathscr{H}om_{\mathcal{O}_M}\left( \mathcal{N}_{\mathcal{F}_0}^*, \Omega_M^1 \right)$ locally on $U_i$ by $F_0'\left( \phi_i, \psi_i, X \right)=\left(D_0'\left(\phi_i, X \right), E_0'\left(\psi_i, X \right) \right)$. We define $F_1'$ is the composition of the natural quotient map with $F_1$. Then we have 
\begin{align} \label{ssc1}
\mathbb{H}^i\left( M, \mathcal{E}_{\mathcal{F}_0}^\bullet \right) \cong \mathbb{H}^i\left( M, \mathcal{F}_0^\bullet \right),\,\,\,\,\,\,\,\,\, i\geq 0
\end{align}
(for the detail on the complex of sheaves $\mathcal{E}_{\mathcal{F}_0}^\bullet$, see Part II). We denote by $\mathcal{A}^{0,p}\left( \mathcal{E}_{\mathcal{F}_0} \right)$ be the sheaf of $C^\infty(0,p)$-forms with coefficients in $\mathcal{E}_{\mathcal{F}_0}$ and denote by $A^0\left( M, \mathcal{E}_{\mathcal{F}_0} \right)$ the global section of $\mathcal{A}^{0,p}\left( \mathcal{E}_{\mathcal{F}_0} \right)$. Then we have the following Dolbeault type bicomplex associated to $\mathcal{E}_{\mathcal{F}_0}^\bullet$:
{\Tiny{\begin{equation}
\begin{CD}
\cdots\\
@A\hat{F}_3AA \\
\frac{A^{0,0}\left(M, \bigwedge^3 \Theta_{\mathcal{F}_0}^*\otimes \Theta_M\right)}{A^{0,0}\left(M,\bigwedge^3 \Theta_{\mathcal{F}_0}^*\otimes \Theta_{\mathcal{F}_0}\right)}\bigoplus A^{0,0}\left(M - S, \left(\mathcal{N}_{\mathcal{F}_0}^*\right)^*\otimes \tilde{\mathcal{S}}^3 \right)\bigoplus A^{0,0}\left(M-S, \left(\mathcal{N}_{\mathcal{F}_0}^* \right)^*\otimes \bigwedge^2 \Theta_\mathcal{F}^* \right) @> - \bar{\partial} >> \cdots \\
@A \hat{F}_2AA \\
\frac{A^{0,0}\left(M, \bigwedge^2 \Theta_{\mathcal{F}_0}^*\otimes \Theta_M\right)}{A^{0,0}\left(M,\bigwedge^2 \Theta_{\mathcal{F}_0}^*\otimes \Theta_{\mathcal{F}_0}\right)}\bigoplus A^{0,0}\left(M - S, \left(\mathcal{N}_{\mathcal{F}_0}^*\right)^*\otimes \tilde{\mathcal{S}}^2 \right)\bigoplus A^{0,0}\left(M, \left(\mathcal{N}_{\mathcal{F}_0}^* \right)^*\otimes \Theta_\mathcal{F}^* \right)@>\bar{\partial}>> \cdots \\
@A \hat{F}_1'AA @A \hat{F}_1'AA \\
A^{0,0}\left(M, \Theta_{\mathcal{F}_0}^*\otimes \Theta_M\right) \bigoplus A^{0,0}\left(M,\left(\mathcal{N}_{\mathcal{F}_0}^*\right)^*\otimes \Omega_M^1\right)@>-\bar{\partial}>> A^{0,1}\left(M, \Theta_{\mathcal{F}_0}^*\otimes \Theta_M\right)\bigoplus A^{0,1}\left(M,\left(\mathcal{N}_{\mathcal{F}_0}^*\right)^*\otimes \Omega_M^1\right)@>\bar{\partial}>>\cdots \\
@A\hat{F}_0' AA @A \hat{F}_0' AA  \\
A^{0,0}\left(M, \mathcal{E}_{\mathcal{F}_0} \right) @>\bar{\partial}>> A^{0,1} \left(M, \mathcal{E}_{\mathcal{F}_0} \right) @>-\bar{\partial}>> \cdots 
\end{CD}
\end{equation}}}
Let us denote the $i$-th cohomology group of the complex associated to the above bicomplex by $\textnormal{F}'^i$. Then $\textnormal{F}'^i\cong \textnormal{F}^i(i\geq 0)$ $\textnormal{(\ref{d26})}$ and so $\textnormal{F}'^i \cong \mathbb{H}^i\left( M, \mathcal{E}_{\mathcal{F}_0}^\bullet \right)$ for $i=0,1,2$.

We reinterpret $(\ref{d27})$ in terms of the above bicomplex associated to $\mathcal{E}_{\mathcal{F}_0}^\bullet$. From $(\ref{d28})$ and $(\ref{d29})$, we have
\begin{align*}
\left( -\xi_\mu, \left\{ \left( \bar{\partial} \Lambda_{i|\mu}, \bar{\partial}B_{i|\mu} \right) \right\} \right) \in A^{0,2}\left( M, \mathcal{E}_{\mathcal{F}_0} \right)
\end{align*}

Then we see that from $(\ref{d27})$
{\Small{\begin{align}\label{sd6}
& \left(\overline{B}_\mu, I_\mu , Z_\mu, \Phi_\mu,  \phi_\mu   , \left(- \xi_\mu , \left\{ \left(  \bar{\partial} \Lambda_{i|\mu}, \bar{\partial} B_{i|\mu} \right\} \right) \right)\right) \\
& \in    \frac{A^{0,0}\left( M, \mathscr{H}om_{\mathcal{O}_M}\left( \bigwedge^2 \Theta_{\mathcal{F}_0} ,  \Theta_M \right) \right)}{A^{0,0}\left( M, \mathscr{H}om_{\mathcal{O}_M}\left( \bigwedge^2 \Theta_\mathcal{F} , \Theta_{\mathcal{F}_0}  \right) \right)}\bigoplus A^{0,0}\left(M- S, \mathscr{H}om_{\mathcal{O}_M}\left( \mathcal{N}_{\mathcal{F}_0}^* , \tilde{S}^2 \right)\right) \bigoplus A^{0,0} \left(M, \mathscr{H}om_{\mathcal{O}_M}\left(\mathcal{N}_{\mathcal{F}_0}^*  , \Theta_{\mathcal{F}_0}^*  \right)  \right) \notag \\
 &\,\,\,\,\,\bigoplus A^{0,1}\left(M, \mathscr{H}om_{\mathcal{O}_M}\left( \Theta_{\mathcal{F}_0} ,  \Theta_M \right)  \right)  \bigoplus  A^{0,1}\left( M, \mathscr{H}om_{\mathcal{O}_M}\left( \mathcal{N}_{\mathcal{F}_0}^*  , \Omega_M^1 \right) \right) \bigoplus A^{0,2}\left(M, \mathcal{E}_{\mathcal{F}_0} \right) \notag
\end{align}}}
defines a $2$-cocycle in the above Dolbeault type bicomplex associated to $\mathcal{E}_{\mathcal{F}_0}^\bullet$. 

We define the \"Holder norms on the sections of $\mathcal{A}^{0,q}\left(\Theta_M \right)$ and $\mathcal{A}^{0,p}\left(\mathscr{H}om_{\mathcal{O}_M}\left( \Theta_{\mathcal{F}_0}, \Theta_M \right) \right)$ and $\mathcal{A}^{0,q}\left( \mathscr{H}om_{\mathcal{O}_M}\left( \bigwedge^2 \Theta_{\mathcal{F}_0}, \Theta_M \right) \right)$ and introduce a harmonic theory as in the proof of convergence of Theorem \ref{tt2}. We also define \" Holder norms on the section of $\mathcal{A}^{0,q}\left( \mathscr{H}om_{\mathcal{O}_M}\left( \mathcal{N}_{\mathcal{F}_0}^*, \Omega_M^1 \right) \right)$ and $\mathcal{A}^{0,q}\left( \mathscr{H}om_{\mathcal{O}_M}\left( \mathcal{N}_{\mathcal{F}_0}^*, \bigwedge^{q+2} \Omega_M^1 \otimes \mathcal{L}_0 \right) \right)$ and introduce a harmonic theory as in the proof of convergence of Theorem \ref{nn2}. 

 We define the H\"older norm $|-|_{k+\alpha}$ (k: an integer $\geq 2$, $0<\alpha< 1$) for sections of $\mathcal{A}^{0,r}\left(\mathscr{H}om_{\mathcal{O}_M}\left( \mathcal{N}_{\mathcal{F}_0}^*, \Theta_{\mathcal{F}_0}^* \right) \right)$ as follows: let $\phi\in \Gamma\left(U_i, \mathcal{A}^{0,r}\left(\mathscr{H}om_{\mathcal{O}_M}\left(  \mathcal{N}_{\mathcal{F}_0}^* ,  \Theta_{\mathcal{F}_0}^* \right) \right) \right)$ and we write 
\begin{align*}
\phi\left(T_i^\beta \right)\left( w_{0i}^\gamma \right)=A_i^{\beta \gamma},\,\,\,\,\,\,\,\,\,\,\,A_i^{\beta\gamma}=\frac{1}{r!}\sum A_{i\mu_1\cdots \mu_r}^{\beta\gamma} (z_i)dz_i^{\mu_1}\wedge \cdots \wedge dz_i^{\mu_r} \in \Gamma\left(U_i, \mathcal{A}^{0,r}\right)
\end{align*}
in terms of local coordinates $\left(z_i^1,..., z_i^n \right)$ and let
\begin{align*}
\left| \phi \right|_{k+\alpha}^{U_i} = \sum_{h=0}^k \sup \left| D_i^h A_{i \mu_1\cdots \mu_r }^{\beta\gamma }(z_i) \right| + \sup \frac{ \left| D_i^kA_{ i \mu_1\cdots \mu_r}^{\beta\gamma }(z_i)- D_i^k A_{i\mu_1\cdots \mu_r}^{\beta\gamma }(y_i) \right| }{\left|z_i-y_i\right|^\alpha}
\end{align*}
where the $``\sup$" is extended over all points $z,y\in U_i$, all indices $\beta, \gamma, \mu_1,...,\mu_r$, and all partial derivatives $D_i^h, D_i^k$ of order $h, k$ with respect to $z_i^1,..., z_i^n, \bar{z}_i^1,..., \bar{z}_i^n$. For $\phi\in A^{0,r}\left( M, \mathscr{H}om_{\mathcal{O}_M}\left( \mathcal{N}_{\mathcal{F}_0}^* , \Theta_{\mathcal{F}_0}^*  \right)     \right)$, we define
\begin{align*}
\left| \phi \right|_{k+\alpha} = \max_i \left| \phi \right|_{k+\alpha}^{U_i} 
\end{align*}

Let $\mathcal{E}_{\mathcal{F}_0}$ be the generalized Atiyah extension associate to $\Theta_{\mathcal{F}_0}$ and $\mathcal{N}_{\mathcal{F}_0}^*$ which is locally of the form
\begin{align*}
\Gamma\left( U_i, \mathcal{E}_{\mathcal{F}_0} \right):= \Gamma\left( U_i, \Theta_M \right) \bigoplus \Gamma\left(  U_i, \mathscr{H}om_{\mathcal{O}_M}\left( \Theta_{\mathcal{F}_0} , \Theta_{\mathcal{F}_0}  \right)\right) \bigoplus \Gamma\left( U_i, \mathscr{H}om_{\mathcal{O}_M}\left( \mathcal{N}_{\mathcal{F}_0}^*, \mathcal{N}_{\mathcal{F}_0}^* \right) \right)
\end{align*}
We define the H\"older norm $|-|_{k+\alpha}$ (k: an integer $\geq 2, 0<\alpha <1$) for sections of $\mathcal{A}^{0,r}\left(\mathcal{E}_{\mathcal{F}_0} \right)$ as follows: let $\xi\in \Gamma(U_i, \mathcal{A}^{0,r}(\mathcal{E}_{\Theta_\mathcal{F}}))$ and we write $\xi=\left(\varphi_i, \phi_i ,\psi_i \right)$ as
\begin{align}
\varphi_i &=\sum_{\lambda=1}^n \varphi_i^\lambda \frac{\partial}{\partial z_i^\lambda},\,\,\,\,\,\,\,\,\varphi_i^\lambda= \frac{1}{r!} \sum \varphi_{i\mu_1\cdots \mu_r}^\lambda(z_i)d\bar{z}_i^{\mu_1}\wedge \cdots \wedge d\bar{z}_i^{\mu_r}\in \Gamma\left(U_i, \mathcal{A}^{0,r}\right)\\
\phi_i\left(T_{0i}^\beta\right)&=\sum_{\gamma=1}^p A_{i}^{\beta\gamma} T_{0i}^\gamma \,\,\,\,\,\,\,\,\,\,\,A_i^{\beta\gamma}=\frac{1}{r!}\sum A_{i\mu_1\cdots \mu_r}^{\beta\gamma} (z_i)dz_i^{\mu_1}\wedge \cdots \wedge dz_i^{\mu_r} \in \Gamma\left(U_i, \mathcal{A}^{0,r}\right) \notag\\
\psi_i\left(w_{0i}^\eta \right)&=\sum_{\gamma=1}^q B_{i}^{\eta \sigma} w_{0i}^\sigma \,\,\,\,\,\,\,\,\,\,\, B_i^{\eta\sigma}=\frac{1}{r!}\sum B_{i\mu_1\cdots \mu_r}^{\eta\sigma} (z_i)dz_i^{\mu_1}\wedge \cdots \wedge dz_i^{\mu_r} \in \Gamma\left(U_i, \mathcal{A}^{0,r}\right) \notag
\end{align} 
in terms of local coordinates $\left( z_i^1,..., z_i^n \right)$ and let
\begin{align}
\left| \xi \right|_{k+\alpha}^{U_i} = \sum_{h=0}^k \sup \left| D_i^h C_{i \mu_1\cdots \mu_r }(z_i) \right| + \sup \frac{ \left| D_i^k C_{i \mu_1\cdots \mu_r}(z_i)- D_i^k C_{\mu_1\cdots \mu_r}(y_i) \right| }{\left|z_i-y_i\right|^\alpha}
\end{align}
where $C_{i\mu_1\cdots \mu_r}= \varphi_{i \mu_1\cdots \mu_r}^\lambda (z_i) $ or $A_{i\mu_1,...\mu_r}^{\beta \gamma}(z_i)$ or $B_{i\mu_1,...\mu_r}^{\eta \sigma}(z_i)$, and the $``\sup$" is extended over all points $z,y\in U_i$, all indices $\lambda, \beta, \gamma, \eta, \sigma, \mu_1,...,\mu_r$, and all partial derivatives $D_i^h, D_i^k$ of order $h, k$ with respect to $z_i^1,..., z_i^n, \bar{z}_i^1,..., \bar{z}_i^n$. For $\xi \in A^{0,r}\left( M, \mathcal{E}_{ \mathcal{F}_0 } \right)$, we define
\begin{align}
\left| \xi \right|_{k+\alpha}=\max_i \left|\xi \right|_{k+\alpha}^{U_i}
\end{align}
We introduce a harmonic theory on $\mathcal{E}_{\mathcal{F}_0}$ as in \cite{Kod05} p.157-p.161.
\begin{align}
\textnormal{We denote by $\tilde{\mathfrak{d}}$ the adjoint operator of $\bar{\partial}$, and $\widetilde{\square}=\tilde{\mathfrak{d}}\bar{\partial} + \bar{\partial} \tilde{\mathfrak{d}}$ and $\tilde{G}$ the Green's operator.}
\end{align}

We recall the notations $(\ref{d30})$ and $(\ref{d31})$. With this preparation, we will show that for a fixed integer $k\geq 2$ and $0<\alpha<1$, the inductive construction of $\varphi , T_i^\alpha, r_{ij}^{\alpha\beta}, g_{i\alpha\beta}^\gamma$ and $h_{ij}^{\alpha\beta}, w_j^\alpha$ in the previous subsection can be carried out in such a way that
\begin{align}\label{dt1}
\left|\varphi \right|_{k+\alpha} \ll A(t) \\
\left| T_i^{\alpha }- T_{0i}^{\alpha} \right|_{k+ \alpha}  \ll  A(t) \\
\left| r_{ij}^{\alpha\beta }- r_{0ij}^{\alpha\beta} \right|_{k+ 1+ \alpha} \ll A(t) \\
\left| g_{i\alpha\beta}^{\gamma} - g_{0i\alpha\beta}^\gamma \right|_{ k -1 +  \alpha} \ll A(t)\\
\left| w_i^{\alpha }- w_{0i}^\alpha \right|_{k+\alpha} \ll A(t)\\
\left| h_{ij}^{\alpha\beta } - h_{0ij}^{\alpha\beta}\right|_{k+1+\alpha} \ll A(t) \label{dt2}
\end{align}
Then it suffices to prove that for $\mu=1,2,3,\cdots$,
\begin{align}
\left| \varphi^{\mu}\right|_{k+\alpha}  \ll A(t) \label{sd1} \\
\left| T_i^{\alpha\mu}- T_{0i}^{\alpha} \right|_{k+ \alpha}  \ll  A(t) \\
\left| r_{ij}^{\alpha\beta \mu}- r_{0ij}^{\alpha\beta} \right|_{k+ 1+ \alpha} \ll A(t) \\
\left| g_{i\alpha\beta}^{\gamma \mu} - g_{0i\alpha\beta}^\gamma \right|_{ k -1+  \alpha} \ll A(t)\\
\left| w_i^{\alpha \mu}- w_{0i}^\alpha \right|_{k+\alpha} \ll A(t)\\
\left| h_{ij}^{\alpha\beta \mu} - h_{0ij}^{\alpha\beta}\right|_{k+1+\alpha} \ll A(t)\label{sd2}
\end{align}
for some proper choice of constants $c>b>0$. We prove $(\ref{sd1})_\mu -(\ref{sd2})_\mu$ by induction on $\mu$. For $\mu=1$ we have $(\ref{st3}),(\ref{st4}),(\ref{st5}),(\ref{st6})$, and $(\ref{ns3}),(\ref{st8}),(\ref{st9})$, and the linear term of $A(t)$ is $\frac{b}{16}\left(t_1+ \cdots  + t_r \right)$. Therefore $(\ref{sd1})_1-(\ref{sd2})_1$ holds if $b$ is sufficiently large.

Now assume that $(\ref{sd1})_{\mu-1}-(\ref{sd2})_{\mu-1}$ are satisfied. We will derive $(\ref{sd1})_\mu - (\ref{sd2})_\mu$. In the following $C_1,C_2,...$ will denote constants which depend only on $k,\alpha, M,\Theta_{\mathcal{F}_0}, \mathcal{N}_{\mathcal{F}_0}^*$. 

By taking $\left( \tilde{\varphi}_\mu, \left\{ \left(\tilde{b}_{i|\mu} , \tilde{\lambda}_{i|\mu} \right)\right\} \right):= \tilde{\mathfrak{d}}\tilde{G} \left( -\xi_\mu, \left\{ \left( \bar{\partial} \Lambda_{i|\mu}, \bar{\partial} B_{i|\mu}\right)\right\} \right)=0$, we have as in the same way with $(\ref{d25})$ and $(\ref{sd5})$
\begin{align}\label{dt3}
\left| \left( \tilde{\varphi}_\mu, \left\{ \left( \tilde{b}_{i|\mu}  , \tilde{\lambda}_{i|\mu}\right) \right\}      \right)   \right|_{k+\alpha} \ll C_1\left| \left( - \xi_\mu,  \left\{ \left( \bar{\partial}\Lambda_{i|\mu}, \bar{\partial} B_{i|\mu} \right)   \right\} \right)  \right|_{k-1+\alpha}
\end{align}
where $C_1$ is a constant independent of $\left(- \xi_\mu, \left\{ \left( \bar{\partial} \Lambda_{i|\mu}, \bar{\partial} B_{i|\mu} \right) \right\} \right)$. As in the same way with $(\ref{d41})$ and $(\ref{ns38})$, we can show that
\begin{align}\label{dt4}
\left( \varphi_\mu'', \left\{\left( b_{i|\mu}'', \lambda_{i|\mu}'' \right) \right\}\right)\in A^{0,1}\left( M, \mathcal{E}_{\mathcal{F}_0} \right),\,\,\,\,\,\,\,\,\mathfrak{d}\varphi_\mu'' = 0,\,\,\,\,\,\,\,\,\bar{\partial}\left( \varphi_\mu'', \left\{ \left( b_{i|\mu}'', \lambda_{i|\mu}''    \right) \right\} \right)=\left( - \xi_\mu, \left\{ \left( \bar{\partial}\Lambda_{i|\mu}, \bar{\partial} B_{i|\mu} \right) \right\}    \right)
\end{align}
and we can show that as in the same way with $(\ref{ds3})$ and $(\ref{ns27})$ 
\begin{align}\label{sd22}
\left| \varphi_\mu'' \right|_{k+\alpha},\,\,\,\,\,\, \left| b_{i|\mu}'' \right|_{k+\alpha}^{U_i}, \,\,\,\,\,\, \left| \lambda_{i|\mu}'' \right|_{k+\alpha}^{U_i} \ll C_2 \left|\left( - \xi_\mu , \left\{ \left( \bar{\partial} \Lambda_{i|\mu} , \bar{\partial} B_{i|\mu}    \right)\right\} \right) \right|_{k-1+\alpha}
\end{align}

From $(\ref{sd6})$, we consider
{\small{\begin{align*}
&0\in A^{0,2}\left( M, \mathcal{E}_{\mathcal{F}_0} \right)\\
&\Phi_\mu'= \left\{ T_{0i}^\alpha \mapsto \Phi_{i|\mu}^\alpha - \bar{\partial} \Gamma_{i|\mu}^\alpha +\sum_{\xi=1}^p \Lambda_{i|\mu}^{\alpha\xi} T_{0i}^\xi +\left[ \varphi_\mu'', T_{0i}^\alpha \right]  - \sum_{\xi=1}^p b_{i|\mu}''^{\alpha  \xi} T_{0i}^\xi   \right\}  \in  A^{0,1}\left(M, \mathscr{H}om_{\mathcal{O}_M}\left(\Theta_{\mathcal{F}_0}, \Theta_M \right) \right) \\
&\phi_\mu'=\left\{ w_{0i}^\alpha \mapsto A_{i|\mu}^\alpha - \bar{\partial} C_{i|\mu}^\alpha + \sum_{\gamma=1}^q B_{i|\mu}^{\alpha \gamma} w_{0i}^\gamma + \mathcal{L}_{\varphi_\mu''}\left( w_{0i}^\alpha \right) - \sum_{\xi=1}^q \lambda_{i|\mu}''^{\alpha \xi} w_{0i}^\xi \right\} \in A^{0,1}\left(M, \mathscr{H}om_{\mathcal{O}_M}\left(  \mathcal{N}_{\mathcal{F}_0}^*, \Omega_M^1 \right) \right) \\
&\overline{B}_\mu-\overline{\left\{ T_{0i}^\alpha \wedge T_{0i}^\beta \mapsto \Pi_{i|\mu}^{\alpha\beta} + \left[ \Gamma_{i|\mu}^\alpha, T_{0i}^\beta \right] - \left[ \Gamma_{i|\mu}^\beta, T_{0i}^\alpha  \right] - \sum_{\gamma=1}^p g_{0i\alpha\beta}^\gamma \Gamma_{i|\mu}^\gamma + \sum_{\gamma=1}^p \Psi_{i|\mu}^{\alpha\beta \gamma} T_{0i}^\gamma    \right\} } \in \frac{ A^{0,0}\left(M, \mathscr{H}om_{\mathcal{O}_M}\left( \bigwedge^2 \Theta_{\mathcal{F}_0}, \Theta_M \right)\right) }{ A^{0,0}\left( M, \mathscr{H}om_{\mathcal{O}_M}\left( \bigwedge^2 \Theta_{\mathcal{F}_0}, \Theta_{\mathcal{F}_0} \right) \right) }\\
&\psi_\mu=\left\{ w_{0i}^\alpha \mapsto  E_{i|\mu}^\alpha + \sum_{\beta=1}^q w_{0i}^1 \wedge \cdots \wedge C_{i|\mu}^\beta \wedge \cdots \wedge w_{0i}^q \wedge dw_{0i}^\alpha + w_{0i}^1 \wedge \cdots \wedge w_{0i}^q \wedge  \partial C_{i|\mu}^\alpha \right\}  \in A^{0,0}\left( M, \mathscr{H}om_{\mathcal{O}_M}\left( \mathcal{N}_{\mathcal{F}_0}^*, \bigwedge^{q+2} \Omega_M^1 \otimes \mathcal{L}_0 \right) \right) \\
& Z_\mu=\left\{  w_{0i}^\beta \mapsto \left( T_{0i}^\alpha \mapsto K_{i|\mu}^{\alpha\beta} + i_{\Gamma_{i|\mu}^\alpha}\left( w_{0i}^\beta\right) + i_{T_{0i}^\alpha}\left( C_{i|\mu}^\beta \right) \right)\right\} \in A^{0,0}\left( M, \mathscr{H}om_{\mathcal{O}_M}\left( \mathcal{N}_{\mathcal{F}_0}^* , \Theta_{\mathcal{F}_0}^\bullet \right) \right)
\end{align*}}}
Then $\left(0, \Phi_\mu',\phi_\mu', \overline{B}_\mu, \psi_\mu|_{M-S}, Z_\mu \right)$ defines a $2$-cocycle in the above Dolbeault type bicomplex associated with $\mathcal{E}_{\mathcal{F}_0}^\bullet$. Since $\mathbb{H}^2\left( M, \mathcal{F}_0^\bullet \right)= \mathbb{H}^2\left( M,  \mathcal{E}_{\mathcal{F}_0}^\bullet \right)=0$ by the assumption of Theorem \ref{ss6}, there exist $\left( \chi_\mu, \left\{ \left( \eta_i^{\chi_\mu}, \delta_i^{\chi_\mu} \right) \right\} \right)\in A^{0,1}\left( M, \mathcal{E}_{\mathcal{F}_0} \right), \sigma_\mu\in A^{0,0}\left( M, \mathscr{H}om_{\mathcal{O}_M}\left( \Theta_{\mathcal{F}_0}, \Theta_M \right) \right)$ and $\Sigma_\mu\in A^{0,0}\left( M, \mathscr{H}om_{\mathcal{O}_M}\left( \mathcal{N}_{\mathcal{F}_0}^* , \Omega_M^1 \right) \right)$ such that we have $-\bar{\partial}\left(\xi_\mu, \left\{\left( \eta_i^{\chi_\mu}, \delta_i^{\chi_\mu} \right) \right\} \right)=0, - \bar{\partial} \left(\sigma_\mu, \Sigma_\mu\right)+ \hat{F}_0'\left( \chi_\mu, \left\{\left(\eta_i^{\chi_\mu}, \delta_i^{\chi_\mu} \right) \right\}\right)=\left( \Phi_\mu', \phi_\mu' \right)$ and $\hat{F}_1'\left( \sigma_\mu, \Sigma_\mu\right)=\left( \overline{B}_\mu , \psi_\mu|_{M-S}, Z_\mu \right)$. By using the following Lemma, we will choose appropriate $\left(\chi_\mu, \left\{ \left( \eta_i^{\chi_\mu}, \delta_i^{\chi_\mu}  \right) \right\} \right), \sigma_\mu, \Sigma_\mu$ in a way that $\varphi, T_i^\alpha, r_{ij}^{\alpha\beta}, g_{i\alpha\beta}^\gamma, w_i^\alpha$ and $h_{ij}^{\alpha\beta}$ converge. By combining Lemma \ref{cvt1} and Lemma \ref{nnc2}, we prove

\begin{lemma} \label{sd18}
Suppose that $\left(\varphi, \left\{ \eta_i\right\}, \left\{ \delta_i \right\} \right) \in A^{0,1}\left( M, \mathcal{E}_{\mathcal{F}_0} \right)$ where $\varphi\in A^{0,1}\left( M, \Theta_M \right), \eta_i \in \Gamma\left( U_i, \mathscr{H}om_{\mathcal{O}_M} \left( \Theta_{\mathcal{F}_0} , \Theta_{\mathcal{F}_0}   \right) \right)$ and $\delta_i \in \Gamma\left( U_i, \mathscr{H}om_{\mathcal{O}_M} \left( \mathcal{N}_{\mathcal{F}_0}^*, \mathcal{N}_{\mathcal{F}_0}^* \right) \right)$ with $\eta_i\left( T_{0i}^\alpha\right)=\eta_j\left( T_{0i}^\alpha \right)+\sum_{\beta=1}^p \left[ \varphi, r_{0ij}^{\alpha\beta}\right] T_{0j}^\beta$ for $\alpha=1,...,p$, and $\delta_i\left(w_{0i}^\alpha \right) = \delta_j \left(w_{0i}^\alpha \right)+ \sum_{\beta=1}^q \left[ \varphi, h_{0ij}^{\alpha\beta} \right] w_{0j}^\beta$ for $\alpha=1,...,q$, and $\Phi\in A^{0,0}\left( M, \mathscr{H}om_{\mathcal{O}_M}\left( \Theta_{\mathcal{F}_0}, \Theta_M \right)\right), V\in A^{0,0}\left( M, \mathscr{H}om_{\mathcal{O}_M} \left( \mathcal{N}_{\mathcal{F}_0}^*, \Omega_M^1 \right) \right)$, and $B\in A^{0,0}\left( M, \mathscr{H}om_{\mathcal{O}_M} \left( \bigwedge^2 \Theta_{\mathcal{F}_0}, \Theta_M \right) \right)$, and  $I \in A^{0,0}\left(M, \mathscr{H}om_{\mathcal{O}_M}\left( \mathcal{N}_{\mathcal{F}_0}^*, \bigwedge^{q+2} \Omega_M^1 \otimes \mathcal{L}_0 \right) \right)$ such that $I_{M-S}\in A^{0,0}\left(M-S, \mathscr{H}om_{\mathcal{O}_M}\left( \mathcal{N}_{\mathcal{F}_0}^*, \tilde{\mathcal{S}}^2 \right) \right)$ and $J\in A^{0,0} \left( M, \mathscr{H}om_{\mathcal{O}_M}\left( \mathcal{N}_{\mathcal{F}_0}^*, \Theta_{\mathcal{F}_0}^* \right) \right)$ such that $\left( 0, \left( \Phi, V\right)+ \hat{F}_1'\left( \varphi, \left\{ \left(\eta_i, \delta_i  \right) \right\} \right) \right.$, $\left. \overline{B}, I_{M-S}, J \right)$ defines a $2$-cocycle in the Dolbeault type bicomplex associated to $\mathcal{E}_{\mathcal{F}_0}^\bullet$, and moreover $(-1)^q\bar{\partial} I + E_1^\sharp(V)=0$. Then we can find $\left( \chi, \left\{ \eta_i^\chi\right\}, \left\{ \delta_i^\chi \right\} \right) \in A^{0,1}\left( M, \mathcal{E}_{\mathcal{F}_0} \right),\sigma\in A^{0,0}\left( M, \mathscr{H}om_{\mathcal{O}_M} \left( \Theta_{\mathcal{F}_0} , \Theta_M \right) \right)$ and the associated element $\pi_{i\sigma}\in \Gamma\left( U_i, \mathscr{H}om_{\mathcal{O}_M}\left( \bigwedge^2 \Theta_{\mathcal{F}_0}, \Theta_{\mathcal{F}_0} \right) \right)$ and $\Sigma\in A^{0,0}\left( M, \mathscr{H}om_{\mathcal{O}_M}\left( \mathcal{N}_{\mathcal{F}_0}^*, \Omega_M^1 \right) \right)$ satisfying $(\ref{tpc13})-(\ref{tpc14})$ and $(\ref{ns11})-(\ref{ns7})$ and in addition
\begin{align}
&i_{\sigma\left( T_{0i}^\alpha \right)}\left( w_{0i}^\beta \right) + i_{T_{0i}^\alpha}\left( \Sigma \left( w_{0i}^\beta \right) \right) = J\left( w_{0i}^\beta\right)\left( T_{0i}^\alpha \right) \label{sd7} \\
\left| \left(\chi, \left\{ \left( \eta_i^\chi, \delta_i^\chi \right) \right\}  \right) \right|_{k+\alpha} \ll C &\left( \left| \left(\varphi, \left\{ \eta_i \right\}, \left\{ \delta_i \right\} \right) \right|_{k+\alpha} + \left| \Phi \right|_{k-1+\alpha} + \left| V \right|_{k-1+\alpha} + \left| B \right|_{k-1+\alpha}+ \left| I \right|_{k-1+\alpha} + \left| J \right|_{k-1+\alpha}    \right)\\
\left| \sigma\right|_{k+\alpha} \ll C &\left( \left| \left(\varphi, \left\{ \eta_i \right\}, \left\{ \delta_i \right\} \right) \right|_{k+\alpha} + \left| \Phi \right|_{k-1+\alpha} + \left| V \right|_{k-1+\alpha} + \left| B \right|_{k-1+\alpha}+ \left| I \right|_{k-1+\alpha} + \left| J \right|_{k-1+\alpha}    \right)\\
\left| \pi_{i\sigma}\right|_{k+\alpha}^{U_i} \ll C &\left( \left| \left(\varphi, \left\{ \eta_i \right\}, \left\{ \delta_i \right\} \right) \right|_{k-1\alpha} + \left| \Phi \right|_{k-1+\alpha} + \left| V \right|_{k-1+\alpha} + \left| B \right|_{k-1+\alpha}+ \left| I \right|_{k-1+\alpha} + \left| J \right|_{k-1+\alpha}    \right)\\
\left| \Sigma  \right|_{k+\alpha} \ll C &\left( \left| \left(\varphi, \left\{ \eta_i \right\}, \left\{ \delta_i \right\} \right) \right|_{k+\alpha} + \left| \Phi \right|_{k-1+\alpha} + \left| V \right|_{k-1+\alpha} + \left| B \right|_{k-1+\alpha}+ \left| I \right|_{k-1+\alpha} + \left| J \right|_{k-1+\alpha}    \right)
\end{align}
where $C$ is a constant which is independent of $\left(\varphi, \left\{ \eta_i \right\}, \left\{ \delta_i \right\} \right), V, I, J$
\end{lemma}

\begin{proof}
For any $\zeta=\left( \left( \varphi, \left\{ \left( \eta_i, \delta_i \right) \right\} \right), V, I, J\right)$ as above, let
\begin{align*}
||\zeta || &= \left| \left(\varphi  , \left\{ \left( \eta_i^\varphi,  \delta_i^\varphi  \right) \right\}\right)\right|_{k+\alpha} + \left| \Phi \right|_{k-1+\alpha} + \left|V \right|_{k-1+\alpha} +\left| B \right|_{k-1+\alpha} + \left| I \right|_{k-1+\alpha} + \left| J \right|_{k-1+\alpha} \\
\iota(\zeta)&= \inf \left( \left|\left(\chi, \left\{ \left( \eta_i^\chi , \delta_i^\chi \right) \right\} \right)\right|_{k+\alpha} +\left| \sigma \right|_{k+\alpha} + \max_i \left| \pi_{i \sigma} \right|_{k-1+\alpha}^{U_i} + \left|\Sigma  \right|_{k+\alpha}  \right)
\end{align*}
where $\inf$ is taken with respect to all solutions $\left( \left(\chi, \left\{\left(\eta_i^\chi , \delta_i^\chi  \right) \right\}\right), \left(\sigma, \left\{ \pi_{i\sigma}\right\}\right), \Sigma \right)$ of the equalities $(\ref{tpc13})-(\ref{tpc14})$ and $(\ref{ns11})-(\ref{ns7})$ and $(\ref{sd7})$. It suffices to prove the existence of $C$ such that
\begin{align*}
\iota\left( \zeta \right)\leq C ||\zeta ||\,\,\,\,\,\,\,\,\,\textnormal{for all $\zeta$}
\end{align*}
Assume that there is no such constant $C$. Then we can find a sequence $\zeta^{(1)},\zeta^{(2)},\cdots, \zeta^{(v)},\cdots$ of quintuple $\zeta^{(v)}=\left( \left( \varphi^{(v)} , \left\{ \eta_i^{(v)} \right\} \right), \Phi^{(v)},  V^{(v)}, I^{(v)}, J^{(v)} \right)$ such that
\begin{align*}
\iota\left(\zeta^{(v)} \right)=1\,\,\,\,\,\,\,\,\textnormal{and}\,\,\,\,\,\,\,\,\left|\left| \zeta^{(v)}\right|\right| < \frac{1}{v}
\end{align*}
The first equality implies the existence of $\left( \chi^{(v)}, \left\{\left( \eta_i^{\chi^{(v)}}, \delta_i^{\chi^{(v)}} \right) \right\} \right)\in A^{0,1}\left( M, \mathcal{E}_{\mathcal{F}_0} \right)$ and $\sigma^{(v)}\in A^{0,0}\left( M, \mathscr{H}om_{\mathcal{O}_M}\left( \Theta_{\mathcal{F}_0}, \Theta_M \right) \right)$ and $\pi_{i\sigma^{(v)}}\in \Gamma\left( U_i, \mathcal{A}^{0,0}\left( M,  \mathscr{H}om_{\mathcal{O}_M}\left(\bigwedge^2 \Theta_{\mathcal{F}_0}, \Theta_{\mathcal{F}_0} \right) \right) \right)$ and $\Sigma^{(v)}\in A^{0,0}\left( M, \mathscr{H}om_{\mathcal{O}_M}\left( \mathcal{N}_{\mathcal{F}_0}^*, \Omega_M^1 \right) \right)$ such that $(\ref{tpc16})-(\ref{sd9})$ and $(\ref{sd10})-(\ref{sd11})$

\begin{align}
&\tilde{\square}\left( \chi^{(v)}, \left\{ \left( \eta_i^{\chi^{(v)}}, \delta_i^{\chi^{(v)}} \right) \right\} \right)=0\\
i_{\sigma^{(v)}\left(T_{0i}^\alpha \right)} & \left( w_{0i}^\beta \right)  + i_{T_{0i}^\alpha }\left( \Sigma^{(v)}\left(w_{0i}^\beta \right)\right) = J^{(v)}\left( w_{0i}^\beta\right)\left( T_{0i}^\alpha \right) \label{sd17}\\
\left| \left( \chi^{(v)}, \left\{ \left( \eta_i^{\chi^{(v)}}, \delta_i^{\chi^{(v)}} \right) \right\} \right)   \right|&_{k+\alpha} +\left| \sigma^{(v)} \right|_{k+\alpha} + \max_i \left| \pi_{i\sigma^{(v)}}\right|_{k-1+\alpha}^{U_i} + \left| \Sigma^{(v)} \right|_{k+\alpha} <2 \label{sd12}
\end{align}
Then $(\ref{sd12})$ implies that we may assume that $ \left( \chi, \left\{ \left( \eta_i^{\chi}, \delta_i^{\chi} \right\} \right) \right)=\lim \left(  \chi^{(v)}, \left\{ \left( \eta_i^{\chi^{(v)}}, \delta_i^{\chi^{(v)}} \right) \right\} \right)$ and $\Sigma= \lim \Sigma^{(v)}$ exist in the norm $|-|_k$, and $\sigma=\lim \sigma^{(v)}$ exists in the norm $\left|- \right|_k$ and $\pi_{i\sigma}=\lim \pi_{i\sigma^{(v)}}$ exists in the norm $|-|_{k-1}$, so that $\left( \chi, \left\{ \left( \eta_i^\chi, \delta_i^\chi \right) \right\} \right)$  and $\sigma, \Sigma$ are of class $C^k$, and $\pi_{i\sigma}$ are of class $C^k$. We note that since $\left( \chi, \left\{ \left( \eta_i^{\chi} , \delta_i^{\chi}  \right) \right\} \right)$ satisfies an elliptic partial differential equation $\widetilde{\square}\left( \chi, \left\{ \left( \eta_i^\chi, \delta_i^\chi \right) \right\} \right)=0$, $\left(\chi, \left\{ \left( \eta_i^\chi, \delta_i^\chi \right) \right\} \right)$ is of $C^\infty$. Moreover since $\widetilde{\square}$ is strongly elliptic, by using `a priori estimate' (see \cite{Kod05} Theorem 4.3 p.436) we can show that $\left( \chi^{(v)}, \left\{ \left( \eta_i^{\chi^{(v)}}, \delta_i^{\chi^{(v)}} \right)  \right\} \right)$ converges to $\left(\chi, \left\{ \left( \eta_i^\chi , \delta_i^\chi \right) \right\} \right)$ in the norm $|-|_{k+\alpha}$ in the same way with $(\ref{ns10})$. On the other hand, as in the proof of Lemma \ref{cvt1} and Lemma \ref{nnc2}, $\sigma$ is $C^\infty$ and $\sigma^{(v)}$ converge to $\sigma$ in $\left|-\right|_{k+\alpha}$, and $\pi_{i\sigma}$ is of $C^\infty$ and $\pi_{i\sigma^{(v)}}$ converge to $\pi_{i\sigma}$ in $\left|-\right|_{k-1+\alpha}$, and $\Sigma$ is $C^\infty$ and $\Sigma^{(v)}$ converge to $\Sigma$ in $\left| - \right|_{k+\alpha}$, and we have $(\ref{sd15}),(\ref{sd16}),(\ref{sd13}),(\ref{sd14})$ and in addition, from $(\ref{sd17})$
\begin{align*}
i_{\left( \sigma^{(v)}-\sigma\right)\left( T_{0i}^\alpha \right) }\left( w_{0i}^\beta \right) + i_{T_{0i}^\alpha}\left( \left(\Sigma^{(v)}- \Sigma \right) \left( w_{0i}^\beta \right) \right) = J^{(v)} \left( w_{0i}^\beta \right) \left( T_{0i}^\alpha \right)
\end{align*}

On the other hand, we have
{\small{\begin{align*}
\left|  \left( \chi^{(v)}, \left\{ \left( \eta_i^{\chi^{(v)}} , \delta_i^{\chi^{(v)}} \right) \right\} \right) - \left( \chi , \left\{ \left( \eta_i^\chi , \delta_i^\chi \right) \right\} \right)\right|_{k+\alpha} <\frac{1}{5},\,\,\,\,\,\,\,\,\,\left| \sigma^{(v)}- \sigma \right|_{k+\alpha} <\frac{1}{5},\,\,\,\,\,\,\,\,\,\,\, \left| \pi_{i\sigma^{(v)}} -  \pi_{i\sigma} \right|_{k-1+\alpha} < \frac{1}{5},\,\,\,\,\,\,\,\,\left| \Sigma^{(v)}- \Sigma \right|_{k+\alpha} < \frac{1}{5}
\end{align*}}}
for sufficiently large integer $v$. This contradicts to $\iota\left( \zeta \right)=1$. This completes the proof of Lemma \ref{sd18}.

\end{proof}

In Lemma \ref{sd18}, we set
\begin{align}\label{sd21}
& \left( \varphi, \left\{ \left( \eta_i  , \delta_i \right) \right\} \right):= \left( \varphi_\mu'', \left\{  \left( b_{i|\mu}'' , \lambda_{i|\mu}'' \right) \right\} \right)\in A^{0,1}\left( M, \mathcal{E}_{ \mathcal{F}_0 } \right)\\
& \Phi:= \left\{ T_{0i}^\alpha \mapsto \Phi_{i|\mu}^\alpha - \bar{\partial} \Gamma_{i|\mu}^\alpha +\sum_{\xi=1}^p \Lambda_{i|\mu}^{\alpha\xi} T_{0i}^\xi   \right\}  \in  A^{0,1}\left(M, \mathscr{H}om_{\mathcal{O}_M}\left(\Theta_{\mathcal{F}_0}, \Theta_M \right) \right) \notag \\
&B:= \left\{ T_{0i}^\alpha \wedge T_{0i}^\beta \mapsto \Pi_{i|\mu}^{\alpha\beta} + \left[ \Gamma_{i|\mu}^\alpha, T_{0i}^\beta \right] - \left[ \Gamma_{i|\mu}^\beta, T_{0i}^\alpha  \right] - \sum_{\gamma=1}^p g_{0i\alpha\beta}^\gamma \Gamma_{i|\mu}^\gamma + \sum_{\gamma=1}^p \Psi_{i|\mu}^{\alpha\beta \gamma} T_{0i}^\gamma    \right\}  \in A^{0,0}\left(M, \mathscr{H}om_{\mathcal{O}_M}\left( \bigwedge^2 \Theta_{\mathcal{F}_0}, \Theta_M \right)\right)  \notag\\
V&:= \phi_\mu'=\left\{ w_{0i}^\alpha \mapsto  A_{i|\mu}^\alpha - \bar{\partial} C_{i|\mu}^\alpha + \sum_{\gamma=1}^q B_{i|\mu}^{\alpha \gamma} w_{0i}^\gamma \right\}\notag \\
I&:= \psi_\mu=\left\{ w_{0i}^\alpha \mapsto E_{i|\mu}^\alpha + \sum_{\beta=1}^q w_{0i}^1 \wedge \cdots \wedge C_{i|\mu}^\beta \wedge \cdots \wedge w_{0i}^q \wedge dw_{0i}^\alpha + w_{0i}^1 \wedge \cdots \wedge w_{0i}^q \wedge \partial C_{i|\mu}^\alpha \right\}. \notag\\
J&:=Z_\mu=\left\{ w_{0i}^\beta \mapsto \left( T_{0i}^\alpha \mapsto K_{i|\mu}^{\alpha\beta}+ i_{\Gamma_{i|\mu}^\alpha}\left( w_{0i}^\beta \right)  + i_{T_{0i}^\alpha} \left( C_{i|\mu}^\beta \right) \right) \right\} \notag
\end{align}
We will estimate $\left(\varphi, \left\{ \left( \eta_i , \delta_i \right) \right\} \right),\Phi,B, V, I, J$. First we will estimate $(\ref{tt3})-(\ref{ds47})$ and $(\ref{nn3})-(\ref{ns14})$ and $(\ref{se6})$. By induction hypothesis $(\ref{sd1})_{\mu-1}-(\ref{sd2})_{\mu-1}$ and from $(\ref{ds22})-(\ref{sd20})$ we may assume that
\begin{align}\label{sd23}
\left| \xi_\mu \right|_{k-1+\alpha} , \left| \Phi_{i|\mu}^\alpha \right|_{k-1+\alpha},\left| \Lambda_{ij|\mu}^{\alpha\beta} \right|_{k+\alpha} ,\left| \Gamma_{ij|\mu}^\alpha \right|_{k+\alpha},  \left| \Pi_{i|\mu}^{\alpha\beta} \right|_{k-1+\alpha}, \left| \lambda_{ijk|\mu}^{\alpha\beta} \right|_{k+1+\alpha}, \left| \Psi_{ij|\mu}^{\alpha\beta \xi} \right|_{k-1+\alpha} \ll C_3 \frac{b}{c} A(t)
\end{align}
and from $(\ref{ns16})-(\ref{ns20})$ we may assume that
\begin{align*}
\left| A_{i|\mu}^\alpha \right|_{k-1+\alpha},  \left| B_{ij|\mu}^{\alpha\eta}\right|_{k+\alpha},  \left| C_{ij|\mu}^\alpha \right|_{k+\alpha}, \left| D_{ijk|\mu}^{\alpha\beta} \right|_{k+1+\alpha} , \left| E_{i|\mu}^\alpha \right|_{k-1+\alpha} \ll C_4 \frac{b}{c} A(t)
\end{align*}
It remains to estimate
\begin{align}
&K_{i|\mu}^{\alpha\beta}= \left[ i_{T_i^{\alpha(\mu-1)}}\left( w_i^{\beta(\mu-1)} \right) \right]_\mu =\left[ i_{ \left( T_i^{\alpha(\mu-1)} - T_{0i}^\alpha + T_{0i}^\alpha \right) } \left( w_i^{\beta(\mu-1)} - w_{0i}^\beta + w_{0i}^\beta \right) \right]_\mu = \left[ i_{ \left( T_i^{\alpha(\mu-1)} - T_{0i}^\alpha \right) } \left( w_i^{\beta(\mu-1)} - w_{0i}^\beta \right) \right]_\mu \notag \\
&\Longrightarrow \left| K_{i|\mu}^{\alpha\beta} \right|_{k+\alpha} \ll C_5 \left| T_i^{\alpha(\mu-1)} - T_{0i}^\alpha \right|_{k+\alpha} \left| w_i^{\beta(\mu-1)} -  w_{0i}^\beta \right|_{k+\alpha} \ll C_6 \frac{b}{c} A(t) \label{sd29}
\end{align}

As in $(\ref{ds12})-(\ref{ds25})$ and $(\ref{ns24})-(\ref{ns29})$ we may assume that
\begin{align}\label{sd24}
\left| \lambda_{ij|\mu} \right|_{k+1+\alpha}, \left| \Lambda_{i|\mu} \right|_{k+\alpha}, \left| \Gamma_{i|\mu} \right|_{k+\alpha} ,  \left| \Psi_{i|\mu} \right|_{k-1+\alpha}, \left| D_{ij|\mu} \right|_{k+1+\alpha}, \left| B_{i|\mu} \right|_{k+\alpha}, \left| C_{i|\mu} \right|_{k+\alpha} \ll C_7 \frac{b}{c} A(t)
\end{align}

With this preparation, now we estimate $(\ref{sd21})$ and then apply Lemma \ref{sd18}. Let us estimate $\left( \varphi_\mu'', \left\{  \left( b_{i|\mu}'', \lambda_{i|\mu}''  \right) \right\} \right)$. First we note that from $(\ref{sd22})$ and $(\ref{sd23})$ and $(\ref{sd24})$
\begin{align}\label{dt10}
\left| \left(  \varphi_\mu'', \left\{ \left( b_{i|\mu}'', \lambda_{i|\mu}'' \right) \right\}  \right)   \right|_{k+\alpha} \ll C_8 \frac{b}{c} A(t)
\end{align}
From $(\ref{sd25})$ and $(\ref{sd26})$ and $(\ref{sd27})$ and $(\ref{sd28})$, we may assume that
\begin{align}\label{dt11}
\left| \Phi \right|_{k-1+\alpha}, \left| B \right|_{k-1+\alpha},  \left| V \right|_{k-1+\alpha} , \left| I \right|_{k-1+\alpha} \ll C_9 \frac{b}{c} A(t)
\end{align}
On the other hand, from $(\ref{sd29})$ and $(\ref{sd24})$, we have
\begin{align}\label{dt12}
\left| J \right|_{k-1+\alpha} \ll C_{10} \frac{b}{c} A(t)
\end{align}
By Lemma \ref{sd18}, we can find $\left( \chi_\mu, \left\{ \left(  \eta_i^{\chi_\mu}, \delta_i^{\chi_\mu} \right) \right\} \right) \in A^{0,1}\left( M, \mathcal{E}_{\mathcal{F}_0} \right)$ and $\sigma_\mu\in A^{0,0}\left( M, \mathscr{H}om_{\mathcal{O}_M}\left( \Theta_{\mathcal{F}_0}, \Theta_M \right) \right)$ and $\pi_{i\sigma_\mu}\in \Gamma\left( U_i, \mathscr{H}om_{\mathcal{O}_M}\left( \bigwedge^2 \Theta_{\mathcal{F}_0}, \Theta_{\mathcal{F}_0} \right) \right)$ and $\Sigma_\mu \in A^{0,0}\left(M, \mathscr{H}om_{\mathcal{O}_M}\left( \mathcal{N}_{\mathcal{F}_0}^*, \Omega_M^1 \right) \right)$ satisfying $(\ref{sd32})-(\ref{ds2})$ and $(\ref{sd30}),(\ref{sd31})$ and
\begin{align} \label{dt15}
\bar{\partial} \chi_\mu=0,\,\,\,\,\, \bar{\partial} \eta_i^{\chi_\mu}\left( T_{0i}^\alpha \right)=0,\,\,\,\,\,\bar{\partial} \delta_i^{\chi_\mu}\left( w_{0i}^\alpha \right) =0 \,\,\,\,\,\,\left( \tilde{\square}\left( \chi_\mu, \left\{ \left( \eta_i^{\chi_\mu}, \delta_i^{\chi_\mu}    \right) \right\} \right)=0     \right)
\end{align}
\begin{align}\label{dt16}
i_{\sigma_\mu\left( T_{0i}^\alpha \right)}\left( w_{0i}^\beta \right) + i_{T_{0i}^\alpha}\left( \Sigma_\mu\left( w_{0i}^\beta \right) \right) = K_{i|\mu}^{\alpha\beta}+ i_{\Gamma_{i|\mu}^\alpha}\left( w_{0i}^\beta \right)  + i_{T_{0i}^\alpha}\left( C_{i|\mu}^\beta \right)
\end{align}
\begin{align}
\left| \left( \chi_\mu, \left\{ \left( \eta_i^{\chi_\mu}, \delta_i^{\chi_\mu} \right) \right\} \right) \right|_{k+\alpha}, \left| \sigma_\mu \right|_{k+\alpha}, \left| \pi_{i\sigma_\mu} \right|_{k-1+\alpha} , \left| \Sigma_\mu \right|_{k+\alpha} \ll C_{11} \frac{b}{c} A(t)
\end{align}
Then as in the proof of Theorem \ref{tt2} and Theorem \ref{nn2}, if we set $(\ref{ds40})-(\ref{ds43})$ and $(\ref{ns46}),(\ref{ns47})$, then we can choose $b$ and $c$ satisfying $(\ref{sd1})_\mu-(\ref{sd2})_\mu$. Here the missing point is to check $(\ref{s13})$: in fact,
\begin{align*}
&i_{T_{i|\mu}'^\alpha}\left( w_{0i}^\beta \right) + i_{T_{0i}^\alpha}\left( w_{i|\mu}'^\beta \right)= i_{\sigma_{i|\mu}^\alpha- \left[ \mathfrak{d} G \chi_\mu , T_{0i}^\alpha \right]+ \sum_{\xi=1}^p S_{i|\mu}^{\alpha \xi} T_{0i}^\xi} \left( w_{0i}^\beta \right) + i_{T_{0i}^\alpha}\left( \Sigma_{i|\mu}^\beta - \mathcal{L}_{\mathfrak{d} G \chi_\mu}\left( w_{0i}^\beta\right) + U_{i|\mu}\left( w_{0i}^\beta \right) \right)   \\
&= i_{\sigma_{i|\mu}^\alpha}\left( w_{0i}^\alpha \right) + i_{T_{0i}^\alpha }\left( \Sigma_{i|\mu}^\beta \right) - i_{\left[ \mathfrak{d} G \chi_\mu, T_{0i}^\alpha \right]}\left( w_{0i}^\beta \right) - i_{T_{0i}^\alpha} \mathcal{L}_{\mathfrak{d} G \chi_\mu}\left( w_{0i}^\beta \right) \\
&=K_{i|\mu}^{\alpha\beta} + i_{\Gamma_{i|\mu}^\alpha} \left( w_{0i}^\beta \right) + i_{T_{0i}^\alpha}\left( C_{i|\mu}^\beta \right)  - \mathcal{L}_{\mathfrak{d} G \chi_\mu} i_{T_{0i}^\beta }
 \left( w_{0i}^\beta \right) = K_{i|\mu}^{\alpha\beta} + i_{\Gamma_{i|\mu}^\alpha}\left( w_{0i}^\beta \right) + i_{T_{0i}^\alpha}\left( C_{i|\mu}^\beta \right)
\end{align*}
Then $\varphi(t), T_i^\alpha, r_{ij}^{\alpha\beta}, g_{i\alpha\beta}^\gamma$ and $w_i^\alpha, h_{ij}^{\alpha\beta}$ converge. Then by the same argument in the proof of Theorem \ref{tt2} and Theorem \ref{nn2}, we can construct a foliated analytic family $\pi:\left( \mathcal{M}, \Theta_{\mathcal{F}}, \mathcal{N}_\mathcal{F}^* \right)\to \Delta_\epsilon$ such that the foliated Kodaira-Spencer map $\varphi_0: T_0\left( \Delta_\epsilon \right)\to \mathbb{H}^1\left( M, \mathcal{F}^\bullet \right)$ is bijective. This completes the proof of Theorem \ref{ss6}.

\end{proof}

\section{Theorem of stability for deformations of foliated complex analytic structures in terms of tangent sheaves}

Our proof of theorem of existence of deformations of foliated complex analytic structures in terms of both tangent sheaves and tangent sheaves (Theorem \ref{ss6}) provides the proof of theorem of stability for deformations of foliated complex analytic structures in terms of both tangent and cotangent sheaves (Theorem \ref{snt1}). We define a complex of sheaves which we truncate the $0$-th term of $\mathcal{F}_0^\bullet$.
\begin{center}
$\mathcal{F}_0^{\bullet \geq 1}:\begin{CD}
\cdots \\
@AF_3AA \\
\mathscr{H}om_{\mathcal{O}_M}\left(  \bigwedge^3 \Theta_{\mathcal{F}_0}, \frac{\Theta_M}{\Theta_{\mathcal{F}_0}} \right) \bigoplus \mathscr{H}om_{\mathcal{O}_M}\left( \mathcal{N}_{\mathcal{F}_0}^* , \tilde{\mathcal{S}}^3 \right) \bigoplus \mathscr{H}om_{\mathcal{O}_M}\left( \mathcal{N}_{\mathcal{F}_0}^*, \bigwedge^2  \Theta_{\mathcal{F}_0}^*  \right) \\
@AF_2AA \\
\mathscr{H}om_{\mathcal{O}_M}\left(  \bigwedge^2 \Theta_{\mathcal{F}_0}, \frac{\Theta_M}{\Theta_{\mathcal{F}_0}} \right) \bigoplus \mathscr{H}om_{\mathcal{O}_M}\left( \mathcal{N}_{\mathcal{F}_0}^* , \tilde{\mathcal{S}}^2 \right) \bigoplus \mathscr{H}om_{\mathcal{O}_M}\left( \mathcal{N}_{\mathcal{F}_0}^*, \Theta_{\mathcal{F}_0}^*  \right) \\
@AF_1AA \\
\mathscr{H}om_{\mathcal{O}_M} \left(  \Theta_{\mathcal{F}_0}, \frac{\Theta_M}{\Theta_{\mathcal{F}_0}}\right) \bigoplus \mathscr{H}om_{\mathcal{O}_M} \left( \mathcal{N}_{\mathcal{F}_0}^*, \frac{\Omega_M^1}{\mathcal{N}_{\mathcal{F}_0}^*} \right)
\end{CD}$
\end{center}
We will denote the $i$-th cohomology group of $\mathcal{F}_0^{\bullet \geq 1}$ by $\mathbb{H}^i\left( M, \mathcal{F}_0^{\bullet \geq 1} \right)$.

\begin{theorem}[Theorem of stability for deformations of foliated complex analytic structures in terms of tangent sheaves] \label{snt1}
Let $\left( M, \Theta_{\mathcal{F}_0} \right)$ be a compact foliated complex manifold with $\Theta_{\mathcal{F}_0}$ locally free. Assume that $\mathbb{H}^1\left( M, \mathcal{F}_0^{\bullet \geq 1} \right)=0$. Then for any complex analytic family $\pi:\mathcal{M}\to B$ of deformations of $\pi^{-1}(0)=M, 0\in B$ in the sense of Kodaira-Spencer, there exists an open neighborhood $N\subset B$ of $0$ and a locally free subsheaf $\Theta_{\mathcal{F}}$ of $\Theta_{\frac{\mathcal{M}|_N}{N}}$ and a locally free subsheaf $\mathcal{N}_{\mathcal{F}}^*$ of $\Omega_{\frac{\mathcal{M}|_N}{N}}^1$ such that $\left(\mathcal{M}|_N, \Theta_{\mathcal{F}} , \mathcal{N}_{\mathcal{F}}^* \right)$ defines a foliated analytic family $\pi|_N: \left(\mathcal{M}|_N , \Theta_{\mathcal{F}} , \mathcal{N}_\mathcal{F}^* \right) \to N$ of deformations of $\left(M, \Theta_{\mathcal{F}_0}  , \mathcal{N}_{\mathcal{F}_0}^* \right)=\pi^{-1}(0)$ in terms of both tangent and cotangent sheaves.
\end{theorem}

\begin{proof}
We copy the proof of Theorem \ref{ss6}. The difference is that we assume that $\varphi(t)$ (determined by $\mathcal{M}$ over some neighborhood $N'\subset B$ of $0$) is given from the beginning and replace $\varphi^\mu$ by $\varphi(t)$ in the proof of Theorem \ref{tt2}. Accordingly we replace $\xi_\mu$ and $\varphi_\mu$ by $0$. Then the obstructions for solving $(\ref{te12})-(\ref{t33})$ are in $\mathbb{H}^1\left( M, \mathcal{F}_0^{\bullet \geq 1} \right)=0$. Hence we can construct $\Theta_{\mathcal{F}}$ and $\mathcal{N}_{\mathcal{F}}^*$ on $\mathcal{M}|_N$ for some neighborhood $N\subset N'$ of $0$. 
\end{proof}

\section{Theorem of completeness of deformations of foliated complex analytic structures in terms of both tangent and cotangent sheaves}

\subsection{Change of parameters}\

Consider a foliated complex analytic family $\left( \mathcal{M}, \Theta_\mathcal{F}, \mathcal{N}_\mathcal{F}, B, \pi \right)$ with $\Theta_{\mathcal{F}}$ and $\mathcal{N}_\mathcal{F}^*$ locally free of simultaneous deformation of $\left( M_t, \Theta_{\mathcal{F}_t}, \mathcal{N}_{\mathcal{F}_t}^* \right)=\pi^{-1}(t),t\in B$ where $B$ is a domain of $\mathbb{C}^m$. Let $D$ be a domain of $\mathbb{C}^{m'}$ and $s:u\to t=s(u), u\in D$, a holomorphic map of $D$ into $B$. Then by changing the parameter from $t$ to $u$, we will construct a foliated analytic family $\left\{ \left( M_{s(u)} , \Theta_{\mathcal{F}_{u(s)}}, \mathcal{N}_{\mathcal{F}_{u(s)}}^*  \right) \right\}$ on the parameter $D$. By subsection \ref{r5} (we keep the notations in subsection \ref{r5}), we have the foliated analytic family $\left( \mathcal{M}\times_B D, p^* \Theta_\mathcal{F}, D, \pi' \right)$ induced from $\left(\mathcal{M}, \Theta_\mathcal{F}, B, \pi \right)$. On the other hand, by subsection \ref{r6} (we keep the notations in subsection \ref{r6}), we have the foliated analytic family $\left( \mathcal{M}\times_B D, p^* \mathcal{N}_\mathcal{F}^* , D, \pi' \right)$ induced from $\left( \mathcal{M}, \mathcal{N}_\mathcal{F}^*, B, \pi \right)$. Moreover, since $\mathcal{N}_\mathcal{F}^*$ vanishes on $\Theta_\mathcal{F}$, $p^* \mathcal{N}_\mathcal{F}^*$ vanishes on $p^* \Theta_\mathcal{F}^*$. Then we have $p^* \Theta_\mathcal{F} \subset \left( \frac{\Omega_{\mathcal{M}\times_B D/D}^1}{p^* \mathcal{N}_\mathcal{F}^*} \right)^*$ and $p^* \mathcal{N}_\mathcal{F}^*\subset \left( \frac{\Theta_{\mathcal{M}\times_B D/D}}{p^* \Theta_\mathcal{F}^*} \right)^*$. Then we have an exact sequence $0\to \left(\frac{\Omega_{\mathcal{M}\times_B D/D}^1}{p^* \mathcal{N}_\mathcal{F}^*} \right)^* \to \Theta_{\mathcal{M}\times_B D} \to \left(p^* \mathcal{N}_\mathcal{F}^* \right)^*$, and for each $u\in D$, by restricting the natural map $\frac{\Theta_{\mathcal{M}\times_B D}}{p^* \Theta_\mathcal{F}}\to \left(p^* \mathcal{N}_\mathcal{F}^* \right)^*$ to $\pi'^{-1}(u)$, we have an injection $\frac{\Theta_{M_{s(u)}}}{\mathcal{N}_{\mathcal{F}_{s(u)}}^*}\to \mathcal{N}_{\mathcal{F}_{s(u)}}^*$. We note that $\frac{\Theta_{\mathcal{M}\times_B D/D}}{p^* \Theta_\mathcal{F}}$ and $\left( p^* \mathcal{N}_\mathcal{F}^* \right)^*$ are flat over $D$. This implies that $\frac{\Theta_{\mathcal{M}\times_B D/D}}{p^* \Theta_\mathcal{F}}\to \left(p^* \mathcal{N}_\mathcal{F}^* \right)^*$ is injective, so that $0 \to p^* \Theta_\mathcal{F} \to \Theta_{\mathcal{M}\times_B D}\to \left(p^* \mathcal{N}_\mathcal{F}^* \right)^*$ is exact. Hence $p^* \Theta_\mathcal{F}=  \left(\frac{\Omega_{\mathcal{M}\times_B D/D}^1}{p^* \mathcal{N}_\mathcal{F}^*} \right)^*$. On the other hand, we have an exact sequence $0\to \left( \frac{\Theta_{\mathcal{M}\times_B D/D}}{p^* \Theta_\mathcal{F}} \right)^* \to \Omega_{\mathcal{M}\times_B D/D}^1 \to \left(p^* \Theta_\mathcal{F} \right)^*$, and for each $u\in D$, by restricting the natural map $\frac{\Omega_{\mathcal{M}\times_B D/D}^1}{p^* \mathcal{N}_\mathcal{F}^*}\to \left(p^* \Theta_\mathcal{F} \right)^*$ to $\pi'^{-1}(u)$, we have an injection $\frac{\Omega_{M_{s(u)}}^1}{\mathcal{N}_{\mathcal{F}_{s(u)}}^*}\to \Theta_{\mathcal{F}_{s(u)}}^*$. We note that $\frac{\Omega_{\mathcal{M}\times_B D/D}^1}{p^* \mathcal{N}_\mathcal{F}^* }$ and $\left( p^* \Theta_\mathcal{F}\right)^*$ are flat over $D$. This implies that $\frac{\Omega_{\mathcal{M}\times_B D/D}^1}{p^* \mathcal{N}_\mathcal{F}^*}\to \left(p^* \Theta_\mathcal{F} \right)^*$ is injective, so that $0\to p^* \mathcal{N}_\mathcal{F}^* \to \Omega_{\mathcal{M}\times_B D/D}^1 \to \left( p^* \Theta_\mathcal{F}^* \right)^*$. Hence $p^*\mathcal{N}_\mathcal{F}^* =\left( \frac{\Theta_{\mathcal{M}\times_B D/D}}{p^* \Theta_\mathcal{F}} \right)^*$. This implies that $\left\{ \left( M_{s(u)}, \Theta_{\mathcal{F}_{u(s)}}, \mathcal{N}_{\mathcal{F}_{u(s)}}^* \right)| u\in D \right\}$ forms a foliated complex analytic family $\left( \mathcal{M}\times_B D , p^* \Theta_\mathcal{F}, p^* \mathcal{N}_\mathcal{F}^*, D, \pi' \right)$.

\begin{definition}
The foliated complex analytic family $\left( \mathcal{M} \times_B D, p^* \Theta_\mathcal{F}, p^* \mathcal{N}_\mathcal{F}^*,  \pi' \right)$ is called the foliated complex analytic family induced from $\left( \mathcal{M}, \Theta_\mathcal{F}, \mathcal{N}_\mathcal{F}^*, B, \pi \right)$ by the holomorphic map $s:D\to B$.
\end{definition}

We note that by combining Theorem \ref{r8} and Theorem \ref{r11}, we have
\begin{theorem}
For any tangent vector $\frac{\partial }{\partial u}=c_1\frac{\partial}{\partial u_1} + \cdots + c_r\frac{\partial}{\partial u_{m'}}\in T_u(D)$, the infinitesimal foliated deformation of $(M_{s(u)}, \Theta_{\mathcal{F}_{s(u)}}, \mathcal{N}_{\mathcal{F}_{s(u)}}^*)$ along $\frac{\partial}{\partial u}$ is given by
\begin{align*}
\frac{\partial \left(M_{s(u)}, \Theta_{\mathcal{F}_{s(u)}}, \mathcal{N}_{\mathcal{F}_{s(u)}}^*\right)}{\partial u}&=\left(\sum_{\gamma=1}^m \frac{\partial t_\gamma}{\partial u}\frac{\partial M_t}{\partial t_\gamma},\left\{T_j^\alpha\left(z_j,s(u)\right) \mapsto - \overline{\sum_{\gamma=1}^m \frac{\partial t_\gamma}{\partial u} \frac{\partial T_j^\alpha}{\partial t_\gamma}} \right\}, \left\{w_j^\alpha\left(z_j, s(u)\right) \mapsto -\overline{\sum_{\gamma=1}^m \frac{\partial t_\gamma}{\partial u}\frac{\partial w_j^\alpha}{\partial t_\gamma}} \right\} \right)\\
            &=\sum_{\gamma=1}^m \frac{\partial t_\gamma}{\partial u} \frac{\partial \left(M_t, \Theta_{\mathcal{F}_t}, \mathcal{N}_{\mathcal{F}_t}^* \right) }{\partial t_\gamma}
\end{align*}
\end{theorem}

\subsection{Theorem of completeness of deformations of foliated complex analytic structures in terms of both tangent and cotangent sheaves}\

\begin{definition}\label{r10}
Let $\left( \mathcal{M}, \Theta_\mathcal{F}, \mathcal{N}_\mathcal{F}^* , B, \pi \right)$ with $\Theta_\mathcal{F}$ and $\mathcal{N}_\mathcal{F}^*$ locally free be a foliated complex analytic family of compact foliated complex manifolds in terms of both tangent and cotangent sheaves, and $t^0\in B$. Then $\left( \mathcal{M}, \Theta_\mathcal{F}, \mathcal{N}_\mathcal{F}^*, B, \pi \right)$ is called complete at $t^0\in B$ if for any foliated complex analytic family $\left( \mathcal{M}', \Theta_{\mathcal{F}'}, \mathcal{N}_{\mathcal{F}'}^*, D, \pi'  \right)$ in terms of both tangent and cotangent sheaves such that $D$ is a domain of $\mathbb{C}^{m'}$ containing $0$ and that $\pi'^{-1}(0)=\pi^{-1}(t^0)=\left( M, \Theta_{\mathcal{F}_0}, \mathcal{N}_{\mathcal{F}_0}^* \right)$, there are a sufficiently small domain $\Delta$ with $0\in \Delta \subset D$, and a holomorphic map $s:u\to t=s(u)$ with $s(0)=t^0$ such that $\left( \mathcal{M}'_\Delta , \Theta_{\mathcal{F}'_\Delta}, \mathcal{N}_{\mathcal{F}'_\Delta}^*,\Delta, \pi' \right)$ is the foliated complex analytic family induced from $\left( \mathcal{M}, \Theta_\mathcal{F}, \mathcal{N}_\mathcal{F}^*, B, \pi \right)$ by $s$ where $\left( \mathcal{M}'_\Delta , \Theta_{\mathcal{F}'_\Delta}, \mathcal{N}_{\mathcal{F}'_\Delta}^*,\Delta, \pi' \right)$ is the restriction of $\left( \mathcal{M}', \Theta_{\mathcal{F}'}, \mathcal{N}_{\mathcal{F}'}^*, D, \pi' \right)$ to $\Delta$.
\end{definition}

We shall prove the following theorem
\begin{theorem}\label{sc1}
Let $\left(\mathcal{M}, \Theta_\mathcal{F}, \mathcal{N}_\mathcal{F}^* ,  B, \pi \right)$ be a foliated analytic family of deformations of a complex foliated complex manifold $(M , \Theta_{\mathcal{F}_0}, \mathcal{N}_{\mathcal{F}_0}^* )=\omega^{-1}(0)$ such that both $\Theta_{\mathcal{F}_0}$ and $\mathcal{N}_{\mathcal{F}_0}^*$ are locally free, and $B$ is a domain of $\mathbb{C}^r$ containing $0$. If the foliated Kodaira-Spence map $\varphi_0: T_0(B)\to \mathbb{H}^1(M, \mathcal{F}_0^\bullet)$ is surjective, the foliated analytic family $\left(\mathcal{M}, \Theta_\mathcal{F}, \mathcal{N}_\mathcal{F}^* ,  B, \pi \right)$ in terms of both tangent and cotangent sheaves is complete at $0 \in B$.
\end{theorem}

\begin{proof}
Let $\left( \mathcal{M}, \Theta_\mathcal{F}, \mathcal{N}_\mathcal{F}^*, B, \pi \right)$ be a foliated analytic family in terms of both tangent and cotangent sheaves which is represented as in Remark \ref{ss17}. We keep the notations in the proof of  Theorem \ref{tc1} and Theorem \ref{nc1}. Then we only consider the condition on $\left(M, \Theta_{\mathcal{F}_0}, \mathcal{N}_{\mathcal{F}_0}^* \right)$: for $\alpha=1,....,p$ and $\beta=1,...,q$,
\begin{align*}
i_{T_{0i}^\alpha}\left( w_{0i}^\beta \right) = 0
\end{align*}

Let $\left(\mathcal{M}', \Theta_{\mathcal{F}'}, \mathcal{N}_{\mathcal{F}'}^*  , D, \pi' \right)$ be an another foliated analytic family such that $\pi'^{-1}(0')=\left(M, \Theta_{\mathcal{F}_0}, \mathcal{N}_{\mathcal{F}_0}^*\right)$. We may assume the following: $(\ref{sc2})-(\ref{sc3})$ in the proof of Theorem \ref{tc1}, and $(\ref{sc4})-(\ref{sc5})$ in the proof of Theorem \ref{nc1} and additionally
\begin{align}\label{cs4}
i_{T_i'^\alpha}\left( w_i'^\beta \right)=0 \,\,\,\,\,\,\,\,\textnormal{for}\,\,\,\alpha=1,...,p,\,\,\,\beta=1,...,q
\end{align}

In order to prove Theorem \ref{sc1}, it suffices to construct holomorphic functions
\begin{align*}
\varphi_i&:\mathcal{U}_i'\to \mathbb{C}^n \\
s &:D \to \mathbb{C}^{m}
\end{align*}
and a matrix function $\left(b_i^{\alpha\beta}(\xi_i, u)\right)_{\alpha,\beta=1,...,p}$ with $b_i^{\alpha\beta}:\mathcal{U}_i'\to \mathbb{C}$ and a matrix function $\left( c_i^{\alpha\beta}(\xi_i, u) \right)_{\alpha,\beta=1,...,q}$ with $c_i^{\alpha\beta}:\mathcal{U}_i' \to \mathbb{C}$ such that $(\ref{tc211})-(\ref{tc4})$ and $(\ref{nc211})-(\ref{nc4})$.

First we prove the existence of formal solution of $(\ref{tc2})-(\ref{tc4})$ and $(\ref{nc2})-(\ref{nc4})$. We recall Notation \ref{te27}. Then $(\ref{tc2})-(\ref{tc4})$ are equivalent to the systems of congruences $(\ref{tc5})_\mu-(\ref{tc7})_\mu$, respectively, and $(\ref{nc2})-(\ref{nc4})$ are equivalent to the systems of congruences $(\ref{nc5})_\mu-(\ref{nc7})_\mu$, respectively for $\mu=1,2,3, \cdots$. We shall construct $\varphi_i^\mu, s^\mu, b_i^{\alpha\beta \mu}$ and $c_i^{\alpha\beta \mu}$ satisfying $(\ref{tc5})_\mu-(\ref{tc7})_\mu$, and $(\ref{nc5})_\mu- (\ref{nc7})_\mu$ by induction on $\mu$. We assume that $\varphi_i^{\mu-1}, s^{\mu-1}, b_i^{\alpha\beta(\mu-1)}$ and $c_i^{\alpha\beta (\mu-1)}$ are already determined. Then we define homogenous polynomials $\Gamma_{ij|\mu}^\alpha, B_{ij|\mu}^{\alpha\beta}, G_{i|\mu}^{\alpha\beta \gamma}$, and $\Pi_{i|\mu}^{\alpha \gamma}$ of degree $\mu$ by the congruences $(\ref{tc35})-(\ref{tc38})$, respectively, and we define homogenous polynomials $A_{ij|\mu}^{\alpha\beta}$ and $W_{i|\mu}^{\alpha \gamma}$ of degree $\mu$ by the congruences $(\ref{nc9}), (\ref{nc10})$, respectively. We set
\begin{align*}
\Gamma_{ij|\mu}:= \sum_{\alpha=1}^n \Gamma_{ij|\mu}^\alpha \frac{\partial}{\partial z_i^\alpha},\,\,\,\,\,\,\,\,\,\, \Pi_{i|\mu}^\alpha:=\sum_{\gamma=1}^n \Pi_{i|\mu}^{\alpha \gamma} \frac{\partial}{\partial z_i^\gamma},\,\,\,\,\,\,\,\,\,\,\, W_{i|\mu}^\alpha := \sum_{\gamma=1}^n W_{i|\mu}^{\alpha \gamma} dz_i^\gamma
\end{align*}

\begin{lemma}\label{sc8}
We have the following equalities$:$ $(\ref{tc32})-(\ref{tc34})$ and $(\ref{nc12})-(\ref{nc14})$ and additionally
\begin{align}\label{sc6}
i_{\Pi_{i|\mu}^\alpha}\left( w_{0i}^\beta \right) - i_{T_{0i}^\alpha}\left( W_{i|\mu}^\beta \right) = 0 \,\,\,\,\,\,\,\,\textnormal{for}\,\,\,\alpha=1,...,p,\,\,\,\beta=1,...,q
\end{align}
\end{lemma}

\begin{proof}
$(\ref{tc32})-(\ref{tc34})$ follows from Lemma \ref{tc39}, and $(\ref{n12})-(\ref{nc14})$ follows from Lemma \ref{nc11}. We prove $(\ref{sc6})$. In fact, we note that
\begin{align*}
i_{T_i^\alpha}\left(w_i^\beta\right)&= \sum_{\gamma=1}^n  T_i^{\alpha\gamma}(z_i,t) w_i^{\beta\gamma}(z_i,t)=0 \\
i_{T_i'^\alpha}\left(w_i'^\beta\right)&= \sum_{\gamma=1}^n  T_i'^{\alpha\gamma}(\xi_i,u) w_i'^{\beta\gamma}(\xi_i, u )=0
\end{align*}
Then we have from $(\ref{tc38})$ and $(\ref{nc10})$
\begin{align*}
0&=\sum_{\gamma=1}^n T_i^{\alpha\gamma}\left(\varphi_i^{\mu-1}, s^{\mu-1}\right) w_i^{\beta\gamma}\left( \varphi_i^{\mu-1}, s^{\mu-1}\right) = \sum_{\gamma=1}^n \left( \Pi_{i|\mu}^{\alpha\gamma}+ \sum_{\eta=1}^p \sum_{\sigma=1}^n b_i^{\alpha \eta(\mu-1)} T_i'^{\eta \sigma}\frac{\partial \varphi_i^{\gamma(\mu-1)}}{\partial \xi_i^\sigma}   \right)w_i^{\beta\gamma}\left( \varphi_i^{\mu-1}, s^{\mu-1}\right)\\
&= \sum_{\gamma=1}^n \Pi_{i|\mu}^{\alpha\gamma} w_{0i}^{\beta\gamma} + \sum_{\sigma=1}^n\sum_{\eta=1}^p b_i^{\alpha\eta(\mu-1)} T_i'^{\eta\sigma} \left( \sum_{\gamma=1}^n w_i^{\beta\gamma}\left(\varphi_i^{\mu-1}, s^{\mu-1}\right) \frac{\partial \varphi_i^{\gamma(\mu-1)}}{\partial \xi_i^\sigma}   \right) \\
&= \sum_{\gamma=1}^n \Pi_{i|\mu}^{\alpha\gamma} w_{0i}^{\beta\gamma} + \sum_{\sigma=1}^n\sum_{\eta=1}^p b_i^{\alpha\eta(\mu-1)} T_i'^{\eta\sigma} \left( \sum_{\delta=1}^q c_i^{\beta\delta (\mu-1)}w_i'^{\delta \sigma} - W_{i|\mu}^{\beta\sigma} \right)\\
&=\sum_{\gamma=1}^n \Pi_{i|\mu}^{\alpha\gamma} w_{0i}^{\beta\gamma}- \sum_{\sigma =1}^n \sum_{\eta=1}^p b_{0i}^{\alpha\eta} T_{0i}'^{\eta \sigma} W_{i|\mu}^{\beta\sigma} + \sum_{\delta=1}^q \sum_{\eta=1}^p b_i^{\alpha\eta(\mu-1)}  c_i^{\beta\delta(\mu-1)} \left( \sum_{\delta=1}^n T_i'^{\eta\sigma} w_i'^{\delta\sigma}     \right)\\
&= \sum_{\gamma=1}^n \Pi_{i|\mu}^{\alpha\gamma} w_{0i}^{\beta\gamma} - \sum_{\sigma=1}^n T_{0i}^{\alpha \sigma} W_{i|\mu}^{\beta\sigma} = i_{\Pi_{i|\mu}^\alpha}\left(w_{0i}^\beta \right)- i_{T_{0i}^\alpha}\left(W_{i|\mu}^\beta \right)
\end{align*}
This completes the proof of Lemma \ref{sc8}.
\end{proof}

Our purpose is to determine $\varphi^\mu= \varphi^{\mu-1}+ \varphi_{i|\mu}, s^\mu= s^{\mu-1}+ s_\mu, b_i^{\alpha\beta(\mu-1)}+ b_{i|\mu}^{\alpha\beta}$, and $c_i^{\alpha\beta(\mu-1)}+ c_{i|\mu}^{\alpha\beta}$ satisfying $(\ref{tc5})_\mu-(\ref{tc7})_\mu$, and $(\ref{nc5})_\mu-(\ref{nc7})_\mu$.

\begin{lemma}\label{sc7}
$(\ref{tc5})_\mu-(\ref{tc7})_\mu$ are equivalent to $(\ref{tc8})-(\ref{tc10})$, respectively, and $(\ref{nc5})_\mu-(\ref{nc7})_\mu$ are equivalent to $(\ref{nc17})-(\ref{nc19})$, respectively.
\end{lemma}

\begin{proof}
Lemma \ref{sc7} follows from Lemma \ref{tc41} and Lemma \ref{nc23}.
\end{proof}

We define an element $\overline{\Pi}_{i|\mu}\in \Gamma\left( U_i, \mathscr{H}om_{\mathcal{O}_M} \left( \Theta_{\mathcal{F}_0}, \frac{\Theta_M}{\Theta_{\mathcal{F}_0}}   \right) \right)$ as in $(\ref{tc213})$ and $\overline{W}_{i|\mu}\in \Gamma\left( U_i , \mathscr{H}om_{\mathcal{O}_M} \left( \mathcal{N}_{\mathcal{F}_0}^*, \frac{\Omega_M^1}{\mathcal{N}_{\mathcal{F}_0}^*}    \right) \right)$ as in $(\ref{nc125})$. Then $(\ref{tc32})-(\ref{tc34})$ and $(\ref{nc12})-(\ref{nc14})$ and $(\ref{sc6})$ imply that
\begin{align}
\left( \left\{ \overline{\Pi}_{i|\mu} \right\}, \left\{ - \overline{W}_{i|\mu} \right\} , \left\{ \Gamma_{ij|\mu} \right\}  \right) \in C^0\left( \mathcal{U}, \mathscr{H}om_{\mathcal{O}_M}\left( \Theta_{\mathcal{F}_0} , \frac{\Theta_M}{\Theta_{\mathcal{F}_0}}   \right) \right) \bigoplus C^0\left( \mathcal{U}, \mathscr{H}om_{\mathcal{O}_M} \left(  \mathcal{N}_{\mathcal{F}_0}^*, \frac{\Omega_M^1}{\mathcal{N}_{\mathcal{F}_0}^*}  \right) \right) \bigoplus C^1\left( \mathcal{U}, \Theta_M \right) \label{ssc2}
\end{align}
defines a $1$-cocycle in the following \v Cech resolution of $\mathcal{F}_0^\bullet$:
{\Tiny{\begin{center}
$\begin{CD}
\cdots \\
@AAA \\
C^0\left( \mathcal{U}, \bigwedge^2 \Theta_{\mathcal{F}_0}^*\otimes \frac{\Theta_M}{\Theta_{\mathcal{F}_0}} \right) \bigoplus C^0\left( \mathcal{U}, \left(\mathcal{N}_{\mathcal{F}_0}^* \right)^* \otimes \tilde{\mathcal{S}}^2 \right) \bigoplus C^0\left( \mathcal{U}, \left( \mathcal{N}_{\mathcal{F}_0}^* \right)^*\otimes \Theta_{\mathcal{F}_0}^* \right) @>-\delta >> \cdots \\
@AAA @AAA  \\
C^0\left(\mathcal{U}, \Theta_{\mathcal{F}_0}^* \otimes \frac{\Theta_M}{\Theta_{\mathcal{F}_0} } \right) \bigoplus C^0\left( \mathcal{U}, \left(\mathcal{N}_{\mathcal{F}_0}^* \right)^* \otimes \frac{\Omega_M^1}{\mathcal{N}_{\mathcal{F}_0}^*} \right) @>\delta>>  C^1\left(\mathcal{U}, \Theta_{\mathcal{F}_0}^* \otimes \frac{\Theta_M}{\Theta_{\mathcal{F}_0} }  \right) \bigoplus C^1\left( \mathcal{U}, \left( \mathcal{N}_{\mathcal{F}_0}^* \right)^* \otimes  \frac{\Omega_M^1}{\mathcal{N}_{\mathcal{F}_0}^*} \right) @>-\delta >> \cdots \\
@AAA @AAA @AAA \\
C^0(\mathcal{U}, \Theta_M) @>-\delta>> C^1(\mathcal{U}, \Theta_M) @>\delta>> C^2(\mathcal{U}, \Theta_M) @>-\delta>> \cdots
\end{CD}$
\end{center}}}
By the hypothesis that the foliated Kodaira-Spencer map $\varphi_0: T_0(B) \to \mathbb{H}^1\left( M, \mathcal{F}_0^\bullet \right)$ is surjective, we can find homogeneous polynomial $s_\mu^\lambda$ such that
\begin{align*}
\varphi_0\left( \sum_{\lambda=1}^r s_\mu^\lambda \frac{\partial}{\partial t_\lambda}   \right) = \left(  \left\{  \overline{\Pi}_{i|\mu}  \right\} , \left\{ - \overline{W}_{i|\mu} \right\}   \left\{ \Gamma_{ij|\mu}  \right\}   \right)
\end{align*}

Since we have
\begin{align*}
\varphi_0\left( \frac{\partial}{\partial t_\lambda} \right) =\left( \{\alpha_{i\lambda}\} =\left\{ T_{0i}^\alpha \mapsto \overline{ -\frac{\partial T_i^\alpha (z_i,t)}{\partial t_\lambda}|_{t=0} }\right\}, \left\{ \beta_{i\lambda}\right\}= \left\{ w_{0i}^\alpha \mapsto \overline{- \frac{\partial w_i^\alpha(z_i,t)}{\partial t_\lambda} |_{t=0} }\right\} ,  \{\rho_{ij\lambda} \} =  \left\{ \sum_{\alpha=1}^n \frac{\partial f_{ij}^\alpha}{\partial t_\lambda}|_{t=0} \frac{\partial}{\partial z_i^\alpha } \right\}   \right),
\end{align*}
there exists $\{\varphi_{i|\mu}\}\in C^0\left(\mathcal{U}, \Theta_M \right)$ and $\left\{ b_{i|\mu}\right\} \in C^0\left( \mathcal{U}, \mathscr{H}om_{\mathcal{O}_M}\left( \Theta_{\mathcal{F}_0}, \Theta_{\mathcal{F}_0} \right) \right)$ and $\left\{ c_{i|\mu} \right\} \in C^0\left( \mathcal{U}, \mathscr{H}om_{\mathcal{O}_M} \left( \mathcal{N}_{\mathcal{F}_0}^*, \mathcal{N}_{\mathcal{F}_0}^* \right) \right)$ satisfying $(\ref{tc8})-(\ref{tc10})$ and $(\ref{nc17})-(\ref{nc19})$ as in the proof of Theorem \ref{tc1} and Theorem \ref{nc1}. This completes the inductive construction of $\varphi_i^\mu,s^\mu, b_i^{\alpha\beta \mu}$ and $c_i^{\alpha \beta \mu}$.

\begin{remark}\label{ssc7}
We recall from $(\ref{ssc1})$ that we have $\mathbb{H}^1\left( M, \mathcal{F}_0^\bullet \right)\cong \mathbb{H}^1\left( M, \mathcal{E}_{\mathcal{F}_0}^\bullet \right)$. We will reinterpret $(\ref{ssc2})$ in $\mathbb{H}^1\left( M, \mathcal{F}_0^\bullet \right)$ in terms of $\mathbb{H}^1\left(  M, \mathcal{E}_{\mathcal{F}_0}^\bullet \right)$ in the following \v Cech resolution of $\mathcal{E}_{\mathcal{F}_0}^\bullet:$
{\Tiny{\begin{equation} \label{ssc12}
\begin{CD}
\cdots \\
@AF_2AA \\
C^0\left(\mathcal{U}, \bigwedge^2 \Theta_{\mathcal{F}_0}\otimes \Theta_M \right) \bigoplus C^0\left( \right) \bigoplus C^0\left( \mathcal{U}, \left( \mathcal{N}_{\mathcal{F}_0}^* \right)^*\otimes \tilde{\mathcal{S}}^2 \right) \bigoplus C^0\left( \mathcal{U}, \left( \mathcal{N}_{\mathcal{F}_0}^* \right)^* \otimes \Theta_{\mathcal{F}_0}^* \right) @>-\delta>> \cdots \\
@AF_1'AA @AF_1'AA \\
C^0\left( \mathcal{U}, \Theta_{\mathcal{F}_0}^*\otimes \Theta_M \right)\bigoplus C^0\left(\mathcal{U}, \left( \mathcal{N}_{\mathcal{F}_0}^* \right)^*\otimes \Omega_M^1 \right) @>\delta>> C^1\left( \mathcal{U}, \Theta_{\mathcal{F}_0}^*\otimes \Theta_M \right)\bigoplus C^1\left(\mathcal{U}, \left( \mathcal{N}_{\mathcal{F}_0}^* \right)^*\otimes \Omega_M^1 \right)  @>-\delta>> \cdots  \\
@AF_0'AA @AF_0'AA @AF_0'AA \\
C^0\left( \mathcal{U}, \mathcal{E}_{\mathcal{F}_0 }  \right) @> - \delta >>  C^1\left( \mathcal{U}, \mathcal{E}_{\mathcal{F}_0} \right) @>\delta>> C^2\left( \mathcal{U}, \mathcal{E}_{\mathcal{F}_0} \right)
\end{CD}
\end{equation}}}
We set $B_{ij|\mu}\in \Gamma\left( U_{ij}, \mathscr{H}om_{\mathcal{O}_M}\left( \Theta_{\mathcal{F}_0}, \Theta_{\mathcal{F}_0} \right) \right)$ as in $(\ref{ssc3})$ and $A_{ij|\mu}\in \Gamma\left( U_{ij}, \mathscr{H}om_{\mathcal{O}_M} \left( \mathcal{N}_{\mathcal{F}_0}^*, \mathcal{N}_{\mathcal{F}_0}^* \right) \right)$ as in $(\ref{ssc4})$. Then $\left\{ \left( \Gamma_{ij|\mu}, B_{ij|\mu}, A_{ij|\mu}     \right) \right\}\in C^1\left( \mathcal{U}, \mathcal{E}_{\mathcal{F}_0} \right)$ and $\delta\left( \left\{ \left( \Gamma_{ij|\mu}, B_{ij|\mu}, A_{ij|\mu} \right)  \right\} \right)=0$. Hence
\begin{align}
\left(\left\{ \Pi_{i|\mu} \right\}, \left\{ - W_{i|\mu} \right\}, \left\{ \left( \Gamma_{ij|\mu}, B_{ij|\mu}, A_{ij|\mu}  \right) \right\} \right) \in C^0\left( \mathcal{U}, \mathscr{H}om_{\mathcal{O}_M}\left( \Theta_{\mathcal{F}_0}, \Theta_M \right) \right) \bigoplus C^0\left( \mathcal{U}, \mathscr{H}om_{\mathcal{O}_M}\left( \mathcal{N}_{\mathcal{F}_0}^*, \Omega_M^1 \right) \right) \bigoplus C^1 \left( \mathcal{U}, \mathcal{E}_{\mathcal{F}_0} \right)
\end{align}
defines a $1$-cocycle in the above \v Cech resolution of $\mathcal{E}_{\mathcal{F}_0}^\bullet$.
\end{remark}

\subsection{Proof of convergence} \label{scc5} \

We prove that we can choose appropriate $\varphi_{i|\mu}, s_\mu, b_{i|\mu}^{\alpha\beta}$ and $c_{i|\mu}^{\alpha \beta}$ satisfying $(\ref{tc8})-(\ref{tc10})$ and $(\ref{nc17})-(\ref{nc19})$ in each inductive step so that
\begin{align*}
s(u)&=s_1(u)+ s_2(u)+\cdots + s_\mu(u)+ \cdots\\
\varphi_i\left( \xi_i, u\right)&= \xi_i+ \varphi_{i|1}\left(\xi_i, u\right) +  \varphi_{i|2}\left( \xi_i, u \right) + \cdots + \varphi_{i|\mu}(\xi_i, u) + \cdots \\
b_i^{\alpha\beta}(\xi_i, u)&= b_{0i}^{\alpha\beta}(\xi_i) + b_{i|1}^{\alpha\beta}(\xi_i,u) + \cdots + b_{i|\mu}^{\alpha\beta} (\xi_i, u) + \cdots \\
c_i^{\alpha\beta}(\xi_i, u)&= c_{0i}^{\alpha\beta}(\xi_i) + c_{i|1}^{\alpha\beta}(\xi_i,u) + \cdots + c_{i|\mu}^{\alpha\beta} (\xi_i, u) + \cdots
\end{align*}
converge absolutely and uniformly for $|u|< \epsilon$ if $\epsilon>0$ is sufficiently small. It suffices to prove the estimates $s(u)\ll A(u), \varphi_i(\xi_i,u)-\xi_i\ll A(u)$ and $b_i^{\alpha\beta}(\xi_i, u)- b_{0i}^{\alpha\beta}(\xi_i) \ll A(u)$ for suitable constants $b$ and $c$ from $(\ref{ncc2})$, equivalently
\begin{align}\label{ssc5}
s^\mu(u)\ll A(u),\,\,\,\,\,\,\,\,\,\varphi_i^\mu(\xi_i, u)- \xi_i\ll A(u),\,\,\,\,\,\,\,\,\,\,\,b_i^{\alpha\beta \mu}(\xi_i, u)- b_{0i}^{\alpha\beta}(\xi_i) \ll A(u),\,\,\,\,\,\,\,\,\,\, c_i^{\alpha\beta}(\xi_i, u)- c_{0i}^{\alpha\beta}(\xi_i) \ll A(u)
\end{align}
for $\mu=1,2,3,\cdots$. We will prove $(\ref{ssc5})$ by induction on $\mu=1,2,3,\cdots$. For $\mu=1$, since the linear term of $A(u)$ is $\frac{b}{16}\left( u_1+ \cdots + u_{r'} \right)$, the estimates $\left(\ref{ssc5} \right)_1$ holds if $b$ is sufficiently large. Let $\mu\geq  2$ and assume that the induction $\left(\ref{ssc5}\right)_{\mu-1}$ holds for $\mu-1$, i.e.
\begin{align*}
s^{\mu-1}(u)\ll A(u),\,\,\,\,\,\,\,\,\,\varphi_i^{\mu-1}(\xi_i, u)- \xi_i \ll A(u),\,\,\,\,\,\,\,\,\,\,\,b_i^{\alpha\beta(\mu-1)}(\xi_i, u) - b_{0i}^{\alpha\beta}(\xi_i) \ll A(u),\,\,\,\,\,\,\,\,\,\, c_i^{\alpha\beta(\mu-1)} - c_{0i}^{\alpha\beta} \ll A(u)
\end{align*}
We will prove that $(\ref{ssc5})_\mu$ holds. We estimate $(\ref{tc35}),(\ref{tc36}), (\ref{tc38})$ and $(\ref{nc9}), (\ref{nc10})$ in the following Lemma.
\begin{lemma}\label{ssc6}
\begin{align}
\Gamma_{ij|\mu}(\xi_i, u) &\ll \left( \frac{K_1}{b}+ \frac{K_2}{c}+ \frac{K_3 b}{c} \right)A(u) \,\,\,\,\,\,\,\,\textnormal{on}\,\,\, U_{ij}  \\
B_{ij|\mu}^{\alpha\beta} (\xi_i, u) &\ll \left( \frac{K_{4}}{b} + \frac{K_{5}}{c} + \frac{K_{6} b}{c} \right) A(u) \,\,\,\,\,\,\,\,\,\,\textnormal{on}\,\,\, U_{ij}  \\
\Pi_{i|\mu}^{\alpha\gamma} (\xi_i, u)  & \ll \left( \frac{K_{7} }{c} + \frac{K_{8} b}{c} \right) A(u) \,\,\,\,\,\,\,\textnormal{on}\,\,\, U_i^\delta\\
A_{ij|\mu}^{\alpha\beta} (\xi_i, u) &\ll \left( \frac{K_{9}}{b} + \frac{K_{10}}{c} + \frac{K_{11} b}{c} \right) A(u) \,\,\,\,\,\,\,\,\,\,\textnormal{on}\,\,\, U_{ij}  \\
W_{i|\mu}^{\alpha\gamma} (\xi_i, u)  & \ll \left( \frac{K_{12} }{c} + \frac{K_{13} b}{c} \right) A(u) \,\,\,\,\,\,\,\textnormal{on}\,\,\, U_i^\delta
\end{align}
where $K_1,K_2,K_3, \cdots$ are constants independent of $\mu$.
\end{lemma}
\begin{proof}
Lemma \ref{ssc6} follows from Lemma \ref{tcc31} and Lemma \ref{ncc20}.
\end{proof}
We recall Remark \ref{ssc7}. For any $\sigma=\left( \Pi, -W, \left( \Gamma, B, A \right) \right)=\left( \left\{ \Pi\right\}, \left\{ - W_i \right\}, \left\{ \left( \Gamma_{ij}, B_{ij}, A_{ij} \right) \right\} \right)\in C^0\left( \mathcal{U}, \mathscr{H}om_{\mathcal{O}_M}\left( \Theta_{\mathcal{F}_0}, \Theta_M \right) \right) $ $ \bigoplus C^0\left( \mathcal{U} , \mathscr{H}om_{\mathcal{O}_M} \left( \mathcal{N}_{\mathcal{F}_0}^*, \Omega_M^1 \right) \right) \bigoplus C^1\left( \mathcal{U}, \mathcal{E}_{\mathcal{F}_0} \right)$ which is a $1$-cocycle in the \v Cech resolution $(\ref{ssc12})$ of $\mathcal{E}_{\mathcal{F}_0}^\bullet$, where $\Pi_i \in \Gamma\left( U_i, \mathscr{H}om_{\mathcal{O}_M}\left( \Theta_{\mathcal{F}_0}, \Theta_M \right) \right)$ as in $(\ref{ssc9})$ and $\Gamma_{ij}\in \Gamma\left( U_{ij}, \Theta_M \right)$ and $B_{ij}\in \Gamma\left( U_{ij}, \mathscr{H}om_{\mathcal{O}_M} \left( \Theta_{\mathcal{F}_0} , \Theta_{\mathcal{F}_0} \right)  \right)$ as in $(\ref{ssc10})$ and $W_i\in \Gamma\left( U_i, \mathscr{H}om_{\mathcal{O}_M}\left( \mathcal{N}_{\mathcal{F}_0}^*, \Omega_M^1 \right) \right)$ as in $(\ref{ssc13})$ and $A_{ij}\in \Gamma\left( U_{ij}, \mathscr{H}om_{\mathcal{O}_M} \left( \mathcal{N}_{\mathcal{F}_0}^* , \mathcal{N}_{\mathcal{F}_0}^* \right) \right)$ as in $(\ref{ssc15})$, we define the norm $\left|\left| \sigma \right| \right|$ by 
\begin{align*}
||\sigma|| = \left|\left| \Gamma \right|\right| + \left|\left|B \right|\right| + \left|\left| \Pi \right|\right| + \left|\left| A \right| \right| + \left| \left| W \right| \right|
\end{align*}
where 
\begin{align*}
\left|\left| \Gamma \right|\right|=\max_{i,j} \sup_{\xi_i\in U_i\cap U_j} & \left| \Gamma_{ij} (\xi_i)\right| , \,\,\,\,\,\,\,\,\,\,\,\,\,\,\left|\left| B\right|\right|=\max_{i,j}\max_{\alpha,\beta}\sup_{\xi_i \in U_i\cap U_j}\left|B_{ij}^{\alpha\beta} \right|\,\,\,\,\,\,\,\,\,\,\,\,\,\,\left|\left| \Pi \right| \right| = \max_{i }\max_{\alpha,\gamma}\sup_{\xi_i\in U_i^\delta } \left| \Pi_i^{\alpha\gamma} (\xi_i) \right|\\
&\left|\left| A \right|\right|=\max_{i,j}\max_{\alpha,\beta}\sup_{\xi_j \in U_i\cap U_j}\left|A_{ij}^{\alpha\beta} \right|\,\,\,\,\,\,\,\,\,\,\,\,\,\,\left|\left| W \right| \right| = \max_{i }\max_{\alpha,\gamma}\sup_{\xi_i\in U_i^\delta } \left| W_i^{\alpha\gamma} (\xi_i) \right|
\end{align*}

\begin{lemma}\label{scc8}
For any quintuple $\sigma=\left( \Pi, -W , \left( \Gamma, B, A \right) \right)=\left(\left\{ \Pi_i \right\}, \left\{ - W_i\right\}, \left\{ \left( \Gamma_{ij}, B_{ij}, A_{ij} \right)\right\} \right)$ which is a $1$-cocycle in the \v Cech resolution of $\mathcal{E}_{\mathcal{F}_0}^\bullet$, we can find $\varphi_i(\xi_i), s^\lambda, b_i^{\alpha \gamma}(\xi_i)$ and $c_i^{\alpha \gamma}(\xi_i)$ satisfying
\begin{align}
\varphi_i- \varphi_j &= - \Gamma_{ij} + \sum_{\lambda=1}^m s^\lambda \rho_{ij\lambda} \\
\left[ \varphi_i, T_{0i}^\alpha \right] &= - \Pi_i^\alpha + \sum_{\lambda=1}^m s^\lambda \alpha_{i\lambda}^\alpha + \sum_{\gamma=1}^p b_{i}^{\alpha \gamma} T_{0i}^\gamma \\
\mathcal{L}_{\varphi_i}\left( w_{0i}^\alpha \right) - \sum_{\beta=1}^q & c_{i|\mu}^{\alpha\beta} w_{0i}'^\beta = W_i^\alpha  - \sum_{\lambda=1}^m s_\mu^\lambda \beta_{i\lambda}^\alpha \\
\left| \varphi_i(\xi_i)\right|\leq K \left| \left| \sigma\right|\right|,\,\,\,\,\,\,\,\,\,\,\,|s| \leq K & ||\sigma|| ,\,\,\,\,\,\,\,\,\,\,\,\, \left|b_i^{\alpha \gamma}(\xi_i)\right| \leq K ||\sigma||,\,\,\,\,\,\,\,\,\,\left| c_i^{\alpha \gamma} (\xi_i) \right| \leq K || \sigma ||
\end{align}
where $K$ is a constant independent of $\sigma$.
\end{lemma}

\begin{proof}
We can prove Lemma \ref{scc8} by combining the proof of Lemma \ref{tcc30} and Lemma \ref{ncc29} in a similar way.
\end{proof}

Then by Lemma \ref{ssc6} and Lemma \ref{scc8} we can choose solutions $\varphi_{i|\mu}(\xi_i,u)$ and $s_{\mu}^\lambda (u)$ and $b_{i|\mu}^{\alpha \beta}(\xi_i, u)$ and $c_{i|\mu}^{\alpha \beta}(\xi_i, u)$ of equations $(\ref{tc8})$ and $(\ref{tc10})$ and $(\ref{nc19})$ such that
\begin{align*}
\varphi_{i|\mu}(\xi_i, u)\ll K K^* A(u),\,\,\,\,\,\,\,\,\,\,\,\,\,s_\mu(u)\ll K K^* A(u),\,\,\,\,\,\,\,\,\,\,\,\,\, b_{i|\mu}^{\alpha \beta} \ll K K^* A(u),\,\,\,\,\,\,\,\,\,\, c_{i|\mu}^{\alpha\beta} \ll K K^* A(u)
\end{align*}
where $K^*=\frac{K_1+ K_{4}+ K_9}{b}+ \frac{K_2 + K_{5}+ K_{7} +K_{10}+ K_{12}}{c} + \frac{( K_3 + K_{6} + K_{8}+ K_{11} +  K_{13} )b}{c}$. We choose $b$ and $c$ in a way that $K K^* <1$. Then we have $\varphi_{i|\mu}(\xi_i, u) \ll A(u), s_\mu (u)\ll A(u)$ and $b_{i|\mu}^{\alpha\beta}(\xi_i, u) \ll A(u)$ and $c_{i|\mu}^{\alpha\beta}(\xi_i, u)\ll A(u)$. Hence the induction $(\ref{ssc5})_\mu$ holds for $\mu$. This completes the proof of Theorem \ref{sc1}.

\end{proof}

\section{Deformations of foliated complex analytic structures in terms of both tangent and cotangent sheaves 2} \label{cs9}

In this section we study deformations of foliated complex analytic structures defined by locally free subsheaves of tangent and cotangent sheaves at the same time in terms of both tangent and cotangent sheaves by using the complex of sheaves $\mathcal{F}_0'^\bullet$ (see Appendix \ref{}) which is isomorphic to $\mathcal{F}_0^\bullet$. We recall Definition \ref{s1}. We keep the notations in section \ref{ss1}. In this case, we replace the integrability condition (6) $d_{\mathcal{A}/B}\left( \mathcal{N}_{\mathcal{F}}^* \right)\subset \mathcal{N}_\mathcal{F}^* \bigwedge \Omega_{\mathcal{M}/B}^1$ at every point of $\mathcal{M}-\bigcup_{t\in B} S_t$, where $t= \textnormal{Sing}\left(\mathcal{F}_t \right)$ in Definition \ref{s1} by the condition $(6)'$ $d_{\mathcal{M}/B}\left(\mathcal{N}_\mathcal{F}^*\right)$ vanishes on $\bigwedge^2 \Theta_\mathcal{F}$. We assume that $\Theta_\mathcal{F}$ and $\mathcal{N}_\mathcal{F}^*$ are both locally free. We keep the notation in Remark \ref{ta1} and Remark \ref{n5} and Remark \ref{ss17}. Let us consider the condition $(6)'$ which implies that for $\alpha,\beta=1,...,p$ and $\gamma=1,...,q$,
\begin{align}\label{ss2}
i_{T_j^\alpha(z_j,t)\wedge T_j^\beta(z_j,t)}\left( dw_j^\gamma(z_j,t) \right)= 0
\end{align}
Then we can write
\begin{align}\label{ss3}
i_{T_j^\alpha(z_j,t)}\left( dw_j^\gamma(z_j,t) \right) = \sum_{\sigma=1}^q a_{j\alpha}^{\gamma \sigma}(z_j,t) w_j^\sigma(z_j,t)\,\,\,\,\,\,\,\,\,\textnormal{for some}\,\,\,a_{j\alpha}^{\gamma\eta}(z_j,t)\in \Gamma\left( \mathcal{U}_j, \mathcal{O}_{\mathcal{M}} \right)
\end{align}
Let $\left(\mathcal{M}, \Theta_\mathcal{F}, \mathcal{N}_\mathcal{F}^*, B, \pi \right)$ be a foliated analytic family with both $\Theta_{\mathcal{F}}$ and $\mathcal{N}_\mathcal{F}^*$ locally free as in Definition \ref{s1} with $(6)'$, so that $\mathcal{F}_t=\left(\Theta_{\mathcal{F}_t}, \mathcal{N}_{\mathcal{F}_t}^* \right)$ defines a (singular) holomorphic foliation on a compact complex manifold $M_t$ with $\Theta_{\mathcal{F}_t}$ and $\mathcal{N}_{\mathcal{F}_t}^*$ locally free. Then we have the complex of sheaves $\mathcal{F}_t'^\bullet$ on $M_t$ associated to $\mathcal{F}_t=\left( \Theta_{\mathcal{F}_t}, \mathcal{N}_{\mathcal{F}_t}^* \right)$ (see Appendix \ref{AppendixA4}).
\begin{center}
$\mathcal{F}_t'^\bullet: \begin{CD}
\cdots \\
@AF_3'^tAA \\
\mathscr{H}om_{\mathcal{O}_{M_t}}\left(\bigwedge^3 \Theta_{\mathcal{F}_t},  \frac{\Theta_{M_t}}{\Theta_{\mathcal{F}_t}} \right) \bigoplus \mathscr{H}om_{\mathcal{O}_{M_t}}\left(\mathcal{N}_{\mathcal{F}_t}^*, \bigwedge^3 \Theta_{\mathcal{F}_t}\right) \bigoplus \mathscr{H}om_{\mathcal{O}_{M_t}}\left(\mathcal{N}_{\mathcal{F}_t}^*,\bigwedge^2 \Theta_{\mathcal{F}_t}^* \right)       \\
@AF_2'^t AA \\
\mathscr{H}om_{\mathcal{O}_{M_t}}\left(\bigwedge^2 \Theta_{\mathcal{F}_t}, \frac{\Theta_{M_t}}{\Theta_{\mathcal{F}_t} }\right) \bigoplus \mathscr{H}om_{\mathcal{O}_{M_t}}\left(\mathcal{N}_{\mathcal{F}_t}^*, \bigwedge^2 \Theta_{\mathcal{F}_t}^* \right) \bigoplus \mathscr{H}om_{\mathcal{O}_{M_t}}\left(\mathcal{N}_{\mathcal{F}_t}^*, \Theta_{\mathcal{F}_t}^* \right)       \\
@AF_1'^t AA \\
\mathscr{H}om_{\mathcal{O}_{M_t}}\left(\Theta_{\mathcal{F}_t}, \frac{\Theta_{M_t}}{\Theta_{\mathcal{F}_t}} \right) \bigoplus \mathscr{H}om_{\mathcal{O}_M}\left( \mathcal{N}_{\mathcal{F}_t}^*, \frac{\Omega_{M_t}^1}{\mathcal{N}_{\mathcal{F}_t}^* } \right)  \\
@AF_0^t AA \\
\Theta_{M_t}
\end{CD}$
\end{center}
We will denote the $i$-th cohomology group by $\mathbb{H}^1\left( M_t,  \mathcal{F}_t'^\bullet \right)$. We can compute $\mathbb{H}^i\left( M_t, \mathcal{F}_t'^\bullet \right)$ by the following \v Cech resolution of $\mathcal{F}_t'^\bullet$. Here $\delta$ is the \v Cech map and $\mathcal{U}_t=\mathcal{U}\cap M_t=\left\{ U_j^t:=\mathcal{U}_j\cap M_t|j=1,2,... \right\}$ is an open covering of $M_t$:
{\tiny{\begin{center}
$\begin{CD}
\cdots \\
@AF_2'^tAA \\
C^0\left( \mathcal{U}_t, \left(\bigwedge^2 \Theta_{\mathcal{F}_t}^*\otimes \frac{\Theta_{M_t}}{\Theta_{\mathcal{F}_t}}\right) \bigoplus \left(\left(\mathcal{N}_{\mathcal{F}_t}^* \right)^*\otimes \bigwedge^2 \Theta_{\mathcal{F}_t}^* \right) \bigoplus \left( \left( \mathcal{N}_{\mathcal{F}_t}^* \right)^*\otimes \Theta_{\mathcal{F}_t}^* \right)   \right) @>-\delta>> \cdots \\
@AF_1'^tAA @AF_1'^tAA \\
C^0\left(\mathcal{U}_t, \left(\Theta_{\mathcal{F}_t}^*\otimes \frac{\Theta_{M_t}}{\Theta_{\mathcal{F}_t}}\right) \bigoplus \left(\left(\mathcal{N}_{\mathcal{F}_t}^* \right)^*\otimes \frac{\Omega_{M_t}^1}{\mathcal{N}_{\mathcal{F}_t}^*} \right) \right) @>\delta>> C^1\left(\mathcal{U}_t, \left(\Theta_{\mathcal{F}_t}^*\otimes \frac{\Theta_{M_t}}{\Theta_{\mathcal{F}_t}}\right) \bigoplus \left(\left(\mathcal{N}_{\mathcal{F}_t}^* \right)^*\otimes \frac{\Omega_{M_t}^1}{\mathcal{N}_{\mathcal{F}_t}^*} \right) \right) @>-\delta>> \cdots \\
@AF_0^tAA @AF_0^tAA  @AF_0^tAA \\
C^0\left(\mathcal{U}_t, \Theta_{M_t} \right) @>-\delta >> C^1\left( \mathcal{U}_t, \Theta_{M_t} \right) @>\delta>> C^2\left(\mathcal{U}_t, \Theta_{M_t} \right)  
\end{CD}$
\end{center}}}

We will relate the first cohomology group $\mathbb{H}^1\left( M_t, \mathcal{F}_t^\bullet \right)$ to infinitesimal foliated deformations of $\pi^{-1}(t)= \left( M_t, \Theta_{\mathcal{F}_t}, \mathcal{N}_{\mathcal{F}_t}'^*  \right)$ in the foliated analytic family $\left( \mathcal{M}, \Theta_\mathcal{F}, \mathcal{N}_\mathcal{F}^*, B, \pi  \right)$ in terms of both tangent sheaves and cotangent sheaves. By taking the derivative of $(\ref{ss2})$ with respect to $t$, we have
\begin{align*}
i_{\frac{\partial T_j^\alpha(z_j,t)}{\partial t} \wedge T_j^\beta(z_j,t)} \left(  dw_j^\gamma(z_j,t)\right) + i_{T_j^\alpha(z_j,t) \wedge \frac{\partial T_j^\beta(z_j,t)}{\partial t}}\left( dw_j^\gamma \right) + i_{T_j^\alpha(z_j,t)\wedge T_j^\beta(z_j,t)}\left( d\left( \frac{\partial w_j^\gamma(z_j,t)}{\partial t}   \right) \right)= 0
\end{align*}
Equivalently we have from $(\ref{ss3})$
\begin{align}\label{ss4}
- \sum_{\sigma=1}^q a_{j\beta}^{\gamma \sigma}(z_j,t) i_{\frac{\partial T_j^\alpha (z_j,t) }{\partial t}} \left( w_j^\gamma(z_j,t) \right) + \sum_{\sigma=1}^q a_{j\alpha}^{\gamma \sigma}(z_j,t) i_{\frac{\partial T_j^\beta(z_j,t)}{\partial t}}\left( w_j^\gamma(z_j,t) \right) + i_{T_j^\alpha(z_j,t) \wedge T_j^\beta(z_j,t)}\left( d\left( \frac{\partial w_j^\gamma(z_j,t)}{\partial t} \right) \right) = 0
\end{align}
Since from $(\ref{ss18})$, we have $i_{\frac{\partial T_j^\alpha(z_j,t)}{\partial t}}\left(  w_j^\gamma(z_j,t) \right) + i_{T_j^\alpha(z_j,t)}\left( \frac{\partial w_j^\gamma(z_j,t)}{\partial t} \right)=0$, $(\ref{ss4})$ is equivalent to
\begin{align*}
 \sum_{\sigma=1}^q a_{j\beta}^{\gamma \sigma}(z_j,t) i_{T_j^\alpha(z_j,t)} \left( \frac{\partial w_j^\gamma(z_j,t) }{\partial t}\right) - \sum_{\sigma=1}^q a_{j\alpha}^{\gamma \sigma}(z_j,t) i_{T_j^\beta(z_j,t) }\left( \frac{\partial w_j^\gamma(z_j,t)}{\partial t} \right) + i_{T_j^\alpha(z_j,t) \wedge T_j^\beta(z_j,t)}\left( d\left( \frac{\partial w_j^\gamma(z_j,t)}{\partial t} \right) \right) = 0
\end{align*}
This implies that $E_1'\left( \beta_j(t) \right)=0$. Hence $\left( \left\{ \alpha_j(t), \beta_j(t) \right\} \right)$ defines a $1$-cocycle in the above \v Cech resolution of $\mathcal{F}_t'^\bullet$, so that an element in $\mathbb{H}^1\left( M_t, \mathcal{F}_t'^\bullet \right)$. Then we can define the foliated Kodaira-Spencer map $\varphi_0:T_0 B \to \mathbb{H}^1\left(  M_t, \mathcal{F}_t'^\bullet \right)$ in terms of both tangent and cotangent sheaves as in Definition \ref{s8}.

\section{Theorem of existence of deformations of foliated complex analytic structures in terms of both tangent and cotangent sheaves 2}

We shall keep the notation in subsection \ref{ss5}. By using the complex of sheaves $\mathcal{F}_0'^\bullet$ instead of $\mathcal{F}_0^\bullet$, we prove
\begin{theorem}[Theorem of existence for deformations of foliated complex analytic structures in terms of both tangent and cotangent sheaves]\label{dt6}
Let $\left( M , \mathcal{F}_0 \right)$ be a compact foliated complex manifold with both $\Theta_{\mathcal{F}_0}$ and $\mathcal{N}_{\mathcal{F}_0}^*$ locally free. Suppose that $\mathbb{H}^2\left( M, \mathcal{F}_0'^\bullet \right)=0$. Then there exists a foliated analytic family $\left( \mathcal{M}, \Theta_{\mathcal{F}}, \mathcal{N}_\mathcal{F}^* , B, \pi \right)$ with $0\in B \subset \mathbb{C}^m$ satisfying the following conditions$:$
\begin{enumerate}
\item $\pi^{-1}(0)=\left( M, \mathcal{F}_0 \right)$
\item The foliated Kodaira-Spencer map $\varphi_0: T_0 B \to \mathbb{H}^1\left( M, \mathcal{F}_0'^{ \bullet} \right)$ in terms of both tangent and cotangent sheaf is an isomorphism.
\end{enumerate}
\end{theorem}

\begin{proof}

We shall keep the notations in the proof of Theorem \ref{tt2} and Theorem \ref{nn2} and Theorem \ref{ss6}.

\subsection{Existence of formal solutions}\

We construct  solutions of $(\ref{t10})-(\ref{te5})$ and $(\ref{n20})-(\ref{n212}), (\ref{n21})$ and $(\ref{se1})$ which are power series in $t$ and additionally
\begin{align}\label{ss19}
i_{T_i^\alpha \wedge T_i^\beta}\left(  \partial w_i^\gamma \right) = 0, \,\,\,\,\,\,\alpha,\beta=1,...,p,\,\,\,\gamma=1,...,q
\end{align}
 Then $(\ref{t10})-(\ref{te5})$ are equivalent to the systems of congruences $(\ref{t13})_\mu-(\ref{te6})_\mu$ for $\mu=1,2,3,\cdots$ and $(\ref{n22})-(\ref{n212}),(\ref{n21})$ are equivalent to the system of congruences $(\ref{n23})_\mu- (\ref{n213})_\mu, (\ref{n24})_\mu$ for $\mu=1,2,3,\cdots$ and $(\ref{se1})$ are equivalent to the systems of congruences $(\ref{se3})_\mu$ for $\mu=1,2,3,...$, and additionally $(\ref{ss19})$ is equivalent to the system of congruences
\begin{align}\label{ss9}
i_{T_i^{\alpha \mu} \wedge T_i^{\beta \mu}} \left(  \partial w_i^{\gamma \mu} \right)\equiv_\mu 0,\,\,\,\,\,\,\,\alpha , \beta=1,...,p,\,\,\,\gamma=1,...,q
\end{align}
for $\mu=1,2,3, \cdots$. We define homogenous polynomials $(\ref{tt3})-(\ref{tt4})$ and $(\ref{nn3})-(\ref{n215})$, and $(\ref{se6})$, and additionally 
\begin{align} \label{ss20}
X_{i\alpha\beta |\mu}^\gamma &\equiv_\mu i_{T_i^{\alpha (\mu-1)} \wedge T_i^{\beta(\mu-1)}}\left( \partial w_{i}^{\gamma(\mu-1)} \right) 
\end{align}

\begin{lemma}
We have the equalities $(\ref{t19})-(\ref{te7})$ and $(\ref{n25})-(\ref{n28})$ and $(\ref{se2}),(\ref{se4})$ and additionally
\begin{align}\label{ss7}
X_{i\alpha\beta|\mu}^\gamma = i_{T_{0i}^\alpha}\partial K_{i|\mu}^{\beta \gamma} - i_{\Pi_{i|\mu}^{\alpha\beta}}\left( w_{0i}^\gamma \right) - \sum_{\xi=1}^p g_{0i\alpha\beta}^\xi K_{i|\mu}^{\xi \gamma} - i_{T_{0i}^\beta}\partial K_{i|\mu}^{\alpha \gamma}
\end{align}

\end{lemma}

\begin{proof}
$(\ref{t19})-(\ref{te7})$ follows from Lemma \ref{te10}, and $(\ref{n25})-(\ref{n28})$ follows from Lemma \ref{ne10} and $(\ref{se2}),(\ref{se4})$ follows from Lemma \ref{ss8}. We prove $(\ref{ss7})$. In fact, from $(\ref{ss20})$ and $(\ref{tt121})$ and $(\ref{se6})$, we have
\begin{align*}
X_{i\alpha\beta |\mu}^\gamma &\equiv_\mu i_{T_i^{\alpha (\mu-1)} \wedge T_i^{\beta(\mu-1)}}\left( \partial w_{i}^{\gamma(\mu-1)} \right) = i_{T_i^{\beta(\mu-1)} } i_{T_i^{\alpha(\mu-1)}}\left(  \partial w_{i}^{\gamma(\mu-1)} \right) \\
& \equiv_\mu i_{T_i^{\beta(\mu-1)}}\mathcal{L}_{T_i^{\alpha(\mu-1)}}\left( w_{i}^{\alpha(\mu-1)} \right)- i_{T_i^{\beta(\mu-1)}} \partial i_{T_i^{\alpha(\mu-1)}} \left( w_i^{\gamma(\mu-1)} \right)\\
&\equiv_\mu \mathcal{L}_{T_i^{\alpha(\mu-1)}} i_{T_i^{\beta(\mu-1)}}\left( w_i^{\gamma(\mu-1)} \right) - i_{\left[ T_i^{\alpha(\mu-1)}, T_i^{\beta(\mu-1)} \right]}\left( w_i^{\gamma(\mu-1)} \right) - i_{T_{0i}^\beta}\partial K_{i|\mu}^{\alpha \gamma}\\
&\equiv_\mu \mathcal{L}_{T_{0i}^\alpha} K_{i|\mu}^{\beta \gamma} - i_{\Pi_{i|\mu}^{\alpha\beta}+ \sum_{\xi=1}^p g_{i\alpha\beta}^{\xi(\mu-1)} T_i^{\xi(\mu-1)}}\left( w_i^{\gamma(\mu-1)} \right)  - i_{T_{0i}^\beta}\partial K_{i|\mu}^{\alpha \gamma}\\
&\equiv_\mu i_{T_{0i}^\alpha}\partial K_{i|\mu}^{\beta \gamma} - i_{\Pi_{i|\mu}^{\alpha\beta}}\left( w_{0i}^\gamma \right) - \sum_{\xi=1}^p g_{0i\alpha\beta}^\xi K_{i|\mu}^{\xi \gamma} - i_{T_{0i}^\beta}\partial K_{i|\mu}^{\alpha \gamma}
\end{align*}
\end{proof}

Our purpose is to construct $\varphi^\mu=\varphi^{\mu-1}+ \varphi_\mu, r_{ij}^{\alpha\beta\mu}= r_{ij}^{\alpha\beta(\mu-1)}+ r_{ij|\mu}^{\alpha\beta}, T_i^{\alpha\mu}=T_i^{\alpha(\mu-1)}+ T_{i|\mu}^\alpha$,  $g_{i\alpha\beta}^{\gamma \mu}= g_{i\alpha\beta}^{\gamma(\mu-1)}+ g_{i\alpha\beta|\mu}^\gamma$ and $h_{ij}^{\alpha\beta \mu}= h_{ij}^{\alpha\beta (\mu-1)} + h_{ij|\mu}^{\alpha\beta}, w_i^{\alpha \mu}=w_i^{\alpha (\mu-1)}+ w_{i|\mu}^\alpha$ satisfying $(\ref{t13})_\mu-(\ref{te6})_\mu$ and $(\ref{n23})_\mu-(\ref{n213})_\mu, (\ref{n24})_\mu$ and $(\ref{se3})_\mu$ and additionally $(\ref{ss9})_\mu$.

\begin{lemma}
$(\ref{t13})_\mu-(\ref{te6})_\mu$ are equivalent to $(\ref{te11})-(\ref{t33})$, and $(\ref{n23})_\mu-(\ref{n213})_\mu, (\ref{n24})_\mu$ are equivalent to $(\ref{n42}) -(\ref{n46}) , (\ref{nn121})$ and $(\ref{se3})_\mu$ is equivalent to $(\ref{se5})$ and additionally $(\ref{ss9}_\mu)$ is equivalent to
\begin{align}\label{ss11}
- X_{i\alpha\beta|\mu}^\gamma =   \sum_{\sigma=1}^q a_{0i\alpha}^{\gamma \sigma} i_{T_{i|\mu}^\beta}\left( w_{0i}^\sigma \right) - \sum_{\sigma=1}^q a_{0i\beta}^{\gamma \sigma} i_{T_{i|\mu}^\alpha}\left( w_{0i}^\sigma \right) + i_{T_{0i}^\alpha \wedge T_{0i}^\beta}\left(\partial w_{i|\mu}^\gamma \right)
\end{align}

\end{lemma}

\begin{proof}
$(\ref{te11})-(\ref{t33})$ follows from Lemma \ref{tt5}, and $(\ref{n42})-(\ref{n46}) , (\ref{nn121})$ follows from Lemma \ref{nn5} and $(\ref{se5})$ follows from Lemma \ref{ss10}. It remains to prove $(\ref{ss11})$. In fact,  
\begin{align*}
&i_{\left( T_i^{\alpha(\mu-1)} + T_{i|\mu}^\alpha \right) \wedge \left( T_i^{\beta(\mu-1)} + T_{i|\mu}^\beta \right) } \left(  \partial \left( w_i^{\gamma(\mu-1)} + w_{i|\mu}^\gamma \right)   \right)\\
& \equiv_\mu X_{i\alpha\beta | \mu }^\gamma + i_{T_{0i}^\alpha \wedge T_{i|\mu}^\beta}\left( d w_{0i}^\gamma \right) + i_{T_{i|\mu}^\alpha \wedge T_{0i}^\beta}\left( d w_{0i}^\gamma \right) + i_{T_{0i}^\alpha \wedge T_{0i}^\beta}\left(\partial w_{i|\mu}^\gamma \right) \\ &\equiv_\mu X_{i\alpha\beta |\mu}^\gamma + \sum_{\sigma=1}^q a_{0i\alpha}^{\gamma \sigma} i_{T_{i|\mu}^\beta}\left( w_{0i}^\sigma \right) - \sum_{\sigma=1}^q a_{0i\beta}^{\gamma \sigma} i_{T_{i|\mu}^\alpha}\left( w_{0i}^\sigma \right) + i_{T_{0i}^\alpha \wedge T_{0i}^\beta}\left(\partial w_{i|\mu}^\gamma \right)= 0
\end{align*}

\end{proof}

\begin{lemma} \label{ss24}
Under the hypothesis $\mathbb{H}^2\left( M, \mathcal{F}_0'^\bullet \right)=0$,  we can find $\varphi_\mu, r_{ij|\mu}^{\alpha\beta}, T_{i|\mu}^\alpha, g_{i\alpha\beta |\mu}^\gamma$ and $h_{ij|\mu}^{\alpha\beta}, w_{i|\mu}^\alpha$ which satisfy $(\ref{te11})-(\ref{t33})$ and $(\ref{n42})-(\ref{n46}),(\ref{nn121})$ and $(\ref{se5})$ and additionally $(\ref{ss11})$.
\end{lemma}

\begin{proof}
We keep the notations in the proof of Lemma \ref{te47} and Lemma \ref{ne25} and Lemma \ref{ss12}. We shall define a global section in $A^{0,0}\left( M, \mathscr{H}om_{\mathcal{O}_M} \left( \mathcal{N}_{\mathcal{F}_0}^*, \bigwedge^2 \Theta_{\mathcal{F}_0}^* \right)  \right)$. First we recall that we can write $i_{T_{0i}^\alpha}\left( w_{0i}^\gamma \right) = \sum_{\eta=1}^q a_{0i\alpha}^{\gamma \eta} w_{0i}^\eta$ for some $a_{0i\alpha}^{\gamma \eta} \in \Gamma\left(U_i, \mathcal{O}_M \right)$. Then we have
\begin{align}
&w_{0i}^\gamma= \sum_{\delta=1}^q  h_{0ij}^{\gamma \delta} w_{0j}^\delta \notag \Longrightarrow dw_{0i}^\gamma= \sum_{\delta=1}^q d h_{0ij}^{\gamma \delta} \wedge w_{0j}^\delta + \sum_{\delta=1}^q h_{0ij}^{\gamma\delta} dw_{0j}^\delta  \notag \\
& \Longrightarrow i_{T_{0i}^\alpha}\left( dw_{0i}^\gamma \right)= \sum_{\delta=1}^q i_{T_{0i}^\alpha}\left(d h_{0ij}^{\gamma \delta} \right) w_{0j}^\delta + \sum_{\delta=1}^q h_{0ij}^{\gamma \delta} i_{T_{0i}^\alpha} \left( dw_{0j}^\delta \right) \notag \\
&\Longrightarrow \sum_{\eta=1}^q a_{0i\alpha}^{\gamma \eta} w_{0i}^\eta = \sum_{\delta=1}^q i_{T_{0i}^\alpha}\left( d h_{0ij}^{\gamma \delta}\right) w_{0j}^\delta + \sum_{\delta, \xi=1}^q \sum_{\beta=1}^p h_{0ij}^{\gamma \delta} r_{0ij}^{\alpha\beta}  a_{0j\beta}^{\delta \xi} w_{0j}^\xi  \notag \\
&\Longrightarrow \sum_{\eta=1}^q a_{0i\alpha}^{\gamma \eta} h_{0ij}^{\eta \xi}  = i_{T_{0i}^\alpha}\left( d h_{0ij}^{\gamma \xi} \right) + \sum_{\delta=1}^q \sum_{\beta=1}^p h_{0ij}^{\gamma \delta} r_{0ij}^{\alpha\beta} a_{0j\beta}^{\delta \xi} \label{ss21}
\end{align}

Then from $(\ref{ss20})$ and $(\ref{tt6})$ and $(\ref{nn6})$ and $(\ref{ss21})$, we have

{\tiny{\begin{align*}
&X_{i\alpha\beta |\mu}^\gamma \equiv_\mu i_{T_i^{\alpha (\mu-1)} \wedge T_i^{\beta(\mu-1)}}\left( \partial w_{i}^{\gamma(\mu-1)} \right) \equiv_\mu i_{ \left(\Gamma_{ij|\mu}^\alpha+ \sum_{\xi=1}^p r_{ij}^{\alpha \xi(\mu-1)} T_j^{\xi(\mu-1)}\right) \wedge \left( \Gamma_{ij|\mu}^\beta + \sum_{\eta=1}^p r_{ij}^{\beta \eta(\mu-1)} T_j^{\eta (\mu-1)}  \right)}\left( \partial\left( C_{ij|\mu}^\gamma + \sum_{\delta=1}^q h_{ij}^{\gamma\delta(\mu-1)} w_j^{\delta(\mu-1)}  \right)    \right)\\
&\equiv_\mu i_{\Gamma_{ij|\mu}^\alpha \wedge T_{0i}^\beta}\left( dw_{0i}^\gamma \right) + i_{T_{0i}^\alpha \wedge \Gamma_{ij|\mu}^\beta}\left( dw_{0i}^\gamma \right) + i_{T_{0i}^\alpha \wedge T_{0i}^\beta}\left( \partial C_{ij|\mu}^\gamma \right) + \sum_{\xi, \eta=1}^p \sum_{\delta=1}^q r_{ij}^{\alpha \xi(\mu-1)} r_{ij}^{\beta \eta(\mu-1)} h_{ij}^{\gamma \delta(\mu-1)} i_{T_i^{\xi(\mu-1)}\wedge T_i^{\eta(\mu-1)}}\left( \partial w_j^{\delta (\mu-1)} \right)\\
&+ \sum_{\xi, \eta=1}^p \sum_{\delta=1}^q r_{ij}^{\alpha \xi(\mu-1)} r_{ij}^{\beta \eta(\mu-1)} i_{T_i^{\xi(\mu-1)}\wedge T_i^{\eta(\mu-1)}}\left(\partial h_{ij}^{\gamma \delta(\mu-1)} \wedge w_j^{\delta(\mu-1)}   \right)\\
&\equiv_\mu i_{\Gamma_{ij|\mu}^\alpha \wedge T_{0i}^\beta}\left( dw_{0i}^\gamma \right) + i_{T_{0i}^\alpha \wedge \Gamma_{ij|\mu}^\beta}\left( dw_{0i}^\gamma \right) + i_{T_{0i}^\alpha \wedge T_{0i}^\beta}\left( \partial C_{ij|\mu}^\gamma \right)  + \sum_{\xi, \eta=1}^p \sum_{\delta=1}^q r_{0ij}^{\alpha \xi} r_{0ij}^{\beta \eta} h_{0ij}^{\gamma \delta} X_{j\xi\eta|\mu}^\delta + \sum_{\xi, \eta=1}^p \sum_{\delta=1}^q r_{0ij}^{\alpha \xi} r_{0ij}^{\beta \eta}  i_{T_{0i}^\xi}\left( dh_{0ij}^{\gamma \delta} \right) K_{j|\mu}^{\eta \delta}\\
&- \sum_{\xi, \eta=1}^p \sum_{\delta=1}^q r_{0ij}^{\alpha \xi} r_{0ij}^{\beta \eta}  i_{T_{0i}^\eta}\left( d h_{0ij}^{\gamma \delta} \right) K_{j|\mu}^{\xi \delta} \\
&\equiv_\mu i_{\left(\sum_{\xi=1}^p r_{0ij}^{\alpha \xi} \Gamma_{j|\mu}^\xi - \Gamma_{i|\mu}^\alpha \right)\wedge T_{0i}^\beta}\left( dw_{0i}^\gamma \right) + i_{T_{0i}^\alpha \wedge \left(\sum_{\eta=1}^p r_{0ij}^{\beta \eta} \Gamma_{j|\mu}^\eta - \Gamma_{i|\mu}^\beta \right)}\left( dw_{0i}^\gamma \right) + i_{T_{0i}^\alpha \wedge T_{0i}^\beta}\left(\partial \left( \sum_{\delta=1}^q h_{0ij}^{\gamma \delta } C_{j|\mu}^\delta \right) - \partial C_{i|\mu}^\gamma \right)\\
&+ \sum_{\xi, \eta=1}^p \sum_{\delta=1}^q r_{0ij}^{\alpha \xi} r_{0ij}^{\beta \eta} h_{0ij}^{\gamma \delta} X_{j\xi\eta|\mu}^\delta + \sum_{\xi, \eta=1}^p \sum_{\delta=1}^q r_{0ij}^{\alpha \xi} r_{0ij}^{\beta \eta}  i_{T_{0j}^\xi}\left( dh_{0ij}^{\gamma \delta} \right) K_{j|\mu}^{\eta \delta}- \sum_{\xi, \eta=1}^p \sum_{\delta=1}^q r_{0ij}^{\alpha \xi} r_{0ij}^{\beta \eta}  i_{T_{0j}^\eta}\left( d h_{0ij}^{\gamma \delta} \right) K_{j|\mu}^{\xi \delta}\\
&\equiv_\mu - \sum_{\xi, \eta=1}^p \sum_{\delta=1}^q r_{0ij}^{\alpha \xi} r_{0ij}^{\beta \eta} i_{\Gamma_{j|\mu}^\xi} i_{T_{0j}^\eta}\left(d h_{0ij}^{\gamma \delta} \wedge w_{0j}^\delta + h_{0ij}^{\gamma \delta} dw_{0j}^\delta \right) + i_{\Gamma_{i|\mu}^\alpha i_{T_{0i}^\beta}}\left( dw_{0i}^\gamma \right) + \sum_{\xi, \eta=1}^p \sum_{\delta=1}^q r_{0ij}^{\alpha \xi} r_{0ij}^{\beta \eta}i_{\Gamma_{j|\mu}^\eta} i_{T_{0j}^\xi} \left( dh_{0ij}^{\gamma \delta} \wedge w_{0j}^\delta + h_{0ij}^{\gamma \delta} dw_{0j}^\delta \right) \\
& - i_{\Gamma_{i|\mu}^\beta}i_{T_{0i}^\alpha}\left( dw_{0i}^\gamma \right) + \sum_{\xi, \eta=1}^p r_{0ij}^{\alpha \xi} r_{0ij}^{\beta \eta} h_{0ij}^{\gamma \delta} i_{T_{0j}^\xi \wedge T_{0j}^\eta}\left( \partial C_{j|\mu}^\delta \right) + \sum_{\delta=1}^q i_{T_{0i}^\alpha \wedge T_{0i}^\beta}\left( dh_{0ij}^{\gamma \delta} \wedge C_{j|\mu}^\delta   \right) - i_{T_{0i}^\alpha \wedge T_{0i}^\beta }\left( \partial C_{i|\mu}^\gamma \right) + \sum_{\xi, \eta=1}^p \sum_{\delta=1}^q r_{0ij}^{\alpha \xi} r_{0ij}^{\beta \eta} h_{0ij}^{\gamma \delta} X_{j\xi\eta|\mu}^\delta\\
&- \sum_{\xi, \eta=1}^p \sum_{\delta,\sigma=1}^q r_{0ij}^{\alpha \xi} r_{0ij}^{\beta \eta} h_{0ij}^{\gamma \sigma}  a_{0j\xi}^{\sigma \delta} K_{j|\mu}^{\eta \delta} +\sum_{\eta=1}^p \sum_{\delta, \sigma=1}^q   a_{0i\alpha}^{\gamma \sigma} r_{0ij}^{\beta \eta} h_{0ij}^{\sigma \delta} K_{j|\mu}^{\eta \delta} + \sum_{\xi, \eta=1}^p \sum_{\delta,\sigma=1}^q r_{0ij}^{\alpha \xi} r_{0ij}^{\beta \eta} h_{0ij}^{\gamma \sigma}  a_{0j\eta}^{\sigma \delta}  K_{j|\mu}^{\xi \delta} - \sum_{\eta=1}^p \sum_{\delta, \sigma=1}^q   a_{0i\beta}^{\gamma \sigma} r_{0ij}^{\alpha \xi} h_{0ij}^{\sigma \delta} K_{j|\mu}^{\xi \delta} \\
&\equiv_\mu  \cancel{\sum_{\xi, \eta=1}^p \sum_{\delta, \sigma=1}^q r_{0ij}^{\alpha \xi} r_{0ij}^{\beta \eta}  h_{0ij}^{\gamma \sigma} a_{0j\eta}^{\sigma \delta}  i_{\Gamma_{j|\mu}^\xi}\left( w_{0j}^\delta \right)} - \sum_{\xi=1}^p \sum_{\delta, \sigma=1}^q  a_{0i\beta}^{\gamma \sigma}  r_{0ij}^{\alpha \xi} h_{0ij}^{\sigma \delta} i_{\Gamma_{j|\mu}^\xi}\left( w_{0i}^\delta \right) - \cancel{\sum_{\xi, \eta=1}^p \sum_{\delta=1}^q r_{0ij}^{\alpha \xi} r_{0ij}^{\beta \eta} h_{0ij}^{\gamma \delta} a_{0j\xi}^{\delta \sigma}   i_{\Gamma_{j|\mu}^\eta }\left( w_{0j}^\sigma \right)} \\
&+ \sum_{\delta=1}^q a_{0i\beta}^{\gamma \delta}  i_{\Gamma_{i|\mu}^\alpha}\left( w_{0i}^\delta \right) - \cancel{\sum_{\xi, \eta=1}^p \sum_{\delta, \sigma=1}^q r_{0ij}^{\alpha \xi} r_{0ij}^{\beta \eta}  h_{0ij}^{\gamma \sigma} a_{0j\xi}^{\sigma \delta}  i_{\Gamma_{j|\mu}^\eta }\left( w_{0j}^\delta \right) } + \sum_{\eta=1}^p \sum_{\delta, \sigma=1}^q  a_{0i\alpha}^{\gamma \sigma}   r_{0ij}^{\beta \eta} h_{0ij}^{\sigma \delta} i_{\Gamma_{j|\mu}^\eta}\left( w_{0i}^\delta \right)\\
& + \cancel{ \sum_{\xi, \eta=1}^p \sum_{\delta=1}^q r_{0ij}^{\alpha \xi} r_{0ij}^{\beta \eta} h_{0ij}^{\gamma \delta} a_{0j\eta}^{\delta \sigma}   i_{\Gamma_{j|\mu}^\xi}\left( w_{0j}^\sigma \right) }  - \sum_{\delta=1}^q  a_{0i \alpha}^{\gamma \delta}  i_{\Gamma_{i|\mu}^\beta}\left( w_{0i}^\delta \right)  +  \sum_{\xi, \eta=1}^p r_{0ij}^{\alpha \xi} r_{0ij}^{\beta \eta} h_{0ij}^{\gamma \delta} i_{T_{0j}^\xi \wedge T_{0j}^\eta}\left( \partial C_{j|\mu}^\delta \right) + \sum_{\xi, \eta=1}^p\sum_{\delta=1}^q r_{0ij}^{\alpha \xi} r_{0ij}^{\beta \eta} i_{T_{0j}^\xi} \left( dh_{0ij}^{\gamma \delta} \right) i_{T_{0j}^\eta}\left( C_{j|\mu}^\delta \right)\\
& -  \sum_{\xi, \eta=1}^p  \sum_{\delta=1}^q r_{0ij}^{\alpha \xi} r_{0ij}^{\beta \eta} i_{T_{0j}^\eta} \left( dh_{0ij}^{\gamma \delta} \right) i_{T_{0j}^\xi}\left( C_{j|\mu}^\delta \right) - i_{T_{0i}^\alpha \wedge T_{0i}^\beta }\left( \partial C_{i|\mu}^\gamma \right) + \sum_{\xi, \eta=1}^p \sum_{\delta=1}^q r_{0ij}^{\alpha \xi} r_{0ij}^{\beta \eta} h_{0ij}^{\gamma \delta} X_{j\xi\eta|\mu}^\delta - \sum_{\xi, \eta=1}^p \sum_{\delta,\sigma=1}^q r_{0ij}^{\alpha \xi} r_{0ij}^{\beta \eta} h_{0ij}^{\gamma \sigma}  a_{0j\xi}^{\sigma \delta}  K_{j|\mu}^{\eta \delta}\\
& + \sum_{\eta=1}^p \sum_{\delta, \sigma=1}^q  a_{0i \alpha}^{\gamma \sigma} r_{0ij}^{\beta \eta} h_{0ij}^{\sigma \delta} K_{j|\mu}^{\eta \delta} + \sum_{\xi, \eta=1}^p \sum_{\delta,\sigma=1}^q r_{0ij}^{\alpha \xi} r_{0ij}^{\beta \eta} h_{0ij}^{\gamma \sigma} a_{0j \eta}^{\sigma \delta}  K_{j|\mu}^{\xi \delta} - \sum_{\eta=1}^p \sum_{\delta, \sigma=1}^q   a_{0i\beta}^{\gamma \sigma}  r_{0ij}^{\alpha \xi} h_{0ij}^{\sigma \delta} K_{j|\mu}^{\xi \delta}\\
&\equiv_\mu- \sum_{\xi=1}^p \sum_{\delta, \sigma=1}^q a_{0i \beta}^{\gamma \sigma}   r_{0ij}^{\alpha \xi} h_{0ij}^{\sigma \delta} i_{\Gamma_{j|\mu}^\xi}\left( w_{0i}^\delta \right) + \sum_{\delta=1}^q a_{0i \beta}^{\gamma \delta} i_{\Gamma_{i|\mu}^\alpha}\left( w_{0i}^\delta \right) - \sum_{\eta=1}^p \sum_{\delta, \sigma=1}^q  a_{0i\alpha}^{\gamma \sigma}  r_{0ij}^{\beta \eta} h_{0ij}^{\sigma \delta} i_{\Gamma_{j|\mu}^\eta}\left( w_{0i}^\delta \right)\\
& - \sum_{\delta=1}^q a_{0i\alpha}^{\gamma \delta} i_{\Gamma_{i|\mu}^\beta}\left( w_{0i}^\delta \right)  +  \sum_{\xi, \eta=1}^p r_{0ij}^{\alpha \xi} r_{0ij}^{\beta \eta} h_{0ij}^{\gamma \delta} i_{T_{0j}^\xi \wedge T_{0j}^\eta}\left( \partial C_{j|\mu}^\delta \right) - \sum_{\xi, \eta=1}^p\sum_{\delta, \sigma=1}^q r_{0ij}^{\alpha \xi} r_{0ij}^{\beta \eta} h_{0ij}^{\gamma \sigma} a_{0j\xi}^{\sigma \delta}   i_{T_{0j}^\eta}\left( C_{j|\mu}^\delta \right)\\
&+ \sum_{\eta=1}^p \sum_{\delta, \sigma=1}^q   a_{0i\alpha}^{\gamma \sigma}  r_{0ij}^{\beta \eta} h_{0ij}^{\sigma \delta} i_{T_{0j}^\eta}\left( C_{j|\mu}^\delta \right) +  \sum_{\xi, \eta=1}^p\sum_{\delta, \sigma=1}^q r_{0ij}^{\alpha \xi} r_{0ij}^{\beta \eta} h_{0ij}^{\gamma \sigma}  a_{0j\eta}^{\sigma \delta}  i_{T_{0j}^\xi}\left( C_{j|\mu}^\delta \right) - \sum_{\xi=1}^p \sum_{\delta, \sigma=1}^q  a_{0i\beta}^{\gamma \sigma}   r_{0ij}^{\alpha \xi} h_{0ij}^{\sigma \delta} i_{T_{0j}^\xi}\left( C_{j|\mu}^\delta \right)\\
&- i_{T_{0i}^\alpha \wedge T_{0i}^\beta }\left( \partial C_{i|\mu}^\gamma \right) + \sum_{\xi, \eta=1}^p \sum_{\delta=1}^q r_{0ij}^{\alpha \xi} r_{0ij}^{\beta \eta} h_{0ij}^{\gamma \delta} X_{j\xi\eta|\mu}^\delta - \sum_{\xi, \eta=1}^p \sum_{\delta,\sigma=1}^q r_{0ij}^{\alpha \xi} r_{0ij}^{\beta \eta} h_{0ij}^{\gamma \sigma}  a_{0j\xi}^{\sigma \delta}  K_{j|\mu}^{\eta \delta}\\
& + \sum_{\eta=1}^p \sum_{\delta, \sigma=1}^q  a_{0i \alpha}^{\gamma \sigma}  r_{0ij}^{\beta \eta} h_{0ij}^{\sigma \delta} K_{j|\mu}^{\eta \delta} + \sum_{\xi, \eta=1}^p \sum_{\delta,\sigma=1}^q r_{0ij}^{\alpha \xi} r_{0ij}^{\beta \eta} h_{0ij}^{\gamma \sigma}  a_{0j\eta}^{\sigma \delta}  K_{j|\mu}^{\xi \delta} - \sum_{\eta=1}^p \sum_{\delta, \sigma=1}^q a_{0i \beta}^{\gamma \sigma}  r_{0ij}^{\alpha \xi} h_{0ij}^{\sigma \delta} K_{j|\mu}^{\xi \delta}
\end{align*}}}

Then we have
{\tiny{\begin{align*}
X_{i\alpha\beta |\mu}^\gamma & = \sum_{\xi, \eta=1}^p \sum_{ \sigma=1}^q r_{0ij}^{\alpha \xi} r_{0ij}^{\beta \eta} h_{0ij}^{\gamma \sigma} \left( X_{j\xi\eta |\mu}^\delta + i_{T_{0j}^\xi \wedge T_{0j}^\eta}\left(\partial C_{j|\mu}^\delta \right) +\sum_{\delta=1}^q \left(-   a_{0j\xi}^{\sigma \delta} i_{T_{0j}^\eta}\left( C_{j|\mu}^\delta \right) +  a_{0j\eta}^{\sigma \delta}  i_{T_{0j}^\xi}\left( C_{j|\mu}^\delta \right) - a_{0j\xi}^{\sigma \delta}  K_{j|\mu}^{\eta \delta} +  a_{0j\eta}^{\sigma \delta}  K_{j|\mu}^{\xi \delta} \right) \right)\\
&-  \sum_{\sigma=1}^q a_{0i\beta}^{\gamma \sigma} \left( \cancel{ i_{\Gamma_{i|\mu}^\alpha } \left( w_{0i}^\sigma \right) } + i_{T_{0i}^\alpha}\left( C_{i|\mu}^\sigma \right) + K_{i|\mu}^{\alpha \sigma}  \right)  +  \sum_{\sigma=1}^q a_{0i \alpha}^{\gamma \sigma} \left( \bcancel{ i_{\Gamma_{i|\mu}^\beta}\left( w_{0i}^\sigma \right) } + i_{T_{0i}^\beta}\left( C_{i|\mu}^\sigma \right) + K_{i|\mu}^{\beta \sigma}  \right) \\
& + \sum_{\delta=1}^q a_{0i \beta}^{\gamma \delta}  \cancel{ i_{\Gamma_{i|\mu}^\alpha} \left( w_{0i}^\delta \right) }  - \sum_{\delta=1}^q  a_{0i\alpha}^{\gamma \delta}  \bcancel{ i_{\Gamma_{i|\mu}^\beta}\left( w_{0i}^\delta \right) } - i_{T_{0i}^\alpha \wedge T_{0i}^\beta}\left( \partial C_{i|\mu}^\gamma \right)
\end{align*}}}

Then we have
{\small{\begin{align}\label{sa1}
&X_{i\alpha\beta |\mu}^\gamma + i_{T_{0i}^\alpha \wedge T_{0i}^\beta}\left( \partial C_{i|\mu}^\gamma \right) + \sum_{\sigma=1}^q\left(  - a_{0i\alpha}^{\gamma \sigma}  i_{T_{0i}^\beta}\left( C_{i|\mu}^\sigma \right) + a_{0i\beta}^{\gamma \sigma} i_{T_{0i}^\alpha}\left( C_{i|\mu}^\sigma \right) - a_{0i \alpha}^{\gamma \sigma}  K_{i|\mu}^{\beta \sigma} +   a_{0i \beta}^{\gamma \sigma}  K_{i|\mu}^{\alpha \sigma} \right)\\
& = \sum_{\xi, \eta=1}^p \sum_{ \sigma=1}^q r_{0ij}^{\alpha \xi} r_{0ij}^{\beta \eta} h_{0ij}^{\gamma \sigma} \left( X_{j\xi\eta |\mu}^\delta + i_{T_{0j}^\xi \wedge T_{0j}^\eta}\left(\partial C_{j|\mu}^\delta \right) +\sum_{\delta=1}^q \left(-  a_{0j\xi}^{\sigma \delta} i_{T_{0j}^\eta}\left( C_{j|\mu}^\delta \right) +   a_{0j\eta}^{\sigma \delta}  i_{T_{0j}^\xi}\left( C_{j|\mu}^\delta \right) -  a_{0j\xi}^{\sigma \delta}  K_{j|\mu}^{\eta \delta} +  a_{0j\eta}^{\sigma \delta}  K_{j|\mu}^{\xi \delta} \right) \right) \notag
\end{align}}}

We define $J_{i|\mu}\in \Gamma\left( U_i, \mathcal{A}^{0,0}\left( \mathscr{H}om_{\mathcal{O}_M}\left( \mathcal{N}_{\mathcal{F}_0}^*, \bigwedge^2 \Theta_{\mathcal{F}_0}^* \right) \right) \right)$ by 
{\small{\begin{align*}
&J_{i|\mu}:\Gamma\left( U_i, \mathcal{N}_{\mathcal{F}_0}^* \right) \to \Gamma\left( U_i, \mathcal{A}^{0,0}\left( \bigwedge^2 \Theta_{\mathcal{F}_0} \right) \right)\\
& w_{0i}^\gamma \mapsto \left( T_{0i}^\alpha \wedge T_{0i}^\beta \mapsto X_{i\alpha\beta |\mu}^\gamma + i_{T_{0i}^\alpha \wedge T_{0i}^\beta}\left( \partial C_{i|\mu}^\gamma \right) + \sum_{\sigma=1}^q\left(  - a_{0i\alpha}^{\gamma \sigma}i_{T_{0i}^\beta}\left( C_{i|\mu}^\sigma \right) +  a_{0i\beta}^{\gamma \sigma}  i_{T_{0i}^\alpha}\left( C_{i|\mu}^\sigma \right) -  a_{0i\alpha}^{\gamma \sigma}  K_{i|\mu}^{\beta \sigma} + a_{0i\beta}^{\gamma \sigma} K_{i|\mu}^{\alpha \sigma} \right)   \right)
\end{align*}}}
Then from $(\ref{sa1})$, we have
\begin{align*}
J_\mu:=\left\{ J_{i|\mu}\right\} \in A^{0,0}\left( M, \mathscr{H}om_{\mathcal{O}_M}\left(  \mathcal{N}_{\mathcal{F}_0}^*, \bigwedge^2 \Theta_{\mathcal{F}_0}^* \right) \right)
\end{align*}
We claim that
{\Small{\begin{align}\label{d33}
 \left(\overline{B}_\mu, J_\mu , Z_\mu, \bar{\Phi}_\mu,  \overline{\phi}_\mu   , - \xi_\mu \right) \in &   \frac{A^{0,0}\left( M, \mathscr{H}om_{\mathcal{O}_M}\left( \bigwedge^2 \Theta_{\mathcal{F}_0} ,  \Theta_M \right) \right)}{A^{0,0}\left( M, \mathscr{H}om_{\mathcal{O}_M}\left( \bigwedge^2 \Theta_\mathcal{F} , \Theta_{\mathcal{F}_0}  \right) \right)}\bigoplus A^{0,0}\left(M, \mathscr{H}om_{\mathcal{O}_M}\left( \mathcal{N}_{\mathcal{F}_0}^* , \bigwedge^2 \Theta_{\mathcal{F}_0}^* \right)\right) \bigoplus A^{0,0} \left(M, \mathscr{H}om_{\mathcal{O}_M}\left(\mathcal{N}_{\mathcal{F}_0}^*  , \Theta_{\mathcal{F}_0}^*  \right)  \right) \\
 &\bigoplus \frac{A^{0,1}\left(M, \mathscr{H}om_{\mathcal{O}_M}\left( \Theta_{\mathcal{F}_0} ,  \Theta_M \right)  \right)}{A^{0,1}\left(M, \mathscr{H}om_{\mathcal{O}_M} \left( \Theta_{\mathcal{F}_0} , \Theta_{\mathcal{F}_0} \right)  \right)} \bigoplus \frac{A^{0,1}\left( M, \mathscr{H}om_{\mathcal{O}_M}\left( \mathcal{N}_{\mathcal{F}_0}^*  , \Omega_M^1 \right) \right)}{A^{0,1}\left( M,  \mathscr{H}om_{\mathcal{O}_M} \left( \mathcal{N}_{\mathcal{F}_0}^* , \mathcal{N}_{\mathcal{F}_0}^* \right) \right)} \bigoplus A^{0,2}\left(M, \Theta_M\right) \notag
\end{align}}} 
defines a $2$-cocycle in the following Dolbeault resolution of $\mathcal{F}_0'^\bullet$:

{\Tiny{\begin{center}
$\begin{CD}
\cdots \\
@A \hat{F}_2'AA \\
\frac{A^{0,0}\left(M, \bigwedge^2 \Theta_{\mathcal{F}_0}^*\otimes \Theta_M\right)}{A^{0,0}\left(M,\bigwedge^2 \Theta_{\mathcal{F}_0}^*\otimes \Theta_{\mathcal{F}_0}\right)}\bigoplus A^{0,0}\left(M , \left(\mathcal{N}_{\mathcal{F}_0}^*\right)^*\otimes \bigwedge^2 \Theta_{\mathcal{F}_0}^* \right)\bigoplus A^{0,0}\left(M, \left(\mathcal{N}_{\mathcal{F}_0}^* \right)^*\otimes \Theta_\mathcal{F}^* \right)@>\bar{\partial}>> \cdots \\
@A \hat{F}_1'AA @A \hat{F}_1'AA \\
\frac{A^{0,0}\left(M, \Theta_{\mathcal{F}_0}^*\otimes \Theta_M\right)}{A^{0,0}\left(M, \Theta_{\mathcal{F}_0}^*\otimes \Theta_{\mathcal{F}_0} \right)}\bigoplus \frac{A^{0,0}\left(M,\left(\mathcal{N}_{\mathcal{F}_0}^*\right)^*\otimes \Omega_M^1\right)}{A^{0,0}\left( M, \left(\mathcal{N}_{\mathcal{F}_0}^* \right)^* \otimes \mathcal{N}_{\mathcal{F}_0}^*\right)} @>-\bar{\partial}>>\frac{A^{0,1}\left(M, \Theta_{\mathcal{F}_0}^*\otimes \Theta_M\right)}{A^{0,1}\left(M, \Theta_{\mathcal{F}_0}^*\otimes \Theta_{\mathcal{F}_0} \right)}\bigoplus \frac{A^{0,1}\left(M,\left(\mathcal{N}_{\mathcal{F}_0}^*\right)^*\otimes \Omega_M^1\right)}{A^{0,1}\left( M, \left(\mathcal{N}_{\mathcal{F}_0}^* \right)^* \otimes \mathcal{N}_{\mathcal{F}_0}^*\right)}@>\bar{\partial}>>\cdots \\
@A\hat{F}_0 AA @A \hat{F}_0 AA @A \hat{F}_0 AA \\
A^{0,0}\left(M, \Theta_M \right) @>\bar{\partial}>> A^{0,1} \left(M, \Theta_M\right) @>-\bar{\partial}>> A^{0,2}\left(\Theta_M\right)  
\end{CD}$
\end{center}}}

From the proof of Lemma \ref{te47} and Lemma \ref{ne25} and Lemma \ref{ss12}, it remains to show that $\hat{E}_2'\left( J_{i|\mu} \right)=0$ and $\bar{\partial} J_{i|\mu}+ \hat{E}_1'\left( \bar{\phi}_{i|\mu} \right)=0$ and $\beta_2\left( \bar{B}_{i|\mu}\right) + \gamma_2\left( J_{i|\mu} \right)- \hat{E}_1'\left( Z_{i|\mu} \right)=0$. on $U_i$. We show $\hat{E}_2'\left( J_{i|\mu} \right)=0$.  In fact, since we have from $(\ref{ss7})$
{\Small{\begin{align*}
&X_{i\alpha\beta |\mu}^\gamma + i_{T_{0i}^\alpha \wedge T_{0i}^\beta}\left( \partial C_{i|\mu}^\gamma \right) + \sum_{\sigma=1}^q\left(  - a_{0i\alpha}^{\gamma \sigma}i_{T_{0i}^\beta}\left( C_{i|\mu}^\sigma \right) +  a_{0i\beta}^{\gamma \sigma}  i_{T_{0i}^\alpha}\left( C_{i|\mu}^\sigma \right) -  a_{0i\alpha}^{\gamma \sigma}  K_{i|\mu}^{\beta \sigma} + a_{0i\beta}^{\gamma \sigma} K_{i|\mu}^{\alpha \sigma} \right) \\
&= i_{T_{0i}^\alpha}\partial K_{i|\mu}^{\beta \gamma} - i_{\Pi_{i|\mu}^{\alpha\beta}}\left( w_{0i}^\gamma \right) - \sum_{\xi=1}^p g_{0i\alpha\beta}^\xi K_{i|\mu}^{\xi \gamma} - i_{T_{0i}^\beta}\partial K_{i|\mu}^{\alpha \gamma} + i_{T_{0i}^\alpha \wedge T_{0i}^\beta}\left( \partial C_{i|\mu}^\gamma \right) + \sum_{\sigma=1}^q\left(  - a_{0i\alpha}^{\gamma \sigma}i_{T_{0i}^\beta}\left( C_{i|\mu}^\sigma \right) +  a_{0i\beta}^{\gamma \sigma}  i_{T_{0i}^\alpha}\left( C_{i|\mu}^\sigma \right) -  a_{0i\alpha}^{\gamma \sigma}  K_{i|\mu}^{\beta \sigma} + a_{0i\beta}^{\gamma \sigma} K_{i|\mu}^{\alpha \sigma} \right)
\end{align*}}}
 it is enough to show that $\hat{E}_2'\left( \tilde{\Pi}_{i|\mu} \right)=0$ where $\tilde{\Pi}_{i|\mu}\in \Gamma\left( U_i, \mathscr{H}om_{\mathcal{O}_M}\left( \mathcal{N}_{\mathcal{F}_0}^*, \bigwedge^2 \Theta_{\mathcal{F}_0}^*  \right) \right)$ is defined by
 \begin{align*}
 \tilde{\Pi}_{i|\mu}:\Gamma\left( U_i, \mathcal{N}_{\mathcal{F}_0}^* \right) &\to \Gamma\left(  U_i, \bigwedge^2 \Theta_{\mathcal{F}_0}^* \right) \\
  w_{0i}^\gamma &\mapsto \left( T_{0i}^\alpha \wedge T_{0i}^\beta \mapsto i_{\Pi_{i|\mu}^{\alpha\beta}}\left(w_{0i}^\gamma \right)    \right)
 \end{align*}
 But since $\hat{D}_2\left( \Pi_{i|\mu} \right)=0$, i.e.  $\hat{D}_2\left( \Pi_{i|\mu}\right)\left( T_{0i}^\alpha \wedge T_{0i}^\beta \right) \in \Theta_{\mathcal{F}_0}$ from $(\ref{te7})$, we see that
 \begin{align*}
 \hat{E}_2'\left(\tilde{\Pi}_{i|\mu} \right)\left(T_{0i}^\alpha \wedge T_{0i}^\beta \right)\left(w_{0i}^\gamma \right) = i_{w_{0i}^\gamma} \left( \hat{D}_2\left(\Pi_{i|\mu} \right)\left( T_{0i}^1\wedge T_{0i}^\beta\right)  \right) = 0
 \end{align*}
 
We show that $\bar{\partial} J_{i|\mu} + \hat{E}_1'\left(\bar{\phi}_{i|\mu} \right)=0$. In fact, from $(\ref{ss7})$ and $(\ref{se4})$ and $(\ref{t28})$, we have
{ \small{\begin{align*}
& \bar{\partial} J_{i\alpha\beta|\mu}^\gamma =\bar{\partial} X_{i\alpha\beta|\mu}^\gamma - i_{T_{0i}^\alpha \wedge T_{0i}^\beta}\left(\partial \bar{\partial} C_{i|\mu}^\gamma \right) + \sum_{\sigma=1}^q\left(  a_{0i\alpha}^{\gamma \sigma} i_{T_{0i}^\beta}\left( \bar{\partial} C_{i|\mu}^\sigma \right)    - a_{0i\beta}^{\gamma \sigma} i_{T_{0i}^\alpha}\left( \bar{\partial} C_{i|\mu}^\sigma \right) - a_{0i\alpha}^{\gamma \sigma} \bar{\partial} K_{i|\mu}^{\beta \sigma} + a_{0i\beta}^{\gamma \sigma} \bar{\partial} K_{i|\mu}^{\alpha \sigma} \right)\\
 &= i_{T_{0i}^\alpha} \partial \bar{\partial} K_{i|\mu}^{\beta \gamma} - i_{\bar{\partial} \Pi_{i|\mu}^{\alpha\beta}}\left( w_{0i}^\gamma \right) - \sum_{\xi=1}^p g_{0i\alpha\beta}^\xi \bar{\partial} K_{i|\mu}^{\xi \gamma} - i_{T_{0i}^\beta}\left(\partial \bar{\partial} K_{i|\mu}^{\alpha \gamma} \right) - i_{T_{0i}^\alpha \wedge T_{0i}^\beta}\left(\partial \bar{\partial} C_{i|\mu}^\gamma \right) \\
 & + \sum_{\sigma=1}^q\left(  a_{0i\alpha}^{\gamma \sigma} i_{T_{0i}^\beta}\left( \bar{\partial} C_{i|\mu}^\sigma \right)    - a_{0i\beta}^{\gamma \sigma} i_{T_{0i}^\alpha}\left( \bar{\partial} C_{i|\mu}^\sigma \right) - a_{0i\alpha}^{\gamma \sigma} \bar{\partial} K_{i|\mu}^{\beta \sigma} + a_{0i\beta}^{\gamma \sigma} \bar{\partial} K_{i|\mu}^{\alpha \sigma} \right)\\
 &= - i_{T_{0i}^\alpha} \partial i_{\Phi_{i|\mu}^\beta}\left( w_{0i}^\gamma \right) + i_{T_{0i}^\alpha}\partial i_{T_{0i}^\beta}\left( A_{i|\mu}^\gamma \right) - i_{-\left[ \Phi_{i|\mu}^\alpha, T_{0i}^\beta\right]- \left[ T_{0i}^\alpha, \Phi_{i|\mu}^\beta \right] +\sum_{\xi=1}^p g_{0i\alpha\beta}^\xi \Phi_{i|\mu}^\xi}\left( w_{0i}^\gamma \right) + \sum_{\xi=1}^p g_{0i\alpha\beta}^\xi i_{\Phi_{i|\mu}^\xi}\left(w_{0i}^\gamma \right) - \sum_{\xi=1}^p g_{0i\alpha\beta}^\xi i_{T_{0i}^\xi}\left( A_{i|\mu}^\gamma \right) \\
 &+ i_{T_{0i}^\beta}\partial i_{\Phi_{i\mu}^\alpha}\left( w_{0i}^\gamma \right) - i_{T_{0i}^\beta}\partial i_{T_{0i}^\alpha}\left( A_{i|\mu}^\gamma \right) - i_{T_{0i}^\alpha \wedge T_{0i}^\beta}\left(\partial \bar{\partial} C_{i|\mu}^\gamma \right) +\sum_{\sigma=1}^q a_{0i\alpha}^{\gamma \sigma} i_{T_{0i}^\beta}\left(\bar{\partial} C_{i|\mu}^\sigma \right) - \sum_{\sigma=1}^q a_{0i\beta}^{\gamma \sigma} i_{T_{0i}^\alpha}\left( \bar{\partial} C_{i|\mu}^\sigma \right)\\
 &+ \sum_{\sigma=1}^q a_{0i\alpha}^{\gamma \sigma}i_{\Phi_{i|\mu}^\beta}\left( w_{0i}^\sigma \right) - \sum_{\sigma=1}^q a_{0i}^{\gamma \sigma} i_{T_{0i}^\beta}\left( A_{i|\mu}^\sigma \right) - \sum_{\sigma=1}^q  a_{0i\beta}^{\gamma \sigma} i_{\Phi_{i|\mu}^{\alpha }}\left( w_{0i}^\sigma \right) + \sum_{\sigma=1}^q a_{0i\beta}^{\gamma \sigma}i_{T_{0i}^\alpha}\left( A_{i|\mu}^\sigma \right)\\
 &=-\cancel{\mathcal{L}_{T_{0i}^\alpha} i_{\Phi_{i|\mu}^\beta}\left( w_{0i}^\gamma \right)} +  \xcancel{\mathcal{L}_{T_{0i}^\alpha} i_{T_{0i}^\beta}\left( A_{i|\mu}^\gamma \right) } - \bcancel{ \mathcal{L}_{T_{0i}^\beta} i_{\Phi_{i|\mu}^\alpha}\left( w_{0i}^\gamma \right) }+ \cancel{ i_{\Phi_{i|\mu}^\alpha}\mathcal{L}_{T_{0i}^\beta}\left(  w_{0i}^\gamma \right) }^1+ \cancel{\mathcal{L}_{T_{0i}^\alpha} i_{\Phi_{i|\mu}^\beta}\left( w_{0i}^\gamma \right)} -\bcancel{ i_{\Phi_{i|\mu}^\beta}\mathcal{L}_{T_{0i}^\alpha}\left( w_{0i}^\gamma \right) }^2 - \xcancel{\sum_{\xi=1}^p g_{0i\alpha\beta}^\xi i_{T_{0i}^\xi}\left( A_{i|\mu}^\gamma \right)} \\
 &+ \bcancel{\mathcal{L}_{T_{0i}^\beta} i_{\Phi_{i|\mu}^\alpha}\left( w_{0i}^\gamma \right) } - \xcancel{ i_{T_{0i}^\beta}\mathcal{L}_{T_{0i}^\alpha}\left( A_{i|\mu}^\gamma \right)} + i_{T_{0i}^\alpha \wedge T_{0i}^\beta}\left(\partial A_{i|\mu}^\gamma \right) - i_{T_{0i}^\alpha \wedge T_{0i}^\beta}\left(\partial \bar{\partial} C_{i|\mu}^\gamma \right) +\sum_{\sigma=1}^q a_{0i\alpha}^{\gamma \sigma} i_{T_{0i}^\beta}\left(\bar{\partial} C_{i|\mu}^\sigma \right) - \sum_{\sigma=1}^q a_{0i\beta}^{\gamma \sigma} i_{T_{0i}^\alpha}\left( \bar{\partial} C_{i|\mu}^\sigma \right)\\
  &+ \bcancel{\sum_{\sigma=1}^q a_{0i\alpha}^{\gamma \sigma}i_{\Phi_{i|\mu}^\beta}\left( w_{0i}^\sigma \right)} - \sum_{\sigma=1}^q a_{0i}^{\gamma \sigma} i_{T_{0i}^\beta}\left( A_{i|\mu}^\sigma \right) - \cancel{ \sum_{\sigma=1}^q  a_{0i\beta}^{\gamma \sigma} i_{\Phi_{i|\mu}^{\alpha }}\left( w_{0i}^\sigma \right) }^1 + \sum_{\sigma=1}^q a_{0i\beta}^{\gamma \sigma}i_{T_{0i}^\alpha}\left( A_{i|\mu}^\sigma \right)\\
  &= i_{T_{0i}^\alpha \wedge T_{0i}^\beta}\left( \partial \left(  A_{i|\mu}^\gamma - \bar{\partial} C_{i|\mu}^\gamma \right) \right) - \sum_{\sigma=1}^q a_{0i}^{\gamma \sigma} i_{T_{0i}^\beta}\left(  A_{i|\mu}^\sigma - \bar{\partial} C_{i|\mu}^\sigma \right)  + \sum_{\sigma=1}^q a_{0i \beta}^{\gamma \sigma} i_{T_{0i}^\alpha}\left( A_{i|\mu}^\sigma - \bar{\partial} C_{i|\mu}^\sigma \right)
 \end{align*}}}

 We show that $\beta_2\left( \bar{B}_{i|\mu} \right) + \gamma_2\left( J_{i|\mu} \right)- \hat{E}_1'\left( Z_{i|\mu} \right)=0$. In fact, from $(\ref{ss7})$
 \begin{align*}
 &i_{\Pi_{i|\mu}^{\alpha\beta}}\left( w_{0i}^\gamma \right) + X_{i\alpha\beta|\mu}^\gamma - \sum_{\sigma=1}^q a_{0i\alpha}^{\gamma \sigma} K_{i|\mu}^{\beta \sigma} + \sum_{\sigma=1}^q a_{0i\beta}^{\gamma \sigma} K_{i|\mu}^{\alpha \sigma} - E_1'\left( K_{i|\mu} \right)\left( T_{0i}^\alpha \wedge T_{0i}^\beta \right) \left( w_{0i}^\gamma \right)\\
 &= i_{\Pi_{i|\mu}^{\alpha\beta} } \left(  w_{0i}^\gamma \right) + i_{T_{0i}^\alpha } \partial K_{i|\mu}^{\beta \gamma} - i_{\Pi_{i|\mu}^{\alpha\beta}}\left( w_{0i}^\gamma \right) - \sum_{\xi=1}^p g_{0i\alpha\beta}^\xi K_{i|\mu}^{\xi \gamma} - i_{T_{0i}^\beta}\partial K_{i|\mu}^{\alpha \gamma} - \sum_{\sigma=1}^q a_{0i\alpha}^{\gamma \sigma} K_{i|\mu}^{\beta \sigma} + \sum_{\sigma=1}^q a_{0i\beta}^{\gamma \sigma} K_{i|\mu}^{\alpha \sigma} \\
 &+ i_{T_{0i}^\beta}\partial K_{i|\mu}^{\alpha \gamma} - \sum_{\eta=1}^q a_{i\beta}^{\alpha \eta} K_{i|\mu}^{\eta \gamma} - i_{T_{0i}^\alpha } \partial K_{i|\mu}^{\beta \gamma} + \sum_{\eta=1}^q a_{0i\alpha}^{\gamma \eta} K_{i|\mu}^{\beta \eta} + \sum_{\xi=1}^p g_{0i\alpha\beta}^\xi K_{i|\mu}^{\xi \gamma} = 0
 \end{align*}
 
 This proves $\left(\bar{B}_\mu, J_\mu, Z_\mu, \bar{\Phi}_\mu, \bar{\phi}_\mu, - \xi_\mu   \right)$ defines a $2$-cocycle in the Dolbeualt resolution of $\mathcal{F}_0'^{\bullet}$. Then by hypothesis $\mathbb{H}^2\left( M, \mathcal{F}_0'^\bullet \right)=0$, there exists
 \begin{align*}
 \left( \overline{T}_\mu' , \overline{w_\mu'}, \varphi_\mu' \right)\in \frac{A^{0,0}\left( M, \mathscr{H}om_{\mathcal{O}_M}\left( \Theta_{\mathcal{F}_0}, \Theta_M \right) \right)}{A^{0,0}\left( M, \mathscr{H}om_{\mathcal{O}_M}\left(  \Theta_{\mathcal{F}_0}, \Theta_{\mathcal{F}_0} \right) \right)}\bigoplus \frac{A^{0,0}\left( M, \mathscr{H}om_{\mathcal{O}_M}\left( \mathcal{N}_{\mathcal{F}_0}^*, \Omega_M^1 \right) \right)}{A^{0,0}\left( M, \mathscr{H}om_{\mathcal{O}_M} \left( \mathcal{N}_{\mathcal{F}_0}^*, \mathcal{N}_{\mathcal{F}_0}^*   \right) \right)} \bigoplus A^{0,1}\left( M, \Theta_M \right)
 \end{align*}
 such that $(\ref{te30}), (\ref{te31}), (\ref{te32})$ and $(\ref{n55})$ and $(\ref{s13})$ and additionally
{\small{\begin{align}\label{ss23}
\hat{E}_1'\left(\overline{w_\mu'} \right)= J_\mu=\left\{ w_{0i}^\alpha \mapsto \left(  T_{0i}^\alpha \wedge T_{0i}^\beta \mapsto  X_{i\alpha\beta |\mu}^\gamma + i_{T_{0i}^\alpha \wedge T_{0i}^\beta}\left( \partial C_{i|\mu}^\gamma \right) + \sum_{\sigma=1}^q\left(  - a_{0i\alpha}^{\gamma \sigma}i_{T_{0i}^\beta}\left( C_{i|\mu}^\sigma \right) +  a_{0i\beta}^{\gamma \sigma}  i_{T_{0i}^\alpha}\left( C_{i|\mu}^\sigma \right) -  a_{0i\alpha}^{\gamma \sigma}  K_{i|\mu}^{\beta \sigma} + a_{0i\beta}^{\gamma \sigma} K_{i|\mu}^{\alpha \sigma} \right)   \right)  \right\}
\end{align}}}
 In the proof of Lemma \ref{te47} and Lemma \ref{ne25} and Lemma \ref{ss12}, we have already shown that $\varphi_\mu, T_{i|\mu}, r_{ij|\mu}^{\alpha\beta}$ and $w_{i|\mu}, h_{ij|\mu}^{\alpha\beta}$ satisfying $(\ref{te11})-(\ref{t33})$ and $(\ref{n42})-(\ref{n46}), (\ref{nn121})$ and $(\ref{se5})$. It remains to show $(\ref{ss11})$. In fact, we recall that $w_{i|\mu}^\alpha= C_{i|\mu}^\alpha - w_{i|\mu}'^\alpha + \sum_{\beta=1}^q P_{i|\mu}^{\alpha\beta} w_{0i}^\beta$ from $(\ref{ne19})$. Then we have from $(\ref{ss23})$
\begin{align*}
& - \cancel{ \sum_{\sigma=1}^q a_{0i\alpha}^{\gamma \sigma} i_{T_{0i}^\beta}\left( w_{i|\mu}'^\sigma \right) } +  \bcancel{ \sum_{\sigma=1}^q a_{0i\beta}^{\gamma \sigma} i_{T_{0i}^\alpha}\left( w_{i|\mu}'^\sigma \right)  } + i_{T_{0i}^\alpha \wedge T_{0i}^\beta}\left( \partial w_{i|\mu}'^\gamma  \right)\\
&= X_{i\alpha\beta|\mu}^\gamma + i_{T_{0i}^\alpha \wedge T_{0i}^\beta}\left(\partial C_{i|\mu}^\gamma \right) - \cancel{ \sum_{\sigma=1}^q a_{0i\alpha}^{\gamma \sigma} i_{T_{0i}^\beta}\left( C_{i|\mu}^\sigma \right) } + \bcancel{ \sum_{\sigma=1}^q a_{0i\beta}^{\gamma \sigma} i_{T_{0i}^\alpha}\left( C_{i|\mu}^\sigma \right) } + \sum_{\sigma=1}^q a_{0i\alpha}^{\gamma \sigma} \left( \cancel{ i_{T_{0i}^\beta}\left( w_{i|\mu}^\sigma \right) } + i_{T_{i|\mu}^\beta} \left( w_{0i}^\sigma \right) \right) \\
&- \sum_{\sigma=1}^q a_{0i\beta}^{\gamma \sigma} \left( \bcancel{ i_{T_{0i}^\alpha}\left( w_{i|\mu}^\sigma \right) } + i_{T_{i|\mu}^\alpha}\left( w_{0i}^\sigma \right)   \right)\\
&\iff - X_{i\alpha\beta |\mu}^\gamma = i_{T_{0i}^\alpha \wedge T_{0i}^\beta}\left( \partial\left( C_{i|\mu}^\gamma - w_{i|\mu}'^\gamma \right) \right) + \sum_{\sigma=1}^q a_{0i\alpha}^{\gamma \sigma} i_{T_{i|\mu}^\beta}\left( w_{0i}^\sigma \right) - \sum_{\sigma=1}^q a_{0i\beta}^{\gamma \sigma} i_{T_{i|\mu}^\alpha}\left( w_{0i}^\sigma \right)
\end{align*}
This completes the proof of Lemma \ref{ss24}.

\end{proof}

It remains to determine $\varphi_1,r_{ij}^{\alpha\beta 1}, g_{i\alpha\beta}^{\gamma 1}$ and  $ w_{i}^{\alpha 1}, h_{ij}^{\alpha\beta 1}$ satisfying $(\ref{t13})_1-(\ref{te6})_1$ and $(\ref{n22})_1-(\ref{n212})_1,(\ref{n21})_1$  and $(\ref{se3})_1$, and additionally $(\ref{ss9})_1$. Given $\dim_\mathbb{C} \mathbb{H}^1\left( M, \mathcal{F}_0'^\bullet \right)=r$, we can find a basis of $\mathbb{H}^1\left( M, \mathcal{F}_0'^\bullet \right)$ by using the Dolbeault resolution of $\mathcal{F}_0'^\bullet$ from $(\ref{})$ and represent the basis by
\begin{align*}
\left( \overline{\psi_{\lambda}}, \overline{\kappa_{\lambda}}, \rho_\lambda \right) \in \frac{A^{0,0}\left( M, \mathscr{H}om_{\mathcal{O}_M}\left( \Theta_{\mathcal{F}_0}, \Theta_M \right)\right)}{A^{0,0}\left(M, \mathscr{H}om_{\mathcal{O}_M}\left( \Theta_{\mathcal{F}_0}, \Theta_{\mathcal{F}_0} \right) \right)} \bigoplus \frac{ A^{0,0}\left( M, \mathscr{H}om_{\mathcal{O}_M}\left( \mathcal{N}_{\mathcal{F}_0}^*, \Omega_M^1 \right) \right) }{A^{0,0}\left( M, \mathscr{H}om_{\mathcal{O}_M}\left( \mathcal{N}_{\mathcal{F}_0}^*, \mathcal{N}_{\mathcal{F}_0}^* \right) \right)} \bigoplus A^{0,1}\left( M , \Theta_M \right)
\end{align*}
where $\psi_\lambda$ from $(\ref{st2})$ satisfying $(\ref{st3})-(\ref{tp3})$ and $\kappa_\lambda$ from $(\ref{st7})$ satisfying $(\ref{nnc3})$. We have also the equality $(\ref{st15})$ and additionally we have
\begin{align}\label{st16}
\sum_{\sigma=1}^q a_{0i\beta}^{\gamma \sigma} i_{T_{0i}^\alpha}\left(\kappa_{i\lambda}^\sigma \right) - \sum_{\sigma=1}^q a_{0j\alpha}^{\gamma \sigma} i_{T_{0j}^\beta}\left(\kappa_{i\lambda}^\sigma \right) + i_{T_{0i}^\alpha\wedge T_{0i}^\beta}\left( \partial \kappa_{i\lambda}^\gamma \right) = 0
\end{align}
Then we set $\varphi_1$ as in $(\ref{st3})$, $T_i^{\alpha1}$ as in $(\ref{st4})$, and $r_{ij}^{\alpha\beta 1}$ as in $(\ref{st5})$, and $g_{i\alpha\beta}^{\gamma 1}$ as in $(\ref{st6})$. We also set $w_i^{\alpha 1}$ as in $(\ref{st8})$, and $h_{ij}^{\alpha\beta 1}$ as in $(\ref{st9})$. Then $(\ref{t13})_1-(\ref{te6})_1$ and $(\ref{n23})_1 - (\ref{n213})_1$, $(\ref{n24})_1$ and $(\ref{se3})_1$ are satisfied. It remains to check $(\ref{ss9})_1$. In fact, from $(\ref{st15})$ and $(\ref{st16})$ we have
\begin{align*}
i_{T_{i}^{\alpha 1} \wedge T_i^{\beta 1}}\left( \partial w_i^{\gamma 1} \right) &= i_{\left( T_{0i}^\alpha + \sum_{\lambda=1}^r t_\lambda \left( \psi_{i\lambda}^\alpha + W_{i\lambda}^\alpha  \right) \right) \wedge \left( T_{0i}^\beta + \sum_{\lambda=1}^r t_\lambda \left( \psi_{i\lambda}^\beta + W_{i\lambda}^\beta \right) \right)}\left(  dw_{0i}^\gamma + \sum_{\lambda=1}^r t_\lambda \left(\partial \kappa_{i\lambda}^\gamma + \partial \mathfrak{K}_{i\lambda}^\gamma \right) \right) \\
 &\equiv_1  i_{T_{0i}^\alpha \wedge \left( \sum_{\lambda=1}^r t_\lambda \psi_{i\lambda}^\beta \right)} \left( dw_{0i}^\gamma \right) - i_{T_{0i}^\beta \wedge \left( \sum_{\lambda=1}^r t_\lambda \psi_{i\lambda}^\alpha \right)}\left( dw_{0i}^\gamma \right) + \sum_{\lambda=1}^r t_\lambda i_{T_{0i}^\alpha \wedge T_{0i}^\beta } \left( \partial \kappa_{i\lambda}^\gamma \right) \\
 &=\sum_{\lambda=1}^r t_\lambda  \left( \sum_{\sigma=1}^q a_{0i\alpha}^{\gamma \sigma} i_{\psi_{i\lambda}^\beta}\left( w_{0i}^\sigma     \right)  - \sum_{\sigma=1}^q a_{0i\beta}^{\gamma \sigma} i_{\psi_{i\lambda}^\alpha}\left( w_{0i}^\sigma \right)   + i_{T_{0i}^\alpha \wedge T_{0i}^\beta}\left( \partial \kappa_{i\lambda}^\gamma \right)         \right) = 0
\end{align*}
This completes the inductive construction of $\varphi, T_i^\alpha, r_{ij}^{\alpha\beta}, g_{i\alpha\beta}^\gamma$ and $w_i^\alpha, h_{ij}^{\alpha\beta}$ satisfying $(\ref{t10})-(\ref{te5})$ and $(\ref{n22})-(\ref{n212}),(\ref{n21})$ and $(\ref{se1})$ and $(\ref{ss19})$.

\subsection{Proof of convergence}\

We will prove that $\varphi=\sum_{\mu=1}^\infty \varphi_\mu, T_i^\alpha = \sum_{\mu=0}^\infty T_{i|\mu}^\alpha, r_{ij}^{\alpha\beta}= \sum_{\mu=0}^\infty r_{ij|\mu}^{\alpha\beta}, g_{i\alpha\beta}^\gamma = \sum_{\mu=0}^\infty g_{i\alpha\beta|\mu}^\gamma$, and $w_i^\alpha= \sum_{\mu=0}^\infty w_{i|\mu}^\alpha, h_{ij}^{\alpha\beta}= \sum_{\mu=0}^\infty h_{ij|\mu}^{\alpha\beta}$ converge. Before proceeding the discussion, we introduce another complex of sheaves which also controls foliated deformations of $\left( M, \Theta_{\mathcal{F}_0}, \mathcal{N}_{\mathcal{F}_0}^* \right)$ but removes the quotient in the degree $1$. We shall define the following complex of sheaves $\mathcal{E}_{\mathcal{F}_0}^\bullet$:
\begin{equation}
\begin{CD}
\cdots \\
@AF_3'AA \\
\mathscr{H}om_{\mathcal{O}_M}\left( \bigwedge^3 \Theta_{\mathcal{F}_0}, \frac{\Theta_M}{\Theta_{\mathcal{F}_0}} \right) \bigoplus \mathscr{H}om_{\mathcal{O}_M}\left( \mathcal{N}_{\mathcal{F}_0}^*, \bigwedge^3 \Theta_{\mathcal{F}_0}^* \right) \bigoplus \mathscr{H}om_{\mathcal{O}_M}\left( \mathcal{N}_{\mathcal{F}_0}^*, \bigwedge^2 \Theta_{\mathcal{F}_0}^* \right) \\
@AF_2'AA \\
\mathscr{H}om_{\mathcal{O}_M}\left( \bigwedge^2 \Theta_{\mathcal{F}_0}, \frac{\Theta_M}{\Theta_{\mathcal{F}_0}} \right) \bigoplus \mathscr{H}om_{\mathcal{O}_M}\left( \mathcal{N}_{\mathcal{F}_0}^*, \bigwedge^2  \Theta_{\mathcal{F}_0}^* \right) \bigoplus \mathscr{H}om_{\mathcal{O}_M}\left( \mathcal{N}_{\mathcal{F}_0}^*, \Theta_{\mathcal{F}_0}^* \right) \\
@AF_1''AA \\
\mathscr{H}om_{\mathcal{O}_M}\left( \Theta_{\mathcal{F}_0}, \Theta_M \right) \bigoplus \mathscr{H}om_{\mathcal{O}_M}\left( \mathcal{N}_{\mathcal{F}_0}^*, \Omega_M^1 \right) \\
@AF_0'AA \\
\mathcal{E}_{\mathcal{F}_0}
\end{CD}
\end{equation}
where $F_0'$ is defined as in $(\ref{d32})$ and $F_1''$ is the composition of the natural quotient map with $F_1'$. We will denote the $i$-th cohomology group of $\mathcal{E}_{\mathcal{F}_0}'^\bullet$ by $\mathbb{H}^i\left( M, \mathcal{E}_{\mathcal{F}_0}'^\bullet \right)$. Then we have 
\begin{align*}
\mathbb{H}^i\left( M,  \mathcal{E}_{\mathcal{F}_0}'^\bullet  \right) \cong \mathbb{H}^i\left( M, \mathcal{F}_0'^\bullet \right),\,\,\,\,\,\,\,i\geq 0
\end{align*}
and we have the following Dolbeault resolution of $\mathcal{E}_{\mathcal{F}_0}'^\bullet$:
{\Tiny{\begin{center}
$\begin{CD}
\cdots \\
@A \hat{F}_2'AA \\
\frac{A^{0,0}\left(M, \bigwedge^2 \Theta_{\mathcal{F}_0}^*\otimes \Theta_M\right)}{A^{0,0}\left(M,\bigwedge^2 \Theta_{\mathcal{F}_0}^*\otimes \Theta_{\mathcal{F}_0}\right)}\bigoplus A^{0,0}\left(M , \left(\mathcal{N}_{\mathcal{F}_0}^*\right)^*\otimes \bigwedge^2 \Theta_{\mathcal{F}_0}^* \right)\bigoplus A^{0,0}\left(M, \left(\mathcal{N}_{\mathcal{F}_0}^* \right)^*\otimes \Theta_\mathcal{F}^* \right)@>\bar{\partial}>> \cdots \\
@A \hat{F}_1'AA @A \hat{F}_1'AA \\
A^{0,0}\left(M, \Theta_{\mathcal{F}_0}^*\otimes \Theta_M\right) \bigoplus A^{0,0}\left(M,\left(\mathcal{N}_{\mathcal{F}_0}^*\right)^*\otimes \Omega_M^1\right) @>-\bar{\partial}>> A^{0,1}\left(M, \Theta_{\mathcal{F}_0}^*\otimes \Theta_M\right)\bigoplus A^{0,1}\left(M,\left(\mathcal{N}_{\mathcal{F}_0}^*\right)^*\otimes \Omega_M^1\right) @>\bar{\partial}>>\cdots \\
@A\hat{F}_0' AA @A \hat{F}_0' AA @A \hat{F}_0' AA \\
A^{0,0}\left(M, \mathcal{E}_{\mathcal{F}_0} \right) @>\bar{\partial}>> A^{0,1} \left(M, \mathcal{E}_{\mathcal{F}_0} \right) @>-\bar{\partial}>> \cdots   
\end{CD}$
\end{center}}}

We reinterpret $(\ref{d33})$ in terms of the above bicomplex associated to $\mathcal{E}_{\mathcal{F}_0}^\bullet$. From $(\ref{d28})$ and $(\ref{d29})$, we have
\begin{align*}
\left( -\xi_\mu, \left\{ \left( \bar{\partial} \Lambda_{i|\mu}, \bar{\partial}B_{i|\mu} \right) \right\} \right) \in A^{0,2}\left( M, \mathcal{E}_{\mathcal{F}_0} \right)
\end{align*}

Then we see that from
{\Small{\begin{align} \label{dt5}
& \left(\overline{B}_\mu, J_\mu , Z_\mu, \Phi_\mu,  \phi_\mu   , \left(- \xi_\mu , \left\{ \left(  \bar{\partial} \Lambda_{i|\mu}, \bar{\partial} B_{i|\mu} \right\} \right) \right)\right) \\
& \in    \frac{A^{0,0}\left( M, \mathscr{H}om_{\mathcal{O}_M}\left( \bigwedge^2 \Theta_{\mathcal{F}_0} ,  \Theta_M \right) \right)}{A^{0,0}\left( M, \mathscr{H}om_{\mathcal{O}_M}\left( \bigwedge^2 \Theta_\mathcal{F} , \Theta_{\mathcal{F}_0}  \right) \right)}\bigoplus A^{0,0}\left(M, \mathscr{H}om_{\mathcal{O}_M}\left( \mathcal{N}_{\mathcal{F}_0}^* , \bigwedge^2 \Theta_{\mathcal{F}_0}^* \right)\right) \bigoplus A^{0,0} \left(M, \mathscr{H}om_{\mathcal{O}_M}\left(\mathcal{N}_{\mathcal{F}_0}^*  , \Theta_{\mathcal{F}_0}^*  \right)  \right) \notag \\
 &\,\,\,\,\,\bigoplus A^{0,1}\left(M, \mathscr{H}om_{\mathcal{O}_M}\left( \Theta_{\mathcal{F}_0} ,  \Theta_M \right)  \right)  \bigoplus  A^{0,1}\left( M, \mathscr{H}om_{\mathcal{O}_M}\left( \mathcal{N}_{\mathcal{F}_0}^*  , \Omega_M^1 \right) \right) \bigoplus A^{0,2}\left(M, \mathcal{E}_{\mathcal{F}_0} \right) \notag
\end{align}}}
defines a $2$-cocycle in the above Dolbeault resolution of $\mathcal{E}_{\mathcal{F}_0}'^\bullet$.

We define the \"Holder norms on the sections of $\mathcal{A}^{0,q}\left(\Theta_M \right)$, $\mathcal{A}^{0,p}\left(\mathscr{H}om_{\mathcal{O}_M}\left( \Theta_{\mathcal{F}_0}, \Theta_M \right) \right)$ and $\mathcal{A}^{0,q}\left( \mathscr{H}om_{\mathcal{O}_M}\left( \bigwedge^2 \Theta_{\mathcal{F}_0}, \Theta_M \right) \right)$ and introduce a harmonic theory as in the proof of convergence of Theorem \ref{tt2}. We also define \" Holder norms on the section of $\mathcal{A}^{0,q}\left( \mathscr{H}om_{\mathcal{O}_M}\left( \mathcal{N}_{\mathcal{F}_0}^*, \Omega_M^1 \right) \right)$ and introduce a harmonic theory as in the proof of convergence of Theorem \ref{nn2}. We define the \" Holder norm on the sections of $\mathcal{A}^{0,r}\left( \mathcal{E}_{\mathcal{F}_0} \right)$ and $\mathcal{A}^{0,r}\left( \mathscr{H}om_{\mathcal{O}_M}\left( \mathcal{N}_{\mathcal{F}_0}^*, \Theta_{\mathcal{F}_0}^* \right) \right)$ as in the proof of convergence of Theorem \ref{ss6}.

 We define the H\"older norm $|-|_{k+\alpha}$ (k: an integer $\geq 2$, $0<\alpha< 1$) for sections of $\mathcal{A}^{0,r}\left(\mathscr{H}om_{\mathcal{O}_M}\left( \mathcal{N}_{\mathcal{F}_0}^*, \bigwedge^2 \Theta_{\mathcal{F}_0}^* \right) \right)$ as follows: let $\phi\in \Gamma\left(U_i, \mathcal{A}^{0,r}\left(\mathscr{H}om_{\mathcal{O}_M}\left(  \mathcal{N}_{\mathcal{F}_0}^* ,  \bigwedge^2 \Theta_{\mathcal{F}_0}^* \right) \right) \right)$ and we write 
\begin{align*}
\phi\left(T_i^\beta \wedge T_i^\gamma \right)\left( w_{0i}^\eta \right)=A_i^{\beta \gamma\eta},\,\,\,\,\,\,\,\,\,\,\,A_i^{\beta\gamma\eta}=\frac{1}{r!}\sum A_{i\mu_1\cdots \mu_r}^{\beta\gamma \eta} (z_i)dz_i^{\mu_1}\wedge \cdots \wedge dz_i^{\mu_r} \in \Gamma\left(U_i, \mathcal{A}^{0,r}\right)
\end{align*}
in terms of local coordinates $\left(z_i^1,..., z_i^n \right)$ and let
\begin{align*}
\left| \phi \right|_{k+\alpha}^{U_i} = \sum_{h=0}^k \sup \left| D_i^h A_{i \mu_1\cdots \mu_r }^{\beta\gamma \eta }(z_i) \right| + \sup \frac{ \left| D_i^kA_{ i \mu_1\cdots \mu_r}^{\beta\gamma }(z_i)- D_i^k A_{i\mu_1\cdots \mu_r}^{\beta\gamma \eta }(y_i) \right| }{\left|z_i-y_i\right|^\alpha}
\end{align*}
where the $``\sup$" is extended over all points $z,y\in U_i$, all indices $\beta, \gamma, \eta, \mu_1,...,\mu_r$, and all partial derivatives $D_i^h, D_i^k$ of order $h, k$ with respect to $z_i^1,..., z_i^n, \bar{z}_i^1,..., \bar{z}_i^n$. For $\phi\in A^{0,r}\left( M, \mathscr{H}om_{\mathcal{O}_M}\left( \mathcal{N}_{\mathcal{F}_0}^* ,  \bigwedge^2 \Theta_{\mathcal{F}_0}^*  \right)     \right)$, we define
\begin{align*}
\left| \phi \right|_{k+\alpha} = \max_i \left| \phi \right|_{k+\alpha}^{U_i} 
\end{align*}

We recall the notations $(\ref{d30})$ and $(\ref{d31})$. With this preparation, we will show that for a fixed integer $k\geq 2$ and $0<\alpha<1$, the inductive construction of $\varphi , T_i^\alpha, r_{ij}^{\alpha\beta}, g_{i\alpha\beta}^\gamma$ and $h_{ij}^{\alpha\beta}, w_j^\alpha$ in the previous subsection can be carried out in such a way that satisfy $(\ref{dt1})-(\ref{dt2})$. Then it suffices to prove the $(\ref{sd1})_\mu-(\ref{sd2})_\mu$ for $\mu=1,2,3,\cdots$ and for some proper choice of constants $c>b>0$. We prove $(\ref{sd1})_\mu -(\ref{sd2})_\mu$ by induction on $\mu$. For $\mu=1$ we have $(\ref{st3}),(\ref{st4}),(\ref{st5}),(\ref{st6})$, and $(\ref{ns3}),(\ref{st8}),(\ref{st9})$, and the linear term of $A(t)$ is $\frac{b}{16}\left(t_1+ \cdots  + t_r \right)$. Therefore $(\ref{sd1})_1-(\ref{sd2})_1$ holds if $b$ is sufficiently large.

Now assume that $(\ref{sd1})_{\mu-1}-(\ref{sd2})_{\mu-1}$ are satisfied. We will derive $(\ref{sd1})_\mu - (\ref{sd2})_\mu$. In the following $D_1,D_2,...$ will denote constants which depend only on $k,\alpha, M,\Theta_{\mathcal{F}_0}, \mathcal{N}_{\mathcal{F}_0}^*$. Then as in $(\ref{dt3})$ and $(\ref{dt4})$ and $(\ref{sd22})$, we may assume that
\begin{align}
\left| \left( \tilde{\varphi}_\mu, \left\{ \left( \tilde{b}_{i|\mu}  , \tilde{\lambda}_{i|\mu}\right) \right\}      \right)   \right|_{k+\alpha} \ll D_1\left| \left( - \xi_\mu,  \left\{ \left( \bar{\partial}\Lambda_{i|\mu}, \bar{\partial} B_{i|\mu} \right)   \right\} \right)  \right|_{k-1+\alpha}
\end{align}
where $D_1$ is a constant independent of $\left(- \xi_\mu, \left\{ \left( \bar{\partial} \Lambda_{i|\mu}, \bar{\partial} B_{i|\mu} \right) \right\} \right)$, and we have
\begin{align}
\left( \varphi_\mu'', \left\{\left( b_{i|\mu}'', \lambda_{i|\mu}'' \right) \right\}\right)\in A^{0,1}\left( M, \mathcal{E}_{\mathcal{F}_0} \right),\,\,\,\,\,\,\,\,\mathfrak{d}\varphi_\mu'' = 0,\,\,\,\,\,\,\,\,\bar{\partial}\left( \varphi_\mu'', \left\{ \left( b_{i|\mu}'', \lambda_{i|\mu}''    \right) \right\} \right)=\left( - \xi_\mu, \left\{ \left( \bar{\partial}\Lambda_{i|\mu}, \bar{\partial} B_{i|\mu} \right) \right\}    \right)
\end{align}
\begin{align}
\left| \varphi_\mu'' \right|_{k+\alpha},\,\,\,\,\,\, \left| b_{i|\mu}'' \right|_{k+\alpha}^{U_i}, \,\,\,\,\,\, \left| \lambda_{i|\mu}'' \right|_{k+\alpha}^{U_i} \ll D_2 \left|\left( - \xi_\mu , \left\{ \left( \bar{\partial} \Lambda_{i|\mu} , \bar{\partial} B_{i|\mu}    \right)\right\} \right) \right|_{k-1+\alpha}
\end{align}

From $(\ref{dt5})$, we consider
{\small{\begin{align*}
&0\in A^{0,2}\left( M, \mathcal{E}_{\mathcal{F}_0} \right)\\
&\Phi_\mu'= \left\{ T_{0i}^\alpha \mapsto \Phi_{i|\mu}^\alpha - \bar{\partial} \Gamma_{i|\mu}^\alpha +\sum_{\xi=1}^p \Lambda_{i|\mu}^{\alpha\xi} T_{0i}^\xi +\left[ \varphi_\mu'', T_{0i}^\alpha \right]  - \sum_{\xi=1}^p b_{i|\mu}''^{\alpha  \xi} T_{0i}^\xi   \right\}  \in  A^{0,1}\left(M, \mathscr{H}om_{\mathcal{O}_M}\left(\Theta_{\mathcal{F}_0}, \Theta_M \right) \right) \\
&\phi_\mu'=\left\{ w_{0i}^\alpha \mapsto A_{i|\mu}^\alpha - \bar{\partial} C_{i|\mu}^\alpha + \sum_{\gamma=1}^q B_{i|\mu}^{\alpha \gamma} w_{0i}^\gamma + \mathcal{L}_{\varphi_\mu''}\left( w_{0i}^\alpha \right) - \sum_{\xi=1}^q \lambda_{i|\mu}''^{\alpha \xi} w_{0i}^\xi \right\} \in A^{0,1}\left(M, \mathscr{H}om_{\mathcal{O}_M}\left(  \mathcal{N}_{\mathcal{F}_0}^*, \Omega_M^1 \right) \right) \\
&\overline{B}_\mu-\overline{\left\{ T_{0i}^\alpha \wedge T_{0i}^\beta \mapsto \Pi_{i|\mu}^{\alpha\beta} + \left[ \Gamma_{i|\mu}^\alpha, T_{0i}^\beta \right] - \left[ \Gamma_{i|\mu}^\beta, T_{0i}^\alpha  \right] - \sum_{\gamma=1}^p g_{0i\alpha\beta}^\gamma \Gamma_{i|\mu}^\gamma + \sum_{\gamma=1}^p \Psi_{i|\mu}^{\alpha\beta \gamma} T_{0i}^\gamma    \right\} } \in \frac{ A^{0,0}\left(M, \mathscr{H}om_{\mathcal{O}_M}\left( \bigwedge^2 \Theta_{\mathcal{F}_0}, \Theta_M \right)\right) }{ A^{0,0}\left( M, \mathscr{H}om_{\mathcal{O}_M}\left( \bigwedge^2 \Theta_{\mathcal{F}_0}, \Theta_{\mathcal{F}_0} \right) \right) }\\
&J_\mu=\left\{ w_{0i}^\gamma \mapsto \left( T_{0i}^\alpha \wedge T_{0i}^\beta \mapsto X_{i\alpha\beta |\mu}^\gamma + i_{T_{0i}^\alpha \wedge T_{0i}^\beta}\left( \partial C_{i|\mu}^\gamma \right) + \sum_{\sigma=1}^q\left(  - a_{0i\alpha}^{\gamma \sigma}i_{T_{0i}^\beta}\left( C_{i|\mu}^\sigma \right) +  a_{0i\beta}^{\gamma \sigma}  i_{T_{0i}^\alpha}\left( C_{i|\mu}^\sigma \right) -  a_{0i\alpha}^{\gamma \sigma}  K_{i|\mu}^{\beta \sigma} + a_{0i\beta}^{\gamma \sigma} K_{i|\mu}^{\alpha \sigma} \right)   \right)\right\} \\
&\,\,\,\,\,\,\,\,\,\,\,\,\,\,\,\, \in A^{0,0}\left( M, \mathscr{H}om_{\mathcal{O}_M}\left( \mathcal{N}_{\mathcal{F}_0}^*, \bigwedge^2 \Theta_{\mathcal{F}_0}^* \right) \right) \\
& Z_\mu=\left\{  w_{0i}^\beta \mapsto \left( T_{0i}^\alpha \mapsto K_{i|\mu}^{\alpha\beta} + i_{\Gamma_{i|\mu}^\alpha}\left( w_{0i}^\beta\right) + i_{T_{0i}^\alpha}\left( C_{i|\mu}^\beta \right) \right)\right\} \in A^{0,0}\left( M, \mathscr{H}om_{\mathcal{O}_M}\left( \mathcal{N}_{\mathcal{F}_0}^* , \Theta_{\mathcal{F}_0}^\bullet \right) \right)
\end{align*}}}
Then $\left( 0, \Phi_\mu', \phi_\mu', \overline{B}_\mu, J_\mu, Z_\mu \right)$ defines a $2$-cocycle in the above Dolbeault resolution of $\mathcal{E}_{\mathcal{F}_0}'^\bullet$. Since $\mathbb{H}^2\left( M, \mathcal{F}_0'^\bullet \right)=\mathbb{H}^2\left( M, \mathcal{E}_{\mathcal{F}_0}'^\bullet \right)=0$ by the assumption of Theorem \ref{dt6}, there exist $\left( \chi_\mu, \left\{ \left( \eta_i^{\chi_\mu}, \delta_i^{\chi_\mu} \right) \right\} \right)\in A^{0,1}\left( M, \mathcal{E}_{\mathcal{F}_0} \right), \sigma_\mu\in A^{0,0}\left( M, \mathscr{H}om_{\mathcal{O}_M}\left( \Theta_{\mathcal{F}_0}, \Theta_M \right) \right)$ and $\Sigma_\mu\in A^{0,0}\left( M, \mathscr{H}om_{\mathcal{O}_M}\left( \mathcal{N}_{\mathcal{F}_0}^* , \Omega_M^1 \right) \right)$ such that we have $-\bar{\partial}\left(\xi_\mu, \left\{\left( \eta_i^{\chi_\mu}, \delta_i^{\chi_\mu} \right) \right\} \right)=0, - \bar{\partial} \left(\sigma_\mu, \Sigma_\mu\right)+ \hat{F}_0'\left( \chi_\mu, \left\{\left(\eta_i^{\chi_\mu}, \delta_i^{\chi_\mu} \right) \right\}\right)=\left( \Phi_\mu', \phi_\mu' \right)$ and $\hat{F}_1''\left( \sigma_\mu, \Sigma_\mu\right)=\left( \overline{B}_\mu , J_\mu , Z_\mu \right)$. By using the following Lemma, we will choose appropriate $\left(\chi_\mu, \left\{ \left( \eta_i^{\chi_\mu}, \delta_i^{\chi_\mu}  \right) \right\} \right), \sigma_\mu, \Sigma_\mu$ in a way that $\varphi, T_i^\alpha, r_{ij}^{\alpha\beta}, g_{i\alpha\beta}^\gamma, w_i^\alpha$ and $h_{ij}^{\alpha\beta}$ converge. By combining Lemma \ref{cvt1} and Lemma \ref{nnc2} and Lemma \ref{sd18}, we can prove 

\begin{lemma} \label{dt8}
Suppose that $\left(\varphi, \left\{ \eta_i\right\}, \left\{ \delta_i \right\} \right) \in A^{0,1}\left( M, \mathcal{E}_{\mathcal{F}_0} \right)$ where $\varphi\in A^{0,1}\left( M, \Theta_M \right), \eta_i \in \Gamma\left( U_i, \mathscr{H}om_{\mathcal{O}_M} \left( \Theta_{\mathcal{F}_0} , \Theta_{\mathcal{F}_0}   \right) \right)$ and $\delta_i \in \Gamma\left( U_i, \mathscr{H}om_{\mathcal{O}_M} \left( \mathcal{N}_{\mathcal{F}_0}^*, \mathcal{N}_{\mathcal{F}_0}^* \right) \right)$ with $\eta_i\left( T_{0i}^\alpha\right)=\eta_j\left( T_{0i}^\alpha \right)+\sum_{\beta=1}^p \left[ \varphi, r_{0ij}^{\alpha\beta}\right] T_{0j}^\beta$ for $\alpha=1,...,p$, and $\delta_i\left(w_{0i}^\alpha \right) = \delta_j \left(w_{0i}^\alpha \right)+ \sum_{\beta=1}^q \left[ \varphi, h_{0ij}^{\alpha\beta} \right] w_{0j}^\beta$ for $\alpha=1,...,q$, and $\Phi\in A^{0,0}\left( M, \mathscr{H}om_{\mathcal{O}_M}\left( \Theta_{\mathcal{F}_0}, \Theta_M \right)\right), V\in A^{0,0}\left( M, \mathscr{H}om_{\mathcal{O}_M} \left( \mathcal{N}_{\mathcal{F}_0}^*, \Omega_M^1 \right) \right)$, and $B\in A^{0,0}\left( M, \mathscr{H}om_{\mathcal{O}_M} \left( \bigwedge^2 \Theta_{\mathcal{F}_0}, \Theta_M \right) \right)$, and  $N \in A^{0,0}\left(M, \mathscr{H}om_{\mathcal{O}_M}\left( \mathcal{N}_{\mathcal{F}_0}^*, \bigwedge^2 \Theta_{\mathcal{F}_0}^* \right) \right)$ and $J\in A^{0,0} \left( M, \mathscr{H}om_{\mathcal{O}_M}\left( \mathcal{N}_{\mathcal{F}_0}^*, \Theta_{\mathcal{F}_0}^* \right) \right)$ such that $\left( 0, \left( \Phi, V\right)+ \hat{F}_1'\left( \varphi, \left\{ \left(\eta_i, \delta_i  \right) \right\} \right) \right.$, $\left. \overline{B}, N, J \right)$ defines a $2$-cocycle in the Dolbeault resolution of $\mathcal{E}_{\mathcal{F}_0}'^\bullet$. Then we can find $\left( \chi, \left\{ \eta_i^\chi\right\}, \left\{ \delta_i^\chi \right\} \right) \in A^{0,1}\left( M, \mathcal{E}_{\mathcal{F}_0} \right),\sigma\in A^{0,0}\left( M, \mathscr{H}om_{\mathcal{O}_M} \left( \Theta_{\mathcal{F}_0} , \Theta_M \right) \right)$ and the associated element $\pi_{i\sigma}\in \Gamma\left( U_i, \mathscr{H}om_{\mathcal{O}_M}\left( \bigwedge^2 \Theta_{\mathcal{F}_0}, \Theta_{\mathcal{F}_0} \right) \right)$ and $\Sigma\in A^{0,0}\left( M, \mathscr{H}om_{\mathcal{O}_M}\left( \mathcal{N}_{\mathcal{F}_0}^*, \Omega_M^1 \right) \right)$ satisfying $(\ref{tpc13})-(\ref{tpc14})$ and $(\ref{ns11})-(\ref{dt7})$ and $(\ref{sd7})$, and in addition
\begin{align}
\sum_{\sigma=1}^q a_{0i\beta}^{\gamma \sigma} i_{T_{0i}^\alpha} \left( \Sigma\left( w_{0i}^\sigma \right) \right) - \sum_{\sigma=1}^q a_{0i\alpha}^{\gamma \sigma} i_{T_{0i}^\beta}\left( \Sigma \left( w_{0i}^\sigma \right) \right)  + i_{T_{0i}^\alpha \wedge T_{0i}^\beta}\left(\partial \left( \Sigma \left( w_{0i}^\gamma \right) \right) \right) =  N\left(T_{0i}^\alpha\wedge T_{0i}^\beta \right)\left( w_{0i}^\gamma \right)
\end{align}
\begin{align}
\left| \left(\chi, \left\{ \left( \eta_i^\chi, \delta_i^\chi \right) \right\}  \right) \right|_{k+\alpha} \ll D &\left( \left| \left(\varphi, \left\{ \eta_i \right\}, \left\{ \delta_i \right\} \right) \right|_{k+\alpha} + \left| \Phi \right|_{k-1+\alpha} + \left| V \right|_{k-1+\alpha} + \left| B \right|_{k-1+\alpha}+ \left| N \right|_{k-1+\alpha} + \left| J \right|_{k-1+\alpha}    \right)\\
\left| \sigma\right|_{k+\alpha} \ll D &\left( \left| \left(\varphi, \left\{ \eta_i \right\}, \left\{ \delta_i \right\} \right) \right|_{k+\alpha} + \left| \Phi \right|_{k-1+\alpha} + \left| V \right|_{k-1+\alpha} + \left| B \right|_{k-1+\alpha}+ \left| N \right|_{k-1+\alpha} + \left| J \right|_{k-1+\alpha}    \right)\\
\left| \pi_{i\sigma}\right|_{k+\alpha}^{U_i} \ll D &\left( \left| \left(\varphi, \left\{ \eta_i \right\}, \left\{ \delta_i \right\} \right) \right|_{k-1\alpha} + \left| \Phi \right|_{k-1+\alpha} + \left| V \right|_{k-1+\alpha} + \left| B \right|_{k-1+\alpha}+ \left| N \right|_{k-1+\alpha} + \left| J \right|_{k-1+\alpha}    \right)\\
\left| \Sigma  \right|_{k+\alpha} \ll D &\left( \left| \left(\varphi, \left\{ \eta_i \right\}, \left\{ \delta_i \right\} \right) \right|_{k+\alpha} + \left| \Phi \right|_{k-1+\alpha} + \left| V \right|_{k-1+\alpha} + \left| B \right|_{k-1+\alpha}+ \left| N \right|_{k-1+\alpha} + \left| J \right|_{k-1+\alpha}    \right)
\end{align}
where $D$ is a constant which is independent of $\left(\varphi, \left\{ \eta_i \right\}, \left\{ \delta_i \right\} \right), V, N, J$

\end{lemma}

In Lemma \ref{dt8}, we set
{\small{\begin{align}\label{dt21}
& \left( \varphi, \left\{ \left( \eta_i  , \delta_i \right) \right\} \right):= \left( \varphi_\mu'', \left\{  \left( b_{i|\mu}'' , \lambda_{i|\mu}'' \right) \right\} \right)\in A^{0,1}\left( M, \mathcal{E}_{ \mathcal{F}_0 } \right)\\
& \Phi:= \left\{ T_{0i}^\alpha \mapsto \Phi_{i|\mu}^\alpha - \bar{\partial} \Gamma_{i|\mu}^\alpha +\sum_{\xi=1}^p \Lambda_{i|\mu}^{\alpha\xi} T_{0i}^\xi   \right\}  \in  A^{0,1}\left(M, \mathscr{H}om_{\mathcal{O}_M}\left(\Theta_{\mathcal{F}_0}, \Theta_M \right) \right) \notag \\
&B:= \left\{ T_{0i}^\alpha \wedge T_{0i}^\beta \mapsto \Pi_{i|\mu}^{\alpha\beta} + \left[ \Gamma_{i|\mu}^\alpha, T_{0i}^\beta \right] - \left[ \Gamma_{i|\mu}^\beta, T_{0i}^\alpha  \right] - \sum_{\gamma=1}^p g_{0i\alpha\beta}^\gamma \Gamma_{i|\mu}^\gamma + \sum_{\gamma=1}^p \Psi_{i|\mu}^{\alpha\beta \gamma} T_{0i}^\gamma    \right\}  \in A^{0,0}\left(M, \mathscr{H}om_{\mathcal{O}_M}\left( \bigwedge^2 \Theta_{\mathcal{F}_0}, \Theta_M \right)\right)  \notag\\
V&:= \phi_\mu'=\left\{ w_{0i}^\alpha \mapsto  A_{i|\mu}^\alpha - \bar{\partial} C_{i|\mu}^\alpha + \sum_{\gamma=1}^q B_{i|\mu}^{\alpha \gamma} w_{0i}^\gamma \right\}\notag \\
N&:= J_\mu=\left\{ w_{0i}^\gamma \mapsto \left( T_{0i}^\alpha \wedge T_{0i}^\beta \mapsto X_{i\alpha\beta |\mu}^\gamma + i_{T_{0i}^\alpha \wedge T_{0i}^\beta}\left( \partial C_{i|\mu}^\gamma \right) + \sum_{\sigma=1}^q\left(  - a_{0i\alpha}^{\gamma \sigma}i_{T_{0i}^\beta}\left( C_{i|\mu}^\sigma \right) +  a_{0i\beta}^{\gamma \sigma}  i_{T_{0i}^\alpha}\left( C_{i|\mu}^\sigma \right) -  a_{0i\alpha}^{\gamma \sigma}  K_{i|\mu}^{\beta \sigma} + a_{0i\beta}^{\gamma \sigma} K_{i|\mu}^{\alpha \sigma} \right)   \right) \right\} \notag\\
J&:=Z_\mu=\left\{ w_{0i}^\beta \mapsto \left( T_{0i}^\alpha \mapsto K_{i|\mu}^{\alpha\beta}+ i_{\Gamma_{i|\mu}^\alpha}\left( w_{0i}^\beta \right)  + i_{T_{0i}^\alpha} \left( C_{i|\mu}^\beta \right) \right) \right\} \notag
\end{align}}}
We will estimate $\left(\varphi, \left\{ \left( \eta_i , \delta_i \right) \right\} \right),\Phi,B, V, N, J$. First we will estimate $(\ref{tt3})-(\ref{ds47})$ and $(\ref{nn3})-(\ref{n215})$ and $(\ref{se6})$ and $(\ref{ss20})$. By induction hypothesis $(\ref{sd1})_{\mu-1}-(\ref{sd2})_{\mu-1}$ and from $(\ref{ds22})-(\ref{sd20})$ we may assume that
\begin{align}\label{sd23}
\left| \xi_\mu \right|_{k-1+\alpha} , \left| \Phi_{i|\mu}^\alpha \right|_{k-1+\alpha},\left| \Lambda_{ij|\mu}^{\alpha\beta} \right|_{k+\alpha} ,\left| \Gamma_{ij|\mu}^\alpha \right|_{k+\alpha},  \left| \Pi_{i|\mu}^{\alpha\beta} \right|_{k-1+\alpha}, \left| \lambda_{ijk|\mu}^{\alpha\beta} \right|_{k+1+\alpha}, \left| \Psi_{ij|\mu}^{\alpha\beta \xi} \right|_{k-1+\alpha} \ll D_3 \frac{b}{c} A(t)
\end{align}
and from $(\ref{ns16})-(\ref{ns19})$ and $(\ref{sd29})$ we may assume that
\begin{align*}
\left| A_{i|\mu}^\alpha \right|_{k-1+\alpha},  \left| B_{ij|\mu}^{\alpha\eta}\right|_{k+\alpha},  \left| C_{ij|\mu}^\alpha \right|_{k+\alpha}, \left| D_{ijk|\mu}^{\alpha\beta} \right|_{k+1+\alpha} , \left| K_{i|\mu}^{\alpha \beta} \right|_{k+\alpha} \ll D_4 \frac{b}{c} A(t)
\end{align*}

It remains to estimate
\begin{align*}
X_{i\alpha\beta|\mu}^\gamma &\equiv_\mu \left[ i_{T_i^{\alpha(\mu-1)} \wedge T_i^{\beta(\mu-1)}}\left(\partial w_i^{\gamma(\mu-1)} \right) \right]_\mu = \left[ i_{\left( T_i^{\alpha(\mu-1)} - T_{0i}^\alpha + T_{0i}^\alpha \right) \wedge \left( T_i^{\beta(\mu-1)} - T_{0i}^\beta + T_{0i}^\beta \right)} \left(\partial \left( w_i^{\gamma(\mu-1)} - w_{0i}^\gamma \right) + \partial w_{0i}^\gamma \right) \right]_\mu \\
&\equiv_\mu \left[ i_{\left( T_i^{\alpha(\mu-1)} - T_{0i}^\alpha \right) \wedge \left( T_i^{\beta(\mu-1)} - T_{0i}^\beta \right)} \left(\partial \left( w_i^{\gamma(\mu-1)} - w_{0i}^\gamma \right) \right) \right]_\mu + \left[ i_{T_{0i}^\alpha \wedge \left( T_i^{\beta(\mu-1)} - T_{0i}^\beta \right)}\left( \partial \left( w_i^{\gamma(\mu-1)} - w_{0i}^\gamma \right) \right) \right]_\mu \\
& +\left[ i_{\left( T_i^{\alpha(\mu-1)} - T_{0i}^\alpha \right) \wedge T_{0i}^\beta }\left( \partial \left( w_i^{\gamma(\mu-1)} - w_{0i}^\gamma \right) \right) \right]_\mu + \left[ i_{\left( T_i^{\alpha(\mu-1)} - T_{0i}^\alpha \right) \wedge \left(  T_i^{\beta(\mu-1)} - T_{0i}^\beta \right) } \left(\partial w_{0i}^\gamma \right) \right]_\mu
\end{align*}

Then we have
\begin{align}\label{dt14}
\left| X_{i\alpha\beta|\mu}^\gamma \right|_{k-1+\alpha}\ll D_5 A(t)^3 + D_6 A(t)^2 +  
D_7 A(t)^2 + D_8 A(t)^2 \ll D_9\left( \frac{b}{c} \right)^2 A(t) +  D_{10}\frac{b}{c} A(t) \ll D_{11} \frac{b}{c} A(t)
\end{align}
As in $(\ref{sd24})$ and $(\ref{sd29})$, we may assume that
\begin{align}\label{dt13}
\left| \lambda_{ij|\mu} \right|_{k+1+\alpha}, \left| \Lambda_{i|\mu} \right|_{k+\alpha}, \left| \Gamma_{i|\mu} \right|_{k+\alpha} ,  \left| \Psi_{i|\mu} \right|_{k-1+\alpha}, \left| D_{ij|\mu} \right|_{k+1+\alpha}, \left| B_{i|\mu} \right|_{k+\alpha}, \left| C_{i|\mu} \right|_{k+\alpha} , \left| K_{i|\mu}^{\alpha\beta} \right|_{k+\alpha} \ll D_{12} \frac{b}{c} A(t)
\end{align}

With this preparation, now we estimate $(\ref{dt21})$ and then apply Lemma \ref{dt8}. As in $(\ref{dt10})$ and $(\ref{dt11})$ and $(\ref{dt12})$, we may assume that 
\begin{align}
\left| \left(  \varphi_\mu'', \left\{ \left( b_{i|\mu}'', \lambda_{i|\mu}'' \right) \right\}  \right)   \right|_{k+\alpha}, \left| \Phi \right|_{k-1+\alpha}, \left| B \right|_{k-1+\alpha},  \left| V \right|_{k-1+\alpha} , \left| J \right|_{k-1+\alpha} \ll D_{13} \frac{b}{c} A(t)
\end{align}
On the other hand, from $(\ref{dt14})$ and $(\ref{dt13})$, we have
\begin{align*}
\left| N \right|_{k-1+\alpha} \ll D_{14} \frac{b}{c} A(t)
\end{align*}

 By Lemma \ref{dt8}, we can find $\left( \chi_\mu, \left\{ \left(  \eta_i^{\chi_\mu}, \delta_i^{\chi_\mu} \right) \right\} \right) \in A^{0,1}\left( M, \mathcal{E}_{\mathcal{F}_0} \right)$ and $\sigma_\mu\in A^{0,0}\left( M, \mathscr{H}om_{\mathcal{O}_M}\left( \Theta_{\mathcal{F}_0}, \Theta_M \right) \right)$ and $\pi_{i\sigma_\mu}\in \Gamma\left( U_i, \mathscr{H}om_{\mathcal{O}_M}\left( \bigwedge^2 \Theta_{\mathcal{F}_0}, \Theta_{\mathcal{F}_0} \right) \right)$ and $\Sigma_\mu \in A^{0,0}\left(M, \mathscr{H}om_{\mathcal{O}_M}\left( \mathcal{N}_{\mathcal{F}_0}^*, \Omega_M^1 \right) \right)$ satisfying $(\ref{sd32})-(\ref{ds2})$ and $(\ref{sd30}),(\ref{dt15})$ and $(\ref{dt16})$ and 
\begin{align}
&\sum_{\sigma=1}^q a_{0i\beta}^{\gamma \sigma} i_{T_{0i}^\alpha} \left( \Sigma\left( w_{0i}^\sigma \right) \right) - \sum_{\sigma=1}^q a_{0i\alpha}^{\gamma \sigma} i_{T_{0i}^\beta}\left( \Sigma \left( w_{0i}^\sigma \right) \right)  + i_{T_{0i}^\alpha \wedge T_{0i}^\beta}\left(\partial \left( \Sigma \left( w_{0i}^\gamma \right) \right) \right) \\
&=  X_{i\alpha\beta |\mu}^\gamma + i_{T_{0i}^\alpha \wedge T_{0i}^\beta}\left( \partial C_{i|\mu}^\gamma \right) + \sum_{\sigma=1}^q\left(  - a_{0i\alpha}^{\gamma \sigma}i_{T_{0i}^\beta}\left( C_{i|\mu}^\sigma \right) +  a_{0i\beta}^{\gamma \sigma}  i_{T_{0i}^\alpha}\left( C_{i|\mu}^\sigma \right) -  a_{0i\alpha}^{\gamma \sigma}  K_{i|\mu}^{\beta \sigma} + a_{0i\beta}^{\gamma \sigma} K_{i|\mu}^{\alpha \sigma} \right) \notag
\end{align}
 \begin{align}
\left| \left( \chi_\mu, \left\{ \left( \eta_i^{\chi_\mu}, \delta_i^{\chi_\mu} \right) \right\} \right) \right|_{k+\alpha}, \left| \sigma_\mu \right|_{k+\alpha}, \left| \pi_{i\sigma_\mu} \right|_{k-1+\alpha} , \left| \Sigma_\mu \right|_{k+\alpha} \ll D_{15} \frac{b}{c} A(t)
\end{align}
Then as in the proof of Theorem \ref{tt2} and Theorem \ref{nn2}, if we set $(\ref{ds40})-(\ref{ds43})$ and $(\ref{ns46}),(\ref{ns47})$, then we can choose $b$ and $c$ satisfying $(\ref{sd1})_\mu-(\ref{sd2})_\mu$. Here the missing point is to check $(\ref{ss23})$: In fact,
{\tiny{\begin{align*}
&\sum_{\sigma=1}^q a_{0i\beta}^{\gamma \sigma} i_{T_{0i}^\alpha} \left( w_{i|\mu}'^\sigma \right) - \sum_{\sigma=1}^q a_{0i\alpha}^{\gamma \sigma} i_{T_{0i}^\beta}\left(  w_{i|\mu}'^\sigma \right) + i_{T_{0i}^\alpha \wedge T_{0i}^\beta}\left(\partial \left(  w_{i|\mu}'^\gamma \right) \right) \\
&= \sum_{\sigma=1}^q a_{0i\beta}^{\gamma \sigma} i_{T_{0i}^\alpha} \left( \Sigma_{i|\mu}^\sigma - \mathcal{L}_{\mathfrak{d}G \chi_\mu}\left( w_{0i}^\sigma \right) + U_{i|\mu}\left( w_{0i}^\sigma \right) \right) - \sum_{\sigma=1}^q a_{0i\alpha}^{\gamma \sigma} i_{T_{0i}^\beta}\left(  \Sigma_{i|\mu}^\sigma - \mathcal{L}_{\mathfrak{d}G \chi_\mu}\left( w_{0i}^\sigma \right) + U_{i|\mu}\left( w_{0i}^\sigma \right)\right) + i_{T_{0i}^\alpha \wedge T_{0i}^\beta}\left(\partial \left( \Sigma_{i|\mu}^\gamma - \mathcal{L}_{\mathfrak{d}G \chi_\mu}\left( w_{0i}^\gamma \right) + U_{i|\mu}\left( w_{0i}^\gamma \right) \right) \right)\\
&=  X_{i\alpha\beta |\mu}^\gamma + i_{T_{0i}^\alpha \wedge T_{0i}^\beta}\left( \partial C_{i|\mu}^\gamma \right) + \sum_{\sigma=1}^q\left(  - a_{0i\alpha}^{\gamma \sigma}i_{T_{0i}^\beta}\left( C_{i|\mu}^\sigma \right) +  a_{0i\beta}^{\gamma \sigma}  i_{T_{0i}^\alpha}\left( C_{i|\mu}^\sigma \right) -  a_{0i\alpha}^{\gamma \sigma}  K_{i|\mu}^{\beta \sigma} + a_{0i\beta}^{\gamma \sigma} K_{i|\mu}^{\alpha \sigma} \right)  
\end{align*}}}
Then $\varphi(t), T_i^\alpha, r_{ij}^{\alpha\beta}, g_{i\alpha\beta}^\gamma$ and $w_i^\alpha, h_{ij}^{\alpha\beta}$ converge. Then by the same argument in the proof of Theorem \ref{tt2} and Theorem \ref{nn2}, we can construct a foliated analytic family $\pi:\left( \mathcal{M}, \Theta_{\mathcal{F}}, \mathcal{N}_\mathcal{F}^* \right)\to \Delta_\epsilon$ such that the foliated Kodaira-Spencer map $\varphi_0: T_0\left( \Delta_\epsilon \right)\to \mathbb{H}^1\left( M, \mathcal{F}^\bullet \right)$ is bijective. This completes the proof of Theorem \ref{dt6}.

\end{proof}

\section{Theorem of completeness of deformations of foliated complex analytic structures in terms of both tangent and cotangent sheaves 2}

We shall prove Theorem \ref{sc1} again by using the complex of sheaves $\mathcal{F}_0'^{\bullet}$ instead of $\mathcal{F}_0^\bullet$. We recall Definition \ref{r10} of completeness of a foliated complex analytic family in terms of both tangent and cotangent sheaves.
\begin{theorem}\label{cs1}
Let $\left(\mathcal{M}, \Theta_\mathcal{F}, \mathcal{N}_\mathcal{F}^* ,  B, \pi \right)$ be a foliated analytic family of deformations of a complex foliated complex manifold $(M , \Theta_{\mathcal{F}_0}, \mathcal{N}_{\mathcal{F}_0}^* )=\omega^{-1}(0)$ such that both $\Theta_{\mathcal{F}_0}$ and $\mathcal{N}_{\mathcal{F}_0}^*$ are locally free, and $B$ is a domain of $\mathbb{C}^r$ containing $0$. If the foliated Kodaira-Spence map $\varphi_0: T_0(B)\to \mathbb{H}^1(M,\mathcal{F}_0'^\bullet)$ is surjective, the foliated analytic family $\left(\mathcal{M}, \Theta_\mathcal{F}, \mathcal{N}_\mathcal{F}^* ,  B, \pi \right)$ in terms of both tangent and cotangent sheaves is complete at $0 \in B$.
\end{theorem}

\begin{proof}

Let $\left( \mathcal{M}, \Theta_\mathcal{F}, \mathcal{N}_\mathcal{F}^*, B, \pi \right)$ be a foliated analytic family in terms of both tangent and cotangent sheaves which is represented as in section \ref{cs9}. We keep the notations in the proof of Theorem \ref{tc1} and Theorem \ref{nc1} and Theorem \ref{sc1}. Then we consider the integrability condition $i_{T_{0i}^\alpha \wedge T_{0i}^\beta}\left( dw_{0i}^\gamma \right)=0$ instead of $w_{0i}^1 \wedge \cdots \wedge w_{0i}^q \wedge dw_{0i}^\alpha=0$ on $\left(M, \Theta_{\mathcal{F}_0},\mathcal{N}_{\mathcal{F}_0}^*\right)$.

Let $\left(\mathcal{M}', \Theta_{\mathcal{F}'}, \mathcal{N}_{\mathcal{F}'}^*  , D, \pi' \right)$ be an another foliated analytic family such that $\pi'^{-1}(0')=\left(M, \Theta_{\mathcal{F}_0}, \mathcal{N}_{\mathcal{F}_0}^*\right)$. We may assume the following: $(\ref{sc2})-(\ref{sc3})$ in the proof of Theorem \ref{tc1}, and $(\ref{sc4})-(\ref{cs2}), (\ref{cs3})-(\ref{sc5})$ in the proof of Theorem \ref{nc1} and $(\ref{cs4})$, and additionally
\begin{align}
i_{T_i'^\alpha\wedge T_i'^\beta}\left( dw_i'^\gamma \right) = 0 \,\,\,\,\,\,\,\,\,\,\textnormal{for}\,\,\, \alpha,\beta=1,...,p\,\,\,\gamma=1,...,q.
\end{align}

In order to prove Theorem \ref{sc1}, it suffices to construct holomorphic functions
\begin{align*}
\varphi_i&:\mathcal{U}_i'\to \mathbb{C}^n \\
s &:D \to \mathbb{C}^{m}
\end{align*}
and a matrix function $\left(b_i^{\alpha\beta}(\xi_i, u)\right)_{\alpha,\beta=1,...,p}$ with $b_i^{\alpha\beta}:\mathcal{U}_i'\to \mathbb{C}$ and a matrix function $\left( c_i^{\alpha\beta}(\xi_i, u) \right)_{\alpha,\beta=1,...,q}$ with $c_i^{\alpha\beta}:\mathcal{U}_i' \to \mathbb{C}$ such that $(\ref{tc211})-(\ref{tc4})$ and $(\ref{nc211})-(\ref{nc4})$.

First we prove the existence of formal solution of $(\ref{tc2})-(\ref{tc4})$ and $(\ref{nc2})-(\ref{nc4})$. We recall Notation \ref{te27}. Then $(\ref{tc2})-(\ref{tc4})$ are equivalent to the systems of congruences $(\ref{tc5})_\mu-(\ref{tc7})_\mu$, respectively, and $(\ref{nc2})-(\ref{nc4})$ are equivalent to the systems of congruences $(\ref{nc5})_\mu-(\ref{nc7})_\mu$, respectively for $\mu=1,2,3, \cdots$. We shall construct $\varphi_i^\mu, s^\mu, b_i^{\alpha\beta \mu}$ and $c_i^{\alpha\beta \mu}$ satisfying $(\ref{tc5})_\mu-(\ref{tc7})_\mu$, and $(\ref{nc5})_\mu- (\ref{nc7})_\mu$ by induction on $\mu$. We assume that $\varphi_i^{\mu-1}, s^{\mu-1}, b_i^{\alpha\beta(\mu-1)}$ and $c_i^{\alpha\beta (\mu-1)}$ are already determined. Then we define homogenous polynomials $\Gamma_{ij|\mu}^\alpha, B_{ij|\mu}^{\alpha\beta}, G_{i|\mu}^{\alpha\beta \gamma}$, and $\Pi_{i|\mu}^{\alpha \gamma}$ of degree $\mu$ by the congruences $(\ref{tc35})-(\ref{tc38})$, respectively, and we define homogenous polynomials $A_{ij|\mu}^{\alpha\beta}$ and $W_{i|\mu}^{\alpha \gamma}$ of degree $\mu$ by the congruences $(\ref{nc9}), (\ref{nc10})$, respectively. We set
\begin{align*}
\Gamma_{ij|\mu}:= \sum_{\alpha=1}^n \Gamma_{ij|\mu}^\alpha \frac{\partial}{\partial z_i^\alpha},\,\,\,\,\,\,\,\,\,\, \Pi_{i|\mu}^\alpha:=\sum_{\gamma=1}^n \Pi_{i|\mu}^{\alpha \gamma} \frac{\partial}{\partial z_i^\gamma},\,\,\,\,\,\,\,\,\,\,\, W_{i|\mu}^\alpha := \sum_{\gamma=1}^n W_{i|\mu}^{\alpha \gamma} dz_i^\gamma
\end{align*}

\begin{lemma}\label{cs5}
We have the following equalities$:$ $(\ref{tc33})-(\ref{tc34})$ and $(\ref{nc13})$ and $(\ref{sc6})$ and additionally
\begin{align}\label{cs6}
  \sum_{\delta=1}^q  a_{0i\alpha}^{\gamma \delta} i_{T_{0i}^\beta}\left(  W_{i|\mu}^{\delta }  \right) - \sum_{\delta=1}^q  a_{0i\beta}^{\gamma \delta} i_{T_{0i}^\alpha} \left( W_{i|\mu}^\delta \right)  - i_{T_{0i}^\alpha\wedge T_{0i}^\beta}\left( d W_{i|\mu}^\gamma \right) = 0 
\end{align}
\end{lemma}

\begin{proof}
$(\ref{tc32})-(\ref{tc34})$ follows from Lemma \ref{tc39}, and $(\ref{nc13})$ follows from Lemma \ref{nc11}, and $(\ref{sc6})$ follows from Lemma \ref{sc8}. We prove $(\ref{sc6})$. In fact, first we note that from $(\ref{ss3})$
\begin{align}
i_{T_i^\alpha}\left( dw_i^\gamma \right) &= \sum_{\delta=1}^q a_{i\alpha}^{\gamma \delta}(z_i,t) w_i^\delta \Longrightarrow i_{T_i^\alpha}\left( \sum_{\delta, \eta=1}^n \frac{\partial w_i^{\gamma \eta}}{\partial z_i^\delta} dz_i^\delta \wedge dz_i^\eta \right) = \sum_{\delta, \eta=1}^n \left(T_i^{\alpha\delta}\frac{\partial w_i^{\gamma \eta}}{\partial z_i^\delta}- T_i^{\alpha\delta}\frac{\partial w_i^{\gamma \delta}}{\partial z_i^\eta} \right)dz_i^\eta \notag \\
&\Longrightarrow \sum_{\delta=1}^n \left(T_i^{\alpha\delta}\frac{\partial w_i^{\gamma \eta}}{\partial z_i^\delta}- T_i^{\alpha\delta}\frac{\partial w_i^{\gamma \delta}}{\partial z_i^\eta} \right) = \sum_{\delta=1}^q a_{i\alpha}^{\gamma \delta} w_i^{\delta \eta} \label{cs11}
\end{align}
Since $i_{T_i^\alpha \wedge T_i^\beta}\left(  dw_i^\gamma \right)=0$, we see that
\begin{align}\label{cs10}
0&=i_{T_i^\beta}\left( i_{T_i^\alpha}\left(dw_i^\gamma \right)   \right)= \sum_{\delta=1}^q \sum_{\eta=1}^n a_{i\alpha}^{\gamma\delta} (z_i,t) T_i^{\beta \eta} (z_i,t) w_i^{\delta \eta}(z_i,t)
\end{align}
We set $a_{0i\alpha}^{\gamma \delta}:= a_{i\alpha}^{\gamma \delta}(z_i,0)$. Then we have from $(\ref{cs10})$ and $(\ref{cs11})$ and $(\ref{tc38})$ and $(\ref{nc10})$,
{\small{\begin{align*}
&0=\sum_{\delta=1}^q \sum_{\eta=1}^n a_{i\alpha}^{\gamma\delta} \left(\varphi_i^{\mu-1}, s^{\mu-1}\right) T_i^{\beta \eta} \left(\varphi_i^{\mu-1}, s^{\mu-1}\right)w_i^{\delta \eta}\left(\varphi_i^{\mu-1}, s^{\mu-1}\right) \\
&\equiv_\mu \sum_{\delta=1}^q \sum_{\eta=1}^n a_{i\alpha}^{\gamma\delta} \left(\varphi_i^{\mu-1}, s^{\mu-1}\right) \left( \Pi_{i|\mu}^{\beta \eta} + \sum_{c=1}^p \sum_{d=1}^n b_i^{\beta c(\mu-1)} T_i'^{c d}\frac{\partial \varphi_i^\eta}{\partial \xi_i^d} \right) w_i^{\delta \eta}\left(\varphi_i^{\mu-1}, s^{\mu-1}\right) \\
&\equiv_\mu \sum_{\delta=1}^q \sum_{\eta=1}^n a_{0i\alpha}^{\gamma \delta} \Pi_{i|\mu}^{\beta \eta} w_{0i}^{\delta \eta} + \sum_{c=1}^p\sum_{\delta=1}^q \sum_{\eta, d=1}^n b_i^{\beta c(\mu-1)} T_i'^{cd}\frac{\partial \varphi_i^{\eta(\mu-1)}}{\partial \xi_i^d}  T_i^{\alpha \delta}\left(\varphi_i^{\mu-1}, s^{\mu-1} \right) \left(\frac{\partial w_i^{\gamma \eta}}{\partial z_i^\delta}\left(\varphi_i^{\mu-1}, s^{\mu-1} \right) - \frac{\partial w_i^{\gamma \delta}}{\partial z_i^\eta}\left( \varphi_i^{\mu-1}, s^{\mu-1} \right)\right)\\
&\equiv_\mu  \sum_{\delta=1}^q \sum_{\eta=1}^n a_{0i\alpha}^{\gamma \delta} \Pi_{i|\mu}^{\beta \eta} w_{0i}^{\delta \eta} +  \sum_{c=1}^p\sum_{\delta=1}^q \sum_{\eta, d=1}^n b_i^{\beta c(\mu-1)} T_i'^{cd}\frac{\partial \varphi_i^{\eta(\mu-1)}}{\partial \xi_i^d} \left( \Pi_{i|\mu}^{\alpha \delta} + \sum_{g=1}^p \sum_{h=1}^n b_i^{\alpha g (\mu-1)} T_i'^{ g h} \frac{\partial \varphi_i^{\delta(\mu-1)}}{\partial \xi_i^h} \right) \frac{\partial w_i^{\gamma \eta}}{\partial z_i^\delta}\left(  \varphi_i^{\mu-1} , s^{\mu-1} \right) \\
& -  \sum_{c=1}^p\sum_{\delta=1}^q \sum_{\eta, d=1}^n b_i^{\beta c(\mu-1)} T_i'^{cd}\frac{\partial \varphi_i^{\eta (\mu-1)} }{\partial \xi_i^d} \left( \Pi_{i|\mu}^{\alpha \delta} + \sum_{g=1}^p \sum_{h=1}^n b_i^{\alpha g (\mu-1)} T_i'^{ g h} \frac{\partial \varphi_i^{\delta(\mu-1)}}{\partial \xi_i^h} \right) \frac{\partial w_i^{\gamma \delta}}{\partial z_i^\eta}\left(  \varphi_i^{\mu-1} , s^{\mu-1} \right) \\
&\equiv_\mu \sum_{\delta=1}^q \sum_{\eta=1}^n a_{0i\alpha}^{\gamma \delta} \Pi_{i|\mu}^{\beta \eta} w_{0i}^{\delta \eta} - \sum_{\delta=1}^q \sum_{\eta=1}^n a_{0i\beta}^{\gamma \delta} \Pi_{i|\mu}^{\alpha \eta} w_{0i}^{\delta \eta} + \sum_{c, g=1}^p \sum_{\eta, d, h=1}^n b_i^{\beta  c (\mu-1)} T_i'^{cd} \frac{\partial \varphi_i^{\eta (\mu-1)} }{\partial \xi_i^d} b_i^{\alpha g (\mu-1)} T_i'^{gh}\frac{\partial}{\partial \xi_i^h}\left( w_i^{\gamma \eta}\left( \varphi_i^{\mu-1}, s^{\mu-1} \right) \right)\\
& - \sum_{c, g=1}^p \sum_{\delta, d, h=1}^n b_i^{\beta  c (\mu-1)} T_i'^{cd}  b_i^{\alpha g (\mu-1)} T_i'^{gh} \frac{\partial \varphi_i^{\delta(\mu-1)}}{\partial \xi_i^h}\frac{\partial}{\partial \xi_i^d}\left( w_i^{\gamma \delta}\left( \varphi_i^{\mu-1}, s^{\mu-1} \right) \right)\\
&\equiv_\mu \sum_{\delta=1}^q \sum_{\eta=1}^n a_{0i\alpha}^{\gamma \delta} \Pi_{i|\mu}^{\beta \eta} w_{0i}^{\delta \eta} - \sum_{\delta=1}^q \sum_{\eta=1}^n a_{0i\beta}^{\gamma \delta} \Pi_{i|\mu}^{\alpha \eta} w_{0i}^{\delta \eta} + \sum_{c,g=1}^p \sum_{ d, h=1}^n b_i^{\beta c(\mu-1)} T_i'^{cd} b_i^{\alpha g(\mu-1)} T_i'^{g h} \frac{\partial}{\partial \xi_i^h }\left( \sum_{\eta=1}^q c_i^{\gamma\eta (\mu-1)} w_i'^{\eta d} - W_{i|\mu}^{\gamma d} \right)\\
&- \sum_{c, g=1}^p \sum_{d, h=1}^n b_i^{\beta c (\mu-1)} T_i'^{cd} b_i^{\alpha g (\mu-1)} T_i'^{gh}\frac{\partial}{\partial \xi_i^d}\left( \sum_{\eta=1}^q c_i^{\gamma \eta(\mu-1)} w_i'^{\eta h} - W_{i|\mu}^{\gamma h} \right) \left(\textnormal{we note that $i_{T_i'^c \wedge T_i'^g}\left( dw_i'^\eta  \right)=0$ and $\sum_{d=1}^n T_i'^{cd} w_i'^{\eta d} =0$ } \right)\\
& =   \sum_{\delta=1}^q \sum_{\eta=1}^n a_{0i\alpha}^{\gamma \delta} T_{0i}^{\beta \eta} W_{i|\mu}^{\delta \eta} - \sum_{\delta=1}^q \sum_{\eta=1}^n a_{0i\beta}^{\gamma \delta} T_{0i}^{\alpha \eta} W_{i|\mu}^{\delta \eta}  - \sum_{d,h=1}^n T_{0i}^{\beta d} T_{0i}^{\alpha h}\left(  \frac{\partial W_{i|\mu}^{\gamma d}}{\partial \xi_i^h} - \frac{\partial W_{i|\mu}^{\gamma h}}{\partial \xi_i^d}\right) \\
& =  \sum_{\delta=1}^q  a_{0i\alpha}^{\gamma \delta} i_{T_{0i}^\beta}\left(  W_{i|\mu}^{\delta }  \right) - \sum_{\delta=1}^q  a_{0i\beta}^{\gamma \delta} i_{T_{0i}^\alpha} \left( W_{i|\mu}^\delta \right)  - i_{T_{0i}^\alpha\wedge T_{0i}^\beta}\left( d W_{i|\mu}^\gamma \right)
\end{align*}}}
This completes the proof of Lemma \ref{cs5}.
\end{proof}

Our purpose is to determine $\varphi^\mu= \varphi^{\mu-1}+ \varphi_{i|\mu}, s^\mu= s^{\mu-1}+ s_\mu, b_i^{\alpha\beta(\mu-1)}+ b_{i|\mu}^{\alpha\beta}$, and $c_i^{\alpha\beta(\mu-1)}+ c_{i|\mu}^{\alpha\beta}$ satisfying $(\ref{tc5})_\mu-(\ref{tc7})_\mu$, and $(\ref{nc5})_\mu-(\ref{nc7})_\mu$.

\begin{lemma}\label{cs8}
$(\ref{tc5})_\mu-(\ref{tc7})_\mu$ are equivalent to $(\ref{tc8})-(\ref{tc10})$, respectively, and $(\ref{nc5})_\mu-(\ref{nc7})_\mu$ are equivalent to $(\ref{nc17})-(\ref{nc19})$, respectively.
\end{lemma}

\begin{proof}
Lemma \ref{cs8} follows from Lemma \ref{sc7}.
\end{proof}

We define an element $\overline{\Pi}_{i|\mu}\in \Gamma\left( U_i, \mathscr{H}om_{\mathcal{O}_M} \left( \Theta_{\mathcal{F}_0}, \frac{\Theta_M}{\Theta_{\mathcal{F}_0}}   \right) \right)$ as in $(\ref{tc213})$ and $\overline{W}_{i|\mu}\in \Gamma\left( U_i , \mathscr{H}om_{\mathcal{O}_M} \left( \mathcal{N}_{\mathcal{F}_0}^*, \frac{\Omega_M^1}{\mathcal{N}_{\mathcal{F}_0}^*}    \right) \right)$ as in $(\ref{nc125})$. Then $(\ref{tc32})-(\ref{tc34})$ and $(\ref{nc13})$ and $(\ref{sc6})$ and $(\ref{cs6})$ imply that
\begin{align*}
\left( \left\{ \overline{\Pi}_{i|\mu} \right\}, \left\{ - \overline{W}_{i|\mu} \right\} , \left\{ \Gamma_{ij|\mu} \right\}  \right) \in C^0\left( \mathcal{U}, \mathscr{H}om_{\mathcal{O}_M}\left( \Theta_{\mathcal{F}_0} , \frac{\Theta_M}{\Theta_{\mathcal{F}_0}}   \right) \right) \bigoplus C^0\left( \mathcal{U}, \mathscr{H}om_{\mathcal{O}_M} \left(  \mathcal{N}_{\mathcal{F}_0}^*, \frac{\Omega_M^1}{\mathcal{N}_{\mathcal{F}_0}^*}  \right) \right) \bigoplus C^1\left( \mathcal{U}, \Theta_M \right)
\end{align*}
defines a $1$-cocycle in the following \v Cech resolution of $\mathcal{F}_0'^\bullet$:
{\Tiny{\begin{center}
$\begin{CD}
\cdots \\
@AAA \\
C^0\left( \mathcal{U}, \bigwedge^2 \Theta_{\mathcal{F}_0}^*\otimes \frac{\Theta_M}{\Theta_{\mathcal{F}_0}} \right) \bigoplus C^0\left( \mathcal{U}, \left(\mathcal{N}_{\mathcal{F}_0}^* \right)^* \otimes \bigwedge^2 \Theta_{\mathcal{F}_0}^* \right) \bigoplus C^0\left( \mathcal{U}, \left( \mathcal{N}_{\mathcal{F}_0}^* \right)^*\otimes \Theta_{\mathcal{F}_0}^* \right) @>-\delta >> \cdots \\
@AAA @AAA  \\
C^0\left(\mathcal{U}, \Theta_{\mathcal{F}_0}^* \otimes \frac{\Theta_M}{\Theta_{\mathcal{F}_0} } \right) \bigoplus C^0\left( \mathcal{U}, \left(\mathcal{N}_{\mathcal{F}_0}^* \right)^* \otimes \frac{\Omega_M^1}{\mathcal{N}_{\mathcal{F}_0}^*} \right) @>\delta>>  C^1\left(\mathcal{U}, \Theta_{\mathcal{F}_0}^* \otimes \frac{\Theta_M}{\Theta_{\mathcal{F}_0} }  \right) \bigoplus C^1\left( \mathcal{U}, \left( \mathcal{N}_{\mathcal{F}_0}^* \right)^* \otimes  \frac{\Omega_M^1}{\mathcal{N}_{\mathcal{F}_0}^*} \right) @>-\delta >> \cdots \\
@AAA @AAA @AAA \\
C^0(\mathcal{U}, \Theta_M) @>-\delta>> C^1(\mathcal{U}, \Theta_M) @>\delta>> C^2(\mathcal{U}, \Theta_M) 
\end{CD}$
\end{center}}}
By the hypothesis that the foliated Kodaira-Spencer map $\varphi_0: T_0(B) \to \mathbb{H}^1\left( M, \mathcal{F}_0^\bullet \right)$ is surjective, we can find homogeneous polynomial $s_\mu^\lambda$ such that
\begin{align*}
\varphi_0\left( \sum_{\lambda=1}^r s_\mu^\lambda \frac{\partial}{\partial t_\lambda}   \right) = \left(  \left\{  \overline{\Pi}_{i|\mu}  \right\} , \left\{ - \overline{W}_{i|\mu} \right\}   \left\{ \Gamma_{ij|\mu}  \right\}   \right)
\end{align*}

Since we have
\begin{align*}
\varphi_0\left( \frac{\partial}{\partial t_\lambda} \right) =\left( \{\alpha_{i\lambda}\} =\left\{ T_{0i}^\alpha \mapsto \overline{ -\frac{\partial T_i^\alpha (z_i,t)}{\partial t_\lambda}|_{t=0} }\right\}, \left\{ \beta_{i\lambda}\right\}= \left\{ w_{0i}^\alpha \mapsto \overline{- \frac{\partial w_i^\alpha(z_i,t)}{\partial t_\lambda} |_{t=0} }\right\} ,  \{\rho_{ij\lambda} \} =  \left\{ \sum_{\alpha=1}^n \frac{\partial f_{ij}^\alpha}{\partial t_\lambda}|_{t=0} \frac{\partial}{\partial z_i^\alpha } \right\}   \right),
\end{align*}
there exists $\{\varphi_{i|\mu}\}\in C^0\left(\mathcal{U}, \Theta_M \right)$ and $\left\{ b_{i|\mu}\right\} \in C^0\left( \mathcal{U}, \mathscr{H}om_{\mathcal{O}_M}\left( \Theta_{\mathcal{F}_0}, \Theta_{\mathcal{F}_0} \right) \right)$ and $\left\{ c_{i|\mu} \right\} \in C^0\left( \mathcal{U}, \mathscr{H}om_{\mathcal{O}_M} \left( \mathcal{N}_{\mathcal{F}_0}^*, \mathcal{N}_{\mathcal{F}_0}^* \right) \right)$ satisfying $(\ref{tc8})-(\ref{tc10})$ and $(\ref{nc17})-(\ref{nc19})$ as in the proof of Theorem \ref{tc1} and Theorem \ref{nc1}. This completes the inductive construction of $\varphi_i^\mu, s^\mu, b_i^{\alpha\beta \mu}$ and $c_i^{\alpha\beta \mu}$ and we can prove the convergence as in the same way with subsection \ref{scc5}. This completes the proof of Theorem \ref{cs1}.

\end{proof}

\appendix

\section{Isomorphism between the leaf complex and the dual leaf complex for regular foliations, and definition of the dual leaf complex of (singular) holomorphic foliations defined by locally free subsheaves of cotangent sheaves}\label{appendixA}

\subsection{The leaf complex of a compact foliated complex manifold $(M, \Theta_{\mathcal{F}_0})$} \label{AppendixA1}\

Let $(M, \Theta_{\mathcal{F}_0})$ be a compact foliated complex manifold. Then we have the leaf complex (see \cite{GM88})
\begin{align*}
\Theta_{\mathcal{F}_0}^\bullet: \Theta_M \xrightarrow{D_0} \mathscr{H}om_{\mathcal{O}_M}\left( \Theta_{\mathcal{F}_0}, \frac{\Theta_M}{\Theta_{\mathcal{F}_0}} \right) \xrightarrow{D_1} \mathscr{H}om_{\mathcal{O}_M}\left( \bigwedge^2 \Theta_{\mathcal{F}_0}, \frac{\Theta_M}{\Theta_{\mathcal{F}_0}} \right) \xrightarrow{D_2} \mathscr{H}om_{\mathcal{O}_M}\left( \bigwedge^3 \Theta_{\mathcal{F}_0}, \frac{\Theta_M}{\Theta_{\mathcal{F}_0}} \right)  \xrightarrow{D_3} \cdots
\end{align*}
where $D_0$ is defined for $X \in \Theta_M$, and $Y\in \Theta_{\mathcal{F}_0}$ by
\begin{align*}
(D_0 X)(Y)=\pi([X,Y]) =\overline{[X, Y]}
\end{align*}
where $\pi$ is the projection $\Theta_M\to \Theta_M/\Theta_{\mathcal{F}_0}$, and for $p>0$, $D_p$ is defined in the following way: for $X_1,...,X_{p+1}\in \Theta_{\mathcal{F}_0}$ and $w\in \mathscr{H}om_{\mathcal{O}_M}\left( \bigwedge^p \Theta_{\mathcal{F}_0}, \frac{\Theta_M}{\Theta_{\mathcal{F}_0}} \right)$,
\begin{align*} 
(D_pw)\left(X_1,...,X_{p+1}\right)=\sum_{i=1}^{p+1} (-1)^{i+1} \overline{\left[X_i, w\left(X_1,...,\hat{X}_i,...,X_{p+1}\right)\right]}+\sum_{i<j}(-1)^{i+j}w\left(\left[X_i,X_j \right],...,\hat{X}_i,...,\hat{X}_j,...,X_{p+1}\right)
\end{align*}
where the first bracket is the map $\overline{[-.-]}:\Theta_\mathcal{F}\times \Theta_M/\Theta_\mathcal{F}\to \Theta_M/\Theta_\mathcal{F}$ induced by the Lie bracket.

\subsection{The dual leaf complex for a compact complex manifold with a regular foliation}\label{AppendixA2}\

Let $\left(M, \mathcal{N}_{\mathcal{F}_0}^* \right)$ be a compact foliated complex manifold such that $\mathcal{N}_{\mathcal{F}_0}^*\subset \Omega_M^1$ defines a codimension $q$ regular foliation on $M$. We will define a complex of sheaves called the dual leaf complex 
\begin{align*}
\mathcal{N}_{\mathcal{F}_0}^{*\bullet}: \Theta_M \xrightarrow{E_0}  \mathscr{H}om_{\mathcal{O}_M}\left( \mathcal{N}_{\mathcal{F}_0}^*, \frac{\Omega_M^1}{\mathcal{N}_{\mathcal{F}_0}^*} \right) \xrightarrow{E_1} \mathscr{H}om_{\mathcal{O}_M}\left( \mathcal{N}_{\mathcal{F}_0}^*, \bigwedge^2 \frac{\Omega_M^1}{\mathcal{N}_{\mathcal{F}_0}^* } \right) \xrightarrow{E_2} \mathscr{H}om_{\mathcal{O}_M}\left( \mathcal{N}_{\mathcal{F}_0}^*, \bigwedge^3 \frac{\Omega_M^1}{\mathcal{N}_{\mathcal{F}_0}^*} \right) \xrightarrow{E_3} \cdots
\end{align*}

Let $\mathcal{U}=\{U_i\}$ be an open covering of $M$ by coordinate neighborhoods such that $\mathcal{N}_{\mathcal{F}_0}^*$ on $U_i$ is generated by $w_i^1, ..., w_i^q\in \Gamma\left(U_i, \Omega_M^1 \right)$ with the relation $w_i^\alpha= \sum_{\beta=1}^q h_{ij}^{\alpha\beta} w_j^\beta$ for $h_{ij}^{\alpha\beta}\in \Gamma(U_{ij}, \mathcal{O}_M)$. Since $\mathcal{N}_{\mathcal{F}_0}^*$ is a regular foliation, the integrability condition $w_i^1\wedge \cdots \wedge w_i^q \wedge dw_i^\alpha=0,\alpha=1,...,q$ implies that we can write $dw_i^\alpha= \sum_{\beta=1}^q a_i^{\alpha\beta} \wedge w_i^\beta, \alpha=1,...,q$ for some $a_i^{\alpha\beta}\in \Gamma\left( U_i, \Omega_M^1 \right)$. Then for $X\in \Theta_M$, $E_0(X)$ is locally defined on $U_i$ by
\begin{align*}
E_0(X)\left(w_i^\alpha \right) = \overline{\mathcal{L}_X\left(w_i^\alpha\right)},\,\,\,\,\,\,\,\,\,\alpha=1,..., q, \,\,\,\,\,\,\,\,\,\textnormal{$\mathcal{L}=$ the Lie derivative and $\overline{\mathcal{L}_X\left(w_i^\alpha\right)}$ is the image of $\mathcal{L}_X \left(w_i^\alpha\right)$ in $\frac{\Omega_M^1}{\mathcal{N}_{\mathcal{F}_0}^*}$}
\end{align*}
and linearly extends to $\Gamma(U_i, \mathcal{N}_{\mathcal{F}_0}^*)$. This is well-defined since 
\begin{align*}
E_0(X)(w_i^\alpha)- E_0(X)\left(\sum_{\beta=1}^q h_{ij}^{\alpha\beta} w_j^\beta\right)= \overline{\mathcal{L}_X(w_i^\alpha)}- \sum_{\beta=1}^q h_{ij}^{\alpha\beta} E_0(X)\left(w_j^\beta\right)= \overline{\mathcal{L}_X(w_i^\alpha)} - \sum_{\beta=1}^q h_{ij}^{\alpha\beta} \overline{\mathcal{L}_X\left( w_j^\beta\right)}=0
\end{align*}
and for $p>0$, $E_p$ is defined in the following way: for $\phi\in \mathscr{H}om_{\mathcal{O}_M}\left( \mathcal{N}_{\mathcal{F}_0}^*, \bigwedge^p \frac{\Omega_M^1}{\mathcal{N}_{\mathcal{F}_0}^*} \right)$ and $E_p(\phi)$ is locally defined on $U_i$ by
\begin{align*}
E_p(\phi)\left(w_i^\alpha \right)= \overline{d\widetilde{\phi \left(w_i^\alpha \right)}-\sum_{\beta=1}^q a_i^{\alpha\beta}\wedge \widetilde{\phi\left(w_i^\beta\right)}}
\end{align*}
and linearly extends to $\Gamma(U_i, \mathcal{N}_{\mathcal{F}_0}^*)$. Here $\widetilde{\phi(w_i^\alpha)}$ is a lifting of $\phi(w_i^\alpha)$ under $\bigwedge^p\Omega_M^1\to \bigwedge^p \frac{\Omega_M^1}{\mathcal{N}_{\mathcal{F}_0}^*} $\footnote{Strictly speaking $\bigwedge^p \Omega_M^1\to \bigwedge^p\frac{\Omega_M^1}{\mathcal{N}_{\mathcal{F}_0}^*}$ is not surjective on $U_i$. More precisely for $x\in U_i$, choose a neighborhood $U_x\subset U_i$ such that $\phi(w_i^\alpha)$ on $U_x$ has a lifting $\widetilde{\phi_x(w_i^\alpha)}\in \Gamma\left(U_x, \bigwedge^p \Omega_M^1\right)$. We compute $\overline{d\widetilde{\phi_x(w_i^\alpha)}- \sum_{\beta=1}^q a_i^{\alpha\beta}\wedge \widetilde{\phi_x(w_i^\beta)}}$ locally on $U_x$  for each $x\in U$, and since they are equal on $U_x\cap U_y\ne \emptyset$ by the above arguments, we can glue together to define $E_p(\phi)(w_i^\alpha)$ on $U_i$.}, and $\overline{A}$ is the image of $A\in \bigwedge^{p+1} \Omega_M^1$ in $\bigwedge^{p+1} \frac{\Omega_M^1}{\mathcal{N}_{\mathcal{F}_0}^*}$. We show that $E_p(\phi)(w_i^\alpha)$ is independent of choices of a lifting of $\phi(w_i^\alpha), \alpha=1,...,p $. In fact, for $A, B^1,..., B^q\in \mathcal{N}_{\mathcal{F}_0}^* \bigwedge \left( \bigwedge^{p-1} \Omega_M^1 \right)$, we have $dA- \sum_{\beta=1}^q a_i^{\alpha\beta} \wedge B^\beta\in \mathcal{N}_{\mathcal{F}_0}^* \bigwedge \left( \bigwedge^p \Omega_M^1 \right)$. We show that $E_p(\phi)(w_i^\alpha)$ is independent of choice of $a_i^{\alpha\beta}$. In fact, we simply note that if $dw_i^\alpha= \sum_{\beta=1}^q a_i^{\alpha\beta} \wedge w_i^\beta= a_i'^{\alpha\beta} \wedge w_i^\beta$, then $\sum_{\beta=1}^q \left(a_i^{\alpha\beta} - a_i'^{\alpha\beta}\right)\wedge w_i^\beta=0$, i.e. $a_i^{\alpha\beta}- a_i'^{\alpha\beta}\in \mathcal{N}_{\mathcal{F}_0}^*$. Next we show that $E_p(\phi)$ is well-defined: in other words, $E_p(\phi)$ on $U_i$ and $E_p(\phi)$ on $U_j$ defines the same homomorphism on $U_i\cap U_j$. We recall that $w_i^\alpha=\sum_{\beta=1}^q h_{ij}^{\alpha\beta} w_j^\beta$. Then we have
\begin{align*}
\sum_{\beta,\gamma=1}^q h_{ij}^{\beta\gamma}a_i^{\alpha\beta} \wedge w_j^\gamma= \sum_{\beta=1}^q a_i^{\alpha\beta} \wedge w_i^\beta=dw_i^\alpha=\sum_{\beta=1}^q dh_{ij}^{\alpha\beta} \wedge  w_j^\beta +\sum_{\beta=1}^q h_{ij}^{\alpha\beta} dw_j^\beta=\sum_{\gamma=1}^q d h_{ij}^{\alpha\gamma} \wedge w_j^\gamma +\sum_{\beta,\gamma=1}^q h_{ij}^{\alpha\beta} a_j^{\beta\gamma}\wedge w_j^\gamma
\end{align*}
Then we see that
\begin{align*}
\sum_{\gamma=1}^q \left(   \sum_{\beta=1}^q h_{ij}^{\beta\gamma} a_i^{\alpha\beta}- dh_{ij}^{\alpha\gamma}-\sum_{\beta=1}^q h_{ij}^{\alpha\beta} a_j^{\beta\gamma}   \right) \wedge w_j^\gamma=0
\end{align*}
so that 
\begin{align*}
\sum_{\beta=1}^q h_{ij}^{\beta\gamma} a_i^{\alpha\beta}- dh_{ij}^{\alpha\gamma}-\sum_{\beta=1}^q h_{ij}^{\alpha\beta} a_j^{\beta\gamma} \in \mathcal{N}_{\mathcal{F}_0}^*
\end{align*}
By abuse of notations, we write $\phi(w_i^\alpha)$ instead of $\widetilde{\phi(w_i^\alpha)}$. Then
\begin{align*}
E_p(\phi)\left(\sum_{\beta=1}^q h_{ij}^{\alpha\beta} w_j^\beta\right)= \sum_{\beta=1}^q h_{ij}^{\alpha\beta} E_p(\phi)\left(w_j^\beta\right)=\sum_{\beta=1}^q h_{ij}^{\alpha\beta} \left( d\phi\left(w_j^\beta\right)-\sum_{\gamma=1}^q a_j^{\beta\gamma} \wedge \phi\left(w_j^\gamma\right) \right)
\end{align*}
On the other hand,
\begin{align*}
E_p(\phi)(w_i^\alpha)&= d\phi(w_i^\alpha)-\sum_{\beta=1}^q a_i^{\alpha\beta} \wedge \phi(w_i^\beta)=d\phi\left(\sum_{\beta=1}^q h_{ij}^{\alpha\beta} w_j^\beta \right)-\sum_{\beta=1}^q a_i^{\alpha\beta}\wedge \phi\left(\sum_{\gamma=1}^q h_{ij}^{\beta\gamma} w_j^\gamma\right)\\
&=\sum_{\gamma=1}^q d h_{ij}^{\alpha\gamma} \wedge \phi(w_j^\gamma) +\sum_{\beta=1}^q h_{ij}^{\alpha\beta} d\phi(w_j^\beta) -\sum_{\beta,\gamma=1}^q h_{ij}^{\beta\gamma} a_i^{\alpha\beta} \wedge \phi(w_j^\gamma)
\end{align*}

Hence we have
\begin{align*}
 E_p(\phi)\left(\sum_{\beta=1}^q h_{ij}^{\alpha\beta} w_j^\beta\right)-E_p(\phi)(w_i^\alpha)\in \mathcal{N}_{\mathcal{F}_0}^* \bigwedge \left( \bigwedge^p \Omega_M^1 \right)
\end{align*}
so that $E_p(\phi)$ is well-defined.

Next we show that $E_{p}\circ E_{p+1}=0$. First we show that $E_1\circ E_0=0$ locally on $U_i$. In fact, for $X\in \Theta_{\mathcal{F}_0}$,
\begin{align*}
E_1(E_0(X))\left(w_i^\alpha\right)&= dE_0(X)\left(w_i^\alpha\right)-\sum_{\beta=1}^q a_i^{\alpha\beta}\wedge E_0(X)\left(w_i^\beta\right)=d\mathcal{L}_X\left(w_i^\alpha\right)-\sum_{\beta=1}^q a_i^{\alpha\beta} \wedge \mathcal{L}_X \left(w_i^\beta\right)\\
            &=\mathcal{L}_X \left(\sum_{\beta=1}^q a_i^{\alpha\beta} \wedge w_i^\beta \right)-\sum_{\beta=1}^q a_i^{\alpha\beta} \wedge \mathcal{L}_X\left(w_i^\beta\right)=\sum_{\beta=1}^q \mathcal{L}_X\left(a_i^{\alpha\beta} \right) \wedge w_i^\beta\in \mathcal{N}_{\mathcal{F}_0}^* \wedge \Omega_M^1
\end{align*}
so that $E_1\circ E_0=0$. For $p>0$, we show that $E_{p+1}\circ E_p=0$ locally on $U_i$. We note that from $dw_i^\alpha=\sum_{\beta=1}^q a_i^{\alpha\beta} \wedge w_i^\beta$, we have $0=d(dw_i^\alpha)=\sum_{\beta=1}^q da_i^{\alpha\beta}\wedge w_i^\alpha -a_i^{\alpha\beta}\wedge dw_i^\beta=\sum_{\beta=1}^q da_i^{\alpha\beta} \wedge w_i^\alpha-\sum_{\beta,\gamma=1}^q a_i^{\alpha\beta} \wedge a_i^{\beta\gamma}\wedge w_i^{\gamma}$, so that $\sum_{\beta=1}^q\left( da_i^{\alpha\beta}  -\sum_{\gamma=1}^q a_i^{\alpha\gamma}\wedge a_i^{\gamma\beta} \right) \wedge w_i^\beta=0$. Hence, for each $\alpha=1,...,p$, $da_i^{\alpha\beta}-\sum_{\gamma=1}^q a_i^{\alpha\gamma}\wedge a_i^{\gamma\beta}\in \mathcal{N}_{\mathcal{F}_0}^* \wedge \Omega_M^1$. Then for $\phi\in \mathscr{H}om_{\mathcal{O}_M}\left( \mathcal{N}_{\mathcal{F}_0}^*, \bigwedge^p \frac{\Omega_M^1}{\mathcal{N}_{\mathcal{F}_0}^*} \right)$, we have
\begin{align*}
E_{p+1}(E_p(\phi))\left(w_i^\alpha\right)&=d\left(d\phi\left(w_i^\alpha\right)-\sum_{\beta=1}^q a_i^{\alpha\beta}\wedge \phi\left(w_i^\beta\right)   \right)-\sum_{\gamma=1}^q a_i^{\alpha\gamma}\wedge\left(d\phi\left(w_i^\gamma\right) -  \sum_{\beta=1}^q  a_i^{\gamma\beta}\wedge \phi \left(w_i^\beta \right) \right)\\
&=-\sum_{\beta=1}^q da_i^{\alpha\beta}\wedge \phi\left(w_i^\beta\right) +\sum_{\beta=1}^q a_i^{\alpha\beta}\wedge d\phi\left(w_i^\beta\right) -\sum_{\beta=1}^q a_i^{\alpha\beta}\wedge d\phi\left(w_i^\beta\right)+\sum_{\beta,\gamma=1}^q a_i^{\alpha\gamma} \wedge a_i^{\gamma\beta} \wedge \phi\left(w_i^\beta\right)\\
&=-\sum_{\beta=1}^q\left( da_i^{\alpha\beta}-\sum_{\gamma=1}^q a_i^{\alpha\gamma}\wedge a_i^{\gamma\beta} \right) \wedge \phi\left(w_i^\beta\right) \in \mathcal{N}_{\mathcal{F}_0}^*\bigwedge \left( \bigwedge^{p+1} \Omega_M^1\right)
\end{align*}
This implies that $E_{p+1}\circ E_p=0$. Hence $\mathcal{N}_{\mathcal{F}_0}^{*\bullet}$ defines a complex of sheaves.

\begin{proposition}\label{appendixA3}
Let $\mathcal{F}_0$ be a regular foliation on a compact complex manifold $M$, so that $\mathcal{N}_{\mathcal{F}_0}^* \cong \left(\Theta_M/ \Theta_{\mathcal{F}_0} \right)^*$ and $\Theta_{\mathcal{F}_0}\cong \left( \Omega_M^1/\mathcal{N}_{\mathcal{F}_0}^* \right)^*$. Then we have an isomorphism between the leaf complex $\Theta_{\mathcal{F}_0}^\bullet$ and the dual leaf complex $\mathcal{N}_{\mathcal{F}_0}^{*\bullet}$$:$
{\small{\begin{center}
$\begin{CD}
\Theta_{\mathcal{F}_0}^\bullet:@. \Theta_M @>D_0>> \mathscr{H}om_{\mathcal{O}_M}\left( \Theta_{\mathcal{F}_0}, \frac{\Theta_M}{\Theta_{\mathcal{F}_0}} \right) @>D_1>> \mathscr{H}om_{\mathcal{O}_M}\left( \bigwedge^2 \Theta_{\mathcal{F}_0}, \frac{\Theta_M}{\Theta_{\mathcal{F}_0}}\right) @>D_2>> \mathscr{H}om_{\mathcal{O}_M}\left( \bigwedge^3 \Theta_{\mathcal{F}_0}, \frac{\Theta_M}{\Theta_{\mathcal{F}_0} }\right) @>D_3>> \cdots\\
 @.@V-\textnormal{Id}VV @V\alpha_1 VV @V \alpha_2 VV @V\alpha_3 VV \\
 \mathcal{N}_{\mathcal{F}_0}^{*\bullet} :@.\Theta_M @>E_0>> \mathscr{H}om_{\mathcal{O}_M}\left(\mathcal{N}_{\mathcal{F}_0}^*,  \frac{\Omega_M^1}{ \mathcal{N}_{\mathcal{F}_0}^* }  \right) @>{E}_1>> \mathscr{H}om_{\mathcal{O}_M}\left(\mathcal{N}_{\mathcal{F}_0}^*, \bigwedge^2 \frac{\Omega_M^1}{\mathcal{N}_{\mathcal{F}_0}^*}  \right) @> {E}_2>> \mathscr{H}om_{\mathcal{O}_M}\left( \mathcal{N}_{\mathcal{F}_0}^*, \bigwedge^3 \frac{\Omega_M^1}{\mathcal{N}_{\mathcal{F}_0}^*} \right) @> {E}_3>> \cdots 
\end{CD}$
\end{center}}}
where $\alpha_p:\mathscr{H}om_{\mathcal{O}_M}\left( \bigwedge^p \Theta_{\mathcal{F}_0} , \frac{\Theta_M}{\Theta_{\mathcal{F}_0}} \right) \to \mathscr{H}om_{\mathcal{O}_M}\left( \mathcal{N}_{\mathcal{F}_0}^*, \bigwedge^p \frac{\Omega_M^1}{\mathcal{N}_{\mathcal{F}_0}^*} \right) $ is the dual map.
\end{proposition}

\begin{proof}
First we show that the following diagram commutes:
\begin{center}
$\begin{CD}
\Theta_M @>D_0>> \mathscr{H}om_{\mathcal{O}_M}\left(\Theta_\mathcal{F}, \frac{\Theta_M}{\Theta_{\mathcal{F}_0}}\right)\\
@V\textnormal{Id} VV @VV-\alpha_1 V\\
\Theta_M @>E_0>> \mathscr{H}om_{\mathcal{O}_M}\left(\mathcal{N}_\mathcal{F}^*, \frac{\Omega_M^1}{ \mathcal{N}_{\mathcal{F}_0}^*} \right)
\end{CD}$
\end{center}

We show the commutativity locally on $U_i$. We keep the notations of subsection \ref{AppendixA1} and subsection \ref{AppendixA2}. For $X\in \Gamma(U_i, \Theta_M)$,
\begin{align*}
\alpha_1\circ D_0(X)(w_i^\alpha)\in \Gamma\left(U_i, \frac{\Omega_M^1}{\mathcal{N}_{\mathcal{F}_0}^*} \right)=\Gamma\left(U_i, \mathscr{H}om_{\mathcal{O}_M}\left(\Theta_{\mathcal{F}_0} , \mathcal{O}_M \right) \right)
\end{align*}
Then for $T\in \Gamma(U_i, \Theta_{\mathcal{F}_0})$,
\begin{align*}
\left(\alpha_1\circ D_0(X)(w_i^\alpha)\right)(T)= i_{[X, T]}(w_i^\alpha)
\end{align*}

On the other hand,
\begin{align*}
E_0(X)(w_0^\alpha)= \overline{\mathcal{L}_X(w_i^\alpha)} \in \Gamma\left( U_i, \frac{\Omega_M^1}{\mathcal{N}_{\mathcal{F}_0}^*}\right) = \Gamma\left( U_i, \mathscr{H}om_{\mathcal{O}_M}\left( \Theta_{\mathcal{F}_0}, \mathcal{O}_M \right)\right)
\end{align*}
Then for $T\in \Gamma\left( U_i, \Theta_{\mathcal{F}_0}\right)$, we have
\begin{align*}
\left(E_0(X)(w_i^\alpha)\right)(T)= i_T\left( \mathcal{L}_X(w_i^\alpha) \right) =\mathcal{L}_X(i_T(w_i^\alpha))-i_{[X,T]}(w_i^\alpha)=-i_{[X,T]}(w_i^\alpha)
\end{align*}
This proves the commutativity of the above diagram.

Next we show that the following diagram commutes:
\begin{center}
$\begin{CD}
\mathscr{H}om_{\mathcal{O}_M}\left(\bigwedge^p\Theta_{\mathcal{F}_0} ,  \frac{\Theta_M}{\Theta_{\mathcal{F}_0} } \right) @>D_p >> \mathscr{H}om_{\mathcal{O}_M} \left( \bigwedge^{p+1} \Theta_\mathcal{F}, \frac{\Theta_M}{\Theta_{\mathcal{F}_0} }\right)\\
@V\alpha_p VV @VV\alpha_{p+1} V\\
\mathscr{H}om_{\mathcal{O}_M}\left(\mathcal{N}_\mathcal{F}^*, \bigwedge^p \frac{ \Omega_M^1}{\mathcal{N}_{\mathcal{F}_0}^*} \right) @>E_p >> \mathscr{H}om_{\mathcal{O}_M}\left(\mathcal{N}_\mathcal{F}^*, \bigwedge^{p+1}  \frac{\Omega_M^1}{ \mathcal{N}_{\mathcal{F}_0}^* } \right)
\end{CD}$
\end{center}
We show the commutativity locally on $U_i$. For $\phi\in \Gamma\left(U_i, \mathscr{H}om_{\mathcal{O}_M}(\bigwedge^p\Theta_{\mathcal{F}_0},  \frac{\Theta_M }{\Theta_{\mathcal{F}_0}})\right)$,
\begin{align*}
\alpha_{p+1}\circ D_p(\phi)(w_i^\alpha)\in \Gamma\left( U_i, \bigwedge^{p+1} \frac{\Omega_M^1}{\mathcal{N}_{\mathcal{F}_0}^*} \right) = \Gamma\left( U_i, \mathscr{H}om_{\mathcal{O}_M} \left( \bigwedge^{p+1} \Theta_{\mathcal{F}_0}, \mathcal{O}_M \right) \right)
\end{align*}
Then for $X_1,,,,, X_{p+1}\in  \Gamma\left( U_i, \Theta_{\mathcal{F}_0}\right)$,
\begin{align*}
&\left(\alpha_{p+1}\circ D_p(\phi)\left(w_i^\alpha\right) \right)\left(X_1,...,X_{p+1}\right)\\
&=i_{\left( \sum_{i=1}^{p+1} (-1)^{i+1} \left[X_i, \phi\left(X_1,...,\hat{X}_i,...,X_{p+1}\right)\right]+\sum_{i<j}(-1)^{i+j}\phi\left(\left[X_i,X_j\right],...,\hat{X}_i,...,\hat{X}_j,...,X_{p+1}\right)
  \right)} \left(w_i^\alpha\right) 
\end{align*}
On the other hand,
\begin{align}\label{ap1}
E_p\circ \alpha_p(\phi)(w_i^\alpha)=d\left(\alpha_p(\phi)(w_i^\alpha)\right)-\sum_{\beta=1}^q a_i^{\alpha\beta} \wedge \left( \alpha_p(\phi)\left(w_i^\beta\right) \right) \in \Gamma\left( U_i, \bigwedge^{p+1}\frac{\Omega_M^1}{\mathcal{N}_{\mathcal{F}_0}^*} \right) = \Gamma\left( U_i, \mathscr{H}om_{\mathcal{O}_M}\left( \bigwedge^{p+1} \Theta_{\mathcal{F}_0}, \mathcal{O}_M \right) \right)
\end{align}
Let us consider the first term of $\textnormal{(\ref{ap1})}$. Then for $X_1,...,X_{p+1}\in \Gamma\left( U_i, \Theta_{\mathcal{F}_0} \right)$, we have (here we recall the formulas: $\mathcal{L}_{X_k}=di_{X_k}+ i_{X_k}d$ and $i_{[X_k,X_l]}=\mathcal{L}_{X_k}i_{X_l}- i_{X_l}\mathcal{L}_{X_k}$)
{\small{\begin{align*}
&(d(\alpha_p(\phi)(w_i^\alpha)))(X_1,...,X_{p+1})=i_{X_1\wedge \cdots \wedge X_{p+1}}(d(\alpha_p(\phi)(w_i^\alpha)))=i_{X_2\wedge \cdots \wedge X_p} i_{X_1} (d(\alpha_p(\phi)(w_i^\alpha)) )\\
&=i_{X_2\wedge \cdots \wedge X_{p+1}}(\mathcal{L}_{X_1}(\alpha_p(\phi)(w_i^\alpha))- di_{X_1}(\alpha_p(\phi)(w_i^\alpha)))=i_{X_3\wedge \cdots \wedge X_{p+1}} i_{X_2}(\mathcal{L}_{X_1}(\alpha_p(\phi)(w_i^\alpha))) - i_{X_3\wedge \cdots \wedge X_{p+1}}i_{X_2}di_{X_1}(\alpha_p(\phi)(w_i^\alpha))\\
&=i_{X_3\wedge \cdots \wedge X_{p+1}}\mathcal{L}_{X_1}i_{X_2}(\alpha_p(\phi)(w_i^\alpha))-i_{X_3\wedge \cdots \wedge X_{p+1}}i_{[X_1,X_2]}(\alpha_p(\phi)(w_i^\alpha))\\
&-i_{X_3\wedge \cdots \wedge X_{p+1}}\mathcal{L}_{X_2}i_{X_1}(\alpha_p(\phi)(w_i^\alpha))+i_{X_3\wedge \cdots \wedge X_{p+1}}d i_{X_2}i_{X_1}(\alpha_p(\phi)(w_i^\alpha))\\
&=i_{X_4\wedge \cdots \wedge X_{p+1}}i_{X_3}\mathcal{L}_{X_1} i_{X_2}(\alpha_p(\phi)(w_i^\alpha))-i_{\phi\left(\left[X_1,X_2\right],...,X_{p+1}\right)}(w_i^\alpha)\\
&-i_{X_4\wedge \cdots \wedge X_{p+1}} i_{X_3}\mathcal{L}_{X_2}i_{X_1}(\alpha_p(\phi)(w_i^\alpha))+i_{X_4\wedge \cdots \wedge X_{p+1}}i_{X_3}d i_{X_2}i_{X_1}(\alpha_p(\phi)(w_i^\alpha))\\
&=\mathcal{L}_{X_1}i_{X_2\wedge \cdots\wedge X_{p+1}}(\alpha_p(\phi)(w_i^\alpha))  - \sum_{k=2}^{p+1}    i_{\phi\left(X_2,...,\left[X_1,X_k\right],...,X_{p+1}\right)}(w_i^\alpha)\\
& - \mathcal{L}_{X_2} i_{X_1\wedge \hat{X}_2\wedge \cdots \wedge X_{p+1}} (\alpha_p(\phi)(w_i^\alpha)) +\sum_{k=3}^{p+1}i_{\phi\left(X_1,..., \left[X_2,X_k\right],...,X_{p+1}\right)}(w_i^\alpha)             \\
& +i_{X_4\wedge \cdots \wedge X_{p+1}}\mathcal{L}_{X_3}i_{X_2}i_{X_1}(\alpha_p(\phi)(w_i^\alpha)) -i_{X_4\wedge \cdots \wedge X_{p+1}}di_{X_3}i_{X_2}i_{X_1}(\alpha_p(\phi)(w_i^\alpha))  \\
&=\sum_{k=1}^{p+1}(-1)^{k-1} \mathcal{L}_{X_k} i_{X_1\wedge \cdots \wedge \hat{X}_k\wedge \cdots X_{p+1}}(\alpha_p(\phi)(w_i^\alpha)) +\sum_{i<j}  (-1)^i i_{\phi\left(X_1,...,\hat{X}_i....,\left[X_i,X_j\right],X_{j+1},...,X_{p+1} \right) }(w_i^\alpha) \\
&= \sum_{k=1}^{p+1}(-1)^{k-1} \mathcal{L}_{X_k} i_{\phi(X_1, \cdots, \hat{X}_k, \cdots X_{p+1})}(w_i^\alpha) +\sum_{i<j}  (-1)^{i+j} i_{\phi \left(\left[X_i,X_j\right],X_1,...,\hat{X}_i....,\hat{X}_j,...,X_{p+1}     \right)}(w_i^\alpha)
\end{align*}}}

On the other hand, let us consider the second term of $\textnormal{(\ref{ap1})}$. We note that from $dw_i^\alpha= \sum_{\beta=1}^q a_i^{\alpha\beta} \wedge w_i^\beta$, we have 
\begin{align*}
\mathcal{L}_{X_k}(w_i^\alpha) = di_{X_k}\left(w_i^\alpha\right)+i_{X_k}\left(dw_i^\alpha\right)=\sum_{\beta=1}^q i_{X_k}\left(a_i^{\alpha\beta}\right)w_i^\beta-\sum_{\beta=1}^qa_i^{\alpha\beta} i_{X_k}\left(w_i^\beta\right)= \sum_{\beta=1}^q i_{X_k}\left( a_i^{\alpha\beta} \right) w_i^\beta\,\,\,\,\,\,\,\textnormal{for $X_k\in \Gamma(U_i, \Theta_{\mathcal{F}_0})$}
\end{align*}
Then for $X_1,..., X_{p+1}\in \Gamma\left( U_i, \Theta_{\mathcal{F}_0} \right)$, we have
\begin{align*}
&\left(\sum_{\beta=1}^q a_i^{\alpha\beta}\wedge \left( \alpha_p(\phi)\left( w_i^\beta \right) \right)\right)\left( T_1 , ... , T_{p+1} \right) =\sum_{\beta=1}^q i_{X_1\wedge \cdots \wedge X_{p+1}}\left(a_i^{\alpha\beta}\wedge \left(\alpha_p(\phi)\left(w_i^\beta\right)   \right) \right)\\
&=\sum_{\beta=1}^q i_{X_2\wedge \cdots \wedge X_{p+1}}\left(i_{X_1}\left(a_i^{\alpha\beta}\right) a_p(\phi)\left(w_i^\beta\right)-a_i^{\alpha\beta}\wedge i_{X_1}\left(\alpha_p(\phi)\left(w_i^\beta\right)\right)\right)\\
&=\sum_{\beta=1}^q i_{X_1}\left(a_i^{\alpha\beta}\right)i_{\phi\left(\hat{X}_1,X_2,...,X_{p+1}\right)} \left(w_i^\beta \right)-\sum_{\beta=1}^qi_{X_3\wedge \cdots \wedge X_{p+1}}\left( i_{X_2}\left(a_i^{\alpha\beta}\right)i_{X_1}\left(\alpha_p(\phi)\left(w_i^\beta\right)\right)- a_i^{\alpha\beta} \wedge i_{X_1\wedge X_2}\left(\alpha_p(\phi)\left(w_i^\alpha\right)\right)  \right)\\
&=i_{\phi\left(\hat{X}_1,X_2,...,X_{p+1}\right)}\left( \sum_{\beta=1}^q i_{X_1}\left(a_i^{\alpha\beta}\right) w_i^\beta \right)-i_{\phi(X_1,\hat{X}_2,...,X_{p+1})}\left(\sum_{\beta=1}^q i_{X_2}\left(a_i^{\alpha\beta}\right) w_i^\beta \right)\\
&+\sum_{\beta=1}^q i_{X_4\wedge \cdots \wedge X_{p+1}}\left(i_{X_3}\left(a_i^{\alpha\beta}\right)\wedge i_{X_1\wedge X_2}\left(\alpha_p(\phi)\left(w_i^\beta\right)\right)-a_i^{\alpha\beta}\wedge i_{X_1\wedge X_2\wedge X_3}\left(\alpha_p(\phi)\left(w_i^\beta\right)\right)\right)\\
&=i_{\phi\left(\hat{X}_1,X_2,...,X_{p+1}\right)}\left( \sum_{\beta=1}^q i_{X_1}\left(a_i^{\alpha\beta}\right) w_i^\beta \right)-i_{\phi\left(X_1,\hat{X}_2,...,X_{p+1}\right)}\left(\sum_{\beta=1}^q i_{X_2}\left(a_i^{\alpha\beta}\right) w_i^\beta \right)\\
&+i_{\phi\left(X_1,X_2,\hat{X}_3,...,X_{p+1}\right)}\left(\sum_{\beta=1}^q i_{X_2}\left(a_i^{\alpha\beta}\right)\left(w_i^\beta\right)  \right)-\sum_{k=1}^q i_{X_5\wedge \cdots \wedge X_{p+1}}\left(i_{X_4}\left(a_i^{\alpha\beta}\right)-a_i^{\alpha\beta}i_{X_1\wedge \cdots \wedge X_4}\left(\alpha_p(\phi)\left(w_i^\beta\right)\right)\right)\\
&=\sum_{k=1}^{p+1} (-1)^{k-1} i_{\phi\left(X_1,...,\hat{X}_k,...,X_{p+1}\right)}\left(\sum_{\beta=1}^q i_{X_k}\left(a_i^{\alpha\beta}\right)w_i^\beta  \right)=\sum_{k=1}^{p+1}(-1)^{k-1} i_{\phi \left( X_1,...,\hat{X}_k,...,X_{p+1} \right) }\mathcal{L}_{X_k}\left(w_i^\alpha\right)
\end{align*}

Hence we have
\begin{align*}
&\left(E_i\circ \alpha_p(\phi)(w_{0i}^\alpha) \right)\left(X_1,...,X_{p+1}\right)=\left(d\alpha_p(\phi)\left(w_i^\alpha\right)-\sum_{\beta=1}^q a_i^{\alpha\beta} \wedge \left( \alpha_p(\phi)\left(w_i^\beta\right) \right)\right)\left(X_1,...,X_{p+1}\right)\\
&=\sum_{k=1}^{p+1}(-1)^{k-1}\left(\mathcal{L}_{X_k}i_{\phi\left(X_1,...,\hat{X}_k,...,X_{p+1}\right)} - i_{\phi\left(X_1,...,\hat{X}_k,...,X_{p+1}\right)}\mathcal{L}_{X_k} \right)(w_i^\alpha)+\sum_{i<j}  (-1)^{i+j} i_{\phi\left(\left[X_i,X_j\right],X_1,...,\hat{X}_i....,\hat{X}_j,...,X_{p+1}     \right)}(w_i^\alpha)\\
&=i_{\sum_{k=1}^{p+1}(-1)^{k+1}\left[X_k,\phi\left(X_1,...,\hat{X}_k,...,X_{p+1}\right)\right]+\sum_{i<j}(-1)^{i+j} \phi\left(\left[X_i,X_j\right],...,\hat{X}_i,...,\hat{X}_j,...,X_{p+1}\right)}\left(w_i^\alpha\right)
\end{align*}
which shows the commutativity of the above diagram.

\end{proof}

\subsection{The dual leaf complex of $\left( M, \mathcal{N}_{\mathcal{F}_0}^*\right)$ with $\mathcal{N}_{\mathcal{F}_0}^*$ locally free} \label{AppendixA5}\

Let $(M, \mathcal{N}_{\mathcal{F}_0}^*)$ be a compact foliated complex manifold with $\mathcal{N}_{\mathcal{F}_0}^*$ locally free of rank $q$. We shall construct the dual leaf complex of $(M, \mathcal{N}_{\mathcal{F}_0}^*)$ which extends the dual leaf complex for a regular foliation on a compact complex manifold in subsection \ref{AppendixA2}.

We note that we have a section $\omega_0\in \Gamma\left( M, \bigwedge^q\Omega_M^1\otimes \mathcal{L}_0 \right)$ where $\mathcal{L}_0= \bigwedge^q \left( \mathcal{N}_{\mathcal{F}_0}^* \right)^*$, which is locally decomposable and satisfy integrability condition. Let $\mathcal{U}=\{U_i\}$ be an open covering of $M$ by coordinate neighborhoods such that $\Gamma(U_i, \mathcal{N}_{\mathcal{F}_0}^*)$  is generated by $w_i^1,..., w_i^q\in \Gamma(U_i, \mathcal{N}_{\mathcal{F}_0}^*)$ with the relation $w_i^\alpha= \sum_{\beta=1}^q h_{ij}^{\alpha\beta}w_j^\beta$ on $U_i\cap U_j$ for $h_{ij}^{\alpha\beta}\in \Gamma(U_{ij}, \mathcal{O}_M)$ and satisfy $w_i^1\wedge \cdots \wedge w_i^q \wedge dw_i^\alpha=0,\alpha=1,...,q$. We set $w_i:= w_i^1\wedge \cdots \wedge w_i^q \in \Gamma\left( U_i, \bigwedge^q \Omega_M^1 \right)$. Then $w_i=\det\left( h_{ij}^{\alpha\beta} \right) w_j$ on $U_i\cap U_j$. Then $w_0$ is defined by $\{w_i\}$. We define the subsheaf $\tilde{\mathcal{S}}^r\subset \bigwedge^{q+r}\Omega_M^1\otimes \mathcal{L}_0$ by the following property that outside of $S=\textnormal{Sing}(\mathcal{F}_0)$ (we recall that codimension of $S$ is at least $2$), it is locally of the form $\omega_0\wedge A$ where $A\in \bigwedge^r \Omega_M^1$. More precisely, for an open set $U$ of $M$, we define for $r\geq 1$,
{\Small{\begin{align*}
&\Gamma(U, \tilde{\mathcal{S}}^r)\\
&=\left\{ a\in \Gamma\left(U, \bigwedge^{q+r}\Omega_M^1\otimes L\right)|\text{for any $x\in U-S$, there exists $U_x\subset U-S$ such that $a=\omega_0\wedge b_x$ on $U_x$ for some $b_x\in \Gamma\left(U_x, \bigwedge^r \Omega_M^1\right)$}        \right\}
\end{align*}}}

Then we claim that $\tilde{S}^r|_{M-S}\cong \bigwedge^r \frac{\Omega_M^1}{\mathcal{N}_{\mathcal{F}_0}^*}|_{M-S}$. For an open set of $U\subset M-S$, we define a homomorphism $\alpha:\bigwedge^r \frac{\Omega_M^1}{\mathcal{N}_{\mathcal{F}_0}^*}|_{M-S}\to \tilde{\mathcal{S}}|_{M-S}$ in the following way:
\begin{align*}
\alpha(U):\Gamma\left( U, \bigwedge^r \frac{\Omega_M^1}{\mathcal{N}_{\mathcal{F}_0}^*}|_{M-S} \right) \to \Gamma\left( U,  \tilde{\mathcal{S}}^r|_{M-S} \right)
\end{align*}
Let $A\in \Gamma\left( U, \bigwedge^r \frac{\Omega_M^1}{\mathcal{N}_{\mathcal{F}_0}^*}|_{M-S} \right)$. Then for any $x\in U$, there exists a neighborhood of $U_x\subset U$ such that $A|_{U_x}$ has a lifting $\tilde{A}_x\in \Gamma\left(U_x, \bigwedge^r \Omega_X^1 \right)$. Then $\omega_0\wedge \tilde{A}_x\in \Gamma\left(U_x, \tilde{\mathcal{S}}^r|_{M-S} \right)$. For another lifting $\tilde{A}_x'$ of $A|_{U_x}$, we have $\tilde{A}_x- \tilde{A}_x'\in \mathcal{N}_{\mathcal{F}_0}^*\bigwedge \left(\bigwedge^{r-1} \Omega_M^1\right)$, so that $\omega_0\wedge \tilde{A}_x= \omega_0 \wedge \tilde{A}'_x$. By the same argument, $\omega_0 \wedge \tilde{A}_x= \omega_0 \wedge \tilde{A}_y$ on $U_x\cap U_y\ne \emptyset$, so that $\left\{\omega_0 \wedge A_x |x\in U \right\}$ glues together to define an element $\alpha(U)(A)\in \Gamma\left(U, \tilde{S}^r|_{M-S} \right)$. Now we show that $\alpha$ is an isomorphism. If $\omega_0 \wedge \tilde{A}_x=0$, then $\tilde{A}_x\in \mathcal{N}_{\mathcal{F}_0}^*\wedge \left(\bigwedge^{r-1} \Omega_M^1 \right)$, and $A|_{U_x}=0$. On the other hand, for any $\omega_0 \wedge B_x\in \Gamma\left(U_x, \tilde{\mathcal{S}}^r \right)$ for $B_x\in \Gamma\left(U_x, \bigwedge^r \Omega_X^1 \right)$. Let $\bar{B}_x$ is the image of $B_x$ in $\Gamma\left(U_x, \bigwedge^r \frac{\Omega_M^1}{\mathcal{N}_{\mathcal{F}_0}^*} \right)$. Then $\alpha(U_x)(\bar{B}_x)=\omega_0 \wedge B_x$.

Now we shall define the dual leaf complex for $(M, \mathcal{N}_{\mathcal{F}_0}^*)$ with $\mathcal{N}_{\mathcal{F}_0}^*$ locally free as above:
\begin{align} \label{ap2}
\mathcal{N}_{\mathcal{F}_0}^{*\bullet}: \Theta_M \xrightarrow{E_0} \mathscr{H}om_{\mathcal{O}_M}\left( \mathcal{N}_{\mathcal{F}_0}^*, \frac{\Omega_M^1}{\mathcal{N}_{\mathcal{F}_0}^*} \right) \xrightarrow{E_1} \mathscr{H}om_{\mathcal{O}_M}\left( \mathcal{N}_{\mathcal{F}_0}^*, \tilde{\mathcal{S}}^2 \right) \xrightarrow{E_2} \mathscr{H}om_{\mathcal{O}_M}\left( \mathcal{N}_{\mathcal{F}_0}^*, \tilde{\mathcal{S}}^3\right) \xrightarrow{E_3} \to \cdots
\end{align}
For $X\in \Theta_M$, $E_0(X)$ is locally defined on $U_i$ by
\begin{align*}
E_0(X)(w_i^\alpha)= \overline{\mathcal{L}_X(w_i^\alpha)}, \,\,\,\,\,\alpha=1,...,q
\end{align*}
as in the dula leaf complex for a regular foliation in subsection \ref{AppendixA2}, and $E_1$ is defined in the following way: for $\phi\in \mathscr{H}om_{\mathcal{O}_M}\left( \mathcal{N}_{\mathcal{F}_0}^*, \frac{\Omega_M^1}{\mathcal{N}_{\mathcal{F}_0}^*} \right)$, $E_1(\phi)$ is locally defined on $U_i$ by
{\small{\begin{align*}
E_1(\phi)(w_i^\alpha)= \sum_{\beta=1}^q w_i^1\wedge \cdots \wedge \widetilde{\phi\left(w_i^\beta\right)} \wedge \cdots \wedge w_i^q \wedge dw_i^\alpha + w_i^1\wedge \cdots \wedge w_i^q \wedge d \widetilde{\phi\left(w_i^\alpha \right)}\in \Gamma\left(U_i, \bigwedge^{q+2}\Omega_M^1\right)\cong \Gamma\left(U_i, \bigwedge^{q+2}\Omega_M^1\otimes \mathcal{L}_0 \right)
\end{align*}}}
for $\alpha=1,...,q$, and linearly extends to $\Gamma\left(U_i, \mathcal{N}_{\mathcal{F}_0}^*\right)$. Here $\widetilde{\phi\left(w_i^\alpha\right)}$ is a lifting of $\phi\left(w_i^\alpha\right)$ under $\Omega_M^1 \to \frac{\Omega_M^1}{\mathcal{N}_{\mathcal{F}_0}^*}$.\footnote{Strictly speaking $\Omega_M^1\to \frac{\Omega_M^1}{\mathcal{N}_{\mathcal{F}_0}^*}$ is not surjective on $U_i$. More precisely for $x\in U_i$, choose a neighborhood $U_x$ of $U_i$ such that $\phi(w_i^\alpha)$ on $U_x$ has a lifting $\widetilde{\phi_x(w_i^\alpha)}$. Then we compute $ \sum_{\beta=1}^q w_i^1\wedge \cdots \wedge \widetilde{\phi_x\left(w_i^\beta\right)} \wedge \cdots \wedge w_i^q \wedge dw_i^q + w_i^1\wedge \cdots \wedge w_i^q \wedge d \widetilde{\phi_x\left(w_i^\alpha \right)}$ locally on $U_x$ for each $x\in U_i$, since they are equal on $U_x\cap U_y\ne \emptyset$ by the above arguments, we can glue together to define $E_p(\phi)(w_i^\alpha)$ on $U_i$} We show that $E_1(\phi)(w_i^\alpha)$ is independent of choices of a lifting of $\phi(w_i^\alpha)$. In fact, for $A_1,...,A_q\in \mathcal{N}_{\mathcal{F}_0}^*$, we have $\sum_{\beta=1}^q w_i^1 \wedge \cdots \wedge A_\beta \wedge w_i^q \wedge dw_i^\alpha + w_i^1 \wedge \cdots \wedge w_i^q \wedge dA_\alpha= 0$ on $U_i-S$, and hence $0$ on $U_i$. We show that $E_1(\phi_i)(w_i^\alpha)$ actually lies in $\Gamma\left(U_i, \tilde{\mathcal{S}}^2 \right)$. For each $x\in U_i-S$,  choose a neighborhood $U_x\subset U_i-S$ such that $dw_i^\alpha= \sum_{\beta=1}^q a_{i_x}^{\alpha\beta} \wedge w_i^\beta$ for some $a_{i_x}^{\alpha\beta}\in \Gamma(U_x, \Omega_M^1)$. Then we see that on $U_x$
\begin{align*}
&\sum_{\beta=1}^q w_i^1\wedge \cdots \wedge \widetilde{\phi\left(w_i^\beta\right)} \wedge \cdots \wedge w_i^q \wedge dw_i^\alpha + w_i^1\wedge \cdots \wedge w_i^q \wedge d \widetilde{\phi\left(w_i^\alpha \right)}\\
&= \sum_{\beta=1}^q w_i^1\wedge \cdots \wedge \widetilde{\phi\left(w_i^\beta\right)} \wedge \cdots \wedge w_i^q \wedge \left(\sum_{\gamma=1}^q a_{i_x}^{\alpha\gamma} \wedge w_i^\gamma \right)+ w_i^1\wedge \cdots \wedge w_i^q \wedge d \widetilde{\phi\left(w_i^\alpha \right)}\\
&= w_i^1\wedge \cdots \wedge w_i^q\wedge \left( d\widetilde{\phi(w_i^\alpha)} - \sum_{\beta=1}^q a_{i_x}^{\alpha \beta} \wedge \widetilde{\phi(w_i^\beta)} \right) \in \Gamma(U_x, \tilde{\mathcal{S}}^r)
\end{align*}

Next we show that $E_1(\phi)$ is well-defined: in other words, $E_1(\phi)$ on $U_i$ and $E_1(\phi)$ on $U_j$ defines the same homomorphism on $U_i\cap U_j$. We recall that $w_i^\alpha= \sum_{\beta=1}^q h_{ij}^{\alpha\beta} w_j^\beta$. By abuse of notation, we write $\phi(w_i^\alpha)$ instead of $\widetilde{\phi(w_i^\alpha)}$. We set $H_{ij}=\det\left( h_{ij}^{\alpha\beta} \right)$.
\begin{align*}
&E_1(\phi)(w_i^\alpha)= \sum_{\beta=1}^q w_i^1\wedge \cdots \wedge \phi\left(w_i^\beta\right) \wedge \cdots \wedge w_i^q \wedge dw_i^\alpha + w_i^1\wedge \cdots \wedge w_i^q \wedge d \phi\left(w_i^\alpha \right)\\
&=\sum_{\beta=1}^q H_{ij} w_j^1 \wedge \cdots \wedge \phi\left(w_j^\beta\right)\wedge \cdots \wedge w_j^q \wedge d\left(\sum_{\gamma=1}^q h_{ij}^{\alpha\gamma} w_j^\gamma    \right) + H_{ij} w_j^1 \wedge \cdots \wedge w_j^q\wedge d\left(   \sum_{\beta=1}^q h_{ij}^{\alpha\beta} \phi\left(w_j^\beta\right)  \right)\\
&=\sum_{\beta=1}^q H_{ij} w_j^1 \wedge \cdots \wedge \phi\left(w_j^\beta\right)\wedge \cdots \wedge w_j^q \wedge \left(\sum_{\gamma=1}^q d h_{ij}^{\alpha\gamma} \wedge w_j^\gamma  + \sum_{\gamma=1}^q h_{ij}^{\alpha \gamma}  dw_j^\gamma  \right)\\
& + H_{ij} w_j^1 \wedge \cdots \wedge w_j^q\wedge \left(   \sum_{\beta=1}^q dh_{ij}^{\alpha\beta} \wedge \phi\left(w_j^\beta\right) + \sum_{\beta=1}^q h_{ij}^{\alpha \beta} d \phi\left(w_j^\beta \right)  \right) \\
&= H_{ij}\sum_{\beta=1}^q h_{ij}^{\alpha\beta}\left(  \sum_{\gamma=1}^q w_j^1 \wedge \cdots \wedge \phi\left( w_j^\gamma \right) \wedge \cdots \wedge w_j^q \wedge d w_j^\beta + w_j^1 \wedge \cdots \wedge w_j^q \wedge d\phi\left(w_j^\beta \right)   \right)= H_{ij}\sum_{\beta=1}^q h_{ij}^{\alpha\beta} E_p(\phi)\left(w_j^\beta\right)
\end{align*}
Hence $E_1$ is well-defined.

For $r>1$, we define $E_r:\mathscr{H}om_{\mathcal{O}_M}\left(\mathcal{N}_{\mathcal{F}_0}^*, \tilde{\mathcal{S}}^r\right) \to \mathscr{H}om_{\mathcal{O}_M}\left( \mathcal{N}_{\mathcal{F}_0}^*, \tilde{\mathcal{S}}^{r+1} \right)$. For $\phi\in \mathscr{H}om_{\mathcal{O}_M}\left(\mathcal{N}_{\mathcal{F}_0}^*, \tilde{\mathcal{S}}^r \right)$, $E_r(\phi)$ is locally defined on $U_i$ in the following way: for any $x\in U_i-S$, choose a neighborhood $U_x\subset U_i-S$, then we can write $dw_i^\alpha= \sum_{\beta=1}^q a_{i_x}^{\alpha\beta} \wedge w_i^\beta$. Since $\phi(w_i^\alpha)\in \Gamma\left(U_i, \tilde{\mathcal{S}}^r\right)$, we have $\phi(w_i^\alpha)= w_i^1\wedge \cdots \wedge w_i^q \wedge  B_{i_x}^\alpha$ for some $B_{i_x}^\alpha \in \Gamma\left(U_{i_x}, \bigwedge^r \Omega_M \right)$. Then we set
\begin{align*}
A_{i_x}^\alpha := w_i^1\wedge \cdots \wedge w_i^q \wedge \left( d B_{i_x}^\alpha  - \sum_{\beta=1}^q a_{i_x}^{\alpha\beta} \wedge B_{i_x}^\beta  \right) \in \Gamma\left( U_x, \tilde{\mathcal{S}}^{r+1} \right)
\end{align*}
As in the construction of the dual leaf complex for a regular foliation, $A_i^x$ is independent of choice of $a_{i_x}^{\alpha\beta}$ and $B_{i_x}^\alpha$. We note that on the intersection $U_x\cap U_y \ne \emptyset$ for $x,y\in U_i-S$, we have $w_i^1\wedge \cdots \wedge w_i^q \wedge B_{i_x}^\alpha = w_i^1 \wedge \cdots \wedge w_i^q \wedge B_{i_y}^\alpha $, so that
$B_{i_x}^\alpha - B_{i_y}^\alpha \in \mathcal{N}_{\mathcal{F}_0}^*$. This implies that $A_{i_x}^\alpha= A_{i_y}^\alpha$, so that we can glue together to define $A_i^\alpha$ on $U_i-S$, and so $A_i^\alpha\in \Gamma\left(U_i, \tilde{\mathcal{S}}^{r+1} \right)\subset \Gamma\left( U_i, \bigwedge^{q+r+1}\Omega_M^1 \otimes \mathcal{L}_0 \right)$. Now we define $E_r(\phi)(w_i^\alpha)= A_i^\alpha$, and linearly extends to $\Gamma\left( U_i, \mathcal{N}_{\mathcal{F}_0}^* \right)$. Next we show that $E_r$ is well-defined: in other words, $E_r(\phi)$ on $U_i$ and $E_r(\phi)$ on $U_j$ defines the same homomorphism on $U_i\cap U_j$. Choose $x\in U_i\cap U_j$ and choose a neighborhood $U_x\subset (U_i\cap U_j)-S$. Then From $\phi(w_i^\alpha)=  \sum_{\beta=1}^q h_{ij}^{\alpha\beta} \phi\left(w_j^\beta \right)$, we see that $w_i^1 \wedge \cdots \wedge w_i^q \wedge B_{i_x}^\alpha= H_{ij} \cdot  \sum_{\beta=1}^q h_{ij}^{\alpha\beta} \left(w_j^1 \wedge \cdots \wedge w_j^q \wedge B_{j_x}^\beta \right)$ on $U_x$. Then $w_i^1\wedge \cdots \wedge w_i^q \wedge \left( B_{i_x}^\alpha - \sum_{\beta=1}^q h_{ij}^{\alpha\beta} B_{j_x}^\beta   \right)=0$, so that we have $B_{i_x}^\alpha - \sum_{\beta=1}^q h_{ij}^{\alpha\beta} B_{j_x}^\beta \in \mathcal{N}_{\mathcal{F}_0}^*$. We also note that 
\begin{align*}
\sum_{\beta, \gamma=1}^q h_{ij}^{\beta \gamma} a_{i_x}^{\alpha\beta} \wedge w_j^\gamma  =\sum_{\beta=1}^q a_{i_x}^{\alpha\beta} \wedge w_i^\beta =  dw_i^\alpha= d\left(\sum_{\beta=1}^q h_{ij}^{\alpha\beta} w_j^\beta \right)= \sum_{\gamma=1}^q d h_{ij}^{\alpha\gamma}  \wedge w_j^\gamma + \sum_{\beta=1}^q h_{ij}^{\alpha\beta} \left( \sum_{\gamma=1}^q a_{j_x}^{\beta \gamma}  \wedge w_j^\gamma    \right)
 \end{align*}
 Then we have $\sum_{\gamma=1}^q \left( \sum_{\beta=1}^q h_{ij}^{\beta \gamma} a_{i_x}^{\alpha\beta} - d h_{ij}^{\alpha \gamma} - \sum_{\beta=1}^q h_{ij}^{\alpha\beta} a_{j_x}^{\beta \gamma}      \right) \wedge w_j^\gamma =0$, so that $ \sum_{\beta=1}^q h_{ij}^{\beta \gamma} a_{i_x}^{\alpha\beta} - d h_{ij}^{\alpha \gamma} - \sum_{\beta=1}^q h_{ij}^{\alpha\beta} a_{j_x}^{\beta \gamma}  \in \mathcal{N}_{\mathcal{F}_0}^*$. Then we have
\begin{align*}
E_r(\phi)(w_i^\alpha) &= A_{i_x}^\alpha= w_i^1 \wedge \cdots\wedge w_i^q \wedge \left( d B_{i_x}^\alpha - \sum_{\beta=1}^q a_{i_x}^{\alpha\beta} \wedge B_{i_x}^\beta    \right) \\
&= w_i^1 \wedge \cdots \wedge w_i^q \wedge \left( d\left(\sum_{\beta=1}^q h_{ij}^{\alpha\beta} B_{j_x}^\beta \right) - \sum_{\beta=1}^q a_{i_x}^{\alpha\beta} \wedge \left(  \sum_{\gamma=1}^q h_{ij}^{\beta \gamma} B_{j_x}^\gamma   \right)    \right) \\
&= H_{ij} w_j^1 \wedge \cdots \wedge w_j^q\wedge \left( \sum_{\beta=1}^q d h_{ij}^{\alpha\beta} \wedge B_{j_x}^\beta + \sum_{\beta=1}^q h_{ij}^{\alpha\beta} d B_{j_x}^\beta - \sum_{\gamma=1}^q \left(  d h_{ij}^{\alpha \gamma} + \sum_{\beta=1}^q h_{ij}^{\alpha\beta}  a_{j_x}^{\beta \gamma}     \right) B_{j_x}^\gamma      \right) \\
&= H_{ij} \sum_{\beta=1}^q h_{ij}^{\alpha\beta}\left( w_j^1 \wedge \cdots \wedge w_j^q \wedge \left( dB_{j_x}^\beta - \sum_{\gamma=1}^q a_{j_x}^{\beta \gamma} B_{j_x}^\gamma \right)    \right) = H_{ij} \sum_{\beta=1}^q h_{ij}^{\alpha\beta} A_{j_x}^\beta = H_{ij} \sum_{\beta=1}^q h_{ij}^{\alpha\beta} E_r(\phi)\left(w_j^\beta\right)
\end{align*}
Hence $E_r(\phi)$ is well-defined.

We show that $E_1\circ E_0=0$. In fact, for $X\in \Theta_M$,
\begin{align*}
E_1(E_0(X))(w_i^\alpha) &= \sum_{\beta=1}^q w_i^1 \wedge \cdots \wedge \mathcal{L}_X\left(w_i^\beta\right) \wedge \cdots \wedge w_i^q \wedge dw_i^\alpha + \sum_{\beta=1}^q w_i^1 \wedge \cdots \wedge w_i^q \wedge d \mathcal{L}_X(w_i^\alpha) \\
 &= \mathcal{L}_X \left( w_i^1 \wedge \cdots \wedge w_i^q \wedge dw_i^\alpha \right) = 0
\end{align*}

We show that for $r\geq 1$,  $E_{r+1}\circ E_r =0$. In fact, it is enough to check locally on $U_x\subset M-S$ for each $x\in M-S$. Let $x\in U_i-S$ for some $i$ and choose $U_x\subset U_i-S$. Then from $dw_i^\alpha= \sum_{\beta=1}^q a_{i_x}^{\alpha\beta} \wedge w_i^\beta$, we have $d a_{i_x}^{\alpha\beta}- \sum_{\gamma=1}^q  a_{i_x}^{\alpha \gamma} \wedge a_{i_x}^{\gamma \beta}\in \mathcal{N}_{\mathcal{F}_0}^* \wedge \Omega_M^1$. Then for  $\phi\in \mathscr{H}om_{\mathcal{O}_M}\left( \mathcal{N}_{\mathcal{F}_0}^*, \frac{\Omega_M^1}{\mathcal{N}_{\mathcal{F}_0}^*} \right)$ or $\mathscr{H}om_{\mathcal{O}_M}\left(\mathcal{N}_{\mathcal{F}_0}^*, \tilde{\mathcal{S}}^r \right) (r \geq 2)$, let $\phi(w_i^\alpha)= w_i^1 \wedge \cdots \wedge w_i^q \wedge A_{i_x}^\alpha$ for some $A_{i_x}^\alpha \in \Gamma\left( U_x, \bigwedge^r \Omega_M^1 \right)$. Then $E_p(\phi)(w_i^\alpha)= w_i^1\wedge \cdots \wedge w_i^q \wedge \left( d A_{i_x}^\alpha - \sum_{\beta=1}^q a_{i_x}^{\alpha\beta} A_{i_x}^\beta \right)$. Then
\begin{align*}
&E_{r+1}(E_r(\phi))(w_i^\alpha) = w_i^1 \wedge \cdots \wedge w_i^q \wedge \left( d\left( dA_{i_x}^\alpha - \sum_{\beta=1}^q a_{i_x}^{\alpha\beta} \wedge A_{i_x}^\beta \right) - \sum_{\gamma=1}^q  a_{i_x}^{\alpha\gamma} \wedge \left( dA_{i_x}^\gamma - \sum_{\beta=1}^q a_{i_x}^{\gamma \beta} \wedge A_{i_x}^{\gamma \beta} \right) \right) \\
&= w_i^1 \wedge \cdots \wedge w_i^q \wedge \left(  - \sum_{\beta=1}^q d a_{i_x}^{\alpha\beta} \wedge A_{i_x}^\beta + \sum_{\beta=1}^q a_{i_x}^{\alpha\beta}  \wedge d A_{i_x}^\beta - \sum_{\gamma=1}^q  a_{i_x}^{\alpha \gamma} \wedge d A_{i_x}^\gamma + \sum_{\beta=1}^q \left( \sum_{\gamma=1}^q a_{i_x}^{\alpha \gamma} \wedge a_{i_x}^{\gamma \beta}  \right) \wedge A_{i_x}^{\gamma \beta}   \right) = 0
\end{align*}
Hence $\textnormal{(\ref{ap2})}$ defines a complex of sheaves. 

\section{The complex of sheaves controlling simultaneous deformations of a compact foliated complex manifolds defined by both locally free subsheaves of tangents sheaf and cotangent sheaf}\label{AppendixA4}

Let $M$ be a compact foliated complex manifold with a singular holomorphic foliation $\mathcal{F}_0$ on $M$ such that the associated subsheaves $\Theta_{\mathcal{F}_0}$ of $\Theta_M$ and $\mathcal{N}_{\mathcal{F}_0}^*$ of $\Omega_M^1$ are both locally free of rank $p$ and rank $q$, respectively.

First we shall construct a complexes of sheaves 
\begin{align*}
\Theta_M \xrightarrow{E_0'} \mathscr{H}om_{\mathcal{O}_M}\left( \Theta_{\mathcal{F}_0}, \left( \mathcal{N}_{\mathcal{F}_0}^* \right)^*  \right) \xrightarrow{E_1'} \mathscr{H}om_{\mathcal{O}_M}\left( \bigwedge^2 \Theta_{\mathcal{F}_0}, \left( \mathcal{N}_{\mathcal{F}_0}^* \right)^* \right) \to \mathscr{H}om_{\mathcal{O}_M}\left( \bigwedge^3 \Theta_{\mathcal{F}_0},  \left( \mathcal{N}_{\mathcal{F}_0}^* \right)^* \right) \xrightarrow{E_3'} \cdots
\end{align*}

Let $\mathcal{U}=\{U_i\}$ be an open covering of $M$ by coordinate neighborhoods such that $\mathcal{N}_{\mathcal{F}_0}^*$ on $U_i$ is generated by $w_i^1,..., w_i^q\in \Gamma\left( U_i, \Omega_M^1 \right)$ with the relation $w_i^\alpha=\sum_{\beta=1}^q h_{ij}^{\alpha\beta} w_j^\beta$ for $h_{ij}^{\alpha\beta}\in \Gamma\left( U_{ij}, \mathcal{O}_M \right)$ and $w_i^1\wedge \cdots \wedge w_i^q \wedge dw_i^\alpha=0, \alpha=1,...,q$. We note that for $T \in \Theta_{\mathcal{F}_0}$,  $i_{T}\left(w_i^\beta \right)=0, \beta=1,...q$. For $X\in \Theta_X$, $E_0'(X)$ is locally defined on $U_i$ by
\begin{align*}
E_0'(X)(T)\left(w_i^\beta\right)= i_{\left[ X, T\right]}\left(w_i^\beta \right)
\end{align*}
and linearly extends to $\Gamma\left(U_i, \mathcal{N}_{\mathcal{F}_0}^* \right)$. This it is well-defined since
\begin{align*}
&E_0'(X)(T)\left(w_i^\beta\right)- E_0'(X)\left( T \right)\left( \sum_{\delta=1}^q h_{ij}^{\beta \delta} w_j^\delta\right) = i_{[X, T]}\left(w_i^\beta\right)- \sum_{\delta=1}^q h_{ij}^{\beta \delta} i_{[X, T]}\left(w_j^\delta\right) = 0
\end{align*}
and for $r>0$, $E_r'$ is defined in the following way: for $\phi \in \mathscr{H}om_{\mathcal{O}_M}\left(\bigwedge^r \Theta_{\mathcal{F}_0}, \left( \mathcal{N}_{\mathcal{F}_0}^* \right)^* \right)$, $E_r(\phi)$ is locally defined on $U_i$ by\footnote{We will define the foliated bracket $\left\{ - ,-\right\}_{\left(\mathcal{N}_{\mathcal{F}_0}^* \right)^*}: \Theta_{\mathcal{F}_0}\times \left( \mathcal{N}_{\mathcal{F}_0}^* \right)^*\to \left( \mathcal{N}_{\mathcal{F}_0}^* \right)^* $  in Part V which gives a natural way to define $E_r'$. However we can ignore this at this point.}
{\Small{\begin{align*}
&E_r'(\phi)\left(T_1,,.., T_{r+1}\right)(w_i^\gamma)\\
&\left(= \sum_{\alpha=1}^{r+1} (-1)^{\alpha+1} \left\{ T_\alpha, \phi\left(T_1,..., \hat{T}_\alpha ,..., T_{r+1} \right)    \right\}_{\left(\mathcal{N}_{\mathcal{F}_0}^*\right)^* }(w_i^\gamma) + \sum_{\alpha <\beta} (-1)^{\alpha+\beta} \phi\left(\left[ T_\alpha, T_\beta \right] , T_1,..., \hat{T}_\alpha, ..., \hat{T}_\beta,..., T_{r+1}            \right)(w_i^\gamma)\right) \\
&= \sum_{\alpha=1}^{r+1} (-1)^{\alpha+1}\left[ T_\alpha, \phi\left(T_1,..., \hat{T}_\alpha ,..., T_{r+1} \right) (w_i^\gamma)\right] - \phi(T_1,..., \hat{T}_\alpha,..., T_{r+1})\left( \mathcal{L}_{T_\alpha}\left(w_i^\gamma \right)     \right) + \sum_{\alpha <\beta} (-1)^{\alpha+\beta} \phi\left(\left[ T_\alpha, T_\beta \right] , T_1,..., \hat{T}_\alpha, ..., \hat{T}_\beta,..., T_{r+1}            \right)(w_i^\gamma)
\end{align*}}}
and then we linearly extends to $\Gamma\left( U_i, \mathcal{N}_{\mathcal{F}_0}^* \right)$. We show that this is well-defined. First we note that for $x\in U_i-S$, choose $U_x\subset U_i-S$ such that $dw_i^\gamma= \sum_{\beta=1}^q a_{i_x}^{\gamma\beta} \wedge w_i^\beta$ for some $a_{i_x}^{\gamma\beta}\in \Gamma\left(U_x, \Omega_M^1 \right)$. Then for $T\in \Theta_{\mathcal{F}_0}$, we have $\mathcal{L}_T (w_i^\gamma)= i_T \left( dw_i^\gamma \right)= \sum_{\beta=1}^q i_{T}\left(a_{i_x}^{\gamma \beta} \right) w_i^\beta \in \Gamma(U_x, \mathcal{N}_{\mathcal{F}_0}^*)$, so that $\mathcal{L}_{T_\alpha}(w_i^\gamma)$ is in $\Gamma(U_i-S, \mathcal{N}_{\mathcal{F}_0}^*)$, i.e. in $\Gamma(U_i, \mathcal{N}_{\mathcal{F}_0}^*)$. Hence the second term is well-defined. We show that $E_r'(\phi)$ on $U_i$ and $E_r'(\phi)$ on $U_j$ define the same homomorphisms: In fact, 
{\Small{\begin{align*}
&E_r'(\phi)(T_1,..., T_{r+1}) \left(w_i^\gamma\right)\\
&= \sum_{\alpha=1}^{r+1} (-1)^{\alpha+1}\left[ T_\alpha, \phi\left(T_1,..., \hat{T}_\alpha ,..., T_{r+1} \right) (w_i^\gamma)\right] - \phi(T_1,..., \hat{T}_\alpha,..., T_{r+1})\left( \mathcal{L}_{T_\alpha}(w_i^\gamma)   \right) + \sum_{\alpha <\beta} (-1)^{\alpha+\beta} \phi\left(\left[ T_\alpha, T_\beta \right] , T_1,..., \hat{T}_\alpha, ..., \hat{T}_\beta,..., T_{r+1}            \right)(w_i^\gamma) \\
&= \sum_{\alpha=1}^{r+1}(-1)^{\alpha+1}\left[ T_\alpha, \sum_{\eta=1}^q h_{ij}^{\gamma \eta}\phi\left( T_1,..., \hat{T}_\alpha, ..., T_{r+1} \right)(w_i^\eta) \right]  - \phi\left( T_1,..., \hat{T}_\alpha, ..., T_{r+1} \right)\left( \sum_{\eta =1}^q \left[ T_\alpha, h_{ij}^{\gamma \eta} \right] w_j^\eta + \sum_{\eta=1}^q h_{ij}^{\gamma \eta} \mathcal{L}_{T_\alpha}\left( w_j^\eta\right)      \right)\\
&+ \sum_{\eta=1}^q h_{ij}^{\gamma \eta} \sum_{\alpha <\beta} (-1)^{\alpha+\beta} \phi\left(\left[ T_\alpha, T_\beta \right] , T_1,..., \hat{T}_\alpha, ..., \hat{T}_\beta,..., T_{r+1}            \right)(w_j^\eta)  = \sum_{\eta=1}^q h_{ij}^{\gamma \eta} E_r'(\phi)(T_1,..., T_{r+1}(w_j^\eta)
\end{align*}}}

We show that we have a commutative diagram
{\small{\begin{center}
$\begin{CD}
\Theta_M @>D_0>> \mathscr{H}om_{\mathcal{O}_M}\left(\Theta_{\mathcal{F}_0}, \frac{\Theta_M}{\Theta_{\mathcal{F}_0}}  \right) @>D_1>> \mathscr{H}om_{\mathcal{O}_M}\left(\bigwedge^2 \Theta_{\mathcal{F}_0},  \frac{\Theta_M}{\Theta_{\mathcal{F}_0}} \right) @>D_2>> \mathscr{H}om_{\mathcal{O}_M}\left( \bigwedge^3 \Theta_{\mathcal{F}_0}, \frac{\Theta_M}{\Theta_{\mathcal{F}_0}} \right) @>D_3>> \cdots \\
@| @V\beta_1VV @V\beta_2VV @V\beta_3VV\\
\Theta_M @>E_0'>>  \mathscr{H}om_{\mathcal{O}_M}\left(\Theta_{\mathcal{F}_0},  \left( \mathcal{N}_{\mathcal{F}_0}^* \right)^* \right) @>E_1'>> \mathscr{H}om_{\mathcal{O}_M}\left( \bigwedge^2 \Theta_{\mathcal{F}_0}, \left( \mathcal{N}_{\mathcal{F}_0}^*\right)^* \right) @>E_2'>> \mathscr{H}om_{\mathcal{O}_M}\left( \bigwedge^3 \Theta_{\mathcal{F}_0}, \left( \mathcal{N}_{\mathcal{F}_0}^* \right)^* \right) @>E_3'>> \cdots
\end{CD}$
\end{center}}}
where $E_0' = D_0\circ \beta_1$, and  $\beta_r (r\geq 1)$ is induced by the inclusion $\frac{\Theta_M}{\Theta_{\mathcal{F}_0}} \to \left( \mathcal{N}_{\mathcal{F}_0}^* \right)^*$. We show that $\beta_{r+1}\circ D_r = E_r'\circ \beta_r$ for $r\geq 1$. In fact, we show the commutativity locally on $U_i$. For $\phi\in \mathscr{H}om_{\mathcal{O}_M}\left(\bigwedge^r \Theta_{\mathcal{F}_0}, \frac{\Theta_M}{\Theta_{\mathcal{F}_0}} \right)$, and $T_1,..., T_{r+1}\in \Theta_{\mathcal{F}_0}$, we have
\begin{align*}
&(\beta_{r+1}\circ D_r)(\phi)(T_1,..., T_{r+1})(w_i^\gamma)\\
&= i_{\sum_{\alpha=1}^{r+1} (-1)^{\alpha+1} \left[ T_\alpha, \phi\left(T_1,..., \hat{T}_\alpha ,..., T_{r+1} \right)    \right] + \sum_{\alpha <\beta} (-1)^{\alpha+\beta} \phi\left(\left[ T_\alpha, T_\beta \right] , T_1,..., \hat{T}_\alpha, ..., \hat{T}_\beta,..., T_{r+1}            \right)  }\left(w_i^\gamma \right)
\end{align*}
On the other hand, we note that $\beta_r(\phi)(T_1,..., T_r)(w_i^\gamma)= i_{\phi(T_1,..., T_r)}(w_i^\gamma) $
{\small{\begin{align*}
&(E_r'\circ \beta_r)(\phi)(T_1,..., T_{r+1})(w_i^\gamma)\\
&= \sum_{\alpha=1}^{r+1} (-1)^{\alpha+1}\left[ T_\alpha, \beta_r(\phi)\left(T_1,..., \hat{T}_\alpha ,..., T_{r+1} \right) (w_i^\gamma)\right] - \beta_r(\phi)(T_1,..., \hat{T}_\alpha,..., T_{r+1})\left( \mathcal{L}_{T_\alpha}\left(w_i^\gamma \right)     \right)  \\
&+ \sum_{\alpha <\beta} (-1)^{\alpha+\beta} \beta_r(\phi)\left(\left[ T_\alpha, T_\beta \right] , T_1,..., \hat{T}_\alpha, ..., \hat{T}_\beta,..., T_{r+1}            \right)(w_i^\gamma)\\
&= \sum_{\alpha=1}^{r+1}(-1)^{\alpha+1} \mathcal{L}_{ T_\alpha} i_{\phi\left(T_1\wedge \cdots \wedge \hat{T}_\alpha \wedge \cdots \wedge T_{r+1} \right)}(w_i^\gamma) -i_{\phi\left( T_1\wedge \cdots \wedge \hat{T}_\alpha \wedge \cdots \wedge T_{r+1} \right)}\left(\mathcal{L}_{T_\alpha}(w_i^\gamma) \right) + \sum_{\alpha < \beta}(-1)^{\alpha+\beta}i_{ \phi\left( \left[ T_\alpha, T_\beta \right], T_1,...,\hat{T}_\alpha, ..., \hat{T}_\beta,... ,T_{r+1} \right) }\left( w_i^\gamma \right)
\end{align*}}}
Hence the diagram commutes. 

We show that we have an isomorphism
{\small{\begin{equation}\label{ad5}
\begin{CD}
\Theta_M @>E_0>> \mathscr{H}om_{\mathcal{O}_M}\left( \mathcal{N}_{\mathcal{F}_0}^*, \frac{\Omega_M^1}{\mathcal{N}_{\mathcal{F}_0}^*} \right) @>E_1>>  \mathscr{H}om_{\mathcal{O}_M}\left( \mathcal{N}_{\mathcal{F}_0}^*, \tilde{\mathcal{S}}^2 \right) @>E_2>> \mathscr{H}om_{\mathcal{O}_M}\left( \mathcal{N}_{\mathcal{F}_0}^*, \tilde{\mathcal{S}}^3 \right) @>E_3>> \cdots \\
@| @| @V\alpha_2VV @V\alpha_3 VV\\
\Theta_M @>E_0>> \mathscr{H}om_{\mathcal{O}_M}\left( \mathcal{N}_{\mathcal{F}_0}^*, \frac{\Omega_M^1}{\mathcal{N}_{\mathcal{F}_0}^*} \right) @>\tilde{E}_1>> \mathscr{H}om_{\mathcal{O}_M}\left( \bigwedge^2 \Theta_{\mathcal{F}_0}, \left(\mathcal{N}_{\mathcal{F}_0}^* \right)^* \right) @>E_2' >> \mathscr{H}om_{\mathcal{O}_M}\left(\bigwedge^3 \Theta_{\mathcal{F}_0} \left( \mathcal{N}_{\mathcal{F}_0}^* \right)^* \right) @>E_3'>> \cdots
\end{CD}
\end{equation}}}
Here $\tilde{E}_1$ is defined by the composition of $\mathscr{H}om_{\mathcal{O}_M}\left( \mathcal{N}_{\mathcal{F}_0}^*, \frac{\Omega_M^1}{\mathcal{N}_{\mathcal{F}_0}^*} \right) \xrightarrow{\gamma_1} \mathscr{H}om_{\mathcal{O}_M}\left( \mathcal{N}_{\mathcal{F}_0}^*, \Theta_{\mathcal{F}_0}^* \right) \xrightarrow{E_1'} \mathscr{H}om_{\mathcal{O}_M}\left( \bigwedge^2 \Theta_{\mathcal{F}_0}, \left(\mathcal{N}_{\mathcal{F}_0}^* \right)^* \right)$ where $\gamma_1$ is induced by the inclusion $\frac{\Omega_M^1}{\mathcal{N}_{\mathcal{F}_0}^*}\to \Theta_{\mathcal{F}_0}^*$. Since $(\gamma_1\circ E_0)(X)(T)(w_i^\alpha)= i_{E_0(X)(T)}(w_i^\alpha)= i_{[X, T]}(w_i^\alpha)= E_0'(X)(T)(w_i^\alpha)$ for $T\in \Theta_{\mathcal{F}_0}$, we have $\tilde{E}_1\circ E_0=0$, and $E_2'\circ \tilde{E}_1=0$.

We define $\alpha_r:\mathscr{H}om_{\mathcal{O}_M}\left( \mathcal{N}_{\mathcal{F}_0}^*, \tilde{\mathcal{S}}^r \right)\to \mathscr{H}om_{\mathcal{O}_M}\left( \bigwedge^r \Theta_{\mathcal{F}_0}, \left( \mathcal{N}_{\mathcal{F}_0}^* \right)^* \right)$ on $U_i$ locally in the following way:  first we note that for $T_1,..., T_r\in \Gamma\left( U_i, \Theta_{\mathcal{F}_0} \right)$, let us consider 
\begin{align*}
 i_{T_1\wedge \cdots \wedge T_r}\left( \phi\left(w_i^\alpha\right) \right) 
\end{align*}
For $x\in U_i-S$, choose $U_x\subset U_i-S$ such that $\phi(w_i^\alpha)= w_i^1\wedge \cdots \wedge w_i^q \wedge A_{i_x}^\alpha$ on $U_{x}$ for some $A_{i_x}\in \Gamma\left( U_x , \bigwedge^r \Omega_M^1 \right)$. Then
\begin{align*}
i_{T_1\wedge \cdots \wedge T_r}\left( \phi\left( w_i^\alpha \right) \right) = i_{T_1\wedge \cdots \wedge T_r}(A_{i_x}) w_i^1\wedge \cdots \wedge w_i^q
\end{align*}
For $y\in U_i-S$ such that $U_x\cap U_y\ne \emptyset$, we have $\left( i_{T_1\wedge \cdots \wedge T_r}\left(A_{i_x}\right)  -  i_{T_1\wedge \cdots \wedge T_r} \left( A_{i_y}^\alpha \right) \right)w_i^1\wedge \cdots \wedge w_i^q=0$ on $U_x\cap U_y$. This implies that $i_{T_1\wedge \cdots \wedge T_r}(A_{i_x}^\alpha)= i_{T_1\wedge \cdots \wedge T_r}(A_{i_y}^\alpha)$. Then $\{i_{T_1\wedge \cdots \wedge T_r}(A_{i_x}^\alpha)|x\in U_i -S \}$ glues together to define a $f_i^\alpha \in \Gamma(U_i-S, \mathcal{O}_M)$, and so $f_i^\alpha \in \Gamma(U_i,\mathcal{O}_M)$. This implies that $i_{T_1\wedge \cdots \wedge T_r}(\phi(w_i^\alpha)) = f_i^\alpha w_i^1\wedge \cdots \wedge w_i^q$. Then we define $\alpha_r(\phi)(T_1,..., T_r)(w_i^\alpha)=f_i^\alpha$. Then we linearly extends to $\Gamma\left( U_i, \mathcal{N}_{\mathcal{F}_0}^* \right)$. This is well-defined since we have (recall that $H_{ij}:= \det\left( h_{ij}^{\alpha\beta} \right)$)
\begin{align*}
i_{T_1\wedge \cdots \wedge T_r}\left( \phi(w_i^\alpha) \right)= H_{ij} \sum_{\beta=1}^q h_{ij}^{\alpha\beta} i_{T_1\wedge \cdots \wedge T_r}\left( \phi\left(w_j^\beta\right) \right) = H_{ij}\sum_{\beta=1}^q h_{ij}^{\alpha\beta} f_j^\beta w_j^1\wedge \cdots \wedge w_j^q= \sum_{\beta=1}^q h_{ij}^{\alpha\beta} f_j^\beta w_i^1\wedge \cdots \wedge w_i^q
\end{align*}
Hence we have $f_i^\alpha= \sum_{\beta=1}^q h_{ij}^{\alpha\beta} f_j^\beta$. In other words, $\alpha_r(\phi)(T_1,..., T_r)(w_i^\alpha)= \sum_{\beta=1}^q h_{ij}^{\alpha\beta} \alpha_r(\phi)(T_1,..., T_r)\left(w_j^\beta \right)$.

We show that $\tilde{E}_1=\alpha_2\circ E_1$. For $\phi\in \mathscr{H}om_{\mathcal{O}_M}\left(\mathcal{N}_{\mathcal{F}_0}^*, \frac{\Omega_M^1}{\mathcal{N}_{\mathcal{F}_0}^*} \right)$.
\begin{align*}
&\tilde{E}_1(\phi)(T_1, T_2)(w_i^\alpha)=( E_1'\circ \gamma_1)(\phi)(T_1,T_2)(w_i^\alpha)\\
&=\left[ T_1, \gamma_1(\phi)(T_2)(w_i^\alpha)\right] - \gamma_1(\phi)(T_2)\left( \mathcal{L}_{T_1}(w_i^\alpha) \right) - \left[ T_2, \gamma_1(\phi)(T_1)(w_i^\alpha) \right] + \gamma_1(\phi)(T_1)\left( \mathcal{L}_{T_2}(w_i^\alpha) \right) - \gamma_1(\phi)\left(\left[ T_1, T_2\right] \right)(w_i^\gamma)\\
&= \left[ T_1, i_{T_2} \left( \phi(w_i^\alpha)\right)\right] - i_{T_2}\left( \phi(\mathcal{L}_{T_1}(w_i^\alpha)) \right) - \left[ T_2, i_{T_1}( \phi( w_i^\alpha)) \right] + i_{T_1}\left( \phi(\mathcal{L}_{T_2}(w_i^\alpha)) \right) - i_{\left[ T_1, T_2\right] }( \phi ( w_i^\alpha ) )
\end{align*}
On the other hand, let us compute $\left(\alpha_2\circ E_1\right)(\phi)(T_1,T_2)(w_i^\alpha) = i_{T_1\wedge T_2}\left(E_1(\phi)(w_i^\alpha) \right)$. For $x\in U_i - S$, choose $U_x \subset U_i - S$ such that $dw_i^\alpha= \sum_{\beta=1}^q a_{i_x}^{\alpha\beta} \wedge w_i^\beta$ for $a_{i_x}^{\alpha\beta} \in \Gamma\left( U_x, \Omega_M^1 \right)$. We note that $\mathcal{L}_T(w_i^\alpha)=i_T \left(dw_i^\alpha \right)=\sum_{\beta=1}^q i_{T}\left( a_{i_x}^{\alpha\beta}\right)w_i^\beta$ for $T\in \Theta_{\mathcal{F}_0}$, so  that $\phi\left( \mathcal{L}_T(w_i^\alpha) \right)= \sum_{\beta=1}^q i_T\left(a_{i_x}^{\alpha\beta} \right) \phi\left(w_i^\beta\right)$. Then on $U_x$
{\small{\begin{align*}
&i_{T_1\wedge T_2}\left(E_1(\phi)(w_i^\alpha)  \right) = i_{T_1\wedge T_2}\left( w_i^1\wedge \cdots \wedge w_i^1 \wedge \left( d\phi(w_i^\alpha) - \sum_{\beta=1}^q a_{i_x}^{\alpha\beta} \wedge \phi \left(w_i^\beta\right)    \right) \right)\\
&= w_i^1\wedge \cdots \wedge w_i^q \wedge \left( i_{T_1\wedge T_2}\left(\left(d\phi(w_i^\alpha) \right) - \sum_{\beta=1}^q a_{i_x}^{\alpha\beta} \wedge \phi\left(w_i^\beta\right) \right) \right)\\
&= w_i^1\wedge \cdots \wedge w_i^q \wedge \left( i_{T_2} \left( \mathcal{L}_{T_1}- d i_{T_1} \right) \left( \phi\left( w_i^\alpha \right)\right)  \right) - w_i^1\wedge \cdots \wedge w_i^q \wedge \left( \sum_{\beta=1}^q i_{T_1}\left(a_{i_x}^{\alpha\beta} \right) i_{T_2}\left(\phi\left(w_i^\beta \right)\right) - i_{T_2}\left(a_{i_x}^{\alpha\beta}\right) i_{T_1}\left( \phi\left(w_i^\beta \right) \right)   \right)\\
&= w_i^1\wedge \cdots \wedge w_i^q\cdot \left(  \mathcal{L}_{T_1}i_{T_2}\left(\phi(w_i^\alpha)\right) - i_{[T_1,T_2]}\left(\phi(w_i^\alpha) \right)    - \mathcal{L}_{T_2}i_{T_1}\left(\phi(w_i^\alpha) \right) - i_{T_2} \phi\left(\mathcal{L}_{T_1}(w_i^\alpha) \right) + i_{T_1} \phi\left( \mathcal{L}_{T_2}(w_i^\alpha) \right)          \right)
\end{align*}}}
Since this holds for any $x\in U_i-S$, this implies $\tilde{E}_1= \alpha_2\circ E_1$.

Next we show that the following diagram commutes
\begin{center}
$\begin{CD}
\mathscr{H}om_{\mathcal{O}_M}\left( \mathcal{N}_{\mathcal{F}_0}^*, \tilde{\mathcal{S}}^r \right) @>E_r>> \mathscr{H}om_{\mathcal{O}_M} \left( \mathcal{N}_{\mathcal{F}_0}^*, \tilde{\mathcal{S}}^{r+1}   \right) \\
@V\alpha_r VV @VV\alpha_{r+1}V\\
\mathscr{H}om_{\mathcal{O}_M}\left( \bigwedge^r \Theta_{\mathcal{F}_0}, \left( \mathcal{N}_{\mathcal{F}_0}^* \right)^* \right) @>E_r'>> \mathscr{H}om_{\mathcal{O}_M}\left( \bigwedge^{r+1} \Theta_{\mathcal{F}_0}, \left( \mathcal{N}_{\mathcal{F}_0}^* \right)^* \right) 
\end{CD}$
\end{center}
Let us check the commutativity locally on $U_i$. Let us compute $(E_r\circ \alpha_r)(\phi)(T_1,..., T_r)\left(w_i^\gamma\right)$.

{\tiny{\begin{align*}
 \sum_{\alpha=1}^{r+1} (-1)^{\alpha+1}\left(\left[ T_\alpha, \alpha_r(\phi)\left(T_1,..., \hat{T}_\alpha ,..., T_{r+1} \right) (w_i^\gamma)\right] - \alpha_r(\phi)(T_1,..., \hat{T}_\alpha,..., T_{r+1})\left( \mathcal{L}_{T_\alpha}\left(w_i^\gamma \right)     \right) \right) + \sum_{\alpha <\beta} (-1)^{\alpha+\beta} \alpha_r(\phi)\left(\left[ T_\alpha, T_\beta \right] , T_1,..., \hat{T}_\alpha, ..., \hat{T}_\beta,..., T_{r+1}            \right)(w_i^\gamma)
\end{align*}}}
For $x\in U_i-S$, choose $U_x \subset U_i - S$ such that $dw_i^\alpha= \sum_{\beta=1}^q a_{i_x}^{\alpha\beta}\wedge w_i^\beta$ for some $a_{i_x}^{\alpha\beta}\in \Gamma\left(U_x, \bigwedge^r \Omega_M^1 \right)$ on $U_x$ and $\phi(w_i^\gamma)= w_i^1\wedge \cdots \wedge w_i^q \wedge A_{i_x}^\gamma$ on $U_x$, and we note that 
\begin{align*}
\phi\left( \mathcal{L}_{T_\alpha}(w_i^\gamma) \right) = \phi\left( i_{T_\alpha}\left(d w_i^\gamma \right) \right) = \phi\left( \sum_{\eta=1}^q i_{T_\alpha}\left(a_{i_x}^{\gamma\eta}\right) w_i^\eta   \right) = w_i^1\wedge \cdots \wedge w_i^q \wedge \left( \sum_{\eta=1}^q i_{T_\alpha}\left(   a_{i_x}^{\gamma\eta} \right)  A_{i_x}^\eta     \right)
\end{align*}

Then $(E_r\circ \alpha_r)(\phi)(T_1,..., T_{r+1})(w_i^\gamma)$ on $U_x$ is
\begin{align*}
&(E_r\circ \alpha_r)(\phi)(T_1,..., T_{r+1})(w_i^\gamma)\\
&= \sum_{\alpha=1}^{r+1} (-1)^{\alpha+1} \left[ T_\alpha, i_{T_1\wedge \cdots \wedge \hat{T}_\alpha \wedge \cdots \wedge T_{r+1}}(A_{i_x}^\gamma)  \right] - \sum_{\alpha=1}^{r+1} (-1)^{\alpha+1} i_{T_1\wedge \cdots \wedge \hat{T}_\alpha \wedge \cdots \wedge T_{r+1}}\left(\sum_{\eta=1}^q i_{T_\alpha}\left( a_{i_x}^{\gamma \eta} \right) A_{i_x}^\eta  \right) \\
&+\sum_{\alpha < \beta}(-1)^{\alpha+\beta} i_{[T_\alpha, T_\beta]\wedge T_1\wedge \cdots \wedge \hat{T}_\alpha \wedge \cdots \wedge \hat{T}_\beta \wedge \cdots \wedge T_{r+1}}(A_{i_x}^\gamma)
\end{align*}

Let us compute $\alpha_{r+1}\circ E_r(\phi)(T_1,..., T_r)(w_i^\gamma)$ on $U_x \subset U_i- S$. Then we note that
\begin{align*}
E_r(\phi)(w_i^\gamma)= w_i^1\wedge \cdots \wedge w_i^q \wedge \left(  dA_{i_x}^\gamma - \sum_{\eta=1}^q a_{i_x}^{\gamma \eta} \wedge A_{i_x}^\eta   \right)
\end{align*}
Then
{\small{\begin{align*}
&\alpha_{r+1}\circ E_r(\phi)(T_1,..., T_{r+1})(w_i^\gamma) = i_{T_1\wedge \cdots \wedge T_{r+1}}\left(d A_{i_x}^\gamma - \sum_{\eta=1}^q a_{i_x}^{\gamma \eta} \wedge A_{i_x}^\eta \right)\\
&= i_{T_{r+1}} \cdots  i_{T_2}\left(i_{T_1} d \right) \left( A_{i_x}^\gamma  \right) - \sum_{\beta=1}^q i_{T_2\wedge \cdots \wedge T_{r+1}}\left(i_{T_1}\left(  a_{i_x}^{\gamma \eta}\right) A_{i_x}^\eta  - a_{i_x}^{\gamma \eta}  \wedge i_{T_1}\left(A_{i_x}^\eta \right)       \right)\\
&= i_{T_{r+1}}\cdots i_{T_2}\mathcal{L}_{T_1} A_{i_x}^\gamma - i_{T_{r+1}} \cdots i_{T_3}\left(i_{T_2}d \right) i_{T_1}\left(A_{i_x}^\gamma \right) - \sum_{\beta=1}^q i_{T_2\wedge \cdots \wedge T_{r+1}}\left(i_{T_1}\left(  a_{i_x}^{\gamma \eta}\right) A_{i_x}^\eta \right) \\
&+ \sum_{\beta=1}^q i_{T_3}\wedge \cdots \wedge i_{T_{r+1}}\left( i_{T_2}\left(a_{i_x}^{\gamma \eta} \right) i_{T_1}\left( A_{i_x}^\eta \right) - a_{i_x}^{\gamma \eta} \wedge i_{T_1\wedge T_2}\left(A_{i_x}^\eta \right)    \right) \\
&= i_{T_{r+1}} \cdots i_{T_3}\left( \mathcal{L}_{T_1} i_{T_2} - i_{[T_1,T_2]} \right)(A_{i_x}^\gamma) - i_{T_{r+1}} \cdots i_{T_3}\left(\mathcal{L}_{T_2}   \right) i_{T_1}(A_{i_x}^\gamma) + i_{T_{r+1}}\cdots i_{T_4}\left(i_{T_3}d \right) i_{T_2}i_{T_1}(A_{i_x}^\gamma)\\
&- \sum_{\alpha=1}^{r+1} (-1)^{\alpha+1} i_{T_1\wedge \cdots \wedge \hat{T}_\alpha \wedge \cdots \wedge T_{r+1}}\left(\sum_{\eta=1}^q i_{T_\alpha}\left( a_{i_x}^{\gamma \eta} \right) A_{i_x}^\eta  \right)\\
&= \sum_{\alpha=1}^{r+1} (-1)^{\alpha+1} \left[ T_\alpha, i_{T_1\wedge \cdots \wedge \hat{T}_\alpha \wedge \cdots \wedge T_{r+1}}(A_{i_x}^\gamma)  \right] - \sum_{\alpha=1}^{r+1} (-1)^{\alpha+1} i_{T_1\wedge \cdots \wedge \hat{T}_\alpha \wedge \cdots \wedge T_{r+1}}\left(\sum_{\eta=1}^q i_{T_\alpha}\left( a_{i_x}^{\gamma \eta} \right) A_{i_x}^\eta  \right) \\
&+\sum_{\alpha < \beta}(-1)^{\alpha+\beta} i_{[T_\alpha, T_\beta]\wedge T_1\wedge \cdots \wedge \hat{T}_\alpha \wedge \cdots \wedge \hat{T}_\beta \wedge \cdots \wedge T_{r+1}}(A_{i_x}^\gamma)\end{align*}}}
Hence the above diagram commutes on $U_x\subset U_i-S$ for any $x\in U_i$. This implies that it commutes on $U_i-S$, so on $U_i$.

Now we define a complex of sheaves $\mathcal{F}_0^\bullet$ by combining the leaf complex and the dual leaf complex in the following way:
{\small{\begin{equation}\label{ad3}
\begin{CD}
\cdots \\
@AF_3AA \\
\mathscr{H}om_{\mathcal{O}_M}\left(\bigwedge^3 \Theta_{\mathcal{F}_0},  \frac{\Theta_M}{\Theta_{\mathcal{F}_0}} \right) \bigoplus \mathscr{H}om_{\mathcal{O}_M}\left( \mathcal{N}_{\mathcal{F}_0}^*, \tilde{\mathcal{S}}^3 \right) \bigoplus \mathscr{H}om_{\mathcal{O}_M}\left(  \bigwedge^2 \Theta_{\mathcal{F}_0} , \left( \mathcal{N}_{\mathcal{F}_0}^* \right)^* \right) \\
@AF_2AA \\
\mathscr{H}om_{\mathcal{O}_M}\left( \bigwedge^2 \Theta_{\mathcal{F}_0}, \frac{\Theta_M}{\Theta_{\mathcal{F}_0}} \right) \bigoplus \mathscr{H}om_{\mathcal{O}_M}\left( \mathcal{N}_{\mathcal{F}_0}^*, \tilde{\mathcal{S}}^2 \right) \bigoplus \mathscr{H}om_{\mathcal{O}_M}\left( \Theta_{\mathcal{F}_0}, \left(\mathcal{N}_{\mathcal{F}_0}^* \right)^* \right) \\
@AF_1AA\\
\mathscr{H}om_{\mathcal{O}_M}\left(\Theta_{\mathcal{F}_0}, \frac{\Theta_M}{\Theta_{\mathcal{F}_0}}   \right) \bigoplus \mathscr{H}om_{\mathcal{O}_M}\left( \mathcal{N}_{\mathcal{F}_0}^* , \frac{\Omega_M^1}{\mathcal{N}_{\mathcal{F}_0}^*} \right)\\
@AF_0AA\\
\Theta_M
\end{CD}
\end{equation}}}
We define $F_r$ in the following way: we define $F_0: \Theta_M\to \mathscr{H}om_{\mathcal{O}_M}\left(\Theta_{\mathcal{F}_0}, \frac{\Theta_M}{\Theta_{\mathcal{F}_0}}   \right) \bigoplus \mathscr{H}om_{\mathcal{O}_M}\left( \mathcal{N}_{\mathcal{F}_0}^* , \frac{\Omega_M^1}{\mathcal{N}_{\mathcal{F}_0}^*} \right)$ by $F_0(a,b)=(D_0(a), E_0(a))$, and  we define
{\Small{\begin{align*}
F_1: \mathscr{H}om_{\mathcal{O}_M}\left(\Theta_{\mathcal{F}_0}, \frac{\Theta_M}{\Theta_{\mathcal{F}_0}}   \right) \bigoplus \mathscr{H}om_{\mathcal{O}_M}\left( \mathcal{N}_{\mathcal{F}_0}^* , \frac{\Omega_M^1}{\mathcal{N}_{\mathcal{F}_0}^*} \right) &\to \mathscr{H}om_{\mathcal{O}_M}\left( \bigwedge^2 \Theta_{\mathcal{F}_0}, \frac{\Theta_M}{\Theta_{\mathcal{F}_0}} \right) \bigoplus \mathscr{H}om_{\mathcal{O}_M}\left( \mathcal{N}_{\mathcal{F}_0}^*, \tilde{\mathcal{S}}^2 \right) \bigoplus \mathscr{H}om_{\mathcal{O}_M}\left( \Theta_{\mathcal{F}_0}, \left(\mathcal{N}_{\mathcal{F}_0}^* \right)^* \right) \\
(a,b)&\mapsto (D_1(a), E_1(a), \beta_1(a) + \gamma_1(b))
\end{align*}}}
We recall that $\gamma_1$ is induced from the inclusion $\frac{\Omega_M^1}{\mathcal{N}_{\mathcal{F}_0}^*}\to \Theta_{\mathcal{F}_0}^*$. In general, for $r >1$,  we define

{\Small{\begin{align*}
F_r&:\mathscr{H}om_{\mathcal{O}_M}\left(\bigwedge^r \Theta_{\mathcal{F}_0}, \frac{\Theta_M}{\Theta_{\mathcal{F}_0} } \right) \bigoplus \mathscr{H}om_{\mathcal{O}_M}\left(\mathcal{N}_{\mathcal{F}_0}^*, \tilde{\mathcal{S}}^r \right) \bigoplus \mathscr{H}om_{\mathcal{O}_M}\left(\mathcal{N}_{\mathcal{F}_0}^*, \bigwedge^{r-1} \Theta_{\mathcal{F}_0}^* \right) \\
&\to \mathscr{H}om_{\mathcal{O}_M}\left(\bigwedge^{r+1} \Theta_{\mathcal{F}_0} , \frac{\Theta_M}{\Theta_{\mathcal{F}_0} } \right) \bigoplus \mathscr{H}om_{\mathcal{O}_M}\left(\mathcal{N}_{\mathcal{F}_0}^*,  \tilde{\mathcal{S}}^{r+1} \right) \bigoplus \mathscr{H}om_{\mathcal{O}_M}\left(\mathcal{N}_{\mathcal{F}_0}^*,\bigwedge^r  \Theta_{\mathcal{F}_0}^* \right)\\
&(a,b,c )\mapsto \left( D_r(a),E_r(b), \beta_r(a)+ \alpha_r(b)- E_{r-1}'(c) \right)
\end{align*}}}

Then we note that since $\left(\beta_1\circ D_0(X) + \gamma_1\circ E_0(X)\right)(T)(w_i^\alpha)=i_{[X, T]}(w_i^\alpha)+ i_T\mathcal{L}_X(w_i^\alpha)=0$ for $T\in \Theta_{\mathcal{F}_0}$,  
\begin{align*}
F_1F_0(X)=\left(D_1D_0(X), E_1E_0(X), \beta_1 \circ D_0(X)+  \gamma_1\circ E_0(X)\right)=0
\end{align*}
and we note that
\begin{align*}
F_2F_1(a,b)&= F_2\left( D_1(a), E_1(b), \beta_1(a)+ \gamma_1(b)    \right) =\left( D_2D_1(a), E_2E_1(a), \beta_2D_1(a)+ \alpha_2E_1(b) - E_1'\left( \beta_1(a)+ \gamma_1(b)\right)     \right) \\
                  &= \left( 0, 0, \left(\beta_2 D_1- E_1'\beta_1 \right)(a) + \left(\alpha_2 E_1 - E_1'\gamma_1 \right)(b)    \right) = \left(0,0, \left( \alpha_2 E_1- \tilde{E}_1\right)(b) \right) = 0
\end{align*}
and we note that for $r>1$,
\begin{align*}
F_{r+1}F_r(a,b,c)&= F_{r+1}\left( D_r(a), D_r(b), \beta_r(a)+ \alpha_r(b) - E_{r-1}'(c)    \right)\\
&= \left(D_{r+1}D_r(a), E_{r+1}E_r(b),  \beta_{r+1} D_r(a)+ \alpha_{r+1} D_r(b) - E_{r}'\left( \beta_r(a)+ \alpha_r(b)- E_{r-1}'(b) \right)    \right) \\
&= \left( 0, 0, \left( \beta_{r+1} D_r- E_r'\beta_r \right)(a)  + \left( \alpha_{r+1} E_r - E_r' \alpha_r \right)(b)      \right) = 0
\end{align*}
This implies that $\mathcal{F}_0^\bullet$ defines a complex of sheaves.

\subsection{Another type of complex of sheaves controlling simultaneous deformations of  compact foliated complex manifolds defined by both locally free subsheaves of tangent sheaves and cotangent sheaves}\

Let $\left( M, \Theta_{\mathcal{F}_0}, \mathcal{N}_{\mathcal{F}_0}^* \right)$ be a compact foliated complex manifold with both $\Theta_{\mathcal{F}_0}$ and $\mathcal{N}_{\mathcal{F}_0}^*$ locally free. We recall that the dual leaf complex $\mathcal{N}_{\mathcal{F}_0}^{*\bullet}$ is isomorphic to $\mathcal{N}_{\mathcal{F}_0}'^{ *\bullet}$ from $\textnormal{(\ref{ad5})}$. Then by replacing tha part $\mathcal{N}_{\mathcal{F}_0}^{*\bullet}$ in $\textnormal{(\ref{ad3})}$ by $\mathcal{N}_{\mathcal{F}_0}'^{*\bullet}$, we have the following complex of sheaves

\begin{equation}\label{ad10}
\begin{CD}
\cdots \\
@AF_3^\sharp AA \\
\mathscr{H}om_{\mathcal{O}_M}\left( \bigwedge^3 \Theta_{\mathcal{F}_0}, \frac{\Theta_M}{\Theta_{\mathcal{F}_0}} \right) \bigoplus \mathscr{H}om_{\mathcal{O}_M}\left( \mathcal{N}_{\mathcal{F}_0}^*, \bigwedge^3 \Theta_{\mathcal{F}_0}^* \right) \bigoplus \mathscr{H}om_{\mathcal{O}_M}\left( \mathcal{N}_{\mathcal{F}_0}^*, \bigwedge^2 \Theta_{\mathcal{F}_0}^* \right) \\
@AF_2^\sharp AA \\
\mathscr{H}om_{\mathcal{O}_M}\left( \bigwedge^2 \Theta_{\mathcal{F}_0}, \frac{\Theta_M}{\Theta_{\mathcal{F}_0}} \right) \bigoplus \mathscr{H}om_{\mathcal{O}_M}\left( \mathcal{N}_{\mathcal{F}_0}^*, \bigwedge^2 \Theta_{\mathcal{F}_0}^* \right) \bigoplus \mathscr{H}om_{\mathcal{O}_M}\left( \mathcal{N}_{\mathcal{F}_0}^*, \Theta_{\mathcal{F}_0}^* \right) \\
@AF_1^\sharp AA \\
\mathscr{H}om_{\mathcal{O}_M}\left( \Theta_{\mathcal{F}_0}, \frac{\Theta_M}{\Theta_{\mathcal{F}_0}} \right) \bigoplus \mathscr{H}om_{\mathcal{O}_M}\left( \mathcal{N}_{\mathcal{F}_0}^*, \frac{\Omega_M^1}{\mathcal{N}_{\mathcal{F}_0}^*} \right) \\
@AF_0AA \\
\Theta_M
\end{CD}
\end{equation}
We define $F_0(a,b)=\left(D_0(a) , E_0(b) \right)$, and we define
{\Small{\begin{align*}
F_1^\sharp: \mathscr{H}om_{\mathcal{O}_M}\left(\Theta_{\mathcal{F}_0}, \frac{\Theta_M}{\Theta_{\mathcal{F}_0}}   \right) \bigoplus \mathscr{H}om_{\mathcal{O}_M}\left( \mathcal{N}_{\mathcal{F}_0}^* , \frac{\Omega_M^1}{\mathcal{N}_{\mathcal{F}_0}^*} \right) &\to \mathscr{H}om_{\mathcal{O}_M}\left( \bigwedge^2 \Theta_{\mathcal{F}_0}, \frac{\Theta_M}{\Theta_{\mathcal{F}_0}} \right) \bigoplus \mathscr{H}om_{\mathcal{O}_M}\left( \mathcal{N}_{\mathcal{F}_0}^*, \bigwedge^2 \Theta_{\mathcal{F}_0}^* \right) \bigoplus \mathscr{H}om_{\mathcal{O}_M}\left( \mathcal{N}_{\mathcal{F}_0}^*, \Theta_{\mathcal{F}_0}^* \right) \\
(a,b)&\mapsto (D_1(a), \tilde{E}_1(a), \beta_1(a) + \gamma_1(b))
\end{align*}}}
In general $r>1$, we define
{\Small{\begin{align*}
F_r^\sharp&:\mathscr{H}om_{\mathcal{O}_M}\left(\bigwedge^r \Theta_{\mathcal{F}_0}, \frac{\Theta_M}{\Theta_{\mathcal{F}_0} } \right) \bigoplus \mathscr{H}om_{\mathcal{O}_M}\left(\mathcal{N}_{\mathcal{F}_0}^*, \bigwedge^r \Theta_{\mathcal{F}_0}^* \right) \bigoplus \mathscr{H}om_{\mathcal{O}_M}\left(\mathcal{N}_{\mathcal{F}_0}^*, \bigwedge^{r-1} \Theta_{\mathcal{F}_0}^* \right) \\
&\to \mathscr{H}om_{\mathcal{O}_M}\left(\bigwedge^{r+1} \Theta_{\mathcal{F}_0} , \frac{\Theta_M}{\Theta_{\mathcal{F}_0} } \right) \bigoplus \mathscr{H}om_{\mathcal{O}_M}\left(\mathcal{N}_{\mathcal{F}_0}^*,  \bigwedge^{r+1} \Theta_{\mathcal{F}_0}^* \right) \bigoplus \mathscr{H}om_{\mathcal{O}_M}\left(\mathcal{N}_{\mathcal{F}_0}^*,\bigwedge^r  \Theta_{\mathcal{F}_0}^* \right)\\
&(a,b,c )\mapsto \left( D_r(a),E_r'(b), \beta_r(a)+ b- E_{r-1}'(c) \right)
\end{align*}}}

\section{Formulas on bracket and Lie derivative}\label{fb1}

We collect definitions and formulas on bracket and Lie derivative which are used in the main body of Part I and Part III: unfoldings. Let $M$ be a compact complex manifold and $z=(z_1,..., z_n)$ be a local coordinate on $M$. For the definition on the bracket on $\bigoplus_p A^{0,p}\left( M, \Theta_M \right)$ and formulas, where $A^{0,p}\left( M, \Theta_M \right)$ is the global section of $\mathcal{A}^{0,p}\left(\Theta_M\right)$, which is the sheaf of germs of $C^\infty(0,p)$-forms with coefficients in $\Theta_M$, we refer to \cite{Kod05} p.266, or more generally, for the definition on the bracket on $\bigoplus_{p,q} A^{0,p}\left( M, \bigwedge^q \Theta_M \right)$ and formulas, we refer to the author's thesis.

Let $T$ be a $C^\infty$-vector of the form $T= \sum_{\alpha=1}^n T^\alpha(z)\frac{\partial}{\partial z_\alpha}$, i.e. a section of $\mathcal{A}^{0,0}\left( \Theta_M \right)$. Then we define the Lie derivative $\mathcal{L}_T$ of a $C^\infty$ $1$-form $w$ of the form $w= \sum_{\beta=1}^n w_\beta dz_\beta + \sum_{\beta=1}^n \overline{w}_\beta d \bar{z}_\beta$, i.e. a section of $\mathcal{A}^{1,0}\bigoplus \mathcal{A}^{0,1}$, where $\mathcal{A}^{p,q}$ is the sheaf of germs of $C^\infty$ sections of $\bigwedge^p \Omega_M^1 \otimes \bigwedge^q \overline{\Omega}_M^1$ in the following way \footnote{In general, the original Lie derivative $\mathcal{L}_T$ for $C^\infty$-vector $T\in \mathcal{A}^{0,0}\left( \Theta_M \bigoplus \overline{\Theta_M} \right)$ is defined by $\mathcal{L}_T= d i_{T} + i_T d = \left( \partial i_T + i_T \partial \right) + \left(\bar{\partial} i_T + i_T \bar{\partial}\right)$. In our case, a $C^\infty$-vector $T$ is always in $\mathcal{A}^{0,0}\left( \Theta_M \right)$ and by abuse of notation we define the Lie derivative $\mathcal{L}_T:= \partial i_T + i_T \partial$. If $w$ is a $C^\infty(p,q)$-form, then $\mathcal{L}_T\left( w\right)$ in our case is the $(p,q)$-component of the original Lie derivative. In other words, we ignore $(p-1, q+1)$-component. Hence it is well-defined.}:
\begin{align*}
\mathcal{L}_T\left( w \right) = \partial i_T\left( w \right) + i_T \partial (w)
\end{align*}
where $\partial w= \sum_{\beta, \gamma=1}^n \frac{\partial w_\beta}{\partial z_\gamma} dz_\gamma \wedge dz_\beta+ \sum_{\beta, \gamma=1}^n \frac{\partial \overline{w}_\beta}{\partial z_\gamma} d z_\gamma \wedge d \bar{z}_\beta$. More generally, for a $C^\infty (0,p)$-vector of the form $T=\sum_{\gamma_1,..., \gamma_p} \sum_{\alpha=1}^n d \bar{z}_{\gamma_1} \wedge \cdots \wedge d\bar{z}_{\gamma_1} T_{\gamma_1,..., \gamma_p}^\alpha \frac{\partial}{\partial z_\alpha}$, i.e. a section of $\mathcal{A}^{0,p}\left(\Theta_M \right)$, and a $C^\infty(p,q)$-form $w$, i.e. a section of $\mathcal{A}^{p,q}$, we define
\begin{align*}
\mathcal{L}_T\left( w \right )&= \sum_{\gamma_1,..., \gamma_p} d \bar{z}_{\gamma_1} \wedge \cdots \wedge d \bar{z}_{\gamma_p} \wedge  \mathcal{L}_{\sum_{\alpha=1}^n T_{\gamma_1,..., \gamma_p}^\alpha \frac{\partial}{\partial z_\alpha}}\left( w \right) \\
& = \sum_{\gamma_1,..., \gamma_p} d\bar{z}_{\gamma_1} \wedge \cdots \wedge d \bar{z}_{\gamma_p}\wedge \left( \partial i_{\sum_{\alpha=1}^n T_{\gamma_1,..., \gamma_p}^\alpha \frac{\partial}{\partial z_\alpha}}\left(  w\right)  + i_{ \sum_{\alpha=1}^n T_{\gamma_1,..., \gamma_p}^\alpha \frac{\partial}{\partial z_\alpha} }\partial \left( w \right)      \right)
\end{align*}

Let us consider a $C^\infty(0,1)$-vector $\varphi=\sum_{\lambda, v=1}^n \varphi_v^\lambda d\bar{z}_v\frac{\partial}{\partial z_\lambda}$ and a a $q$-form $w_1 \wedge \cdots \wedge w_q$. Then we have
\begin{align*}
\mathcal{L}_\varphi \left(w_1\wedge \cdots \wedge w_q \right) &=\sum_{\lambda, v=1}^n d\bar{z}_v\wedge \mathcal{L}_{\varphi_v^\lambda\frac{\partial}{\partial z_\lambda}}\left(w_1\wedge \cdots \wedge w_1^q \right)=\sum_{\lambda, v=1}^n\sum_{\alpha=1}^q d\bar{z}_v\wedge w_1\wedge \cdots \wedge \left(\mathcal{L}_{\varphi_v^\lambda \frac{\partial}{\partial z_\lambda}} \left(w_\alpha \right) \right) \wedge \cdots \wedge w_q\\
&=\sum_{\alpha=1}^q (-1)^{\alpha-1} w_1\wedge \cdots \wedge\left(  \sum_{\lambda,v=1}^n d\bar{z}_v \wedge \mathcal{L}_{\varphi_v^\lambda \frac{\partial}{\partial z_\lambda}} \left(w_\alpha \right) \right) \wedge \cdots \wedge w_q\\
&=\sum_{\alpha=1}^q  (-1)^{\alpha-1} w_1\wedge \cdots \wedge \mathcal{L}_\varphi \left(w_\alpha \right)\wedge \cdots \wedge w_q
\end{align*}

Let us consider a $C^\infty(0,1)$-vector $\varphi=\sum_{\lambda, v=1}^n \varphi_v^\lambda d\bar{z}_v \frac{\partial}{\partial z_\lambda}$ and $1$-form $w=\sum_{\alpha=1}^n w_\alpha dz_\alpha$. Then we claim that
\begin{align}\label{bf1}
\bar{\partial} \mathcal{L}_{\varphi}\left(w\right)= \mathcal{L}_{\bar{\partial} \varphi} \left(w\right)- \mathcal{L}_{\varphi}\left(\bar{\partial } w \right)
\end{align}
In fact, we compute
\begin{align}\label{bf2}
\bar{\partial} \mathcal{L}_\varphi \left(w\right)&=\bar{\partial}\left( \sum_{\lambda, v,\alpha=1}^n d\bar{z}_v \wedge \mathcal{L}_{\varphi_v^\lambda\frac{\partial}{\partial z_\lambda}}(w_\alpha dz_\alpha) \right)=\bar{\partial}\left( \sum_{\lambda, v,\alpha=1}^n d\bar{z}_v \wedge\left(\varphi_v^\lambda \frac{\partial w_\alpha}{\partial z_\lambda} dz_\alpha \right) + \sum_{\lambda, v,\alpha=1}^n d\bar{z}_v \wedge w_\alpha \frac{\partial \varphi_v^\alpha }{\partial z_\lambda} dz_\lambda \right)\\
&= \bar{\partial}\left( \sum_{\lambda, v,\alpha=1}^n \left( \varphi_v^\lambda \frac{\partial w_\alpha}{\partial z_\lambda} + w_\lambda \frac{\partial \varphi_v^\lambda}{\partial z_\alpha} \right) d\bar{z}_v \wedge dz_\alpha \right) \notag
\end{align}
We note that $\bar{\partial}\varphi=\sum_{\beta=1}^n \frac{\partial \varphi_v^\lambda}{\partial \bar{z}_\beta} d\bar{z}_\beta \wedge d\bar{z}_v \frac{\partial}{\partial z_\lambda}$. Then we have
\begin{align}\label{bf3}
\mathcal{L}_{\bar{\partial} \varphi}\left( w\right)&=\sum_{\beta, v, \lambda, \alpha=1}^n d\bar{z}_\beta \wedge d\bar{z}_v\wedge \mathcal{L}_{\frac{\partial \varphi_v^\lambda}{\partial \bar{z}_\beta} \frac{\partial}{\partial z_\lambda}} \left(w_\alpha dz_\alpha \right)=\sum_{\beta, v,\lambda, \alpha=1}^n d\bar{z}_\beta \wedge d\bar{z}_v\wedge \left(\frac{\partial \varphi_v^\lambda}{\partial \bar{z}_\beta}\frac{\partial w_\alpha}{\partial z_\lambda}dz_\alpha + w_\alpha \frac{\partial^2 \varphi_v^\alpha}{\partial{z_\lambda}\partial \bar{z}_\beta}dz_\lambda \right)\\
 &=\sum_{\beta, v, \lambda, \alpha=1}^n d\bar{z}_\beta\wedge d\bar{z}_v \wedge \left( \frac{\partial \varphi_v^\lambda}{\partial \bar{z}_\beta}\frac{\partial w_\alpha}{\partial z_\lambda}dz_\alpha + w_\lambda \frac{\partial^2 \varphi_v^\lambda}{\partial{z_\alpha}\partial \bar{z}_\beta}dz_\alpha \right)\notag
\end{align}
On the other hand, we note that $\bar{\partial}w= \sum_{\alpha,\beta=1}^n \frac{\partial w_\alpha}{\partial \bar{z}_\beta}d\bar{z}_\beta \wedge dz_\alpha$. Then we have
\begin{align} \label{bf4}
\mathcal{L}_{\varphi}\left( \bar{\partial} w \right) &=\sum_{\lambda, v=1}^n d\bar{z}_v \wedge \mathcal{L}_{\varphi_v^\lambda \frac{\partial}{\partial z_\lambda}}\left( \sum_{\alpha, \beta=1}^n \frac{\partial w_\alpha}{\partial \bar{z}_\beta}d\bar{z}_\beta \wedge dz_\alpha \right)\\
& =\sum_{\lambda, v=1}^n d\bar{z}_v\wedge \left( \sum_{\alpha, \beta=1}^n \mathcal{L}_{\varphi_v^\lambda \frac{\partial}{\partial z_\lambda}}    \left(\frac{\partial w_\alpha}{\partial \bar{z}_\beta} d\bar{z}_\beta   \right) \wedge dz_\alpha + \sum_{\alpha, \beta=1}^n \frac{\partial w_\alpha}{\partial \bar{z}_\beta} d\bar{z}_\beta \wedge \mathcal{L}_{\varphi_v^\lambda \frac{\partial}{\partial z_\lambda}}(dz_\alpha)\right)\notag\\
&=\sum_{\lambda, v, \alpha, \beta =1}^n d\bar{z}_v \wedge\left(\varphi_v^\lambda \frac{\partial^2 w_\alpha}{\partial z_\lambda \partial \bar{z}_\beta}d\bar{z}_\beta \wedge dz_\alpha +\sum_{\alpha,\beta=1}^n \frac{\partial w_\alpha}{\partial \bar{z}_\beta} d\bar{z}_\beta \wedge \frac{\partial \varphi_v^\alpha}{\partial z_\lambda} dz_\lambda \right)\notag\\
&=-\sum_{\lambda, v,\alpha, \beta=1}^n d\bar{z}_\beta \wedge d\bar{z}_v\wedge \left( \varphi_v^\lambda \frac{\partial^2 w_\alpha}{\partial z_\lambda \partial \bar{z}_\beta} dz_\alpha + \frac{\partial w_\lambda}{\partial \bar{z}_\beta}\frac{\partial \varphi_v^\lambda}{\partial z_\alpha}  dz_\alpha \right) \notag
\end{align}
Then $\textnormal{(\ref{bf2}),(\ref{bf3})}$ and $\textnormal{(\ref{bf4})}$ imply $\textnormal{(\ref{bf1})}$.

Let us consider a $C^\infty(0,1)$-vector $\varphi= \sum_{\lambda, v=1}^n \varphi_v^\lambda d\bar{z}_v \frac{\partial}{\partial z_\lambda}$ and a $C^\infty (0,1)$-vector $\psi= \sum_{\eta, \delta =1}^n \psi_\delta^\eta d \bar{z}_\delta \frac{\partial}{\partial z_\eta}$ and $C^\infty$ $1$-form $w= \sum_{\alpha=1}^n w_\alpha dz_\alpha$. Then we claim that
\begin{align}\label{bf7}
\mathcal{L}_{\left[ \varphi, \psi \right]}(w)= \mathcal{L}_\varphi \mathcal{L}_\psi(w) + \mathcal{L}_\psi \mathcal{L}_\varphi (w)
\end{align}
In fact, we note that we have
\begin{align*}
\sum_{v,\lambda, \eta, \delta=1}^n \left[\varphi_v^\lambda d\bar{z}_v\frac{\partial}{\partial z_\lambda},  \psi_\eta^\delta d\bar{z}_\eta\frac{\partial}{\partial z_\delta} \right]=\sum_{v, \lambda, \eta, \delta=1}^n (-1)^{1\cdot 2} d\bar{z}_v\wedge d\bar{z}_\eta\left[ \varphi_v^\lambda \frac{\partial}{\partial z_\lambda}, \psi_\eta^\delta \frac{\partial}{\partial z_\delta} \right]
\end{align*}
Then we have
\begin{align}\label{bf5}
\mathcal{L}_{\left[\varphi, \psi\right]}(w)= \sum_{v, \lambda, \eta, \delta=1}^n d\bar{z}_v\wedge d\bar{z}_\eta \wedge \mathcal{L}_{\left[ \varphi_v^\lambda\frac{\partial}{\partial z_\lambda}, \psi_\eta^\delta \frac{\partial}{\partial z_\delta}\right]}(w) =\sum_{v, \lambda, \eta, \delta=1}^n d\bar{z}_v \wedge d\bar{z}_\eta\wedge \left(L_{\varphi_v^\lambda \frac{\partial}{\partial z_\lambda}} \mathcal{L}_{\psi_\eta^\delta\frac{\partial}{\partial z_\delta}}(w) - \mathcal{L}_{\psi_\eta^\delta\frac{\partial}{\partial z_\delta}} \mathcal{L}_{\varphi_v^\lambda \frac{\partial}{\partial z_\lambda}}(w) \right)
\end{align}
On the other hand, we have
{\Small{\begin{align}\label{bf6}
\mathcal{L}_{\varphi}\mathcal{L}_\psi(w)= \sum_{v, \lambda, \eta, \delta=1}^n \mathcal{L}_{\varphi_v^\lambda d\bar{z}_v \frac{\partial}{\partial z_\lambda}} \mathcal{L}_{\psi_\eta^\delta d\bar{z}_\eta \frac{\partial}{\partial z_\delta}}(w)= \sum_{v,\lambda, \eta, \delta=1}^n d\bar{z}_v \wedge \left( \mathcal{L}_{\varphi_v^\lambda \frac{\partial}{\partial z_\lambda}} \left(d\bar{z}_\eta\wedge \mathcal{L}_{\psi_\eta^\delta\frac{\partial}{\partial z_\delta}}(w) \right)\right)=\sum_{v, \lambda, \eta, \delta=1}^n d\bar{z}_v \wedge d\bar{z}_\eta \wedge \mathcal{L}_{\varphi_v^\lambda \frac{\partial}{\partial z_\lambda} } \mathcal{L}_{\psi_\eta^\delta \frac{\partial}{\partial z_\delta} } (w)
\end{align}}}
Then $\textnormal{(\ref{bf5})}$ and $\textnormal{(\ref{bf6})}$ imply $\textnormal{(\ref{bf7})}$. Then in particular if $\varphi= \psi$, then we have
\begin{align*}
\frac{1}{2}\mathcal{L}_{\left[ \varphi, \varphi \right]}\left( w \right)= \mathcal{L}_\varphi \mathcal{L}_\varphi(w)
\end{align*}

Next for a $C^\infty(0,p)$-vector $T=\sum_{\gamma_1,..., \gamma_p}\sum_{\alpha=1}^n d \bar{z}_{\gamma_1} \wedge \cdots \wedge d \bar{z}_{\gamma_p} T_{\gamma_1,..., \gamma_p}^\alpha \frac{\partial}{\partial z_\alpha}$ and $C^\infty(p,q)$-form $w$, we define the interior multiplication by
\begin{align*}
i_T\left( w \right)= \sum_{\gamma_1,..., \gamma_p} d\bar{z}_{\gamma_1} \wedge \cdots \wedge d \bar{z}_{\gamma_p} \wedge i_{\sum_{\alpha=1}^n T_{\gamma_1,..,\gamma_p}^\alpha \frac{\partial}{\partial z_\alpha}}\left( w \right)
\end{align*}
Then we can write
\begin{align*}
\mathcal{L}_T(w)= (-1)^p \partial i_T \left( w\right)+ i_T \partial(w)
\end{align*}
In particular, for a $C^\infty(0,1)$-vector $\varphi$ and a $C^\infty$ $1$-form $w$, we have
\begin{align*}
\mathcal{L}_{\varphi}(w)=  - \partial i_{\varphi}(w) + i_{\varphi}\left( \partial w \right)
\end{align*}
For a $C^\infty (0,1)$-vector $\varphi= \sum_{v, \lambda=1}^n \varphi_v^\lambda d \bar{z}_v \frac{\partial}{\partial z_\lambda}$ and $C^\infty$ 1-forms $w_1,..., w_q$, we claim that
\begin{align*}
i_\varphi \left(w_1\wedge \cdots \wedge w_q\right)= \sum_{\alpha=1}^q w_1 \wedge \cdots \wedge i_\varphi \left(w_\alpha\right) \wedge \cdots \wedge w_q
\end{align*}
In fact, we have
\begin{align*}
i_\varphi(w_1\wedge \cdots \wedge w_q)&= \sum_{\lambda, v=1}^n d\bar{z}_v \wedge i_{ \varphi_v^\lambda \frac{\partial}{\partial z_\lambda}}\left(w_1\wedge \cdots \wedge w_q \right) = \sum_{\lambda, v=1}^n d\bar{z}_v \wedge\left( \sum_{\alpha=1}^q (-1)^{\alpha-1} w_1 \wedge \cdots \wedge i_{\varphi_v^\lambda\frac{\partial}{\partial z_\lambda}}\left(w_\alpha\right)\wedge \cdots \wedge w_q   \right) \\
&= \sum_{\alpha=1}^q w_1 \wedge \cdots \wedge \left( \sum_{\lambda, v=1}^n d\bar{z}_v \wedge i_{\varphi_v^\lambda \frac{\partial}{\partial z_\lambda}}\left(w_\alpha\right)   \right) \wedge \cdots \wedge w_q=\sum_{\alpha=1}^q w_1 \wedge \cdots \wedge i_\varphi \left(w_\alpha\right) \wedge \cdots \wedge w_q
\end{align*}

Let us consider a $C^\infty\left(0,1 \right)$-vector $\varphi= \sum_{\lambda , v=1}^n \varphi_v^\lambda d\bar{z}_v \frac{\partial}{\partial z_\lambda}$ and $C^\infty$ 1-form $w=\sum_{\alpha=1}^n w_\alpha dz_\alpha$. Then we claim that
\begin{align*}
\mathcal{L}_\varphi\left( \partial w \right) = - \partial \mathcal{L}_\varphi\left( w \right)
\end{align*}
In fact, we note that $\mathcal{L}_\varphi\left( \partial w \right)= - \partial i_\varphi\left(\partial w\right) + i_\varphi \partial \left( \partial w \right)= - \partial \left(\mathcal{L}_\varphi(w)- (-1)\partial i_\varphi (w) \right)= - \partial \mathcal{L}_\varphi(w) $. 

Let us consider a $C^\infty$-vector of the form $T=\sum_{\alpha=1}^n  T^\alpha \frac{\partial}{\partial z_\alpha}$ and a $C^\infty$-form $w=\sum_{\beta=1}^n w^\beta dz_\beta$. Then we claim that
\begin{align}\label{bf11}
\bar{\partial} i_T\left( w \right) = i_{\bar{\partial T}}\left( w \right) - i_T \left( \bar{\partial} w \right)
\end{align}
In fact, $\textnormal{(\ref{bf11})}$ follows from 
\begin{align*}
\bar{\partial} i_{T}(w)&=\bar{\partial}\left( \sum_{\alpha=1}^n T^\alpha w^\alpha \right)=\sum_{\alpha,\beta=1}^n \frac{\partial T^\alpha}{\partial \bar{z}_\alpha} w^\alpha d\bar{z}_\beta + \sum_{\alpha,\beta=1}^n T^\alpha \frac{\partial w^\alpha}{\partial \bar{z}_\beta} d\bar{z}_\beta\\
i_{\bar{\partial} T}(w)&= i_{\sum_{\alpha ,\beta=1}^n \frac{\partial T^\alpha}{\partial \bar{z}_\beta} d\bar{z}_\beta \frac{\partial}{\partial z_\alpha}}\left(\sum_{\gamma=1}^n w^\gamma dz_\gamma \right) = \sum_{\alpha, \beta=1}^n \frac{\partial T^\alpha}{\partial \bar{z}_\alpha} w^\alpha d \bar{z}_\beta\\
i_T\left(\bar{\partial} w \right)&= i_{\sum_{\alpha=1}^n T^\alpha \frac{\partial}{\partial z_\alpha} } \left( \sum_{\beta,\gamma=1}^n \frac{\partial w^\beta}{\partial \bar{z}_\gamma} d\bar{z}_\gamma \wedge dz _\beta \right)= - \sum_{\alpha,\beta=1}^n T^\alpha \frac{\partial w^\alpha}{\partial \bar{z}_\beta} d\bar{z}_\beta
\end{align*}

Let us consider a $C^\infty(0,1)$-form $\varphi=\sum_{\lambda, v=1}^n \varphi_v^\lambda d\bar{z}_v \frac{\partial}{\partial z_\lambda}$ and a $C^\infty$-vector of the form $T=\sum_{\alpha=1}^n T_\alpha \frac{\partial}{\partial z_\alpha}$ and $C^\infty$ 1-form $w$, we claim that
\begin{align*}
i_{\left[ \varphi, T \right]}(w) = \mathcal{L}_\varphi i_T(w)+i_T\mathcal{L}_\varphi(w)
\end{align*}
In fact, we have
\begin{align*}
i_{\left[ \varphi, T \right]}(w)= \sum_{v, \lambda=1}^n d\bar{z}_v\wedge i_{\left[\varphi_v^\lambda \frac{\partial}{\partial z_\lambda}, T \right]}(w) =\sum_{v,\lambda=1}^n d\bar{z}_v \wedge\left( \mathcal{L}_{\varphi_v^\lambda \frac{\partial}{\partial z_\lambda}} i_T(w) -i_T \mathcal{L}_{\varphi_v^\lambda \frac{\partial}{\partial z_\lambda}}(w)    \right) = \mathcal{L}_\varphi i_T(w)+i_T\mathcal{L}_\varphi(w)
\end{align*}

Let us consider $C^{\infty}(0,1)$-vector $\varphi=\sum_{\lambda, v=1}^n \varphi_v^\lambda d\bar{z}_v\frac{\partial}{\partial z_\lambda}$, and $C^\infty$-form  $w=\sum_{\alpha=1}^n w_\alpha dz_\alpha$. Then we claim that
\begin{align} \label{bf12}
\bar{\partial}i_\varphi(w)=i_{\bar{\partial} \varphi}(w)+  i_\varphi (\bar{\partial} w)
\end{align}
In fact, $\textnormal{(\ref{bf12})}$ follows from
\begin{align*}
\bar{\partial} i_{\varphi}\left(w \right) &=\bar{\partial}\left( \sum_{\lambda, v=1}^n \varphi_v^\lambda w_\lambda d\bar{z}_v \right) \\
i_{\bar{\partial}\varphi}(w)&= i_{\sum_{\lambda, v,\eta=1}^n \frac{\partial \varphi_v^\lambda}{\partial \bar{z}_\eta}d\bar{z}_\eta\wedge d\bar{z}_v\frac{\partial}{\partial z_\lambda} } \left(\sum_{\alpha=1}^n w_\alpha dz_\alpha \right) = \sum_{\lambda, v , \eta=1}^n \frac{\partial \varphi_v^\lambda}{\partial \bar{z}_\eta} w_\lambda d\bar{z}_\eta\wedge d\bar{z}_v\\
i_{\varphi}( \bar{\partial} w) &= i_{\sum_{\lambda, v=1}^n \varphi_v^\lambda d\bar{z}_v\frac{\partial}{\partial z_\lambda}} \left( \sum_{\alpha, \eta=1}^n \frac{\partial w_\alpha}{\partial \bar{z}_\eta} d\bar{z}_\eta \wedge dz_\alpha    \right)= - \sum_{\lambda, v, \eta=1}^n \varphi_v^\lambda  \frac{\partial w_\lambda}{\partial \bar{z}_\eta} d\bar{z}_v\wedge d\bar{z}_\eta= \sum_{\lambda, v, \eta=1}^n \varphi_v^\lambda  \frac{\partial w_\lambda}{\partial \bar{z}_\eta} d\bar{z}_\eta\wedge d\bar{z}_v 
\end{align*}

For a $C^\infty(0,1)$-vector $\varphi= \sum_{ v, \lambda=1}^n \varphi_v^\lambda d \bar{z}_v \frac{\partial}{\partial z_\lambda}$ and a $C^\infty (0,1)$-vector $\psi= \sum_{\eta, \delta=1}^n \psi_\eta^\delta d \bar{z}_\eta \frac{\partial}{\partial z_\delta}$ and a $C^\infty$-1 form $w$, we claim that
\begin{align}\label{bf13}
i_{\left[ \varphi, \psi \right]} (w) = \mathcal{L}_\varphi i_{\psi}(w) - i_{\psi}\mathcal{L}_\varphi(w)
\end{align}
In fact, $\textnormal{(\ref{bf13})}$ follows from
\begin{align*}
&i_{\left[ \varphi, \psi \right]}(w)= \sum_{v, \lambda, \eta, \delta=1}^n d\bar{z}_v \wedge d\bar{z}_\eta  i_{\left[ \varphi_v^\lambda \frac{\partial}{\partial z_\lambda}, \psi_\eta^\delta \frac{\partial}{\partial z_\delta} \right]}(w) = \sum_{ v, \lambda, \eta, \delta=1}^n d\bar{z}_v\wedge d\bar{z}_\eta \left( \mathcal{L}_{\varphi_v^\lambda \frac{\partial}{\partial z_\lambda}} i_{\psi_\eta^\delta \frac{\partial}{\partial z_\delta}}   - i_{\psi_\eta^\delta \frac{\partial}{\partial z_\delta}}\mathcal{L}_{\varphi_v^\lambda \frac{\partial}{\partial z_\lambda}}  \right)(w)\\
&\mathcal{L}_\varphi i_\psi(w)= \mathcal{L}_\varphi\left( \sum_{\eta, \delta=1}^n d\bar{z}_\eta i_{\psi_\eta^\delta \frac{\partial}{\partial z_\delta}}(w)   \right) = \sum_{v, \lambda, \eta, \delta=1}^n d\bar{z}_v \wedge \mathcal{L}_{\varphi_v^\lambda \frac{\partial}{\partial z_\lambda}}\left( i_{\psi_\eta^\delta \frac{\partial}{\partial z_\delta}}(w)  d\bar{z}_\eta\right) = \sum_{ v, \lambda, \eta , \delta =1}^n  d\bar{z}_v \wedge d\bar{z}_\eta \mathcal{L}_{\varphi_v^\lambda \frac{\partial}{\partial z_\lambda}}\left( i_{\psi_\eta^\delta \frac{\partial}{\partial z_\delta}}(w) \right)\\
&i_\psi \mathcal{L}_\varphi(w)= i_\psi\left(  \sum_{v, \lambda=1}^n d\bar{z}_v \wedge \mathcal{L}_{\varphi_v^\lambda \frac{\partial}{\partial z_\lambda}}(w)          \right) = \sum_{v, \lambda, \eta, \delta=1}^n d\bar{z}_\eta \wedge i_{\psi_\eta^\delta\frac{\partial}{\partial z_\delta}}\left( d\bar{z}_v \wedge \mathcal{L}_{\varphi_v^\lambda \frac{\partial}{\partial z_\lambda}}(w)     \right) =- \sum_{v, \lambda, \eta, \delta=1}^n d\bar{z}_\eta \wedge d\bar{z}_v i_{\psi_\eta^\delta\frac{\partial}{\partial z_\delta}} \mathcal{L}_{\varphi_v^\lambda \frac{\partial}{\partial z_\lambda}}(w) 
\end{align*}

\section{Dolbeault resolution of $\Theta_{\mathcal{F}_0}^{\bullet}$ and Dolbeault type bicomplex associated to $\mathcal{N}_{\mathcal{F}_0}^{*\bullet}$}\label{app3}

\subsection{Dolbeault resolution of the leaf complex $\Theta_{\mathcal{F}_0}^\bullet$}\label{ad1}\

Let $\left(M, \Theta_{\mathcal{F}_0} \right)$ be a compact foliated complex manifold with $\Theta_{\mathcal{F}_0}$ locally free. We will describe the Dolbeault resolution of the leaf complex $\Theta_{\mathcal{F}_0}^\bullet$:
\begin{align*}
\Theta_{\mathcal{F}_0}^\bullet : \Theta_M \xrightarrow{D_0} \mathscr{H}om_{\mathcal{O}_M}\left( \Theta_{\mathcal{F}_0}, \frac{\Theta_M}{\Theta_{\mathcal{F}_0} } \right)  \xrightarrow{D_1} \mathscr{H}om_{\mathcal{O}_M}\left( \bigwedge^2 \Theta_{\mathcal{F}_0}, \frac{\Theta_M}{\Theta_{\mathcal{F}_0}} \right) \xrightarrow{D_2} \mathscr{H}om_{\mathcal{O}_M}\left(  \bigwedge^3 \Theta_{\mathcal{F}_0}, \frac{\Theta_M}{\Theta_{\mathcal{F}_0} } \right) \xrightarrow{D_3} \cdots
\end{align*}
Let us denote by $\mathcal{A}^{0,p}\left( \Theta_M \right)$ the sheaf of germs of $C^\infty(0,p)$-forms with coefficients in $\Theta_M$ and denote its global section by $A^{0,p}\left(M, \Theta_M \right)$. Then we have the Dolbeualt resolution $A^{0,0}\left(M, \Theta_M \right)\xrightarrow{\bar{\partial}} A^{0,1}\left( M, \Theta_M \right) \xrightarrow{\bar{\partial}} A^{0,2}\left( M, \Theta_M \right) \xrightarrow{\bar{\partial}} \cdots $ of $\Theta_M$. On the other hand, let us describe the Dolbeault resolution of $\mathscr{H}om_{\mathcal{O}_M}\left( \bigwedge^r \Theta_{\mathcal{F}_0}, \frac{\Theta_M}{\Theta_{\mathcal{F}_0} } \right)$. Let us denote by $\mathcal{A}^{0,p}\left( \mathscr{H}om_{\mathcal{O}_M}\left(\bigwedge^r \Theta_{\mathcal{F}_0}, \Theta_{\mathcal{F}_0} \right)\right)$ the sheaf of germs of $C^\infty(0,p)$-forms with coefficients in the locally free sheaf $\mathscr{H}om_{\mathcal{O}_M}\left( \bigwedge^r \Theta_{\mathcal{F}_0} ,  \Theta_{\mathcal{F}_0} \right)$ and denote the global section by $A^{0,p}\left( M, \mathscr{H}om_{\mathcal{O}_M}\left( \bigwedge^r \Theta_{\mathcal{F}_0}, \Theta_{\mathcal{F}_0} \right) \right)$. On the other hand, let us denote by $\mathcal{A}^{0,p}\left( \mathscr{H}om_{\mathcal{O}_M}\left(\bigwedge^r \Theta_{\mathcal{F}_0}, \Theta_M \right)\right)$ the sheaf of germs of $C^\infty(0,p)$-forms with coefficients in $\mathscr{H}om_{\mathcal{O}_M}\left( \bigwedge^r \Theta_{\mathcal{F}_0} , \Theta_M \right)$ and denote the global section by $A^{0,p}\left( M, \mathscr{H}om_{\mathcal{O}_M}\left( \bigwedge^r \Theta_{\mathcal{F}_0}, \Theta_M \right) \right)$. Then we have the following exact sequence of complexes of sheaves such that each vertical sequence is exact. 
{\tiny{\begin{center}
$\begin{CD}
@. \cdots @.\cdots @. \cdots \\
@. @A\bar{\partial}AA @A\bar{\partial}AA @AAA\\
0@>>> \mathcal{A}^{0,3}\left(\mathscr{H}om_{\mathcal{O}_M} \left(\bigwedge^r \Theta_{\mathcal{F}_0} , \Theta_{\mathcal{F}_0} \right)\right)@>>> \mathcal{A}^{0,3}\left(\mathscr{H}om_{\mathcal{O}_M} \left(\bigwedge^r \Theta_F, \Theta_M\right)\right)@>>>  \frac{\mathcal{A}^{0,3}\left(\mathscr{H}om_{\mathcal{O}_M} \left(\bigwedge^r \Theta_F, \Theta_M\right)\right)}{\mathcal{A}^{0,3}\left(\mathscr{H}om_{\mathcal{O}_M} \left(\bigwedge^r \Theta_{\mathcal{F}_0}, \Theta_{\mathcal{F}_0} \right)\right)}@>>> 0 \\
@. @A\bar{\partial}AA @A\bar{\partial}AA @AAA\\
0@>>> \mathcal{A}^{0,2}\left(\mathscr{H}om_{\mathcal{O}_M} \left(\bigwedge^r \Theta_{\mathcal{F}_0}, \Theta_{\mathcal{F}_0} \right)\right)@>>> \mathcal{A}^{0,2}\left(\mathscr{H}om_{\mathcal{O}_M} \left(\bigwedge^r \Theta_{\mathcal{F}_0} , \Theta_M\right)\right)@>>>  \frac{\mathcal{A}^{0,2}\left(\mathscr{H}om_{\mathcal{O}_M} \left(\bigwedge^r \Theta_{\mathcal{F}_0}, \Theta_M\right)\right)}{\mathcal{A}^{0,2} \left(\mathscr{H}om_{\mathcal{O}_M}\left(\bigwedge^r \Theta_{\mathcal{F}_0}, \Theta_{\mathcal{F}_0} \right)\right)}@>>> 0 \\
@. @A\bar{\partial}AA @A\bar{\partial}AA @AAA\\
0@>>> \mathcal{A}^{0,1}\left(\mathscr{H}om_{\mathcal{O}_M} \left(\bigwedge^r \Theta_{\mathcal{F}_0}, \Theta_{\mathcal{F}_0}\right)\right)@>>> \mathcal{A}^{0,1}\left(\mathscr{H}om_{\mathcal{O}_M}\left(\bigwedge^r \Theta_{\mathcal{F}_0}, \Theta_M\right)\right)@>>> \frac{ \mathcal{A}^{0,1}\left(\mathscr{H}om_{\mathcal{O}_M} \left(\bigwedge^r \Theta_{\mathcal{F}_0}, \Theta_M\right)\right)}{\mathcal{A}^{0,1}\left(\mathscr{H}om_{\mathcal{O}_M} \left(\bigwedge^r \Theta_{\mathcal{F}_0} , \Theta_{\mathcal{F}_0} \right)\right)}@>>> 0 \\
@. @A\bar{\partial}AA @A\bar{\partial}AA @AAA\\
0@>>> \mathcal{A}^{0,0}\left(\mathscr{H}om_{\mathcal{O}_M} \left(\bigwedge^r \Theta_{\mathcal{F}_0} ,  \Theta_{\mathcal{F}_0} \right)\right)@>>> \mathcal{A}^{0,0}\left(\mathscr{H}om_{\mathcal{O}_M} \left(\bigwedge^r \Theta_{\mathcal{F}_0} , \Theta_M\right)\right)@>>> \frac{\mathcal{A}^{0,0}\left(\mathscr{H}om_{\mathcal{O}_M} \left(\bigwedge^r \Theta_{\mathcal{F}_0}, \Theta_M\right)\right)}{\mathcal{A}^{0,0}\left(\mathscr{H}om_{\mathcal{O}_M} \left(\bigwedge^r \Theta_{\mathcal{F}_0}, \Theta_{\mathcal{F}_0} \right)\right)}@>>> 0 \\
@. @AAA @AAA @AAA\\
0@>>> \mathscr{H}om_{\mathcal{O}_M} \left(\bigwedge^r \Theta_{\mathcal{F}_0},  \Theta_{\mathcal{F}_0} \right) @>>> \mathscr{H}om_{\mathcal{O}_M}\left(\bigwedge^r \Theta_{\mathcal{F}_0}, \Theta_M\right) @>>> \mathscr{H}om_{\mathcal{O}_M} \left(\bigwedge^r \Theta_{\mathcal{F}_0} , \frac{ \Theta_M}{\Theta_{\mathcal{F}_0}}\right)@>>> 0
\end{CD}$
\end{center}}}
We note that since $H^i\left(M,  \mathcal{A}^{0,p}\left(\mathscr{H}om_{\mathcal{O}_M} \left(\bigwedge^r \Theta_{\mathcal{F}_0}, \Theta_{\mathcal{F}_0} \right)\right) \right)=0$ and $H^i\left(M, \mathcal{A}^{0,p}\left(\mathscr{H}om_{\mathcal{O}_M}\left(\bigwedge^r \Theta_{\mathcal{F}_0}, \Theta_M \right)\right)\right)=0$ for $i\geq 1$, we have from the exact sequence above
\begin{align*}
H^i\left(M,\frac{\mathcal{A}^{0,p}\left(\mathscr{H}om_{\mathcal{O}_M}\left(\bigwedge^r \Theta_{\mathcal{F}_0}, \Theta_M\right)\right)}{\mathcal{A}^{0,p}\left(\mathscr{H}om_{\mathcal{O}_M} \left(\bigwedge^r \Theta_{\mathcal{F}_0}, \Theta_{\mathcal{F}_0} \right)\right)}  \right)=0\,\,\,\,\,\textnormal{for}\,\,\, i\geq 1.
\end{align*}
and we have
{\small{\begin{align*}
H^0\left(M, \frac{\mathcal{A}^{0,p}\left(\mathscr{H}om_{\mathcal{O}_M}\left(\bigwedge^r \Theta_{\mathcal{F}_0} , \Theta_M  \right)\right)}{\mathcal{A}^{0,p}\left(\mathscr{H}om_{\mathcal{O}_M} \left(\bigwedge^r \Theta_{\mathcal{F}_0} , \Theta_{\mathcal{F}_0} \right)\right) }\right)\cong \frac{A^{0,p}\left(M,\mathscr{H}om_{\mathcal{O}_M} \left(\bigwedge^r \Theta_{\mathcal{F}_0}, \Theta_M\right) \right)}{A^{0,p}\left(M, \mathscr{H}om_{\mathcal{O}_M}\left(\bigwedge^r \Theta_{\mathcal{F}_0} ,  \Theta_{\mathcal{F}_0}\right)\right)}
\end{align*}}}
Hence we have the following Dolbeault resolution of $\Theta_{\mathcal{F}_0}^\bullet$.
\begin{equation}
\begin{CD}
\cdots \\
@A\hat{D}_3AA \\
\frac{A^{0,0}\left( M,  \bigwedge^3 \Theta_{\mathcal{F}_0}^*\otimes \Theta_M\right)}{A^{0,0}\left(M, \bigwedge^3 \Theta_{\mathcal{F}_0}^*\otimes \Theta_{\mathcal{F}_0} \right)}@>-\bar{\partial}>> \cdots\\
@A\hat{D}_2AA @A\hat{D}_2AA \\
\frac{A^{0,0}\left(M,  \bigwedge^2 \Theta_{\mathcal{F}_0}^*\otimes \Theta_M\right)}{A^{0,0}\left(M, \bigwedge^2 \Theta_{\mathcal{F}_0}^* \otimes \Theta_{\mathcal{F}_0} \right)} @>\bar{\partial}>> \frac{A^{0,1}\left( M, \bigwedge^2 \Theta_{\mathcal{F}_0}^*\otimes \Theta_M \right)}{A^{0,1}\left(M, \bigwedge^2 \Theta_{\mathcal{F}_0}^*\otimes \Theta_{\mathcal{F}_0} \right)}@>-\bar{\partial}>> \cdots \\
@A\hat{D}_1AA @A\hat{D}_1AA @A\hat{D}_1AA\\
\frac{A^{0,0}\left(M, \Theta_{\mathcal{F}_0}^*\otimes \Theta_M\right)}{A^{0,0}\left(M, \Theta_{\mathcal{F}_0}^*\otimes \Theta_{\mathcal{F}_0} \right)} @>-\bar{\partial}>> \frac{A^{0,1}\left(M, \Theta_{\mathcal{F}_0}^*\otimes \Theta_M\right)}{A^{0,1}\left(M, \Theta_{\mathcal{F}_0}^*\otimes \Theta_{\mathcal{F}_0} \right)} @>\bar{\partial}>> \frac{A^{0,2}\left(M, \Theta_{\mathcal{F}_0}^*\otimes \Theta_M\right)}{A^{0,2}\left(M, \Theta_{\mathcal{F}_0}^*\otimes \Theta_{\mathcal{F}_0}\right)} @>-\bar{\partial}>>\cdots \\
@A\hat{D}_0AA @A\hat{D}_0AA @A\hat{D}_0AA @A\hat{D}_0AA\\
A^{0,0}(M, \Theta_M) @>\bar{\partial} >> A^{0,1}(M, \Theta_M) @>-\bar{\partial}>> A^{0,2}(M, \Theta_M)  @>\bar{\partial}>> A^{0,3}(M, \Theta_M) @>-\bar{\partial}>> \cdots\\
\end{CD}
\end{equation}
We will explicitly describe $\hat{D}_0$ and $\hat{D}_1$ and $\hat{D}_2$ that we use in the main body of the paper. Let $\mathcal{U}=\left\{U_i\right\}$ be a Stein open covering of $M$ by coordinate neighborhoods with local coordinates $\left(z_i^1,...,z_i^n \right)$ on $U_i$ such that $\Gamma\left(U_j , \Theta_{\mathcal{F}_0}\right)$ is generated by $T_{j}^1,..., T_{j}^p\in \Gamma\left( U_j, \Theta_M \right)$ with the relation $T_{j}^\alpha= \sum_{\beta=1}^p r_{jk}^{\alpha\beta} T_{j}^\beta$ for some $r_{jk}^{\alpha\beta}\in \Gamma\left( U_j\cap U_k, \mathcal{O}_M \right)$ and $\left[ T_{j}^\alpha, T_{j}^\beta \right]= \sum_{\gamma=1}^p g_{j\alpha\beta}^\gamma T_{j}^\gamma$ for some $g_{j\alpha\beta}^\gamma \in \Gamma\left( U_j, \mathcal{O}_M \right)$.

\subsubsection{\textnormal{Description of $\hat{D}_0$}}\

We describe $\hat{D}_0: A^{0,p}\left( M, \Theta_M \right)\to \frac{A^{0,p}\left( M, \mathscr{H}om_{\mathcal{O}_M}\left( \Theta_{\mathcal{F}_0}, \Theta_M \right) \right)}{A^{0,p}\left( M, \mathscr{H}om_{\mathcal{O}_M}\left( \Theta_{\mathcal{F}_0}, \Theta_{\mathcal{F}_0}  \right) \right)}$. For $\phi\in A^{0,p}\left( M, \Theta_M \right)$, we define an element $\eta_{ij}\in \Gamma\left(U_{ij}, \mathcal{A}^{0,p}\left(\mathscr{H}om_{\mathcal{O}_M}\left( \Theta_{\mathcal{F}_0}, \Theta_{\mathcal{F}_0} \right) \right) \right)$ by
\begin{align*}
\eta_{ij}: \Gamma\left( U_{ij}, \Theta_{\mathcal{F}_0} \right) &\to \Gamma\left( U_{ij}, \mathcal{A}^{0,p}\left( \Theta_{\mathcal{F}_0} \right) \right) \\
 T_i^\alpha &\mapsto \left[ \phi, T_i^\alpha \right] - \sum_{\beta=1}^p r_{ij}^{\alpha\beta}\left[ \phi, T_j^\beta \right]= \sum_{\beta=1}^p \left[ \phi, r_{ij}^{\alpha \beta} \right] T_j^\beta
\end{align*}
and linearly extends to $\Gamma\left( U_{ij}, \Theta_{\mathcal{F}_0} \right)$. Then we see that
\begin{align*}
\left(\eta_{ij}- \eta_{ik} + \eta_{jk} \right)\left( T_i^\alpha \right)&= \eta_{ij}\left(T_i^\alpha \right)- \eta_{ik}\left( T_i^\alpha \right) + \sum_{\beta=1}^p r_{ij}^{\alpha\beta} \eta_{jk}\left( T_j^\beta \right) \\
&=\left[ \phi, T_i^\alpha \right] - \sum_{\beta=1}^p r_{ij}^{\alpha\beta}\left[ \phi, T_j^\beta \right] - \left[ \phi, T_i^\alpha \right] + \sum_{\beta=1}^p r_{ik}^{\alpha\beta}\left[ \phi, T_k^\beta \right] + \sum_{\beta=1}^p r_{ij}^{\alpha\beta}\left( \left[ \phi, T_i^\beta \right] - \sum_{\gamma=1}^p r_{jk}^{\beta \gamma} \left[ \phi, T_k^\gamma \right] \right)= 0
\end{align*}
This implies that there exist $\left\{ \eta_i \right\}$ with $\eta_i\in \Gamma\left( U_i, \mathcal{A}^{0,p}\left( \mathscr{H}om_{\mathcal{O}_M}\left( \Theta_{\mathcal{F}_0}, \Theta_{\mathcal{F}_0} \right) \right) \right)$ such that $\eta_j- \eta_i= \eta_{ij}$ where
\begin{align*}
\eta_i: \Gamma\left(  U_i, \Theta_{\mathcal{F}_0} \right) &\to \Gamma\left( U_i, \mathcal{A}^{0,p}\left( \Theta_{\mathcal{F}_0} \right) \right) \\
 T_i^\alpha &\mapsto W_i^\alpha
\end{align*}
and then linearly extends to $\Gamma\left(U_i, \Theta_{\mathcal{F}_0} \right)$. We define $\hat{D}_{0i}' \left( \phi\right) \in \Gamma\left( U_i, \mathcal{A}^{0,p}\left( \mathscr{H}om_{\mathcal{O}_M}\left( \Theta_{\mathcal{F}_0}, \Theta_M \right) \right) \right)$ by $\hat{D}_{0i}' (\phi ) \left(T_i^\alpha \right)=\left[ \phi, T_i^\alpha \right] + W_i^\alpha$. Then we have
\begin{align*}
\left[ \phi, T_i^\alpha \right] + W_i^\alpha - \sum_{\beta=1}^p r_{ij}^{\alpha\beta}\left( \left[ \phi, T_j^\beta \right] + W_j^\beta \right)=0
\end{align*}
This implies that $\hat{D}_0'(\phi):=\left\{ \hat{D}_{0i}' (\phi)  \right\}\in A^{0,p}\left( M, \mathscr{H}om_{\mathcal{O}_M}\left( \Theta_{\mathcal{F}_0}, \Theta_M \right) \right)$. Then we define $\hat{D}_0\left( \phi \right)$ to be the image of $\hat{D}_0'(\phi)$ in $\frac{A^{0,p}\left( M, \mathscr{H}om_{\mathcal{O}_M}\left( \Theta_{\mathcal{F}_0}, \Theta_M \right) \right)}{A^{0,p}\left(  M, \mathscr{H}om_{\mathcal{O}_M}\left( \Theta_{\mathcal{F}_0}, \Theta_{\mathcal{F}_0} \right) \right)}$. We note that it is independent of choices of $W_i^\alpha$. Let $\eta_i'$ be another choice such that $\eta_j'- \eta_i'= \eta_{ij}$ where $\eta_i':\Gamma\left( U_i, \Theta_{\mathcal{F}_0} \right) \to \Gamma\left( U_i, \mathcal{A}^{0,p}\left( \Theta_{\mathcal{F}_0} \right) \right), T_i^\alpha \mapsto W_i'^\alpha$. Then $W_i^\alpha - W_i'^\alpha  = \sum_{\beta=1}^p r_{ij}^{\alpha\beta} \left( W_j^\beta - W_j'^\beta\right)$, so that $\sigma=\left\{\eta_i-\eta_i' \right\} \in A^{0,p}\left( M, \mathscr{H}om_{\mathcal{O}_M}\left( \Theta_{\mathcal{F}_0}, \Theta_M \right) \right)$ and $\left[ \phi, T_i^\alpha \right] + W_i^\alpha - \left(\left[ \phi, T_i^\alpha \right] + W_i'^\alpha  \right)= W_i^\alpha - W_i'^\alpha=\sigma\left( T_i^\alpha \right)$, so that $\hat{D}_0(\phi)$ is independent of choices of $W_i^\alpha$.

\subsubsection{\textnormal{Description of $\hat{D}_1$}}\

We describe $\hat{D}_1: \frac{A^{0,p}\left( M, \mathscr{H}om_{\mathcal{O}_M}\left( \Theta_{\mathcal{F}_0}, \Theta_M \right) \right) }{A^{0,p}\left( M , \mathscr{H}om_{\mathcal{O}_M}\left( \Theta_{\mathcal{F}_0}, \Theta_{\mathcal{F}_0} \right) \right)} \to \frac{A^{0,p}\left( M, \mathscr{H}om_{\mathcal{O}_M}\left( \bigwedge^2 \Theta_{\mathcal{F}_0} , \Theta_M    \right) \right)}{A^{0,p}\left( M, \mathscr{H}om_{\mathcal{O}_M} \left(  \bigwedge^2 \Theta_{\mathcal{F}_0}, \Theta_{\mathcal{F}_0}  \right) \right)}$. For $\overline{\phi} \in \frac{A^{0,p}\left( M, \mathscr{H}om_{\mathcal{O}_M}\left( \Theta_{\mathcal{F}_0}, \Theta_M \right) \right) }{A^{0,p}\left( M , \mathscr{H}om_{\mathcal{O}_M}\left( \Theta_{\mathcal{F}_0}, \Theta_{\mathcal{F}_0} \right) \right)} $, we define an element $\eta_{ij}\in \Gamma\left( U_{ij}, \mathcal{A}^{0,p} \left( \bigwedge^2 \Theta_{\mathcal{F}_0}, \Theta_{\mathcal{F}_0} \right) \right)$ \footnote{We note that $\phi\left( T_i^\alpha \right)$ is of the form $\phi\left( T_i^\alpha \right) = \sum_{\eta=1}^n \sum_{\gamma_1,..., \gamma_p} d \bar{z}_i^{\gamma_1}\wedge \cdots \wedge d \bar{z}_i^{\gamma_p} A_{\gamma_1,,..,\gamma_p}^\eta \frac{\partial}{\partial z_i^\eta}$. Then we note that 
\begin{align*} 
\left[ T_i^\beta, \phi\left( T_i^\alpha \right)\right] = \left[ T_i^\beta,  \sum_{\gamma_1,..., \gamma_p}  \sum_{\eta=1}^n d\bar{z}_i^{\gamma_1}\wedge \cdots \wedge d \bar{z}_i^{\gamma_p} A_{\gamma_1,...,\gamma_p}^\eta \frac{\partial}{\partial z_i^\beta} \right] = \sum_{\gamma_1,..., \gamma_p} d\bar{z}_i^{\gamma_1}\wedge \cdots \wedge d \bar{z}_i^{\gamma_p} \left[ T_i^\beta, \sum_{\eta=1}^n  A_{\gamma_1,...,\gamma_p}^\eta \frac{\partial}{\partial z_i^\eta} \right]
\end{align*}
Hence we do note need to care for the sign, i.e the order of $(0,p)$-form in the local trivilization.} by
{\small{\begin{align*}
\eta_{ij}: \Gamma\left( U_{ij}, \bigwedge^2 \Theta_{\mathcal{F}_0} \right) &\to \Gamma\left( U_{ij}, \mathcal{A}^{0,p}\left( \Theta_{\mathcal{F}_0} \right) \right) \\
 T_i^\alpha \wedge T_i^\beta &\mapsto  \left[ T_i^\alpha, \phi\left( T_i^\beta \right) \right] - \left[ T_i^\beta, \phi\left( T_i^\alpha \right) \right] - \phi \left( \left[ T_i^\alpha, T_i^\beta \right] \right) - \sum_{\gamma, \eta=1}^p r_{ij}^{\alpha \gamma} r_{ij}^{\beta \eta} \left(  \left[T_j^\gamma, \phi\left(T_j^\eta\right)\right]-\left[T_j^\eta,  \phi\left(T_j^\gamma\right)\right]- \phi\left(\left[T_j^\gamma, T_j^\eta \right]\right)     \right)
\end{align*}}}
We will show that the image of $\eta_{ij}$ actually lies in $\Gamma\left(U_{ij}, \mathcal{A}^{0,p}\left(\Theta_{\mathcal{F}_0} \right) \right)$. In fact, we note that
\begin{align*}
&\left[T_i^\alpha, \phi\left(T_i^\beta\right)\right]-\left[T_i^\beta,  \phi\left(T_i^\alpha\right)\right]- \phi\left( \left[T_i^\alpha, T_i^\beta \right] \right) \\
&=\left[\sum_{\gamma=1}^p r_{ij}^{\alpha \gamma} T_j^{\gamma} ,  \sum_{\eta=1}^p r_{ij}^{\beta\eta}\phi\left(T_j^\eta\right)\right] -\left[ \sum_{\eta=1}^p r_{ij}^{\beta\eta} T_j^\eta , \sum_{\gamma=1}^p r_{ij}^{\alpha \gamma} \phi \left(T_j^\gamma \right)    \right] - \sum_{\xi,\eta=1}^p g_{i\alpha\beta}^\xi r_{ij}^{\xi\eta} \phi\left( T_j^\eta  \right)\\
&=\sum_{\gamma,\eta=1}^p r_{ij}^{\alpha\gamma} r_{ij}^{\beta\eta} \left( \left[T_j^\gamma, \phi \left(T_j^\eta \right) \right] - \left[T_j^\eta, \phi \left(T_j^\gamma \right) \right]    \right) +\sum_{\gamma,\eta=1}^p r_{ij}^{\alpha\gamma} T_j^\gamma \left(r_{ij}^{\beta\eta}\right)\phi \left(T_j^\eta \right) - \sum_{\gamma,\eta=1}^p r_{ij}^{\beta\eta}\left[\phi\left(T_j^\eta\right), r_{ij}^{\alpha\gamma}\right] T_j^\gamma\\
& - \sum_{\eta,\gamma=1}^p r_{ij}^{\beta\eta} T_j^\eta \left(r_{ij}^{\alpha\gamma} \right)\phi \left(T_j^\gamma \right) + \sum_{\eta,\gamma=1}^p r_{ij}^{\alpha\gamma}\left[\phi\left(T_j^\gamma\right), r_{ij}^{\beta\eta} \right] T_j^\eta - \sum_{\xi,\eta=1}^p g_{i\alpha\beta}^\xi r_{ij}^{\xi\eta}\phi \left(T_j^\eta\right)
\end{align*}
On the other hand,
\begin{align*}
 &\sum_{\gamma,\eta=1}^p r_{ij}^{\alpha\gamma}r_{ij}^{\beta \eta}\left(  \left[T_j^\gamma, \phi(T_j^\eta)\right]-\left[T_j^\eta,  \phi(T_j^\gamma)\right]- \phi(\left[T_j^\gamma, T_j^\eta \right]) \right) \\
 & =\sum_{\gamma,\eta=1}^p r_{ij}^{\alpha\gamma} r_{ij}^{\beta\eta}\left( \left[T_j^\gamma, \phi \left(T_j^\eta \right) \right]- \left[ T_j^\eta, \phi \left(T_j^\gamma \right)\right]  \right) - \sum_{\gamma,\eta,\xi=1}^p r_{ij}^{\alpha\gamma} r_{ij}^{\beta\eta}g_{j\gamma\eta}^\xi\phi \left(T_j^\xi \right)
\end{align*}
Then from $(\ref{te19})$, we have
\begin{align*}
\eta_{ij}: \Gamma\left( U_{ij}, \bigwedge^2 \Theta_{\mathcal{F}_0} \right) &\to \Gamma\left( U_{ij}, \mathcal{A}^{0,p}\left( \Theta_{\mathcal{F}_0} \right) \right) \\
T_i^\alpha \wedge T_i^\beta &\mapsto   \sum_{\eta=1}^p \left[\phi \left(T_i^\alpha \right), r_{ij}^{\beta\eta}\right]T_j^\eta - \sum_{\gamma=1}^p \left[ \phi\left(T_i^\beta\right), r_{ij}^{\alpha\gamma} \right]T_j^\gamma
\end{align*}
Then $\eta_{jk}- \eta_{ik} + \eta_{ij}=0$. This implies that there exist $\left\{ \eta_i \right\}$ with $\eta_i \in \Gamma\left( U_i, \mathcal{A}^{0,p}\left( \mathscr{H}om_{\mathcal{O}_M} \left(  \Theta_{\mathcal{F}_0},  \Theta_{\mathcal{F}_0}   \right) \right) \right)$ such that $\eta_j- \eta_i= \eta_{ij}$ where 
\begin{align*}
\eta_i : \Gamma\left( U_i, \bigwedge^2 \Theta_{\mathcal{F}_0} \right) &\to \Gamma\left( U_i, \mathcal{A}^{0,p} \left( \Theta_{\mathcal{F}_0} \right)   \right) \\
 T_i^\alpha \wedge T_i^\beta &\mapsto W_i^{\alpha\beta}
\end{align*}
Then we define $\hat{D}_{1i}' (\phi)\in \Gamma\left( U_i, \mathcal{A}^{0,p}\left( \mathscr{H}om_{\mathcal{O}_M}\left( \bigwedge^2 \Theta_{\mathcal{F}_0}, \Theta_M  \right) \right) \right)$ by $\hat{D}_{1i}'(\phi)\left( T_i^\alpha\wedge T_i^\beta \right)=  \left[ T_i^\alpha, \phi\left( T_i^\beta \right) \right]-\left[ T_i^\beta, \phi\left( T_i^\alpha \right) \right]- \phi\left( \left[ T_i^\alpha, T_i^\beta \right] \right) + W_i^{\alpha\beta}$. Then $\hat{D}_1'(\phi):= \left\{ \hat{D}_{1i}'(\phi) \right\}\in A^{0,p}\left( M , \mathscr{H}om_{\mathcal{O}_M} \left( \bigwedge^2 \Theta_{\mathcal{F}_0} , \Theta_M \right) \right)$. Then we define $\hat{D}_1\left( \overline{\phi} \right)$ to be the image of $\hat{D}_1'(\phi)$ in $\frac{A^{0,p}\left( M , \mathscr{H}om_{\mathcal{O}_M}\left( \bigwedge^2  \Theta_{\mathcal{F}_0} ,\Theta_M \right) \right)}{A^{0,p}\left( M, \mathscr{H}om_{\mathcal{O}_M} \left(  \bigwedge^2 \Theta_{\mathcal{F}_0} , \Theta_{\mathcal{F}_0}   \right) \right)}$. We note that it is independent of choices of $W_i^{\alpha\beta}$. Moreover if $\phi \in A^{0,p}\left( M, \mathscr{H}om_{\mathcal{O}_M}\left( \Theta_{\mathcal{F}_0}, \Theta_{\mathcal{F}_0} \right) \right)$, then $\hat{D}_1'(\phi)\in A^{0,p}\left( M, \mathscr{H}om_{\mathcal{O}_M}\left( \bigwedge^2 \Theta_{\mathcal{F}_0} , \Theta_{\mathcal{F}_0} \right) \right)$.  Hence $\hat{D}_0\left(\overline{\phi} \right)$ is independent of choice of $\phi$.

\subsubsection{\textnormal{Description of $\hat{D}_2$}}\

We describe $\hat{D}_2: \frac{A^{0,p}\left( M, \mathscr{H}om_{\mathcal{O}_M}\left( \bigwedge^2 \Theta_{\mathcal{F}_0}, \Theta_M \right) \right)}{A^{0,p}\left( M, \mathscr{H}om_{\mathcal{O}_M} \left( \bigwedge^2 \Theta_{\mathcal{F}_0}, \Theta_{\mathcal{F}_0}   \right) \right)} \to \frac{A^{0,p}\left( M, \mathscr{H}om_{\mathcal{O}_M}\left( \bigwedge^3 \Theta_{\mathcal{F}_0}, \Theta_M \right) \right)}{A^{0,p}\left( M, \mathscr{H}om_{\mathcal{O}_M} \left( \bigwedge^3 \Theta_{\mathcal{F}_0}, \Theta_{\mathcal{F}_0}   \right) \right)} $. In a similar way to define $\hat{D}_1$ we define $\hat{D}_2$. For $\overline{B}\in \frac{A^{0,p}\left( M, \mathscr{H}om_{\mathcal{O}_M} \left( \bigwedge^2 \Theta_{\mathcal{F}_0}, \Theta_M \right) \right)}{A^{0,p}\left( M, \mathscr{H}om_{\mathcal{O}_M} \left( \bigwedge^2 \Theta_{\mathcal{F}_0} , \Theta_{\mathcal{F}_0}  \right) \right)}$, we define an element $\eta_{ij} \in \Gamma\left( U_{ij},  \mathscr{H}om_{\mathcal{O}_M} \left( \bigwedge^3 \Theta_{\mathcal{F}_0}, \Theta_{\mathcal{F}_0}  \right) \right)$ 
{\Small{\begin{align*}
&\eta_{ij}: \Gamma\left( U_{ij}, \bigwedge^3 \Theta_{\mathcal{F}_0} \right) \to  \Gamma\left(  U_{ij} , \mathcal{A}^{0,p}\left(\Theta_{\mathcal{F}_0} \right) \right) \\
              &  T_i^\alpha \wedge T_i^\beta \wedge T_i^\gamma \\
               & \mapsto  \left[T_i^\alpha, B\left(T_i^\beta, T_i^\gamma\right)\right]- \left[T_i^\beta, B\left(T_i^\alpha, T_i^\gamma\right)\right]+ \left[T_i^\gamma, B \left(T_i^\alpha ,  T_i^\beta \right) \right] - B\left( \left[T_i^\alpha, T_i^\beta \right], T_i^\gamma \right) + B\left(\left[T_i^\alpha, T_i^\gamma \right], T_i^\beta \right) -B\left( \left[T_i^\beta, T_i^\gamma  \right], T_i^\alpha \right)\\
               &- \sum_{\eta,\delta, \xi=1}^p r_{ij}^{\alpha \eta} r_{ij}^{\beta \delta} r_{ij}^{\gamma \xi}\left(\left[T_j^\eta, B \left(T_j^\delta, T_j^\xi \right)\right]- \left[T_j^\delta, B\left(T_j^\eta, T_i^\xi \right)\right]+ \left[T_j^\xi, B \left(T_j^\eta ,  T_j^\delta \right) \right] - B\left( \left[T_j^\eta, T_j^\delta \right], T_j^\xi \right) + B\left(\left[T_j^\eta, T_j^\xi \right], T_j^\delta \right) -B\left( \left[T_j^\delta, T_j^\xi  \right], T_j^\eta \right) \right)
\end{align*}}}
Then $\eta_{ij}- \eta_{ik} + \eta_{jk}=0$. This implies that there exits $\left\{ \eta_i \right\}$ with $\eta_i \in \Gamma\left( U_i, \mathcal{A}^{0,p}\left( \mathscr{H}om_{\mathcal{O}_M}\left( \bigwedge^3 \Theta_{\mathcal{F}_0}, \Theta_{\mathcal{F}_0} \right) \right)   \right)$ such that $\eta_j- \eta_i= \eta_{ij}$ where $\eta_{ij}:\Gamma\left(U_i, \bigwedge^3 \Theta_{\mathcal{F}_0} \right) \to \Gamma\left( U_i, \mathcal{A}^{0,p} \left( \Theta_{\mathcal{F}_0} \right) \right), T_i^\alpha \wedge T_i^\beta \wedge T_i^\gamma \mapsto W_i^{\alpha\beta\gamma}$. Then we define $\hat{D}_{2i}'(B)\in \Gamma\left( U_i, \mathcal{A}^{0,p}\left( \mathscr{H}om_{\mathcal{O}_M} \left( \bigwedge^3 \Theta_{\mathcal{F}_0}, \Theta_M \right) \right)\right)$ by
{\tiny{\begin{align*}
\hat{D}_{2i}'(B)\left(T_i^\alpha\wedge T_i^\beta \wedge T_i^\gamma \right)=  \left[T_i^\alpha, B\left(T_i^\beta, T_i^\gamma\right)\right]- \left[T_i^\beta, B\left(T_i^\alpha, T_i^\gamma\right)\right]+ \left[T_i^\gamma, B \left(T_i^\alpha ,  T_i^\beta \right) \right] - B\left( \left[T_i^\alpha, T_i^\beta \right], T_i^\gamma \right) + B\left(\left[T_i^\alpha, T_i^\gamma \right], T_i^\beta \right) -B\left( \left[T_i^\beta, T_i^\gamma  \right], T_i^\alpha \right) + W_i^{\alpha\beta\gamma}
\end{align*}}}
Then $\hat{D}_2'(B):=\left\{ \hat{D}_{2i}'(B) \right\}\in A^{0,p}\left( M , \mathscr{H}om_{\mathcal{O}_M}\left( \bigwedge^3 \Theta_{\mathcal{F}_0}, \Theta_M   \right) \right)$ and we define $\hat{D}_2\left( \overline{B} \right)$ to be the image of $\hat{D}_2'(B)$ in $\frac{A^{0,p}\left( M, \mathscr{H}om_{\mathcal{O}_M}\left( \bigwedge^3 \Theta_{\mathcal{F}_0}, \Theta_M \right) \right)}{A^{0,p}\left( M, \mathscr{H}om_{\mathcal{O}_M} \left( \bigwedge^3 \Theta_{\mathcal{F}_0}, \Theta_{\mathcal{F}_0}    \right) \right)}$.

In practice in the main body of paper, we only need to show that $\hat{D}_2\left( \overline{B}\right)=0$ for some $\overline{B}\in \frac{A^{0,0}\left( M, \mathscr{H}om_{\mathcal{O}_M}\left( \bigwedge^2 \Theta_{\mathcal{F}_0}, \Theta_M \right) \right)}{A^{0,0}\left( M, \mathscr{H}om_{\mathcal{O}_M} \left( \bigwedge^2 \Theta_{\mathcal{F}_0}, \Theta_{\mathcal{F}_0}  \right) \right)} $. In this case, we don't need to find a global section $\hat{D}_2' (B)$ in $A^{0,0}\left( M, \mathscr{H}om_{\mathcal{O}_M}\left( \bigwedge^3 \Theta_{\mathcal{F}_0}, \Theta_M \right) \right)$. It is enough to check that locally on $U_i$, 
{\small{\begin{align*}
 \left[T_i^\alpha, B\left(T_i^\beta, T_i^\gamma\right)\right]- \left[T_i^\beta, B\left(T_i^\alpha, T_i^\gamma\right)\right]+ \left[T_i^\gamma, B \left(T_i^\alpha ,  T_i^\beta \right) \right] - B\left( \left[T_i^\alpha, T_i^\beta \right], T_i^\gamma \right) + B\left(\left[T_i^\alpha, T_i^\gamma \right], T_i^\beta \right) -B\left( \left[T_i^\beta, T_i^\gamma  \right], T_i^\alpha \right) \in \Theta_{\mathcal{F}_0}
\end{align*}}}

\subsection{Dolbeault type bicomplex associated to the dual leaf complex $\mathcal{N}_{\mathcal{F}_0}^{* \bullet}$} \label{ad2}\

Let $\left( M, \mathcal{N}_{\mathcal{F}_0}^* \right)$ be a compact foliated complex manifold with $\mathcal{N}_{\mathcal{F}_0}^*$ locally free. We will construct a Dolbeault-type bicomplex associated the the dual leaf complex $\mathcal{N}_{\mathcal{F}_0}^{* \bullet}$:
\begin{align*}
\Theta_M \xrightarrow{E_0} \mathscr{H}om_{\mathcal{O}_M}\left( \mathcal{N}_{\mathcal{F}_0}^* ,  \frac{\Omega_M^1}{\mathcal{N}_{\mathcal{F}_0}^*} \right)  \xrightarrow{E_1} \mathscr{H}om_{\mathcal{O}_M}\left( \mathcal{N}_{\mathcal{F}_0}^*, \tilde{\mathcal{S}}^2 \right) \xrightarrow{E_2} \mathscr{H}om_{\mathcal{O}_M}\left( \mathcal{N}_{\mathcal{F}_0}^*, \tilde{\mathcal{S}}^3 \right) \xrightarrow{E_3} \cdots
\end{align*}

Since $\tilde{\mathcal{S}}^r, r\geq 2$ is not necessarily locally free, we do note have the Dolbeault resolution of $\tilde{\mathcal{S}}^r$ on $M$. However, since $\tilde{\mathcal{S}}^r$ is locally free on $M-S$, we have the Dolbeault resolution of $\tilde{\mathcal{S}}^r$ on $M-S$. Let us denote $\mathcal{A}^{0,p}\left( \mathscr{H}om_{\mathcal{O}_M}\left(\mathcal{N}_{\mathcal{F}_0}^*, \tilde{\mathcal{S}}^r \right)|_{M-S}\right)$ be the sheaf of germs of $C^\infty(0,p)$-forms with coefficient in $\mathscr{H}om_{\mathcal{O}_M}\left(\mathcal{N}_{\mathcal{F}_0}^*, \tilde{\mathcal{S}}^r \right)|_{M-S}$ on $M-S$ and denote the global section by $A^{0,p}\left( M-S, \mathscr{H}om_{\mathcal{O}_M}\left(\mathcal{N}_{\mathcal{F}_0}^*, \tilde{\mathcal{S}}^r \right) \right)$. Then we have the well-defined complex
\begin{align*}
A^{0,0}\left(M-S, \mathscr{H}om_{\mathcal{O}_M}\left( \mathcal{N}_{\mathcal{F}_0}^* , \tilde{\mathcal{S}}^r \right) \right) \xrightarrow{\bar{\partial}} A^{0,1}\left( M-S, \mathscr{H}om_{\mathcal{O}_M}\left( \mathcal{N}_{\mathcal{F}_0}^*, \tilde{\mathcal{S}}^r \right) \right) \xrightarrow{\bar{\partial}} A^{0,2}\left( M, \mathscr{H}om_{\mathcal{O}_M}\left( \mathcal{N}_{\mathcal{F}_0}^*, \tilde{\mathcal{S}}^r \right) \right) \xrightarrow{\bar{\partial}} \cdots
\end{align*}
Let $\mathcal{A}^{0,p}\left( \mathscr{H}om_{\mathcal{O}_M}\left( \mathcal{N}_{\mathcal{F}_0}^*, \mathcal{N}_{\mathcal{F}_0}^* \right) \right)$ be the sheaf of germs of $C^\infty(0,p)$-forms with coefficients in $\mathscr{H}om_{\mathcal{O}_M}\left( \mathcal{N}_{\mathcal{F}_0}^*, \mathcal{N}_{\mathcal{F}_0}^* \right)$ on $M$ and denote the global section by $A^{0,p}\left( M, \mathscr{H}om_{\mathcal{O}_M}\left( \mathcal{N}_{\mathcal{F}_0}^*, \mathcal{N}_{\mathcal{F}_0}^* \right) \right)$. On the other hand, let $\mathcal{A}^{0,p}\left( \mathscr{H}om_{\mathcal{O}_M}\left( \mathcal{N}_{\mathcal{F}_0}^*, \Omega_M^1 \right) \right)$ be the sheaf of germs of $C^\infty(0,p)$-forms with coefficients in $\mathscr{H}om_{\mathcal{O}_M}\left( \mathcal{N}_{\mathcal{F}_0}^*, \Omega_M^1 \right)$ on $M$ and denote the global section by $A^{0,p}\left( M, \mathscr{H}om_{\mathcal{O}_M}\left( \mathcal{N}_{\mathcal{F}_0}^*, \Omega_M^1 \right) \right)$. Then as in the same arguments with $\mathscr{H}om_{\mathcal{O}_M}\left( \Theta_{\mathcal{F}_0}, \frac{\Theta_M}{\Theta_{\mathcal{F}_0}} \right)$ in subsection \ref{ad1}, we have the complex
\begin{align*}
\frac{A^{0,0}\left(M, \mathscr{H}om_{\mathcal{O}_M}\left( \mathcal{N}_{\mathcal{F}_0}^*, \Omega_M^1  \right) \right)}{A^{0,0}\left( M, \mathscr{H}om_{\mathcal{O}_M}\left( \mathcal{N}_{\mathcal{F}_0}^*, \mathcal{N}_{\mathcal{F}_0}^* \right) \right)}\xrightarrow{\bar{\partial}}\frac{A^{0,1}\left(M, \mathscr{H}om_{\mathcal{O}_M}\left( \mathcal{N}_{\mathcal{F}_0}^*, \Omega_M^1  \right) \right)}{A^{0,1}\left( M, \mathscr{H}om_{\mathcal{O}_M}\left( \mathcal{N}_{\mathcal{F}_0}^*, \mathcal{N}_{\mathcal{F}_0}^* \right) \right)}\xrightarrow{\bar{\partial}}  \frac{A^{0,2}\left(M, \mathscr{H}om_{\mathcal{O}_M}\left( \mathcal{N}_{\mathcal{F}_0}^*, \Omega_M^1  \right) \right)}{A^{0,2}\left( M, \mathscr{H}om_{\mathcal{O}_M}\left( \mathcal{N}_{\mathcal{F}_0}^*, \mathcal{N}_{\mathcal{F}_0}^* \right) \right)}\xrightarrow{\bar{\partial}}  \cdots
\end{align*}

Then we have the following Dolbeault type bicomplex associated to $\mathcal{N}_{\mathcal{F}_0}^{*\bullet}$
{\small{\begin{equation}
\begin{CD}
\cdots  \\
@A\hat{E}_3AA   \\
A^{0,0}\left( M-S, \left(\mathcal{N}_{\mathcal{F}_0}^*\right)^*\otimes \tilde{\mathcal{S}}^3 \right)@>-(-1)^q \bar{\partial}>> \cdots \\
@A\hat{E}_2AA @A\hat{E}_2AA \\
A^{0,0}\left( M-S, \left(\mathcal{N}_{\mathcal{F}_0}^*\right)^*\otimes \tilde{\mathcal{S}}^2\right) @>(-1)^q\bar{\partial}>> A^{0,1}\left(M-S, \left(\mathcal{N}_{\mathcal{F}_0}^*\right)^*\otimes \tilde{\mathcal{S}}^2  \right) @>-(-1)^q\bar{\partial}>> \cdots  \\
@A\hat{E}_1AA @A \hat{E}_1AA @A \hat{E}_1AA \\
\frac{A^{0,0}\left(M, \left(\mathcal{N}_{\mathcal{F}_0}^* \right)^* \otimes \Omega_M^1 \right)}{A^{0,0}\left(M,\left(\mathcal{N}_{\mathcal{F}_0}^* \right)^*\otimes \mathcal{N}_{\mathcal{F}_0}^*\right)} @>-\bar{\partial}>>\frac{A^{0,1}\left(M, \left(\mathcal{N}_{\mathcal{F}_0}^*\right)^*\otimes \Omega_M^1 \right)}{A^{0,1}\left(M,\left(\mathcal{N}_{\mathcal{F}_0}^* \right)^*\otimes \mathcal{N}_{\mathcal{F}_0}^*\right)}@>\bar{\partial}>> \frac{A^{0,2}\left(M, \left(\mathcal{N}_{\mathcal{F}_0}^*\right)^*\otimes \Omega_M^1 \right)}{A^{0,2}\left(M,\left(\mathcal{N}_{\mathcal{F}_0}^* \right)^*\otimes \mathcal{N}_{\mathcal{F}_0}^*\right)}@>-\bar{\partial}>> \cdots \\
@A\hat{E}_0AA @A\hat{E}_0AA @A\hat{E}_0AA @A\hat{E}_0AA \\
A^{0,0}\left(M, \Theta_M\right) @>\bar{\partial}>> A^{0,1}\left(M, \Theta_M\right) @>-\bar{\partial}>> A^{0,2}\left(M, \Theta_M \right) @>\bar{\partial}>> A^{0,3}\left(M, \Theta_M\right) @>-\bar{\partial}>> \cdots
\end{CD}
\end{equation}}}
where $\hat{E}_1: \frac{A^{0,p}\left( M, \mathscr{H}om_{\mathcal{O}_M}\left( \mathcal{N}_{\mathcal{F}_0}^* ,   \Omega_M^1 \right) \right)}{A^{0,p}\left( M, \mathscr{H}om_{\mathcal{O}_M}\left( \mathcal{N}_{\mathcal{F}_0}^* ,  \mathcal{N}_{\mathcal{F}_0}^* \right) \right)} \to A^{0,p}\left( M-S, \mathscr{H}om_{\mathcal{O}_M}\left( \mathcal{N}_{\mathcal{F}_0}^* , \tilde{\mathcal{S}}^2 \right) \right)$ is defined in the following way: first we note that we have the composition $\mathscr{H}om_{\mathcal{O}_M}\left( \mathcal{N}_{\mathcal{F}_0}^*, \frac{\Omega_M^1}{\mathcal{N}_{\mathcal{F}_0}^*} \right)\xrightarrow{E_1} \mathscr{H}om_{\mathcal{O}_M}\left( \mathcal{N}_{\mathcal{F}_0}^*, \tilde{\mathcal{S}}^2\right) \hookrightarrow \mathscr{H}om_{\mathcal{O}_M}\left( \mathcal{N}_{\mathcal{F}_0}^*, \bigwedge^{q+2} \Omega_M^1\otimes \mathcal{L}_0 \right)$ where $\mathcal{L}_0:= \bigwedge^q \left(\mathcal{N}_{\mathcal{F}_0}^*\right)^*$ which induces 
\begin{align*}
E_1^\sharp:\frac{A^{0,p}\left( M, \mathscr{H}om_{\mathcal{O}_M}\left( \mathcal{N}_{\mathcal{F}_0}^* ,   \Omega_M^1 \right) \right)}{A^{0,p}\left( M, \mathscr{H}om_{\mathcal{O}_M}\left( \mathcal{N}_{\mathcal{F}_0}^* ,  \mathcal{N}_{\mathcal{F}_0}^* \right) \right)}  \to  A^{0,p}\left( M, \mathscr{H}om_{\mathcal{O}_M}\left( \mathcal{N}_{\mathcal{F}_0}^*, \bigwedge^{q+2} \Omega_M^1 \otimes \mathcal{L}_0 \right) \right)
\end{align*}
We note that when we restrict the image of $E_1^\sharp$ to $M-S$, the restriction lies in $A^{0,p}\left(M-S, \mathscr{H}om_{\mathcal{O}_M}\left( \mathcal{N}_{\mathcal{F}_0}^*, \tilde{\mathcal{S}}^2 \right) \right)\subset A^{0,p}\left(M-S, \mathscr{H}om_{\mathcal{O}_M}\left( \mathcal{N}_{\mathcal{F}_0}^* , \bigwedge^{q+2} \Omega_M^1 \otimes \mathcal{L}_0 \right) \right)$. Then we define $\hat{E}_1$ is the composition of $E_1^\sharp$ and the restriction on $M-S$.

We will denote the $i$-th cohomology group of the complex associated to the above bicomplex by $\textnormal{S}^i$. Then we will show in Part IV-1 that for $i=0,1,2$,
\begin{align}\label{d5}
\textnormal{S}^i \cong \mathbb{H}^i\left( M, \mathcal{N}_{\mathcal{F}_0}^{*\bullet} \right)
\end{align}
by explicitly constructing isomorphisms between $\textnormal{S}^i$ and $i$-th cohomolgoy group from the \v Cech resolution of $\mathcal{N}_{\mathcal{F}_0}^{*\bullet}$ for $i=0,1,2$.

We will describe $\hat{E}_0$ and $\hat{E}_1$ and $\hat{E}_2$ explicitly which we use in the main body of the paper. Let $\mathcal{U}=\left\{U_i\right\}$ be a Stein open covering of $M$ by coordinate neighborhoods with local coordinates $\left( z_i^1,..., z_i^n \right)$ on $U_i$ such that $\Gamma\left( U_j, \mathcal{N}_{\mathcal{F}_0}^* \right)$ is generated by $w_i^1,..., w_i^q\in \Gamma\left( U_j, \Omega_M^1 \right)$ with the relation $w_j^\alpha= \sum_{\beta=1}^q h_{jk}^{\alpha\beta} w_j^\beta$ for some $h_{jk}^{\alpha\beta}\in \Gamma\left( U_j \cap U_k, \mathcal{O}_M \right)$ and $w_j^1 \wedge \cdots \wedge w_j^q \wedge dw_i^\alpha=0$ for $\alpha= 1,...,q$.

\subsubsection{\textnormal{Description of $\hat{E}_0$}}\ 

We describe $\hat{E}_0: A^{0,p}\left( M, \Theta_M \right) \to \frac{A^{0,p}\left( M, \mathscr{H}om_{\mathcal{O}_M}\left( \mathcal{N}_{\mathcal{F}_0}^*, \Omega_M^1 \right) \right)}{A^{0,p}\left( M, \mathscr{H}om_{\mathcal{O}_M}\left( \mathcal{N}_{\mathcal{F}_0}^* , \mathcal{N}_{\mathcal{F}_0}^* \right) \right)}$.  For $\phi \in A^{0,p}\left( M, \Theta_M \right)$, we define an element $\eta_{ij}\in \Gamma\left(U_{ij}, \mathcal{A}^{0,p}\left( \mathscr{H}om_{\mathcal{O}_M}\left( \mathcal{N}_{\mathcal{F}_0}^*, \mathcal{N}_{\mathcal{F}_0}^* \right) \right) \right)$ by
\begin{align*}
\eta_{ij} : \Gamma\left( U_{ij}, \mathcal{N}_{\mathcal{F}_0}^* \right) &\to \Gamma\left(  U_{ij}, \mathcal{A}^{0,p}\left( \mathcal{N}_{\mathcal{F}_0}^* \right)\right) \\
 w_i^\alpha &\mapsto \mathcal{L}_{\phi}\left( w_i^\alpha \right) - \sum_{\beta=1}^q h_{ij}^{\alpha\beta} \mathcal{L}_\phi\left( w_j^\beta \right) =  \sum_{\beta=1}^q \left[ \phi, h_{ij}^{\alpha\beta} \right] \wedge w_j^\beta
\end{align*}
Then $\eta_{ij}- \eta_{ik} + \eta_{jk}=0$. This implies that there exists $\left\{ \eta_i \right\}$ with $\eta_i \in \Gamma\left( U_i, \mathscr{H}om_{\mathcal{O}_M} \left( \mathcal{N}_{\mathcal{F}_0}^*, \mathcal{N}_{\mathcal{F}_0}^*   \right) \right)$ such that $\eta_j- \eta_i=\eta_{ij}$, where $\eta_i:\Gamma\left(U_i, \mathcal{N}_{\mathcal{F}_0}^* \right)\to \Gamma\left( U_i , \mathcal{A}^{0,p}\left( \mathcal{N}_{\mathcal{F}_0}^* \right) \right), w_i^\alpha \mapsto V_i^\alpha$. Then we define $\hat{E}_{0i}' (\phi)\in \Gamma\left( U_i, \mathcal{A}^{0,p}\left( \mathscr{H}om_{\mathcal{O}_M}\left( \mathcal{N}_{\mathcal{F}_0}^*, \Omega_M^1 \right) \right) \right)$ by $\hat{E}_{0i}' (\phi)\left( w_i^\alpha \right)= \mathcal{L}_\phi\left( w_i^\alpha \right) + V_i^\alpha$. Then $\hat{E}_0'\left( \phi \right):=\left\{ \hat{D}_{0i}'(\phi) \right\} \in A^{0,p}\left( M, \mathscr{H}om_{\mathcal{O}_M}\left( \mathcal{N}_{\mathcal{F}_0}^*, \Omega_M^1 \right) \right)$. Then we define $\hat{E}_0(\phi)$ to be the image of $\hat{E}_0'(\phi)$ in $\frac{A^{0,p}\left( M, \mathscr{H}om_{\mathcal{O}_M}\left( \mathcal{N}_{\mathcal{F}_0}^*, \Omega_M^1 \right) \right)}{A^{0,p}\left(  M, \mathscr{H}om_{\mathcal{O}_M}\left( \mathcal{N}_{\mathcal{F}_0}^*, \mathcal{N}_{\mathcal{F}_0}^*   \right) \right)}$.

\subsubsection{\textnormal{Description of $\hat{E}_1$}}\

We describe $\hat{E}_1: \frac{A^{0,p}\left( M, \mathscr{H}om_{\mathcal{O}_M}\left( \mathcal{N}_{\mathcal{F}_0}^*, \Omega_M^1 \right) \right)}{A^{0,p}\left( M, \mathscr{H}om_{\mathcal{O}_M}\left( \mathcal{N}_{\mathcal{F}_0}^*, \mathcal{N}_{\mathcal{F}_0}^* \right) \right)} \to A^{0,p}\left( M-S, \mathscr{H}om_{\mathcal{O}_M}\left( \mathcal{N}_{\mathcal{F}_0}^*, \tilde{\mathcal{S}}^2 \right) \right) $. It is sufficient to describe $\hat{E}_1^\sharp: \frac{A^{0,p}\left(  M, \mathscr{H}om_{\mathcal{O}_M}\left( \mathcal{N}_{\mathcal{F}_0}^*, \Omega_M^1  \right) \right)}{A^{0,p}\left( M, \mathscr{H}om_{\mathcal{O}_M} \left( \mathcal{N}_{\mathcal{F}_0}^*, \mathcal{N}_{\mathcal{F}_0}^* \right) \right)} \to A^{0,p}\left( M, \mathscr{H}om_{\mathcal{O}_M}\left( \mathcal{N}_{\mathcal{F}_0}^* , \bigwedge^{q+2} \Omega_M^1 \otimes \mathcal{L}_0 \right) \right)$. For $\overline{\phi} \in \frac{A^{0,p}\left(M, \mathscr{H}om_{\mathcal{O}_M}\left( \mathcal{N}_{\mathcal{F}_0}^*, \Omega_M^1 \right) \right)}{A^{0,p}\left( M, \mathscr{H}om_{\mathcal{O}_M}\left( \mathcal{N}_{\mathcal{F}_0}^*, \mathcal{N}_{\mathcal{F}_0}^* \right) \right)}$, let $\phi\left(w_i^\alpha  \right)= \sum_{\gamma_1< \cdots < \gamma_r} d\bar{z}_i^{\gamma_1}\wedge \cdots d \bar{z}_i^{\gamma_p} \wedge A_{\gamma_1,..., \gamma_p}^\alpha$ on $U_i$ with $A_{\gamma_1,..., \gamma_p}^\alpha \in \Gamma\left( U_i,  \mathcal{A}^{0,0}\left(\Omega_M^1 \right) \right)$. Then we have
\begin{align*}
\hat{E}_1^\sharp\left( \phi \right)\left( w_i^\alpha \right)=  (-1)^{pq} \sum_{\gamma_1<\cdots <\gamma_p} d \bar{z}_i^{\gamma_1}\wedge \cdots \wedge d \bar{z}_i^{\gamma_p}\wedge\left( \sum_{\beta=1}^q w_i^1 \wedge \cdots \wedge A_{\gamma_1,..., \gamma_p}^\beta \wedge \cdots \wedge w_i^q \wedge dw_i^\alpha + w_i^1 \wedge \cdots \wedge w_i^q \wedge \partial\left( A_{\gamma_1,...,\gamma_p}^\alpha \right)        \right)
\end{align*}
In particular, for $p=1$, we have
\begin{align*}
\hat{E}_1\left( \overline{\phi} \right)(w_i^\alpha) &= (-1)^q \left( \sum_{\beta=1}^q (-1)^{\beta-1} w_i^1 \wedge \cdots  \wedge \phi\left( w_i^\beta \right) \wedge \cdots \wedge w_i^q \wedge d w_i^\alpha - (-1)^q w_i^1 \wedge \cdots \wedge w_i^q \wedge \partial\left( \phi\left( w_i^\alpha \right) \right) \right)
\end{align*}

\subsubsection{\textnormal{Description of $\hat{E}_2$}}\

We describe $\hat{E}_2: A^{0,p}\left( M-S, \mathscr{H}om_{\mathcal{O}_M}\left( \mathcal{N}_{\mathcal{F}_0}^*, \tilde{\mathcal{S}}^2 \right) \right) \to A^{0,p}\left(M-S, \mathscr{H}om_{\mathcal{O}_M}\left( \mathcal{N}_{\mathcal{F}_0}^*,  \tilde{\mathcal{S}}^3 \right) \right)$. For $x\in U_x\subset U_i-S$, we locally compute on $U_x$. We note that we have $dw_i^\alpha- \sum_{\beta=1}^q a_{i_x}^{\alpha\beta} \wedge w_i^\beta$ for some $a_{i_x}^{\alpha\beta}\in \Gamma\left( U_x , \Omega_M^1 \right)$.  For $C\in A^{0,p}\left(M-S, \mathscr{H}om_{\mathcal{O}_M} \left( \mathcal{N}_{\mathcal{F}_0}^*, \tilde{\mathcal{S}}^2 \right) \right)$, we can write $C\left( w_i^\alpha \right)= w_{i}^1 \wedge \cdots \wedge w_i^q \wedge C_{i_x}^\alpha $ on $U_x$, where $C_{i_x}^\alpha=  \sum_{\gamma_1< \cdots < \gamma_p} d\bar{z}_i^{\gamma_1}\wedge \cdots \wedge d \bar{z}_i^{\gamma_p} \wedge C_{\gamma_1,..., \gamma_p}^\alpha$ with $C_{\gamma_1,..,\gamma_p}^\alpha\in \Gamma\left( U_x, \mathcal{A}^{0,0}\left(\bigwedge^2 \Omega_M^1\right) \right)$. Then we define
\begin{align*}
\hat{E}_2\left(C \right)\left( w_i^\alpha \right) &= w_i^1 \wedge \cdots \wedge w_i^q \wedge \left(  \sum_{\gamma_1< \cdots < \gamma_p} d\bar{z}_i^{\gamma_1}\wedge \cdots \wedge d \bar{z}_i^{\gamma_p} \wedge \left( \partial \left( C_{\gamma_1,..., \gamma_p}^\alpha \right)- \sum_{\beta=1}^q a_{i_x}^{\alpha\beta} \wedge C_{\gamma_1,...,\gamma_p}^\beta  \right)  \right) \\
&=(-1)^p w_i^1 \wedge \cdots \wedge w_i^q \wedge \left( \partial \left( C_{i_x}^\alpha \right) - \sum_{\beta=1}^q a_{i_x}^{\alpha\beta} \wedge C_{i_x}^\beta     \right)
\end{align*}

\subsection{Dolbeault type bicomplex associated to $\mathcal{F}^\bullet$} \

Let $\left( M, \Theta_{\mathcal{F}_0}, \mathcal{N}_{\mathcal{F}_0}^* \right)$ be a compact  foliated complex manifold with both $\Theta_{\mathcal{F}_0}$ and $\mathcal{N}_{\mathcal{F}_0}^*$ locally free. Then as in a similar way to subsection \ref{ad1} and subsection \ref{ad2}. We have a Dolbeault type bicomplex associated to $\mathcal{F}_0^\bullet$ from $\textnormal{(\ref{ad3})}$.
{\Tiny{\begin{equation}
\begin{CD}
\cdots\\
@A\hat{F}_3AA \\
\frac{A^{0,0}\left(M, \bigwedge^3 \Theta_{\mathcal{F}_0}^*\otimes \Theta_M\right)}{A^{0,0}\left(M,\bigwedge^3 \Theta_{\mathcal{F}_0}^*\otimes \Theta_{\mathcal{F}_0}\right)}\bigoplus A^{0,0}\left(M - S, \left(\mathcal{N}_{\mathcal{F}_0}^*\right)^*\otimes \tilde{\mathcal{S}}^3 \right)\bigoplus A^{0,0}\left(M-S, \left(\mathcal{N}_{\mathcal{F}_0}^* \right)^*\otimes \bigwedge^2 \Theta_\mathcal{F}^* \right) @> - \left(\bar{\partial}, (-1)^q \bar{\partial}, \bar{\partial}  \right)>> \cdots \\
@A \hat{F}_2AA \\
\frac{A^{0,0}\left(M, \bigwedge^2 \Theta_{\mathcal{F}_0}^*\otimes \Theta_M\right)}{A^{0,0}\left(M,\bigwedge^2 \Theta_{\mathcal{F}_0}^*\otimes \Theta_{\mathcal{F}_0}\right)}\bigoplus A^{0,0}\left(M - S, \left(\mathcal{N}_{\mathcal{F}_0}^*\right)^*\otimes \tilde{\mathcal{S}}^2 \right)\bigoplus A^{0,0}\left(M, \left(\mathcal{N}_{\mathcal{F}_0}^* \right)^*\otimes \Theta_\mathcal{F}^* \right)@>\bar{\partial}, (-1)^q \bar{\partial}, \bar{\partial} >> \cdots \\
@A \hat{F}_1AA @A \hat{F}_1AA \\
\frac{A^{0,0}\left(M, \Theta_{\mathcal{F}_0}^*\otimes \Theta_M\right)}{A^{0,0}\left(M, \Theta_{\mathcal{F}_0}^*\otimes \Theta_{\mathcal{F}_0} \right)}\bigoplus \frac{A^{0,0}\left(M,\left(\mathcal{N}_{\mathcal{F}_0}^*\right)^*\otimes \Omega_M^1\right)}{A^{0,0}\left( M, \left(\mathcal{N}_{\mathcal{F}_0}^* \right)^* \otimes \mathcal{N}_{\mathcal{F}_0}^*\right)} @>-\bar{\partial}>>\frac{A^{0,1}\left(M, \Theta_{\mathcal{F}_0}^*\otimes \Theta_M\right)}{A^{0,1}\left(M, \Theta_{\mathcal{F}_0}^*\otimes \Theta_{\mathcal{F}_0} \right)}\bigoplus \frac{A^{0,1}\left(M,\left(\mathcal{N}_{\mathcal{F}_0}^*\right)^*\otimes \Omega_M^1\right)}{A^{0,1}\left( M, \left(\mathcal{N}_{\mathcal{F}_0}^* \right)^* \otimes \mathcal{N}_{\mathcal{F}_0}^*\right)}@>\bar{\partial}>>\cdots \\
@A\hat{F}_0 AA @A \hat{F}_0 AA  \\
A^{0,0}\left(M, \Theta_M \right) @>\bar{\partial}>> A^{0,1} \left(M, \Theta_M\right) @>-\bar{\partial}>> \cdots 
\end{CD}
\end{equation}}}
We will denote the $i$-th cohomology group of the complex associated to the bicomplex by $\textnormal{F}^i$. Then we will show in Part IV-1 that for $i=0,1,2$.
\begin{align}\label{d26}
\textnormal{F}^i \cong \mathbb{H}^i\left( M, \mathcal{F}_0^\bullet \right)
\end{align}
by explicitly constructing isomorphisms between $\textnormal{F}^i$ and $i$-th cohomology group from the \v Cech resolution of $\mathcal{F}_0^\bullet$ for $i=0,1,2$.

We will describe $\hat{F}_0$ and $\hat{F}_1$ and $\hat{F}_2$ explicitly which we use in the main body of the paper. Let $\mathcal{U}= \left\{ U_i \right\}$ be a Stein open covering of $M$ by coordinate neighborhoods with local coordinates $\left(z_i^1,..., z_i^n  \right)$ on $U_i$ such that $\Gamma\left(U_j, \Theta_{\mathcal{F}_0} \right)$ is generated by $T_j^1,..., T_j^p\in \Gamma\left( U_j, \Theta_M \right)$ with the relation $T_j^\alpha = \sum_{\beta=1}^q r_{jk}^{\alpha\beta} T_j^\beta$ for some $r_{jk}^{\alpha\beta}\in \Gamma\left( U_j \cap U_k, \mathcal{O}_M \right)$ and $\left[ T_j^\alpha, T_i^\beta \right]=\sum_{\gamma=1}^p g_{j\alpha\beta}^\gamma T_i^\gamma$ for some $g_{j\alpha\beta}^\gamma \in \Gamma\left( U_j , \mathcal{O}_M \right)$. On the other hand, $\Gamma\left( U_j, \mathcal{N}_{\mathcal{F}_0}^* \right)$ is generated by $w_j^1,..., w_j^q \in \Gamma\left( U_j, \Omega_M^1 \right)$ with the relation $w_j^\alpha = \sum_{\beta=1}^q h_{jk}^{\alpha\beta} w_j^\beta \in \Gamma\left( U_j \cap U_j, \mathcal{O}_M \right)$ and $w_j^1\wedge \cdots \wedge w_j^q \wedge d w_j^\alpha=0$ for $\alpha=1,...,q$.

\subsubsection{\textnormal{Description of $\hat{F}_0$} } \

We describe $\hat{F}_0: A^{0,p}\left( M, \Theta_M \right) \to \frac{A^{0,p}\left( M, \mathscr{H}om_{\mathcal{O}_M} \left( \Theta_{\mathcal{F}_0}, \Theta_M \right) \right)}{ A^{0,p}\left( M, \mathscr{H}om_{\mathcal{O}_M} \left( \Theta_{\mathcal{F}_0}, \Theta_{\mathcal{F}_0}  \right)  \right)} \bigoplus \frac{A^{0,0}\left( M, \mathscr{H}om_{\mathcal{O}_M}\left( \mathcal{N}_{\mathcal{F}_0}^*, \Omega_M^1 \right) \right)}{A^{0,0}\left( M, \mathscr{H}om_{\mathcal{O}_M}\left( \mathcal{N}_{\mathcal{F}_0}^*, \mathcal{N}_{\mathcal{F}_0}^* \right) \right)} $. For $\phi\in A^{0,0}\left( M, \Theta_M \right)$, we define $\hat{F}_0(\phi):=\left( \hat{D}_0(\phi), \hat{E}_0 (\phi)  \right)$.

\subsubsection{\textnormal{Description of $\hat{F}_1$} } \

We describe
{\small{\begin{align*}
& \frac{A^{0,p}\left( M, \mathscr{H}om_{\mathcal{O}_M} \left( \Theta_{\mathcal{F}_0}, \Theta_M \right) \right)}{ A^{0,p}\left( M, \mathscr{H}om_{\mathcal{O}_M} \left( \Theta_{\mathcal{F}_0}, \Theta_{\mathcal{F}_0}  \right)  \right)} \bigoplus \frac{A^{0,p}\left( M, \mathscr{H}om_{\mathcal{O}_M}\left( \mathcal{N}_{\mathcal{F}_0}^*, \Omega_M^1 \right) \right)}{A^{0,p}\left( M, \mathscr{H}om_{\mathcal{O}_M}\left( \mathcal{N}_{\mathcal{F}_0}^*, \mathcal{N}_{\mathcal{F}_0}^* \right) \right)} \\
& \xrightarrow{\hat{F}_1}  \frac{A^{0,p}\left( M, \mathscr{H}om_{\mathcal{O}_M} \left( \bigwedge^2  \Theta_{\mathcal{F}_0}, \Theta_M \right) \right)}{ A^{0,p}\left( M, \mathscr{H}om_{\mathcal{O}_M} \left( \bigwedge^2 \Theta_{\mathcal{F}_0}, \Theta_{\mathcal{F}_0}  \right)  \right)} \bigoplus A^{0,p}\left(M-S, \mathscr{H}om_{\mathcal{O}_M}\left( \mathcal{N}_{\mathcal{F}_0}^*, \tilde{\mathcal{S}}^2 \right) \right) \bigoplus A^{0,p}\left( M, \mathscr{H}om_{\mathcal{O}_M}\left( \mathcal{N}_{\mathcal{F}_0}^*, \Theta_{\mathcal{F}_0}^* \right) \right)
\end{align*}}}
We define $\hat{F}_1\left(\overline{\phi}, \overline{\psi} \right)=\left( \hat{D}_1\left(\overline{\phi} \right), \hat{E}_1\left(\overline{\psi} \right) , \hat{E}_0^* \left(\overline{\phi}, \overline{\psi} \right) \right)$ where $\hat{E}_0^*\left( \overline{\phi}, \overline{\psi} \right)$ is defined in the following way:  Locally on $U_i$, we can write $\phi\left( T_i^\alpha \right) =\sum_{\gamma_1<\cdots <\gamma_p} d\bar{z}_i^{\gamma_1}\wedge \cdots d\bar{z}_i^{\gamma_p} A_{\gamma_1,...,\gamma_p}^\alpha$ where $A_{\gamma_1,...,\gamma_p}^\alpha\in \Gamma\left( U_i, \mathcal{A}^{0,0}\left(\Theta_M \right) \right)$ and we can write $\psi\left(w_i^\beta \right)=\sum_{\gamma_1<\cdots < \gamma_p} d\bar{z}_i^{\gamma_1} \wedge \cdots \wedge d \bar{z}_i^{\gamma_p} C_{\gamma_1,...,\gamma_p}^\beta$, where $C_{\gamma_1,...,\gamma_p}^\beta\in \Gamma\left( U_i, \mathcal{A}^{0,0}\left( \Omega_M^1 \right) \right)$. Then we define
\begin{align*}
\hat{E}_0^*\left( \overline{\phi}, \overline{\psi}\right)\left( T_i^\alpha \right)\left( w_i^\beta\right)= \sum_{\gamma_1<\cdots < \gamma_p} d\bar{z}_i^{\gamma_1} \wedge \cdots \wedge d \bar{z}_i^{\gamma_p} \left(i_{T_i^\alpha}\left( C_{\gamma_1,...,\gamma_p}^\beta \right)  + i_{A_{\gamma_1,...,\gamma_p}^\alpha}\left(  w_i^\beta \right) \right)
\end{align*}

\subsubsection{\textnormal{Description of $\hat{F}_2$}} \

We describe 
{\small{\begin{align*}
&  \frac{A^{0,p}\left( M, \mathscr{H}om_{\mathcal{O}_M} \left( \bigwedge^2  \Theta_{\mathcal{F}_0}, \Theta_M \right) \right)}{ A^{0,p}\left( M, \mathscr{H}om_{\mathcal{O}_M} \left( \bigwedge^2 \Theta_{\mathcal{F}_0}, \Theta_{\mathcal{F}_0}  \right)  \right)} \bigoplus A^{0,p}\left(M-S, \mathscr{H}om_{\mathcal{O}_M}\left( \mathcal{N}_{\mathcal{F}_0}^*, \tilde{\mathcal{S}}^2 \right) \right) \bigoplus A^{0,p}\left( M, \mathscr{H}om_{\mathcal{O}_M}\left( \mathcal{N}_{\mathcal{F}_0}^*, \Theta_{\mathcal{F}_0}^* \right) \right)\\
&\xrightarrow{\hat{F}_2}  \frac{A^{0,p}\left( M, \mathscr{H}om_{\mathcal{O}_M} \left( \bigwedge^3  \Theta_{\mathcal{F}_0}, \Theta_M \right) \right)}{ A^{0,p}\left( M, \mathscr{H}om_{\mathcal{O}_M} \left( \bigwedge^3 \Theta_{\mathcal{F}_0}, \Theta_{\mathcal{F}_0}  \right)  \right)} \bigoplus A^{0,p}\left(M-S, \mathscr{H}om_{\mathcal{O}_M}\left( \mathcal{N}_{\mathcal{F}_0}^*, \tilde{\mathcal{S}}^3 \right) \right) \bigoplus A^{0,p}\left( M-S, \mathscr{H}om_{\mathcal{O}_M}\left( \mathcal{N}_{\mathcal{F}_0}^*, \bigwedge^2 \Theta_{\mathcal{F}_0}^* \right) \right)
\end{align*}}}
We define $\hat{F}_2\left(\overline{B}, C, M \right)=\left(\hat{D}_2\left( \overline{B} \right), \hat{E}_2\left(C \right), \hat{E}_1^*\left(\overline{B}, C, M \right) \right)$ where $\hat{E}_1^*$ is defined in the following way: first we note that we have an isomorphism
\begin{align*}
A^{0,p}\left( M-S, \mathscr{H}om_{\mathcal{O}_M}\left( \mathcal{N}_{\mathcal{F}_0}^*, \tilde{\mathcal{S}}^2 \right) \right) \xrightarrow{\alpha_2} A^{0,p}\left( M-S, \mathscr{H}om_{\mathcal{O}_M}\left( \mathcal{N}_{\mathcal{F}_0}^*,  \bigwedge^2 \Theta_{\mathcal{F}_0}^* \right) \right)
\end{align*} 
which is defined in the following way: for $C\in A^{0,p}\left(M-S, \mathscr{H}om_{\mathcal{O}_M}\left( \mathcal{N}_{\mathcal{F}_0}^*, \tilde{\mathcal{S}}^2 \right) \right)$, we choose $x\in U_x \subset U_i-S$ such that $C\left( w_i^\gamma \right)= w_i^1 \wedge \cdots \wedge w_i^q \wedge C_{i_x}^\gamma $ where $C_{i_x}^\gamma \in \Gamma\left(U_x, \mathcal{A}^{0,p} \left( \bigwedge^2 \Omega_M^1 \right) \right)$. Then we have
\begin{align*}
i_{T_i^\alpha \wedge T_i^\beta}\left( C\left(w_i^\gamma\right) \right) = i_{T_i^\alpha \wedge T_i^\beta}\left( C_{i_x}^\gamma \right)\cdot w_i^1 \wedge \cdots \wedge w_i^q
\end{align*}
For $y\in U_i-S$ such that $U_x\cap U_y \ne \emptyset$, we have $i_{T_i^\alpha\wedge T_i^\beta}\left( C_{i_x}^{\gamma} \right)= i_{T_i^\alpha \wedge T_i^\beta} \left( C_{i_y}^\gamma \right)$. This implies that $\left\{ i_{T_i^\alpha \wedge T_i^\beta}\left( C_{i_x}^\gamma \right) \right\}$ glues together to define $f_{i\alpha\beta}^\gamma \in \Gamma\left(U_i-S, \mathcal{A}^{0,p} \right)$. Then $\alpha_2\left(C \right) \in A^{0,p}\left( M-S, \mathscr{H}om_{\mathcal{O}_M}\left( \mathcal{N}_{\mathcal{F}_0}^*, \bigwedge^2 \Theta_{\mathcal{F}_0}^* \right) \right)$ is defined locally on $U_i-S$ by $\alpha_2\left(C\right)\left(w_i^\gamma \right)\left( T_i^\alpha\wedge T_i^\beta \right)= f_{i\alpha\beta}^\gamma$. On the other hand, locally on $U_i$, we have $B\left( T_i^\alpha \wedge T_i^\beta \right) \in \Gamma\left( U_i, \mathcal{A}^{0,p} \left( \Theta_M \right) \right)$. Let $\hat{E}_0'(M)\in A^{0,p}\left( M, \mathscr{H}om_{\mathcal{O}_M}\left( \mathcal{N}_{\mathcal{F}_0}^*, \bigwedge^2 \Theta_{\mathcal{F}_0}^* \right) \right)$ be defined locally on $U_i$ in the following way: let $M\left( T_i^\alpha \right)\left( w_i^\gamma \right)= \sum_{\gamma_1 <\cdots <\gamma_p} d \bar{z}_i^{\gamma_1}\wedge \cdots \wedge d\bar{z}_i^{\gamma_p}  M_{\gamma_1,..., \gamma_p}^{\alpha \gamma} $ on $U_i$ where $M_{\gamma_1,...,\gamma_p}^{\alpha \gamma}\in \Gamma\left( U_i, \mathcal{A}^{0,0} \right)$. Then
{\tiny{\begin{align*}
\hat{E}_0'(M)\left( T_i^\alpha \wedge T_i^\beta \right)\left(  w_i^\gamma \right) = \sum_{\gamma_1 < \cdots < \gamma_p} d\bar{z}_i^{\gamma_1} \wedge \cdots \wedge d \bar{z}_i^{\gamma_p} \left( \left[ T_i^\alpha, M_{\gamma_1,...,\gamma_p}^{\beta \gamma} \right]  - \sum_{\delta=1}^q a_{i\alpha}^{\gamma \delta} M_{\gamma_1,..., \gamma_p}^{\beta \delta} - \left[ T_i^\beta , M_{\gamma_1,..., \gamma_p}^{\alpha \gamma}\right] + \sum_{\delta=1}^q a_{i\beta}^{\gamma \delta} M_{\gamma_1,..., \gamma_p}^{\alpha \delta}  - \sum_{\eta=1}^p g_{i\alpha\beta}^\eta M_{\gamma_1,..., \gamma_p}^{\eta \gamma}    \right)
\end{align*}}}
where $\mathcal{L}_{T_i^\alpha}\left( w_i^\gamma \right)= \sum_{\delta=1}^q a_{i\alpha}^{\gamma \delta} w_i^\delta$.

Then we define $\hat{E}_1^*\left( \overline{B}, C, M \right)$ locally on $U_i-S$ by
\begin{align*}
\hat{E}_1^* \left( \overline{B}, C, M \right)\left( T_i^\alpha \wedge T_i^\beta \right)\left( w_i^\gamma \right) = i_{B\left(T_i^\alpha \wedge T_i^\beta \right) }\left( w_i^\gamma \right) + \alpha_2 (C)\left( T_i^\alpha \wedge T_i^\beta \right) \left( w_i^\gamma \right)  - \hat{E}_0' (M)\left( T_i^\alpha \wedge T_i^\beta \right)\left( w_i^\gamma \right)
\end{align*}

\subsection{Dolbeault resolution of $\mathcal{F}_0'^\bullet$} \

Let $\left(M, \Theta_{\mathcal{F}_0}, \mathcal{N}_{\mathcal{F}_0}^* \right)$ be a compact foliated complex manifold with both $\Theta_{\mathcal{F}_0}$ and $\mathcal{N}_{\mathcal{F}_0}^*$ locally free. Then as in the a similar way to subsection \ref{ad1}, we have the Dolbeualt resolution of $\mathcal{F}_0'^\bullet$ from $\textnormal{(\ref{ad10})}$.
{\Tiny{\begin{center}
$\begin{CD}
\cdots \\
@A \hat{F}_2'AA \\
\frac{A^{0,0}\left(M, \bigwedge^2 \Theta_{\mathcal{F}_0}^*\otimes \Theta_M\right)}{A^{0,0}\left(M,\bigwedge^2 \Theta_{\mathcal{F}_0}^*\otimes \Theta_{\mathcal{F}_0}\right)}\bigoplus A^{0,0}\left(M , \left(\mathcal{N}_{\mathcal{F}_0}^*\right)^*\otimes \bigwedge^2 \Theta_{\mathcal{F}_0}^* \right)\bigoplus A^{0,0}\left(M, \left(\mathcal{N}_{\mathcal{F}_0}^* \right)^*\otimes \Theta_\mathcal{F}^* \right)@>\bar{\partial}>> \cdots \\
@A \hat{F}_1'AA @A \hat{F}_1'AA \\
\frac{A^{0,0}\left(M, \Theta_{\mathcal{F}_0}^*\otimes \Theta_M\right)}{A^{0,0}\left(M, \Theta_{\mathcal{F}_0}^*\otimes \Theta_{\mathcal{F}_0} \right)}\bigoplus \frac{A^{0,0}\left(M,\left(\mathcal{N}_{\mathcal{F}_0}^*\right)^*\otimes \Omega_M^1\right)}{A^{0,0}\left( M, \left(\mathcal{N}_{\mathcal{F}_0}^* \right)^* \otimes \mathcal{N}_{\mathcal{F}_0}^*\right)} @>-\bar{\partial}>>\frac{A^{0,1}\left(M, \Theta_{\mathcal{F}_0}^*\otimes \Theta_M\right)}{A^{0,1}\left(M, \Theta_{\mathcal{F}_0}^*\otimes \Theta_{\mathcal{F}_0} \right)}\bigoplus \frac{A^{0,1}\left(M,\left(\mathcal{N}_{\mathcal{F}_0}^*\right)^*\otimes \Omega_M^1\right)}{A^{0,1}\left( M, \left(\mathcal{N}_{\mathcal{F}_0}^* \right)^* \otimes \mathcal{N}_{\mathcal{F}_0}^*\right)}@>\bar{\partial}>>\cdots \\
@A\hat{F}_0 AA @A \hat{F}_0 AA @A \hat{F}_0 AA \\
A^{0,0}\left(M, \Theta_M \right) @>\bar{\partial}>> A^{0,1} \left(M, \Theta_M\right) @>-\bar{\partial}>> A^{0,2}\left(\Theta_M\right)  
\end{CD}$
\end{center}}}

We will explicitly describe $\hat{F}_0$ and $\hat{F}_1'$ and $\hat{F}_2'$ that we use in the main body of the paper. Let $\mathcal{U}= \left\{ U_i \right\}$ be a Stein open covering of $M$ by coordinate neighborhoods with local coordinates $\left(z_i^1,..., z_i^n  \right)$ on $U_i$ such that $\Gamma\left(U_j, \Theta_{\mathcal{F}_0} \right)$ is generated by $T_j^1,..., T_j^p\in \Gamma\left( U_j, \Theta_M \right)$ with the relation $T_j^\alpha = \sum_{\beta=1}^q r_{jk}^{\alpha\beta} T_j^\beta$ for some $r_{jk}^{\alpha\beta}\in \Gamma\left( U_j \cap U_k, \mathcal{O}_M \right)$ and $\left[ T_j^\alpha, T_i^\beta \right]=\sum_{\gamma=1}^p g_{j\alpha\beta}^\gamma T_i^\gamma$ for some $g_{j\alpha\beta}^\gamma \in \Gamma\left( U_j , \mathcal{O}_M \right)$. On the other hand, $\Gamma\left( U_j, \mathcal{N}_{\mathcal{F}_0}^* \right)$ is generated by $w_j^1,..., w_j^q \in \Gamma\left( U_j, \Omega_M^1 \right)$ with the relation $w_j^\alpha = \sum_{\beta=1}^q h_{jk}^{\alpha\beta} w_j^\beta \in \Gamma\left( U_j \cap U_j, \mathcal{O}_M \right)$ and $w_j^1\wedge \cdots \wedge w_j^q \wedge d w_j^\alpha=0$ for $\alpha=1,...,q$.

\subsubsection{\textnormal{Description of $\hat{F}_0$} } \

We describe $\hat{F}_0: A^{0,p}\left( M, \Theta_M \right) \to \frac{A^{0,p}\left( M, \mathscr{H}om_{\mathcal{O}_M} \left( \Theta_{\mathcal{F}_0}, \Theta_M \right) \right)}{ A^{0,p}\left( M, \mathscr{H}om_{\mathcal{O}_M} \left( \Theta_{\mathcal{F}_0}, \Theta_{\mathcal{F}_0}  \right)  \right)} \bigoplus \frac{A^{0,0}\left( M, \mathscr{H}om_{\mathcal{O}_M}\left( \mathcal{N}_{\mathcal{F}_0}^*, \Omega_M^1 \right) \right)}{A^{0,0}\left( M, \mathscr{H}om_{\mathcal{O}_M}\left( \mathcal{N}_{\mathcal{F}_0}^*, \mathcal{N}_{\mathcal{F}_0}^* \right) \right)} $. For $\phi\in A^{0,0}\left( M, \Theta_M \right)$, we define $\hat{F}_0(\phi):=\left( \hat{D}_0(\phi), \hat{E}_0 (\phi)  \right)$.

\subsubsection{\textnormal{Description of $\hat{F}_1'$}}\

We describe
{\small{\begin{align*}
& \frac{A^{0,p}\left( M, \mathscr{H}om_{\mathcal{O}_M} \left( \Theta_{\mathcal{F}_0}, \Theta_M \right) \right)}{ A^{0,p}\left( M, \mathscr{H}om_{\mathcal{O}_M} \left( \Theta_{\mathcal{F}_0}, \Theta_{\mathcal{F}_0}  \right)  \right)} \bigoplus \frac{A^{0,p}\left( M, \mathscr{H}om_{\mathcal{O}_M}\left( \mathcal{N}_{\mathcal{F}_0}^*, \Omega_M^1 \right) \right)}{A^{0,p}\left( M, \mathscr{H}om_{\mathcal{O}_M}\left( \mathcal{N}_{\mathcal{F}_0}^*, \mathcal{N}_{\mathcal{F}_0}^* \right) \right)} \\
& \xrightarrow{\hat{F}_1'}  \frac{A^{0,p}\left( M, \mathscr{H}om_{\mathcal{O}_M} \left( \bigwedge^2  \Theta_{\mathcal{F}_0}, \Theta_M \right) \right)}{ A^{0,p}\left( M, \mathscr{H}om_{\mathcal{O}_M} \left( \bigwedge^2 \Theta_{\mathcal{F}_0}, \Theta_{\mathcal{F}_0}  \right)  \right)} \bigoplus A^{0,p}\left( M , \mathscr{H}om_{\mathcal{O}_M}\left( \mathcal{N}_{\mathcal{F}_0}^*, \bigwedge^2 \Theta_{\mathcal{F}_0}^* \right) \right) \bigoplus A^{0,p}\left( M, \mathscr{H}om_{\mathcal{O}_M}\left( \mathcal{N}_{\mathcal{F}_0}^*, \Theta_{\mathcal{F}_0}^* \right) \right)
\end{align*}}}
We define $\hat{F}_2'\left(\overline{\phi},  \overline{\psi} \right)=\left(\hat{D}_1\left(\overline{\phi} \right) , \hat{E}_1'\left( \overline{\psi} \right) , \hat{E}_0^* \left(\overline{\phi}, \overline{\psi} \right) \right)$, where $\hat{E}_1'\left(\overline{\psi} \right)$ is defined in the following way: locally on $U_i$, we can write $\psi\left( w_i^\beta \right)=\sum_{\gamma_1< \cdots < \gamma_p} d \bar{z}_i^{\gamma_1} \wedge \cdots \wedge d \bar{z}_i^{\gamma_p}  \wedge C_{\gamma_1,..., \gamma_p}^\beta$, where $C_{\gamma_1,..., \gamma_p}^\beta \in \Gamma \left(U_i, \mathcal{A}^{0,0}\left( \Omega_M^1 \right) \right)$. Then we define
{\Tiny{\begin{align*}
&\hat{E}_1'\left( \overline{\psi} \right)\left( T_i^\alpha \wedge T_i^\beta \right)\left( w_i^\gamma \right)\\
& = \sum_{\gamma_1< \cdots < \gamma_p} d \bar{z}_i^{\gamma_1} \wedge \cdots \wedge d \bar{z}_i^{\gamma_p} \left( \left[ T_i^\alpha, i_{T_i^\beta}\left( C_{\gamma_1,..., \gamma_p}^\gamma \right) \right] - \sum_{\delta=1}^q a_{i\alpha}^{\gamma \delta}  i_{T_i^\beta}\left( C_{\gamma_1,..., \gamma_p}^\delta \right) - \left[ T_i^\beta,  i_{T_i^\alpha}\left( C_{\gamma_1,..., \gamma_p}^\gamma \right) \right]  + \sum_{\delta=1}^q a_{i\beta}^{\gamma \delta} i_{T_i^\alpha}\left( C_{\gamma_1,..., \gamma_p}^{\delta}  \right) - \sum_{\eta=1}^p g_{i\alpha\beta}^{\eta} i_{T_i^\eta}\left( C_{\gamma_1,..., \gamma_p}^\gamma \right) \right)
\end{align*}}}

\subsubsection{\textnormal{Description of $\hat{F}_2'$}} \

We describe 
{\small{\begin{align*}
&  \frac{A^{0,p}\left( M, \mathscr{H}om_{\mathcal{O}_M} \left( \bigwedge^2  \Theta_{\mathcal{F}_0}, \Theta_M \right) \right)}{ A^{0,p}\left( M, \mathscr{H}om_{\mathcal{O}_M} \left( \bigwedge^2 \Theta_{\mathcal{F}_0}, \Theta_{\mathcal{F}_0}  \right)  \right)} \bigoplus A^{0,p}\left(M, \mathscr{H}om_{\mathcal{O}_M}\left( \mathcal{N}_{\mathcal{F}_0}^*, \bigwedge^2 \Theta_{\mathcal{F}_0}^* \right) \right) \bigoplus A^{0,p}\left( M, \mathscr{H}om_{\mathcal{O}_M}\left( \mathcal{N}_{\mathcal{F}_0}^*, \Theta_{\mathcal{F}_0}^* \right) \right)\\
&\xrightarrow{\hat{F}_2}  \frac{A^{0,p}\left( M, \mathscr{H}om_{\mathcal{O}_M} \left( \bigwedge^3  \Theta_{\mathcal{F}_0}, \Theta_M \right) \right)}{ A^{0,p}\left( M, \mathscr{H}om_{\mathcal{O}_M} \left( \bigwedge^3 \Theta_{\mathcal{F}_0}, \Theta_{\mathcal{F}_0}  \right)  \right)} \bigoplus A^{0,p}\left(M , \mathscr{H}om_{\mathcal{O}_M}\left( \mathcal{N}_{\mathcal{F}_0}^*, \bigwedge^3 \Theta_{\mathcal{F}_0}^* \right) \right) \bigoplus A^{0,p}\left( M, \mathscr{H}om_{\mathcal{O}_M}\left( \mathcal{N}_{\mathcal{F}_0}^*, \bigwedge^2 \Theta_{\mathcal{F}_0}^* \right) \right)
\end{align*}}}
We define $\hat{F}_2\left( \overline{B}, C, M \right) =  \left( \hat{D}_2\left( \overline{B} \right) , \hat{E}_2' \left( C \right),  \hat{E}_1^{**}\left(\overline{B} , C, M \right)\right)$, where $\hat{E}_1^{**}\left( \overline{B}, C, M \right)$ is defined locally on $U_i$:
\begin{align*}
\hat{E}_1^{**}\left(\overline{B}, C, M \right) \left( T_i^\alpha \wedge T_i^\beta \right) \left( w_i^\gamma  \right) = i_{B\left( T_i^\alpha \wedge T_i^\beta \right)}\left( w_i^\gamma \right) + C\left( T_i^\alpha \wedge T_i^\beta \right) \left( w_i^\gamma \right) - \hat{E}_1'(M)\left(T_i^\alpha \wedge T_i^\beta \right) \left( w_i^\gamma \right)
\end{align*} 
On the other hand, $\hat{E}_2'(C)$ is defined locally in the following way: locally on $U_i$, 
\begin{align*}
C \left( T_i^\alpha \wedge T_i^\beta \right)\left( w_i^\gamma \right)= \sum_{\gamma_1< \cdots <\gamma _p} d\bar{z}_i^{\gamma_1} \wedge \cdots \wedge d \bar{z}_i^{\gamma_p}  C_{\gamma_1,..., \gamma_p}^{\alpha\beta \gamma}
\end{align*}
Then we define locally on $U_i$
\begin{align*}
\hat{E}_2'\left( C \right)\left( T_i^\alpha \wedge T_i^\beta \wedge T_i^\gamma \right)\left( w_i^\delta \right)= \sum_{\gamma_1< \cdots <\gamma_p} d\bar{z}_i^{\gamma_1} \wedge \cdots \wedge d\bar{z}_i^{\gamma_p} N_{\gamma_1,...,\gamma_p}^{\alpha\beta\gamma \delta}
\end{align*}
where
\begin{align*}
N_{\gamma_1,.., \gamma_p}^{\alpha\beta \gamma \delta}:&= \left[ T_i^\alpha, C_{\gamma_1,..,\gamma_p}^{\beta \gamma \delta}     \right] - \sum_{\eta=1}^q a_{i\alpha}^{\delta \eta} C_{\gamma_1,..., \gamma_p}^{\beta \gamma \eta} - \left[ T_i^\beta, C_{\gamma_1,..., \gamma_p}^{\alpha \gamma \delta} \right]+ \sum_{\eta=1}^q a_{i\beta}^{\delta \eta} C_{\gamma_1,..., \gamma_p}^{\alpha\gamma \eta} + \left[ T_i^\gamma, C_{\gamma_1,..., \gamma_p}^{\alpha\beta \delta} \right] - \sum_{\eta=1}^q a_{i\gamma}^{\delta \eta} C_{\gamma_1,..., \gamma_p}^{\alpha \beta \eta}\\
& - \sum_{\xi=1}^p g_{i\alpha\beta}^\xi C_{\gamma_1,..., \gamma_p}^{\xi \gamma \delta} + \sum_{\xi=1}^p g_{i\alpha \gamma}^\xi C_{\gamma_1,..., \gamma_p}^{\xi \beta \delta} - \sum_{\xi=1}^p g_{i\beta \gamma}^\xi C_{\gamma_1,..., \gamma_p}^{\xi \alpha \delta}
\end{align*}

\bibliographystyle{amsalpha}
\bibliography{References-Rev9}

\end{document}